%% file: 1-19.tex
\documentclass{article}
\usepackage{amsmath,amssymb,amscd}

\input{notation}
\input{new_theorem}

\begin{document}

\title{Rescalability of integrable mixed twistor $D$-modules}
\author{Takuro Mochizuki}
\date{}
\maketitle

\begin{abstract}
We study the rescalability of 
integrable mixed twistor $D$-modules.
We prove some basic functoriality of
the rescalability and the associated irregular Hodge filtrations.
We also observe that rescalable integrable mixed twistor $D$-modules
are equivalent to exponential Hodge modules. 

\end{abstract}

\section{Introduction}

\subsection{Background}

\subsubsection{Mixed Hodge modules and mixed twistor $\nbigd$-modules}

A branch of the Hodge theory is
the study of the functoriality of Hodge structures.
It is achieved by the theory of mixed Hodge modules
due to Morihiko Saito \cite{saito1,saito2}.
A Hodge structure on a regular holonomic $\nbigd$-module
is formulated as a mixed Hodge module,
and it is established that
the Hodge structure is functorial with respect to
the standard operations for $\nbigd$-modules.

In \cite{s3}, Simpson introduced the notion of mixed twistor structure
as a generalization of mixed Hodge structure,
that is a holomorphic vector bundle on $\proj^1$.
He proposed a principle called Simpson's Meta Theorem,
which roughly says,
a theory concerning Hodge structures
should be generalized to a theory concerning twistor structures.
As a twistor version of Hodge modules,
the theory of twistor $\nbigd$-modules was developed
by Sabbah \cite{Sabbah-pure-twistor, Sabbah-wild-twistor}
and the author \cite{mochi2, Mochizuki-wild, Mochizuki-MTM}.
A mixed twistor structure on a holonomic $\nbigd$-module
is formulated as a mixed twistor $\nbigd$-module,
and it is established that
the twistor structure is functorial with respect to
the standard operations for $\nbigd$-modules.

\subsubsection{Irregular Hodge theory}

The Hodge numbers are associated with a mixed Hodge structure
as the dimensions of the graded pieces
of the Hodge filtration and the weight filtration.
They are important numerical invariants
which are useful for applications.
In the theory of mixed twistor $\nbigd$-modules,
there does not exist the numerical invariant
corresponding to the Hodge numbers, a priori.
To improve the situation, Sabbah pursued
the idea of irregular Hodge theory due to Deligne,
by partially collaborating with Yu.

Deligne \cite{Deligne-Malgrange-Ramis}
introduced a filtration on the $\nbigd$-module $(\nbigo,d+df)$
associated with 
an algebraic function $f$ on a complex algebraic curve,
and he proved the $E_1$-degeneration
of the spectral sequence of the associated filtered de Rham complex,
which defines a filtration
on the twisted de Rham cohomology.
The filtrations are now called the irregular Hodge filtrations.
In \cite{Sabbah-Fourier-LaplaceII},
Sabbah generalized the construction of Deligne
to the case of
a Hodge module twisted by an algebraic function
on an algebraic curve,
and established the desired $E_1$-degeneration property.
In the proof, he related the irregular Hodge filtration
with the $V$-filtration of
the rescaling of integrable twistor $\nbigd$-modules.
As noted in \cite{Sabbah-irregular-Hodge},
the rescaling procedure in the context of
irregular Hodge theory goes back to
the study of Hertling and Sevenheck
\cite{Hertling, Hertling-Sevenheck, Hertling-Sevenheck2}.
Independently, Yu \cite{Yu} generalized the irregular Hodge filtration
to the case of twisted de Rham complex
associated with any algebraic function on any quasi-projective variety.
He conjectured the $E_1$-degeneration of the associated spectral sequence,
and he proved it in some interesting cases.
The filtered twisted de Rham complex of Yu
is related with complexes introduced by Kontsevich.
Esnault, Sabbah and Yu \cite{Esnault-Sabbah-Yu},
established the $E_1$-degeneration property of the complexes.
Sabbah and Yu \cite{Sabbah-Yu}
introduced the irregular Hodge filtration
for a mixed Hodge module twisted by an algebraic function,
and studied some further properties of Kontsevich complexes.
Later, the author \cite{Mochizuki-Kontsevich-complexes}
studied the Kontsevich complexes and the associated irregular Hodge filtration
on the basis of the theory of mixed twistor $\nbigd$-modules.

In \cite{Sabbah-irregular-Hodge},
Sabbah introduced the notion of
``irregular mixed Hodge modules''
as a general framework
to study the irregular Hodge filtrations.
In particular,
the meaning of an irregular Hodge filtration
was clarified in terms of the $V$-filtrations of
the rescaling of integrable mixed twistor $\nbigd$-modules.

There are some remarkable and interesting
studies on the irregular Hodge theory.
Fres\'{a}n, Sabbah and Yu \cite{Fresan-Sabbah-Yu}
applied the irregular Hodge theory to
prove the analytic continuation and the functional equation
of some $L$-functions, conjectured by Broadhurst and Roberts.
They computed some numerical invariants
associated with Kloosterman connections
by using the irregular Hodge theory,
and with the aid of the theory of Patrikis and Taylor \cite{Patrikis-Taylor},
they obtained the desired property of the $L$-functions.
Casta\~{n}o Dom\'{\i}nguez, Martin, Reichelt,
Sabbah, Sevenheck and Yu computed 
the irregular Hodge numbers of confluent hypergeometric systems
in \cite{CDS, CDSR, Martin, Sabbah-irregular-Hodge, Sabbah-Yu2}.
(See also the work of Fedorov \cite{Fedorov}
in the regular singular case.)
Takahiro Saito described the irregular Hodge filtration of
the Fourier transform of a monodromic mixed Hodge module
\cite{Takahiro-Saito, Takahiro-Saito2}.

\subsection{Purposes of this paper}

In this paper, we study the basic part of
the theory of integrable mixed twistor $\nbigd$-modules
and the irregular Hodge filtrations
though a considerable part has been already established
in \cite{Sabbah-irregular-Hodge}.

First, we shall explain that the $\nbigd$-module
underlying an integrable mixed twistor $\nbigd$-module $\nbigt$
is equipped with the irregular Hodge filtration
if $\nbigt$ is {\em rescalable}.
We shall prove that the rescalability condition
is preserved by standard operations for
integrable mixed twistor $\nbigd$-modules.
In particular, in the algebraic setting,
the rescalability is preserved by
the standard $6$-operations.
Moreover,
for the projective push-forward,
the duality, the external tensor product,
and the non-characteristic inverse image,
the irregular Hodge filtrations are transformed
as in the case of Hodge modules.
For the localization,
the nearby cycle functor
and the vanishing cycle functor,
under the additional assumption of the regularity,
the irregular Hodge filtrations are transformed
as in the case of Hodge modules.
(See \cite[Remark 0.6]{Sabbah-irregular-Hodge}.)
See \S\ref{subsection;22.7.8.1} for some more details.

For those purposes,
we also revisit the basic part of the theory of
integrable mixed twistor $\nbigd$-modules.
In particular, we first study the Malgrange extension
of integrable mixed twistor $\nbigd$-modules,
which is helpful to clarify the conditions
for the existence of irregular Hodge filtrations.
We also study the non-characteristic inverse image
of integrable mixed twistor $\nbigd$-modules.
See \S\ref{subsection;22.7.8.2} for some more details.

\begin{rem}
In {\rm\cite{Sabbah-irregular-Hodge}}
the functoriality is preserved by
the pull back via smooth morphisms,
the push-forward via projective morphisms
and the external tensor product.
See also {\rm\cite{Esnault-Sabbah-Yu, Sabbah-Yu}}
for the push-forward. 
The duality was studied in {\rm\cite{Yu}},
and the external tensor product was studied in
{\rm\cite{Chen-Yu}}
in some special but interesting cases.
\hfill\qed
\end{rem}

\begin{rem}
We do not use the terminology
``irregular mixed Hodge modules''
in {\rm\cite{Sabbah-irregular-Hodge}}.
\hfill\qed
\end{rem}

It is another purpose to clarify the relation
between rescalable integrable mixed twistor $\nbigd$-modules
and exponential mixed Hodge modules.
In \cite{Kontsevich-Soibelman-exponential-Hodge},
Kontsevich and Soibelman introduced
the notion of exponential mixed Hodge structure
to study the cohomology of an algebraic variety
equipped with an algebraic function.
In \cite{Sabbah-irregular-Hodge},
it is proved that the Fourier transform
induces a fully faithful functor
from the category of exponential mixed Hodge structures
to the category of ``irregular mixed Hodge structures''.
We refine it to an equivalence
between exponential mixed Hodge modules and
integrable mixed twistor $\nbigd$-modules
which are rescalable along both $0$ and $\infty$.
See \S\ref{subsection;22.7.8.10}
for some more details.

\subsection{Some basic operations for integrable mixed twistor $\nbigd$-modules}
\label{subsection;22.7.8.2}

\subsubsection{Notation}

Let $X$ be a complex manifold.
Let $\nbigd_X$ denote the sheaf of holomorphic linear differential operators
on $X$.
We set $\nbigx:=\cnum\times X$.
Let $p_X:\nbigx\lrarr X$ denote the projection.
Let $\nbigr_X\subset\nbigd_{\nbigx}$
denote the sheaf of subalgebras
generated by $\lambda\cdot p_X^{\ast}\Theta_X$
over $\nbigo_{\nbigx}$,
where $\Theta_X$ denotes the tangent sheaf of $X$.
Let $\nbigrtilde_X:=\nbigr_X\langle\lambda^2\del_{\lambda}\rangle
\subset\nbigd_{\nbigx}$.
Let $\nbigc(X)$ denote the category of
$\nbigrtilde_X$-modules
underlying integrable mixed twistor $\nbigd_X$-modules.

We set $\gbigx:=\proj^1\times X$
and $\gbigx^{\infty}:=\{\infty\}\times X$
and $\gbigx^0:=\{0\}\times X$.
Let $p_{1,X}:\gbigx\lrarr X$ denote the projection.
Let $\nbigo_{\gbigx}(\ast\gbigx^{\infty})$
denote the sheaf of meromorphic functions on $\gbigx$
which may allow poles along $\gbigx^{\infty}$.
We set
$\nbigd_{\gbigx}(\ast\gbigx^{\infty}):=
\nbigo_{\gbigx}(\ast\gbigx^{\infty})\otimes_{\nbigo_{\gbigx}}
\nbigd_{\gbigx}$.
Let $\gbigr_X\subset\nbigd_{\gbigx}(\ast\gbigx^{\infty})$ denote
the sheaf of subalgebras generated by
$\lambda\cdot p_{1,X}^{\ast}\Theta_X$
over $\nbigo_{\gbigx}(\ast\gbigx^{\infty})$.
We set
$\gbigrtilde_X=\gbigr_X\langle\lambda^2\del_{\lambda}\rangle
\subset\nbigd_{\gbigx}(\ast\gbigx^{\infty})$.
Let $\gbigc(X)$ denote the full subcategory of
$\gbigrtilde_X$-modules $\gbigm$
such that
$\gbigm_{|\nbigx}\in\nbigc(X)$
and that
$\gbigm_{|\gbigx\setminus\gbigx^0}$
is a holonomic $\nbigd_{\gbigx\setminus\gbigx^0}$-module.

\subsubsection{Malgrange extension}

We introduce a strong regularity condition
of a holonomic $\nbigd$-module along a smooth hypersurface.
(See Definition \ref{df;21.4.14.5}.)
Then, we shall prove the following theorem.

\begin{thm}[Theorem
\ref{thm;21.1.23.1}]
There exists a unique functor
$\Upsilon:\nbigc(X)\lrarr\gbigc(X)$,
called the Malgrange extension,
such that the following holds
for any $\nbigm\in\nbigc(X)$:
\begin{itemize}
 \item $\Upsilon(\nbigm)_{|\nbigx}=\nbigm$.
 \item $\Upsilon(\nbigm)_{|\gbigx\setminus\gbigx^0}$
       is strongly regular along $\gbigx^{\infty}$.
\end{itemize}
\end{thm}
In \S\ref{subsection;21.6.29.1},
we shall see that the Malgrange extension $\Upsilon$ is compatible with
standard operations such as
the direct image,
the localization,
the Beilinson functors, 
the nearby cycle functors,
and the vanishing cycle functors.

\subsubsection{Inverse image functors}

Let $H$ be a hypersurface of $X$.
For any $\nbigrtilde_X$-module $\nbigm$,
we obtain $\nbigm(\ast H)=\nbigo_{\nbigx}(\ast\nbigh)\otimes\nbigm$.
This induces a functor from $\nbigc(X)$
to the category of $\nbigrtilde_X(\ast H)$-module.
Let $\nbigc(X;H)$ denote the essential image.
As the complement to \cite{Mochizuki-MTM},
we study the inverse image of objects
of $\nbigc(X;H)$ via a morphism of complex manifolds
$f:Y\lrarr X$. We set $H_Y:=f^{-1}(H)$.

We naturally obtain the $\nbigrtilde_Y(\ast H_Y)$-module
$f^{\ast}(\nbigm)
=\nbigo_{\nbigy}(\ast\nbigh_Y)
\otimes_{f^{-1}(\nbigo_{\nbigx}(\ast\nbigh))}f^{-1}(\nbigm)$.
(See \S\ref{subsection;21.6.29.11}.)

\begin{thm}[Theorem
\ref{thm;21.4.7.50}, Proposition \ref{prop;21.6.29.12}]
If $f$ is strictly non-characteristic for $\nbigm\in\nbigc(X;H)$,
i.e.,
$f_{|Y\setminus H_Y}$ is strictly non-characteristic for
$\nbigm_{|X\setminus H_X}$ (see {\rm\cite{Sabbah-pure-twistor}}),
then 
$f^{\ast}(\nbigm)$ is an object of $\nbigc(Y;H_Y)$. 
We also have
$f^{\ast}(\Upsilon(\nbigm))
=\Upsilon(f^{\ast}(\nbigm))$.
\end{thm}

Let us consider the case where
$f$ is a closed immersion but
not necessarily non-characteristic for $\nbigm$.
We impose an additional assumption that
there exists a finite tuple of hypersurfaces
$H^{(j)}$ $(j=1,\ldots,m)$ of $X$ such that
$f(Y)\cup H_X=\bigcap_{j=1}^m H^{(j)}\cup H_X$.
Then, we can define the inverse images
$(\lefttop{T}f^{\star})^j(\nbigm)$ $(j\in\seisuu,\,\,\star=!,\ast)$
as in the case of mixed Hodge modules \cite{saito2}.
(See \S\ref{subsection;21.6.22.2}.)
We shall see that the inverse image functors
are compatible with the Malgrange extensions
(Proposition \ref{prop;21.7.25.1}).

We can apply the inverse image functors
to define the tensor products
as in \S\ref{subsection;21.6.22.11}.
Together with the consideration for
non-characteristic inverse image,
we shall observe that 
the localization functors can be described
in terms of the tensor products.

\subsection{Rescalable objects and the irregular Hodge filtrations}
\label{subsection;22.7.8.1}

We set $\lefttop{\tau}X:=\cnum_{\tau}\times X$
and $\lefttop{\tau}X_0:=\{0\}\times X$.
For any $\nbigrtilde_X$-module $\nbigm$,
as in \cite{Sabbah-irregular-Hodge},
we naturally obtain
an $\nbigrtilde_{\lefttop{\tau}X\setminus\lefttop{\tau}X_0}$-module
$\lefttop{\tau}\nbigm$ as the pull back by
\[
\cnum_{\lambda}\times\cnum_{\tau}^{\ast}\times X
\lrarr \cnum_{\lambda}\times X,
\quad
(\lambda,\tau, x)\longmapsto (\lambda\tau^{-1},x).
\]
Similarly,
for any $\gbigrtilde_X$-module $\gbigm$,
we naturally obtain
an $\gbigrtilde_{\lefttop{\tau}X}(\ast\tau)$-module
$\lefttop{\tau}\gbigm$
(see \S\ref{subsection;21.6.29.20}).
Let $\pi:\lefttop{\tau}X\lrarr X$ denote the projection.

We say that $\nbigm\in\nbigc(X)$ is rescalable
if there exists $\nbigm_0\in\nbigc(\lefttop{\tau}X)$
such that
$\nbigm_{0|\lefttop{\tau}\nbigx\setminus \lefttop{\tau}\nbigx_0}
=\lefttop{\tau}\nbigm$.

\begin{thm}[Theorem
\ref{thm;21.3.12.12}, Theorem \ref{thm;21.3.15.2},
Proposition \ref{prop;21.4.16.11}]
\label{thm;22.7.28.10}
Suppose that $\nbigm\in\nbigc(X)$ is rescalable.
\begin{itemize}
 \item 
If $\nbigm_0\in\nbigc(\lefttop{\tau}X)$ satisfies
$\nbigm_{0|\lefttop{\tau}\nbigx\setminus\lefttop{\tau}\nbigx_0}
=\lefttop{\tau}\nbigm$,
we have
$\nbigm_0(\ast\tau)
=\lefttop{\tau}(\Upsilon(\nbigm))_{|\lefttop{\tau}\nbigx}$.
\item
$\lefttop{\tau}(\Upsilon(\nbigm))_{|\lefttop{\tau}\nbigx}$
is regular along $\tau$
in the sense that
each $V_a\bigl(
\lefttop{\tau}(\Upsilon(\nbigm))_{|\lefttop{\tau}\nbigx}
\bigr)$
is coherent over $\pi^{\ast}\nbigr_X$.
Moreover,
each $V_a\bigl(
\lefttop{\tau}(\Upsilon(\nbigm))_{|\lefttop{\tau}\nbigx}
\bigr)$
is flat over $\nbigo_{\cnum_{\lambda}\times\cnum_{\tau}}$.
\end{itemize}
 \end{thm}

Let $\iota_1:\{1\}\times X\to \nbigx=\cnum_{\lambda}\times X$
denote the inclusion.
For any $\nbigm\in \nbigc(X)$,
we naturally regard
$\Xi_{\DR}(\nbigm):=\iota_1^{\ast}(\nbigm)$
as a $\nbigd_X$-module.
Let us explain that Theorem \ref{thm;22.7.28.10}
implies the existence of a coherent filtration of
$\Xi_{\DR}(\nbigm)$ if $\nbigm$ is rescalable.

In general, a $\nbigd_X$-module $M$ with a coherent filtration $F$
induces a coherent $\nbigr_X$-module $\Rtilde_F(M)$
as the analytification of the Rees module
$R_F(M)=\sum_{j\in\seisuu} F_j(M)\lambda^j$.
There exists a natural monomorphism
$\Rtilde_F(M)\lrarr p_X^{\ast}(M)(\ast\lambda)$,
which induces
$\Rtilde_F(M)(\ast\lambda)=p_X^{\ast}(M)(\ast\lambda)$.
Let us consider the $\cnum^{\ast}$-action on
$\cnum_{\lambda}\times X$
defined by $b(\lambda,x)=(b\lambda,x)$ $(b\in\cnum^{\ast})$.
By the construction,
$p_X^{\ast}(M)(\ast\lambda)$ is naturally
$\cnum^{\ast}$-equivariant,
and $\Rtilde_F(M)$ is
a $\cnum^{\ast}$-invariant subsheaf of $p_X^{\ast}(M)(\ast\lambda)$.
It is more or less standard to observe that
a coherent filtration of $M$ bijectively corresponds to
a coherent $\nbigr_X$-submodule $\nbigf$ of
$p_X^{\ast}(M)(\ast\lambda)$
such that
(i) $\nbigf(\ast\lambda)=p_X^{\ast}(M)(\ast\lambda)$,
(ii) $\nbigf$ is $\cnum^{\ast}$-invariant.
(See Proposition \ref{prop;21.4.17.1}.)

Let $\iota_{\lambda=\tau}\colon\nbigx\lrarr \lefttop{\tau}\nbigx$
be the morphism induced by the diagonal embedding
$\cnum_{\lambda}\lrarr \cnum_{\lambda}\times\cnum_{\tau}$,
$\lambda\longmapsto(\lambda,\lambda)$.
For any rescalable $\nbigm\in\nbigc(X)$,
$\iota_{\lambda=\tau}^{\ast}
\lefttop{\tau}\Upsilon(\nbigm)$
is naturally identified with
$p_X^{\ast}(\Xi_{\DR}(\nbigm))(\ast\lambda)$
in a way compatible with the $\cnum^{\ast}$-actions.
For each $a\in\real$,
we obtain the coherent $\nbigr_X$-module
$\iota_{\lambda=\tau}^{\ast}V_a\bigl(
\lefttop{\tau}\Upsilon(\nbigm)
\bigr)$,
which is naturally a $\cnum^{\ast}$-invariant subsheaf of
$\iota_{\lambda=\tau}^{\ast}
\lefttop{\tau}\Upsilon(\nbigm)$
such that
$\iota_{\lambda=\tau}^{\ast}V_a\bigl(
\lefttop{\tau}\Upsilon(\nbigm)
\bigr)(\ast\lambda)
=\iota_{\lambda=\tau}^{\ast}
\lefttop{\tau}\Upsilon(\nbigm)$.

\begin{cor}[Proposition
\ref{prop;21.6.29.20}]
For any rescalable object $\nbigm\in\nbigc(X)$,
the underlying $\nbigd_X$-module
$\Xi_{\DR}(\nbigm)$
is equipped with  
the unique coherent filtration $F^{\irr}_{a+\bullet}$
such that
$\Rtilde_{F^{\irr}_{a+\bullet}}\Xi_{\DR}(\nbigm)
=\iota_{\lambda=\tau}^{\ast}
V_a\bigl(
\lefttop{\tau}\Upsilon(\nbigm)
 \bigr)$
under the natural identification
$\iota_{\lambda=\tau}^{\ast}
\lefttop{\tau}\Upsilon(\nbigm)=
p_X^{\ast}(\Xi_{\DR}(\nbigm))(\ast\lambda)$.
It is called the irregular Hodge filtration of $\nbigm$.
\end{cor}

We have the following functoriality
(see 
Proposition \ref{prop;21.6.29.21}, Theorem \ref{thm;21.3.25.20},
Theorem \ref{thm;21.3.26.11},
Corollary \ref{cor;21.6.29.22} and
Proposition \ref{prop;21.6.29.23}).

 \begin{thm}
Rescalable objects are preserved by
(i) the direct image via any projective morphism,
(ii) the duality functor,
(iii) the external product, and
(iv)  the non-characteristic inverse image.
In these cases,
the irregular Hodge filtrations are functorial
in the natural ways.
\end{thm}

We also have the following functoriality.
\begin{thm}[Proposition
\ref{prop;21.6.29.24},
Theorem \ref{thm;21.3.17.20}]
Rescalable objects are preserved by
the localization along any hypersurfaces,
and by the Beilinson functors and the vanishing cycle functor
along any holomorphic function. 
If $\nbigm\in\nbigc(X)$ is regular,
the irregular Hodge filtration is functorial
in natural ways. 
 \end{thm}
 
We study
the rescalability at $\infty$
in \S\ref{section;21.7.3.1}.

\subsection{Partial Fourier transforms}
\label{subsection;22.7.8.10}

\subsubsection{$\nbigd$-modules}

On $\proj^1_t\times\proj^1_{\tau}$,
we set
$H_{\infty}=
(\proj^1_t\times\{\infty\})\cup(\{\infty\}\times\proj^1_{\tau})$.
Let $\nbigo_{\proj^1_t\times\proj^1_{\tau}}(\ast H_{\infty})$
denote the sheaf of meromorphic functions on
$\proj^1_t\times\proj^1_{\tau}$
which may allow poles along $H_{\infty}$.
We have the meromorphic flat bundles
$L(t\tau)=(\nbigo_{\proj^1_t\times\proj^1_{\tau}}(\ast H_{\infty}),
d+d(t\tau))$
and 
$L(-t\tau)=(\nbigo_{\proj^1_t\times\proj^1_{\tau}}(\ast H_{\infty}),
d-d(t\tau))$.
The pull back of $L(\pm t\tau)$
by $\proj^1_t\times\proj^1_{\tau}\times X\lrarr \proj^1_t\times
\proj^1_{\tau}$
are also denoted by $L(\pm t\tau)$.

Let $p_{t,X}:\proj^1_t\times\proj^1_{\tau}\times X\lrarr
\proj^1_t\times X$ and
$p_{\tau,X}:\proj^1_{t}\times\proj^1_{\tau}\times X\lrarr
\proj^1_{\tau}\times X$
denote the projections.
For any holonomic
$\nbigd_{\proj^1_{t}\times X}$-module $M$
such that
$\nbigo(\ast(\{\infty\}\times X))\otimes_{\nbigo}M\simeq M$,
we obtain its partial Fourier transforms $\FT_{\pm}(M)$
as
\begin{equation}
\label{eq;21.7.17.1}
 \FT_{\pm}(M)=
 p_{\tau,X\dagger}^0
 \bigl(
 p_{t,X}^{\ast}(M)\otimes L(\pm t\tau)
 \bigr).
\end{equation}
We have the standard inversion formula
$\FT_{\mp}\circ \FT_{\pm}(M)\simeq M$.
Let $\pi_X:\proj^1_t\times X\lrarr X$ denote the projection.
Then, we have $\pi_{X\dagger}^j(M)=0$ $(j\in\seisuu)$
if and only 
$\nbigo_{\proj^1_{\tau}\times X}(\ast(\{0\}\times X))
\otimes
\FT_{\pm}(M)=\FT_{\pm}(M)$.

\subsubsection{$\nbigrtilde$-modules
and integrable mixed twistor $\nbigd$-modules}

We set
$\nbigctilde(\cnum\times X):=
\nbigc(\proj^1\times X;\{\infty\}\times X)$.
We naturally obtain
the functors
$\FT_{\pm}:\nbigctilde(\cnum\times X)
\lrarr \nbigctilde(\cnum\times X)$
which are enhancement of (\ref{eq;21.7.17.1}).
(See \S\ref{subsection;21.7.25.2}.)
We have the inversion 
$\FT_{\mp}\circ\FT_{\pm}(\nbigm)
\simeq
\lambda^{-1}\nbigm$
for $\nbigm\in\nbigctilde(\cnum\times X)$.
Let $\nbigctilde(\cnum\times X)_{\ast}\subset\nbigctilde(\cnum\times X)$
denote the full subcategory of
$\nbigm\in\nbigctilde(\cnum\times X)$
such that $\pi^j_{X\dagger}(\nbigm[\ast (\{\infty\}\times X)])=0$
$(j\in\seisuu)$.
Let $\nbigctilde(\cnum\times X,[\ast 0])\subset
\nbigctilde(\cnum\times X)$
denote the full subcategory of
$\nbigm\in\nbigctilde(\cnum\times X)$
such that $\nbigm[\ast (\{0\}\times X)]\simeq\nbigm$.
Then, $\FT_{\pm}$ induce equivalences
\[
 \nbigctilde(\cnum\times X)_{\ast}
 \simeq
 \nbigctilde(\cnum\times X,[\ast 0]).
\]

For $\vecn=(n_1,n_2)\in\seisuu^2$
such that $\gcd(n_1,n_2)=1$ and $n_1\geq n_2\geq 0$,
we consider the $\cnum^{\ast}$-action
$\rho_{\vecn}:\cnum^{\ast}\times\cnum_{\lambda}\times \cnum_t\times X
\lrarr\cnum_{\lambda}\times\cnum_t\times X$
determined by
$\rho_{\vecn}(a,\lambda,t,x)=(a^{n_1}\lambda,a^{n_2}t,x)$.
An $\nbigrtilde_{\cnum_t\times X}$-module $\nbigm$
is called homogeneous with respect to $\rho_{\vecn}$
if there exists an isomorphism
$\rho_{\vecn}^{\ast}\nbigm\simeq\nbigm$
satisfying a cocycle condition.
Let $\nbigctilde_{\vecn}(\cnum\times X)\subset\nbigctilde(\cnum\times X)$
denote the full subcategory of
the homogeneous objects with respect to $\rho_{\vecn}$.
We obtain the subcategories
$\nbigctilde_{\vecn}(\cnum\times X)_{\ast}$
and
$\nbigctilde_{\vecn}(\cnum\times X,[\ast 0])$
naturally.
For $\vecn$ as above,
we set $\FT(\vecn)=(n_1,n_1-n_2)$.
Then, $\FT_{\pm}$ induce equivalences
\[
 \nbigctilde_{\vecn}(\cnum\times X)_{\ast}
 \simeq
 \nbigctilde_{\FT(\vecn)}(\cnum\times X,[\ast 0]).
\]

Let $\nbigc_{\res,0,\infty}(X)\subset\nbigc(X)$
denote the full subcategory of the objects
which are rescalable along both $0$ and $\infty$.
We can observe that
the restriction to $\{1\}\times X\subset\cnum\times X$
induces
$\nbigctilde_{(1,1)}(\cnum\times X,[\ast 0])
\simeq
\nbigc_{\res,0,\infty}(X)$
(Proposition \ref{prop;21.7.11.31}).
Let $\nbigctilde_{\Hod}(\cnum\times X)_{\ast}$
denote the category of filtered
$\nbigd_{\proj^1\times X}(\ast(\{\infty\}\times X))$-modules
$(M,F)$ on $\proj^1\times X$
such that
the analytification of the Rees module of $(M,F)$
are objects of $\nbigctilde_{(1,0)}(\cnum\times X)_{\ast}$.
We have
$\nbigctilde_{\Hod}(\cnum\times X)_{\ast}
\simeq
\nbigctilde_{(1,0)}(\cnum\times X)_{\ast}$.
We obtain the following theorem.
\begin{thm}[Theorem
\ref{thm;21.7.14.2}, see also \S\ref{subsection;21.7.26.12}]
\label{thm;21.7.17.3}
$\nbigctilde_{\Hod}(\cnum\times X)_{\ast}
\simeq
\nbigc_{\res,0,\infty}(X)$.
\end{thm}

The equivalences are compatible with
other basic functors
as explained in \S\ref{subsection;21.7.25.3}
and \S\ref{subsection;21.7.26.12}.
Moreover, we shall enhance Theorem \ref{thm;21.7.17.3}
to equivalences between
rescalable mixed twistor $\nbigd$-modules
and exponential $\real$-Hodge modules
(see \S\ref{subsection;21.7.16.71}),
which are also compatible with other basic functors.
\begin{thm}[Theorem
\ref{thm;21.7.16.3}]
We have
$\MHMtilde(\cnum\times X,\real)_{\ast}
\simeq
\MTM^{\integral}_{\res}(X,\real)$.
\end{thm}

\paragraph{Acknowledgement}

I thank Claude Sabbah for many useful discussions and his kindness.
This study grew out of
our discussions on irregular Hodge filtrations,
and my effort to understand \cite{Sabbah-irregular-Hodge}.
I am partially motivated to understand the work \cite{Fresan-Sabbah-Yu}
by Javier Fres\'{a}n, Claude Sabbah and Jeng-Daw Yu,
and to understand the notion of exponential mixed Hodge structure
\cite{Kontsevich-Soibelman-exponential-Hodge} due to
Maxim Kontsevich and Yan Soibelman.
I thank Claus Hertling and Christian Sevenheck
for discussions on TERP structures, integrable twistor structures
and rescalable objects.
I appreciate Takahiro Saito for stimulating discussions
on monodromic Hodge modules and irregular Hodge filtrations,
which are helpful to keep my interest in twistor $\nbigd$-modules.
I also thank him for his comments to improve this paper.
I am grateful to Yota Shamoto for his comments and questions
to improve this manuscript.
I thank Akira Ishii and Yoshifumi Tsuchimoto for constant encouragements.

I am partially supported by
the Grant-in-Aid for Scientific Research (S) (No. 17H06127),
the Grant-in-Aid for Scientific Research (S) (No. 16H06335),
the Grant-in-Aid for Scientific Research (A) (No. 21H04429),
the Grant-in-Aid for Scientific Research (C) (No. 15K04843),
and
the Grant-in-Aid for Scientific Research (C) (No. 20K03609),
Japan Society for the Promotion of Science.

\section{Preliminary}

\subsection{A condition for goodness of meromorphic flat bundles}

Let $X$ be a complex manifold
with a simple normal crossing hypersurface $H$.
Let $\vecnbigi$ be a good system of irregular values
on $(X,H)$.
(See \cite[Definition 2.4.2]{Mochizuki-wild}
for the notion of good system of irregular values.)
Let $\nbigv$ be a coherent reflexive
$\nbigo_X(\ast H)$-module
with a flat connection $\nabla$.
We set $H_{[2]}:=\bigcup_{i\neq j}(H_i\cap H_j)$.
We shall prove the following proposition
in \S\ref{subsection;21.3.8.1}
after the preliminaries in
\S\ref{subsection;21.3.8.2}--\ref{subsection;21.3.8.3}.
\begin{prop}
\label{prop;21.1.22.2}
 Assume that $\nbigv_{|X\setminus H_{[2]}}$
is unramifiedly good
on $(X\setminus H_{[2]},H\setminus H_{[2]})$
whose good system of irregular values
is contained in $\vecnbigi_{|H\setminus H_{[2]}}$.
Then, $\nbigv$ is unramifiedly good on $(X,H)$
whose good system of irregular values
is contained in $\vecnbigi$.
\end{prop}

\subsubsection{The formal structure at the intersection of two hypersurfaces}
\label{subsection;21.3.8.2}

We set $\nbiga_0:=\cnum[\![z_3,\ldots,z_n]\!]$.
We set $\nbiga:=\nbiga_0[\![z_1,z_2]\!]$
and $\nbigr:=\nbiga[z_1^{-1},z_2^{-1}]$.
We set
$\gbigk_1:=\nbiga_0((z_1))$
and $\gbigk_2:=\nbiga_0((z_2))$.
We set
$\gbigr_1:=\gbigk_1[\![z_2]\!]$
and
$\gbigr_2:=\gbigk_2[\![z_1]\!]$.

Let $L$ be a free $\nbiga$-module
equipped with a meromorphic flat connection
\[
 \nabla:
 L\lrarr
 (L\otimes_{\nbiga}\nbigr)\otimes_{\nbigr}
 \Omega^1_{\nbigr/\cnum}.
\]
Let $\nbigi$ be a good set of irregular values
contained in $\nbigr/\nbiga$.
We may assume that the natural map
$\nbigi\lrarr \nbigr/\nbiga[z_2^{-1}]$ is injective
by exchanging $(z_1,z_2)$ if necessary.
Let $\nbigi^{(2)}$ denote the image of
the natural map
$\nbigi\lrarr\nbigr/\nbiga[z_1^{-1}]$.
We assume the following conditions.
\begin{itemize}
 \item There exists a decomposition
\begin{equation}
\label{eq;21.2.11.4}
       (L,\nabla)\otimes\gbigr_1
       =\bigoplus_{\gminib\in\nbigi^{(2)}}
       (L^{(2)}_{\gminib},\nabla^{(2)}_{\gminib})
\end{equation}
       such that
       $\nablatilde^{(2)}_{\gminib}:=\nabla^{(2)}_{\gminib}-d\gminib\id$
       is logarithmic with respect to $L^{(2)}_{\gminib}$
       in the sense that
       $\nablatilde^{(2)}_{\gminib}(\del_{z_i})L^{(2)}_{\gminib}
       \subset L^{(2)}_{\gminib}$ $(i\neq 2)$
       and $\nablatilde^{(2)}_{\gminib}(z_2\del_{z_2})L^{(2)}_{\gminib}
       \subset L^{(2)}_{\gminib}$.
\item There exists a decomposition
\begin{equation}
\label{eq;21.2.11.3}
       (L,\nabla)\otimes\gbigr_2
       =\bigoplus_{\gminia\in\nbigi}
       (L^{(1)}_{\gminia},\nabla^{(1)}_{\gminia})
\end{equation}
       such that
       $\nablatilde^{(1)}_{\gminia}:=\nabla^{(1)}_{\gminia}-d\gminia\id$
      is logarithmic with respect to $L^{(1)}_{\gminia}$
      in the sense that
       $\nablatilde^{(1)}_{\gminia}(\del_{z_i})L^{(1)}_{\gminia}
       \subset L^{(1)}_{\gminia}$ $(i\neq 1)$
       and $\nablatilde^{(1)}_{\gminia}(z_1\del_{z_1})L^{(1)}_{\gminia}
       \subset L^{(1)}_{\gminia}$.
\end{itemize}

\begin{lem}
\label{lem;21.2.11.11}
 There exists a decomposition
\[
       (L,\nabla)
       =\bigoplus_{\gminia\in\nbigi}
       (L_{\gminia},\nabla_{\gminia})
\]
       such that
       $\nabla_{\gminia}-d\gminia\id$
       is logarithmic with respect to $L_{\gminia}$.
\end{lem}
\pf
Let us consider the case $\nbigi=\{0\}$.
Let $\vecv$ be a frame of $L$.
Let $A_i\in M_r(\nbigr)$ $(i=1,2,\ldots,n)$ be determined by
$z_i\nabla_{z_i}\vecv=\vecv A_i$ $(i=1,2)$
and
$\nabla_{z_i}\vecv=\vecv A_i$ $(i=3,\ldots,n)$.
By the condition,
we obtain
$A_i\in M_r(\gbigr_1)\cap M_r(\gbigr_2)$,
which implies
$A_i\in M_r(\nbiga)$.
Hence, $\nabla$ is logarithmic.

\vspace{.1in}
Because $\nbigi$ is a good set of irregular values,
the partial order $\leq_{\seisuu^2}$
induces a total order on
$\bigl\{
\ord(\gminia-\gminib)\,\big|\,
\gminia,\gminib\in\nbigi
\bigr\}$.
Namely, we can order the elements of
$\bigl\{
\ord(\gminia-\gminib)\,\big|\,
\gminia,\gminib\in\nbigi
\bigr\}$
as
$\vecm(1)>_{\seisuu^2}\vecm(2)>_{\seisuu^2}
\cdots>_{\seisuu^2}\vecm(k)$ in $\seisuu^2_{\leq 0}$.
We use an induction on $k$.
Let $q:\nbigi\lrarr \cnum \vecz^{\vecm(k)}$
denote the projection.
We may assume that the image of $q$ contains
at least two elements.

We consider the map
\[
 z^{-\vecm(k)}z_1\nabla_{z_1}:L\lrarr L\otimes\nbigr.
\]
By the conditions, we obtain
$z^{-\vecm(k)}z_1\nabla_{z_1}\bigl(
 L\otimes\gbigr_i\bigr)
 \subset
 L\otimes\gbigr_i$,
 and hence
 $z^{-\vecm(k)}z_1\nabla_{z_1}(L)
 \subset L$.
It is easy to check that
$z^{-\vecm(k)}z_1\nabla_{z_1}\bigl(z_iL\bigr)
 \subset z_iL$
$(i=1,2)$.
It is standard that the induced morphism
\[
 z^{-\vecm(k)}z_1\nabla_{z_1}:
 L/z_1L
 \lrarr
 L/z_1L
\]
is an $\nbiga/z_1\nbiga$-homomorphism.
Because of the existence of the decomposition
(\ref{eq;21.2.11.3}),
there exists a decomposition
\begin{equation}
\label{eq;21.2.11.1}
 (L/z_1L)\otimes\gminik_2
 =\bigoplus_{\alpha\in\cnum}
 \bigl(
  (L/z_1L)\otimes\gminik_2
 \bigr)_{\alpha}
\end{equation}
such that the following holds.
\begin{itemize}
 \item The decomposition (\ref{eq;21.2.11.1})
       is preserved by
       $z^{-\vecm(k)}z_1\nabla_{z_1}$.
       On  $\bigl(
  (L/z_1L)\otimes\gminik_2
       \bigr)_{\alpha}$,
       the eigenvalues are
       contained in $\nbiga/z_1\nbiga$,
       and equal to
       $m(k)_1^{-1}\alpha$ modulo
       $(z_2,\ldots,z_n)$.
\end{itemize}
Then, the decomposition (\ref{eq;21.2.11.1})
extends to a decomposition
\begin{equation}
\label{eq;21.2.11.2}
  L/z_1L
 =\bigoplus_{\alpha\in\cnum}
 \bigl(
  L/z_1L
 \bigr)_{\alpha}
\end{equation}
such that the following holds.
\begin{itemize}
 \item The decomposition (\ref{eq;21.2.11.1})
       is preserved by $z^{-\vecm(k)}z_1\nabla_{z_1}$.
 On $\bigl(L/z_1L\bigr)_{\alpha}$,
       the eigenvalues are
       $m(k)_1^{-1}\alpha$ modulo $(z_2,\ldots,z_n)$.
\end{itemize}
It is standard that there exists a decomposition
\begin{equation}
 \label{eq;21.2.11.10}
  L=\bigoplus_{\alpha\in\cnum}
  L_{\alpha}
\end{equation}
such that the following holds
(for example, see \cite[\S2.1.5]{Mochizuki-Stokes}).
\begin{itemize}
 \item  The decomposition (\ref{eq;21.2.11.10})
	is preserved by
	$z^{-\vecm(k)}z_1\nabla_{z_1}$.
 \item $L_{\alpha}/z_1L_{\alpha}
       = \bigl(L/z_1L\bigr)_{\alpha}$.
\end{itemize}
By the construction,
there exists the following decomposition.
\begin{equation}
 L_{\alpha}\otimes\gbigr_2
  =\bigoplus_{\substack{
  \gminia\in\nbigi\\
 q(\gminia)=\alpha z^{\vecm(k)}}}
  L^{(1)}_{\gminia}.
\end{equation}
It particularly implies 
$\nabla_{z_j}L_{\alpha}
 \subset
 L_{\alpha}\otimes\nbigr$
 $(j=1,\ldots,n)$.

If $m(k)_2=0$,
$L_{\alpha}\otimes\gbigr_1$ is logarithmic with respect to
$\nabla$.
If $m(k)_2<0$, the map
\[
 z^{-\vecm(k)}z_1\nabla_{z_1}:
 L/z_2L\lrarr
 L/z_2L
\]
is an $\nbiga/z_2\nbiga$-homomorphism.
The decomposition (\ref{eq;21.2.11.10})
induces a decomposition
\begin{equation}
  L/z_2L=\bigoplus_{\alpha\in\cnum}
  L_{\alpha}/z_2L_{\alpha},
\end{equation}
which is preserved by
$z^{-\vecm(k)}z_1\nabla_{z_1}$.
The eigenvalues on 
$L_{\alpha}/z_2L_{\alpha}$ are
$m(k)_1^{-1}\alpha$ modulo $(z_1,z_3,\ldots,z_n)$.
Hence, we obtain
\begin{equation}
 L_{\alpha}\otimes\gbigr_1
  =\bigoplus_{\substack{\gminib\in\nbigi^{(2)}\\
 q(\gminib)=\alpha\vecz^{\vecm(k)}
 }}
  L^{(2)}_{\gminib}.
\end{equation}
Therefore, we can apply the assumption of the induction to
each $(L_{\alpha},\nabla)$.
\hfill\qed

\subsubsection{The extension of formal sections of reflexive sheaves}
\label{subsection;21.3.8.3}

Let $Y$ be a complex manifold.
Let $H_Y$ be a smooth hypersurface of $Y$.
Let $\nbige$ be a coherent reflexive $\nbigo_Y$-module.
Let $Z\subset H_Y$ be a complex analytic closed subset
with $\dim Z\leq\dim H_Y-2$.
\begin{lem}
\label{lem;21.2.15.1}
A section of $\nbige_{|\Hhat_Y\setminus Z}$
uniquely extends to a section of $\nbige_{|\Hhat_Y}$.
\end{lem}
\pf
It is enough to study the claim locally around any point of $Y$.
Because $\nbige$ is reflexive,
there exists a morphism of locally free $\nbigo_Y$-modules
$\nbigetilde_1\lrarr\nbigetilde_2$
such that the kernel is isomorphic to $\nbige$.
By the faithful flatness of the formal completion,
we obtain the exact sequence
$0\lrarr \nbige_{|\Hhat_Y}
\lrarr \nbigetilde_{1|\Hhat_Y}
\lrarr \nbigetilde_{2|\Hhat_Y}$.

Let $U$ be an open subset of $H_Y$.
Let $s$ be any local section of
$\nbige_{|\Hhat_Y\setminus Z}$
on $U\setminus Z$.
We may regard $s$
as a section of
$\nbigetilde_{1|\Hhat_Y\setminus Z}$
on $U\setminus Z$.
By Hartogs theorem,
$s$ uniquely extends to
a section $\stilde$ of $\nbigetilde_{2|\Hhat_Y}$ on $U$.
Because $\nbigetilde_{2|\Hhat_Y}$ is locally free,
the induced local section of
$\nbigetilde_{2|\Hhat_Y}$ is $0$.
Hence, $\stilde$ is a section of
$\nbige_{|\Hhat_Y}$.
\hfill\qed

\subsubsection{Proof of Proposition
\ref{prop;21.1.22.2}}
\label{subsection;21.3.8.1}

It is enough to prove the case where
$X$ is a neighbourhood of
$(0,\ldots,0)$ in $\cnum^n$
and $H=X\cap\bigl(\bigcup_{i=1}^{\ell}\{z_i=0\}\bigr)$.
Let $\nbigi$ be a good set of irregular values on $(X,H)$,
which induces a good system of irregular values
$\vecnbigi$ on $(X,H)$.
Let $\nbigv$ be a meromorphic flat bundle on $(X,H)$
as in Proposition \ref{prop;21.1.22.2}
for $\vecnbigi$.
Let $\nbigv^{\DM}$ denote the Deligne-Malgrange lattice of
$(\nbigv,\nabla)$.

Let us consider the case $\ell=2$.
There exists a complex analytic closed subset
$Z\subset H_1\cap H_2$
such that
(i) $\dim Z\leq n-3$,
(ii) $\nbigv^{\DM}$
is locally free on $X\setminus Z$.
We take any $P\in (H_1\cap H_2)\setminus Z$.
By applying Lemma \ref{lem;21.2.11.11}
to the formal completion of $\nbigv^{\DM}$ at $P$,
we obtain that
$(\nbigv,\nabla)$ is unramifiedly good at $P$,
and the good set of irregular values is contained in $\nbigi$.
Hence, there exists the decomposition
\[
(\nbigv^{\DM},\nabla)_{|\Hhat_1\setminus Z}
=\bigoplus_{\gminia\in \nbigi}
 (\nbigv^{\DM}_{\Hhat_1\setminus Z,\gminia},\nabla_{\gminia})
\]
such that $\nabla_{\gminia}-d\gminia\id$ are logarithmic.
Let $U\subset H_1\setminus Z$ be an open subset.
There exist locally free
$\nbigo_{\Hhat_1}$-modules $\nbigl_{\gminia}$
equipped with a logarithmic flat connection
$\nabla^{\reg}_{\gminia}$
such that
$(\nbigl_\gminia,\nabla^{\reg}_{\gminia})_{|\Uhat}
\simeq
 (\nbigv^{\DM}_{\Hhat_1\setminus Z,\gminia},\nabla_{\gminia}-d\gminia)_{|\Uhat}$.
Because $\dim Z\leq \dim H_1-2$,
the isomorphisms extend to isomorphisms
$(\nbigl_\gminia,\nabla^{\reg}_{\gminia})_{\Hhat_1\setminus Z}
\simeq
 (\nbigv^{\DM}_{\Hhat_1\setminus Z,\gminia},\nabla_{\gminia}-d\gminia)$.
By Lemma \ref{lem;21.2.15.1},
we obtain the isomorphism
$\bigoplus_{\gminia\in\nbigi}\nbigl_{\gminia}
\simeq
\nbigv^{\DM}_{|\Hhat_1}$,
under which
$\bigoplus_{\gminia\in\nbigi} (\nabla^{\reg}_{\gminia}+d\gminia\id)
=\nabla$.
Thus, the claim of Proposition \ref{prop;21.1.22.2}
is proved in the case $\ell=2$.

\vspace{.1in}

Let us consider the case $\ell\geq 3$.
Let $H_{[3]}$ denote the union of
the intersections of
three mutually distinct components of $H$.
There exists a decomposition
$(\nbigv^{\DM},\nabla)_{|\Hhat_1\setminus H_{[3]}}=\bigoplus_{\gminia\in\nbigi}
(\nbigv^{\DM}_{\Hhat_1\setminus H_{[3]},\gminia},\nabla_{\gminia})$
such that
$\nabla_{\gminia}-\sum d\gminia\pi_{\gminia}$
are logarithmic.
By a similar argument,
we can prove that
it extends to a decomposition of
locally free $\nbigo_{\Hhat_1}$-modules
such that
$\nabla_{\gminia}-\sum d\gminia\pi_{\gminia}$
are logarithmic on the direct summand.
Then, the claim of Proposition \ref{prop;21.1.22.2} follows.
\hfill\qed

\subsection{An auxiliary resolution for a set of irregular values}
\label{subsection;21.1.22.1}

Let $X$ be a neighbourhood of $(0,\ldots,0)$ in $\cnum^n$.
We set
$H_i=X\cap\{z_i=0\}$ and 
$H=\bigcup_{i=1}^{\ell}H_i$ for some $1\leq \ell\leq n$.
We set
$\Xtilde:=\cnum_t\times X$,
$\Htilde_i:=\cnum_t\times H_i$
and
$\Htilde:=\cnum_t\times H$.
We set
$\Xtilde^{0}:=\{0\}\times X$,
$\Htilde_i^{0}:=\{0\}\times H_i$
and
$\Htilde^{0}:=\{0\}\times H$.
We put
$H_{[2]}:=\bigcup_{1\leq i\neq j\leq \ell}H_i\cap H_j$,
$\Htilde_{[2]}:=\cnum_t\times H_{[2]}$
and 
$\Htilde^{0}_{[2]}:=\{0\}\times H_{[2]}$.

\subsubsection{An auxiliary resolution associated with a totally ordered sequence}

Let $\leq_{\seisuu^{\ell}}$ denote the partial order
on $\seisuu^{\ell}$ defined as
$(a_1,\ldots,a_{\ell})\leq_{\seisuu^{\ell}}
(b_1,\ldots,b_{\ell})
\stackrel{\rm def}{\Longleftrightarrow}
a_i\leq b_i\,\,(\forall i)$.
Let $\nbigs\subset\seisuu_{\leq 0}^{\ell}$
be a finite subset
which is totally ordered with respect to
$\leq_{\seisuu^{\ell}}$,
i.e., we have the ordering
$\nbigs=\{\vecm(1),\vecm(2),\ldots,\vecm(k)\}$
such that
$(0,\ldots,0)>_{\seisuu^{\ell}}\vecm(1)>_{\seisuu^{\ell}}
\vecm(2)>_{\seisuu^{\ell}}
>\cdots>_{\seisuu^{\ell}}\vecm(k)$
in $\seisuu_{\leq 0}^{\ell}$.

Let us construct a projective morphism of complex manifolds
$f_{1}:\Xtilde_{\vecm(1)}\lrarr \Xtilde$.
First, let 
$\Xtilde_{\vecm(1),\ell}\lrarr \Xtilde$
be the projective morphism of complex manifolds
obtained as the $-m(1)_{\ell}$-successive blow up
along the intersection of
the proper transform of $\Xtilde^0$
and the inverse image of $\Htilde_{\ell}$.
Let
$\Xtilde_{\vecm(1),\ell-1}\lrarr \Xtilde_{\vecm(1),\ell}$
be the projective morphism of complex manifolds
obtained as the $-m(1)_{\ell-1}$-successive blow up
along the intersection of
the proper transform of $\Xtilde_0$
and the inverse image of $\Htilde_{\ell-1}$.
Inductively,
we construct morphisms
\[
 \Xtilde_{\vecm(1),1}
 \lrarr
  \Xtilde_{\vecm(1),2}\lrarr
  \cdots\lrarr\Xtilde_{\vecm(1),\ell}\lrarr\Xtilde,
\]
where $\Xtilde_{\vecm(1),k-1}$
is the $-m(1)_{k-1}$-successive blow up
along the intersection of the proper transform of
$\Xtilde^0$
and the inverse image of
$\Htilde_{k-1}$.
We set $\Xtilde_{\vecm(1)}:=\Xtilde_{\vecm(1),1}$,
and let $f_1:\Xtilde_{\vecm(1)}\lrarr\Xtilde$
denote the induced morphism.

We construct
$f_2:\Xtilde_{\vecm(2)}\lrarr\Xtilde_{\vecm(1)}$
as follows.
Let $\Xtilde_{\vecm(2),\ell}\lrarr\Xtilde_{\vecm(1)}$
be the projective morphism of complex manifolds
obtained as the $-m(2)_{\ell}+m(1)_{\ell}$-successive blow up
along the intersection of
the proper transform of $\Xtilde^0$
and the inverse image of $\Htilde_{\ell}$.
We obtain the projective morphism of complex manifolds
$\Xtilde_{\vecm(2),\ell-1}\lrarr \Xtilde_{\vecm(2),\ell}$
as the $-m(2)_{\ell-1}+m(1)_{\ell-1}$-successive blow up
along the intersection of
the proper transform of $\Xtilde^{0}$
and the inverse image of $\Htilde_{\ell-1}$.
Inductively,
we construct the following morphisms:
\[
 \Xtilde_{\vecm(2),1}
 \lrarr
  \Xtilde_{\vecm(2),2}\lrarr
  \cdots\lrarr\Xtilde_{\vecm(2),\ell}\lrarr\Xtilde_{\vecm(1)}.
\]
We set $\Xtilde_{\vecm(2)}:=\Xtilde_{\vecm(2),1}$,
and let $f_2:\Xtilde_{\vecm(2)}\lrarr\Xtilde_{\vecm(1)}$
denote the induced morphism.

Similarly and inductively,
we construct the following morphisms:
\[
\begin{CD}
 \Xtilde_{\vecm(k)}
 @>{f_k}>>
 \Xtilde_{\vecm(k-1)}
 @>{f_{k-1}}>>
 \cdots
 @>{f_2}>>\Xtilde_{\vecm(1)}
 @>{f_1}>>\Xtilde.
\end{CD}
\]
We set
$\Xtilde_{\nbigs}:=\Xtilde_{\vecm(k)}$,
and let $F_{\nbigs}:\Xtilde_{\nbigs}\lrarr \Xtilde$
denote the induced morphism.

\begin{lem}
\mbox{{}}\label{lem;21.3.8.3}
\begin{itemize}
 \item $\Htilde_{\nbigs}:=F_{\nbigs}^{-1}(\Htilde\cup\Xtilde^0)$
       is normal crossing.
 \item $F_{\nbigs}$ induces
       an isomorphism
       $\Xtilde_{\nbigs}\setminus F_{\nbigs}^{-1}(\Htilde^0)
       \simeq
       \Xtilde\setminus\Htilde^0$.
 \item The natural morphism
       $\Theta_{\Xtilde_{\nbigs}}(\log\Htilde_{\nbigs})
       \lrarr
       F_{\nbigs}^{\ast}
       \bigl(
       \Theta_{\Xtilde}(\log(\Htilde\cup\Xtilde^0))
       \bigr)$ is an isomorphism.
 \item For $1\leq i\neq j\leq \ell$,
       we obtain
       $\dim\bigl(
       F_{\nbigs}^{-1}(\Htilde_i)
       \cap
       F_{\nbigs}^{-1}(\Htilde_j)
       \bigr)\leq \dim X-2$.
\end{itemize}
\end{lem}
\pf
Let $U_1$ be an open neighbourhood of
$(0,0)$ in $\cnum^2=\{(x_1,x_2)\}$,
and let $U_2$ be an open neighbourhood of
$(0,\ldots,0)$ in $\cnum^{n-2}=\{(y_1,\ldots,y_{n-2})\}$.
We set
$U:=U_1\times U_2$,
$H_{x,i}:=\{x_i=0\}\times U_2$ $(i=1,2)$
and $H_{y,j}:=U_1\times\{y_j=0\}$ $(j=1,\ldots,\ell-2)$.
Let $\pi:\Utilde\lrarr U$ be the blowing up
along $H_{x,1}\cap H_{x,2}$.
Let $\Utilde'\subset \Utilde$ be the complement
of the proper transform of $H_{x,2}$.
On $\Utilde'$,
we set
$u_1:=x_1/x_2$ and $u_2:=x_2$.
Then, $(u_1,u_2,y_1,\ldots,y_{n-2})$
induces a holomorphic coordinate system
on $\Utilde'$.
We have
$\pi^{-1}(H_{x,2})\cap \Utilde'=\{u_2=0\}$
and
$\pi^{-1}(H_{y,j})\cap\Utilde'=\{y_j=0\}$.
The proper transform of $H_{x,1}$ is contained in $\Utilde'$,
and it is equal to $\{u_1=0\}$.

Because each step of the construction of $F_{\nbigs}$
is locally described as above,
we can prove the claims of the proposition
by an induction.
\hfill\qed

\begin{rem}
The construction of $F_{\nbigs}$
depends on the ordering
$z_1,\ldots,z_{\ell}$.
Namely, for any
automorphism $\sigma$ of
$\{1,\ldots,\ell\}$,
we set $y_{i}=z_{\sigma(i)}$ $(i=1,\ldots,\ell)$
and $y_{i}=z_i$ $(i=\ell+1,\ldots,n)$.
If $\sigma$ is not the identity map,
we obtain a different manifolds
$\Xtilde'_{\nbigs}$
with a morphism $F_{\nbigs}':\Xtilde'_{\nbigs}\lrarr \Xtilde$.
\hfill\qed
\end{rem}

\subsubsection{An auxiliary resolution associated with a good set of irregular values}
\label{subsection;21.3.15.1}

Let $\nbigi$ be a good set of irregular values on $(X,H)$.
(See \cite[Definition 2.1.2]{Mochizuki-wild}.)
For any positive integer $m$,
and for any $\Ptilde\in \Htilde$,
we obtain the induced tuple
\[
 \nbigitilde_{m,\Ptilde}=
 \bigl\{
  t^m\gminia\,\big|\,\gminia\in\nbigi
  \bigr\}
\subset
   \nbigo_{\Xtilde}\bigl(\ast(\Htilde\cup\Xtilde^0)\bigr)_{\Ptilde}
   \big/
   \nbigo_{\Xtilde,\Ptilde}.
\]
If $\Ptilde\in\Htilde\setminus\Htilde^0$,
it is a good set of irregular values at $\Ptilde$.
For any $\Ptilde\in\Xtilde^0\setminus\Htilde$,
we set
\[
 \nbigitilde_{m,\Ptilde}=\bigl\{0\bigr\}
 \subset
   \nbigo_{\Xtilde}(\ast \Xtilde^0)_{\Ptilde}
   \big/
   \nbigo_{\Xtilde,\Ptilde}.
\]
We obtain the family
$\vecnbigitilde_m=
(\nbigitilde_{m,\Ptilde}\,|\,\Ptilde\in \Htilde\cup\Xtilde^0)$.
If $m=1$,
we omit to denote the subscript $m$,
i.e.,
$\nbigitilde_{\Ptilde}=\nbigitilde_{1,\Ptilde}$
and
$\vecnbigitilde=\vecnbigitilde_1$.
Note that for any holomorphic map
$G_Y:Y\lrarr \Xtilde$
such that
$H_Y:=G_Y^{-1}(\Htilde\cup\Xtilde^0)$ is normal crossing,
we obtain the induced family of tuples
$G_Y^{\ast}(\vecnbigitilde_m)
=\bigl(
 G_Y^{\ast}(\vecnbigitilde_m)_Q\,\big|\,Q\in H_Y
 \bigr)$
where
$G_Y^{\ast}(\vecnbigitilde_m)_Q$
denotes
the image of
$\nbigitilde_{m,G_Y(Q)}$
via
$\nbigo_{\Xtilde}(\ast (\Htilde\cup\Xtilde^0))_{G_Y(Q)}
\big/
\nbigo_{\Xtilde,G_Y(Q)}
\lrarr
\nbigo_{Y}(\ast H_Y)_Q\big/
\nbigo_{Y,Q}$.

\vspace{.1in}
Let $\nbigs(\nbigi)\subset\seisuu_{\leq 0}^{\ell}\setminus\{(0,\ldots,0)\}$
denote the set of
$\ord(\gminia-\gminib)$
for distinct $\gminia,\gminib\in\nbigi$.
Because $\nbigi$ is assumed to be a good set of irregular values,
$\nbigs(\nbigi)$ is totally ordered with respect to $\leq_{\seisuu^{\ell}}$.

Suppose that
$\nbigs(\nbigi)\subset (m\seisuu_{\leq 0})^{\ell}$
for a positive integer $m$.
We obtain
the totally ordered set
$\nbigs(\nbigi)/m:=
\bigl\{
m^{-1}\vecn\,\big|\,
\vecn\in\nbigs
\bigr\}\subset \seisuu_{\leq 0}^{\ell}$.
We obtain the following lemma 
by the construction of 
the projective morphism of complex manifolds
$F_{\nbigs(\nbigi)/m}:\Xtilde_{\nbigs(\nbigi)/m}\lrarr \Xtilde$.

\begin{lem}
\label{lem;21.2.15.10}
 $F_{\nbigs(\nbigi)/m}^{\ast}(\vecnbigitilde_{m})$
 is a good set of irregular values
 on $(\Xtilde_{\nbigs(\nbigi)/m},\Htilde_{\nbigs(\nbigi)/m})$.
\hfill\qed
\end{lem}

\section{Strong regularity along a smooth hypersurface}

\subsection{Regularity along a smooth hypersurface}

\subsubsection{Regularity of holonomic $\nbigd$-modules
along a coordinate function}
\label{subsection;22.7.30.1}

Let $X$ be a complex manifold.
Let $U$ be an open subset of $X$.
A holomorphic function $f$ on $U$ is called
a coordinate function if $df$ is nowhere vanishing on $U$.
For any such $f$,
we obtain the subbundle
$\Ker(df)\subset\Theta_{U}$.
Let $\nbigd_{U/f}\subset
\nbigd_{U}$
denote the sheaf of subalgebras
generated by
$\Ker(df)$ over $\nbigo_{U}$.
We also obtain the subsheaf
$\Theta_{U}(\log f)\subset \Theta_{U}$
of holomorphic vector fields $v$ which are logarithmic
with respect to $f$,
i.e., $v(f)\in f\nbigo_{U}$.
Let 
$V_f\nbigd_{U}\subset\nbigd_{U}$
denote the sheaf of subalgebras generated by
$\Theta_{U}(\log f)$
over $\nbigo_{U}$.

If $f_1$ is another coordinate function on $U$
such that $f^{-1}(0)=f_1^{-1}(0)$,
we have
$\Theta_{U}(\log f)=\Theta_{U}(\log f_1)$
but $\Ker(df)\neq \Ker(df_1)$, in general.
Hence,
$V_f\nbigd_{U}$ depends only on $f^{-1}(0)$,
but 
$\nbigd_{U/f}$ is not determined by $f^{-1}(0)$.

Let $M$ be a holonomic $\nbigd_X$-module.
We fix a total order $\leq_{\cnum}$
such that
(i) the restriction of $\leq_{\cnum}$ to $\seisuu$
is equal to the standard order,
(ii) $a_1\leq_{\cnum}a_2$
implies that
$a_1+n\leq_{\cnum}a_2+n$
for any $n\in\seisuu$.
Then, $M_{|U}$ has a $V$-filtration
$V_{\bullet}(M_{|U})$ along $f$
indexed by $(\cnum,\leq_{\cnum})$,
i.e.,
(i) each $V_{a}(M_{|U})$
is $V_f\nbigd_U$-coherent,
(ii) $\bigcup V_a(M_{|U})=M_{|U}$
and $V_{a}(M_{|U})=\bigcap_{b>_{\cnum}a}V_{b}(M_{|U})$,
(iii) we have $fV_a(M_{|U})\subset V_{a-1}(M_{|U})$
for any $a\in\real$
and $fV_a(M_{|U})=V_{a-1}(M_{|U})$ for $a<_{\cnum}0$,
(iv) for a coordinate system
$(f,x_1,\ldots,x_{n-1})$ around any point of $U$,
we have $\del_fV_a(M_{|U})\subset V_{a+1}(M_{|U})$,
and $-\del_ff-a$ are locally nilpotent on $\Gr^V_{a}(M_{|U})$.
We recall that
$M$ is called regular along $f$
if each $V_a(M_{|U})$ is $\nbigd_{U/f}$-coherent.
(See \cite[\S3.1.d]{Sabbah-pure-twistor}.)

\begin{lem}
$M$ is regular along $f$
if and only if $M_{|U}(\ast f)$ is regular along $f$.
\end{lem}
\pf
We may assume that $X=U$.
There exists the natural morphism
$\rho:M\lrarr M(\ast f)$.
The supports of $\Ker(\rho)$ and
$\Cok(\rho)$ are contained in $f^{-1}(0)$.
There exists the exact sequence
\[
 0\lrarr
 V_a(\Ker\rho)
 \lrarr
 V_a(M)\lrarr V_a(M(\ast f))
 \lrarr V_a(\Cok\rho)\lrarr 0.
\]
It is easy to see that
$V_a(\Ker\rho)$
and $V_a(\Cok\rho)$
are $\nbigd_{U/f}$-coherent.
Hence, the claim of the lemma follows.
\hfill\qed

\begin{lem}
\label{lem;21.3.9.30}
 $M$ is regular along $f$
if and only if
there exists a
$V\nbigd_{U}$-submodule
$L\subset M_{|U}$
such that
(i) $L(\ast f)=M_{|U}(\ast f)$,
(ii) $L$ is coherent over
 $\nbigd_{U/f}$.
\end{lem}
\pf
We may assume that $X=U$ and $M=M(\ast f)$.
The only if part of the claim is clear.
Assume that there exists 
a $V\nbigd_{U}$-submodule $L\subset M$
such that
(i) $L(\ast f)=M$,
(ii) $L$ is coherent over $\nbigd_{X/f}$.
Let us prove that
$V_a(M)$ is
$\nbigd_{X/f}$-coherent for each $a\in\cnum$.
Because it is enough to prove the claim locally around any point
of $X$,
we may assume that
$V_a(M)$ is finitely generated over
$V_f\nbigd_X$.
There exists an integer $\ell$ such that
$V_a(M)\subset f^{-\ell}L$.
Note that $f^{-\ell}L$ is coherent over
$\nbigd_{X/f}$,
and that the sheaf of algebras
$\nbigd_{X/f}$ is Noetherian.
Because
$V_a(M)$ is the sum of
coherent $\nbigd_{X/f}$-submodules,
we obtain that $V_a(M)$ is coherent over
$\nbigd_{X/f}$.
(See \cite[Theorem A.29]{kashiwara_text}.)

\hfill\qed

\begin{lem}
\label{lem;21.3.9.20}
 Let $0\lrarr M'\lrarr M\lrarr M''\lrarr 0$
be an exact sequence of holonomic $\nbigd_{X}$-modules.
Then, $M$ is regular along $f$
if and only if $M'$ and $M''$ are regular along $f$.
\end{lem}
\pf
It follows from the exactness of
$0\lrarr V_a(M'_{m,U_X})\lrarr
 V_{a}(M_{m,U_X})\lrarr
  V_a(M''_{m,U_X})\lrarr 0$.
\hfill\qed

\begin{cor}
\label{cor;21.3.9.21}
 Let $g:M_1\lrarr M_2$ be a morphism of holonomic
$\nbigd_{X}$-modules.
If $M_i$ are regular along $f$,
 $\Ker(g)$, $\Image(g)$ and $\Cok(g)$
 are also regular along $f$. 
\hfill\qed
\end{cor}

Let $\iota_f:U\lrarr U\times\cnum_t$ denote the embedding
defined by $\iota_f(x)=(x,f(x))$.

\begin{lem}
$M$ is regular along $f$
if and only if
$\iota_{f\dagger}(M_{|U})$  is regular along $t$.
\end{lem}
\pf
We may assume that $X=U$.
We may also assume that $M(\ast f)=M$.
It is enough to prove the claim locally around any point of $f^{-1}(0)$.
Hence, we may assume that $X$ is the product of
a complex manifold $X_0$
and a neighbourhood $Z$ of $0$ in $\cnum_z$,
and that $f$ is the projection onto $Z$,
i.e., $f=z$.
Let $\pi_0$ denote the projection of $X$ onto $X_0$.
Let $\pi$ denote the projection of $X\times\cnum_t$
onto $X$.

Recall that
$\iota_{z\dagger}(M)
=\iota_{z\ast}(M)\otimes\cnum[\del_t]$
on which we have the following formula:
\[
 P(m\otimes \del_t^j)=(Pm)\otimes \del_t^j
 \quad (P\in \pi_0^{\ast}\nbigd_{X_0}),
 \quad
 \del_z(m\otimes \del_t^j)=
 (\del_zm)\otimes \del_t^j
 -m\otimes\del_t^{j+1},
\]
\[
 t(m\otimes\del_t^j)=(zm)\otimes\del_t^j
 -jm\otimes\del_t^{j-1},
 \quad
 \del_t(m\otimes\del_t^j)=m\otimes\del_t^{j+1}.
\]
We have
$(\del_tt+\del_zz)(m\otimes\del_t^j)
=(\del_zzm)\otimes\del_t^j-jm\otimes\del_t^j$.
We note that $(t-z)(m\otimes 1)=0$.
For $j\geq 1$, the following holds.
\begin{equation}
\label{eq;21.6.29.3}
 (t-z)\del_z\Bigl(
 \del_z^j(m\otimes 1)
 \Bigr)
 =(j+1)\del_z^j(m\otimes 1).
\end{equation}
Let $V_{\bullet}(M)$ denote the $V$-filtration of $M$ along $z$.
We set
$V_{c}(\iota_{z\dagger}M)
=\pi^{\ast}\nbigd_X\Bigl(
\iota_{z\ast}V_c(M)\otimes 1
\Bigr)$
for any $c\in\real$.
Then, it is easy to see that
$V_{\bullet}(\iota_{z\dagger}M)$
is the $V$-filtration of
$\iota_{z\dagger}M$ along $t$.
 
Suppose that $\iota_{z\dagger}M$ is regular along $t$.
Let us prove that $M$ is regular along $z$.
It is enough to study it around any point of $z^{-1}(0)$.
Hence, we may assume that
$V_c(\iota_{z\dagger}M)$
is generated by
$a_j\otimes 1$ $(j=1,\ldots,N)$
over $\pi^{\ast}\nbigd_{X}$.
For any $m\in V_c(M)$,
there exist $N'\in\seisuu_{>0}$ and
$P_{k,j}\in \pi_0^{\ast}\nbigd_{X_0}$
$(1\leq j\leq N,\,\,0\leq k\leq N')$
such that
\[
 m\otimes 1=
 \sum_{k,j}P_{k,j}\del_z^k(a_j\otimes 1).
\]
By using (\ref{eq;21.6.29.3}),
we obtain
\[
 m\otimes 1
 =\sum_{j}(P_{0,j}a_j)\otimes 1.
\]
We obtain that
$V_c(M)$ is generated by $a_j$.
Hence, we obtain the regularity of $M$.
The converse is also easy to prove.
\hfill\qed

\subsubsection{Regularity along smooth hypersurfaces}

Let $H$ be a smooth hypersurface of $X$.
Let $U$ be an open subset of $X$.
A holomorphic function $f$ on $U$ is called
a coordinate function defining $H\cap U$
if $f$ is a coordinate function on $U$ such that
$f^{-1}(0)=H\cap U$.

A holonomic $\nbigd_X$-module $M$ is called
regular along $H$
if the following holds.
\begin{itemize}
 \item For any open subset $U\subset X$
       and any coordinate function $f$ on $U$ defining $H\cap U$,
       $M$ is regular along $f$.     
\end{itemize}

\begin{lem}
\mbox{{}}\label{lem;21.6.29.10}
\begin{itemize}
 \item   
$M$  is regular along $H$ if and only if $M(\ast H)$ is regular along $H$.
 \item Let $0\lrarr M'\lrarr M\lrarr M''\lrarr 0$
       be an exact sequence of holonomic $\nbigd_{X}$-modules.
       Then, $M$ is regular along $H$
       if and only if $M'$ and $M''$ are regular along $H$.
 \item Let $g:M_1\lrarr M_2$ be a morphism of
       holonomic $\nbigd_{X}$-modules.
       If $M_i$ are regular along $H$,
       $\Ker(g)$, $\Image(g)$ and $\Cok(g)$
       are also regular along $H$.       \hfill\qed
\end{itemize}
 \end{lem}

Let $F:X\lrarr Y$ be any proper morphism of complex manifolds.
Let $H_Y$ be a smooth hypersurface.
Assume that $H_X=X\times_YH_Y$ is smooth and reduced,
i.e., $F$ is transversal with $H_Y$.

\begin{lem}
\label{lem;21.3.9.22}
Let $M$ be a holonomic $\nbigd_{X}$-module.
If $M$ is regular along $H_X$,
$F_{\dagger}^i(M)$ are regular along $H_Y$.
\end{lem}
\pf
We have only to prove the claim locally
around any point of $H_Y$.
Hence,
we may assume that
$Y$ is the product of a complex manifold $Y_0$
and a neighbourhood $Z$ of $0$ in $\cnum_z$,
and that $H_Y=Y\times \{0\}$.
By the assumption,
$f_X:=F^{\ast}(z)$ is a coordinate function
defining $H_X$.
Let $\iota_z:Y\lrarr Y\times\cnum_t$
and $\iota_{f_X}:X\lrarr X\times\cnum_t$
denote the graph embeddings.
Then, $M$ is regular along $f_X$
if and only if $\iota_{f_X\dagger}M$ is regular along $t$,
and $F_{\dagger}^j(M)$ are regular along $z$
if and only if
$(F\times\id_{\cnum_t})_{\dagger}^j(\iota_{f_X\dagger}M)$
is regular along $t$.
Hence, we may assume that
(i) $X=X_1\times\cnum_t$
and $Y=Y_1\times\cnum_t$,
(ii) $F$ is induced by a morphism $F_1:X_1\lrarr Y_1$
and the identity map of $\cnum_t$,
(iii) $H_Y=Y_1\times\{0\}$ and $H_X=X_1\times\{0\}$.

Let $p_X:X\lrarr X_1$ and $p_Y:Y\lrarr Y_1$
denote the projections.
There exist natural isomorphisms
$\nbigd_{X/t}=p_X^{\ast}(\nbigd_{X_1})$
and
$\nbigd_{Y/t}=p_Y^{\ast}(\nbigd_{Y_1})$.
There exist natural isomorphisms
(see \cite{Mebkhout-Sabbah}):
\[
 V_a\Bigl(
 F^i_{\dagger}(M)
 \Bigr)
\simeq R^iF_{\ast}\Bigl(
 p_X^{\ast}\nbigd_{Y_1\larr X_1}\otimes^{L}_{p_X^{\ast}(\nbigd_{X_1})}
V_a(M)
 \Bigr).
\]
Because $V_{a}(M)$ are
coherent over $\nbigd_{X/t}$,
we can prove that
$V_aF^i_{\dagger}(M)$ are
$\nbigd_{Y/t}$-coherent
by the standard argument in the theory
of $\nbigd$-modules.
(See \cite[Theorem 4.25]{kashiwara_text}.)
\hfill\qed

\subsection{Strong regularity of holonomic $\nbigd$-modules
along a smooth hypersurface}

Let $\mu$ be a standard coordinate of $\cnum$.
For any positive integer $m$,
let $\varphi_m:\cnum\lrarr \cnum$ be the morphism
defined by $\varphi_m(\mu)=\mu^m$.

Let $X$ be a complex manifold.
Let $H$ be a smooth hypersurface of $X$.
Let $U$ be an open subset of $X$
with a coordinate function $f$ defining $H\cap U$.
We define $U_{f,m}$ as the fiber product of
$\varphi_m:\cnum\lrarr\cnum$ and $f:U\lrarr \cnum$.
Note that $U_{f,m}$ is naturally a complex manifold.
Let $\varphi_{f,m}:U_{f,m}\lrarr U$
denote the induced morphism.
We set $H_{f,m}:=\varphi_{f,m}^{-1}(H\cap U)$.

\begin{lem}
Let $M$ be a holonomic $\nbigd_X$-module.
Let $U$ be an open subset of $X$
with coordinate functions $f_i$ $(i=1,2)$
defining $H\cap U$.
We obtain the induced $\nbigd_{U_{f_i,m}}$-modules
$M_{f_i,m}:=(\varphi_{f_i,m})^{\ast}\bigl(M_{|U}(\ast f_i)\bigr)$.
Then,
$M_{f_1,m}$ is regular along $H_{f_1,m}$
if and only if
$M_{f_2,m}$ is regular along $H_{f_2,m}$.
\end{lem}
\pf
It is enough to prove the claim locally around any point of $U$.
We may assume that $U$ is simply connected.
We set $h:=f_2/f_1$,
which is a nowhere vanishing holomorphic function on $U$.
There exists a holomorphic function $h^{1/m}$ on $U$
such that $(h^{1/m})^m=h$.

We set $g_1=f_1\circ\varphi_{f_1,m}:U_{f_1,m}\lrarr \cnum$.
Let $f_1^{1/m}:U_{f_1,m}\lrarr \cnum$
be the morphism obtained as the projection,
i.e.,
$\varphi_m\circ f_1^{1/m}=f_1$.
By the construction,
we have $g_1=(f_1^{1/m})^m$.

We set $g_2:=f_2\circ\varphi_{f_1,m}:U_{f_1,m}\lrarr \cnum$.
We also obtain holomorphic functions
$h_0:=h\circ\varphi_{f_1,m}$
and $h_0^{1/m}:=h^{1/m}\circ\varphi_{f_1,m}$.
We have
$g_2=h_0\cdot g_1=(h_0^{1/m}\cdot g_1^{1/m})^m$.
We define $k:U_{f_1,m}\lrarr \cnum$
by $k_m=h_0^{1/m}\cdot g_1^{1/m}$.
Then, we obtain
$f_2\circ\varphi_{f_1,m}
=\varphi_{m}\circ k_m$.
Hence, $\varphi_{f_1,m}$ and $k_m$
induce a morphism
$\rho_{1.m}:U_{f_1,m}\lrarr U_{f_2,m}$.
Similarly, we obtain a morphism
$\rho_{2,m}:U_{f_2,m}\lrarr U_{f_1,m}$.
The morphisms
$\rho_{1,m}$ and $\rho_{2,m}$ are mutually inverse,
and we have
$\rho_{1,m}^{-1}(H_{f_2,m})=H_{f_1,m}$.
Then, the claim of the lemma is clear.
\hfill\qed

\begin{df}
 \label{df;21.4.14.5}
We say that a holonomic $\nbigd_{X}$-module $M$
is strongly regular along $H$
if the following condition is satisfied 
for any open subset $U\subset X$,
any coordinate function $f$ on $U$
defining $H\cap U$
and a positive integer $m$.
\begin{itemize}
 \item The induced $\nbigd_{U_{f,m}}$-module
       $M_{f,m}$
       is regular along $H_{f,m}$.
\end{itemize} 
If there exists a holomorphic function $g$ on $X$
such that $H=g^{-1}(0)$,
we also say that $M$ is strongly regular along $g$.
\hfill\qed
\end{df}

By definition,
$M$ is strongly regular along $\mu$
if and only if $M(\ast \mu)$ is strongly regular along $\mu$.
We obtain the following lemma from Lemma \ref{lem;21.6.29.10}.
\begin{lem}
\mbox{{}}\label{lem;21.3.9.33}
\begin{itemize}
 \item 
 Let $0\lrarr M'\lrarr M\lrarr M''\lrarr 0$
be an exact sequence of holonomic $\nbigd_{X}$-modules.
Then, $M$ is strongly regular along $H$
if and only if $M'$ and $M''$ are strongly regular along $H$.
\item
 Let $g:M_1\lrarr M_2$ be a morphism of holonomic
$\nbigd_{Z\times X}$-modules.
If $M_i$ are strongly regular along $\mu$,
 $\Ker(g)$, $\Image(g)$ and $\Cok(g)$
 are also strongly regular along $\mu$. 
     \hfill\qed
\end{itemize}
\end{lem}

Let $F:X\lrarr Y$ be any proper morphism of complex manifolds.
Let $H_Y$ be a smooth hypersurface of $Y$
such that $F$ is transversal with $H_Y$.
We obtain the smooth hypersurface $H_X:=F^{-1}(H_Y)$.

\begin{prop}
\label{prop;21.1.23.2}
Let $M$ be a holonomic $\nbigd_{X}$-module
which is strongly regular along $H_X$.
Then,
$F_{\dagger}^i(M)$ are strongly regular along $H_Y$.
\end{prop}
\pf
It is enough to prove the claim locally around any point of $Y$.
Hence, we may assume that $H_Y$ is defined by
a coordinate function $f_Y$ of $Y$.
We set $f_X:=F^{\ast}(f_Y)$.
Take any $m\in\seisuu_{>0}$.
We obtain
$\varphi_{f_Y,m}\colon Y_{f_Y,m}\lrarr Y$
and
$\varphi_{f_X,m}\colon X_{f_X,m}\lrarr X$.
There exists the induced morphism
$F_{m}:X_{f_{X},m}\lrarr Y_{f_Y,m}$.
There exists a natural isomorphism
\[
 \varphi_{f_Y,m}^{\ast}\bigl(
 F^j_{\dagger}(M)(\ast H_Y)
 \bigr)
\simeq
 F_{m\dagger}^j\bigl(
 \varphi_{f_X,m}^{\ast}M(\ast H_X)
 \bigr)
\]
Then, the claim follows from Lemma \ref{lem;21.3.9.22}.
\hfill\qed

\subsubsection{Regular singularity of meromorphic flat bundles}

\begin{prop}
\label{prop;21.2.16.10}
Suppose that $M$ is a meromorphic flat bundle
on $(Z\times X,\{0\}\times X)$.
Then, $M$ is strongly regular along $\mu$
if and only if
$M$ is a regular singular meromorphic flat bundle
on $(Z\times X,\{0\}\times X)$.
\end{prop}
\pf
The ``if'' part of the claim is clear.
Let us prove that the ``only if'' part.
Suppose that $M$ is not regular singular,
and we shall deduce a contradiction.
According to Malgrange \cite{malgrange}
(see also \cite[Proposition 2.7.6]{Mochizuki-wild}),
there exists a closed complex analytic subset
$A\subset \{0\}\times X$
with $\dim A<\dim X$
and that
$M_{|(Z\times X)\setminus A}$
is a good meromorphic flat bundle
on $((Z\times X)\setminus A,(\{0\}\times X)\setminus A)$.
(See \cite[\S2]{Mochizuki-Stokes} for a survey of
good meromorphic flat bundles.)
For any $(0,P)\in \{0\}\times X$,
there exist a positive integer $m$
and a neighbourhood $X_P$ of $P$ in $X$
such that
$(\varphi_m\times\id)^{\ast}(M)$
is an unramifiedly good meromorphic flat bundle
on $(Z^{(m)}\times X_P,\{0\}\times X_P)$.
There exist
a good set of irregular values
$\nbigi\subset \mu^{-1}\nbigo_{X_P}[\mu^{-1}]$
and the decomposition
\[
 (\varphi_m\times\id)^{\ast}(M,\nabla)_{|\widehat{\{0\}\times X_P}}
 =\bigoplus_{\gminia\in\nbigi}
 (\Mhat_{P,\gminia},\nabla_{\gminia})
\]
such that
$\nabla_{\gminia}-d\gminia\id$ are regular singular.
Note that
\[
 V_a\Bigl(
 (\varphi_m\times\id)^{\ast}(M)_{|\widehat{\{0\}\times X_P}}
 \Bigr)
 =\bigoplus_{\gminia\in\nbigi}
 V_a(\Mhat_{P,\gminia}).
\]
If $\gminia\neq 0$,
we obtain $V_a(\Mhat_{P,\gminia})=\Mhat_{P,\gminia}$ for any $a\in\cnum$.

If $M$ is not regular singular around $(0,P)$,
there exists a non-zero $\gminia\in\nbigi$.
By shrinking $X_P$ and $Z$,
there exists a coordinate function $f_m$
defining $\{0\}\times X_P$ on $Z^{(m)}\times X_P$
such that $\gminia=f_m^{-\ell}$ for a positive integer $\ell$.
It is easy to see that
$V_a(\Mhat_{P,\gminia})$ is not coherent
over $\nbigd_{(Z^{(m)}\times X_P)/f_m|\widehat{\{0\}\times X_P}}$.
Hence,
$V_a
 (\varphi_m\times\id)^{\ast}(M)$
is not coherent over $\nbigd_{(Z^{m}\times X_P)/f_m}$.
Therefore, we can conclude that
$M$ is regular singular around
any point of $(\{0\}\times X)\setminus A$.

There exists a regular singular meromorphic flat bundle
$M'$ on $(Z\times X,0\times X)$
with an isomorphism $M'_{|(Z\setminus\{0\})\times X}
\simeq M_{|(Z\setminus\{0\})\times X}$
which extends to an isomorphism
$M'_{|(Z\times X)\setminus A}
\simeq M_{|(Z\times X)\setminus A}$.
By the Hartogs theorem,
we obtain
$M'\simeq M$.
Hence, we obtain that $M$ is a regular singular
meromorphic flat bundle on $(Z\times X,\{0\}\times X)$.
\hfill\qed

\subsection{Strongly regular extensions of some $\nbigd$-modules}

We use the notation in \S\ref{subsection;21.1.22.1}.
Let $(\nbigv,\nabla)$ be a meromorphic flat bundle on
$(\Xtilde\setminus\Xtilde^0,\Htilde\setminus\Htilde^0)$
which is good over $\vecnbigitilde_{|\Xtilde\setminus\Xtilde^0}$.
\begin{prop}
\label{prop;21.1.23.10}
 There exists a coherent reflexive $\nbigo_{\Xtilde}(\ast \Htilde)$-module
 $\nbigvtilde$ with an integrable connection $\nablatilde$
 such that 
 (i) $(\nbigvtilde,\nablatilde)_{|\Xtilde\setminus\Xtilde^0}
      \simeq(\nbigv,\nabla)$,
 (ii) $(\nbigvtilde,\nablatilde)_{|\Xtilde\setminus \Htilde}$
 is regular singular.
Such $(\nbigvtilde,\nablatilde)$ is unique up to canonical isomorphisms.
Moreover,
$F_{\nbigs(\nbigi)}^{\ast}(\nbigvtilde,\nablatilde)$
on $(\Xtilde_{\nbigs(\nbigi)},\Htilde_{\nbigs(\nbigi)})$ is good over
$F_{\nbigs(\nbigi)}^{\ast}(\vecnbigitilde)$.
\end{prop}
\pf
The uniqueness of such $(\nbigvtilde,\nablatilde)$
follows from
the uniqueness of regular singular extension and
the Hartogs theorem.
Let us study the existence
in the case $\ell=1$.
We may naturally regard
$\Xtilde\setminus \Xtilde^0$
as an open subset of $\Xtilde_{\nbigs(\nbigi)}$.
By using \cite[Proposition 14.4.9]{Mochizuki-MTM},
we extend
$\nbigv$ to a good meromorphic flat bundle
$\nbigvtilde'$ on $(\Xtilde_{\nbigs(\nbigi)},\Htilde_{\nbigs(\nbigi)})$
over $F_{\nbigs(\nbigi)}^{\ast}(\vecnbigitilde)$.
Then, $\nbigvtilde=(F_{\nbigs(\nbigi)})_{\ast}(\nbigvtilde')$
is the desired extension.
Let us study the case $\ell\geq 2$.
By the existence in the case $\ell=1$,
we obtain $\nbigvtilde_1$
on $\bigl(\Xtilde\setminus \Htilde_{[2]}^0,(\Htilde\cup\Xtilde^0)
\setminus\Htilde_{[2]}^0\bigr)$
with the desired property.
There exists the Deligne-Malgrange lattice
$\nbigvtilde_1^{\DM}$ of $\nbigvtilde_1$.
Because the codimension of
$\Htilde_{[2]}^0$ in $\Xtilde$ is $3$,
$\nbigvtilde_1^{\DM}$ extends to
a coherent reflexive $\nbigo_{\Xtilde}$-module
$\nbigvtilde^{\DM}$
by an extension theorem of Siu \cite{Siu-extension}.
Then, $\nbigvtilde=
 \nbigvtilde^{\DM}\otimes
 \nbigo_{\Xtilde}(\ast(\Htilde\cup\Xtilde^0))$
is the desired one.
By the construction in the case $\ell=1$,
the restriction of
$F_{\nbigs(\nbigi)}^{\ast}(\nbigvtilde)$
to $\Xtilde_{\nbigs(\nbigi)}\setminus
F_{\nbigs(\nbigi)}^{-1}(\Htilde_{[2]}^0)$
is unramifiedly good
whose good system of irregular values
is contained in 
$F_{\nbigs(\nbigi)}^{\ast}(\vecnbigitilde)
 _{|\Xtilde_{\nbigs(\nbigi)}\setminus F_{\nbigs(\nbigi)}^{-1}(\Htilde_{[2]}^0)}$.
We also note that
$\dim F_{\nbigi}^{-1}(\Htilde_{[2]})\leq \dim \Xtilde-2$.
Then, the last claim follows from
Proposition \ref{prop;21.1.22.2}.
\hfill\qed

\vspace{.1in}

Let $\{1,\ldots,\ell\}=I\sqcup J$ be a decomposition.
Let $\Htilde(I)=\bigcup_{i\in I}\Htilde_i$
and $\Htilde(J)=\bigcup_{i\in J}\Htilde_i$.
We shall prove the following proposition
in \S\ref{subsection;21.3.9.1}.

\begin{prop}
 \label{prop;21.1.23.11}
The holonomic $\nbigd_{\Xtilde}$-module
$\bigl(
\nbigvtilde(!\Htilde(J))
\bigr)(\ast \Xtilde^0)$
is strongly regular along $t$.
\end{prop}

We also have the following complement.
\begin{lem}
\label{lem;21.6.26.1}
 Let $g$ be a holomorphic function on $X$ such that $g^{-1}(0)\subset H$.
Let $\gtilde$ be the induced holomorphic function on $\Xtilde$.
Let $\iota_{\gtilde}:
\Xtilde\lrarr \Xtilde\times \cnum_s$
 denote the graph embedding.
Let
$\lefttop{s}V\nbigd_{\Xtilde\times\cnum_s/t}(\ast \Xtilde^0)
\subset\nbigd_{\Xtilde\times\cnum_s}(\ast \Xtilde^0)$
denote the sheaf of subalgebras
generated by $\del_{z_i}$ and $s\del_s$
over $\nbigo_{\Xtilde\times\cnum_s}(\ast \Xtilde^0)$.
Let $V_{\bullet}\Bigl(
\iota_{\gtilde\dagger}
\bigl(
\nbigvtilde(!\Htilde(J))
\bigr)(\ast \Xtilde^0)
\Bigr)$
denote the $V$-filtration of 
$\iota_{\gtilde\dagger}
\bigl(
\nbigvtilde(!\Htilde(J))
\bigr)(\ast \Xtilde^0)$
along $s$.
Then,
$V_{a}\Bigl(
\iota_{\gtilde\dagger}
\bigl(
\nbigvtilde(!\Htilde(J))
\bigr)(\ast \Xtilde^0)
\Bigr)(\ast \Xtilde^0)$
are coherent over
$\lefttop{s}V\nbigd_{\Xtilde\times\cnum_s/t}(\ast \Xtilde^0)$.
\end{lem}
\pf
We have the formula for 
the $V$-filtration of
$\iota_{\gtilde\dagger}
\bigl(
\nbigvtilde(!\Htilde(J))
\bigr)(\ast \Xtilde_0)_{|(\Xtilde\setminus\Xtilde^0)\times\cnum_s}$
as in 
\cite[(5.20), (5.21)]{Mochizuki-MTM},
according to which
they are coherent over
$\lefttop{s}V\nbigd_{(\Xtilde\setminus\Xtilde^0)\times\cnum_s/t}$.
It naturally extends to
coherent
$V\nbigd_{\Xtilde\times\cnum_s/t}(\ast \Xtilde^0)$-submodules
of 
$\iota_{\gtilde\dagger}
\bigl(
\nbigvtilde(!\Htilde(J))
\bigr)(\ast \Xtilde_0)$.
It is easy to check that
they are equal to 
$V_{\bullet}\Bigl(
\iota_{\gtilde\dagger}
\bigl(
\nbigvtilde(!\Htilde(J))
\bigr)(\ast \Xtilde_0)
\Bigr)(\ast \Xtilde^0)$.

\hfill\qed

\vspace{.1in}

Let $(\nbigv_i,\nabla)$ $(i=1,2)$ be 
meromorphic flat bundles
on $(\Xtilde\setminus\Xtilde^0,\Htilde\setminus\Htilde^0)$
which are good over $\vecnbigitilde_{|\Xtilde\setminus\Xtilde^0}$.
Let $(\nbigvtilde_i,\nablatilde)$
be the extension as in Proposition \ref{prop;21.1.23.10}.

\begin{prop}
\label{prop;21.2.17.3}
Any morphism
$g:\nbigv_1\lrarr\nbigv_2$
of $\nbigd_{\Xtilde\setminus\Xtilde^0}$-modules
 uniquely extends to
a morphism 
 $\nbigvtilde_1\lrarr\nbigvtilde_2$
 of $\nbigd_{\Xtilde}$-modules.
\end{prop}
\pf
The uniqueness is clear.
By the regularity of $\nbigv_i$
along $\Xtilde^0\setminus\Htilde^0$,
$g$ extends to a morphism
$\nbigv_{1|\Xtilde\setminus\Htilde^0}
\lrarr
\nbigv_{2|\Xtilde\setminus\Htilde^0}$.
Let $\nbigv^{\DM}_{i}$ denote the Deligne-Malgrange lattice
of $\nbigv_i$.
We obtain
$\nbigv^{\DM}_{1|\Xtilde\setminus\Htilde^0}
\lrarr
\nbigv^{\DM}_{2|\Xtilde\setminus\Htilde^0}$
induced by $g$.
By the Hartogs theorem,
it extends to a morphism
$\nbigv^{\DM}_1
\lrarr
\nbigv^{\DM}_2$,
and hence
$g$ extends to
$\nbigv_1\lrarr\nbigv_2$.
\hfill\qed

\vspace{.1in}
Let us mention an easy example.
Let $(V_0,\nabla)$ be a regular singular meromorphic flat bundle on $(X,H)$.
Let $f$ be a meromorphic function on $(X,H)$.
For simplicity,
we assume that the zero and the pole of $f$
does not intersect.
Let $p:\Xtilde\lrarr X$ denote the projection.
We obtain  the meromorphic flat bundle
$\nbigvtilde:=p^{\ast}(V_0)(\ast\Xtilde^0)$
with $\nablatilde:=p^{\ast}(\nabla)+d(tf)$
on $(\Xtilde,\Xtilde^0\cup\Htilde)$.
The following is a corollary of Proposition \ref{prop;21.1.23.11}.
\begin{cor}
$\bigl(
 \nbigvtilde(!\Htilde(J))\bigr)(\ast \Xtilde^0)$
 is strongly regular along $t=0$.
\hfill\qed
\end{cor}

\subsubsection{Ramified coverings}

Let $m\in\seisuu_{>0}$.
We set $\Xtilde_m:=\cnum_{t_m}\times X$.
There exists the ramified covering
$\varphi_m:\Xtilde_m\lrarr \Xtilde$
induced by $t_m\longmapsto t_m^m=t$.
Let $\rho'_m:\cnum^n=\{(z_{m,1},\ldots,z_{m,\ell},z_{\ell+1},\ldots,z_n)\in\cnum^n\}
\lrarr \cnum^n=\{(z_1,\ldots,z_n)\in\cnum^n\}$
denote the ramified covering defined by
$\rho'_m(z_{m,1},\ldots,z_{m,\ell},z_{\ell+1},\ldots,z_n)
=(z_{m,1}^m,\ldots,z_{m,\ell}^m,z_{\ell+1},\ldots,z_n)$.
We set $Y_m:=(\rho_m')^{-1}(X)$.
The induced map $Y_m\lrarr X$ is also denoted by $\rho_m'$.
We set $\Ytilde_m:=\cnum_{t_m}\times Y_m$,
and let $\rho_m$ denote 
the induced ramified covering $\Ytilde_m\lrarr \Xtilde_m$.
The composition
$\kappa_m:=\varphi_m\circ\rho_m:
\Ytilde_{m}\lrarr \Xtilde$ is the ramified covering
induced by
\[
\kappa_m
(t_m,z_{m,1},\ldots,z_{m,\ell},z_{\ell+1},\ldots,z_n)
=(t_m^m,z_{m,1}^m,\ldots,z_{m,\ell}^m,z_{\ell+1},\ldots,z_n).
\]
We set $\Ytilde_m^0:=\{t_m=0\}\cap \Ytilde_m$,
$\Htilde_{Y,m,i}:=\{z_{m,i}=0\}\cap \Ytilde_m$
and $\Htilde_{Y,m}:=\bigcup_{i=1}^{\ell}\Htilde_{Y,m,i}$.
We obtain the meromorphic flat bundle
$\nbigvtilde_{Y,m}=
\kappa_m^{\ast}(\nbigvtilde)$
on $(\Ytilde_m,\Htilde_{Y,m}\cup\Ytilde_m^0)$.
We set
$\Htilde_{Y,m}(I):=\bigcup_{i\in I}\Htilde_{Y,m,i}$
and
$\Htilde_{Y,m}(J):=\bigcup_{i\in J}\Htilde_{Y,m,i}$.
We also set
\[
 \nbigm_{Y,m}(I,J):=
\Bigl(
 \nbigvtilde_{Y,m}(! \Htilde_{Y,m}(J))
\Bigr)(\ast \Ytilde_m^0).
\]

We set
$\nbigi_{Y,m}:=(\rho_m')^{\ast}(\nbigi)$
and
$\vecnbigitilde_{Y,m}:=
\kappa_m^{\ast}(\vecnbigitilde)$.
Note that
$(\vecnbigitilde_{Y,m})_{\Ptilde}
=\bigl\{
 t_m^m\gminib\,\big|\,
 \gminib\in\nbigi_{Y,m}
 \bigr\}$
 if $\Ptilde\in\Htilde_{Y,m}$,
 and 
 $(\vecnbigitilde_{Y,m})_{\Ptilde}
 =\{0\}$
if $\Ptilde\in\Ytilde_m^0\setminus\Htilde_{Y,m}$.
Moreover,
$\nbigs(\nbigi_{Y,m})
=\{m\cdot\vecn\,|\,\vecn\in\nbigs(\nbigi)\}$.

\subsubsection{Good filtered meromorphic flat bundles on the resolution}

By setting
$\Ytilde'_{m}:=
 (\Ytilde_{m})_{\nbigs(\nbigi_{Y,m})/m}$
and 
$F:=
 F_{\nbigs(\nbigi_{Y,m})/m}$,
we obtain a projective morphism of complex manifolds
$F:\Ytilde_{m}'\lrarr \Ytilde_m$.
We also set
$\Htilde'_{Y,m}:=
F^{-1}(\Htilde_{Y,m}\cup \Ytilde_m^0)$.
We obtain the meromorphic flat bundle
$\nbigvtilde'_{Y,m}:=
F^{\ast}(\nbigvtilde_{Y,m})$
on $(\Ytilde'_{m},\Htilde'_{Y,m})$.

\begin{lem}
$\nbigvtilde'_{Y,m}$
is an unramifiedly good meromorphic flat bundle
 on $(\Ytilde_{m}',\Htilde_{Y,m}')$
 over $F^{\ast}(\vecnbigitilde_{Y,m})$.
\end{lem}
\pf
Because there exists a morphism
$\kappa_{m}':\Ytilde_{m}'\lrarr \Xtilde_{\nbigs(\nbigi)}$
such that
$F_{\nbigs(\nbigi)}\circ\kappa_{m}'
=\kappa_{m}\circ F$,
the claim of the lemma follows.
\hfill\qed

\vspace{.1in}

Let $\Htilde'_{Y,m}=
\bigcup_{j\in\Lambda(\nbigi)}\Htilde'_{Y,m,j}$
denote the irreducible decomposition.
We obtain the associated good Deligne-Malgrange filtered bundle
$\nbigp^{\DM}_{\ast}\nbigvtilde'_{Y,m}$
on $(\Ytilde'_{m},\Htilde'_{Y,m})$
indexed by
$\real^{\Lambda(\nbigi)}$
(see \cite[\S2.3]{Mochizuki-DM-lattice} or
\cite[\S2.7.1]{Mochizuki-wild}).

\subsubsection{Submodules}

For $i=1,\ldots,\ell$,
let $\Lambda(\nbigi,i)\subset\Lambda(\nbigi)$
be the subset determined by
$F^{-1}(\Htilde_{Y,m,i})
=\bigcup_{j\in\Lambda(\nbigi,i)}
 \Htilde'_{Y,m,j}$.
There exists $[t_m]\in \Lambda(\nbigi)$
such that
$\Htilde'_{Y,m,[t_m]}$
is the proper transform of $\Ytilde^0_m$.
We have
$\Lambda(\nbigi)
=\{[t_m]\}\sqcup
\bigsqcup_{i=1}^{\ell}\Lambda(\nbigi,i)$.
Let $\epsilon$ be a sufficiently small positive number.
Let $\veca(I,J)\in \real^{\Lambda(\nbigi)}$
be determined by
$\veca(I,J)_j=1$ if
$j\in \{[t_m]\}\sqcup \bigsqcup_{i\in I}\Lambda(\nbigi,i)$,
and
$\veca(I,J)_j=1-\epsilon$ if
$j\in \bigsqcup_{i\in J}\Lambda(\nbigi,i)$.
Note that
$\nbigp^{\DM}_{\veca(I,J)}\nbigvtilde_{Y,m}'$
is independent of any sufficiently small $\epsilon$.

Let $V_{\Htilde'_{Y,m}}\nbigd_{\Ytilde'_{m}}
\subset\nbigd_{\Ytilde'_{m}}$
denote the subsheaf of algebras
generated by
$\Theta_{\Ytilde'_{m}}(\log \Htilde'_{Y,m})$
over $\nbigo_{\Ytilde'_{m}}$.
We obtain the following
$V_{\Htilde'_{Y,m}}\nbigd_{\Ytilde'_{m}}$-submodule
\begin{equation}
\nbign(I,J):=
V_{\Htilde_{Y,m}'}\nbigd_{\Ytilde'_{m}}\cdot
\nbigp^{\DM}_{\veca(I,J)}\nbigvtilde'_{Y,m}
\subset
\nbigvtilde'_{Y,m}.
\end{equation}
Let $V_{\Htilde_{Y,m}\cup \Ytilde_m^0}\nbigd_{\Ytilde_m}
\subset\nbigd_{\Ytilde_m}$
denote the sheaf of subalgebras
generated by
$\Theta_{\Ytilde_m}(\log (\Htilde_{Y,m}\cup \Ytilde_m^0))$
over $\nbigo_{\Ytilde_m}$.
Because
\[
 F^{\ast}\Theta_{\Ytilde_m}(\log (\Htilde_{Y,m}\cup \Ytilde_m^0))
\simeq
\Theta_{\Ytilde'_{m}}(\log \Htilde'_{Y,m}),
\]
we obtain
$V_{\Htilde_{Y,m}\cup \Ytilde_m^0}\nbigd_{\Ytilde_m}$-module
$F_{\ast}(\nbign(I,J))$.

\begin{lem}
\label{lem;21.3.9.10}
There exists a natural isomorphism
\begin{equation}
\label{eq;21.3.10.2}
\nbigm_{Y,m}(I,J)
 \simeq
 \nbigd_{\Ytilde_m}(\ast \Ytilde_m^0)
 \otimes_{V_{\Htilde_{Y,m}}\nbigd_{\Ytilde_m}(\ast\Ytilde_m^0)}
 \Bigl(
 F_{\ast}\bigl(\nbign(I,J)\bigr)(\ast \Ytilde_m^0)
 \Bigr).
\end{equation}
\end{lem}
\pf
We set
$\nbigm'_{Y,m}(I,J)
:=\nbigd_{\Ytilde_m}
\otimes_{V_{\Htilde_{Y,m}}\nbigd_{\Ytilde_m}}
 F_{\ast}(\nbign(I,J))$,
 which is a coherent $\nbigd_{\Ytilde_m}$-module.
The right hand side of (\ref{eq;21.3.10.2})
is naturally isomorphic to
$\nbigm'_{Y,m}(I,J)(\ast\Ytilde^0_m)$.
The natural inclusion
$F_{\ast}(\nbign(I,J))
\lrarr
\nbigvtilde_{Y,m}$
induces
\begin{equation}
 \label{eq;21.3.10.1}
 \nbigm'_{Y,m}(I,J)
\bigl(\ast(\Htilde_{Y,m}\cup\Ytilde_m^0)\bigr)
\lrarr\nbigvtilde_{Y,m}.
\end{equation}
Because the restriction of (\ref{eq;21.3.10.1})
to $\Ytilde_m\setminus\Ytilde_m^0$
is an isomorphism,
we obtain that (\ref{eq;21.3.10.1}) is an isomorphism.
Let $\DD$ denote the duality functor
of $\nbigd_{\Ytilde_m}$-modules,
then there exists a natural isomorphism
$\DD(\nbigm'_{Y,m}(I,J))\bigl(\ast(\Htilde_{Y,m}\cup\Ytilde^0_m)\bigr)
\simeq
 \DD(\nbigvtilde_{Y,m})\bigl(\ast(\Htilde_{Y,m}\cup\Ytilde^0_m)\bigr)$.
Hence, 
there exists the natural morphism
of $\nbigd_{\Ytilde_m}$-modules:
\[
 \nbigvtilde_{Y,m}\bigl(
 !(\Htilde_{Y,m}\cup \Ytilde_m^0)
 \bigr)
 \lrarr
 \nbigm'_{Y,m}(I,J).
\]
We set
$\Htilde_{Y,m}(I):=\bigcup_{i\in I}\Htilde_{Y,m,i}$.
We obtain the following morphisms:
\begin{multline}
\label{eq;21.3.9.2}
 \nbigm_{Y,m}(I,J)
 \llarr
  \nbigvtilde_{Y,m}\bigl(
 !(\Htilde_{Y,m}\cup \Ytilde_m^0)
 \bigr)\bigl(\ast (\Htilde_{Y,m}(I)\cup \Ytilde_m^0)\bigr)
 \lrarr \\
 \nbigm'_{Y,m}(I,J)
 \bigl(\ast (\Htilde_{Y,m}(I)\cup \Ytilde_m^0)\bigr)
 \llarr
 \nbigm'_{Y,m}(I,J)(\ast \Ytilde_m^0).
\end{multline}
The restriction of the morphisms in (\ref{eq;21.3.9.2})
to $\Ytilde_m\setminus\Ytilde_m^0$ are isomorphisms.
(See \cite[\S5.3.3]{Mochizuki-MTM},
where a similar issue is studied for
$\nbigr$-modules.)
Because
$\nbigm_{Y,m}(I,J)(\ast\Ytilde_m^0)
=\nbigm_{Y,m}(I,J)$,
the morphisms in (\ref{eq;21.3.9.2})
are isomorphisms.
\hfill\qed

\vspace{.1in}

Let
$V\nbigd_{\Ytilde_m}:=
\nbigo\langle
 t_{m}\del_{t_{m}},
 \del_{z_{m,1}},\ldots,\del_{z_{m,\ell}},
 \del_{z_{\ell+1}},\ldots,\del_{z_{n}}\rangle$.
 Let $\nbigl(I,J)\subset
 \nbigm_{Y,m}(I,J)$
 denote the $V\nbigd_{\Ytilde_m}$-submodule
 generated by the image of
 $F_{\ast}(\nbign(I,J))
 \lrarr \nbigm_{Y,m}(I,J)$.
By Lemma \ref{lem;21.3.9.10},
we have
\begin{equation}
\label{eq;21.3.9.11}
\nbigl(I,J)(\ast t_m)=\nbigm_{Y,m}(I,J).
\end{equation}

\subsubsection{Regularity}

Let $\nbigu_{Y,m}$ be a neighbourhood of
$(0,\ldots,0)$ in $\Ytilde_m$.
Let $f_{Y,m}$ be a coordinate function on $\nbigu_{Y,m}$
defining $\{t_m=0\}\cap\nbigu_{Y,m}$.

\begin{lem}
 \label{lem;21.3.8.12}
$\nbigl(I,J)_{|\nbigu_{Y,m}}$ is coherent over
$\nbigd_{\nbigu_{Y,m}/f_{Y,m}}$.
\end{lem}
\pf
By reordering $(z_1,\ldots,z_{\ell})$,
we may assume that 
for any $\vecn\in\nbigs(\nbigi)$
$n_i<0$ implies $n_j<0$ for any $j\leq i$.
We obtain the coordinate system
$(f_{Y,m},\ztilde_{m,1},\ldots,\ztilde_{m,\ell},\ztilde_{\ell+1},\ldots,\ztilde_{n})$
by setting
$\ztilde_{m,i}=z_{m,i}$ $(i=1,\ldots,\ell)$
and $\ztilde_i=z_i$ $(i=\ell+1,\ldots,n)$.
There exists a holomorphic function $B$ on $\nbigu_{Y,m}$
such that
$v:=\ztilde_{m,1}\del_{\ztilde_{m,1}}
=z_{m,1}\del_{z_{m,1}}
+Bz_{m,1}t_m\del_{t_m}$.

We may naturally regard $\nbigi$ as a subset of
$H^0(X,\nbigo_X(\ast H))
 \big/
 H^0(X,\nbigo_X)$.
There exists a finite subset
$\nbigi_{1}\subset H^0(X,\nbigo_X(\ast H))$
such that the projection
$H^0(X,\nbigo_X(\ast H))\lrarr
H^0(X,\nbigo_X(\ast H))/H^0(X,\nbigo_X)$
induces $\nbigi_{1}\simeq\nbigi$.
For each $\gminia\in\nbigi$,
let $\gminia_1\in\nbigi_1$ denote the corresponding element.
It is enough to study the case where
$v\bigl(\kappa_m^{\ast}(t\gminia_1)\bigr)\big/
 \kappa_m^{\ast}(t\gminia_1)$
are invertible on $\nbigu$
for any non-zero element $\gminia\in\nbigi$.
We obtain the sections
$F^{\ast}(v)$
of $\Theta_{\Ytilde_{m}'}(\log\Htilde'_{Y,m})$.
Let $\nbiga\subset V_{\Htilde'_{Y,m}}\nbigd_{\Ytilde'_{m}}$
denote the sheaf of subalgebras
generated by $F^{\ast}(v)$
over $\nbigo_{\Ytilde_m'}$.
We have the natural inclusion
$\nbiga\cdot\nbigp^{\DM}_{\veca(I,J)}(\nbigvtilde'_{Y,m})
\subset
 \nbign(I,J)$.
 
\begin{lem}
\label{lem;21.3.8.11}
$\nbign(I,J)=\nbiga\cdot \nbigp^{\DM}_{\veca(I,J)}(\nbigvtilde'_{Y,m})$
on $F^{-1}(\nbigu_{Y,m})$.
\end{lem}
\pf
Let $Q$ be any point of
$F^{-1}(\nbigu_{Y,m})\cap \Htilde_{Y,m}'$.
Because the formal completion is faithfully flat,
it is enough to prove that
$(\nbiga\cdot\nbigp^{\DM}_{\veca(I,J)}\nbigvtilde'_{Y,m})_{|\Qhat}
=\nbign(I,J)_{|\Qhat}$.
Note that
$\nbign(I,J)_{|\Qhat}
=V_{\Htilde'_{Y,m}}\nbigd_{\Ytilde'_{m}|\Qhat}\cdot
 \nbigp^{\DM}_{\veca(I,J)}(\nbigvtilde'_{Y,m})_{|\Qhat}$
and
$(\nbiga\cdot\nbigp^{\DM}_{\veca(I,J)}\nbigvtilde'_{Y,m})_{|\Qhat}
=\nbiga_{|\Qhat}\cdot
 \nbigp^{\DM}_{\veca(I,J)}\nbigvtilde'_{Y,m|\Qhat}$
which also follows from the faithful flatness
of the formal completion.

We obtain
$\nbigitilde'_{Y,m,Q}\subset
\nbigo_{\Ytilde'_{m}}(\ast\Htilde'_{Y,m})_Q
 \big/\nbigo_{\Ytilde'_{m},Q}$
as the inverse image of $\{t\gminia\,|\,\gminia\in\nbigi\}$
by
$\kappa_m\circ F$.
For each non-zero element
$\gminib\in\nbigitilde'_{Y,m,Q}$,
we choose $\gminia\in\nbigi$ such that
$\gminib=(\kappa_m\circ F)^{\ast}(t\gminia)$
in $\nbigo_{\Ytilde'_{m}}(\ast\Htilde'_{Y,m})_Q
 /\nbigo_{\Ytilde'_{m},Q}$,
and we set
$\gminib_1:=(\kappa_m\circ F)^{\ast}(t\gminia_1)$.
Note that
$F^{\ast}(v)(\gminib_1)/\gminib_1$ is invertible around $Q$.
Let $(\gminib)_{Q,\infty}$
denote the pole divisor of $\gminib_1$ around $Q$

There exists a decomposition
$(\nbigp^{\DM}_{\veca(I,J)}\nbigvtilde'_{Y,m},\nablatilde)_{|\Qhat}
=\bigoplus_{\gminib\in\nbigitilde'_{Y,m,Q}}
 (\nbigp^{\DM}_{\veca(I,J)}\nbigvtilde'_{\Qhat,\gminib},
 \nablatilde_{\Qhat,\gminib})$
 such that
 $\nablatilde_{\Qhat,\gminib}-d\gminib_1$ are logarithmic.
It is easy to check
\[
\nbiga_{|\Qhat}\cdot
\nbigp^{\DM}_{\veca(I,J)}\nbigvtilde'_{\Qhat,\gminib}
=\nbigp^{\DM}_{\veca(I,J)}\nbigvtilde'_{\Qhat,\gminib}(\ast (\gminib)_{Q,\infty})
=V_{\Htilde'_{Y,m}}\nbigd_{\Ytilde'_{m}|\Qhat} \cdot
 \nbigp^{\DM}_{\veca(I,J)}\nbigvtilde'_{\Qhat,\gminib}.
\]
Then, we obtain the claim of Lemma \ref{lem;21.3.8.11}.
\hfill\qed
 
\vspace{.1in}

Because $\nbign(I,J)$ is good as an $\nbiga$-module,
we can prove that
$F_{\ast}\nbign(I,J)$ is coherent over
$\nbigo_{\nbigu_{Y,m}}\langle v\rangle$
by the standard argument in the theory of $\nbigd$-modules
(see \cite[Theorem 4.25]{kashiwara_text}).
Note that
$\nbigl(I,J)_{|\nbigu_{Y,m}}$
is equal to the $\nbigd_{\nbigu_{Y,m}/f_{Y,m}}$-submodule
of $\nbigm_{Y,m}(I,J)_{|\nbigu_{Y,m}}$
generated by the image of
$F_{\ast}(\nbign(I,J))_{|\nbigu_{Y,m}}
\lrarr \nbigm_{Y,m}(I,J)_{|\nbigu_{Y,m}}$.
Then, we obtain
the claim of Lemma \ref{lem;21.3.8.12}.
\hfill\qed

\vspace{.1in}
We obtain the following lemma from Lemma \ref{lem;21.3.9.30}
and Lemma \ref{lem;21.3.8.12}.

\begin{lem}
$\nbigm_{Y,m}(I,J)_{|\nbigu_{Y,m}}$ is
regular along $f_{Y,m}$.
\hfill\qed
\end{lem}

\subsubsection{Proof of Proposition \ref{prop;21.1.23.11}}
\label{subsection;21.3.9.1}

Let $\nbigu_{X,m}$ be a neighbourhood of
$(0,\ldots,0)$ in $\Xtilde_m$.
Let $f_{X,m}$ be a coordinate function on $\nbigu_{X,m}$
defining $\{t_m=0\}\cap\nbigu_{X,m}$.
Note that
\[
 \bigl(
 \rho_{m\ast}(\nbigvtilde_{Y,m})
 (!\Htilde_{X,m}(J))
 \bigr)(\ast \Xtilde_m^0)
 \simeq
 \rho_{m\dagger}(\nbigm_{Y,m}(I,J)).
\]
Because $\nbigm_{Y,m}(I,J)$
is regular along $\rho_m^{\ast}(f_{X,m})$,
we obtain that
$\bigl(
 \rho_{m\ast}(\nbigvtilde_{Y,m})
 (!\Htilde_{X,m}(J))
 \bigr)(\ast \Xtilde_m^0)$
 is regular along $f_{X,m}$
by Lemma \ref{lem;21.3.9.22}.
There exists the natural inclusion
$\varphi_m^{\ast}(\nbigvtilde)
 \subset
 \rho_{m\ast}(\nbigvtilde_{Y,m})$,
which induces the following inclusion:
\[
 \varphi_m^{\ast}(\nbigvtilde)(!\Htilde_{X,m}(J))
 (\ast \Xtilde_m^0)
 \subset
 \rho_{m\ast}(\nbigvtilde_{Y,m})
 (!\Htilde_{X,m}(J))
 (\ast \Xtilde_m^0).
\]
There exists the following natural isomorphism:
\[
 \varphi_m^{\ast}(\nbigvtilde)(! \Htilde_{X,m}(J))
 (\ast \Xtilde_m^0)
 \simeq
 \varphi_m^{\ast}
 \Bigl(
 \bigl(
 \nbigvtilde(! \Htilde(J))
 \bigr)(\ast \Xtilde^0)
 \Bigr).
\]
Hence, 
$\varphi_m^{\ast}
 \Bigl(
 \bigl(
 \nbigvtilde(! \Htilde(J))
 \bigr)(\ast \Xtilde^0)
 \Bigr)$
is regular along $f_{X,m}$
by Lemma \ref{lem;21.3.9.20}.
Thus, Proposition \ref{prop;21.1.23.11}
is proved.
\hfill\qed

\section{$\nbigrtilde$-modules and $\gbigrtilde$-modules}

\subsection{Basic operations for $\nbigr$-modules
and $\nbigrtilde$-modules}
Let $X$ be a complex manifold.
Let $\Theta_X$ denote the tangent sheaf of $X$.
We set
$\nbigx=\cnum_{\lambda}\times X$.
Let $p_X:\nbigx\lrarr X$ denote the projection.
Let $\nbigr_X\subset\nbigd_{\nbigx}$ denote
the sheaf of subalgebras generated by
$\lambda p_X^{\ast}\Theta_X$ over $\nbigo_{\nbigx}$.
We set $\nbigrtilde_X:=\nbigr_X\langle\lambda^2\del_{\lambda}\rangle
\subset\nbigd_{\nbigx}$.

For any open subset $\nbigu\subset \cnum_{\lambda}\times X$,
the restriction of an $\nbigo_{\nbigx}$-module $\nbigm$ to
$\nbigu$ is denoted by $\nbigm_{|\nbigu}$.
If $\nbigu=\cnum_{\lambda}\times U$ for an open subset $U\subset X$,
the restriction is also denoted by $\nbigm_{|U}$.

Let $H\subset X$ be a complex hypersurface.
For an $\nbigo_{\nbigx}$-module $\nbigm$,
we set
$\nbigm(\ast H):=\nbigm\otimes_{\nbigo_{\nbigx}}
 \nbigo_{\nbigx}(\ast\nbigh)$.
We set
$\nbigr_{X(\ast H)}:=\nbigr_X(\ast H)$
and
$\nbigrtilde_{X(\ast H)}:=\nbigrtilde_X(\ast H)$.
If $H=f^{-1}(0)$ for a holomorphic function $f$ on $X$,
$\nbigm(\ast H)$ is also denoted by $\nbigm(\ast f)$.

We shall recall some basic operations
for $\nbigrtilde_{X(\ast H)}$-modules
from \cite{Sabbah-pure-twistor} and \cite{Mochizuki-MTM}.

\begin{rem}
We set $\nbigx^{\circ}:=(\cnum_{\lambda}\times\{0\})\times X$
and $\nbigh^{\circ}:=(\cnum_{\lambda}\times\{0\})\times H$.
Because $\nbigrtilde_{X(\ast H)|\nbigx^{\circ}}=
\nbigd_{\nbigx^{\circ}}(\ast\nbigh^{\circ})$,
any $\nbigrtilde_{X(\ast H)}$-module $\nbigm$
induces a $\nbigd_{\nbigx^{\circ}}(\ast\nbigh^{\circ})$-module
$\nbigm(\ast H)_{|\nbigx^{\circ}}$ by the restriction.
\hfill\qed
\end{rem}

\subsubsection{Direct images by proper morphisms}

Let $F:X\lrarr Y$ be a proper morphism of complex manifolds.
We set $\omega_{\nbigx}:=\lambda^{-\dim X}p_X^{\ast}\omega_{X}$,
where $\omega_X$ denotes the canonical line bundle of $X$.
Similarly, $\omega_{\nbigy}:=\lambda^{-\dim Y}p_Y^{\ast}\omega_Y$.
For a hypersurface $H_Y$ of $Y$,
we set $H_X:=F^{-1}(H_Y)$.
We set
$\nbigr_{Y\larr X}:=
 \omega_{\nbigx}\otimes_{F^{-1}(\nbigo_{\nbigy})}
 F^{-1}(\nbigr_{Y}\otimes\omega_{\nbigy}^{-1})$. 
For any $\nbigr_{X(\ast H_X)}$-module $\nbigm$,
we obtain the following $\nbigr_{Y(\ast H_Y)}$-modules:
\[
 F^i_{\dagger}(\nbigm):=
 R^i(\id_{\cnum}\times F)_{\ast}
 \bigl(
  \nbigr_{Y\larr X}(\ast H_X)\otimes^L_{\nbigr_{X(\ast H_X)}}
  \nbigm
 \bigr).
\]
If $\nbigm$ is a good $\nbigr_{X(\ast H_X)}$-module,
then $F_{\dagger}^i(\nbigm)$ are also good.
(See \cite[\S1.4]{Sabbah-pure-twistor}.)
If $\nbigm$ is an $\nbigrtilde_{X(\ast H_X)}$-module,
$F^i_{\dagger}(\nbigm)$ are naturally
$\nbigrtilde_{Y(\ast H_Y)}$-modules.
The restriction
$F^i_{\dagger}(\nbigm)_{|\nbigy^{\circ}}$
is naturally identified with the direct image
$F^i_{\dagger}(\nbigm_{|\nbigx^{\circ}})$
of the $\nbigd_{\nbigx^{\circ}}(\ast\nbigh^{\circ})$-module
$\nbigm_{|\nbigx^{\circ}}$.

\subsubsection{Strict specializability along a coordinate function}
\label{subsection;21.4.12.40}

Let us consider the case where $X$ is an open subset of
$\proj^1\times X_0$ for a complex manifold $X_0$.
Let $t$ denote the standard coordinate of $\cnum\subset\proj^1$.
Let $H$ be a hypersurface such that
$\{\infty\}\times X\subset H$
and $\dim (H\cap \{t=0\})<\dim H$.
Let $V\nbigr_{X(\ast H)}\subset \nbigr_{X(\ast H)}$ denote
the sheaf of subalgebras generated by
$\lambda p_X^{\ast}\bigl(\Theta_X(\log t)\bigr)$
over $\nbigo_{\nbigx}(\ast H)$.
We set
$V\nbigrtilde_{X(\ast H)}=
V\nbigr_{X(\ast H)}\langle\lambda^2\del_{\lambda}\rangle
\subset\nbigrtilde_{X(\ast H)}$.
We set $\deldel_t:=\lambda\del_t$ on $\nbigx$.

Let $\nbigm$ be an $\nbigrtilde_{X(\ast H)}$-module
which is coherent over $\nbigr_{X(\ast H)}$.
We say that $\nbigm$ is strictly specializable along $t$
if there exists an increasing filtration
$V_{\bullet}(\nbigm)$ by
$V\nbigrtilde_{X(\ast H)}$-submodules indexed by $\real$
satisfying the following conditions.
\begin{itemize}
 \item $V_a(\nbigm)(\ast t)=\nbigm(\ast t)$.
 \item For any $a\in\real$ and any compact subset $K\subset X$,
       there exists a positive number $\epsilon$
       such that
       $V_a(\nbigm)_{|\cnum_{\lambda}\times K}=
       V_{a+\epsilon}(\nbigm)_{|\cnum_{\lambda}\times K}$.
 \item $\Gr^V_a(\nbigm):=V_{a}(\nbigm)/V_{<a}(\nbigm)$
       is strict,
       i.e., flat over $\nbigo_{\cnum_{\lambda}}$.
 \item We have $tV_a(\nbigm)\subset V_{a-1}(\nbigm)$ for any $a$.
       If $a<0$,
       then $tV_a(\nbigm)=V_{a-1}(\nbigm)$.
 \item We have $\deldel_tV_a(\nbigm)\subset V_{a+1}(\nbigm)$ for any $a$.
       If $a>-1$,
       the induced morphism
       $\Gr^V_{a}(\nbigm)\lrarr
       \Gr^V_{a+1}(\nbigm)$
       is an isomorphism of sheaves.
 \item For each $a\in\real$, $-\deldel_tt-a\lambda$ is locally nilpotent
       on $\Gr^V_a(\nbigm)$.
 \item Each $V_a(\nbigm)$ is coherent over $V\nbigr_{X(\ast H)}$.
\end{itemize}
For such $\nbigm$,
we set
$\nbigm[\ast t]:=
\nbigr_{X(\ast H)}\otimes_{V\nbigr_{X(\ast H)}}V_0\nbigm$
and
$\nbigm[! t]:=
\nbigr_{X(\ast H)}\otimes_{V\nbigr_{X(\ast H)}}V_{<0}\nbigm$.
They are also strictly specializable.

\begin{rem}
If $H=\emptyset$,
the restriction $V_{\bullet}(\nbigm)_{|\nbigx^{\circ}}$
is the $V$-filtration of
the $\nbigd_{\nbigx^{\circ}}$-module $\nbigm_{|\nbigx^{\circ}}$.
 Note that it satisfies a stronger condition
 than the ordinary $V$-filtrations,
i.e., $V_{a}(\nbigm)_{|\nbigx^{\circ}}$
are coherent over $(V\nbigr_X)_{|\nbigx^{\circ}}$,
 not only over $(V\nbigrtilde_X)_{|\nbigx^{\circ}}$.
\hfill\qed
\end{rem}

As a variant,
for an $\nbigrtilde_{X(\ast H)}(\ast t)$-module $\nbigm$
which is coherent over $\nbigr_{X(\ast H)}(\ast t)$,
we say that $\nbigm$ is strictly specializable
if it satisfies the above conditions
replacing the fourth and fifth conditions
with the following conditions.
\begin{itemize}
 \item We have $tV_a(\nbigm)=V_{a-1}(\nbigm)$ for any $a$.
 \item We have $\deldel_tV_a(\nbigm)\subset V_{a+1}(\nbigm)$ for any $a$.
\end{itemize}
Even for such $\nbigm$,
we obtain the coherent $\nbigr_{X(\ast H)}$-modules
$\nbigm[\ast t]:=
\nbigr_{X(\ast H)}\otimes_{V\nbigr_{X(\ast H)}}V_0\nbigm$
and
$\nbigm[! t]:=
\nbigr_{X(\ast H)}\otimes_{V\nbigr_{X(\ast H)}}V_{<0}\nbigm$.

\subsubsection{Strict specializability and localizability
   along a function}

Let $X$ be any complex manifold with a hypersurface $H$.
Let $f$ be a meromorphic function on $X$.
The zero divisor and the pole divisor are denoted by
$(f)_0$ and $(f)_{\infty}$,
respectively.
The supports of the divisors are denoted by
$|(f)_0|$ and $|(f)_{\infty}|$,
respectively.
We say that $f$ is meromorphic on $(X,H)$
if $|(f)_{\infty}|\subset H$.

Assume $|(f)_0|\cap|(f)_{\infty}|=\emptyset$.
Let $\iota_f:X\lrarr \Xtilde:=X\times\proj^1$ denote the graph of $f$.
We set $\Htilde:=(H\times\proj^1)\cup (X\times\{\infty\})$.
An $\nbigrtilde_{X(\ast H)}$-module $\nbigm$
which is coherent over $\nbigr_{X(\ast H)}$
induces an $\nbigrtilde_{\Xtilde(\ast\Htilde)}$-module
$\iota_{f\dagger}(\nbigm)$
which is coherent over $\nbigr_{\Xtilde(\ast \Htilde)}$.
It is called strictly specializable along $f$ modulo $H$
if $\iota_{f\dagger}(\nbigm)$ is strictly specializable along $t$
as an $\nbigrtilde_{X(\ast H)}$-module.
We say that $\nbigm$ is localizable along $f$ modulo $H$
if moreover there exist $\nbigrtilde_{X(\ast H)}$-modules
$\nbigm[!f]$ and $\nbigm[\ast f]$ which are coherent over $\nbigr_{X(\ast H)}$
such that
$\iota_{f\dagger}(\nbigm[\star f])
\simeq
 (\iota_{f\dagger}\nbigm)[\star t]$ for $\star=!,\ast$.
Such $\nbigm[!f]$ and $\nbigm[\ast f]$
are uniquely determined.

Similarly,
an $\nbigrtilde_{X(\ast H)}(\ast f)$-module $\nbigm$
which is coherent over $\nbigr_{X(\ast H)}(\ast f)$
induces an $\nbigrtilde_{\Xtilde(\ast\Htilde)}(\ast t)$-module
$\iota_{f\dagger}(\nbigm)$
which is coherent over $\nbigr_{\Xtilde(\ast \Htilde)}(\ast t)$.
It is called strictly specializable along $f$ modulo $H$
if $\iota_{f\dagger}(\nbigm)$ is strictly specializable along $t$
as an $\nbigrtilde_{X(\ast H)}(\ast t)$-module.
We say that $\nbigm$ is localizable along $f$ modulo $H$
if moreover there exist $\nbigrtilde_{X(\ast H)}$-modules
$\nbigm[!f]$ and $\nbigm[\ast f]$ which are coherent over $\nbigr_{X(\ast H)}$
such that
$\iota_{f\dagger}(\nbigm[\star f])
\simeq
 (\iota_{f\dagger}\nbigm)[\star t]$ for $\star=!,\ast$.
Such $\nbigm[!f]$ and $\nbigm[\ast f]$
are uniquely determined.

\begin{rem}
In the above definitions,
$\nbigm[\star f]$ depend on the choice of $H$. 
For example, we set $H_1=H\cup|(f)_0|$
and $\Htilde_1=(H_1\times\proj^1)\cup(X\times\{\infty\})$.
Note that $\iota_{f\dagger}(\nbigm(\ast f))$
is coherent over
$V\nbigr_{\Xtilde(\ast\Htilde_1)}$.
The $V$-filtration of
$\iota_{f\dagger}(\nbigm(\ast f))$
is trivial as an $\nbigr_{\Xtilde(\ast\Htilde_1)}$-module.
Hence,
$\nbigm(\ast f)$ is always localizable along $f$
modulo $H_1$,
and we have $\nbigm[\star f]=\nbigm(\ast f)$ $(\star=\ast,!)$.
When we emphasize the dependence on $H$,
we use the notation
$\nbigm[\star f](\ast H)$.
\hfill\qed
\end{rem}

Let us consider the case where 
$|(f)_0|\cap|(f)_{\infty}|\neq\emptyset$.
Let $\nbigm$ be an $\nbigrtilde_{X(\ast H)}$-module
which is coherent over $\nbigr_{X(\ast H)}$.
We say that
$\nbigm$ is strictly specializable (resp. localizable) along $f$ modulo $H$
if there exists a proper morphism
$\rho:X'\lrarr X$ such that
(i) $\rho$ induces $X'\setminus \rho^{-1}(H)\simeq X\setminus H$,
(ii) $|(\rho^{\ast}f)_0|\cap|(\rho^{\ast}f)_{\infty}|=\emptyset$,
(iii) $\rho^{\ast}(\nbigm)$ is strictly specializable
(resp. localizable) along $\rho^{\ast}(f)$.
If $\nbigm$ is localizable along $f$ modulo $H$
in the above sense,
we obtain
the $\nbigrtilde_{X'(\ast H')}$-module
$\rho^{\ast}(\nbigm)[\star \rho^{\ast}(f)]$
which is coherent over $\nbigr_{X'(\ast H')}$,
where $H'=\rho^{-1}(H)$.
We obtain the $\nbigrtilde_{X(\ast H)}$-module
$\nbigm[\star f]:=
\rho_{\ast}\Bigl(
\rho^{\ast}(\nbigm)[\star \rho^{\ast}(f)]
\Bigr)$
which is coherent over $\nbigr_{X(\ast H)}$,
and independent of the choice of $\rho$.
Similar concepts are defined for
$\nbigrtilde_{X(\ast H)}(\ast f)$-modules.

\subsubsection{Beilinson functors along a function}

For $a<b$,
we set
\[
 \II^{a,b}_f:=\bigoplus_{a\leq j<b}
 \nbigo_{\nbigx}\bigl(\ast ((f)_0\cup\nbigh)\bigr)(\lambda s)^j.
\]
Here, $s$ denotes a formal variable.
It is naturally an $\nbigrtilde_X$-module
by the meromorphic flat connection $\nabla$
determined by
$\nabla(s^j)=s^{j+1}df/f$.

Let $\nbigm$ be an $\nbigrtilde_{X(\ast H)}(\ast f)$-module
which is a coherent over $\nbigr_{X(\ast H)}(\ast f)$.
We obtain the $\nbigrtilde_{X(\ast H)}(\ast f)$-module
$\Pi^{a,b}_f(\nbigm):=\II^{a,b}_f\otimes_{\nbigo_{\nbigx}(\ast H)}\nbigm$
which is coherent over $\nbigr_{X(\ast H)}(\ast f)$.
If $\nbigm$ is strictly specializable along $f$ modulo $H$,
each $\Pi^{a,b}_f(\nbigm)$ is also strictly specializable along $f$ modulo $H$.
Assume that $\Pi^{a,b}_f(\nbigm)$
are localizable along $f$ modulo $H$
for any $a,b$.
We set
$\Pi^{a,b}_{f,\star}(\nbigm):=
\Pi^{a,b}_f(\nbigm)[\star f]$ for $\star=!,\ast$.
Following \cite{beilinson2},
we define the $\nbigrtilde_{X(\ast H)}$-module
$\Pi^{a,b}_{f,\ast!}(\nbigm)$
as
\[
\Pi^{a,b}_{f,\ast!}(\nbigm):=
 \varprojlim_{N\to\infty}
 \Cok\Bigl(\Pi^{b,N}_{f!}(\nbigm)
 \lrarr
 \Pi^{a,N}_{f\ast}(\nbigm)\Bigr).
\]
In particular,
we set
$\psi^{(a)}_f(\nbigm):=
 \Pi^{a,a}_{f,\ast !}(\nbigm)$
and
$\Xi^{(a)}_f(\nbigm):=
 \Pi^{a,a+1}_{f,\ast !}(\nbigm)$.
 (See \cite[\S4.1]{Mochizuki-MTM} for more details.)

\vspace{.1in}
Let $\nbigm$ be an $\nbigrtilde_{X(\ast H)}$-module
which is coherent over $\nbigr_{X(\ast H)}$.
Suppose that $\nbigm$ is strictly specializable along $f$
and that $\Pi^{a,b}_f(\nbigm(\ast f))$
are localizable along $f$ for any $a,b$.
We set
$\Pi^{a,b}_f(\nbigm):=
\Pi^{a,b}_f(\nbigm(\ast f))$,
$\Pi^{a,b}_{f\star}(\nbigm):=
\Pi^{a,b}_{f\star}(\nbigm(\ast f))$
and
$\Pi^{a,b}_{f,\ast !}(\nbigm):=
\Pi^{a,b}_{f,\ast !}(\nbigm(\ast f))$.
In particular,
we set
$\psi^{(a)}_f(\nbigm):=
\psi^{(a)}_f(\nbigm(\ast f))$
and
$\Xi^{(a)}_f(\nbigm):=\Xi^{(a)}_f(\nbigm(\ast f))$.
There exists the following naturally defined complex of
$\nbigrtilde_{X(\ast H)}$-modules:
\begin{equation}
\label{eq;21.2.16.50}
 \nbigm[!f]
 \lrarr
 \Xi^{(0)}_f(\nbigm)
 \oplus
 \nbigm
 \lrarr
 \nbigm[\ast f].
\end{equation}
We define
the $\nbigrtilde_{X(\ast H)}$-module $\phi^{(0)}_f(\nbigm)$
as the cohomology of the complex (\ref{eq;21.2.16.50}).
There exist the natural morphisms
$\psi^{(1)}_f(\nbigm)\lrarr\phi^{(0)}_f(\nbigm)
\lrarr\psi^{(0)}_f(\nbigm)$.
Then,
as explained in \cite{beilinson2},
$\nbigm$ is naturally isomorphic to
the cohomology of the following complex
as an $\nbigrtilde_{X(\ast H)}$-module:
\begin{equation}
 \psi^{(1)}_f(\nbigm)
  \lrarr
  \phi^{(0)}_f(\nbigm)
  \oplus
  \Xi^{(0)}_f(\nbigm)
  \lrarr
  \psi^{(0)}_f(\nbigm).
\end{equation}
When we emphasize the dependence on $H$,
we use the notation
$\psi^{(a)}_f(\nbigm)(\ast H)$,
$\Xi_f^{(a)}(\nbigm)(\ast H)$,
$\phi^{(0)}_f(\nbigm)(\ast H)$,
etc.

\subsubsection{Localizability along hypersurfaces}

Let $H^{(1)}$ be a hypersurface of $X$.
Let $\nbigm$ be an $\nbigrtilde_{X(\ast H)}$-module
which is coherent over $\nbigr_{X(\ast H)}$.
We say that it is localizable along $H^{(1)}$ modulo $H$
if the following holds.
\begin{itemize}
 \item Let $Y$ be any open subset of $X$.
       Let $f$ be a meromorphic function on $(Y,H\cap Y)$
       such that
       $|(f)_0|\cup (H\cap Y)=(Y\cap H^{(1)})\cup(H\cap Y)$.
       Then, $\nbigm_{|Y}$ is strictly specializable
      and localizable along $f$ modulo $H$.
\end{itemize}
If $\nbigm$ is localizable along $H^{(1)}$ modulo $H$,
there exist $\nbigrtilde_{X(\ast H)}$-modules
$\nbigm[\star H^{(1)}]$ $(\star=!,\ast)$
which are coherent over $\nbigr_{X(\ast H)}$,
such that
$\nbigm[\star H^{(1)}]_{|Y}=\nbigm_{|Y}[\star f]$
for any $(Y,f)$ as above.
When we emphasize the dependence on $H$,
we use the notation
$\nbigm[\star H^{(1)}](\ast H)$.
If $H^{(1)}\subset H$,
the condition is trivial,
and we have $\nbigm[\star H^{(1)}](\ast H)=\nbigm$.

\subsection{$\nbigrtilde$-modules underlying integrable mixed twistor $\nbigd$-modules}

Let $\MTM^{\integral}(X)$ denote the category of
integrable mixed twistor $\nbigd_X$-modules.
Recall that an integrable mixed twistor $\nbigd_X$-module
is defined to be
an integrable $\nbigr_X$-triple $(\nbigm',\nbigm'',C)$
with a weight filtration $W$ satisfying some conditions.
(See \cite[\S2.1.5]{Mochizuki-MTM} for integrable
$\nbigr_X$-triples,
which originally goes back to \cite{Sabbah-pure-twistor}.
See \cite[\S7.2.3]{Mochizuki-MTM} for
integrable mixed twistor $\nbigd$-modules.)
In the following,
an integrable $\nbigr_X$-triple is called
$\nbigrtilde_X$-triple.

Let $\nbigc(X)$ denote the full subcategory of $\nbigrtilde_X$-modules
underlying integrable mixed twistor $\nbigd_X$-modules,
i.e.,
an $\nbigrtilde_X$-module $\nbigm''$ is an object of $\nbigc(X)$
if and only if
there exists
$((\nbigm',\nbigm'',C),W)\in\MTM^{\integral}(X)$.

Let $H$ be a hypersurface of $X$.
Let $\MTM^{\integral}(X;H)$ denote the essential image of
the natural functor from
$\MTM^{\integral}(X)$
to the category of filtered $\nbigrtilde_{X(\ast H)}$-triples.
An object in $\MTM^{\integral}(X;H)$ is called
an integrable mixed twistor $\nbigd_{X(\ast H)}$-module.
Similarly, let $\nbigc(X;H)$ denote the essential image of
the natural functor
from $\nbigc(X)$ to the category of $\nbigrtilde_{X}(\ast H)$-modules.
For any open subset $Y\subset X$,
the restriction induces the functors
$\MTM^{\integral}(X;H)\lrarr \MTM^{\integral}(Y;H_Y)$
and 
$\nbigc(X;H)\lrarr\nbigc(Y;H_Y)$,
where $H_Y:=H\cap Y$.

Because $\MTM^{\integral}(X;H)$ is an abelian category,
the following lemma is obvious.
\begin{lem}
Let $g:\nbigm_1\lrarr\nbigm_2$ be a morphism in $\nbigc(X;H)$.
Suppose that $g$ is induced by a morphism
in $\MTM^{\integral}(X;H)$, 
i.e., there exist
$(\nbigt_i,W)=((\nbigm'_i,\nbigm_i,C_i),W)\in\MTM^{\integral}(X;H)$
and
a morphism
$(g',g):(\nbigt_1,W)\lrarr(\nbigt_2,W)$ in $\MTM^{\integral}(X;H)$.
Then, 
$\Ker(g)$, $\Image(g)$ and $\Cok(g)$
are also objects in $\nbigc(X;H)$. 
\hfill\qed
\end{lem}

\begin{rem}
\label{rem;21.4.13.1}
A morphism $g:\nbigm_1\lrarr\nbigm_2$ in $\nbigc(X;H)$
is not necessarily induced by
a morphism in the category $\MTM^{\integral}(X;H)$.
If not,
$\Ker(g)$, $\Image(g)$ and $\Cok(g)$
are not necessarily objects in $\nbigc(X;H)$.
For example,
for any $\nbigm\in \nbigc(X;H)$,
the natural inclusion $\iota:\nbigm\lrarr \lambda^{-1}\nbigm$
is a morphism in $\nbigc(X;H)$,
which is not induced by a morphism
in $\MTM^{\integral}(X;H)$,
and $\Cok(\iota)$ is not an object in $\nbigc(X;H)$. 
In particular, the category $\nbigc(X;H)$ is not abelian.
\hfill\qed
\end{rem}

We recall that any $\nbigm\in\nbigc(X)$
is localizable along $H$.
Let $\nbigc(X,[\ast H])\subset\nbigc(X)$ denote the full subcategory of
objects $\nbigm\in\nbigc(X)$ such that $\nbigm[\ast H]\simeq\nbigm$.
Similarly,
$\nbigc(X,[!H])\subset\nbigc(X)$ denote the full subcategory of
objects $\nbigm\in\nbigc(X)$ such that $\nbigm[!H]\simeq\nbigm$.
\begin{lem}
The naturally induced functors $\nbigc(X,[\star H])\lrarr \nbigc(X;H)$ $(\star=!,\ast)$
are equivalent. 
\end{lem}
\pf
For any $\nbigm\in\nbigc(X;H)$,
there exists $\nbigmtilde\in\nbigc(X)$ such that
$\nbigmtilde(\ast H)=\nbigm$.
According to \cite[Proposition 11.2.1]{Mochizuki-MTM},
we have $\nbigmtilde[\ast H]\in\nbigc(X)$
which satisfies $(\nbigmtilde[\ast H])(\ast H)=\nbigmtilde(\ast H)=\nbigm$.
Hence, we obtain the essential surjectivity.
See \cite[Lemma 3.3.19]{Mochizuki-MTM}
for the fully faithfulness.
\hfill\qed

\subsubsection{Strictly specializability and localizability}

\begin{prop}
Let $\nbigm\in\nbigc(X;H)$.
\begin{itemize}
 \item $\nbigm$ is strictly specializable and localizable along
       any meromorphic function $f$ on $(X,H)$ modulo $H$,
       and the $\nbigrtilde_{X(\ast H)}$-module
       $\nbigm[\star f]$ $(\star=\ast,!)$
       are objects of $\nbigc(X;H)$.
       Moreover,
       $\Pi^{a,b}_f(\nbigm)$ are localizable
       along $f$ modulo $H$ for any $a<b$,
       and the induced $\nbigrtilde_{X(\ast H)}$-modules
       $\Pi^{a,b}_{f,\star}(\nbigm)$ $(\star=!,\ast)$
       $\Pi^{a,b}_{f,\ast!}(\nbigm)$,
       $\Xi^{(a)}_f(\nbigm)$,
       $\psi^{(a)}_f(\nbigm)$
       and $\phi_f^{(a)}(\nbigm)$
       are also objects in $\nbigc(X;H)$.
 \item $\nbigm$ is localizable along any hypersurface $H^{(1)}$ of $X$ modulo $H$,
       and $\nbigm[\star H^{(1)}]$  $(\star=\ast,!)$
       are objects of $\nbigc(X;H)$.
\end{itemize}
\end{prop}
\pf
For any hypersurface $H^{(1)}$ and $\nbigm\in\nbigc(X;H)$,
we take $\nbigmtilde\in\nbigc(X,[\ast H])$
such that $\nbigmtilde(\ast H)=\nbigm$,
then we have
$\nbigm[\star H^{(1)}]
=(\nbigmtilde[\star H^{(1)}])(\ast H)
\in\nbigc(X;H)$.
Thus, we obtain the second claim.
We obtain the first claim similarly.
\hfill\qed

\subsubsection{Direct image}
\label{subsection;21.6.16.1}

Let $F:X\lrarr Y$ be a projective morphism of complex manifolds.
Let $H_Y$ be a hypersurface of $Y$,
and we set $H_X:=F^{-1}(H_Y)$. 
We obtain the following proposition
from the functorial properties of
integrable mixed twistor $\nbigd$-modules \cite{Mochizuki-MTM}.
\begin{prop}
\label{prop;21.6.22.16}
Let $\nbigm\in\nbigc(X;H_X)$.
\begin{itemize}
 \item $F_{\dagger}^j(\nbigm)\in\nbigc(Y;H_Y)$.
 \item For any meromorphic function $f$ on $(Y,H_Y)$,
       there exist natural isomorphisms
\[
       \Pi^{a,b}_{f,\star}(F^j_{\dagger}(\nbigm))
       \simeq
       F^j_{\dagger}(\Pi^{a,b}_{F^{\ast}(f),\star}(\nbigm))
       \quad
       (\star=\ast,!).
\]
       They induce the following isomorphisms:
\[
       \Pi^{a,b}_{f,\ast!}(F^j_{\dagger}(\nbigm))
       \simeq
       F^j_{\dagger}(\Pi^{a,b}_{F^{\ast}(f),\ast!}(\nbigm)),\quad
       \Xi^{(a)}_{f}(F^j_{\dagger}(\nbigm))
       \simeq
       F^j_{\dagger}(\Xi^{(a)}_{F^{\ast}(f)}(\nbigm)),\quad
\]       
\[
       \psi^{(a)}_{f}(F^j_{\dagger}(\nbigm))
       \simeq
       F^j_{\dagger}(\psi^{(a)}_{F^{\ast}(f)}(\nbigm)),\quad
       \phi^{(a)}_{f}(F^j_{\dagger}(\nbigm))
       \simeq
       F^j_{\dagger}(\phi^{(a)}_{F^{\ast}(f)}(\nbigm)).
\]
 \item For any hypersurface $H^{(1)}_Y$ of $Y$,
       there exist natural isomorphisms
       $F_{\dagger}^j(\nbigm)[\star H_Y^{(1)}]
       \simeq
       F_{\dagger}^j\bigl(\nbigm[\star F^{-1}(H_Y^{(1)})]\bigr)$.
\item If $F$ induces an isomorphism
$X\setminus H_X\simeq Y\setminus H_Y$,
then we have
$F_{\dagger}^j(\nbigm)=0$ $(j\neq 0)$
for any $\nbigm\in\nbigc(X;H_X)$.
Moreover,
$F^0_{\dagger}$ induces an equivalence between
$\nbigc(X;H_X)\simeq \nbigc(Y;H_Y)$. 
Note that
$F_{\dagger}^0(\nbigm)=F_{\ast}(\nbigm)$ in this case.
A quasi-inverse $\nbigc(Y;H_Y)\simeq \nbigc(X;H_X)$
is given by the correspondence
$\nbign\longmapsto F^{\ast}(\nbign)$.
\hfill\qed
\end{itemize}
\end{prop}

\subsubsection{Kashiwara equivalence}
\label{subsection;22.7.31.1}

We continue to use the notation in \S\ref{subsection;21.6.16.1}.
Let $\nbigc_{F(X)}(Y;H_Y)\subset\nbigc(X;H_X)$
denote the full subcategory of
$\nbigm\in\nbigc(Y;H_Y)$ such that
the support of $\nbigm$ is contained in
the closed complex analytic subset $F(X)\subset Y$.
\begin{prop}
\label{prop;21.4.9.21}
Suppose that $F_{|X\setminus H_X}$ is a closed embedding
of $X\setminus H_X$ into $Y\setminus H_Y$.
 We obtain $F^j_{\dagger}(\nbigm)=0$ $(j\neq 0)$
for any $\nbigm\in\nbigc(X;H_X)$.
Moreover, $F^0_{\dagger}$ induces an equivalence
between $\nbigc(X;H_X)$ and $\nbigc_{F(X)}(Y;H_Y)$.
\end{prop}
\pf
Let us consider the case where $F$ is a closed immersion.
We have $F_{\dagger}^j(\nbigm)=0$ $(j\neq 0)$ for any $\nbigm\in\nbigc(X)$.
It implies that $F_{\dagger}^j(\nbigm)=0$ $(j\neq 0)$
for any $\nbigm\in\nbigc(X;H_X)$.
Let $\nbigm\in\nbigc_{F(X)}(Y;H_Y)$.
We obtain $\nbigm[\ast H_Y]\in\nbigc_{F(X)}(Y)$.
By \cite[Proposition 7.2.8]{Mochizuki-MTM},
there exists $\nbign\in\nbigc(X)$
such that $F_{\dagger}^0(\nbign)=\nbigm[\ast H_Y]$.
Because
$\nbigm=F_{\dagger}^0(\nbign)(\ast H_Y)
=F_{\dagger}^0(\nbign(\ast H_Y))$,
we obtain the essential surjectivity
of $F_{\dagger}^0:\nbigc(X;H_X)\lrarr\nbigc(Y;H_Y)$.
It is enough to check the full faithfulness
locally around any point $P$ of $F(X)$.
We may assume that
$Y$ is a neighbourhood of $(0,\ldots,0)$ in $\cnum^n$,
and $X=\bigcap_{i=1}^{\ell}\{z_i=0\}$.
By an induction,
it is enough to study the case $\ell=1$.
Let $g:F_{\dagger}^0(\nbign_1)\lrarr F_{\dagger}^0(\nbign_2)$
be a morphism
in $\nbigc_{F(X)}(Y;H_Y)$ for $\nbign_i\in\nbigc(X;H_X)$.
We have $\nbign_i=\psi_{z_1,-1}F_{\dagger}^0(\nbign_i)$,
and we obtain 
the induced morphism $g_0=\psi_{z_1,-1}(g):\nbign_1\lrarr\nbign_2$
in $\nbigc(X;H_X)$.
It is easy to see that $g=F_{\dagger}^0(g_0)$.

\vspace{.1in}

Let us consider the general case.
There exists a projective morphism of complex manifolds
$\rho_Y:Y'\lrarr Y$
such that
(i) $\rho_Y$ induces an isomorphism
$Y'\setminus \rho_Y^{-1}(H_Y)\simeq Y\setminus H_Y$,
(ii) the proper transform $Z$ of $F(X)$ is a complex submanifold of $Y'$.
(See \cite{Wlodarczyk}.)
We set $H_{Y'}:=\rho_Y^{-1}(H_Y)$.
We set $H_{Z}=H_{Y'}\cap Z$.
The following lemma is obvious because
$\rho_{Y\dagger}^0$ induces an equivalence
$\nbigc(Y';H_{Y'})
\simeq
\nbigc(Y;H_{Y})$.
\begin{lem}
$\rho_{Y\dagger}^0$ induces an equivalence
 $\nbigc_{Z}(Y';H_{Y'})
 \simeq
 \nbigc_{F(X)}(Y;H_Y)$.
\hfill\qed
\end{lem}

We obtain the complex analytic space
$X_1:=X\times_{F(X)}Z$ as the fiber product.
We set $H_{X_1}:=H_X\times_{F(X)}Z\subset X_1$.
The natural morphism $X_1\lrarr X$ induces
$X_1\setminus H_{X_1}\simeq X\setminus H_X$.
There exists a projective birational morphism
$\nu:X'\lrarr X_1$
such that
(i) $X'$ is a complex manifold,
(ii) $X'\setminus \nu^{-1}(H_{X_1})\simeq X_1\setminus H_{X_1}$.
We set $H_{X'}:=\nu^{-1}(H_{X_1})$.
For the induced morphisms
$\nu_X:X'\lrarr X$
and $\nu_Z:X'\lrarr Z$,
we have
$\nu_{X\dagger}^j(\nbigm)=0$ $(j\neq 0)$
and
$\nu_{Z\dagger}^j(\nbigm)=0$ $(j\neq 0)$
for any $\nbigm\in\nbigc(X';H_{X'})$,
and $\nu_{X\dagger}^0$
and $\nu_{Z\dagger}^0$
induce equivalences
$\nbigc(X';H_{X'})
\simeq
\nbigc(X;H_X)$
and 
$\nbigc(X';H_{X'})
\simeq
\nbigc(Z;H_Z)$.
Because
$F_{\dagger}^0\circ \nu_{X\dagger}^0
=\rho_{\dagger}^0\circ \nu_{Z\dagger}^0$,
we obtain that $F_{\dagger}^0$
induces an equivalence
$\nbigc(X;H_X)\simeq\nbigc_{F(X)}(Y;H_Y)$.
\hfill\qed

\vspace{.1in}
Let us study a minor refinement.
Let $H_Y=H_Y^{(1)}\cup H_Y^{(2)}$
be a decomposition of hypersurfaces
such that $\dim H_Y^{(1)}\cap H_Y^{(2)}<\dim H_Y$.
We set $H_X=f^{-1}(H_Y)$
and $H_X^{(2)}=f^{-1}(H_Y^{(2)})$.
We have the decomposition
$H_X=H_X^{(1)}\cup H_X^{(2)}$
such that $\dim H_X^{(1)}\cap H_X^{(2)}<\dim H_X$.
Let $\nbigc(Y,[\star H^{(1)}_Y];H_Y^{(2)})$
denote the essential image of
the natural functor
$\nbigc(Y,[\star H^{(1)}_Y])\lrarr \nbigc(Y;H_Y^{(2)})$.
We use the notation $\nbigc(X,[\star H^{(1)}_X];H_X^{(2)})$
in a similar meaning.
Because
$\nbigc(X,[\star H^{(1)}_X];H_X^{(2)})
\lrarr\nbigc(X;H_X)$
and
$\nbigc(Y,[\star H^{(1)}_Y];H_Y^{(2)})
\lrarr\nbigc(Y;H_Y)$
are equivalence, we obtain the following.

\begin{cor}
\label{cor;22.7.31.2}
In Proposition {\rm\ref{prop;21.4.9.21}},
$F_{\dagger}^0$ induces an equivalence
 $\nbigc(X,[\star H^{(1)}_X];H_X^{(2)})
 \lrarr
\nbigc(Y,[\star H^{(1)}_Y];H_Y^{(2)})$.
\hfill\qed
\end{cor}

\subsubsection{External tensor product}
\label{subsection;21.4.14.3}

Let $X_i$ $(i=1,2)$ be complex manifolds
with a hypersurface $H_i$.
Let $\nbigm_i$ be $\nbigr_{X_i(\ast H_i)}$-modules.
Let $p_i:\cnum_{\lambda}\times(X_1\times X_2)\lrarr\nbigx_i$
denote the projections.
We set $\Htilde=(H_1\times X_2)\cup (X_1\times H_2)$.
We obtain
the $\nbigr_{X_1\times X_2(\ast\Htilde)}$-module
\[
 \nbigm_1\boxtimes\nbigm_2
=p_1^{\ast}(\nbigm_1)\otimes_{\nbigo_{\cnum\times X_1\times X_2}}
 p_2^{\ast}(\nbigm_2).
\]
For $\nbigm_i\in\nbigc(X_i;H_i)$,
we obtain $\nbigm_1\boxtimes\nbigm_2\in\nbigc(X_1\times X_2)(\ast\Htilde)$.
Note that
for $\nbigm_i\in\nbigc(X_i;H_i)$,
we have
$\nbigm_1\boxtimes\nbigm_2\simeq
p_1^{\ast}(\nbigm_1)\otimes^L_{\nbigo_{\cnum\times X_1\times X_2}}
p_2^{\ast}(\nbigm_2)$.

\subsection{Duality}

Let $d_X:=\dim X$.
For any coherent $\nbigr_{X(\ast H)}$-module $\nbigm$,
we obtain the following objects
in the derived category of complexes
of $\nbigr_{X(\ast H)}$-modules
(see \cite[\S13.1]{Mochizuki-MTM}):
\[
 \DD_{X(\ast H)}\nbigm:=
 \nrhom_{\nbigr_{X(\ast H)}}\bigl(
 \nbigm,
 \nbigr_{X(\ast H)}\otimes
 (\lambda^{d_X}\omega_{\nbigx}^{-1})
 \bigr)[d_X].
\]
If $\nbigm$ is an $\nbigrtilde_{X(\ast H)}$-module,
then $\DD_{X(\ast H)}\nbigm$ is naturally an object
of the derived category of
complexes of $\nbigrtilde_{X(\ast H)}$-modules
\cite[\S13.1.6]{Mochizuki-MTM}.
For $\nbigm\in\nbigc(X;H)$,
we obtain $\DD_{X(\ast H)}(\nbigm)\in\nbigc(X;H)$.
For $\nbigm\in\nbigc(X)$,
there exists a natural isomorphism
$\DD_{X(\ast H)}(\nbigm(\ast H))
\simeq
\bigl(
 \DD_X(\nbigm)
 \bigr)(\ast H)$.

\subsubsection{Compatibility with the duality of $\nbigd$-modules
on $\nbigx^{\circ}$}
 
Let $\DDD_{\nbigx^{\circ}}(\nbign)$
denote the duality functor for
$\nbigd_{\nbigx^{\circ}}$-modules.
For a holonomic $\nbigd_{\nbigx^{\circ}}$-module $\nbign$
such that
$\nbign(\ast\nbigh^{\circ})=\nbign$,
we set
$\DDD_{\nbigx^{\circ}(\ast\nbigh^{\circ})}(\nbign):=
 \DDD_{\nbigx^{\circ}}(\nbign)(\ast\nbigh^{\circ})$.
 
\begin{prop}
\label{prop;21.3.23.10}
For $\nbigm\in\nbigc(X;H)$,
$\DD_{X(\ast H)}(\nbigm)_{|\nbigx^{\circ}}$
is naturally isomorphic to
$\DDD_{\nbigx^{\circ}(\ast\nbigh^{\circ})}(\nbigm_{|\nbigx^{\circ}})$.
\end{prop}
\pf
We indicate only an outline.
The natural left action of $\nbigr_{X(\ast H)}$
on $\nbigr_{X(\ast H)}
 \otimes_{\nbigo_{\nbigx}}
 (\lambda^{d_X}\omega_{\nbigx}^{-1})$
by $\ell$.
The right action of $\nbigr_{X(\ast H)}$ on $\lambda^{d_X}\nbigr_{X(\ast H)}$
induces a left $\nbigr_{X(\ast H)}$-action on
$\nbigr_{X(\ast H)}\otimes\lambda^{d_X}\omega_{\nbigx}^{-1}$,
which is denoted by $r$.
There exists the natural action of the differential operator
$\lambda^2\del_{\lambda}$
on 
$\nbigr^{\otimes 2}_{X(\ast H)}:=
\nbigr_{X(\ast H)}
\otimes_{\nbigo_{\cnum_{\lambda}}}
\nbigr_{X(\ast H)}$.
We set
\[
\nbigr_{X(\ast H)}^{\otimes 2}
\langle\lambda^2\del_{\lambda}\rangle
=\bigoplus_{j=0}^{\infty}
(\lambda^2\del_{\lambda})^j\otimes
\nbigr_{X(\ast H)}^{\otimes 2}.
\]
With the multiplication induced by
$(1\otimes f)\cdot
(\lambda^2\del_{\lambda}\otimes 1)
=
\lambda^2\del_{\lambda}\otimes f
-1\otimes(\lambda^2\del_{\lambda}f)$
and
$(\lambda^2\del_{\lambda}\otimes 1)(1\otimes f)
=\lambda^2\del_{\lambda}\otimes f$
for sections $f$ of $\nbigr_{X(\ast H)}^{\otimes 2}$,
it is a sheaf of algebras.
The 
$\nbigr_{X(\ast H)}^{\otimes 2}$-action
$(\ell,r)$
on $\nbigr_{X(\ast H)}\otimes\lambda^{d_X}\omega_{\nbigx}^{-1}$
naturally extends to the action of
$\nbigr_{X(\ast H)}^{\otimes 2}
\langle\lambda^2\del_{\lambda}\rangle$.
Let $(\nbigg^{\bullet},\ell,r)$ be an injective
$\nbigr_{X(\ast H)}^{\otimes 2}
\langle\lambda^2\del_{\lambda}\rangle$-resolution of
$\nbigr_{X(\ast H)}\otimes\lambda^{d_X}\omega_{\nbigx}^{-1}$.
Recall that
$\DD_{X(\ast H)}(\nbigm)=
 \nhom_{\nbigr_{X(\ast H)}}\Bigl(
 \nbigm,
 \nbigg^{\bullet,\ell,r}
 \Bigr)[d_X]$,
 where
 $\nhom_{\nbigr_{X(\ast H)}}\Bigl(
 \nbigm,(\nbigg^j)^{\ell,r}
 \Bigr)$
 denote the sheaf of
$\nbigr_{X(\ast H)}$-homomorphisms from
$\nbigm$ to $(\nbigg^j,\ell)$,
which is naturally an $\nbigrtilde_{X(\ast H)}$-module
induced by the action $r$.

We set $\nbigd_{\nbigx^{\circ}(\ast\nbigh^{\circ})}:=
\nbigd_{\nbigx^{\circ}}(\ast\nbigh^{{\circ}})$.
Let $\Omega_{\nbigx^{\circ}}$ denote the sheaf of
holomorphic $(d_X+1)$-forms on $\nbigx^{\circ}$.
The natural left
$\nbigd_{\nbigx^{\circ}(\ast\nbigh^{\circ})}$-action on
$\nbigd_{\nbigx^{\circ}(\ast\nbigh^{\circ})}
 \otimes_{\nbigo_{\nbigx^{\circ}}}
 \Omega_{\nbigx^{\circ}}^{-1}$
is denoted by $\ell$.
The natural right $\nbigd_{\nbigx^{\circ}(\ast\nbigh^{\circ})}$-action
on $\nbigd_{\nbigx^{\circ}(\ast\nbigh^{\circ})}$
induces a left $\nbigd_{\nbigx^{\circ}(\ast\nbigh^{\circ})}$-action
on $\nbigd_{\nbigx^{\circ}(\ast\nbigh^{\circ})}
\otimes\Omega_{\nbigx^{\circ}}^{-1}$
which is denoted by $r$.
Let $(\nbigl^{\bullet},\ell,r)$ be an injective
$\nbigd_{\nbigx^{\circ}(\ast\nbigh^{\circ})}
\otimes_{\cnum}
\nbigd_{\nbigx^{\circ}(\ast\nbigh^{\circ})}$-resolution of
$\nbigd_{\nbigx^{\circ}(\ast\nbigh^{\circ})}
\otimes\Omega_{\nbigx^{\circ}}^{-1}$.
We set $\nbigm^{\circ}:=\nbigm_{|\nbigx^{\circ}}$.
\begin{lem}
\label{lem;21.4.13.4}
 We have
$\DDD_{\nbigx^{\circ}(\ast\nbigh^{\circ})}(\nbigm^{\circ})=
 \nhom_{\nbigd_{\nbigx^{\circ}}(\ast\nbigh^{\circ})}\Bigl(
 \nbigm^{\circ},
 \nbigl^{\bullet,\ell,r}
 \Bigr)[d_X+1]$.
\end{lem}
\pf
Let $(\nbigl_1^{\bullet},\ell,r)$
denote an injective
$(\nbigd_{\nbigx^{\circ}}\otimes_{\cnum}
\nbigd_{\nbigx^{\circ}})$-resolution
of $\nbigd_{\nbigx^{\circ}}\otimes\Omega_{\nbigx^{\circ}}^{-1}$.
We have
$\DDD_{\nbigx^{\circ}}(\nbigm^{\circ})
=\nhom_{\nbigd_{\nbigx^{\circ}}}
\Bigl(
\nbigm^{\circ},\nbigl_1^{\bullet,\ell,r}
\Bigr)[d_X+1]$.
By using the $\nbigo_{\nbigx^{\circ}}$-action
on $\nbigl^{\circ}_1$ induced by $r$,
we define
\[
 \nbigl_2^{\bullet}:=
  \nbigo_{\nbigx^{\circ}}(\ast\nbigh^{\circ})
 \otimes_{\nbigo_{\nbigx^{\circ}}}^{r}
\nbigl_1^{\bullet}.
\]
By using the $\nbigo_{\nbigx^{\circ}}$-action
on $\nbigl^{\circ}_1$ induced by $\ell$,
we define
\[
 \nbigl_3^{\bullet}:=
\nbigo_{\nbigx^{\circ}}(\ast\nbigh^{\circ})
\otimes_{\nbigo_{\nbigx}^{\circ}}^{\ell}
\nbigl_2^{\bullet}.
\]
They are resolutions
of $\nbigd_{\nbigx^{\circ}(\ast\nbigh^{\circ})}
\otimes\Omega_{\nbigx^{\circ}}^{-1}$.
There exists a  quasi-isomorphism
$(\nbigl_3^{\bullet},\ell,r)
\lrarr(\nbigl^{\bullet},\ell,r)$.
There exist the following natural morphisms:
\begin{multline}
\label{eq;21.4.13.2}
 \DDD_{\nbigx^{\circ}}(\nbigm^{\circ})(\ast \nbigh^{\circ})
 =
 \nhom_{\nbigd_{\nbigx^{\circ}}}
 \bigl(
 \nbigm^{\circ},\nbigl_2^{\bullet,\ell,r}
 \bigr)[d_X+1]
 \lrarr
  \nhom_{\nbigd_{\nbigx^{\circ}(\ast\nbigh^{\circ})}}
 \bigl(
 \nbigm^{\circ},\nbigl_3^{\bullet,\ell,r}
 \bigr)[d_X+1] \\
 \lrarr
 \nhom_{\nbigd_{\nbigx^{\circ}(\ast\nbigh^{\circ})}}
 \bigl(
 \nbigm^{\circ},
 \nbigl^{\bullet,\ell,r}
 \bigr)[d_X+1].
\end{multline}
Let us check that the composition of (\ref{eq;21.4.13.2})
is a quasi-isomorphism.
We may assume that there exists a coherent free $\nbigd_{\nbigx^{\circ}}$-resolution
$\nbigp^{\bullet}$ of $\nbigm$.
We obtain the induced free
$\nbigd_{\nbigx^{\circ}(\ast\nbigh^{\circ})}$-resolution
$\nbigp^{\bullet}(\ast\nbigh^{\circ})$ of $\nbigm$.
It is enough to check that the induced morphism
\[
\nhom_{\nbigd_{\nbigx^{\circ}}}(\nbigp^{\bullet},\nbigl_2^{\bullet,\ell,r})[d_X+1]
\lrarr
\nhom_{\nbigd_{\nbigx^{\circ}(\ast\nbigh^{\circ})}}(\nbigp^{\bullet}(\ast\nbigh^{\circ}),
 \nbigl^{\bullet,\ell,r})[d_X+1]
\]
is a quasi-isomorphism,
which is easy to see.
Thus, we obtain Lemma \ref{lem;21.4.13.4}.
\hfill\qed

\vspace{.1in}

Let $(\nbign,\rho)$
be a $\nbigd_{\nbigx^{\circ}(\ast\nbigh^{\circ})}$-module.
We set
$\nbign\langle\del_{\lambda}\rangle
:=\bigoplus_{j=0}^{\infty}
\del_{\lambda}^j\otimes\nbign$.
The $\nbigd_{\nbigx^{\circ}(\ast\nbigh^{\circ})}$-action $\rho$
on $\nbign$ extends to
a $\nbigd_{\nbigx^{\circ}(\ast\nbigh^{\circ})}$-action $\rho$
on $\nbign\langle\del_{\lambda}\rangle$
as follows:
\[
f\bullet_{\rho} (\del_{\lambda}^j\otimes m)
=\sum_{k=0}^j(-1)^kC(j,k)\del_{\lambda}^{j-k}\otimes
 (\del_{\lambda}^kf\bullet_{\rho} m),
\quad
v\bullet_{\rho} (\del_{\lambda}^j\otimes m)
=\del_{\lambda}^j\otimes(v\bullet_{\rho} m),
\quad
\del_{\lambda}\bullet_{\rho} (\del_{\lambda}^j\otimes m)
=\del_{\lambda}^{j+1}\otimes m.
\]
Here, $C(j,k)$ denote the binomial coefficients,
and $v$ denote local sections of $\nbigd_{\nbigx^{\circ}/\lambda}(\ast\nbigh^{\circ})$.
We set
$C^i(\nbign):=\nbign\langle\del_{\lambda}\rangle$ $(i=0,-1)$.
Let $C^{-1}(\nbign)\lrarr C^{0}(\nbign)$
be the $\nbigd_{\nbigx^{\circ}(\ast\nbigh^{\circ})}$-morphism induced by
\[
 \sum \del_{\lambda}^j\otimes m_j
\longmapsto
-\sum \del_{\lambda}^{j+1}\otimes m_j
+\sum \del_{\lambda}^j\otimes (\del_{\lambda}\bullet_{\rho}m_j).
\]
Thus, we obtain a $\nbigd_{\nbigx^{\circ}(\ast\nbigh^{\circ})}$-complex
$C^{\bullet}(\nbign)$.
Let $C^0(\nbign)\lrarr \nbign$
be the $\nbigd_{\nbigx^{\circ}(\ast\nbigh^{\circ})}$-homomorphism
defined by
$\sum \del_{\lambda}^j\otimes m_j\longmapsto
\sum \del_{\lambda}^j\bullet_{\rho} m_j$,
which induces a quasi-isomorphism
$C^{\bullet}(\nbign)\simeq \nbign$.

We set
$\nbigd_{\nbigx^{\circ}/\lambda}(\ast\nbigh^{\circ})^{\otimes 2}:=
\nbigd_{\nbigx^{\circ}/\lambda}(\ast\nbigh^{\circ})
 \otimes_{\nbigo_{\cnum^{\ast}_{\lambda}}}
 \nbigd_{\nbigx^{\circ}/\lambda}(\ast\nbigh^{\circ})$.
(See \S\ref{subsection;22.7.30.1}
for $\nbigd_{\nbigx^{\circ}/\lambda}$.)
 As in the case of
$\nbigr_{X(\ast H)}^{\otimes 2}\langle\lambda^2\del_{\lambda}\rangle$,
$\nbigd_{\nbigx^{\circ}/\lambda}(\ast\nbigh^{\circ})^{\otimes 2}
\langle\del_{\lambda}\rangle$
is naturally a sheaf of algebras.
Let $(\nbign,\rho_1,\rho_2)$ be a
$\nbigd_{\nbigx^{\circ}/\lambda}(\ast\nbigh^{\circ})^{\otimes 2}
\langle\del_{\lambda}\rangle$-module,
where $\rho_1$ and $\rho_2$
denote the underlying two commuting $\nbigd_{\nbigx^{\circ}/\lambda}$-actions
on $\nbign$.
 We set
$\nbign\langle \del_{\lambda}\rangle
=\bigoplus_{j=0}^{\infty}
 \del_{\lambda}^j\otimes\nbign$.
The $\nbigd_{\nbigx^{\circ}/\lambda}(\ast\nbigh^{\circ})$-action
$\rho_1$ extends
to a $\nbigd_{\nbigx^{\circ}(\ast\nbigh^{\circ})}$-action $\rho_1$
on $\nbign\langle\del_{\lambda}\rangle$
as follows:
\[
 f\bullet_{\rho_1} (\del_{\lambda}^j\otimes m)
=\sum_{k=0}^j(-1)^kC(j,k)
 \del_{\lambda}^{j-k}\otimes(\del_{\lambda}^kf\bullet_{\rho_1}m),
 \quad
 v\bullet_{\rho_1}(\del_{\lambda}^j\otimes m)
 =\del_{\lambda}^j\otimes(v\bullet_{\rho_1}m),
 \quad
 \del_{\lambda}\bullet_{\rho_1}(\del_{\lambda}^j\otimes m)
 =
\del_{\lambda}^{j+1}\otimes m.
\]
The $\nbigd_{\nbigx^{\circ}/\cnum^{\ast}}(\ast\nbigh^{\circ})$-action $\rho_2$
extends to a $\nbigd_{\nbigx^{\circ}(\ast\nbigh^{\circ})}$-action $\rho_2$
on $\nbign\langle\del_{\lambda}\rangle$
as follows:
\[
 f\bullet_{\rho_2} (\del_{\lambda}^j\otimes m)
=\del_{\lambda}^j\otimes (f\bullet_{\rho_2}m),
 \quad
 v\bullet_{\rho_2}(\del_{\lambda}^j\otimes m)
=\del_{\lambda}^j\otimes(v\bullet_{\rho_2}m),
 \quad
 \del_{\lambda}\bullet_{\rho_2}(\del_{\lambda}^j\otimes m)
=-\del_{\lambda}^{j+1}\otimes m
 +\del_{\lambda}^j\otimes (\del_{\lambda}\bullet m).
\]
Here, $\del_{\lambda}\bullet m$
is induced by the original 
$\nbigd_{\nbigx^{\circ}/\lambda}(\ast\nbigh^{\circ})^{\otimes 2}
\langle\del_{\lambda}\rangle$-action on $\nbign$.
The two $\nbigd_{\nbigx^{\circ}}(\ast\nbigh^{\circ})$-actions
are commutative.

Let $\nbigg_0^{\bullet}:=\nbigg^{\bullet}_{|\nbigx^{\circ}}$.
By the multiplication of  $d\lambda$,
we identify
$\omega_{\nbigx|\nbigx^{\circ}}$
with $\Omega_{\nbigx^{\circ}}$.
Note that there exists a natural isomorphism
\[
 \bigl(
 (\nbigr_{X(\ast H)}\otimes\lambda^{d_X}\omega_{\nbigx}^{-1})
  _{|\nbigx^{\circ}}\langle\del_{\lambda}\rangle,
  \ell,r
 \bigr)
\simeq
 \bigl(\nbigd_{\nbigx^{\circ}(\ast\nbigh^{\circ})}
 \otimes\Omega_{\nbigx^{\circ}}^{-1},
 \ell,r
 \bigr).
\]
Hence,
$(\nbigg_0^{\bullet}\langle\del_{\lambda}\rangle,\ell,r)$
is 
a $(\nbigd_{\nbigx^{\circ}(\ast\nbigh^{\circ})}
 \otimes_{\cnum}\nbigd_{\nbigx^{\circ}(\ast\nbigh^{\circ})})$-resolution
 of $\nbigd_{\nbigx^{\circ}(\ast\nbigh^{\circ})}\otimes
 \Omega_{\nbigx^{\circ}}^{-1}$.
There exists a quasi-isomorphism
$(\nbigg_0^{\bullet}\langle\del_{\lambda}\rangle,\ell,r)
\lrarr
(\nbigl^{\bullet},\ell,r)$.

There exists the natural quasi-isomorphism
\[
 \nhom_{\nbigd_{\nbigx^{\circ}}(\ast\nbigh^{\circ})}
 \bigl(
 \nbigm^{\circ},\nbigl^{\bullet,\ell,r}
 \bigr)[d_X+1]
\lrarr
 \nhom_{\nbigd_{\nbigx^{\circ}}(\ast\nbigh^{\circ})}
 \bigl(
 C^{\bullet}(\nbigm^{\circ}),\nbigl^{\bullet,\ell,r}
 \bigr)[d_X+1].
\]
There exists the following natural morphism
\begin{equation}
\label{eq;21.3.23.1}
 \nhom_{\nbigd_{\nbigx^{\circ}}(\ast\nbigh^{\circ})}\bigl(
  C^{\bullet}(\nbigm^{\circ}),
  \nbigg_0^{\bullet,\ell,r}\langle\del_{\lambda}\rangle
  \bigr)[d_X+1]
  \lrarr
  \nhom_{\nbigd_{\nbigx^{\circ}}(\ast\nbigh^{\circ})}
 \bigl(
 C^{\bullet}(\nbigm^{\circ}),\nbigl^{\bullet,\ell,r}
 \bigr)[d_X+1].
\end{equation}
There exist the following natural isomorphism
of $\nbigd_{\nbigx^{\circ}(\ast\nbigh^{\circ})}$-modules for $i=0,-1$.
\[
 \nhom_{\nbigd_{\nbigx^{\circ}}(\ast\nbigh^{\circ})}
 \bigl(
 C^{i}(\nbigm^{\circ}),
\nbigg_0^{\bullet,\ell,r}\langle\del_{\lambda}\rangle
 \bigr)
=\nhom_{\nbigd_{\nbigx^{\circ}}(\ast\nbigh^{\circ})}
 \bigl(
 \nbigm^{\circ}\langle\del_{\lambda}\rangle,
\nbigg_0^{\bullet,\ell,r}\langle\del_{\lambda}\rangle
 \bigr)
 \simeq
  \nhom_{\nbigd_{\nbigx^{\circ}/\lambda}(\ast\nbigh^{\circ})}
 \bigl(
 \nbigm^{\circ},
\nbigg_0^{\bullet,\ell,r}\langle\del_{\lambda}\rangle
 \bigr).
\]
The morphism
\[
 \nhom_{\nbigd_{\nbigx^{\circ}}(\ast\nbigh^{\circ})}
 \bigl(
 C^{0}(\nbigm^{\circ}),
\nbigg_0^{\bullet,\ell,r}\langle\del_{\lambda}\rangle
 \bigr)
 \lrarr
 \nhom_{\nbigd_{\nbigx^{\circ}}(\ast\nbigh^{\circ})}
 \bigl(
 C^{-1}(\nbigm^{\circ}),
\nbigg_0^{\bullet,\ell,r}\langle\del_{\lambda}\rangle
 \bigr)
\]
is identified with the morphism
\[
\Phi:
  \nhom_{\nbigd_{\nbigx^{\circ}/\lambda}(\ast\nbigh^{\circ})}
 \bigl(
 \nbigm^{\circ},
\nbigg_0^{\bullet,\ell,r}\langle\del_{\lambda}\rangle
 \bigr)
 \lrarr
   \nhom_{\nbigd_{\nbigx^{\circ}/\lambda}(\ast\nbigh^{\circ})}
 \bigl(
 \nbigm^{\circ},
\nbigg_0^{\bullet,\ell,r}\langle\del_{\lambda}\rangle
 \bigr)
\]
given as
\[
 \Phi(F)(m)
=\del_{\lambda}\bullet_{\ell}(F(m))
-F\bigl(\del_{\lambda}m\bigr)
\]
for $F\in \nhom_{\nbigd_{\nbigx^{\circ}/\lambda}(\ast\nbigh^{\circ})}
 \bigl(
 \nbigm^{\circ},
\nbigg_0^{\bullet,\ell,r}\langle\del_{\lambda}\rangle
 \bigr)$
 and $m\in \nbigm^{\circ}$.
(See \cite[Remark 1.8.11]{Kashiwara-Schapira}
for the signature.)
Because $\nbigm^{\circ}$ is
coherent over $\nbigd_{\nbigx^{\circ}/\lambda}(\ast\nbigh^{\circ})$,
we have the natural isomorphism
\begin{equation}
\label{eq;21.4.13.5}
 \nhom_{\nbigd_{\nbigx^{\circ}/\lambda}(\ast\nbigh^{\circ})}
 \bigl(
 \nbigm^{\circ},
\nbigg_0^{\bullet,\ell,r}\langle\del_{\lambda}\rangle
 \bigr)
=
 \nhom_{\nbigd_{\nbigx^{\circ}/\lambda}(\ast\nbigh^{\circ})}
 \bigl(
 \nbigm^{\circ},
\nbigg_0^{\bullet,\ell,r}
 \bigr)
 \langle\del_{\lambda}\rangle.
\end{equation}

\begin{lem}
\label{lem;21.4.13.3}
 The morphism {\rm(\ref{eq;21.3.23.1})}
is a quasi-isomorphism. 
\end{lem}
\pf
It is enough to prove the claim
on a neighbourhood $U$ of any point $P$ of $\nbigx^{\circ}$.
We may assume that there exists
a coherent free $\nbigr_{X(\ast H)|U}$-resolution $\nbigp^{\bullet}\lrarr\nbigm_{|U}$.
We set $\nbigh_U^{\circ}:=U\cap\nbigh^{\circ}$.
We obtain the following commutative diagram:
\[
 \begin{CD}
   \nhom_{\nbigd_{\nbigx^{\circ}}(\ast\nbigh^{\circ})}\bigl(
  C^{\bullet}(\nbigm^{\circ}),
  \nbigg_0^{\bullet,\ell,r}\langle\del_{\lambda}\rangle
  \bigr)[d_X+1]_{|U}
@>{a_0}>>
  \nhom_{\nbigd_{\nbigx^{\circ}}(\ast\nbigh^{\circ})}
 \bigl(
 C^{\bullet}(\nbigm^{\circ}),\nbigl^{\bullet,\ell,r}
  \bigr)[d_X+1]_{|U}
  \\
  @V{a_1}VV @V{a_2}VV \\
   \nhom_{\nbigd_{U}(\ast\nbigh^{\circ}_U)}\bigl(
  C^{\bullet}(\nbigp^{\bullet}),
  \nbigg_0^{\bullet,\ell,r}\langle\del_{\lambda}\rangle_{|U}
  \bigr)[d_X+1]
@>{a_3}>>
  \nhom_{\nbigd_{U}(\ast\nbigh^{\circ}_U)}
 \bigl(
 C^{\bullet}(\nbigp^{\bullet}),\nbigl^{\bullet,\ell,r}_{|U}
  \bigr)[d_X+1].
 \end{CD}
\]
Note that $C^{\bullet}(\nbigp^{\bullet})$
is a complex of coherent free $\nbigd_{U}$-modules.
Because
\[
 \nhom_{\nbigd_{\nbigx^{\circ}/\lambda}(\ast\nbigh^{\circ})}
 \bigl(
 \nbigm^{\circ},
\nbigg_0^{\bullet,\ell,r}
 \bigr)_{|U}
 \lrarr
  \nhom_{\nbigd_{U/\lambda}(\ast\nbigh^{\circ}_U)}
 \bigl(
\nbigp^{\bullet},
\nbigg^{\bullet,\ell,r}_{0|U}
 \bigr)
\]
is a quasi-isomorphism,
we obtain that $a_1$ is a quasi-isomorphism
by using (\ref{eq;21.4.13.5}).
Because $a_2$ and $a_3$ are quasi-isomorphism,
we obtain that $a_0$ is a quasi-isomorphism,
which is the claim of Lemma \ref{lem;21.4.13.3}.
\hfill\qed

\vspace{.1in}

The cokernel of $\Phi$
as a $\nbigd_{\nbigx^{\circ}/\lambda}(\ast\nbigh^{\circ})$-module
is naturally isomorphic to
$\DD(\nbigm)_{|\nbigx^{\circ}}$.
We can check that it is compatible
with the induced actions of $\del_{\lambda}$.
Thus, we obtain Proposition \ref{prop;21.3.23.10}.
\hfill\qed

\subsubsection{Duality and regularity}

Let $Z$ be an open neighbourhood of $0$ in $\cnum_t$.
Suppose $X=Z\times X_1$.
Let $\pi:X\lrarr X_1$ be the projection.
We set
$\nbiga_0:=\pi^{\ast}(\nbigr_{X_1})(\ast H)
\subset \nbigr_{X(\ast H)}$.
We obtain a sheaf of algebras
$\nbiga_1:=
\bigl(\nbiga_0\otimes_{\nbigo_{\nbigz}}\nbiga_0
\bigr)
\langle
\lambda^2\del_{\lambda},\deldel_t
\rangle$
as in the case of
$\nbigr_{X(\ast H)}^{\otimes 2}\langle\lambda^2\del_{\lambda}\rangle$.
We set 
$\nbigb:=\pi^{\ast}(\nbigr_{X_1}\otimes\omega_{\nbigx_1}^{-1})(\ast H)$.
Then, $\nbigb(\ast t)$ is naturally
an $\nbiga_1(\ast t)$-module.
Let $\nbigg^{\bullet}_1$ be
an $\nbiga_1(\ast t)$-injective resolution of $\nbigb(\ast t)$.
For any $\nbigrtilde_{X(\ast H)}(\ast t)$-module $\nbign$
which is coherent over $\nbiga_0(\ast t)$,
we obtain the following complex of
$\nbigrtilde_{X(\ast H)}(\ast t)$-modules:
\[
 \nrhom_{\nbiga_0(\ast t)}
 \Bigl(
 \nbign,
 \lambda^{d_{X_1}}\nbigb(\ast t)[d_{X_1}]
  \Bigr)
:=
 \nhom_{\nbiga_0(\ast t)}
 \Bigl(
 \nbign,\lambda^{d_{X_1}}\nbigg_1^{\bullet}[d_{X_1}]
 \Bigr).
\]
\begin{lem}
\label{lem;21.3.30.1}
Let $\nbign$ be an $\nbigrtilde_{X(\ast H)}$-module
such that
 (i) $\nbign$ is coherent over $\nbigr_{X(\ast H)}$,
 (ii) $\nbign(\ast t)$ is coherent over
$\nbiga_0(\ast t)$.
There exists the following natural quasi-isomorphism
of $\nbigrtilde_{X}(\ast t)$-complexes:
 \[
 \DD_{X(\ast H)}(\nbign)(\ast t)
 \simeq
 \lambda\cdot
 \nrhom_{\nbiga_0(\ast t)}
 \Bigl(
 \nbign(\ast t),\,
 \lambda^{d_{X_1}}\nbigb(\ast t)
 [d_{X_1}]
 \Bigr).
\]
\end{lem}
\pf
Similar to Proposition \ref{prop;21.3.23.10}.
\hfill\qed

\vspace{.1in}

Let $\nbigm\in\nbigc(X;H)$.
Suppose that $\nbigm$ is regular along $t$,
i.e.,
each $V_a(\nbigm)$ is $\nbiga_0$-coherent,
which implies that $\nbigm(\ast t)=V_a(\nbigm)(\ast t)$ is
$\nbiga_0(\ast t)$-coherent.
According to Lemma \ref{lem;21.3.30.1},
there exists the following natural isomorphism:
\[
 \DD_{X(\ast H)}(\nbigm)(\ast t)
 \simeq
 \lambda\cdot
 \nrhom_{\nbiga_0(\ast t)}
 \bigl(
 \nbigm(\ast t),
 \lambda^{d_{X_1}}\nbigb(\ast t)
 \bigr)[d_{X_1}].
\]

Let $\nbigg_2^{\bullet}$
be an $\nbiga_1$-injective resolution of
$\nbigb$.
We obtain the following complex
of $V\nbigrtilde_{X(\ast H)}$-modules:
\[
  \nrhom_{\nbiga_0}
 \Bigl(
 V_a(\nbigm),\,
 \lambda^{d_{X_1}}
 \nbigb[d_{X_1}]
 \Bigr)
 :=
 \nhom_{\nbiga_0}
 \Bigl(
 V_a(\nbigm),\,
 \lambda^{d_{X_1}}
 \nbigg_2^{\bullet}[d_{X_1}]
 \Bigr).
\]
There exists the natural morphism:
\[
 \lambda\cdot\nrhom_{\nbiga_0}
 \Bigl(
 V_a(\nbigm),\,
 \lambda^{d_{X_1}}
 \nbigb[d_{X_1}]
 \Bigr)
 \lrarr
 \DD_{X(\ast H)}(\nbigm)(\ast t).
\]

\begin{prop}
\label{prop;21.3.30.3}
If $\nbigm$ is regular along $t$,
$\DD_{X(\ast H)}(\nbigm)(\ast t)$ is regular along $t$.
There exists the following natural isomorphism
for any $a\in\real$:
\[
 V_a\bigl(\DD_{X(\ast H)}\nbigm(\ast t) \bigr)\simeq
 \lambda\cdot \nrhom_{\nbiga_0}
 \Bigl(
 V_{<-1-a}(\nbigm),\,
 \lambda^{d_{X_1}}\nbigb[d_{X_1}]
 \Bigr).
\]
\end{prop}
\pf
It is easy to check the following lemma.
\begin{lem}
Let $\nbigl$ be a coherent $\nbiga_0$-module.
Then, the natural morphism
\[
 \nrhom_{\nbiga_0}
  \bigl(
 \nbigl,
 \nbigb(\ast t)
  \bigr)
\lrarr 
 \nrhom_{\nbiga_0(\ast t)}
  \bigl(
 \nbigl(\ast t),
 \nbigb(\ast t)
  \bigr)
\]
is a quasi-isomorphism.
\hfill\qed
\end{lem}

We set $\nbign:=\nbigm(\ast t)$.
There exists the following quasi-isomorphism:
\[
 \nrhom_{\nbiga_0}
 \bigl(
 V_a(\nbign),
 \nbigb(\ast t)
 \bigr)[d_{X_1}]
\simeq
 \nrhom_{\nbiga_0(\ast t)}
 \bigl(
 \nbign,
 \nbigb(\ast t)
 \bigr)[d_{X_1}].
\]
We have already observed that
\[
 \nbigh^j\Bigl(
 \lambda\nrhom_{\nbiga_0(\ast t)}
 \bigl(
 \nbign,
 \lambda^{d_{X_1}}
 \nbigb(\ast t)
 \bigr)[d_{X_1}]
 \Bigr)
 =\left\{
\begin{array}{ll}
 \DD_X(\nbigm)(\ast t)
  & (j=0)
  \\
 0 & (j\neq 0).
\end{array}
 \right.
\]

Let $\iota_0:X_1\simeq \{0\}\times X_1\lrarr X$
denote the inclusion.
We set $H_1:=\iota_0^{-1}(H)$.

\begin{lem}
\label{lem;22.7.2.1}
\[
 \nbigh^j\Bigl(
 \nrhom_{\nbiga_0}
 \bigl(
 \Gr^V_a(\nbign),\,
 \lambda^{d_{X_1}}\nbigb [d_{X_1}]
 \bigr)
 =\left\{
\begin{array}{ll}
 \iota_{0\ast}
 \DD_{X_1(\ast H_1)}(\iota_0^{\ast}\Gr^V_a\nbign)
 & (j=1)\\
 0 &(j\neq 1).
\end{array}
 \right.
\]
\end{lem}
\pf
There exists $\nbign_{a,0}\in \nbigc(X_1)$
such that
$\iota_{0\ast}\bigl(
\nbign_{a,0}(\ast H_1)
\bigr)=
\Gr^{V}_a(\nbign)$.
There exists the following natural isomorphism
\[
 \nrhom_{\pi^{\ast}\nbigr_{X_1}}
 \bigl(
 \pi^{\ast}(\nbign_{a,0}),
 \pi^{\ast}(\lambda^{d_{X_1}}\nbigr_{X_1}\otimes\omega_{\nbigx_1}^{-1})[d_{X_1}]
 \bigr)
 \simeq
 \pi^{\ast}\bigl(
 \DD_{X_1}(\nbign_{a,0})
 \bigr).
\]
There exists the $\pi^{\ast}\nbigr_{X_1}$-submodule
$t\cdot\pi^{\ast}(\nbign_{a,0})
\subset
 \pi^{\ast}(\nbign_{a,0})$,
and the quotient is isomorphic to
$\iota_{0\ast}(\nbign_{a,0})$.
Hence, we obtain
\[
\nbigh^j\Bigl(
  \nrhom_{\pi^{\ast}\nbigr_{X_1}}
 \bigl(
 \iota_{0\ast}(\nbign_{a,0}),\,
 \pi^{\ast}(\lambda^{d_{X_1}}\nbigr_{X_1}\otimes\omega_{\nbigx_1}^{-1})[d_{X_1}]
 \bigr)
 \Bigr)
 \Bigr)
 \simeq
 \left\{
 \begin{array}{ll}
 \iota_{0\ast}\DD_{X_1}(\nbign_{a,0})
  & (j=1) \\
  0 & (j\neq 1).
 \end{array}
 \right.
\]
Then, we obtain the claim of Lemma \ref{lem;22.7.2.1}.
\hfill\qed

\begin{lem}
\label{lem;22.7.2.2}
For $j\neq 0$,
 \[
 \nbigh^j\Bigl(
 \nrhom_{\nbiga_0}
 \bigl(
 V_a(\nbign),\,\,
 (\nbigb/t\nbigb)
 \bigr)
 [d_{X_1}]
 \Bigr)
=0.
\]
As a result,
we obtain the following vanishing for $j\neq 0$:
\[
\nbigh^j\Bigl(
 \nrhom_{\nbiga_0}
\bigl(
 V_a(\nbign),\,\,
 \nbigb(\ast t)/\nbigb
\bigr)[d_{X_1}]
\Bigr)
=0.
\]
\end{lem}
\pf
We have
$\iota_{0\ast}(\iota_0^{\ast}\nbigb)=\nbigb/t\nbigb$.
There exist the following morphisms:
\begin{multline}
\label{eq;21.3.30.2}
 \iota_{0\ast}\Bigl(
 \nrhom_{\iota_0^{\ast}\nbiga_0}\bigl(
 \iota_0^{\ast}V_a(\nbign),\,
 \iota_0^{\ast}\nbigb[d_{X_1}]
 \bigr)
 \Bigr)
 \lrarr
 \nrhom_{\iota_{0\ast}\iota_0^{\ast}\nbiga_{0}}\bigl(
 \iota_{0\ast}\iota_0^{\ast}V_a(\nbign),\,
  (\nbigb/t\nbigb) [d_{X_1}]
 \bigr)
 \\
\lrarr
 \nrhom_{\nbiga_0}\bigl(
 V_a(\nbign),\,
 (\nbigb/t\nbigb) [d_{X_1}]
 \bigr).
\end{multline}
We can check that the composite of (\ref{eq;21.3.30.2})
is a quasi-isomorphism
locally around any point of $X$
by using an $\nbiga_0$-free resolution of $V_a(\nbign)$.
We have
$\iota_0^{\ast}V_a(\nbign)=
\iota_0^{\ast}(V_{a}(\nbign)/V_{a-1}(\nbign))$.
Then, we obtain Lemma \ref{lem;22.7.2.2}.
\hfill\qed

\begin{lem}
We have the following vanishing for $j\neq 0$:
\[
 \nbigh^j\Bigl(
 \nrhom_{\nbiga_0}\bigl(
 V_a(\nbign),\,\,
 \nbigb
 \bigr)[d_{X_1}]
 \Bigr)=0.
\]
\end{lem}
\pf
From the exact sequence
\[
 0\lrarr\nbigb\stackrel{t}{\lrarr}\nbigb
 \lrarr \nbigb/t\nbigb
 \lrarr 0,
\]
and the previous lemma,
the induced morphism
\[
 \nbigh^j\Bigl(
 \nrhom_{\nbiga_0}\bigl(
 V_a(\nbign),\,
 \nbigb
 \bigr)[d_{X_1}]
 \Bigr)
\stackrel{t}{\lrarr}
  \nbigh^j\bigl(
 \nrhom_{\nbiga_0}\Bigl(
 V_a(\nbign),
 \nbigb
 \bigr)[d_{X_1}]
 \Bigr) 
\]
is an epimorphism unless $j=0$.
We also remark that
the support of
the coherent $\nbiga_0$-modules
\[
\nbigh^j\Bigl(
 \nrhom_{\nbiga_0}
 \bigl(
 V_a(\nbign),
 \nbigb
 \bigr)[d_{X_1}]
 \Bigr)
 \quad (j\neq 0) 
\]
are contained in $\{t=0\}$.
Hence, for any local section $m$,
there exists $N(m)\in\seisuu_{>0}$ such that $t^{N(m)}m=0$.
If
$\nbigh^j\Bigl(
\nrhom_{\nbiga_0}
\bigl(
V_a(\nbign),
\nbigb
\bigr)[d_{X_1}]
\Bigr)\neq 0$
for some $j\neq 0$,
there exists a non-zero section $s_0$ such that
$t s_0=0$.
We can choose $s_i$ $(i=1,2,\ldots)$ such that
$t s_i=s_{i-1}$, inductively.
Let $\nbign'\subset\nbign$
be the $\nbiga_0$-submodule
generated by $\{s_i\}$.
By the Noetherian property of $\nbiga_0$,
$\nbign'$ is finitely generated.
Hence, there exists $N\in\seisuu_{>0}$ such that
$t^N\nbign'=0$,
which contradicts that
$s_{N}\in\nbign'$
for which $t^Ns_{N}=s_0\neq 0$.
\hfill\qed

\vspace{.1in}
The natural morphism
\[
 \nbigh^0\Bigl(
  \nrhom_{\nbiga_0}
 \bigl(
 V_a(\nbign),
 \nbigb
 \bigr)[d_{X_1}]
 \Bigr)
 \lrarr
\nbigh^0\Bigl(
  \nrhom_{\nbiga_0}
 \Bigl(
 V_a(\nbign),\,
 \nbigb(\ast t)
 \Bigr)[d_{X_1}]
 \Bigr)
\]
is a monomorphism.
The tuple
\[
U_{a}\bigl(\DD_{X(\ast H)}(\nbigm)(\ast t)\bigr)
:=\nbigh^0\Bigl(
  \nrhom_{\nbiga_0}
 \bigl(
 V_{<-1-a}(\nbign),\,
 \nbigb
 \bigr)[d_{X_1}]
 \Bigr)
 \quad (a\in\real)
\]
induces a filtration
by $\nbiga_0$-coherent submodules
such that
\[
U_{a}\bigl(\DD_{X(\ast H)}(\nbigm)(\ast t)\bigr)
 \big/
  U_{<a}\bigl(\DD_{X(\ast H)}(\nbigm)(\ast t) \bigr)
 \simeq
 \iota_{0\ast}\Bigl(
 \DD_{X_1(\ast H_1)}\bigl(
 \iota_0^{\ast}\Gr^V_{-1-a}(\nbign)
 \bigr)
 \Bigr),
\]
which are strict for any $a\in\real$.

Let $g$ be a section of
$U_{a,j}:=
 \nhom_{\nbiga_0}\Bigl(
 V_{<-1-a}(\nbign),\nbigg_2^{j}
 \Bigr)$.
Then, $t\cdot g$ is naturally a section of
$U_{a-1,j}=
  \nhom_{\nbiga_0}\Bigl(
 V_{<-1-(a-1)}(\nbign),\nbigg_2^{j}
 \Bigr)$.
Thus, we obtain
$t\cdot U_a\bigl(\DD_{X(\ast H)}(\nbigm)(\ast t)\bigr)
=U_{a-1}\bigl(\DD_{X(\ast H)}(\nbigm)(\ast t)\bigr)$.
Because
$[\deldel_{t},g]$ is naturally a section of
$U_{a+1,j}$,
 we obtain
 $\deldel_t U_a\bigl(\DD_{X(\ast H)}(\nbigm)(\ast t)\bigr)
\subset U_{a+1}\bigl(\DD_{X(\ast H)}(\nbigm)(\ast t)\bigr)$.
For a section $m$ of
$V_{<-1-a}(\nbign)$,
we have
\[
\Bigl(
\bigl(
-\deldel_{t}t-a\lambda
\bigr)\bullet g\Bigr)(m)
=
 [-\deldel_{t},t g](m)
-\lambda a g(m)
=-\deldel_{t}(g(t m))
+g\bigl((\deldel_{t}t-(1+a)\lambda) m\bigr).
\]
If $N\in\seisuu_{>0}$ is sufficiently large,
for any $m\in V_{-1-a}(\nbign)$,
we obtain $tm\in V_{<-1-a}(\nbign)$
and
$(\deldel_tt-(1+a)\lambda)^Nm\in V_{<-1-a}(\nbign)$.
Hence,
$\bigl(
-\deldel_{t}t-(1+a)\lambda
\bigr)^N\bullet g$
induces a section of
$U_{<a,j}=
 \nhom_{\pi^{\ast}\nbigr_{X_1}}\Bigl(
 V_{-1-a}(\nbign),\nbigg_2^{j}
 \Bigr)$.
Thus, we obtain that
$U_{\bullet}\bigl(\DD_X(\nbigm)(\ast t)\bigr)$
is a $V$-filtration of
$\DD_{X(\ast H)}(\nbigm)(\ast t)$.
\hfill\qed

\subsection{Non-characteristic inverse image}
\label{subsection;21.6.29.11}

Let $F:X\lrarr Y$ be a morphism of complex manifolds.
Let $H_Y$ be a hypersurface of $Y$
such that $H_X:=F^{-1}(H_Y)$ is a hypersurface of $X$.
We set
$\nbigr_{X\rarr Y}:=\nbigo_{\nbigx}\otimes_{F^{-1}\nbigo_{\nbigy}}F^{-1}\nbigr_Y$,
which is naturally
a left $\nbigr_X$-module
and a right $F^{-1}\nbigr_Y$-module.
For an $\nbigr_{Y(\ast H_Y)}$-module $\nbigm$,
we set
$LF^{\ast}(\nbigm):=
\nbigr_{X\rarr Y}(\ast H_X)
\otimes^L_{F^{-1}\nbigr_{Y(\ast H_Y)}}F^{-1}\nbigm$
in the derived category of $\nbigrtilde_{X(\ast H_X)}$-modules.

\begin{df}
\label{df;21.4.12.30}
We say that $F$ is strictly non-characteristic for
$\nbigm\in\nbigc(Y;H_Y)$
if $F_{|X\setminus H_X}$ is strictly non-characteristic for
$\nbigm_{|Y\setminus H_Y}$ (see {\rm\cite[\S3.7]{Sabbah-pure-twistor}}).
If $F$ is an embedding,
we also say that $X$ is non-characteristic for $\nbigm$.
\hfill\qed
\end{df}

Let $\nbigm\in\nbigc(Y;H_Y)$.
There exists
$((\nbigm',\nbigm'',C),W)\in\MTM^{\integral}(Y)$
such that $\nbigm=\nbigm''(\ast H_Y)$.

\begin{lem}
If $F$ is strictly non-characteristic for $\nbigm$,
then $F$ is also strictly non-characteristic for
the $\nbigrtilde_{Y(\ast H_Y)}$-module
$\nbigm'(\ast H_Y)\in\nbigc(Y;H_Y)$.
 \end{lem}
\pf
Because $\Gr^W_w(\nbigm',\nbigm'',C)$ are polarizable for any $w\in\seisuu$,
there exist isomorphisms
$\Gr^W_w(\nbigm')\simeq\Gr^W_{-w}(\nbigm'')$.
Hence,
we have
\[
\Ch(\nbigm'_{|Y\setminus H_Y})=
\bigcup_w\Ch(\Gr^W_w(\nbigm'_{|Y\setminus H_Y}))
=
\bigcup_w\Ch(\Gr^W_w(\nbigm''_{|Y\setminus H_Y}))
=\Ch(\nbigm''_{|Y\setminus H_Y}).
\]
Thus, we obtain the claim of the lemma.
\hfill\qed

\begin{lem}
Suppose that $\nbigm\in\nbigc(Y;H_Y)$
and that $F$ is strictly non-characteristic for $\nbigm$.
Then, we have
\begin{equation}
\label{eq;21.4.5.1}
LF^{\ast}(\nbigm)=
F^{\ast}(\nbigm)
=\nbigr_{X\rarr Y}(\ast H_X)\otimes_{F^{-1}\nbigr_{Y(\ast H_Y)}}F^{-1}\nbigm
=\nbigo_{\nbigx}\otimes_{F^{-1}(\nbigo_{\nbigy})}F^{-1}(\nbigm).
\end{equation}
\end{lem}
\pf
There exists the natural morphism
\begin{equation}
\label{eq;21.4.5.2}
 LF^{\ast}(\nbigm)
 \lrarr
 F^{\ast}(\nbigm)
=
 \nbigr_{X\rarr Y}(\ast H_X)\otimes_{F^{-1}\nbigr_{Y(\ast H_Y)}}
 F^{-1}(\nbigm).
\end{equation}
Because $F_{|X\setminus H_X}$ is strictly non-characteristic for
$\nbigm_{|Y\setminus H_Y}$,
we obtain
$\nbigh^j\bigl(
LF^{\ast}(\nbigm)
\bigr)_{|X\setminus H_X}=0$ $(j\neq 0)$
(see \cite[\S3.7]{Sabbah-pure-twistor}).
It is enough to check the claim locally around any point of $Y$.
We may assume the existence of 
an $\nbigr_Y(\ast H_Y)$-free resolution of $\nbigm$,
and we can check that
$\nbigh^j\bigl(
LF^{\ast}(\nbigm)
\bigr)$
are good $\nbigo_{\nbigx}$-modules
in the sense of \cite[Definition 4.22]{kashiwara_text}.
They are also $\nbigr_X(\ast H_X)$-modules.
Hence, we obtain that
$\nbigh^j\bigl(
LF^{\ast}(\nbigm)
\bigr)=0$ $(j\neq 0)$.
Thus, we obtain (\ref{eq;21.4.5.1}).
\hfill\qed

\begin{lem}
For any $\nbigm\in\nbigc(Y;H_Y)$,
we have 
 $\Ch(\DD_{Y\setminus H_Y}(\nbigm_{|Y\setminus H_Y}))
 =\Ch(\nbigm_{|Y\setminus H_Y})$.
As a result,
$F$ is strictly non-characteristic for $\nbigm\in\nbigc(Y;H_Y)$
if and only if
$F$ is strictly non-characteristic for $\DD_{Y(\ast H_Y)}(\nbigm)$.
\end{lem}
\pf
We obtain 
$\Ch(\DD_{Y\setminus H_Y}(\nbigm_{|Y\setminus H_Y}))
\subset\Ch(\nbigm_{|Y\setminus H_Y})$
by using the argument in \cite[Theorem 2.18]{kashiwara_text}.
Because
$\DD_{Y(\ast H_Y)}(\DD_{Y(\ast H_Y)}\nbigm)=\nbigm$,
we obtain
$\Ch(\DD_{Y\setminus H_Y}(\nbigm_{|Y\setminus H_Y}))
=\Ch(\nbigm_{|Y\setminus H_Y})$.
\hfill\qed

\vspace{.1in}
We shall prove the following theorem
in \S\ref{subsection;21.4.8.1}
after the preliminaries
in \S\ref{subsection;21.4.14.1}--\ref{subsection;21.4.14.2}.

\begin{thm}
\label{thm;21.4.7.50}
If $F$ is strictly non-characteristic for $\nbigm\in\nbigc(Y;H_Y)$,
then $F^{\ast}(\nbigm)$ is an object of $\nbigc(X;H_X)$.
Moreover,
there exists a natural isomorphism
$\lambda^{-\dim X}\DD_{X(\ast H)}F^{\ast}(\nbigm)
\simeq
\lambda^{-\dim Y}F^{\ast}(\DD_{Y(\ast H)}\nbigm)$.
\end{thm}

\subsubsection{Non-characteristic tensor product}

We state a consequence of Theorem \ref{thm;21.4.7.50}.
Let $X$ be a complex manifold
with a closed complex hypersurface $H$.
As in \cite[\S4.4]{kashiwara_text},
we introduce the following definition.

\begin{df}
We say that
$\nbigm_i\in\nbigc(X;H)$ $(i=1,2)$
are non-characteristic
if $\Ch(\nbigm_{1|X\setminus H})
\cap
 \Ch(\nbigm_{2|X\setminus H})$
 is contained in the $0$-section of
 $\cnum\times T^{\ast}(X\setminus H)$.
\hfill\qed
\end{df}

\begin{prop}
If $\nbigm_i\in\nbigc(X;H)$ $(i=1,2)$ are non-characteristic,
then
$\nbigm_1\otimes^L_{\nbigo_{\nbigx}}\nbigm_2
=\nbigm_1\otimes_{\nbigo_{\nbigx}}\nbigm_2$
is an object of $\nbigc(X;H)$.
\end{prop}
\pf
We set $\Xtilde=X\times X$
and $\Htilde=(H\times X)\cup(X\times H)$.
Let $\Delta_X:X\lrarr X\times X$ denote the diagonal embedding.
We have
$\nbigm_1\boxtimes\nbigm_2\in\nbigc(\Xtilde;\Htilde)$.
Because $\Delta_X$ is non-characteristic to
$\nbigm_1\boxtimes\nbigm_2$,
we obtain the claim of the proposition from
Theorem \ref{thm;21.4.7.50}.
\hfill\qed

\begin{cor}
If $\nbigm_i\in\nbigc(X;H)$ $(i=1,2)$ are non-characteristic,
then
$(\DD_{X(\ast H)}\nbigm_1)\otimes_{\nbigo_{\nbigx}}\nbigm_2$
is an object of $\nbigc(X;H)$.
\hfill\qed
\end{cor}
\pf
Because
$\Ch(\nbigm_{1|X\setminus H})=\Ch(\DD_{X\setminus H}(\nbigm_{1|X\setminus H}))$,
the claim follows from the previous proposition.
\hfill\qed

\begin{cor}
Let $\nbigm_i\in\nbigc(X;H)$ $(i=1,2)$.
Suppose that $\nbigm_{1|X\setminus H}$
is a locally free $\nbigo_{\nbigx\setminus\nbigh}$-module,
then  
$\nbigm_1\otimes\nbigm_2\in\nbigc(X;H)$.
\hfill\qed
\end{cor}

\subsubsection{Complexes associated with a family of hypersurfaces}
\label{subsection;21.4.14.1}

Let $Y$ be a complex manifold.
Let $H_Y$ be a hypersurface of $Y$.
Let $\nbigs(Y)$ denote the set of closed complex hypersurfaces of $Y$.
Let $\Gamma$ be any finite set
with a map $\gbigk:\Gamma\lrarr \nbigs(Y)$.
For any non-empty subset $I\subset\Gamma$,
let $\gbigk(I)$ denote the hypersurface of $Y$
obtained as the union of $\gbigk(i)$ $(i\in I)$.
We set $\gbigk(\emptyset):=\emptyset$.
Let $B(\Gamma,\gbigk)$
denote the closure of 
$\bigcap_{i\in\Gamma}(\gbigk(i)\setminus H_Y)$
in $Y$.
We introduce complexes associated with objects of
$\nbigc(Y;H_Y)$
by following \cite{saito2}.

Let $\cnum_{\Gamma}$ denote the complex vector space
generated by $\Gamma$.
For each $i\in\Gamma$,
let $v_i\in\cnum_{\Gamma}$ denote the corresponding vector.
For any non-empty finite subset $I=\{i_1,\ldots,i_k\}\subset \Gamma$,
let $V(I)$ denote the one-dimensional subspace of
$\bigwedge^{\bullet} \cnum_{\Gamma}$
generated by $v_{i_1}\wedge\cdots\wedge v_{i_k}$.
For $I=\emptyset$,
we set $V(\emptyset)=\cnum=\bigwedge^0\cnum_{\Gamma}
\subset\bigwedge^{\bullet}\cnum_{\Gamma}$.
For $i\in\Gamma\setminus I$,
we obtain
$e(i):V(I)\lrarr V(I\cup\{i\})$
induced by the exterior product of $v_i$ from the left.
Let $h_{\Gamma}$ denote the Hermitian metric of $\cnum_{\Gamma}$
for which $v_i$ $(i\in\Gamma)$ are orthonormal.
It induces a Hermitian metric $\bigwedge^{\bullet}\cnum_{\Gamma}$.
We obtain the adjoint
$e(i)^{\dagger}:V(I\cup\{i\})\lrarr V(I)$.

Let $\nbigm\in\nbigc(Y;H_Y)$.
For any $I\subset\Gamma$,
we obtain
$\nbigm[\star \gbigk(I)]\in\nbigc(Y;H_Y)$ $(\star=\ast,!)$.
For $i\in\Gamma\setminus I$,
the natural morphism
$\nbigm[\ast\gbigk(I)]\lrarr\nbigm[\ast \gbigk(I\sqcup\{i\})]$
and $e(i):V(I)\lrarr V(I\sqcup\{i\})$ induce
\begin{equation}
\label{eq;21.4.7.1}
 \nbigm[\ast\gbigk(I)]\otimes V(I)
 \lrarr
 \nbigm[\ast\gbigk(I\sqcup\{i\})]\otimes V(I\sqcup\{i\}),
\end{equation}
and the natural morphism
$\nbigm[!\gbigk(I\sqcup\{i\})]\lrarr\nbigm[!\gbigk(I)]$
and $e(i)^{\dagger}:V(I\sqcup\{i\})\lrarr V(I)$ induce
\begin{equation}
\label{eq;21.4.7.2}
  \nbigm[!\gbigk(I\sqcup\{i\})]\otimes V(I\sqcup\{i\})
  \lrarr
  \nbigm[!\gbigk(I)]\otimes V(I).
\end{equation}

For any $k\in\seisuu_{\geq 0}$,
we set
\[
 C^k_{\ast}(\nbigm,\Gamma,\gbigk):=
 \bigoplus_{\substack{I\subset\Gamma\\ |I|=k}}
 \nbigm[\ast \gbigk(I)]\otimes V(I),
 \quad\quad
 C^{-k}_!(\nbigm,\Gamma,\gbigk):=
 \bigoplus_{\substack{I\subset\Gamma\\ |I|=k}}
 \nbigm[!\gbigk(I)]\otimes V(I).
\]
We formally set $C^k_{\ast}(\nbigm,\Gamma,\gbigk)=0$
and $C^{-k}_!(\nbigm,\Gamma,\gbigk)=0$ for $k<0$.
The morphisms (\ref{eq;21.4.7.1}) induce
$C^k_{\ast}(\nbigm,\Gamma,\gbigk)
\lrarr C^{k+1}_{\ast}(\nbigm,\Gamma,\gbigk)$,
and $C^{\bullet}_{\ast}(\nbigm,\Gamma,\gbigk)$ is a complex
in $\nbigc(Y;H_Y)$.
The morphisms (\ref{eq;21.4.7.2}) induce
$C^{-k-1}_{!}(\nbigm,\Gamma,\gbigk)
\lrarr C^{-k}_{!}(\nbigm,\Gamma,\gbigk)$,
and $C^{\bullet}_{!}(\nbigm,\Gamma,\gbigk)$ is a complex
in $\nbigc(Y;H_Y)$.
For $\star=\ast,!$,
let $\nbigh^k_{\star}(\nbigm,\Gamma,\gbigk)$
denote the $\nbigrtilde_Y(\ast H_Y)$-modules
obtained as the $k$-th cohomology of
the complexes $C^{\bullet}_{\star}(\nbigm,\Gamma,\gbigk)$.

\begin{lem}
\label{lem;21.4.8.2}
 $\nbigh^k_{\star}(\nbigm,\Gamma,\gbigk)$ are objects of
$\nbigc(Y;H_Y)$.
The supports of 
$\nbigh^k_{\star}(\nbigm,\Gamma,\gbigk)$ are
contained in $B(\Gamma,\gbigk)$.
\end{lem}
\pf
Let $((\nbigm',\nbigm'',C),W)\in\MTM^{\integral}(Y;H_Y)$.
In the following, we shall omit to denote the induced weight filtration $W$.
The $\nbigrtilde_{Y(\ast H_Y)}$-triple $(\nbigm',\nbigm'',C)$
is denoted by $\nbigt$.
We obtain 
\[
 \nbigt[\ast\gbigk(I)]
=(\nbigm'[!\gbigk(I)],\nbigm''[\ast\gbigk(I)],C[\ast\gbigk(I)])
\in \MTM^{\integral}(Y;H_Y)
\]
for any $I\subset\Gamma$.
Let $C_I$ denote the Hermitian metric of $V(I)$
induced by $h_{\Gamma}$.
We obtain the following integrable mixed twistor $\nbigd_{Y(\ast H_Y)}$-module:
\[
 \nbigt[\ast \gbigk(I)]\otimes T(I)
 :=\Bigl(
  \nbigm'[!\gbigk(I)]\otimes V(I),\,
  \nbigm''[\ast \gbigk(I)]\otimes V(I),\,
  C[\ast \gbigk(I)]\otimes C_I
 \Bigr).
\]
For $I\sqcup\{i\}\subset \Gamma$,
the morphisms
$\nbigm''[\ast\gbigk(I)]\otimes V(I)
\lrarr \nbigm''[\ast\gbigk(I\sqcup\{i\})]\otimes V(I\sqcup\{i\})$
and
$\nbigm'[!\gbigk(I\sqcup\{i\})]\otimes V(I\sqcup\{i\})
\lrarr\nbigm'[!\gbigk(I)]\otimes V(I)$
induce a morphism
$\nbigt[\ast\gbigk(I)]\otimes T(I)
\lrarr
\nbigt[\ast\gbigk(I\sqcup\{i\})]\otimes T(I\sqcup\{i\})$
in $\MTM^{\integral}(Y;H_Y)$.
For $k\geq 0$, we set
\[
 C^k(\nbigt,\Gamma,\gbigk):=
 \bigoplus_{\substack{I\subset\Gamma\\ |I|=k}}
 \nbigt[\ast\gbigk(I)]\otimes T(I).
\]
Then, we obtain a complex
$C^{\bullet}(\nbigt,\Gamma,\gbigk)$
in the abelian category
$\MTM^{\integral}(Y;H_Y)$.
We obtain
$\nbigh^k(\nbigt,\Gamma,\gbigk)
\in\MTM^{\integral}(Y;H_Y)$
as the $k$-th cohomology of the complex
$C^{\bullet}(\nbigt,\Gamma,\gbigk)$,
which consists of 
$\nbigh^{-k}_{!}(\nbigm',\Gamma,\gbigk)$
and
$\nbigh^{k}_{\ast}(\nbigm'',\Gamma,\gbigk)$
with the induced sesqui-linear pairing
and the weight filtration.
Then, the first claim of the lemma immediately follows.
We can prove the second claim
by an argument similar to \cite[Lemma 14.1.18]{Mochizuki-MTM}.
\hfill\qed

\vspace{.1in}

Let $\Gamma'\subset\Gamma$  be a subset.
Let $\gbigk':\Gamma'\lrarr\nbigs(Y)$ denote the map
obtained as the restriction of $\gbigk$.
There exist the natural projections
$C^k_{\ast}(\nbigm,\Gamma,\gbigk)\lrarr
C^k_{\ast}(\nbigm,\Gamma',\gbigk')$
which induce a morphism of complexes
\begin{equation}
\label{eq;21.4.7.3}
C^{\bullet}_{\ast}(\nbigm,\Gamma,\gbigk)
\lrarr
C^{\bullet}_{\ast}(\nbigm,\Gamma',\gbigk').
\end{equation}
There exist the natural inclusions
$C^k_{!}(\nbigm,\Gamma',\gbigk')\lrarr
C^k_{!}(\nbigm,\Gamma,\gbigk)$,
which induce a morphism of complexes
\begin{equation}
\label{eq;21.4.7.4}
C^{\bullet}_{!}(\nbigm,\Gamma',\gbigk')
\lrarr
C^{\bullet}_!(\nbigm,\Gamma,\gbigk).
\end{equation}
\begin{lem}
If $B(\Gamma,\gbigk)=B(\Gamma',\gbigk')$,
 then the induced morphisms
 \[
\nbigh^k_{\ast}(\nbigm,\Gamma,\gbigk)
 \lrarr\nbigh^k_{\ast}(\nbigm,\Gamma',\gbigk'),
 \quad\quad
\nbigh^k_{!}(\nbigm,\Gamma',\gbigk')
\lrarr\nbigh^k_{!}(\nbigm,\Gamma,\gbigk)
 \]
are isomorphisms.
\end{lem}
\pf
Let $(\nbigt,W)\in\MTM^{\integral}(Y;H_Y)$.
We use the notation in the proof of Lemma \ref{lem;21.4.8.2}.
There exists a naturally induced morphism of
the complexes of integrable mixed twistor $\nbigd_{Y(\ast H_Y)}$-modules
$a:C^{\bullet}(\nbigt,\Gamma,\gbigk)
\lrarr
C^{\bullet}(\nbigt,\Gamma',\gbigk')$.
For any $i\in\Gamma\setminus\Gamma'$,
we have
$B(\Gamma',\gbigk')=B(\Gamma,\gbigk)\subset\gbigk(i)$.
Because the localization functor
$\nbigt\longmapsto \nbigt[\ast \gbigk(i)]$ is exact,
we have
\[
 \nbigh^k(\nbigt[\ast \gbigk(i)],\Gamma',\gbigk')
 =\nbigh^k(\nbigt,\Gamma',\gbigk')[\ast\gbigk(i)]=0.
\]
Then, it is easy to check that
$\Ker(a)$ is acyclic.
\hfill\qed

\vspace{.1in}

Let $(\Gamma_i,\gbigk_i)$ $(i=1,2)$ be tuples of
a finite set $\Gamma_i$
and a map $\gbigk_i:\Gamma_i\lrarr \nbigs(Y)$.
If there exists an injection
$\varphi:\Gamma_1\lrarr \Gamma_2$
such that
$\gbigk_2\circ\varphi=\gbigk_1$,
we obtain morphisms of complexes
\[
 \varphi^{\ast}:
 C^{\bullet}_{\ast}(\nbigm,\Gamma_2,\gbigk_2)
 \lrarr
 C^{\bullet}_{\ast}(\nbigm,\Gamma_1,\gbigk_1),
\]
\[
 \varphi_!:
 C^{\bullet}_{!}(\nbigm,\Gamma_1,\gbigk_1)
 \lrarr
 C^{\bullet}_{!}(\nbigm,\Gamma_2,\gbigk_2)
\]
as in (\ref{eq;21.4.7.3}) and (\ref{eq;21.4.7.4}),
by identifying $\Gamma_1$ as $\varphi(\Gamma_1)\subset\Gamma_2$.
We obtain the induced morphisms
$\varphi^{\ast}:
\nbigh^k_{\ast}(\nbigm,\Gamma_2,\gbigk_2)
\lrarr
\nbigh^k_{\ast}(\nbigm,\Gamma_1,\gbigk_1)$
and 
$\varphi_!:
\nbigh^k_{!}(\nbigm,\Gamma_1,\gbigk_1)
\lrarr
\nbigh^k_{!}(\nbigm,\Gamma_2,\gbigk_2)$.

\vspace{.1in}
Let $(\Gamma_i,\gbigk_i)$ $(i=1,2)$ be tuples of
a finite set $\Gamma_i$
and a map $\gbigk_i:\Gamma_i\lrarr \nbigs(Y)$.
Suppose that
$B(\Gamma_2,\gbigk)\subset B(\Gamma_1,\gbigk)$.
We set
$\Gamma_3:=\Gamma_1\sqcup\Gamma_2$.
We obtain the map $\gbigk_3:\Gamma_3\lrarr\nbigs(Y)$
induced by $\gbigk_1$ and $\gbigk_2$.
There exist the natural inclusions
$\iota_i:\Gamma_i\lrarr\Gamma_3$ $(i=1,2)$.
We obtain the following morphisms:
\begin{equation}
\label{eq;21.4.7.5}
\begin{CD}
 \nbigh_{\ast}^k(\nbigm,\Gamma_1,\gbigk_1)
 @<{\iota_1^{\ast}}<<
 \nbigh_{\ast}^k(\nbigm,\Gamma_3,\gbigk_3)
 @>{\iota_2^{\ast}}>{\simeq}>
 \nbigh_{\ast}^k(\nbigm,\Gamma_2,\gbigk_2).
\end{CD}
\end{equation}
We also obtain the following morphisms:
\begin{equation}
\label{eq;21.4.7.6}
\begin{CD}
 \nbigh_!^k(\nbigm,\Gamma_1,\gbigk_1)
 @>{\iota_{1!}}>>
 \nbigh_!^k(\nbigm,\Gamma_3,\gbigk_3)
 @<{\iota_{2!}}<{\simeq}<
 \nbigh_!^k(\nbigm,\Gamma_2,\gbigk_2).
\end{CD}
\end{equation}
We obtain the morphisms
\[
g_{(\Gamma_1,\gbigk_1),(\Gamma_2,\gbigk_2),\ast}:
\nbigh_{\ast}^k(\nbigm,\Gamma_2,\gbigk_2)
\lrarr
\nbigh_{\ast}^k(\nbigm,\Gamma_1,\gbigk_1)
\]
\[
g_{(\Gamma_2,\gbigk_2),(\Gamma_1,\gbigk_1),!}:
\nbigh_{!}^k(\nbigm,\Gamma_1,\gbigk_1)
\lrarr
\nbigh_{!}^k(\nbigm,\Gamma_2,\gbigk_2)
\]
from (\ref{eq;21.4.7.5})
and (\ref{eq;21.4.7.6}), respectively.

\begin{prop}
\label{prop;21.4.9.11}
Let $(\Gamma_4,\gbigk_4)$ be a tuple of
a finite set $\Gamma_4$
and a map $\gbigk_4:\Gamma_4\lrarr\nbigs(Y)$
such that $B(\Gamma_4,\gbigk_4)=B(\Gamma_2,\gbigk_2)$.
Suppose that there exist injections
$\varphi_i:\Gamma_i\lrarr\Gamma_4$ $(i=1,2)$
such that $\gbigk_i=\gbigk_4\circ\varphi_i$.
Then, we have
$\varphi_{1}^{\ast}\circ(\varphi_2^{\ast})^{-1}
 =g_{(\Gamma_1,\gbigk_1),(\Gamma_2,\gbigk_2),\ast}$
and
$(\varphi_{2!})^{-1}\circ \varphi_{1!}
=g_{(\Gamma_1,\gbigk_1),(\Gamma_2,\gbigk_2),!}$.
\end{prop}
\pf
We set $\Gamma_2':=\Gamma_2$
and $\gbigk_2':=\gbigk_2$.
We set
$\Gammatilde_2:=\Gamma_2\sqcup\Gamma_2'$
which is equipped with a map
$\gbigktilde_2:\Gammatilde_2\lrarr\nbigs(Y)$
induced by $\gbigk_2$ and $\gbigk_2'$.
For $j\in\Gamma_2$, let $j'\in\Gamma_2$
denote the corresponding element.
For any finite subset $I\subset\Gamma_2$,
let $I'\subset\Gamma_2'$ denote the corresponding subset.
Let $\iota_2:\Gamma_2\lrarr \Gammatilde_2$
and $\iota_2':\Gamma_2'\lrarr\Gammatilde_2$
denote the inclusions.
We obtain the isomorphisms
\[
 (\iota_2')^{\ast}\circ
 (\iota_2^{\ast})^{-1}:
 \nbigh^k(\nbigm,\Gamma_2,\gbigk_2)
 \simeq
 \nbigh^k(\nbigm,\Gamma'_2,\gbigk'_2).
\]
\[
 (\iota_{2!})^{-1}\circ\iota_{2!}':
  \nbigh^k(\nbigm,\Gamma'_2,\gbigk'_2)\simeq
  \nbigh^k(\nbigm,\Gamma_2,\gbigk_2).
\]
The natural bijection $\id:\Gamma'_2\simeq\Gamma_2$
induces the isomorphisms
\[
 \id^{\ast}:
  \nbigh^k(\nbigm,\Gamma_2,\gbigk_2)
 \simeq
 \nbigh^k(\nbigm,\Gamma'_2,\gbigk'_2),
\]
\[
 \id_{!}:
  \nbigh^k(\nbigm,\Gamma'_2,\gbigk'_2)\simeq
  \nbigh^k(\nbigm,\Gamma_2,\gbigk_2).
\]

\begin{lem}
\label{lem;21.4.7.20}
 We have
$\id^{\ast}= (\iota_2')^{\ast}\circ
(\iota_2^{\ast})^{-1}$
and 
$\id_{!}= (\iota_{2!})^{-1}\circ\iota_{2!}'$.
\end{lem}
\pf
Let $I\subset \Gamma_2$.
For a decomposition $I=I_1\sqcup I_2$,
we obtain
$I_1\sqcup I_2'\subset\Gammatilde_2$.
The natural bijections
$I_1\simeq I_1$ and $I_2\simeq I_2'$
induce $\cnum_I\lrarr \cnum_{\Gammatilde}$,
which induces
$V(I)\simeq V(I_1\sqcup I_2')$.
We have the natural isomorphisms
\[
 a_{\star}(I,I_1,I_2):
 \nbigm(\star \gbigk(I))\otimes V(I)
 \simeq
 \nbigm(\star \gbigktilde(I_1\sqcup I_2'))
 \otimes V(I_1\sqcup I_2').
\]
We set
$b_{\star}(I,I_1,I_2):=
\frac{|I|!}{|I_1|!|I_2|!}
 a_{\star}(I,I_1,I_2)$.
For any $I\subset \Gamma_2$,
we obtain
\begin{equation}
\label{eq;21.4.7.10}
 \nbigm[\ast \gbigk_2(I)]\otimes V(I)
 \lrarr
 \bigoplus_{I_1\sqcup I_2=I}
 \nbigm[\ast\gbigktilde_2(I_1\sqcup I'_2)]\otimes V(I_1\sqcup I_2')
\end{equation}
by $b_{\ast}(I,I_1,I_2)$.
We obtain
\begin{equation}
 \label{eq;21.4.7.11}
 \bigoplus_{I_1\sqcup I_2=I}
 \nbigm[!\gbigktilde_2(I_1\sqcup I'_2)]\otimes V(I_1\sqcup I_2')
 \lrarr
 \nbigm[! \gbigk_2(I)]\otimes V(I)
\end{equation}
by $b_{!}(I,I_1,I_2)$.
We obtain
\begin{equation}
 \label{eq;21.4.7.12}
 C^k_{\ast}(\nbigm,\Gamma_2,\gbigk_2)
 \lrarr
 C^k_{\ast}(\nbigm,\Gammatilde_2,\gbigktilde_2)
\end{equation}
from (\ref{eq;21.4.7.10})
for any $I\subset \Gamma_2$ with $|I|=k$.
Similarly, we obtain
\begin{equation}
 \label{eq;21.4.7.13}
 C^k_{!}(\nbigm,\Gammatilde_2,\gbigktilde_2)
 \lrarr
 C^k_{!}(\nbigm,\Gamma_2,\gbigk_2)
\end{equation}
from (\ref{eq;21.4.7.11})
for any $I\subset \Gamma_2$ with $|I|=k$.
We can check that
(\ref{eq;21.4.7.12}) and (\ref{eq;21.4.7.13})
induce morphisms of complexes
$C_{\ast}^{\bullet}(\nbigm,\Gamma_2,\gbigk_2)
 \lrarr
 C_{\ast}^{\bullet}(\nbigm,\Gammatilde_2,\gbigktilde_2)$
and
$C_{!}^{\bullet}(\nbigm,\Gammatilde_2,\gbigktilde_2)\lrarr
 C_{!}^{\bullet}(\nbigm,\Gamma_2,\gbigk_2)$.
The composition of the morphisms
\[
 C_{\ast}^{\bullet}(\nbigm,\Gamma_2,\gbigk_2)
 \lrarr
 C_{\ast}^{\bullet}(\nbigm,\Gammatilde_2,\gbigktilde_2)
 \lrarr
 C_{\ast}^{\bullet}(\nbigm,\Gamma_2,\gbigk_2)
\]
is the identity.
The composition of the morphisms
\[
 C_{\ast}^{\bullet}(\nbigm,\Gamma_2,\gbigk_2)
 \lrarr
 C_{\ast}^{\bullet}(\nbigm,\Gammatilde_2,\gbigktilde_2)
 \lrarr
 C_{\ast}^{\bullet}(\nbigm,\Gamma'_2,\gbigk'_2)
\]
is the morphism induced by $\id:\Gamma_2'\simeq\Gamma_2$.
Hence, we obtain the 
$\id^{\ast}= (\iota_2')^{\ast}\circ
(\iota_2^{\ast})^{-1}$.
We obtain the other identity similarly.
\hfill\qed

\vspace{.1in}

We set
$\Gamma_5:=\Gamma_4\sqcup\Gamma'_2$,
which is equipped with
$\gbigk_5:\Gamma_5\lrarr\nbigs(Y)$
induced by $\gbigk_4$ and $\gbigk_2'$.
There exist the natural inclusions
$j_4:\Gamma_4\lrarr\Gamma_5$
and $j'_2:\Gamma'_2\lrarr \Gamma_5$.
We can replace
$(\Gamma_4,\gbigk_4,\varphi_i)$
with
$(\Gamma_5,\gbigk_5,j_4\circ\varphi_i)$.
We obtain the injection
$\psi:\Gammatilde_2\lrarr \Gamma_5$
from $j_4\circ\varphi_2$ and $j'_2$.
By using Lemma \ref{lem;21.4.7.20},
we can replace
$(\Gamma_2,j_4\circ\varphi_2)$
with $(\Gamma'_2,j_2')$.
We identify $\Gamma_3$ with $\Gamma_1\sqcup \Gamma_2'$
by identifying $\Gamma_2$ and $\Gamma_2'$.
We obtain the map
$\kappa:\Gamma_3\lrarr \Gamma_5$
by $j_4\circ\varphi_1$ and $j_2'$,
from which we obtain the claim of the proposition.
\hfill\qed

\begin{cor}
$g_{(\Gamma,\gbigk),(\Gamma,\gbigk),\star}$ $(\star=!,\ast)$
are the identity maps.
If
$B(\Gamma_{10},\gbigk_{10})
=B(\Gamma_{11},\gbigk_{11})
=B(\Gamma_{12},\gbigk_{12})$,
then we have
\[
 g_{(\Gamma_{10},\gbigk_{10}),(\Gamma_{11},\gbigk_{11}),\ast}
 \circ
 g_{(\Gamma_{11},\gbigk_{11}),(\Gamma_{12},\gbigk_{12}),\ast}
 =g_{(\Gamma_{10},\gbigk_{10}),(\Gamma_{12},\gbigk_{12}),\ast}
\]
\[
  g_{(\Gamma_{12},\gbigk_{12}),(\Gamma_{11},\gbigk_{11}),!}
 \circ
 g_{(\Gamma_{11},\gbigk_{11}),(\Gamma_{10},\gbigk_{10}),!}
 =g_{(\Gamma_{12},\gbigk_{12}),(\Gamma_{10},\gbigk_{10}),!}
\]
In this sense,
$\nbigh^k_{\star}(\nbigm,\Gamma,\gbigk)$ $(\star=!,\ast)$
depend only on $B(\Gamma,\gbigk)$
up to canonical isomorphisms. 
\hfill\qed
\end{cor}

Let $(\nbigt,W)\in\MTM^{\integral}(Y;H_Y)$.
We construct a complex
$C^{\bullet}(\nbigt,\Gamma,\gbigk)$
in $\MTM^{\integral}(Y;H_Y)$
as in the proof of Lemma \ref{lem;21.4.8.2}.
We obtain
$\nbigh^k(\nbigt,\Gamma,\gbigk)
\in \MTM^{\integral}(Y;H_Y)$.
For any $(\Gamma_i,\gbigk_i)$ $(i=1,2)$
such that $B(\Gamma_1,\gbigk_1)=B(\Gamma_2,\gbigk_2)$,
there exists the following canonical isomorphism
induced by the inclusions
$\Gamma_i\subset\Gamma_1\sqcup\Gamma_2$
as in (\ref{eq;21.4.7.5}) and (\ref{eq;21.4.7.6}):
\[
 g_{(\Gamma_2,\gbigk_2),(\Gamma_1,\gbigk_1)}:
 \nbigh^k(\nbigt,\Gamma_1,\gbigk_1)
 \simeq
 \nbigh^k(\nbigt,\Gamma_2,\gbigk_2).
\]
\begin{cor}
$g_{(\Gamma,\gbigk),(\Gamma,\gbigk)}$ is the identity map.
If $B(\Gamma_1,\gbigk_1)=B(\Gamma_2,\gbigk_2)=B(\Gamma_3,\gbigk_3)$,
we have
$g_{(\Gamma_3,\gbigk_3),(\Gamma_2,\gbigk_2)}\circ
g_{(\Gamma_2,\gbigk_2),(\Gamma_1,\gbigk_1)}
=g_{(\Gamma_3,\gbigk_3),(\Gamma_1,\gbigk_1)}$.
In this sense,
$\nbigh^k(\nbigt,\Gamma,\gbigk)$
depend only on $B(\Gamma,\gbigk)$
up to canonical isomorphisms.
\hfill\qed
\end{cor}

We identify $\cnum_{\Gamma}$ with its dual space
by a natural bilinear form for which
the base $v_i$ $(i\in\Gamma)$ is orthonormal.
By the construction,
we obtain the following lemma.
\begin{lem}
There exists the natural isomorphism
$\DD_{Y(\ast H_Y)}C_!^{\bullet}\bigl(\nbigm,\Gamma,\gbigk\bigr)
 \simeq
 C_{\ast}^{\bullet}\bigl(\DD_{Y(\ast H_Y)}\nbigm,\Gamma,\gbigk\bigr)$,
 which induces
\[
 \DD_{Y(\ast H_Y)}\nbigh^{-\ell}_{!}(\nbigm,\Gamma,\gbigk)
 \simeq
 \nbigh^{\ell}\bigl(
 \DD_{Y(\ast H_Y)}(\nbigm),\Gamma,\gbigk
 \bigr).
\]
Under the isomorphism,
we obtain
$\DD\bigl(
 g_{(\Gamma_2,\gbigk_2),(\Gamma_1,\gbigk_1),!}
 \bigr)
 =g_{(\Gamma_1,\gbigk_1),(\Gamma_2,\gbigk_2),\ast}$,
 where
 $g_{(\Gamma_2,\gbigk_2),(\Gamma_1,\gbigk_1),!}$ in the left hand side
 is defined for $\nbigm$,
 and 
 $g_{(\Gamma_1,\gbigk_1),(\Gamma_2,\gbigk_2),\ast}$
 in the right hand side is defined for
 $\DD_{Y(\ast H_Y)}(\nbigm)$.
\hfill\qed
\end{lem}

\subsubsection{Coordinate hypersurfaces}

Let $Y$ be a complex manifold
with a holomorphic coordinate system $\vecy=(y_1,\ldots,y_n)$.
Let $H_Y$ be a closed complex hypersurface of $Y$.
For any $c\in\cnum$, we set $X_c:=\{y_1=c\}$,
and let $i_c:X_c\lrarr Y$ denote the inclusion.
We set $X:=X_0$ and $i_X:=i_0$.
We also set $H_X:=H_Y\cap X$.

We set $\Gamma:=\{1\}$
and $\gbigk_{\vecy}(1):=\{y_1=0\}$.
We have $B(\Gamma,\gbigk_{\vecy})=X$.
Let $\nbigm\in\nbigc(Y;H_Y)$
such that $i_c$ $(c\in\cnum)$
are strictly non-characteristic for $\nbigm$.

\begin{lem}
The $V$-filtration of
the $\nbigrtilde_{Y(\ast H_Y)}$-module $\nbigm$ along $y_1$
is given by $V_a(\nbigm)=y_1^{[-1-a]}\nbigm$ for $a\leq -1$,
 and $V_a(\nbigm)=\nbigm$ for $a>-1$.
\end{lem}
\pf
As in \cite[\S3.7]{Sabbah-pure-twistor},
$\nbigm_{|Y\setminus H_Y}$ is coherent over
$\nbigo_{\nbigy\setminus\nbigh_Y}\langle
 \deldel_{y_2},\ldots,\deldel_{y_n}
 \rangle$.
There exists $\nbigm_0\in\nbigc(Y)$
such that $\nbigm_0(\ast H_Y)=\nbigm$.
There exists the $V$-filtration $V(\nbigm_0)$ along $y_1$.
We set $V_a(\nbigm):=V_a(\nbigm_0)(\ast H_Y)$,
which is a $V$-filtration of
the $\nbigrtilde_Y(\ast H_Y)$-module $\nbigm$ along $y_1$.
We set
$V'_a(\nbigm)=y_1^{[-a-1]}\nbigm$ for $a\leq -1$,
and $V'_a(\nbigm)=\nbigm$ for $a>-1$.
As in \cite[\S3.7]{Sabbah-pure-twistor},
$V'_{\bullet}(\nbigm)_{|Y\setminus H_Y}$
is the $V$-filtration of $\nbigm_{|Y\setminus H_Y}$ along $y_1$,
and hence
we have
$V'_{\bullet}(\nbigm)_{|Y\setminus H_Y}
=V_a(\nbigm)_{|Y\setminus H_Y}$.
Because both $V'_a(\nbigm)$ and $V_a(\nbigm)$
are good $\nbigo_{\nbigy}$-submodules
of a good $\nbigo_{\nbigy}$-module $\nbigm$
satisfying
$V'_a(\nbigm)(\ast H_Y)=V'_a(\nbigm)$
and $V_a(\nbigm)(\ast H_Y)=V_a(\nbigm)$,
we obtain $V'_a(\nbigm)=V_a(\nbigm)$ for any $a\in\real$.
\hfill\qed

\begin{lem}
We have
$\nbigh^k_{\ast}(\nbigm,\Gamma,\gbigk_{\vecy})=0$
and
$\nbigh^{-k}_{!}(\nbigm,\Gamma,\gbigk_{\vecy})=0$
unless $k=1$. 
\end{lem}
\pf
Because $\Gr^V_0(\nbigm)=0$,
the natural morphism $\nbigm\lrarr\nbigm[\ast y_1]$
is a monomorphism,
and the natural morphism
$\nbigm[!y_1]\lrarr\nbigm$ is an epimorphism.
It implies the claim of the lemma.
\hfill\qed

\vspace{.1in}

Let $\lefttop{1}V\nbigr_Y\subset\nbigr_Y$
denote the sheaf of subalgebras
generated by
$\lambda p_Y^{\ast}\Theta_Y(\log y_1)$.
We obtain
\[
 \nbigm[\ast y_1]=
 \nbigr_Y(\ast H_Y)\otimes_{\lefttop{1}V\nbigr_Y(\ast H_Y)}
 (y_1^{-1}\nbigm)
 \simeq
 \nbigo_{\nbigy}\langle\deldel_{y_1}\rangle
 \otimes_{\nbigo_{\nbigy}\langle y_1\deldel_{y_1}\rangle}
 (y_1^{-1}\nbigm).
\]
\[
 \nbigm[!y_1]=
 \nbigr_Y(\ast H_Y)\otimes_{\lefttop{1}V\nbigr_Y(\ast H_Y)}
 \nbigm
 \simeq
  \nbigo_{\nbigy}\langle\deldel_{y_1}\rangle
 \otimes_{\nbigo_{\nbigy}\langle y_1\deldel_{y_1}\rangle}
 \nbigm.
\]

For a section $m$ of $\nbigm$,
we obtain a section $y_1^{-1}m$ of $\nbigm[\ast y_1]$,
which induces a section
$[y_1^{-1}m]$ of $\nbigm[\ast y_1]/\nbigm$.
This correspondence induces the following morphism
\[
\lambda^{-1}\cdot
 i_{X\ast}\bigl(
 i_X^{\ast}(\nbigm)
 \otimes (dy_1/\lambda)^{-1}
 \bigr)
 \lrarr
 \nbigm[\ast y_1]/\nbigm
\]
by $i_{X\ast}\bigl(
i_X^{\ast}(m)(dy_1)^{-1}\bigr)
\longmapsto
[y_1^{-1}m]$.
It induces an isomorphism
\[
\rho_{\vecy,\ast}:
\lambda^{-1} i_{X\dagger}i_X^{\ast}\nbigm
\simeq
\nbigm[\ast y_1]/\nbigm
\simeq
\nbigh^1_{\ast}(\nbigm,\Gamma,\gbigk_{\vecy})
\]

For a section $m$ of $\nbigm$,
we obtain a section $\deldel_{y_1}\otimes m$ of
$\nbigm[!y_1]$
and a section $\deldel_{y_1}m$ of $\nbigm$,
which induces a section $1\otimes\deldel_{y_1}m$
of $\nbigm[!y_1]$.
We obtain a section
$\deldel_{y_1}\otimes m-1\otimes(\deldel_{y_1}m)$
of the kernel of
$\nbigm[!y_1]\lrarr\nbigm$.
This procedure induces the following morphism
\[
 i_{X\ast}\bigl(
 i_X^{\ast}(\nbigm)\otimes (dy_1/\lambda)^{-1}
 \bigr)
 \lrarr
 \Ker(\nbigm[!y_1]\lrarr\nbigm)
\]
by
$i_{X\ast}(i_X^{\ast}(m)\otimes (dy_1/\lambda)^{-1})
\longmapsto
\deldel_{y_1}\otimes m-
1\otimes(\deldel_{y_1}m)$.
It induces the isomorphism
\[
 \rho_{\vecy,!}:
 i_{X\dagger}\bigl(
  i_X^{\ast}(\nbigm)
  \bigr)
  \simeq
  \Ker(\nbigm[!y_1]\lrarr\nbigm)
  \simeq \nbigh^{-1}_!(\nbigm,\Gamma,\gbigk_{\vecy}).
\]

By using the isomorphisms
$\rho_{\vecy,\ast}$ or $\rho_{\vecy,!}$,
we obtain the following lemma.
\begin{lem}
\label{lem;21.4.9.3}
 $i_X^{\ast}\nbigm$ is an object of
 $\nbigc(X;H_X)$.
\hfill\qed
\end{lem}

Let $(z_1,\ldots,z_n)$ be another holomorphic coordinate system
of $Y$ such that $X=\{z_1=0\}$.
\begin{lem}
\label{lem;21.4.7.40}
We have
$\rho_{\vecy,\star}=\rho_{\vecz,\star}$ $(\star=!,\ast)$.
Namely,
the isomorphisms
$\rho_{\vecy,\star}$ $(\star=\ast,!)$
are independent of the choice of $(y_1,\ldots,y_n)$. 
 \end{lem}
\pf
We have
$dz_{1|X}=\bigl(
 \del_{y_1}(z_1)\cdot dy_1
 \bigr)_{|X}$.
We have
$i_X^{\ast}(m)(dy_1)^{-1}_{|X}
 =i_X^{\ast}(m\cdot \del_{y_1}z_1)\cdot
  (dz_1)^{-1}_{|X}$.
Because
$y_1^{-1}m-z_1^{-1}(\del_yz_1)m$ is a section of $\nbigm$,
we obtain $\rho_{\vecy,\ast}=\rho_{\vecz,\ast}$.
We have
$\del_{y_1}=\sum_{j=1}^n (\del_{y_1}z_j)\del_{z_j}$,
and
\begin{multline}
 \del_{y_1}\otimes m-1\otimes \del_{y_1}m
 =(\del_{y_1}z_1)\del_{z_1}\otimes m
 -1\otimes (\del_{y_1}z_1)\del_{z_1}m
  \\
 =\del_{z_1}\otimes (\del_{y_1}z_1)m
 -1\otimes \del_{z_1}(\del_{y_1}z_1)\cdot m
 -1\otimes (\del_{y_1}z_1)\del_{z_1}m
 =\del_{z_1}\otimes\bigl(
 (\del_{y_1}z_1)m\bigr)
 -1\otimes \del_{z_1}\bigl((\del_{y_1}z_1)m\bigr).
\end{multline}
Hence, we obtain $\rho_{\vecy,!}=\rho_{\vecz,!}$.
\hfill\qed

\begin{cor}
\label{cor;21.4.12.22}
We have
$\rho_{\vecy,\ast}=
g_{(\Gamma,\gbigk_{\vecy}),(\Gamma,\gbigk_{\vecz}),\ast}\circ
\rho_{\vecz,\ast}$
and
 $\rho_{\vecy,!}=
 \rho_{\vecz,!}\circ g_{(\Gamma,\gbigk_{\vecz}),(\Gamma,\gbigk_{\vecy}),!}$.
\hfill\qed
\end{cor}

Let $(i_X)_{\pi}:
 \cnum_{\lambda}\times(X\times_YT^{\ast}Y)_{|X\setminus H_Y}
 \lrarr
\cnum_{\lambda}\times(T^{\ast}Y)_{|Y\setminus H_Y}$
denote the natural inclusion.
We obtain the natural morphism
$(i_X)_d:
\cnum_{\lambda}\times(X\times_YT^{\ast}Y)_{|X\setminus H_Y}
\lrarr
\cnum_{\lambda}\times(TX^{\ast})_{|X\setminus H_X}$
from the tangent map of $i_X$.
\begin{lem}
$\Ch(i_X^{\ast}(\nbigm)_{|X\setminus H_X})
 \subset
 (i_X)_d\bigl(
 (i_X)_{\pi}^{-1}(\Ch(\nbigm_{|Y\setminus H_Y}))
 \bigr)$.
\end{lem}
\pf
It is enough to prove the claim locally around
any point $P$ of $\cnum_{\lambda}\times(X\setminus H_X)$.
We may assume $H_X=\emptyset$ from the beginning.
Let $F\nbigr_Y$ denote the filtration defined by the orders of
differential operators.
Let $\eta_i$ denote the sections of $\Gr^F\nbigr_Y$
induced by $\deldel_{y_i}$.
Let $\pi:\cnum_{\lambda}\times T^{\ast}Y\lrarr \cnum\times Y$
denote the projection.
Let $\nbigu$ denote
a neighbourhood of $P$ in $\cnum\times Y$.
We obtain the induced coordinate system
$(\lambda,y_1,\ldots,y_n,\eta_1,\ldots,\eta_n)$
on $\pi^{-1}(\nbigu)$.

Let $F_{\bullet}(\nbigm_{|\nbigu})$ be
a coherent filtration of $\nbigm_{|\nbigu}$.
We obtain the associated coherent
$\nbigo_{\nbigu}\langle\eta_1,\ldots,\eta_{n}\rangle$-module
$\Gr^F(\nbigm_{|\nbigu})$ on $\nbigu$,
which induces a coherent
$\nbigo_{\pi^{-1}(\nbigu)}$-module
$\Gr^F(\nbigm_{|\nbigu})^{\sim}$.
Let $F_j(i_X^{\ast}\nbigm)$ denote the image of
$i_X^{\ast}(F_j\nbigm)\lrarr i_X^{\ast}(\nbigm)$.
Because the natural morphism
$i_X^{\ast}F_j(\nbigm)\big/
i_X^{\ast}F_j(\nbigm)
\lrarr
i_X^{\ast}\Gr^F_j(\nbigm)$
is an isomorphism,
we obtain the epimorphism
$i_X^{\ast}\Gr^F(\nbigm)
 \lrarr
 \Gr^Fi_X^{\ast}(\nbigm)$.
Therefore, there exists
the natural epimorphism
\[
 (i_X)_{d\ast}
 (i_X)_{\pi}^{\ast}
 \Gr^F(\nbigm)^{\sim}
 \lrarr
 \Gr^F(i_X^{\ast}\nbigm)^{\sim}.
\]
Thus, we obtain the claim of the lemma.
\hfill\qed

\vspace{.1in}

We also have
$\Ch(\nbigm[\star y_1]_{|Y\setminus H_Y})
=\Ch(\nbigm)_{|Y\setminus H_Y}
 \cup
 (i_X)_{\pi}\bigl(
 (i_{X})_d^{-1}
 (\Ch(i_X^{\ast}\nbigm)_{|X\setminus H_X})
 \bigr)$.

\subsubsection{Tuples of coordinate functions}
\label{subsection;21.4.14.2}

Let $Y$ be a complex manifold
with a holomorphic coordinate system $(y_1,\ldots,y_n)$.
Let $H_Y$ be a closed complex hypersurface of $Y$.
For any $\vecc\in\cnum^{\ell}$,
we set $X_{\vecc}:=\bigcap_{i=1}^{\ell}\{y_i=c_i\}$,
and let $i_{\vecc}:X_{\vecc}\lrarr Y$ denote the inclusion.
We set $X:=X_{0,\ldots,0}$
and $i_X:=i_{0,\ldots,0}$.
We also set $H_X:=H_X\cap Y$.
We set $\Gamma:=\{1,\ldots,\ell\}$,
and let $\gbigk_{\vecy}:\Gamma\lrarr\nbigs(Y)$
be defined by $\gbigk_{\vecy}(i):=\{y_i=0\}$.
We have $B(\Gamma,\gbigk_{\vecy})=X$.
For any $1\leq j\leq \ell$,
we set $\Gamma_{\geq j}:=\{j,\ldots,\ell\}$
and $\Gamma_{\leq j}:=\{1,\ldots,j\}$.
Formally, we set $\Gamma_{\leq 0}:=\emptyset$.
Let $\gbigk_{\vecy,\geq j}:\Gamma_{\geq j}\lrarr\nbigs(Y)$
and $\gbigk_{\vecy,\leq j}:\Gamma_{\leq j}\lrarr\nbigs(Y)$
denote the induced maps.
We set $Y_{\leq j}:=\gbigk_{\vecy}(\Gamma_{\leq j})$
and $H_{Y_{\leq j}}:=Y_{\leq j}\cap H_Y$.
We have $Y_{\leq\ell}=X$.
Let $\iota_j:Y_{\leq j}\lrarr Y$
and $\iota_{j-1,j}:Y_{\leq j}\lrarr Y_{\leq j-1}$
denote the inclusions.
Let $\nbigm\in\nbigc(Y;H_Y)$.
Assume the following condition.
\begin{condition}
\label{condition;21.4.12.10}
 $X_{\vecc}$ $(\vecc\in\cnum^{\ell})$
are strictly non-characteristic for $\nbigm$. 
\hfill\qed
\end{condition}

Note that Condition \ref{condition;21.4.12.10} implies that
$\iota_j$ is strictly non-characteristic for $\nbigm$,
and $\iota_j^{\ast}\nbigm$ satisfies Condition \ref{condition;21.4.12.10}
on $Y_{\leq j}$
for $\bigcap_{i=j+1}^{\ell}\{y_j=c_j\}$
$(\vecc\in\cnum^{\Gamma_{\geq j+1}})$
with respect to the coordinate system
$(y_{j+1},\ldots,y_n)$.

\begin{lem}
$\iota_j^{\ast}(\nbigm)$ are objects
of $\nbigc(Y_{\leq j};H_{Y_{\leq j}})$.
\end{lem}
\pf
It follows from Lemma \ref{lem;21.4.9.3}.
\hfill\qed

\vspace{.1in}

Let us construct a quasi-isomorphism
\begin{equation}
\label{eq;21.4.9.1}
  C_{\ast}^{\bullet}(\nbigm,\Gamma,\gbigk_{\vecy})
  \lrarr
  \lambda^{-1}\iota_{1\dagger}
  C_{\ast}^{\bullet}(\iota_1^{\ast}\nbigm,\Gamma_{\geq 2},\gbigk_{\vecy,\geq 2})[-1].
\end{equation}
If $1\in I$,
there exists the morphism
\[
\nbigm[\ast \gbigk_{\vecy}(I)]\otimes V(I)
\lrarr
\Bigl(
\nbigm[\ast \gbigk_{\vecy}(I)]\big/
\nbigm[\ast \gbigk_{\vecy}(I\setminus\{1\})]
\Bigr)\otimes V(I\setminus\{1\})
\]
induced by the natural projection
$\nbigm[\ast \gbigk_{\vecy}(I)]\lrarr
\nbigm[\ast\gbigk_{\vecy}(I)]/\nbigm[\ast\gbigk_{\vecy}(I\setminus\{1\})]$
and $e(1)^{\dagger}$.
If $1\not\in I$,
we have the trivial map
$\nbigm[\ast \gbigk_{\vecy}(I)]\otimes V(I)\lrarr 0$.
They induce a quasi-isomorphism of complexes
\[
 C_{\ast}^{\bullet}(\nbigm,\Gamma,\gbigk_{\vecy})
\lrarr
 C_{\ast}^{\bullet}\bigl(\nbigm[\ast\gbigk_{\vecy}(\Gamma_{\leq 1})]/\nbigm,
 \Gamma_{\geq 2},\gbigk_{\vecy,\geq 2}
 \bigr)[-1].
\]
There exists the isomorphism
\[
 \nbigm[\ast \gbigk_{\vecy}(\Gamma_{\leq 1})]/\nbigm
 \simeq
 \lambda^{-1}\iota_{1\dagger}\iota_1^{\ast}(\nbigm)
\]
induced by
$[m/y_1]\longleftrightarrow
\lambda^{-1}\iota_{1\ast}\iota_1^{\ast}\bigl(m(dy_1/\lambda)^{-1}\bigr)$.
Thus, we obtain the quasi-isomorphism (\ref{eq;21.4.9.1}).
In the same way,
we obtain the following quasi-isomorphisms:
\begin{equation}
\label{eq;21.4.9.4}
C^{\bullet}_{\ast}\bigl(
\iota_j^{\ast}(\nbigm),\,\,
\Gamma_{\geq j+1},\gbigk_{\vecy,\geq j+1}
\bigr)
\lrarr
\lambda^{-1}
\iota_{j,j+1\dagger}
C^{\bullet}_{\ast}\bigl(
\iota_{j+1}^{\ast}(\nbigm),\,\,
\Gamma_{\geq j+2},\gbigk_{\vecy,\geq j+2}
\bigr)[-1].
\end{equation}

Let us construct a quasi-isomorphism
\begin{equation}
 \label{eq;21.4.9.2}
  \iota_{1\dagger}
  C_{!}^{\bullet}(\iota_1^{\ast}\nbigm,\Gamma_{\geq 2},\gbigk_{\vecy,\geq 2})[1]
  \lrarr
  C_{!}^{\bullet}(\nbigm,\Gamma,\gbigk_{\vecy}).
\end{equation}
For $I\subset\Gamma_{\geq 2}$,
there exists the morphism
\[
\Ker
\Bigl(
\nbigm[!\gbigk_{\vecy}(I\sqcup\{1\})]
\lrarr
\nbigm[!\gbigk_{\vecy}(I)]
\Bigr)\otimes V(I)
\lrarr
\nbigm[!\gbigk_{\vecy}(I\sqcup\{1\})]\otimes V(I\sqcup\{1\})
\]
induced by the natural inclusion
$\Ker\Bigl(
\nbigm[!\gbigk_{\vecy}(I\sqcup\{1\})]
\lrarr
\nbigm[!\gbigk_{\vecy}(I)]
\Bigr)
\lrarr
\nbigm[!\gbigk_{\vecy}(I\sqcup\{1\})]$
and $e(1)$.
They induce a quasi-isomorphism of complexes
\[
 C_{!}^{\bullet}\Bigl(
 \Ker\Bigl(
 \nbigm[!\gbigk_{\vecy}(\Gamma_{\leq 1})]
 \lrarr
 \nbigm
 \Bigr),\,
 \Gamma_{\geq 2},\gbigk_{\vecy,\geq 2}
 \Bigr)[1]
\lrarr
 C_{!}^{\bullet}(\nbigm,\Gamma,\gbigk_{\vecy})
\]
There exists the isomorphism
\[
 \Ker\Bigl(
 \nbigm[!\gbigk_{\vecy}(\Gamma_{\leq 1})]
 \lrarr
 \nbigm
 \Bigr)
 \simeq
 \iota_{1\dagger}\iota_1^{\ast}(\nbigm)
\]
induced by
$\deldel_{y_1}\otimes m-1\otimes(\deldel_{y_1}m)
\longleftrightarrow
\iota_{1\ast}\iota_1^{\ast}\bigl(m(dy_1/\lambda)^{-1}\bigr)$.
Thus, we obtain the quasi-isomorphism (\ref{eq;21.4.9.2}).
In the same way,
we obtain the following quasi-isomorphisms:
\begin{equation}
\label{eq;21.4.9.5}
\iota_{j,j+1\dagger}
C^{\bullet}_{!}\bigl(
\iota_{j+1}^{\ast}(\nbigm),\,\,
\Gamma_{\geq j+2},\gbigk_{\vecy,\geq j+2}
\bigr)[1]
\lrarr  
C^{\bullet}_{!}\bigl(
\iota_j^{\ast}(\nbigm),\,\,
\Gamma_{\geq j+1},\gbigk_{\vecy,\geq j+1}
\bigr).
\end{equation}

By an inductive argument,
we obtain the following lemma.
\begin{lem}
We have
$\nbigh^k_{\ast}(\nbigm,\Gamma,\gbigk_{\vecy})=0$
and
$\nbigh^{-k}_!(\nbigm,\Gamma,\gbigk_{\vecy})=0$
unless $k=\ell$.
\hfill\qed
\end{lem}

From (\ref{eq;21.4.9.4}),
we obtain the following isomorphism:
\begin{equation}
\label{eq;21.4.7.30}
\rho_{\vecy,\ast}:
 \lambda^{-\ell}i_{X\dagger}i_X^{\ast}(\nbigm)
 \simeq
 \nbigh^{\ell}_{\ast}\bigl(\nbigm,\Gamma,\gbigk_{\vecy}\bigr).
\end{equation}
From (\ref{eq;21.4.9.5}),
we obtain the following isomorphism:
\begin{equation}
\label{eq;21.4.7.31}
\rho_{\vecy,!}:
 \nbigh^{-\ell}_{!}(\nbigm,\Gamma,\gbigk_{\vecy})
 \simeq
 i_{X\dagger}i_X^{\ast}(\nbigm).
\end{equation}
Let us describe the isomorphisms more explicitly.

For $I\subset \Gamma=\{1,\ldots,\ell\}$,
let $\lefttop{I}V\nbigr_Y\subset\nbigr_Y$ denote
the sheaf of subalgebras generated by
$\lambda p_Y^{\ast}\Theta_Y(\log \gbigk_{\vecy}(I))$.

\begin{lem}
For any $I\subset \Gamma$, we have
\begin{equation}
\label{eq;21.4.7.35}
 \nbigm[\ast \gbigk_{\vecy}(I)](\ast H_Y)
=\nbigr_Y(\ast H_Y)\otimes_{\lefttop{I}V\nbigr_Y(\ast H_Y)}
\bigl(
\nbigm(\gbigk_{\vecy}(I))
\bigr),
\end{equation}
\begin{equation}
\label{eq;21.4.7.32}
 \nbigm[!\gbigk_{\vecy}(I)](\ast H_Y)
 =\nbigr_Y(\ast H_Y)\otimes_{\lefttop{I}V\nbigr_Y(\ast H_Y)}
 \nbigm.
\end{equation}
Moreover,
$\bigcap_{i\in\Gamma\setminus I}\{y_i=b_i\}$
is strictly non-characteristic for
$\nbigm[\star \gbigk_{\vecy}(I)](\ast H_Y)$
for any $\vecb\in\cnum^{\Gamma\setminus I}$.
For $k\in I$,
there exist the following exact sequences:
\begin{equation}
\label{eq;21.4.7.33}
 0\lrarr
 \nbigm[\ast\gbigk_{\vecy}(I\setminus\{k\})](\ast H_Y)
 \lrarr
 \nbigm[\ast\gbigk_{\vecy}(I)](\ast H_Y)
 \lrarr
 \lambda^{-1}j_{k\dagger}j_k^{\ast}
 (\nbigm[\ast\gbigk_{\vecy}(I\setminus\{k\})](\ast H_Y))
 \lrarr 0
\end{equation}
\begin{equation}
\label{eq;21.4.7.34}
 0\lrarr
 j_{k\dagger}j_k^{\ast}
 (\nbigm[\ast\gbigk_{\vecy}(I\setminus\{k\})](\ast H_Y))
 \lrarr
 \nbigm[\ast\gbigk_{\vecy}(I)](\ast H_Y)
 \lrarr
 \nbigm[\ast\gbigk_{\vecy}(I\setminus\{k\})](\ast H_Y)
 \lrarr 0
\end{equation}
Here, $j_k$ denotes the inclusion $\{y_k=0\}\lrarr Y$.
\end{lem}
\pf
It is enough to study the case
$I=\Gamma_{\leq j_0}\subset\Gamma$.
We set $I_0:=\Gamma_{\leq j_0-1}$.
Suppose that we have already proved
the claims for $I_0$.
The $V$-filtration of
the $\nbigrtilde_{Y}(\ast H_Y)$-module
$\nbigm[\ast \gbigk_{\vecy}(I_0)](\ast H_Y)$
along $y_{j_0}$
is given as follows:
\begin{equation}
V_{a}(\nbigm[\ast \gbigk_{\vecy}(I_0)](\ast H_Y))
 =\left\{
   \begin{array}{ll}
 y_{j_0}^{[-a-1]}
  \nbigm[\ast \gbigk_{\vecy}(I_0)](\ast H_Y)
      & (a\leq -1) \\
     \nbigm[\ast \gbigk_{\vecy}(I_0)](\ast H_Y)
 & (a>-1).
    \end{array}  \right.
\end{equation}
We obtain
\begin{equation}
\label{eq;21.4.9.10}
\nbigm[\ast \gbigk_{\vecy}(I)](\ast H_Y)
=\nbigr_Y(\ast H_Y)\otimes_{\lefttop{j_0}V\nbigr_Y(\ast H_Y)}
  \nbigm[\ast \gbigk_{\vecy}(I_0)](\ast H_Y),
\end{equation}
from which we obtain (\ref{eq;21.4.7.35}) for $I$.
We also obtain the exact sequence (\ref{eq;21.4.7.33})
for $I$ and $j_0$,
from which we obtain that
$\bigcap_{i\in\Gamma\setminus I}\{y_i=b_i\}$ $(\vecb\in\cnum^{\Gamma\setminus I})$
are non-characteristic for
$\nbigm[\ast \gbigk_{\vecy}(I)](\ast H_Y)$.
We obtain the claims for 
$\nbigm[!\gbigk_{\vecy}(I)](\ast H_Y)$
similarly.
\hfill\qed

\vspace{.1in}

Let $m$ be a section of $\nbigm$.
We have the section
$y_1^{-1}\cdots y_{\ell}^{-1}m\otimes v_1\wedge\cdots\wedge v_{\ell}$
of $\nbigm[\ast \gbigk_{\vecy}(\Gamma)]\otimes V(\Gamma)$.
The induced section of
$\nbigh_{\ast}^{\ell}(\nbigm,\Gamma,\gbigk)$
is denoted by
$\bigl[
y_1^{-1}\cdots y_{\ell}^{-1}m\otimes(v_{1}\wedge\cdots \wedge v_{\ell})
\bigr]$.
Then, we have
\begin{equation}
\label{eq;21.4.12.20}
 \rho_{\vecy,\ast}\bigl(
  i_{X\ast}i_{X}^{\ast}m\otimes (dy_{\ell}\wedge\cdots \wedge dy_1)^{-1}_{|X}
  \bigr)
  =\bigl[y_1^{-1}\cdots y_{\ell}^{-1}m\otimes
  (v_{1}\wedge\cdots \wedge v_{\ell})
  \bigr].
\end{equation}

For a section $m$ of $\nbigm$,
we obtain the section
\[
 m(\Gamma,\gbigk_{\vecy}):=
 \sum_{I\sqcup J=\Gamma}
 (-1)^{|J|}
 \prod_{i\in I}\deldel_{y_i}
 \otimes
 \prod_{i\in J}\deldel_{y_i}m
\]
of $\nbigm[!\gbigk_{\vecy}(\Gamma)]$.
It is contained in the intersections of
the kernel of
$\nbigm[!\gbigk_{\vecy}(\Gamma)]
\lrarr
\nbigm[!\gbigk_{\vecy}(\Gamma\setminus\{j\})]$
for any $j\in\Gamma$.
It induces the section
$\bigl[m(\Gamma,\gbigk_{\vecy})\otimes (v_1\wedge\cdots \wedge v_{\ell})
 \bigr]$ of
 $\nbigh^{-\ell}_!(\nbigm,\Gamma,\gbigk_{\vecy})$.
We have
\begin{equation}
\label{eq;21.4.12.21}
\rho_{\vecy,!}\Bigl(
\bigl[m(\Gamma,\gbigk_{\vecy})\otimes (v_1\wedge\cdots \wedge v_{\ell})
 \bigr]
 \Bigr)
 = i_{X\ast}i_X^{\ast}(m)
 \bigl((dy_{\ell}/\lambda)\wedge\cdots \wedge(dy_{1}/\lambda)\bigr)^{-1}_{|X}.
\end{equation}
As a result of the formulas (\ref{eq;21.4.12.20})
and (\ref{eq;21.4.12.21}),
we obtain the following lemma.
\begin{lem}
$\rho_{\vecy,\star}$ $(\star=\ast,!)$ are
independent of the choice of
the order
of $y_1,\ldots,y_{\ell}$.
\hfill\qed
\end{lem}

Let $(z_1,\ldots,z_n)$ be another coordinate system of $Y$
such that $X=\bigcap_{i=1}^{\ell}\{z_i=0\}$
and that $\bigcap_{i=1}^{\ell}\{z_i=c_i\}$ $(\vecc\in\cnum^{\ell})$
are non-characteristic for $\nbigm$.
We obtain $\gbigk_{\vecz}:\Gamma\lrarr \nbigs(Y)$
by setting $\gbigk_{\vecz}(i)=\{z_i=0\}$.
We obtain the isomorphisms
$\rho_{\vecz,\ast}:
 \lambda^{-\ell}i_{X\dagger}i_X^{\ast}\nbigm
 \simeq
 \nbigh^{\ell}_{\ast}(\nbigm,\Gamma,\gbigk_{\vecz})$
and
$\rho_{\vecz,!}:
 \nbigh^{-\ell}_!(\nbigm,\Gamma,\gbigk_{\vecz})
 \simeq
 i_{X\dagger}i_X^{\ast}\nbigm$.

\begin{lem}
\label{lem;21.4.9.13}
 We have
$\rho_{\vecy,\ast}=
g_{(\Gamma,\gbigk_{\vecy}),(\Gamma,\gbigk_{\vecz}),\ast}\circ
\rho_{\vecz,\ast}$
and
 $\rho_{\vecy,!}=
 \rho_{\vecz,!}\circ g_{(\Gamma,\gbigk_{\vecz}),(\Gamma,\gbigk_{\vecy}),!}$.
\end{lem}
\pf
We explain
$\rho_{\vecy,\ast}=
g_{(\Gamma,\gbigk_{\vecy}),(\Gamma,\gbigk_{\vecz}),\ast}\circ
\rho_{\vecz,\ast}$,
and the other can be argued similarly.
It is enough to consider the case where
$y_i=z_i$ for $i\neq \ell$.
We set $\Gammatilde:=\{1,\ldots,\ell+1\}$,
$\Gamma_1:=\Gammatilde\setminus\{\ell+1\}$
and $\Gamma_2:=\Gammatilde\setminus\{\ell\}$.
Let $\varphi_i:\Gamma_i\lrarr \Gammatilde$
denote the natural inclusions.
We define $\gbigktilde:\Gammatilde\lrarr\nbigs(Y)$
by $\gbigktilde(i)=\{y_i=0\}$ $(1\leq i\leq \ell)$
and $\gbigktilde(\ell+1)=\{z_{\ell}=0\}$.
We obtain $\gbigk_i:\Gamma_i\lrarr\nbigs(Y)$
as the restriction of $\gbigktilde$.
We have 
$(\Gamma_1,\gbigk_1)=(\Gamma,\gbigk_{\vecy})$.
By the isomorphism
$\rho:\Gamma_2\simeq \Gamma$
defined as
$\rho(i)=i$ $(1\leq i\leq \ell-1)$
and $\rho(\ell+1)=\ell$,
we can identify
$(\Gamma_2,\gbigk_2)$ with
$(\Gamma,\gbigk_{\vecz})$.

We set
$\Gammatilde_{\geq \ell}:=\{\ell,\ell+1\}$.
We have the natural inclusions
$\varphi_3:\{\ell\}\lrarr\Gammatilde_{\geq \ell}$
and
$\varphi_4:\{\ell+1\}\lrarr\Gammatilde_{\geq\ell+1}$.
We set $x_{\ell}=y_{\ell|Y_{\leq \ell-1}}$
and $x_{\ell+1}:=z_{\ell|Y_{\leq \ell-1}}$.
Note that $\{x_{\ell}=0\}=\{x_{\ell+1}=0\}=X\subset Y_{\leq\ell-1}$.
We define $\gbigktilde':\Gammatilde_{\geq \ell}\lrarr\nbigs(Y_{\leq\ell-1})$
by
$\gbigktilde'(\ell)=\gbigktilde'(\ell+1)=X$.
We set $\nbign:=\iota_{\ell-1}^{\ast}\nbigm
\in\nbigc(Y_{\leq \ell-1};H_{Y_{\leq \ell-1}})$.
There exists the following diagram of the isomorphisms:
\begin{equation}
\label{eq;21.4.9.12}
 \begin{CD}
  \nbigh^{\ell}_{\ast}(\nbigm,\Gamma_1,\gbigk_1)
  @<{\varphi^{\ast}_1}<<
  \nbigh^{\ell}_{\ast}(\nbigm,\Gammatilde,\gbigktilde)
  @>{\varphi^{\ast}_2}>>
  \nbigh^{\ell}_{\ast}(\nbigm,\Gamma_2,\gbigk_2)\\
  @V{a_1}VV @VVV @V{a_2}VV \\
  \iota_{\ell-1\dagger}
  \nbigh^1_{\ast}(\nbign,
  \{\ell\},x_{\ell})
  @<{\varphi_3^{\ast}}<<
  \iota_{\ell-1\dagger}
  \nbigh^1_{\ast}\bigl(
  \nbign,\Gammatilde_{\geq \ell},\gbigktilde_{\geq \ell}
  \bigr)
  @>{\varphi_4^{\ast}}>>
  \iota_{\ell-1\dagger}
  \nbigh^1_{\ast}(\nbign,
  \{\ell+1\},x_{\ell+1}).
 \end{CD}
\end{equation}
We have the isomorphism
$b_1:i_{\ell-1,\ell\dagger}(i_X^{\ast}\nbigm)
\simeq
\nbigh^1\bigl(
\nbign,\{\ell\},x_{\ell}
\bigr)$
induced by
$(y_{\ell},y_{\ell+1},\ldots,y_n)_{|Y_{\leq \ell-1}}$
for which
$a_1^{-1}\circ b_1=\rho_{\vecy,\ast}$.
Similarly,
we have the isomorphism
$b_2:i_{\ell-1,\ell\dagger}(i_X^{\ast}\nbigm)
\simeq
\nbigh^1\bigl(
\nbign,\{\ell\},x_{\ell+1}
\bigr)$
induced by
$(z_{\ell},z_{\ell+1},\ldots,z_n)_{|Y_{\leq \ell-1}}$
for which
$a_2^{-1}\circ b_2=\rho_{\vecz,\ast}$.
By Corollary \ref{cor;21.4.12.22},
Proposition \ref{prop;21.4.9.11}
and the commutativity of the diagram (\ref{eq;21.4.9.12}),
we obtain the claim of Lemma \ref{lem;21.4.9.13}.
\hfill\qed

\begin{cor}
The isomorphisms
$\rho_{\vecy,\star}$
are independent of the choice of 
 a holomorphic coordinate system
 $(y_1,\ldots,y_n)$ such that
 $X=\bigcap_{i=1}^{\ell}\{y_i=0\}$
and that $\bigcap_{i=1}^{\ell}\{y_i=a_i\}$ $(\veca\in\cnum^{\ell})$
are non-characteristic for $\nbigm$.
\hfill\qed 
\end{cor}

We obtain the following isomorphisms:
\begin{multline}
 i_{X\dagger}(\DD_{X(\ast H_X)}i_X^{\ast}\nbigm)
 =\DD_{Y(\ast H_Y)}(i_{X\dagger}i_X^{\ast}\nbigm)
 \stackrel{b_1}{\simeq}
 \DD_{Y(\ast H_Y)}\nbigh^{-\ell}_!(\nbigm,\Gamma,\gbigk)
 \\
 \stackrel{b_2}{\simeq}
 \nbigh^{\ell}_{\ast}(\DD_{Y(\ast H_Y)}\nbigm,\gbigk)
 \stackrel{b_3}{\simeq}
 \lambda^{-\ell}i_{X\dagger}i_X^{\ast}(\DD_{Y(\ast H_Y)}\nbigm).
\end{multline}
It is induced by an isomorphism
\begin{equation}
\label{eq;21.4.14.4}
 \DD_{X(\ast H_X)}(i_X^{\ast}\nbigm)
 \simeq
\lambda^{-\ell}i_X^{\ast}(\DD_{Y(\ast H_Y)}\nbigm).
\end{equation}
It is independent of the choice of a coordinate system
$(y_1,\ldots,y_n)$
satisfying Condition \ref{condition;21.4.12.10}.

\vspace{.1in}

Let $(\nbigt,W)=((\nbigm',\nbigm'',C),W)\in\MTM^{\integral}(Y;H_Y)$.
Suppose that
$X_{\veca}$ $(\veca\in\cnum^{\ell})$
are non-characteristic for $\nbigt$.
We obtain
$\nbigh^{\ell}(\nbigt,\Gamma,\gbigk_{\vecy})
\in\MTM^{\integral}(Y;H_Y)$,
which consists of
the $\nbigrtilde_Y(\ast H_Y)$-modules
$\nbigh^{-\ell}_!(\nbigm',\Gamma,\gbigk_{\vecy})$
and
$\nbigh^{\ell}_{\ast}(\nbigm'',\Gamma,\gbigk_{\vecy})$
with the induced sesqui-linear pairing
$\nbigh^{\ell}(C,\Gamma,\gbigk_{\vecy})$
and the weight filtration.
There exists a unique sesqui-linear pairing
$i_X^{\ast}(C)$
of $i_X^{\ast}(\nbigm')$
and $\lambda^{-\ell}i_X^{\ast}(\nbigm'')$
such that
$i_{X\dagger}(i_X^{\ast}(C))
=\nbigh^{\ell}(C,\Gamma,\gbigk_{\vecy})$
under the isomorphisms
$\rho_{\vecy,!}$ for $\nbigm'$
and $\rho_{\vecy,\ast}$ for $\nbigm''$.
Thus, we obtain
the filtered $\nbigrtilde_X(\ast H_X)$-triple
\[
 (i_X^{\ast}\nbigt,W)
=((i_X^{\ast}\nbigm',\lambda^{-\ell}i_X^{\ast}\nbigm'',i_X^{\ast}C),W)
\]
with an isomorphism
$i_{X\dagger}(i_X^{\ast}\nbigt,W)\simeq
\nbigh^{\ell}(\nbigt,\Gamma,\gbigk_{\vecy})$.
Note that
$(i_{X}^{\ast}(\nbigt),W)$ is uniquely determined
up to canonical isomorphisms,
independently from
a holomorphic coordinate system
$(y_1,\ldots,y_n)$ such that $X=\bigcap_{i=1}^{\ell}\{y_i=0\}$
and that $\bigcap_{i=1}^{\ell}\{y_i=a_i\}$ $(\veca\in\cnum^{\ell})$
are strictly non-characteristic for $\nbigt$.
The following lemma is obvious by the construction
and Proposition \ref{prop;21.4.9.21}.
\begin{lem}
 $(i_X^{\ast}\nbigt,W)\in\MTM^{\integral}(X;H_X)$.
\hfill\qed
\end{lem}

Let $(\nbigt_1,W)\lrarr(\nbigt_2,W)$ be any morphism 
in $\MTM^{\integral}(Y;H_Y)$.
Assume that $\bigcap_{i=1}^{\ell}\{y_i=a_i\}$ $(\veca\in\cnum^{\ell})$
are strictly non-characteristic for $\nbigt_i$.
We obtain the induced morphism
$(i_X^{\ast}(\nbigt_1),W)\lrarr (i_X^{\ast}(\nbigt_2),W)$,
which is independent of the choice of $(y_1,\ldots,y_n)$.

Let $\nbigt_0\in\MTM^{\integral}(Y)$ be pure of weight $w$
such that $\bigcap_{i=1}^{\ell}\{y_i=a_i\}$ $(\veca\in\cnum^{\ell})$
are strictly non-characteristic for $\nbigt_0(\ast H_Y)$.
Let $\nbigs_0:\nbigt_0\lrarr \nbigt_0^{\ast}(-w)$
be a polarization.
We obtain
\[
i_X^{\ast}(\nbigs_0):
i_X^{\ast}(\nbigt_0(\ast H_Y))
\lrarr
i_X^{\ast}\bigl(\nbigt_0^{\ast}(-w)(\ast H_Y)\bigr)
=i_X^{\ast}\bigl(\nbigt_0(\ast H_Y)\bigr)^{\ast}(-w-\ell),
\]
which is independent of $(y_1,\ldots,y_n)$.
The following lemma is obvious by the construction.
 \begin{lem}
  Let $\nbigt_X$ denote the image of
$i_X^{\ast}(\nbigt_0(\ast H_Y))[!H_Y]
\lrarr
i_X^{\ast}(\nbigt_0(\ast H_Y))[\ast H_Y]$.  
Then $\nbigt_X$
is pure and polarizable of weight $w+\ell$.
It is equipped with a unique polarization
 $\nbigs_X$ which induces
 $i_X^{\ast}(\nbigs_0)$
 on $\nbigt_X(\ast H_X)=i_X^{\ast}(\nbigt_0(\ast H_Y))$.
\hfill\qed
 \end{lem}

\subsubsection{Proof of Theorem \ref{thm;21.4.7.50}}
\label{subsection;21.4.8.1}

If $F$ is the projection of
$X=Y\times Z$ onto $Y$ for a complex manifold $Z$,
then $F^{\ast}(\nbigm)$ is an object of $\nbigc(X;H_X)$
by the functoriality
with respect to the external tensor product
in \S\ref{subsection;21.4.14.3}.
We also obtain
$\DD_{X(\ast H_X)}(F^{\ast}\nbigm)\simeq
\lambda^{\dim Z}
F^{\ast}(\DD_{Y(\ast H_Y)}\nbigm)$
as in \cite[Proposition 13.3.9]{Mochizuki-MTM}.
Hence, it is enough to study the case where
$F$ is a closed embedding.
We set $\ell=\dim Y-\dim X$.

\vspace{.1in}

Let $(\nbigt,W)=((\nbigm',\nbigm,C),W)\in\MTM^{\integral}(Y;H_Y)$.
Assume that $F$ is non-characteristic for $\nbigt(\ast H_Y)$.
For any $P\in F(X)$,
there exists a holomorphic coordinate neighbourhood
$(Y_P,z_1,\ldots,z_n)$
such that $F(X)\cap Y_P=\bigcap_{i=1}^{\ell}\{z_i=0\}$.
We may assume that Condition \ref{condition;21.4.12.10}
is satisfied for $\nbigt_{|Y_P}$
with respect to the coordinate system $(z_1,\ldots,z_n)$.
We set $X_P:=F^{-1}(Y_P)$ and $H_{X_P}:=F^{-1}(H_Y)$,
and let $F_P:X_P\lrarr Y_P$ denote the induced morphism.
We obtain
$(F_P^{\ast}(\nbigt_{|Y_P}),W)=\bigl(
(F_P^{\ast}(\nbigm'_{|Y_P}),
\lambda^{-\ell}F_P^{\ast}(\nbigm_{|Y_P}),
F_P^{\ast}(C_{|Y_P})),W
\bigr)
\in\MTM^{\integral}(X_P,H_P)$,
which is independent of
the coordinate system $(z_1,\ldots,z_n)$.
There exist natural isomorphisms
$F_P^{\ast}(\nbigm'_{|Y_P})
\simeq
 F^{\ast}(\nbigm')_{|X_P}$
and
$F_P^{\ast}(\nbigm_{|Y_P})
\simeq
 F^{\ast}(\nbigm)_{|X_P}$.
By varying $P\in F(X)$,
and patching $(F_P^{\ast}(\nbigt_{|X_P}),W)$,
we obtain
the filtered $\nbigrtilde_X(\ast H_X)$-triple
$(F^{\ast}(\nbigt),W)
=((F^{\ast}\nbigm',F^{\ast}\nbigm,F^{\ast}C),W)$.

Suppose that there exists
$\nbigt_0\in \MTM^{\integral}(Y)$
such that
(i) $\nbigt_0$ is pure of weight $w$
with a polarization $\nbigs_0$,
(ii) $\nbigt_0(\ast H_Y)=\nbigt$.
We obtain the induced morphism
$F^{\ast}\nbigs_0:
F^{\ast}(\nbigt)\simeq
F^{\ast}(\nbigt)^{\ast}(-w-\ell)$.

\begin{lem}
There exists $\nbigt_X\in\MTM^{\integral}(X)$
such that
(i) $\nbigt_X$ is pure of weight $w+\ell$
 with a polarization $\nbigs_X$,
 (ii) $\nbigt_X(\ast H_X)\simeq F^{\ast}(\nbigt)$
 and $\nbigs_X=F^{\ast}(\nbigs)$ under the isomorphism.
\end{lem}
\pf
For each $P\in F(X)$,
let $\nbigt_{X_P}\in\MTM^{\integral}(X_P)$
be the image of
$F_P^{\ast}(\nbigt_{|Y_P})[! H_{X,P}]
\lrarr
F_P^{\ast}(\nbigt_{|Y_P})[\ast H_{X,P}]$.
Then, it is pure of weight $w+\ell$,
and there exists a polarization $\nbigs_{X_P}$
of $\nbigt_{X_P}$
which induces $F_P^{\ast}(\nbigs_{|F_P})$
under the isomorphism
$\nbigt_{X_P}(\ast H_{X_P})
\simeq
F_P^{\ast}(\nbigt_{|Y_P})$.
By patching $(\nbigt_{X_P},\nbigs_{X_P})$,
we obtain an $\nbigrtilde_{X}$-triple
$\nbigt_X$ with a morphism
$\nbigs_X:\nbigt_X\simeq \nbigt_X^{\ast}\otimes\newTate(-w-\ell)$.
It is easy to check that
$(\nbigt_X,\nbigs_X)$ is a polarized integrable pure twistor $\nbigd_X$-module
of weight $w+\ell$.
\hfill\qed

\vspace{.1in}

In general,
note that 
$W_{w}F^{\ast}(\nbigt):=F^{\ast}(W_{w+\ell}\nbigt)$
by the construction.
By using Proposition \ref{prop;21.4.9.20},
we can prove that
$(F^{\ast}\nbigt,W)\in\MTM^{\integral}(X;H_X)$.
It implies
$F^{\ast}(\nbigm)\in \nbigc(X;H_X)$,
which is the first claim of Theorem \ref{thm;21.4.7.50}.

\vspace{.1in}
By patching the isomorphisms (\ref{eq;21.4.14.4})
for $F_P^{\ast}(\nbigm_{|Y_P})$ for varied $P$,
we obtain the desired isomorphism
\[
 \DD_{X(\ast H_X)}(F^{\ast}\nbigm)
\simeq
 \lambda^{-\dim Y+\dim X}F^{\ast}(\DD_{Y(\ast H_Y)}\nbigm).
\]
Thus, we obtain Theorem \ref{thm;21.4.7.50}.
\hfill\qed

\subsubsection{Appendix}

Let $X$ be any complex manifold
with a closed complex hypersurface $H_X$.
Let $\nbigt$ be an $\nbigrtilde_{X(\ast H_X)}$-triple
equipped with a filtration $W$
and morphisms
       $\nbigs_w:\Gr^W_w(\nbigt)
       \lrarr\Gr^W_w(\nbigt)^{\ast}(-w)$ $(w\in\seisuu)$
such that the following holds.
 \begin{itemize}
 \item $\Gr^W_w(\nbigt)_{|X\setminus H_X}$ are
	polarizable pure twistor $\nbigd_{X\setminus H_X}$-modules
       of weight $w$,
       and
       $\nbigs_{w|X\setminus H_X}$
       are polarizations of
       $\Gr^W_w(\nbigt)_{|X\setminus H_X}$.
 \item Each $P\in H_X$ has 
       a neighbourhood $X_P$ in $X$
       such that
       $(\nbigt_P,W):=(\nbigt,W)_{|X_P}\in \MTM^{\integral}(X_P;H_{X,P})$,
       where $H_{X,P}:=H_X\cap X_P$.
       Moreover,
       there exist polarizations
       $\nbigs_{P,w}$ of
       \[
       \Gr^W_{w}(\nbigt_P)_{!\ast}:=
       \Image\Bigl(
       \Gr^W_w(\nbigt_{P})[!H_{X,P}]
       \lrarr
       \Gr^W_{w}(\nbigt_P)[\ast H_{X,P}]
       \Bigr)
       \]
       which induce
       $\nbigs_{w|X_P}$
       under the isomorphism
       $\Gr^W_{w}(\nbigt_P)_{!\ast}(\ast H_{X,P})
       =\Gr^W(\nbigt)_{|X_P}$.
\end{itemize}

\begin{prop}
\label{prop;21.4.9.20}
 $(\nbigt,W)\in\MTM^{\integral}(X;H_X)$,
i.e.,
there exists $(\nbigt_0,W)\in\MTM^{\integral}(X)$
such that
$(\nbigt,W)=(\nbigt_0,W)(\ast H_X)$.
\end{prop}
\pf
Let us consider the case where
there exists $w_0$ such that
$\Gr^W_w(\nbigt)=0$ unless $w=w_0$.
For each $P\in H_X$,
let $\nbigt_{P,!\ast}$ denote the image of
$\nbigt_{P}[!H_{X,P}]\lrarr
\nbigt_P[\ast H_{X,P}]$,
which is a pure twistor $\nbigd_{X_P}$-module
of weight $w_0$
equipped with the polarization $\nbigs_{P}$
as in the assumption.
By patching
$\nbigt_{|X\setminus H_X}$
and $\nbigt_{P,!\ast}$
for varied $P\in H_X$,
we obtain an $\nbigrtilde_X$-triple $\nbigt_{!\ast}$
equipped with a morphism
$\nbigs_{w_0,!\ast}:
\nbigt_{!\ast}\lrarr\nbigt_{!\ast}(-w)$.
It is easy to check that
$\nbigt_{!\ast}$ is a pure twistor $\nbigd_X$-module
of weight $w_0$
with a polarization $\nbigs_{w_0,!\ast}$.

Before considering the mixed case,
we recall that,
by regarding the hypersurface $H_X$ with the reduced divisor
as an effective divisor of $X$,
there exists the following graded $\nbigrtilde_X$-triple
(see \cite[\S4.4.3]{Mochizuki-MTM}):
\[
 \Gr^W\psi^{(0)}_{H_X}(\nbigt_{!\ast})
=\bigoplus_w\Gr^W_{w}\psi^{(0)}_{H_X}(\nbigt_{!\ast}).
\]
Here, we adjust that
$\Gr^W_w\psi^{(0)}_{H_X}(\nbigt_{!\ast})$
are pure of weight $w$.
Note that there exists the polarizations
$\nbigs_{w_0,w,H_X}$ on
$\Gr^W_w\psi^{(0)}_{H_X}(\nbigt_{!\ast})$
induced by $\nbigs_{!\ast}$.
(See \cite[\S7.1.3.4]{Mochizuki-MTM}.)

Let us consider the mixed case.
For each $P\in H_X$,
we obtain
$(\nbigt_P,W)[\ast H_{X,P}]
\in\MTM^{\integral}(X_P)$.
By patching
$(\nbigt,W)_{|X\setminus H_X}$
and 
$(\nbigt_P,W)[\ast H_{X,P}]$ $(P\in H_X)$,
we obtain
the filtered $\nbigrtilde_X$-triple
$(\nbigt_0,W)$.
To prove that
$(\nbigt_0,W)\in\MTM^{\integral}(X)$,
it is enough to prove that
$\Gr^W_w(\nbigt_0)$ are pure and polarizable of weight $w$
because the other conditions can be checked locally around
any points of $X$,
which follows from the properties of
$(\nbigt,W)_{|X\setminus H_X}$
and $(\nbigt_P,W)[\ast H_{X,P}]$.
Because
$\Gr^W_w(\nbigt_0)_{|X\setminus H_X}
=\Gr^W_w(\nbigt)_{|X\setminus H_X}$
and
$\Gr^W_w(\nbigt_0)_{|X_P\setminus H_{X,P}}
=\Gr^W_w(\nbigt_P)_{|X_P\setminus H_{X,P}}$,
$\Gr^W_w(\nbigt_0)$ is strictly $S$-decomposable.
There exists the following decomposition:
\[
 \Gr^W_w(\nbigt_0)
=\Gr^W_w(\nbigt)_{!\ast}
 \oplus
 \nbigp_{w,1},
\]
where the strict support of $\nbigp_{w,1}$
is contained in $H_X$.
We have already observed that
$\Gr^W_w(\nbigt)_{!\ast}$ are pure and polarizable of weight $w$.
We also note that
$\nbigp_{w,1|X_P}$ are pure and polarizable of weight $w$
for any $P\in H_X$.
Hence, it is easy to see that
$\nbigp_{w,1}$ is pure of weight $w$.
It remains to show the existence of a global polarization.

We have the graded $\nbigrtilde_{X}$-triple
\[
\Gr^W\psi^{(0)}_{H_{X}}\bigl(
 \nbigt_0
 \bigr)
=\bigoplus \Gr^W_w\psi^{(0)}_{H_{X}}\bigl(
 \nbigt
 \bigr).
\]
We have the canonical decomposition
\[
 \Gr^W_w\psi^{(0)}_{H_{X}}\bigl(\nbigt\bigr)
=\bigoplus_{w_0}
  \Gr^W_w\psi^{(0)}_{H_X}\bigl(
   \Gr^W_{w_0}(\nbigt)_{!\ast}
  \bigr).
\]
Hence, $\Gr^W_w\psi^{(0)}_{H_X}(\nbigt)$
is equipped with a polarization
$\bigoplus_{w_0}\nbigs_{w_0,w,H_X}$.
There exists the natural monomorphism
$\nbigp_{1,w}\lrarr
 \Gr^W_w\psi^{(0)}_{H_X}(\nbigt)$.
Because $\nbigp_{1,w}$ is pure of weight $w$,
we obtain that
$\nbigp_{1,w}$ is also polarizable.
\hfill\qed

\subsection{Inverse image by closed embeddings under an additional assumption}
\label{subsection;21.6.22.2}

Let $f:X\lrarr Y$ be a morphism of complex manifolds.
Let $H_Y$ be a hypersurface of $Y$
with a decomposition
$H_Y=H_Y^{(1)}\cup H_Y^{(2)}$
such that $\dim H_Y^{(1)}\cap H_Y^{(2)}<\dim H_Y$.
We assume the following condition.
\begin{condition}
\label{condition;21.6.21.1}
The restriction $f_{|X\setminus H_X}:X\setminus H_X\to Y\setminus H_Y$
is a closed embedding.
There exists a finite set $\Gamma$
with a map $\gbigk:\Gamma\lrarr\nbigs(Y)$
such that
$f(X)\cup H_Y=
B(\Gamma,\gbigk)\cup H_Y$.
\hfill\qed
\end{condition}

We set $H_X=f^{-1}(H_Y)$
and $H_X^{(2)}=f^{-1}(H_Y^{(2)})$.
We have the decomposition
$H_X=H_X^{(1)}\cup H_X^{(2)}$
such that $\dim H_X^{(1)}\cap H_X^{(2)}<\dim H_X$.

For $\nbigm\in\nbigc(Y,[\star H_Y^{(1)}];H^{(2)}_Y)$,
we obtain the complex of
$\nbigrtilde_{Y(\ast H^{(2)}_Y)}$-modules
$C_{\star}(\nbigm,\Gamma,\gbigk)$
as in \S\ref{subsection;21.4.14.1}.
(See \S\ref{subsection;22.7.31.1}
for $\nbigc(Y,[\star H_Y^{(1)}];H^{(2)}_Y)$.)
The following lemma is similar to Lemma \ref{lem;21.4.8.2}.
\begin{lem}
For $\nbigm\in\nbigc(Y,[\star H^{(1)}_Y];H_Y^{(2)})$ $(\star=!,\ast)$,
we have
$\nbigh^i
C_{\star}(\nbigm,\Gamma,\gbigk)
\in \nbigc(Y,[\star H^{(1)}_Y];H_Y^{(2)})$,
and 
the support of
$\nbigh^i
C_{\star}(\nbigm,\Gamma,\gbigk)$
are contained in $\cnum_{\lambda}\times F(X)$.
\hfill\qed
\end{lem}

For $\nbigm\in\nbigc(Y,[!H_Y^{(1)}];H^{(2)}_Y)$,
there exist
$(\lefttop{T}f^{\ast})^i(\nbigm)
\in\nbigc(X,[!H_X^{(1)}];H^{(2)}_X)$
with an isomorphism
\[
 f_{\dagger}\bigl(
 (\lefttop{T}f^{\ast})^i(\nbigm)
 \bigr)
 \simeq
 \nbigh^i
  C_!(\nbigm,\Gamma,\gbigk).
\]
For $\nbigm\in\nbigc(Y,[\ast H_Y^{(1)}];H^{(2)}_Y)$,
there exist
$(\lefttop{T}f^{!})^i(\nbigm)
\in\nbigc(X,[\ast H_X^{(1)}];H_X^{(2)})$
with an isomorphism
\[
 f_{\dagger}\bigl(
 (\lefttop{T}f^{!})^i(\nbigm)
 \bigr)
 \simeq
 \nbigh^i
  C_{\ast}(\nbigm,\Gamma,\gbigk).
\]
They are independent of the choice of $(\Gamma,\gbigk)$.

Note that for $\nbigm\in\nbigc(Y;H_Y)$,
we have $\nbigm[\star H_Y^{(1)}]:=(\nbigm[\star H_Y])(\ast H_Y^{(2)})
\in\nbigc(Y,[\star H_Y^{(1)}];H_Y^{(2)})$.
Similarly, for $\nbign\in\nbigc(X;H_X)$,
we have
$\nbign[\star H_X^{(1)}]:=(\nbign[\star H_X])(\ast H_X^{(2)})$.
By the construction, the following lemma is clear.
\begin{lem}
\label{lem;21.6.22.10}
We have 
$(\lefttop{T}f^{!})^i\bigl(\nbigm[\ast H^{(1)}_Y]\bigr)
=(\lefttop{T}f^{!})^i(\nbigm)[\ast H^{(1)}_X]$,
and 
$(\lefttop{T}f^{\ast})^i\bigl(\nbigm[! H^{(1)}_Y]\bigr)
=(\lefttop{T}f^{\ast})^i(\nbigm)[!H^{(1)}_X]$.
\hfill\qed
\end{lem}

\begin{cor}
\label{cor;21.4.21.1}
Suppose that $f$ is non-characteristic to
$\nbigm\in\nbigc(Y;H_Y)$.
\begin{itemize}
 \item 
We have
$(\lefttop{T}f^!)^i(\nbigm[\ast H_Y^{(1)}])=0$
unless $i=\dim Y-\dim X$,
and
\[
(\lefttop{T}f^!)^{\dim Y-\dim X}(\nbigm[\ast H_Y^{(1)}])
=\lambda^{-\dim Y+\dim X}f^{\ast}(\nbigm)[\ast H^{(1)}_X]
\]
 \item
We have
$(\lefttop{T}f^{\ast})^i(\nbigm[!H_Y^{(1)}])=0$
unless $i=\dim X-\dim Y$,
      and
\[
(\lefttop{T}f^{\ast})^{\dim X-\dim Y}(\nbigm[!H_Y^{(1)}])
=f^{\ast}(\nbigm)[!H^{(1)}_X].
\]
\hfill\qed
\end{itemize}
\end{cor}

\subsubsection{Tensor products}
\label{subsection;21.6.22.11}

Let $X$ be a complex manifold with a hypersurface $H$.
We set $\Htilde:=(H\times X)\cup (X\times H)\subset X\times X$.
Let $\Delta_X:X\lrarr X\times X$ denote the diagonal embedding.
\begin{condition}
There exists a finite set $\Gamma$
with a map $\gbigk:\Gamma\lrarr\nbigs(X\times X)$
such that
$\Delta_X(X)\cup\Htilde
=B(\Gamma,\gbigk)\cup\Htilde$.

\hfill\qed
\end{condition}

For $\nbigm_i\in\nbigc(X;H)$ $(i=1,2)$
and for $\star=!,\ast$,
by using
$\nbigm_1\boxtimes\nbigm_2\in\nbigc(\Xtilde;\Htilde)$,
we define
\[
 \nbigh^k\bigl(
 \nbigm_1\otimes^{\star}\nbigm_2
 \bigr):=
 (\lefttop{T}\Delta^{\star})^k
 \bigl(
 \nbigm_1\boxtimes\nbigm_2
 \bigr)
 \in\nbigc(X;H).
\]

Let $H^{(1)}$ be a hypersurface of $X$.
For $\nbigm\in\nbigc(X;H\cup H^{(1)})$,
we have
$\nbigm[\star H^{(1)}]=\bigl(
\nbigm[\star (H^{(1)}\cup H)]
\bigr)(\ast H)\in\nbigc(X;H)$.

\begin{prop}
\label{prop;21.6.22.12}
Let $\nbigm_i\in\nbigc(X;H)$ $(i=1,2)$.
Suppose that $\nbigm_i(\ast H^{(1)})\in\nbigc(X;H\cup H^{(1)})$ $(i=1,2)$
are non-characteristic.
Then, we have
$\nbigh^{\ell}\bigl(
 \nbigm_1
 \otimes^!
 (\nbigm_2[\ast H^{(1)}])\bigr)=0$
unless $\ell=\dim X$,
and 
\[
 \nbigh^{\dim X}\bigl(
 \nbigm_1
 \otimes^!
 (\nbigm_2[\ast H^{(1)}])
 \bigr)
=\lambda^{-\dim X}
 \Bigl(
 \nbigm_1(\ast H^{(1)})\otimes_{\nbigo_{\nbigx}}\nbigm_2(\ast H^{(1)})
 \Bigr)[\ast H^{(1)}].
\] 
We also have
$\nbigh^{\ell}\bigl(
 \nbigm_1
 \otimes^{\ast}
 (\nbigm_2[!H^{(1)}])\bigr)=0$
unless $\ell=-\dim X$,
and 
\[
 \nbigh^{-\dim X}\bigl(
 \nbigm_1
 \otimes^{\ast}
 (\nbigm_2[!H^{(1)}])
 \bigr)
=\Bigl(
 \nbigm_1(\ast H^{(1)})\otimes_{\nbigo_{\nbigx}}\nbigm_2(\ast H^{(1)})
 \Bigr)[!H^{(1)}].
\] 
\end{prop}
\pf
We set $\Htilde^{(1)}:=X\times H^{(1)}$.
By the assumption,
$\Delta_X$ is non-characteristic for 
$(\nbigm_1\boxtimes\nbigm_2)(\ast \Htilde^{(1)})
\in\nbigc(X\times X;\Htilde\cup\Htilde^{(1)})$.
We also have
$(\nbigm_1\boxtimes\nbigm_2)[\ast \Htilde^{(1)}]
=\nbigm_1\boxtimes(\nbigm_2[\ast H^{(1)}])$
and
$(\nbigm_1\boxtimes\nbigm_2)(\ast \Htilde^{(1)})
=\nbigm_1\boxtimes(\nbigm_2(\ast H^{(1)}))$.
By Corollary \ref{cor;21.4.21.1},
we have
$\nbigh^{\ell}\bigl(
 \nbigm_1
 \otimes^!
 (\nbigm_2[\ast H^{(1)}])\bigr)=0$
unless $\ell=\dim X$,
and 
\begin{multline}
 \nbigh^{\dim X}\bigl(
 \nbigm_1
 \otimes^!
 (\nbigm_2[\ast H^{(1)}])
 \bigr)
=\lambda^{-\dim X}
 \Bigl(
 \nbigm_1\otimes_{\nbigo_{\nbigx}}\nbigm_2(\ast H^{(1)})
 \Bigr)[\ast H^{(1)}]
 \\
 \simeq
 \lambda^{-\dim X}
 \Bigl(
 \nbigm_1(\ast H^{(1)})\otimes_{\nbigo_{\nbigx}}\nbigm_2(\ast H^{(1)})
 \Bigr)[\ast H^{(1)}]
\end{multline}
We obtain the other claim similarly.
\hfill\qed

\begin{cor}
\label{cor;21.6.22.13}
For any $\nbigm\in\nbigc(X;H)$,
We have $\nbigh^{\ell}(\nbigm\otimes^!\nbigo_{\nbigx}[\ast H^{(1)}])=0$
unless $\ell=\dim X$,
and there exists a natural isomorphism
\[
       \nbigh^{\dim X}(\nbigm\otimes^!\nbigo_{\nbigx}[\ast H^{(1)}])
       \simeq
       \lambda^{-\dim X}\nbigm[\ast H^{(1)}].
\]       
We also have $\nbigh^{\ell}(\nbigm\otimes^{\ast}\nbigo_{\nbigx}[!H^{(1)}])=0$
unless $\ell=-\dim X$,
and there exists a natural isomorphism
\[
       \nbigh^{-\dim X}(\nbigm\otimes^{\ast}\nbigo_{\nbigx}[!H^{(1)}])
       \simeq
       \nbigm[!H^{(1)}].
\]         
\hfill\qed
\end{cor}

\subsection{$\gbigr$-modules and $\gbigrtilde$-modules}

We set
$\gbigx:=\proj^1\times X$.
Let $p_{1,X}:\gbigx\lrarr X$ denote the projection.
We set
$\gbigx^{\infty}=\{\infty\}\times X\subset\gbigx$
and
$\gbigx^{\lambda}=\{\lambda\}\times X\subset\gbigx$
($\lambda\in\cnum$).
Let $\gbigr_X\subset\nbigd_{\gbigx}(\ast\gbigx^{\infty})$ denote
the sheaf of subalgebras generated by
$\lambda p_{1,X}^{\ast}\Theta_X$
over $\nbigo_{\gbigx}(\ast\gbigx^{\infty})$.
We set
$\gbigrtilde_X:=
\gbigr_X\langle\lambda^2\del_{\lambda}\rangle
\subset\nbigd_{\gbigx}(\ast\gbigx^{\infty})$.
We set $\gbigx^{\circ}:=\gbigx\setminus\gbigx^0$.

For any open subset $\gbigu\subset\gbigx$,
the restriction of an $\gbigr_X$-module $\gbigm$ to $\gbigu$
is denoted by
$\gbigm_{|\gbigu}$.
If $\gbigu=\proj^1\times U$ for an open subset $U\subset X$,
the restriction is also denoted as $\gbigm_{|U}$.

Let $F:X\lrarr Y$ be any morphism of complex manifolds.
For any $\nbigo_{\gbigx}$-module $\gbigm$,
the push-forward
$(\id\times F)_{\ast}(\gbigm)$ is often denoted by
$F_{\ast}(\gbigm)$.
For any $\nbigo_{\gbigy}$-module $\gbign$,
the pull back
$(\id\times F)^{\ast}(\gbign)$
is often denoted by
$F^{\ast}(\gbign)$.

Let $H$ be any closed complex hypersurface of $X$.
For any $\nbigo_{\gbigx}$-module $\gbign$,
we set $\gbign(\ast H):=
\gbign\otimes\nbigo_{\gbigx}\bigl(\ast\gbigh\bigr)$.
We set $\gbigr_{X(\ast H)}:=\gbigr_X(\ast H)$
and $\gbigrtilde_{X(\ast H)}:=\gbigrtilde_X(\ast H)$.
If $H=f^{-1}(0)$ for a holomorphic function $f$,
$\gbign(\ast H)$ is also denoted by $\gbign(\ast f)$.

\subsubsection{Extension of $\nbigrtilde_X$-modules underlying integrable mixed
twistor $\nbigd$-modules}

Let $\gbigc(X)$ denote the full subcategory of
$\gbigrtilde_X$-modules $\gbigm$
satisfying the following condition.
\begin{itemize}
 \item
        $\gbigm_{|\nbigx}$ is an object in $\nbigc(X)$,
  and 
$\gbigm_{|\gbigx^{\circ}}$ is a holonomic
      $\nbigd_{\gbigx^{\circ}}$-module.
\end{itemize}

For a closed complex hypersurface $H$,
there exists the natural functor from
$\gbigc(X)$ to the category of
$\gbigrtilde_{X(\ast H)}$-modules.
Let $\gbigc(X;H)$
denote the essential image.
For any open subset $U\subset X$,
the restriction induces
$\gbigc(X;H)\lrarr\gbigc(U;H\cap U)$.
The following lemma is obvious.
\begin{lem}
Let $g:\gbigm_1\lrarr\gbigm_2$ be a morphism
in $\gbigc(X;H)$.
If $g_{|\nbigx}:\gbigm_{1|\nbigx}\lrarr\gbigm_{2|\nbigx}$
is induced by a morphism in
$\MTM^{\integral}(X;H)$,
$\Ker(g)$, $\Image(g)$ and
$\Cok(g)$ are objects in
$\gbigc(X;H)$.
(See Remark {\rm\ref{rem;21.4.13.1}}.)
\hfill\qed
\end{lem}

We remark the following coherence property.

\begin{lem}
For any $\gbigm\in\gbigc(X;H)$,
there exists a directed family of coherent $\nbigo_{\gbigx}$-submodules
$\{G_i\}_{i\in\Lambda}$ of $\gbigm$
such that
 $\sum_{i\in\Lambda} G_i=\gbigm$.
Here, a directed family means
that for any $i,i''\in\Lambda$ there exists $i''\in\Lambda$
such that $G_{i}+G_{i'}\subset G_{i''}$.
In particular, $\gbigm$
is a good $\nbigo_{\gbigx}$-module
in the sense of {\rm\cite[Definition 4.22]{kashiwara_text}}.
There also exists a coherent filtration
of $\gbigm$ as an $\gbigrtilde_{\gbigx}$-module
over the filtered ring $\gbigrtilde_X$
in the sense of {\rm\cite[Definition A.19]{kashiwara_text}}.
\end{lem}  
\pf
It is enough to consider the case $H=\emptyset$.
Note that
$\gbigm(\ast\lambda)$
is a holonomic $\nbigd_{\gbigx}$-module.
Hence, according to \cite[Theorem 3.1]{malgrange-holonomic-D-modules},
there exists a good filtration
$F_{\bullet}\bigl(
\gbigm(\ast\lambda)
\bigr)$.
Then, we obtain a family of
coherent $\nbigo_{\gbigx}$-submodules
with the desired property
by setting
$G_j:=
F_{j}\bigl(
\gbigm(\ast\lambda)
\bigr)
\cap
\gbigm$.
We can easily construct
a coherent filtration of $\gbigm$
from $\{G_j\}$.
\hfill\qed

\vspace{.1in}

We obtain various operations for objects $\gbigm$ of $\gbigc(X;H)$
obtained as the gluing of operations for objects
$\gbigm_{|\nbigx}\in\nbigc(X;H)$
and 
the corresponding operations for holonomic
$\nbigd_{\gbigx^{\circ}}\bigl(\ast(\gbigx^{\infty}\cup\gbigh^{\circ})\bigr)$-modules
$\gbigm_{|\gbigx^{\circ}}$,
which we shall explain below.

\subsubsection{Direct image}

Let $F:X\lrarr Y$ be a projective morphism of complex manifolds.
We set
$\omega_{\gbigx}:=
\lambda^{-\dim X}p_{1,X}^{\ast}(\omega_{X})(\ast\gbigx^{\infty})$.
Similarly,
$\omega_{\gbigy}:=
\lambda^{-\dim Y}p_{1,Y}^{\ast}(\omega_Y)(\ast\gbigy^{\infty})$.
We set
$\gbigr_{Y\larr X}:=
 \omega_{\gbigx}\otimes_{F^{-1}(\nbigo_{\gbigy})}
 F^{-1}(\gbigr_Y\otimes\omega_{\gbigy}^{-1})$.
Let $H_Y$ be a closed complex hypersurface of $Y$.
We set $H_X:=F^{-1}(H_Y)$.
For any $\gbigr_{X(\ast H_X)}$-module $\gbigm$,
we set
\[
 F^i_{\dagger}(\gbigm):=
 R^i(\id_{\proj^1}\times F)_{\ast}
 \bigl(
  \gbigr_{Y\larr X}(\ast H_X)\otimes^L_{\gbigr_X(\ast H_X)}
  \gbigm
 \bigr).
\]
If $\gbigm$ is an $\gbigrtilde_{X(\ast H_X)}$-module,
$F^i_{\dagger}(\gbigm)$ is naturally
$\gbigrtilde_{Y(\ast H_Y)}$-module.
Note that the restriction to
$\gbigy^{\circ}$
is equal to the direct image of $\nbigd$-modules.
The following proposition is obvious.

\begin{prop}
\label{prop;21.6.22.15}
If $\gbigm\in\gbigc(X;H_X)$,
then $F_{\dagger}^j(\gbigm)\in\gbigc(Y;H_Y)$.
Thus, we obtain
$F_{\dagger}^j:
\gbigc(X;H_X)\lrarr\gbigc(Y;H_Y)$.
If $F$ induces an isomorphism
$X\setminus H_X\simeq Y\setminus H_Y$,
we have
$F_{\dagger}^j=0$ $(j\neq 0)$,
and $F_{\dagger}^0$ induces an equivalence
 $\gbigc(X;H_X)\simeq
 \gbigc(Y;H_Y)$.
 In this case,
 we have $F_{\dagger}^0(\gbigm)=F_{\ast}(\gbigm)$.
A quasi-inverse $\gbigc(Y;H_Y)\simeq \gbigc(X;H_X)$
is given by the correspondence
$\gbign\longmapsto F^{\ast}(\gbign)$.
\hfill\qed
\end{prop}

Suppose that $F_{|X\setminus H_X}$ is a closed embedding
of $X\setminus H_X$ into $Y\setminus H_Y$.
Let $\gbigc_{F(X)}(Y;H_Y)\subset\gbigc(Y;H_Y)$ denote
the full subcategory of $\gbigm\in\gbigc(Y;H_Y)$
such that the support of $\gbigm$ is contained in
$\proj^1\times F(X)$.
We obtain the following from Proposition \ref{prop;21.4.9.21}
and the Kashiwara equivalence for $\nbigd$-modules.
\begin{prop}
For any $\gbigm\in\gbigc(X;H_X)$,
we have $F_{\dagger}^j(\gbigm)=0$ $(j\neq 0)$.
The functor
 $F_{\dagger}^0:\gbigc(X;H_X)\lrarr
 \gbigc_{F(X)}(Y;H_Y)$  is an equivalence.
It is refined as in Corollary {\rm\ref{cor;22.7.31.2}}.
\hfill\qed 
\end{prop}

\subsubsection{Strict specializability along a coordinate function}
\label{subsection;21.3.13.20}

Let us consider the case where
$X$ is an open subset of $\proj^1\times X_0$
for a complex manifold $X_0$.
We use the notation in \S\ref{subsection;21.4.12.40}.
Let $V\gbigr_{X(\ast H)}\subset \gbigr_{X(\ast H)}$ denote
the sheaf of subalgebras generated by
$\lambda p_{1,X}^{\ast}\Theta_X(\log t)$
over $\nbigo_{\gbigx}\bigl(\ast(\gbigx^{\infty}\cup\gbigh)\bigr)$.
We set
$V\gbigrtilde_{X(\ast H)}=
V\gbigr_{X(\ast H)}\langle\lambda^2\del_{\lambda}\rangle$.

\begin{lem}
\label{lem;21.3.13.21}
For $\gbigm\in\gbigc(X;H)$,
 there exists an increasing filtration
 $V_{\bullet}(\gbigm)$ by
       coherent $V\gbigrtilde_{X(\ast H)}$-submodules indexed by $\real$
       such that
       $V_{\bullet}(\gbigm)_{|\nbigx}$
       is equal to
       the $V$-filtration of $\gbigm_{|\nbigx}$.
 It satisfies the following conditions.
\begin{itemize}
 \item $V_a(\gbigm)(\ast t)=\gbigm(\ast t)$.
 \item For $a\in\real$ and any compact subset $K\subset X$,
       there exists $\epsilon>0$
       such that
       $V_a(\gbigm)_{|\proj^1\times K}=
       V_{a+\epsilon}(\gbigm)_{|\proj^1\times K}$.
 \item $\Gr^V_a(\gbigm):=V_{a}(\gbigm)/V_{<a}(\nbigm)
       \in\gbigc\bigl(X\cap(\{0\}\times X_0)\bigr)$.
       In particular, it is strict,
       i.e., flat over $\nbigo_{\proj^1_{\lambda}}(\ast\infty)$.
 \item $tV_a(\gbigm)\subset V_{a-1}(\gbigm)$ for any $a$.
       If $a<0$,
       then $tV_a(\gbigm)=V_{a-1}(\gbigm)$.
 \item $\deldel_tV_a(\gbigm)\subset V_{a+1}(\gbigm)$ for any $a$.
       If $a>-1$,
       the induced morphism
       $\Gr^V_{a}(\gbigm)\lrarr
       \Gr^V_{a+1}(\gbigm)$
       is an isomorphism of sheaves.
 \item $-\deldel_tt-a\lambda$ is locally nilpotent
       on $\Gr^V_a(\gbigm)$.
\end{itemize}
\end{lem}
\pf
There exists a $V$-filtration
$V'_{\bullet}(\gbigm_{|\gbigx^{\circ}})$
of $\gbigm_{|\gbigx^{\circ}}$
along $t$ as a holonomic $\nbigd_{\gbigx^{\circ}}$-module.
We set
$V_{\bullet}(\gbigm_{|\gbigx^{\circ}}):=
V'_{\bullet}(\gbigm_{|\gbigx^{\circ}})\bigl(\ast(\gbigx^{\infty}\cup\gbigh)\bigr)$.
We obtain a filtration
$V_{\bullet}(\gbigm)$
by gluing
$V_{\bullet}(\gbigm_{|\gbigx^{\circ}})$
and
$V_{\bullet}(\gbigm_{|\nbigx})$.
Then, it satisfies the desired conditions.
\hfill\qed

\begin{rem}
Though $V_a\gbigm_{|\nbigx}$ is coherent over
$V\nbigr_{X(\ast H)}$,
we do not impose the coherence
of $V_a\gbigm$ over $V\gbigr_{X(\ast H)}$.
But, see Proposition {\rm\ref{prop;21.6.26.2}}.
\hfill\qed
\end{rem}

For $\gbigm\in\gbigc(X;H)$,
we set
$\gbigm[\ast t]:=
\gbigr_{X(\ast H)}\otimes_{V\gbigr_{X(\ast H)}}V_0\gbigm$
and
$\gbigm[! t]:=
\bigl(
\gbigr_{X(\ast H)}\otimes_{V\gbigr_{X(\ast H)}}V_{<0}\gbigm
\bigr)$.
They are also objects in $\gbigc(X;H_X)$.

As a variant,
for $\gbigm\in\gbigc\bigl(X;H\cup(\{0\}\times X_0)\bigr)$,
there exists an increasing filtration
$V_{\bullet}(\gbigm)$ by
coherent $V\gbigrtilde_{X(\ast H)}$-submodules indexed by $\real$
such that
$V_{\bullet}(\gbigm)_{|\nbigx}$
is equal to
the $V$-filtration of $\gbigm_{|\nbigx}$.
It satisfies the conditions in Lemma \ref{lem;21.3.13.21}
replacing the fourth and fifth conditions by the following.
\begin{itemize}
 \item $tV_a(\gbigm)=V_{a-1}(\gbigm)$ for any $a\in\real$.
 \item $\deldel_tV_a(\gbigm)\subset V_{a+1}(\gbigm)$ for any $a\in\real$.
\end{itemize}

\subsubsection{Strict specializability, localizability
and Beilinson functors along a function}

Let $X$ be any complex manifold with a closed complex hypersurface $H$.
Let $f$ be a meromorphic function on $(X,H)$.

\begin{lem}
Suppose $|(f)_0|\cap |(f)_{\infty}|=\emptyset$.
For $\gbigm\in\gbigc(X;H)$,
there exist $\gbigm[\star f]\in\gbigc(X;H)$  $(\star=\ast,!)$
such that
 $\iota_{f\dagger}(\gbigm[\star f])
 \simeq
 \iota_{f\dagger}(\gbigm)[\star t]$,
where
$\iota_f:X\lrarr X\times\proj^1$ denotes the graph of $f$.
 \end{lem}
\pf
For $\star=\ast,!$,
we obtain $\gbigrtilde_X$-modules $\gbigm[\star f]$
by gluing 
the holonomic $\nbigd_{\gbigx^{\circ}}$-module
$\Bigl(
 (\gbigm_{|\gbigx^{\circ}})(\star f)
 \Bigr)\bigl(\ast(\gbigx^{\infty}\cup\gbigh^{\circ})\bigr)$
 and 
$\gbigm_{|\nbigx}[\star f]\in\nbigc(X;H)$.
They satisfy the desired condition.
\hfill\qed

\vspace{.1in}
Let us consider the case $|(f)_0|\cap |(f)_{\infty}|\neq \emptyset$.
There exists a projective morphism
$\rho:X'\lrarr X$
such that
$|(\rho^{\ast}(f))_{0}|\cap|(\rho^{\ast}(f))_{\infty}|=\emptyset$,
and that $\rho$ induces
$X'\setminus \rho^{-1}(H)\simeq X\setminus H$.
We put $H':=\rho^{-1}(H)$.
We set $\gbigm':=\rho^{\ast}(\gbigm)$.
It is easy to see that
$\gbigm'$ is an object of $\gbigc(X';H')$,
and that $\rho_{\ast}(\gbigm')\simeq\gbigm$.
We set
$\gbigm[\star f]:=
 \rho_{\ast}\bigl(\gbigm'[\star \rho^{\ast}(f)]\bigr)
 \in\gbigc(X;H)$.

\vspace{.1in}

For $a<b$,
we set
\[
 \II^{a,b}_{\gbigx,f}:=\bigoplus_{a\leq j<b}
 \nbigo_{\gbigx}(\ast\gbigx^{\infty})\bigl(\ast ((f)_{0}\cup\gbigh)\bigr)(\lambda s)^j.
\]
It is naturally an $\gbigrtilde_{X(\ast H)}(\ast f)$-module
with the meromorphic flat connection $\nabla$
defined by $\nabla(s^j)=s^{j+1}df/f$.

Let $\gbigm\in\gbigc(X;H)$.
For $a<b$,
we obtain $\Pi^{a,b}_{f,\star}(\gbigm)\in\gbigc(X;H)$
by gluing
$\Pi^{a,b}_{f,\star}(\gbigm_{|\nbigx})$
and
\[ 
\Bigl(
(\II^{a,b}_{\gbigx,f}\otimes\gbigm)_{|\gbigx^{\circ}}(\star f)
\Bigr)\bigl(\ast (\gbigx^{\infty}\cup\gbigh^{\circ})\bigr).
\]
The Beilinson functor
$\Pi^{a,b}_{f,\ast!}(\gbigm)\in\gbigc(X;H_X)$
is defined as
\[
\Pi^{a,b}_{f,\ast!}(\gbigm):=
 \varprojlim_{N\to\infty}
 \Cok\Bigl(\Pi^{b,N}_{f!}(\gbigm)
 \lrarr
 \Pi^{a,N}_{f\ast}(\gbigm)\Bigr).
\]
In particular,
we set
$\psi^{(a)}_f(\gbigm):=
 \Pi^{a,a}_{f,\ast !}(\gbigm)$
and
$\Xi^{(a)}_f(\gbigm):=
 \Pi^{a,a+1}_{f,\ast !}(\gbigm)$.
There exists the following naturally defined complex
of $\gbigrtilde_X(\ast H)$-modules:
\begin{equation}
\label{eq;21.2.16.60}
 \gbigm[!f]
 \lrarr
 \Xi^{(0)}_f(\gbigm)
 \oplus
 \gbigm
 \lrarr
 \gbigm[\ast f].
\end{equation}
We define $\phi^{(0)}_f(\gbigm)$
as the cohomology of the complex (\ref{eq;21.2.16.60}).
It is easy to see that
$\phi^{(0)}_f(\gbigm)\in\gbigc(X;H)$.
There exist the natural morphisms
$\psi^{(1)}_f(\gbigm)\lrarr\phi^{(0)}_f(\gbigm)
\lrarr\psi^{(0)}_f(\gbigm)$.
Then, as in \cite{beilinson2},
$\gbigm$ is naturally isomorphic to
the cohomology of the following complex:
\begin{equation}
 \psi^{(1)}_f(\gbigm)
  \lrarr
  \phi^{(0)}_f(\gbigm)
  \oplus
  \Xi^{(0)}_f(\gbigm)
  \lrarr
  \psi^{(0)}_f(\gbigm).
\end{equation}

\subsubsection{Localizability along hypersurfaces}

Let $H^{(1)}$ be a hypersurface of $X$.
Let $\gbigm\in\gbigc(X;H)$.
We obtain
$\gbigm[\star H^{(1)}]\in\gbigc(X;H)$ $(\star =!,\ast)$
by gluing
$\bigl(
\gbigm_{|\gbigx^{\circ}}(\star H^{(1)})
\bigr)\bigl(\ast(\gbigx^{\infty}\cup \gbigh^{(1)\circ})\bigr)$
and
$\gbigm_{|\nbigx}[\star H^{(1)}]$.
They satisfy the following condition.
\begin{itemize}
 \item Let $U$ be any open subset of $X$.
       Let $f$ be a meromorphic function on
       $(U,H\cap U)$
       such that
       $|(f)_0|\cup (U\cap H)=(U\cap H^{(1)})\cup (U\cap H)$.
       Then,
       $\gbigm_{|U}[\star f]=(\gbigm[\star H^{(1)}])_{|U}$.
\end{itemize}

\subsubsection{Some compatibility}

Let $F:X\lrarr Y$ be a projective morphism of complex manifolds.
Let $H_Y$ be a closed complex hypersurface of $Y$.
We set $H_X:=F^{-1}(H_Y)$.
Let $f_Y$ be any meromorphic function on $(Y,H_Y)$.
We set $f_X:=F^{\ast}(f_Y)$.
We easily obtain the following proposition
from the property of mixed twistor $\nbigd$-modules
and holonomic $\nbigd$-modules.
\begin{prop}  
Let $\gbigm\in\gbigc(X;H_X)$.
\begin{itemize}
 \item
If $|(f_Y)_0|\cap|(f_Y)_{\infty}|=\emptyset$,
there exists a natural isomorphism:
\[
 V_a\bigl(\iota_{f_Y\dagger}F_{\dagger}^j(\gbigm)\bigr)
\simeq
 R^j(\id_{\proj^1}\times F)_{\ast}
 \Bigl(
 \pi^{\ast}\gbigr_{Y\rarr X}\otimes^L_{\pi^{\ast}\gbigr_X}
 V_a\bigl(
  \iota_{f_X\dagger}(\gbigm)
 \bigr)
 \Bigr).
\]
Here,
$\iota_f$ denotes the graph embedding of $f$,
and 
$\pi:\proj^1_{\lambda}\times \proj^1_t\times X
\lrarr \proj^1_{\lambda}\times X$ denotes the projection.
\item
 There exist natural isomorphisms
 $F_{\dagger}^j(\gbigm[\star f_X])
 \simeq
 F_{\dagger}^j(\gbigm)[\star f_Y]$.
 \item
      There exist natural isomorphisms
      $\Pi^{a,b}_{f_Y,\ast !}(F_{\dagger}^j(\gbigm))
      \simeq
      F_{\dagger}^j\bigl(
      \Pi^{a,b}_{f_X,\ast !}(\gbigm)\bigr)$.
      In particular,
      there exist natural isomorphisms
      $\Xi^{(a)}_{f_Y}(F_{\dagger}^j(\gbigm))
      \simeq
      F_{\dagger}^j\bigl(
      \Xi^{(a)}_{f_X}(\gbigm)\bigr)$
      and 
      $\psi^{(a)}_{f_Y}(F_{\dagger}^j(\gbigm))
      \simeq
      F_{\dagger}^j\bigl(
      \psi^{(a)}_{f_X}(\gbigm)\bigr)$.
      Moreover, we have
      $\phi^{(0)}_{f_Y}(F_{\dagger}^j(\gbigm))
      \simeq
      F_{\dagger}^j\bigl(
      \phi^{(0)}_{f_X}(\gbigm)\bigr)$.
      \hfill\qed
\end{itemize}
\end{prop}

Let $H^{(1)}_Y$ be a complex analytic hypersurface of $Y$.
We obtain $H^{(1)}_X:=F^{-1}(H^{(1)}_Y)$.
\begin{prop}
 For $\gbigm\in\gbigc(X;H_X)$,
 there exist natural isomorphisms
 $F_{\dagger}^j(\nbigm)[\star H^{(1)}_Y]
 \simeq
 F_{\dagger}^j(\nbigm[\star H^{(1)}_X])$
$(\star=!,\ast)$.
\hfill\qed
\end{prop}

\subsubsection{External tensor product}

Let $X_i$ $(i=1,2)$ be complex manifolds
equipped with a closed complex hypersurface $H_i$.
We set
$\Htilde:=(H_1\times X_2)\cup(X_1\times H_2)$.
Let $\gbigm_i\in\gbigc(X_i;H_i)$ $(i=1,2)$.
We obtain an $\gbigrtilde_{X_1\times X_2}(\ast \Htilde)$-module
$\gbigm_1\boxtimes\gbigm_2$
as the gluing of
$\gbigm_{1|\nbigx_1}\boxtimes
 \gbigm_{2|\nbigx_2}\in\nbigc(X_1\times X_2;\Htilde)$
and
the $\nbigd_{(\proj^1\setminus\{0\})\times (X_1\times X_2)}
\bigl(
\ast(\{\infty\}\times(X_1\times X_2))
\bigr)$-module
\[
 p_{1,\gbigx_1^{\circ}}^{\ast}\bigl(
  \gbigm_{1|\gbigx_1^{\circ}}
  \bigr)
  \otimes_{\nbigo_{(\proj^1\setminus\{0\})\times X_1\times X_2}}
  p_{2,\gbigx_2^{\circ}}^{\ast}\bigl(
  \gbigm_{2|\gbigx_2^{\circ}}
  \bigr).
\]
Here,
$p_{i,\gbigx_i^{\circ}}$ denote the projections
$(\proj^1\setminus\{0\})\times (X_1\times X_2)\lrarr\gbigx_i^{\circ}$.
Thus, we obtain
$\gbigm_1\boxtimes\gbigm_2\in\gbigc(X_1\times X_2;\Htilde)$.
 
\subsubsection{Duality}

Let $\gbigm\in\gbigc(X;H)$.
Note that 
$\DD_{X(\ast H)}(\gbigm_{|\nbigx})_{|\nbigx^{\circ}}$
and
$\DDD_{\gbigx^{\circ}}(\gbigm_{|\gbigx^{\circ}})(\ast\gbigh^{\circ})_{|\nbigx^{\circ}}$
are naturally isomorphic,
according to Proposition \ref{prop;21.3.23.10}.
Hence, we obtain
$\DD_{X(\ast H)}(\gbigm)\in\gbigc(X;H)$
by gluing 
$\DD_{X(\ast H)}(\gbigm_{|\nbigx})$
and 
$\DDD_{\gbigx^{\circ}}(\gbigm_{|\gbigx^{\circ}})
\bigl(\ast(\gbigx^{\infty}\cup\gbigh^{\circ})\bigr)$.

We set
$\omega_{\gbigx}:=\lambda^{-d_X}p_{1,X}^{\ast}(\omega_X)(\ast\gbigx^{\infty})$,
where $d_X=\dim X$.
By Proposition \ref{prop;21.3.23.10},
there exists the following natural isomorphism
for any $\gbigm\in\gbigc(X;H)$:
\[
 \DD_{X(\ast H)}(\gbigm)\simeq
 \lambda^{d_X}
 \nrhom_{\gbigr_{X(\ast H)}}\Bigl(
 \gbigm,
 \gbigr_{X(\ast H)}\otimes
 \omega^{-1}_{\gbigx}
 \Bigr)[d_X].
\]

\subsubsection{Non-characteristic inverse image}

Let $f:X\lrarr Y$ be a morphism of complex manifolds.
We put
$\gbigr_{X\rarr Y}:=
 \nbigo_{\gbigx}(\ast\gbigx^{\infty})
 \otimes_{f^{-1}\nbigo_{\gbigy}(\ast\gbigy^{\infty})}
 \gbigr_Y$.
Let $H_Y$ be a closed complex hypersurface of $Y$.
We set $H_X:=f^{-1}(H_Y)$.
For $\gbigrtilde_Y(\ast H_Y)$-module $\gbigm$,
we set
$Lf^{\ast}\gbigm=
\gbigr_{X\rarr Y}(\ast H_X)\otimes^L_{f^{-1}\gbigr_Y(\ast H_Y)}
f^{-1}(\gbigm)$.

\begin{prop}
Let $\gbigm\in\gbigc(Y;H_Y)$.
If $f$ is strictly non-characteristic for $\gbigm_{|\nbigy}$
in the sense of Definition  {\rm\ref{df;21.4.12.30}},
we have
\[
 Lf^{\ast}(\gbigm)
 =f^{\ast}(\gbigm)
=\gbigr_{X\rarr Y}\otimes_{f^{-1}\gbigr_Y}f^{-1}(\gbigm)
 =\nbigo_{\gbigx}(\ast\gbigx^{\infty})
 \otimes_{f^{-1}\nbigo_{\gbigy}(\ast\gbigy^{\infty})}
 f^{-1}(\gbigm)
 \in\gbigc(X;H_X).
\]
Moreover, there exists a natural isomorphism
 $\lambda^{-\dim X}\DD_{X(\ast H_X)}f^{\ast}\gbigm
 \simeq
 \lambda^{-\dim Y}f^{\ast}\DD_{Y(\ast H_Y)}\gbigm$. 
\end{prop}
\pf
For any $k$,
$L^kf^{\ast}(\gbigm)$
are good $\nbigo_{\gbigx}$-modules,
and
$L^kf^{\ast}(\gbigm)(\ast\gbigx^{\infty})
=L^kf^{\ast}(\gbigm)$.
Because
$L^kf^{\ast}(\gbigm)_{|\nbigx}=0$ $(k\neq 0)$
we obtain
$L^kf^{\ast}(\gbigm)=0$ $(k\neq 0)$.
By Theorem \ref{thm;21.4.7.50},
we obtain
$f^{\ast}(\gbigm)\in\gbigc(X)$.
It is easy to see that
the isomorphism
 $\lambda^{-\dim X}\DD_{X(\ast H_X)}f^{\ast}(\nbigm_{|\nbigy})
 \simeq
 \lambda^{-\dim Y}f^{\ast}\DD_{Y(\ast H_Y)}(\nbigm_{|\nbigy})$
 extends to 
 $\lambda^{-\dim X}\DD_{X(\ast H_X)}f^{\ast}(\nbigm)
 \simeq
 \lambda^{-\dim Y}f^{\ast}\DD_{Y(\ast H_Y)}(\nbigm)$.
\hfill\qed

\subsubsection{Inverse image by closed embeddings under an additional condition}
\label{subsection;21.6.22.20}

We use the notation in \S\ref{subsection;21.6.22.2}.
Let $\gbigc(Y,[\star H_Y^{(1)}])\subset\gbigc(Y)$
denote the full subcategory of
objects $\gbigm\in\gbigc(Y)$
such that $\gbigm\simeq\gbigm[\star H_Y^{(1)}]$.
Let $\gbigc(Y,[\star H_Y^{(1)}];H_Y^{(2)})$
denote the essential image of
$\gbigc(Y,[\star H_Y^{(1)}])
\lrarr \gbigc(Y;H_Y^{(2)})$.
For $\gbigm\in\gbigc(Y,[\star H_Y^{(1)}];H^{(2)}_Y)$,
we obtain the complex of
$\gbigrtilde_{Y(\ast H^{(2)}_Y)}$-modules
$C_{\star}(\gbigm,\Gamma,\gbigk)$
in the same way as in \S\ref{subsection;21.4.14.1}.
We have
$\nbigh^i
C_{\star}(\gbigm,\Gamma,\gbigk)
\in \gbigc(Y,[\star H^{(1)}_Y];H_Y^{(2)})$,
and 
the support of
$\nbigh^iC_{\star}(\gbigm,\Gamma,\gbigk)$
are contained in $\proj^1_{\lambda}\times F(X)$.

For $\gbigm\in\gbigc(Y,[!H_Y^{(1)}];H^{(2)}_Y)$,
there exist
$(\lefttop{T}f^{\ast})^i(\gbigm)
\in\gbigc(X,[!H_X^{(1)}];H^{(2)}_X)$
with an isomorphism
\[
 f_{\dagger}\bigl(
 (\lefttop{T}f^{\ast})^i(\gbigm)
 \bigr)
 \simeq
 \nbigh^i
  C_!(\gbigm,\Gamma,\gbigk).
\]
For $\gbigm\in\gbigc(Y,[\ast H_Y^{(1)}];H^{(2)}_Y)$,
there exist
$(\lefttop{T}f^{!})^i(\gbigm)
\in\gbigc(X,[\ast H_X^{(1)}];H_X^{(2)})$
with an isomorphism
\[
 f_{\dagger}\bigl(
 (\lefttop{T}f^{!})^i(\gbigm)
 \bigr)
 \simeq
 \nbigh^i
  C_{\ast}(\gbigm,\Gamma,\gbigk).
\]
They are independent of the choice of $(\Gamma,\gbigk)$.
We have similar formulas
as in Lemma \ref{lem;21.6.22.10}
and Corollary \ref{cor;21.4.21.1}.

\subsubsection{Tensor product}
\label{subsection;21.6.22.18}

We use the notation in \S\ref{subsection;21.6.22.11}.
For $\gbigm_i\in\gbigc(X;H)$ $(i=1,2)$
and for $\star=!,\ast$,
by using
$\gbigm_1\boxtimes\gbigm_2\in\gbigc(\Xtilde;\Htilde)$,
we define
\[
 \nbigh^k\bigl(
 \gbigm_1\otimes^{\star}\gbigm_2
 \bigr):=
 (\lefttop{T}\Delta^{\star})^k
 \bigl(
 \gbigm_1\boxtimes\gbigm_2
 \bigr)
 \in\gbigc(X;H).
\]
We have similar formulas as in
Proposition \ref{prop;21.6.22.12}
and Corollary \ref{cor;21.6.22.13}.

\section{Malgrange extension of integrable mixed twistor $\nbigd$-modules}

\subsection{Existence and uniqueness of the Malgrange extension}
\label{subsection;21.4.14.30}

Let $X$ be a complex manifold.
Let $\iota:\nbigx\lrarr\gbigx$ denote the inclusion.
There exists a natural inclusion
$\gbigrtilde_{X}\subset\iota_{\ast}\nbigrtilde_{X}$.
Any $\nbigrtilde_{X}$-module $\nbigm$ induces
an $\gbigrtilde_{X}$-module $\iota_{\ast}\nbigm$
though it is too large.

\begin{df}
Let $\nbigm$ be an $\nbigrtilde_{X}$-module
such that 
$\nbigm_{|\nbigx\setminus\nbigx^0}$
is a holonomic $\nbigd_{\nbigx\setminus\nbigx^0}$-module.
An $\gbigrtilde_X$-submodule $\gbigm\subset\iota_{\ast}(\nbigm)$
is called a Malgrange extension of $\nbigm$
if the following holds.
\begin{itemize}
 \item $\iota^{-1}\gbigm=\nbigm$.
 \item $\gbigm_{|\gbigx\setminus\gbigx^0}$ is a holonomic
$\nbigd_{\gbigx\setminus\gbigx^0}$-module
which is strongly regular along $\mu=\lambda^{-1}$. 
\hfill\qed
\end{itemize}
\end{df}

We shall prove the following theorem
in \S\ref{subsection;21.4.14.10}
after the preliminary in \S\ref{subsection;21.4.14.40}.
\begin{thm}
\label{thm;21.1.23.1}
For any $\nbigm\in\nbigc(X)$,
there uniquely exists a Malgrange extension
$\Upsilon(\nbigm)\in\gbigc(X)$.
This extension procedure induces a fully faithful functor
$\Upsilon:\nbigc(X)\lrarr \gbigc(X)$.
\end{thm}

\begin{rem}
In {\rm\cite[Theorem 14.4.8]{Mochizuki-MTM}},
we should impose the strong regularity condition
instead of the regularity condition.
\hfill\qed
\end{rem}

\begin{notation}
Let $\gbigc_{\Malg}(X)$ denote the essential image of
$\Upsilon:\nbigc(X)\lrarr \gbigc(X)$.
\hfill\qed
\end{notation}

We shall also observe the following separation property
in \S\ref{subsection;21.4.16.1}
as a complement of Theorem \ref{thm;21.1.23.1}.
\begin{prop}
\label{prop;21.4.16.2}
For the $V$-filtration $\lefttop{\infty}V_{\bullet}\Upsilon(\nbigm)$
of $\Upsilon(\nbigm)$ along $\mu$,
we have
$\bigcap_{a\in\real}(\lefttop{\infty}V_a\Upsilon(\nbigm))=0$. 
\end{prop}

We shall also prove the following complement
in \S\ref{subsection;22.7.27.1}.

\begin{prop}
\label{prop;21.6.26.2}
Let $g$ be any holomorphic function on $X$.
Let $\iota_g:X\to X\times\cnum_t$ denote the induced embedding.
Let $\nbigm\in\nbigc(X)$.
Let 
$V_{\bullet}\Bigl(
\iota_{g\dagger}\Upsilon(\nbigm)\Bigr)$
denote the $V$-filtration of
$\iota_{g\dagger}\Upsilon(\nbigm)$ along $t$
in the sense of Lemma {\rm\ref{lem;21.3.13.21}}.
Then, 
$V_{a}\Bigl(
\iota_{g\dagger}\Upsilon(\nbigm)\Bigr)$
are coherent over $V\gbigr_{X\times\cnum_t}$
not only over $V\gbigrtilde_{X\times\cnum_t}$.
\end{prop}

\subsection{Malgrange extension of $\nbigrtilde$-modules
 in the good-KMS case}
\label{subsection;21.4.14.40}
Let $X$ be a complex manifold
with a simple normal crossing hypersurface
$H=\bigcup_{i\in \Lambda}H_i$.
Let $\nbigv$ be a good-KMS smooth $\nbigr_{X(\ast H)}$-module.
(See \cite[\S5.1.1]{Mochizuki-MTM}.)
Moreover, we assume that
the $\nbigr_{X(\ast H)}$-action on $\nbigv$
extends to an $\nbigrtilde_{X(\ast H)}$-module on $\nbigv$.
In this case,
the eigenvalues of $\Res(\nablatilde)$ along each $H_i$
are real numbers,
and that
there exists a global filtered bundle
$\nbigq_{\ast}\nbigv
=\bigl(\nbigq_{\veca}\nbigv\,\big|\,\veca\in\real^{\Lambda}\bigr)$
over $\nbigv$
whose restriction to $\nbigx\setminus\nbigx^0$
is equal to the Deligne-Malgrange filtered bundle
of the meromorphic flat bundle $\nbigv_{|\nbigx\setminus\nbigx^0}$.

By Proposition \ref{prop;21.1.23.10},
$\nbigv_{|\nbigx\setminus\nbigx^0}$
extends to a meromorphic flat bundle
$\widetilde{(\nbigv_{|\nbigx\setminus\nbigx^0})}$
on $(\gbigx\setminus\gbigx^0,(\gbigh\cup\gbigx^{\infty})\setminus\gbigh^0)$,
which is strongly regular along $\mu=\lambda^{-1}$.
By gluing
$\nbigv$ and $\widetilde{(\nbigv_{|\nbigx\setminus\nbigx^0})}$,
we obtain an $\gbigrtilde_{X(\ast H)}$-module
$\gbigv$, which is a Malgrange extension of $\nbigv$.

Recall that
we constructed the holonomic $\nbigr_X$-module
$\nbigv[!I\ast J]$ for any decomposition $\Lambda=I\sqcup J$,
which is naturally an $\nbigrtilde_X$-module.
(See \cite[\S5.3]{Mochizuki-MTM}.)

\begin{prop}
There uniquely exists
a Malgrange extension
$\Upsilon(\nbigv[!I\ast J])$
of $\nbigv[!I\ast J]$.
\end{prop}
\pf
It is enough to prove the claim locally around
any point of $H$.
Hence, we may assume that
$X$ is a neighbourhood of
$(0,\ldots,0)$ in $\cnum^n$,
and $H=\bigcup_{i=1}^{\ell}\{z_i=0\}$.

For any $e\in\seisuu_{>0}$,
let $\varphi_e:\cnum^n\lrarr \cnum^n$
be the map defined by
\[
 \varphi_e(z_{e,1},\ldots,z_{e,\ell},z_{\ell+1},\ldots,z_{n})
=(z_{e,1}^e,\ldots,z_{e,\ell}^e,z_{\ell+1},\ldots,z_n).
\]
We set $X^{(e)}:=\varphi_e^{-1}(X)$
and $H^{(e)}:=\varphi_e^{-1}(H)$.
For an appropriate choice of $e$,
$\nbigv^{(e)}:=\varphi_e^{\ast}(\nbigv)$
is unramifiedly good-KMS.
By using Proposition \ref{prop;21.1.23.10},
we uniquely extend $\nbigv^{(e)}$
to a $\gbigrtilde_{\gbigx^{(e)}}(\ast\gbigh^{(e)})$-module
$\nbigvtilde^{(e)}$
such that
$\nbigvtilde^{(e)}_{|\gbigx^{(e)}\setminus(\gbigx^{(e)})^0}$
is a meromorphic flat bundle
on
$\bigl(\gbigx^{(e)}\setminus(\gbigx^{(e)})^0,
 (\gbigx^{(e)})^{\infty}\cup
 (\gbigh^{(e)}\setminus(\gbigh^{(e)})^0)
 \bigr)$
such that
 $\nbigvtilde^{(e)}_{|\gbigx^{(e)}\setminus
 ((\gbigx^{(e)})^0\cup\gbigh^{(e)})}$
is regular singular.
By the uniqueness,
it is naturally equivariant
with respect to the action of the Galois group
of the ramified covering $\varphi_e$.
Let $\nbigvtilde$ denote the descent of
$\nbigvtilde^{(e)}$.
The natural inclusion $\nbigvtilde\subset
\varphi_{e\ast}(\nbigvtilde^{(e)})$ induces
\[
\bigl(
 \nbigvtilde_{|\gbigx\setminus\gbigx^{0}}
 \bigr)
 \bigl(!\gbigh(I)\ast(\gbigh(J)\cup\gbigx^{\infty})
 \bigr)
\subset
\varphi_{e\dagger}\Bigl(
\bigl(
\nbigvtilde^{(e)}_{|\gbigx^{(e)}\setminus(\gbigx^{(e)})^0}
\bigr)
\bigl(!\gbigh^{(e)}(I)\ast(\gbigh^{(e)}(J)\cup(\gbigx^{(e)})^{\infty})
 \bigr)
\Bigr).
\]
Hence, by Proposition \ref{prop;21.1.23.11},
Lemma \ref{lem;21.3.9.33}
and Proposition \ref{prop;21.1.23.2},
we obtain a Malgrange extension of
$\nbigv[!I\ast J]$
by gluing $\nbigv[!I\ast J]$
and 
$\bigl(
 \nbigvtilde_{|\gbigx\setminus\gbigx^{0}}
 \bigr)
 \bigl(!\gbigh(I)\ast(\gbigh(J)\cup\gbigx^{\infty})
 \bigr)$.

Let $\Upsilon'(\nbigv[!I\ast J])$
be another Malgrange extension.
Note that
 $\Upsilon'(\nbigv[!I\ast J])_{|\gbigx\setminus(\gbigh\cup\gbigx^0)}$
is regular singular. 
Hence, there exists a natural isomorphism
\[
 \Upsilon'(\nbigv[!I\ast J])_{|\gbigx\setminus\gbigh^{\infty}}
 \simeq
 \Upsilon(\nbigv[!I\ast J])_{|\gbigx\setminus\gbigh^{\infty}}.
\]
It induces the following isomorphism
\[
\Bigl(
 \Upsilon'(\nbigv[!I\ast J])
 (\ast\gbigh)
 \Bigr)
 _{|\gbigx\setminus\gbigh^{\infty}}
 \simeq
 \Bigl(
 \Upsilon(\nbigv[!I\ast J])
 (\ast\gbigh)
 \Bigr)_{|\gbigx\setminus\gbigh^{\infty}}.
\]
By the Hartogs theorem,
it extends to the isomorphism
\begin{equation}
\label{eq;21.2.12.10}
 \Bigl(
 \Upsilon'(\nbigv[!I\ast J])
 (\ast\gbigh)
 \Bigr)
 \simeq
 \Bigl(
 \Upsilon(\nbigv[!I\ast J])
 (\ast\gbigh)
 \Bigr).
\end{equation}
Because
$\Upsilon'(\nbigv[!I\ast J])_{|\gbigx\setminus\gbigx^0}
 =\Bigl(
 \bigl(\Upsilon'(\ast \gbigh)\bigr)
 (!\gbigh(I))
 \Bigr)(\ast\gbigx^{\infty})$,
 (\ref{eq;21.2.12.10})
induces an isomorphism
$\Upsilon'(\nbigv[!I\ast J])
\simeq
\Upsilon(\nbigv[!I\ast J])$.
Thus, we obtain the uniqueness.
\hfill\qed

\vspace{.1in}
We easily obtain the following lemma
by using Proposition \ref{prop;21.2.17.3}.
 
\begin{lem}
Let $\nbigv_i$ $(i=1,2)$
be good-KMS $\nbigrtilde_{X(\ast H)}$-module.
Any morphism  $f:\nbigv_1[!I\ast J]\lrarr\nbigv_2[!I\ast J]$ of
$\nbigrtilde_{X(\ast H)}$-modules
uniquely extends to
a morphism
 $\Upsilon(f):
 \Upsilon(\nbigv_1[!I\ast J])
 \lrarr
 \Upsilon(\nbigv_2[!I\ast J])$.
\hfill\qed
\end{lem}

\subsubsection{Duality}

For a smooth good-KMS $\nbigrtilde_{X(\ast H)}$-module $\nbigv$,
we set
$\nbigv^{\lor}:=
\nhom_{\nbigo_{\nbigx(\ast\nbigh)}}\bigl(
\nbigv,\nbigo_{\nbigx}(\ast\nbigh)
\bigr)$,
which is naturally a smooth good-KMS
$\nbigrtilde_{X(\ast H)}$-module.
There exists a natural isomorphism
$\DD_{X}(\nbigv)(\ast H)
\simeq \lambda^{d_X}\nbigv^{\lor}$.

\begin{lem}
\label{lem;21.3.23.20}
 We have
 $\DD\Bigl(
 \Upsilon(\nbigv[!I\ast J])
 \Bigr)
 \simeq
 \lambda^{d_X}\Upsilon(\nbigv^{\lor}[!J\ast I])$.
\end{lem}
\pf
Note that
 $\DD\Bigl(
 \nbigv[!I\ast J]
 \Bigr)
 \simeq
 \lambda^{d_X}\nbigv^{\lor}[!J\ast I]$.
By using the Hartogs theorem,
we obtain 
$\DD\Bigl(
 \Upsilon(\nbigv[!I\ast J])
 \Bigr)(\ast H)
 \simeq
 \lambda^{d_X}\Upsilon(\nbigv^{\lor})$.
The isomorphism 
  $\DD\Bigl(
 \nbigv[!I\ast J]
 \Bigr)
 \simeq
 \lambda^{d_X}\nbigv^{\lor}[!J\ast I]$
extends to 
  $\DD\Bigl(
 \Upsilon(\nbigv[!I\ast J])
 \Bigr)
 \simeq
 \lambda^{d_X}\Upsilon(\nbigv^{\lor}[!J\ast I])$.
\hfill\qed

\subsubsection{External tensor product}

Let $X^{(i)}$ $(i=1,2)$ be complex manifolds
with simple normal crossing hypersurfaces $H^{(i)}$.
Let $\nbigv_i$ be smooth good-KMS $\nbigrtilde_{X_i(\ast H_i)}$-modules.
Suppose that $\nbigv_i$ underlie
graded polarizable admissible
integrable variations of mixed twistor structure on $(X_i,H_i)$.
Let $H^{(i)}=\bigcup_{j\in\Lambda^{(i)}}H^{(i)}_j$
be the irreducible decomposition.
Let $I^{(i)}\sqcup J^{(i)}=\Lambda^{(i)}$ be decompositions.

\begin{lem}
\label{lem;21.3.23.11}
 There exists a unique Malgrange extension of 
 $\nbigv_1[!I^{(1)}\ast J^{(1)}]
 \boxtimes
  \nbigv_2[!I^{(2)}\ast J^{(2)}]$,
which is equal to
 $\Upsilon(\nbigv_1[!I^{(1)}\ast J^{(1)}])
 \boxtimes
 \Upsilon(\nbigv_2[!I^{(2)}\ast J^{(2)}])$.
\end{lem}
\pf
It is enough to check the claim locally around any point of
$\Xtilde:=X^{(1)}\times X^{(2)}$.
We set $\Htilde:=(X^{(1)}\times H^{(2)})\cup (H^{(1)}\times X^{(2)})$
and $\nbigvtilde=\nbigv_1\boxtimes\nbigv_2$.
There exists a projective morphism
$F:\Xtilde'\lrarr \Xtilde$
such that
(i) $\Htilde'=F^{-1}(\Htilde)$ is simply normal crossing,
(ii) $\Xtilde'\setminus\Htilde'\simeq \Xtilde\setminus \Htilde$,
(iii) $\nbigvtilde'=F^{\ast}(\nbigvtilde)$ is smooth good-KMS.

By setting $\Lambdatilde=\Lambda^{(1)}\sqcup\Lambda^{(2)}$,
we obtain the induced decomposition
$\Htilde=\bigcup_{k\in\Lambdatilde}\Htilde_k$,
where
$\Htilde_k=H^{(1)}\times X^{(2)}$ if $k\in\Lambda^{(1)}$,
and
$\Htilde_k=X^{(1)}\times H^{(2)}$ if $k\in\Lambda^{(2)}$.
We set $\Itilde:=I^{(1)}\sqcup I^{(2)}$
and $\Jtilde:=J^{(1)}\sqcup J^{(2)}$.
There exists the decomposition
$\Htilde'=\bigcup_{i\in\Lambdatilde'}\Htilde'_i$.
Let $\Itilde'\subset\Lambdatilde'$
be determined by
$\bigcup_{i\in \Itilde}\varphi^{-1}(\Htilde_i)
=\bigcup_{j\in\Itilde'}\Htilde'_j$.
We set $\Jtilde':=\Lambdatilde'\setminus\Itilde'$.
We obtain
$\Upsilon(\nbigvtilde'[!\Itilde'\ast \Jtilde'])$.
Note that
$\varphi_{\dagger}\bigl(
\nbigvtilde'[!\Itilde'\ast \Jtilde']
\bigr)
\simeq
\nbigvtilde[!\Itilde\ast \Jtilde]$
as in the proof of
\cite[Proposition 11.4.6]{Mochizuki-MTM}.
Hence, we obtain a Malgrange extension
$\varphi_{\dagger}
\Upsilon(\nbigvtilde'[!\Itilde'\ast \Jtilde'])$
of
$\nbigv_1[!I_1\ast J_1]\boxtimes
\nbigv_2[!I_2\ast J_2]$.
Then, it is standard to obtain that
$\varphi_{\dagger}
\Upsilon(\nbigvtilde'[!\Itilde'\ast \Jtilde'])$
is equal to
 $\Upsilon(\nbigv_1[!I^{(1)}\ast J^{(1)}])
 \boxtimes
 \Upsilon(\nbigv_2[!I^{(2)}\ast J^{(2)}])$.
\hfill\qed

\subsection{Proof of Theorem \ref{thm;21.1.23.1}}
\label{subsection;21.4.14.10}

Let $\nbigm\in\nbigc(X)$.
Let $P$ be any point of $X$.
There exist a neighbourhood $U_P$ of $P$ in $X$,
and a holomorphic function $f_P$ on $U_P$
such that
(i) $(\Supp(\nbigm)\cap U_P)\setminus f_P^{-1}(0)$
is a closed submanifold of
$U_P\setminus f_P^{-1}(0)$,
(ii)
$\dim(\Supp(\nbigm)\cap f_P^{-1}(0))
<\dim(\Supp(\nbigm)\cap U_P\setminus f_P^{-1}(0))$.
We set $\nbigm_P:=\nbigm_{|U_P}$.

\begin{lem}
\label{lem;21.2.17.1}
 There exists a Malgrange extension
$\Upsilon(\nbigm_P)\in\gbigc(U_P)$ of $\nbigm_P$. 
\end{lem}
\pf
Let $Z_P$ denote the closure of
$(\Supp(\nbigm)\cap U_P)\setminus f_P^{-1}(0)$
in $U_P$.
According to
\cite[Proposition 11.1.1, Proposition 11.1.2, Lemma 11.1.3]{Mochizuki-wild},
there exists a tuple as follows.
\begin{itemize}
 \item A smooth complex manifold $Z'_P$
       with a projective morphism
       $\varphi_P:Z_P'\lrarr U_P$
       such that
       (i) $H_P':=\varphi_P^{\ast}(f_P)^{-1}(0)$ is
       a simple normal crossing hypersurface,
       (ii) $\varphi_P$ induces
       $Z_P'\setminus H_P'\simeq Z_P\setminus f_P^{-1}(0)$.
 \item A good-KMS $\nbigrtilde_{Z'_P(\ast H_P')}$-module $\nbigv_P$
       with an isomorphism of $\nbigrtilde_{U_P}(\ast f_P)$-modules
       $\varphi_{P\dagger}(\nbigv_P)\simeq
       \nbigm_{P}(\ast f_P)$.
       Note that $\nbigv_P$ underlies a good-KMS variation of
       mixed twistor structure.
\end{itemize}

We obtain
$\Pi^{a,b}_{\varphi_P^{\ast}(f),\star}(\nbigv_P)
\in\nbigc(Z_P')$ $(\star=\ast,!)$
for any $a<b$.
By Proposition \ref{prop;21.1.23.11},
there uniquely exist the Malgrange extensions
$\Upsilon(\Pi^{a,b}_{\varphi_P^{\ast}(f),\star}(\nbigv_P))
\in\gbigc(Z_P')$.
We obtain
\[
 \Upsilon\bigl(
 \Pi^{a,b}_{f_P,\star}(\nbigm_{P})
 \bigr):=
(\varphi_P)_{\dagger}\Bigl(
 \Upsilon\bigl(\Pi^{a,b}_{\varphi_P^{\ast}(f),\star}(\nbigv_P)
  \bigr)
  \Bigr)
  \in\gbigc(U_P)
\]
 which are Malgrange extensions of
 $\Pi^{a,b}_{f_P,\star}(\nbigm_{P})\in\nbigc(U_P)$
 by Proposition \ref{prop;21.1.23.2}.
Because the Beilinson functor
$\Pi^{a,b}_{f_P,\ast!}(\nbigm_{P})$
is obtained as
the cokernel of
$\Pi^{b,N}_{f_P,!}(\nbigm_{P})
\lrarr
 \Pi^{a,N}_{f_P,\ast}(\nbigm_{P})$
for a sufficiently large $N$,
we obtain a Malgrange extension
$\Upsilon(\Pi^{a,b}_{f_P,\ast!}(\nbigm_{P}))$
as the cokernel of
$\Upsilon\bigl(\Pi^{b,N}_{f_P,!}(\nbigm_{P})\bigr)
\lrarr
 \Upsilon\bigl(\Pi^{a,N}_{f_P,\ast}(\nbigm_{P})\bigr)$.
Recall that
$\Xi^{(a)}_{f_P}(\nbigm_{P})
=\Pi^{a,a+1}_{f_P,\ast!}(\nbigm_{P})$
and
$\psi^{(a)}_{f_P}(\nbigm_{P})
=\Pi^{a,a}_{f_P,\ast!}(\nbigm_{P})$.
Hence, we obtain Malgrange extensions
$\Upsilon(\Xi^{(a)}_{f_P}(\nbigm_{P}))$ and
$\Upsilon(\psi^{(a)}_{f_P}(\nbigm_{P}))$
of $\Xi^{(a)}_{f_P}(\nbigm_{P})$
and $\psi^{(a)}_{f_P}(\nbigm_{P})$,
respectively.
We also obtain the induced morphisms
\begin{equation}
\label{eq;21.2.16.1}
\Upsilon( \psi^{(1)}_{f_P}(\nbigm_{P}))
\stackrel{\gamma_0}{\lrarr}
\Upsilon(\Xi^{(0)}_{f_P}(\nbigm_{P}))
\stackrel{\beta_0}{\lrarr}
\Upsilon( \psi^{(0)}_{f_P}(\nbigm_{P})),
\end{equation}
which is the extension of the natural morphisms
$\psi^{(1)}_{f_P}(\nbigm_{P})
\stackrel{\gamma_0}{\lrarr}
 \Xi^{(0)}_{f_P}(\nbigm_{P})
\stackrel{\beta_0}{\lrarr}
 \psi^{(0)}_{f_P}(\nbigm_{P})$.
By the assumption of the induction on the dimension of the support,
there uniquely exist the morphisms
\begin{equation}
\label{eq;21.2.16.2}
 \Upsilon(\psi^{(1)}_{f_P}(\nbigm_{P}))
\stackrel{\can^{(0)}}{\lrarr}
 \Upsilon(\phi^{(0)}_{f_P}(\nbigm_{P}))
\stackrel{\var^{(0)}}{\lrarr}
 \Upsilon(\psi^{(0)}_{f_P}(\nbigm_{P}))
\end{equation}
which is the extension of the natural morphisms
$\psi^{(1)}_{f_P}(\nbigm_{P})
\stackrel{\can^{(0)}}{\lrarr}
 \phi^{(0)}_{f_P}(\nbigm_{P})
\stackrel{\var^{(0)}}{\lrarr}
 \psi^{(0)}_{f_P}(\nbigm_{P})$.
 Note that the compositions of
 (\ref{eq;21.2.16.1}) and (\ref{eq;21.2.16.2})
 are equal
 by the assumption of the induction on the dimension of the support.
Hence, we obtain the following complex
\begin{equation}
\label{eq;21.2.16.3}
 \Upsilon(\psi^{(1)}_{f_P}(\nbigm_P))
 \stackrel{\gamma_0+\can^{(0)}}{\lrarr}
 \Upsilon(\Xi^{(0)}_{f_P}(\nbigm_P))\oplus
 \Upsilon(\phi^{(0)}_{f_P}(\nbigm_P))
\stackrel{-\beta_0+\var^{(0)}}{\lrarr}
 \Upsilon(\psi^{(0)}_{f_P}(\nbigm_P)).
\end{equation}
Because $\nbigm_{P}$ is naturally isomorphic
to the cohomology of
\[
\psi^{(1)}_{f_P}(\nbigm_{P})
 \stackrel{\gamma_0+\can^{(0)}}{\lrarr}
 \Xi^{(0)}_{f_P}(\nbigm_{P})\oplus
 \phi^{(0)}_{f_P}(\nbigm_{P})
 \stackrel{-\beta_0+\var^{(0)}}{\lrarr}
 \psi^{(0)}_{f_P}(\nbigm_{P}),
\]
we obtain a Malgrange extension
$\Upsilon(\nbigm_{P})$ as the cohomology of (\ref{eq;21.2.16.3}).
Thus, we obtain Lemma \ref{lem;21.2.17.1}.
\hfill\qed

\vspace{.1in}

Let $\gbigm_{P,1}$ be any Malgrange extensions of $\nbigm_P$. 
Let us prove $\gbigm_{P,1}=\Upsilon(\nbigm_P)=:\gbigm_{P,2}$.
According to \cite[Theorem 2.0.2]{Wlodarczyk},
there exits a projective morphism of complex manifolds
$\rho_P:U^{(1)}_P\lrarr U_P$ 
such that the following holds.
\begin{itemize}
 \item $\rho_P^{\ast}(f_P)^{-1}(0)$
       is a normal crossing hypersurface of $U_P^{(1)}$.
 \item The proper transform $Z_P^{(1)}$ of $Z_P$
       is smooth, and it intersects with
       $\rho_P^{\ast}(f_P)^{-1}(0)$
       in the normal crossing way.
 \item $\rho_P$ induces an isomorphism
       $Z_P^{(1)}\setminus \rho_P^{\ast}(f_P)^{-1}(0)
       \simeq
       Z_P\setminus f_P^{-1}(0)$.
\end{itemize}
We set $H^{(1)}_P:=\rho_P^{\ast}(f_P)^{-1}(0)$.
There uniquely exists an $\gbigrtilde_{U^{(1)}_P(\ast H_P^{(1)})}$-module
$\gbigm^{(1)}_{P,i}$ such that
$\rho_{P\dagger}(\gbigm^{(1)}_{P,i})=\gbigm_{P,i}(\ast f_P)$.
Note that there exists a natural isomorphism
\begin{equation}
\label{eq;21.2.16.12}
 \gbigm^{(1)}_{P,1|(\proj^1\setminus\{\infty\})\times U_P^{(1)}}
 \simeq
 \gbigm^{(1)}_{P,2|(\proj^1\setminus\{\infty\})\times U_P^{(1)}}.
\end{equation}

\begin{lem}
\label{lem;21.2.17.2}
 The isomorphism {\rm(\ref{eq;21.2.16.12})}
naturally extends to an isomorphism
$\gbigm^{(1)}_{P,1}\simeq
\gbigm^{(1)}_{P,2}$.
\end{lem}
\pf
The restriction of
$\gbigm^{(1)}_{P,i}$ to
$(\gbigu_P^{(1)})^{\circ}$ is naturally
a holonomic
$\nbigd_{(\gbigu_P^{(1)})^{\circ}}$-module.
We set $H_{Z,P}^{(1)}:=H_{P}^{(1)}\cap Z_P^{(1)}$.
By the Kashiwara's equivalence for $\nbigd$-modules on submanifolds,
there exist meromorphic flat bundles
$\nbigv_{P,i}$ on
\[
\Bigl(
(\proj^1\setminus\{0\})\times
 Z^{(1)}_P,\,\,
 (\{\infty\}\times Z^{(1)}_P)
 \cup
 \bigl(
 (\proj^1\setminus\{0\})
 \times
 H_{Z,P}^{(1)}
 \bigr)
\Bigr)
\]
such that
$\gbigm^{(1)}_{P,i}$
are the direct image of $\nbigv_{P,i}$.
By the construction,
there exists a natural isomorphism
\begin{equation}
\label{eq;21.2.16.11}
 \nbigv_{P,1|(\proj^1\setminus\{0,\infty\})\times Z_P^{(1)}}
  \simeq
 \nbigv_{P,2|(\proj^1\setminus\{0,\infty\})\times Z_P^{(1)}}.
\end{equation}
Note that the restriction of $\nbigv_{P,i}$
to $(\proj^1\setminus\{0\})\times (Z_P^{(1)}\setminus H_P^{(1)})$
is strongly regular along $\{\infty\}\times(Z_P^{(1)}\setminus H_P^{(1)})$,
and hence they are regular singular
as in Proposition \ref{prop;21.2.16.10}.
Therefore,
the isomorphism (\ref{eq;21.2.16.11})
uniquely extends to an isomorphism
on
$\bigl(
 (\proj^1\setminus\{0\})\times Z_P^{(1)}
 \bigr)
\setminus
(\{\infty\}\times H_{Z,P}^{(1)})$.
By the Hartogs theorem,
it extends to an isomorphism
$\nbigv_{P,1}\simeq\nbigv_{P,2}$
on $(\proj^1\setminus\{0\})\times Z_P^{(1)}$.
As a result,
the natural isomorphism
$\gbigm^{(1)}_{P,1|(\proj^1\setminus\{\infty\})\times U_P^{(1)}}
\simeq
\gbigm^{(1)}_{P,2|(\proj^1\setminus\{\infty\})\times U_P^{(1)}}$
uniquely extends to an isomorphism
$\gbigm^{(1)}_{P,1}\simeq\gbigm^{(1)}_{P,2}$.
\hfill\qed

\vspace{.1in}

By Lemma \ref{lem;21.2.17.2}, we obtain
$\gbigm_{P,1}(\ast f_P)
=\gbigm_{P,2}(\ast f_P)$.
We obtain isomorphisms of the $\gbigrtilde_{U_P}(\ast f_P)$-modules:
\[
 \Pi_{f_P}^{a,b}(\gbigm_{P,1})
 \simeq
 \Pi_{f_P}^{a,b}(\gbigm_{P,2}).
\]
We obtain the following natural complexes of
$\gbigrtilde_{U_P}$-modules:
\begin{equation}
\label{eq;21.2.16.20}
 \gbigm_{P,i}[! f_P]
 \lrarr
 \gbigm_{P,i}
 \oplus
 \Xi^{(0)}_{f_P}(\gbigm_{P,i})
 \lrarr
 \gbigm_{P,i}[\ast f_P].
\end{equation}
They are the extensions of the following complex:
\begin{equation}
\label{eq;21.2.16.22}
 \nbigm_{P}[!f_P]
  \lrarr
  \nbigm_{P}
  \oplus
  \Xi^{(0)}_{f_P}(\nbigm_P)
  \lrarr
  \nbigm_{P}[\ast f_P].
\end{equation}
Let $\phi^{(0)}_{f_P}(\gbigm_{P,i})$
be the cohomology of (\ref{eq;21.2.16.20}).
We recover $\gbigm_{P,i}$ as the cohomology of
the following naturally defined complexes:
\begin{equation}
 \label{eq;21.2.16.21}
  \psi^{(1)}_{f_P}(\gbigm_{P,i})
  \lrarr
  \phi^{(0)}_{f_P}(\gbigm_{P,i})
  \oplus
  \Xi^{(0)}_{f_P}(\gbigm_{P,i})
  \lrarr
  \psi^{(0)}_{f_P}(\gbigm_{P,i}).
\end{equation}
They are the extensions of the following complex:
\begin{equation}
 \label{eq;21.2.16.24}
  \psi^{(1)}_{f_P}(\nbigm_{P})
  \lrarr
  \phi^{(0)}_{f_P}(\nbigm_{P})
  \oplus
  \Xi^{(0)}_{f_P}(\nbigm_{P})
  \lrarr
  \psi^{(0)}_{f_P}(\nbigm_{P}).
\end{equation}
We have already proved the following claim
in the construction of $\gbigm_2=\Upsilon(\nbigm)$.
\begin{lem}
\label{lem;21.2.16.25}
 $\Xi^{(0)}_{f_P}(\gbigm_{P,2})$,
 $\psi^{(i)}_{f_P}(\gbigm_{P,2})$
 and
 $\phi^{(0)}_{f_P}(\gbigm_{P,2})$
 are Malgrange extensions of
  $\Xi^{(0)}_{f_P}(\nbigm_{P})$,
 $\psi^{(i)}_{f_P}(\nbigm_{P})$
 and
 $\phi^{(0)}_{f_P}(\nbigm_{P})$,
 respectively. 
\hfill\qed
\end{lem}

\begin{lem}
 $\Xi^{(0)}_{f_P}(\gbigm_{P,1})$,
 $\psi^{(i)}_{f_P}(\gbigm_{P,1})$
 and
 $\phi^{(0)}_{f_P}(\gbigm_{P,1})$
 are Malgrange extensions of
  $\Xi^{(0)}_{f_P}(\nbigm_{P})$,
 $\psi^{(i)}_{f_P}(\nbigm_{P})$
 and
 $\phi^{(0)}_{f_P}(\nbigm_{P})$,
 respectively.  
\end{lem}
\pf
Because $\gbigm_{P,1}(\ast f_P)=\gbigm_{P,2}(\ast f_P)$,
we obtain the claims for
$\Xi^{(0)}_{f_P}(\gbigm_{P,1})$
and
$\psi^{(i)}_{f_P}(\gbigm_{P,1})$
from Lemma \ref{lem;21.2.16.25}.
Because $\gbigm_{P,1}$ is assumed to be
an Malgrange extension of $\nbigm_P$,
$\phi^{(0)}_{f_P}(\gbigm_{P,1})$
is also a Malgrange extension of
$\phi^{(0)}_{f_P}(\nbigm_P)$.
\hfill\qed

\vspace{.1in}
By the assumption of the induction on the dimension of the support,
we obtain
$\phi^{(0)}_{f_P}(\gbigm_{P,1})
=\phi^{(0)}_{f_P}(\gbigm_{P,2})$.
Because $\gbigm_{P,i}$ are recovered
as the cohomology of (\ref{eq;21.2.16.21}),
we obtain that $\gbigm_{P,1}=\gbigm_{P,2}$,
as desired.

By varying $P\in X$
and gluing $\Upsilon(\nbigm_P)$,
we obtain the existence and the uniqueness of
a Malgrange extension $\Upsilon(\nbigm)$ of $\nbigm$.

\vspace{.1in}

Let $g:\nbigm_1\lrarr\nbigm_2$
be any morphism in $\nbigc(X)$.
We set
$Z=\Supp(\nbigm_1)\cup\Supp(\nbigm_2)$.
Let $P$ be any point of $X$.
There exist a neighbourhood $U_P$ of $P$ in $X$,
and a holomorphic function $f_P$ on $U_P$
such that
(i) $(Z\cap U_P)\setminus f_P^{-1}(0)$
is a closed submanifold of $U_P\setminus f_P^{-1}(0)$,
(ii) $\dim (Z\cap f_P^{-1}(0))<\dim ((Z\cap U_P)\setminus f_P^{-1}(0))$.
We set
$\nbigm_{i,P}:=\nbigm_{i|U_P}$.
\begin{lem}
\label{lem;21.3.11.1}
 $g_{|U_P}$ uniquely extends to
a morphism
 $\Upsilon(g_{|U_P}):
 \Upsilon(\nbigm_{1,P})
 \lrarr
 \Upsilon(\nbigm_{2,P})$.
\end{lem}
\pf
The uniqueness is clear.
Each $\nbigm_{i,P}$ is expressed as the cohomology of
\[
 \psi_{f_P}^{(0)}(\nbigm_{i,P})
 \lrarr
 \Xi_{f_P}^{(0)}(\nbigm_{i,P})\oplus
 \phi_{f_P}^{(0)}(\nbigm_{i,P})
 \lrarr
 \psi_{f_P}^{(1)}(\nbigm_{i,P}).
\]
There exist the induced morphisms
\begin{equation}
\label{eq;21.2.16.30}
\psi_{f_P}^{(i)}(g):
 \psi_{f_P}^{(i)}(\nbigm_{1,P})
 \lrarr
  \psi_{f_P}^{(i)}(\nbigm_{2,P}),
\end{equation}
\begin{equation}
\label{eq;21.2.16.31}
 \Xi_{f_P}^{(i)}(g):
 \Xi_{f_P}^{(i)}(\nbigm_{1,P})
 \lrarr
  \Xi_{f_P}^{(i)}(\nbigm_{2,P}),
\end{equation}
\begin{equation}
\label{eq;21.2.16.32}
 \phi_{f_P}^{(i)}(g):
 \phi_{f_P}^{(i)}(\nbigm_{1,P})
 \lrarr
  \phi_{f_P}^{(i)}(\nbigm_{2,P}).
\end{equation}
By the assumption of the induction on the dimension of the support,
(\ref{eq;21.2.16.30})
and 
(\ref{eq;21.2.16.32})
uniquely extend to morphisms
\begin{equation}
\label{eq;21.2.16.33}
\Upsilon(\psi_{f_P}^{(i)}(g)):
 \Upsilon(\psi_{f_P}^{(i)}(\nbigm_{1,P}))
 \lrarr
 \Upsilon(\psi_{f_P}^{(i)}(\nbigm_{2,P})),
\end{equation}
\begin{equation}
\label{eq;21.2.16.34}
\Upsilon(\phi_{f_P}^{(i)}(g)):
 \Upsilon(\phi_{f_P}^{(i)}(\nbigm_{1,P}))
 \lrarr
 \Upsilon(\phi_{f_P}^{(i)}(\nbigm_{2,P})).
\end{equation}
By using
Proposition \ref{prop;21.2.17.3}
and the description of
$\nbigm_{1,P}(\ast f_P)$
as the direct image
of good-KMS $\nbigr_{Z'_P}(\ast H_P')$-module
for some $(Z_P',H_P')$,
we obtain extensions
\begin{equation}
\label{eq;21.2.16.35}
\Upsilon(\psi_{f_P}^{(i)}(g)):
 \Upsilon(\psi_{f_P}^{(i)}(\nbigm_{1,P}))
 \lrarr
 \Upsilon(\psi_{f_P}^{(i)}(\nbigm_{2,P})),
\end{equation}
\begin{equation}
\label{eq;21.2.16.36}
\Upsilon(\Xi_{f_P}^{(i)}(g)):
 \Upsilon(\Xi_{f_P}^{(i)}(\nbigm_{1,P}))
 \lrarr
 \Upsilon(\Xi_{f_P}^{(i)}(\nbigm_{2,P})).
\end{equation}
Note that (\ref{eq;21.2.16.33})
and (\ref{eq;21.2.16.35}) are equal by the uniqueness.
Hence, the above morphisms induce the following commutative diagram:
\begin{equation}
\label{eq;21.2.16.40}
 \begin{CD}
  \Upsilon(\psi_{f_P}^{(1)}(\nbigm_{1,P}))
  @>>>
  \Upsilon(\Xi_{f_P}^{(0)}(\nbigm_{1,P}))
  \oplus
  \Upsilon(\phi_{f_P}^{(0)}(\nbigm_{1,P}))
  @>>>
  \Upsilon(\psi_{f_P}^{(0)}(\nbigm_{1,P}))
  \\
  @VVV @VVV @VVV\\
  \Upsilon(\psi_{f_P}^{(1)}(\nbigm_{2,P}))
  @>>>
  \Upsilon(\Xi_{f_P}^{(0)}(\nbigm_{2,P}))
  \oplus
  \Upsilon(\phi_{f_P}^{(0)}(\nbigm_{2,P}))
  @>>>
  \Upsilon(\psi_{f_P}^{(0)}(\nbigm_{2,P})).
 \end{CD}
\end{equation}
Thus, we obtain the desired morphism
$\Upsilon(g_{|U_P})$,
and Lemma \ref{lem;21.3.11.1} is proved.
\hfill\qed

\vspace{.1in}

Then, it is easy to prove the existence and uniqueness of
a morphism
$\Upsilon(g):\Upsilon(\nbigm_1)\lrarr\Upsilon(\nbigm_2)$
whose restriction to
$(\proj^1\setminus\{\infty\})\times X$
is equal to $g$.
The fully faithfulness is obvious.
The proof of Theorem \ref{thm;21.1.23.1}
is completed.
\hfill\qed

\subsubsection{Proof of Proposition \ref{prop;21.4.16.2}}
\label{subsection;21.4.16.1}

It is enough to prove the claim locally around
any point $P$ of $X$.
We use the induction on the dimension of the support.
Let $U_P$, $f_P$ and $\nbigm_P$ be as above.
There exists the morphism
$\kappa:\nbigm_P\lrarr \nbigm_P[\ast f_P]$.
Note that $\Ker\kappa\in \nbigc(X)$
and that the dimension of the support of $\Ker\kappa$
is strictly smaller than
the dimension of the support of $\nbigm_P$.
By the assumption of the induction,
we obtain
$\bigcap_{a\in\real}
\Bigl(\lefttop{\infty}V_a\bigl(\Upsilon(\Ker\kappa)\bigr)\Bigr)=0$.

\begin{lem}
\label{lem;21.4.16.10}
 $\bigcap_{a\in\real}\Bigl(
 \lefttop{\infty}V_a\bigl(\Upsilon(\nbigm_P[\ast f_P])\bigr)
 \Bigr)=0$.
\end{lem}
\pf
Let $Q$ be any point of
$\Supp(\nbigm_P)\setminus f_P^{-1}(0)$.
Let $U_Q$ be a neighbourhood of $Q$ in $X$
such that $U_Q\cap f_P^{-1}(0)=\emptyset$.
We set $Z_Q:=\Supp(\nbigm_P)\cap U_Q$,
which is a closed complex submanifold of
$U_Q$.
Then, there exists a regular meromorphic flat bundle
$\nbigv_Q$ on
$\bigl(
(\proj^1\setminus\{0\})\times Z_Q,\infty\times Z_Q
\bigr)$
such that 
$\Upsilon(\nbigm_{P|U_Q})
_{|(\proj^1\setminus\{0\})\times U_Q}$
is the direct image of
$\nbigv_Q$.
Hence, it is easy to check
$\bigcap_{a\in\real}\bigl(
\lefttop{\infty}V_a\bigl(\Upsilon(\nbigm_P)\bigr)_{|U_Q}
\bigr)=0$.

Let $s$ be a section of
$\bigcap_{a\in\real}\bigl(
\lefttop{\infty}
V_a(\Upsilon(\nbigm[\ast f_P]))\bigr)$
on $\nbigu\subset\proj^1\times U_P$.
By the previous consideration,
for any $Q\in \Supp(\nbigm_P)\setminus f_P^{-1}(0)$,
we obtain $s_{|(\proj^1\times U_Q)\cap\nbigu}=0$
for a neighbourhood $U_Q$ of $Q$.
It implies $s=0$.
\hfill\qed

\vspace{.1in}
Let $s$ be a section of
$\bigcap_{a\in\real}
\bigl(
\lefttop{\infty}
V_a(\Upsilon(\nbigm_P))\bigr)$.
Because $\kappa(s)$ is a section of
$\bigcap_{a\in\real}
\bigl(\lefttop{\infty}V_a(\Upsilon(\nbigm_P[\ast f_P]))\bigr)$,
we obtain $\kappa(s)=0$ by Lemma \ref{lem;21.4.16.10}.
Because
$\lefttop{\infty}V_a(\Upsilon(\Ker\kappa))=
\lefttop{\infty}V_a(\Upsilon(\nbigm_P))\cap
\Upsilon(\Ker\kappa)$,
we obtain $s=0$.
Thus, we obtain Proposition \ref{prop;21.4.16.2}.
\hfill\qed

\subsubsection{Proof of Proposition \ref{prop;21.6.26.2}}
\label{subsection;22.7.27.1}

We use the induction on the dimension of the support.
The claim is trivial in the $0$-dimensional case.
It is enough to study the claim locally around any point of $X$.
Let $P$ be any point of $X$.
Let $U_P$, $f_P$, $\nbigm_P$ and $Z_P$ be as in \S\ref{subsection;21.4.14.10}.
We may assume that $g^{-1}(0)\cap U_P\subset f_P^{-1}(0)$.
There exists $\varphi_P$ and $\nbigv_P$
as in the proof of Lemma \ref{lem;21.2.17.1}.
By the standard argument,
it is enough to prove the claim of Proposition \ref{prop;21.6.26.2}
for $\Upsilon(\nbigv_P[\star \varphi_P^{\ast}(f_P)])$
with the function $\varphi_P^{\ast}(g)$
which follows from Lemma \ref{lem;21.6.26.1}.
\hfill\qed

\subsection{Basic functoriality of the Malgrange extension}
\label{subsection;21.6.29.1}
 
We obtain the following proposition
from Lemma \ref{lem;21.3.9.33}.

\begin{prop}
Let $g:\nbigm_1\lrarr\nbigm_2$ be a morphism in $\nbigc(X)$.
Assume that
$\Ker(g)$, $\Image(g)$ and $\Cok(g)$ are also objects
in $\nbigc(X)$.
Then, we obtain
$\Upsilon(\Ker g)= \Ker\Upsilon(g)$,
$\Upsilon(\Image g)=\Image\Upsilon(g)$
and
$\Upsilon(\Cok(g))=\Cok\Upsilon(g)$.
\hfill\qed
\end{prop}

We obtain the following proposition from Proposition \ref{prop;21.1.23.2}.
\begin{prop}
\label{prop;21.3.13.10}
 Let $F:X\lrarr Y$ be a projective morphism.
 For $\nbigm\in\nbigc(X)$,
 we obtain
$\Upsilon (F^j_{\dagger}(\nbigm))
=F^j_{\dagger}(\Upsilon(\nbigm))$.
As a result,
for any $\gbigm\in\gbigc_{\Malg}(X)$,
$F^j_{\dagger}(\gbigm)$ are objects
of $\gbigc_{\Malg}(Y)$.
\hfill\qed
\end{prop}

Let $H$ be any closed complex hypersurface of $X$.
\begin{prop}
\label{prop;21.2.18.1}
For any $\nbigm\in\nbigc(X)$,
we have
 $\Upsilon(\nbigm)[\star H]
=\Upsilon(\nbigm[\star H])$. 
As a result,
for any $\gbigm\in\gbigc_{\Malg}(X)$,
$\gbigc[\star H]$  are objects of $\gbigc_{\Malg}(X)$.
\end{prop}
\pf
We use an induction on the dimension of the support of $\nbigm$.
It is enough to study the claim locally around
any point of $X$.
Therefore, we may assume that
there exist a holomorphic function $f$ on $X$ and
a projective morphism of complex manifolds $\varphi:Z\lrarr X$
such that the following holds.
\begin{itemize}
 \item $H_Z:=\varphi^{-1}(H\cup f^{-1}(0))$
       is a simple normal crossing hypersurface of $Z$.
       We set $H_{Z,1}:=\varphi^{-1}(f^{-1}(0))$.
 \item $\varphi(Z)\setminus f^{-1}(0)=\Supp(\nbigm)\setminus f^{-1}(0)$,
       and $\dim(\Supp(\nbigm)\cap f^{-1}(0))<\dim Z$.
 \item There exists a good-KMS
       $\nbigr_{Z(\ast H_{Z,1})}$-module
       $\nbigv$ with an isomorphism
       $\varphi_{\dagger}(\nbigv)
       \simeq
        \nbigm(\ast f)$.
 \end{itemize}
Then, $\nbigm[\star H]$ is isomorphic to the cohomology of the following complex:
\[
 \psi^{(1)}_f(\nbigm)[\star H]
 \lrarr
 \Xi^{(0)}_f(\nbigm)[\star H]
 \oplus
 \phi^{(0)}_f(\nbigm)[\star H]
 \lrarr
 \psi^{(0)}_f(\nbigm)[\star H].
\]
Under the isomorphism,
$\Upsilon(\nbigm[\star H])$ is equal to
the cohomology of the following complex:
\[
\Upsilon\bigl(
 \psi^{(1)}_f(\nbigm)[\star H]
 \bigr)
 \lrarr
\Upsilon\bigl(
 \Xi^{(0)}_f(\nbigm)[\star H]
 \bigr)
 \oplus
\Upsilon\bigl(
 \phi^{(0)}_f(\nbigm)[\star H]
\bigr)
 \lrarr
\Upsilon\bigl(
 \psi^{(0)}_f(\nbigm)[\star H]
\bigr).
\]
By the assumption of the induction,
we have
$\Upsilon\bigl(
 \psi^{(a)}_f(\nbigm)[\star H]
 \bigr)
 =\Upsilon\bigl(
 \psi^{(a)}_f(\nbigm)
 \bigr)[\star H]$
and
$\Upsilon\bigl(
  \phi^{(a)}_f(\nbigm)[\star H]
 \bigr)
 =\Upsilon\bigl(
 \phi^{(a)}_f(\nbigm)
 \bigr)[\star H]$.
 By setting $f_Z:=\varphi^{\ast}(f)$
 and $H_{Z,2}:=\varphi^{-1}(H)$,
we obtain
\[
\Upsilon\bigl(
 \Xi_{f}^{(a)}(\nbigm)[\star H]
\bigr)
=\varphi_{\dagger}\Bigl(
\Upsilon\Bigl(
 \Xi_{f_Z}^{(a)}(\nbigv)[\star H_{Z,2}]
 \Bigr)
 \Bigr).
\]
Recall that
$\Upsilon\Bigl(
 \Xi_{f_Z}^{(a)}(\nbigv)[\star H_{Z,2}]
 \Bigr)$
 is isomorphic to the cokernel of
\[
 \Upsilon\Bigl(
 \Pi_{f_Z,!}^{a+1,N}(\nbigv)
 [\star H_{Z,2}]
  \Bigr)
  \lrarr
\Upsilon\Bigl(
 \Pi_{f_Z,\ast}^{a,N}(\nbigv)
 [\star H_{Z,2}]
 \Bigr)
\] 
for large $N$.
Because
$\Upsilon\Bigl(
 \Pi^{a,b}_{f_Z,\star_1}(\nbigv)[\star_2 H_{Z,2}]
 \Bigr)
=\Upsilon\Bigl(
 \Pi^{a,b}_{f_Z,\star_1}(\nbigv)
 \Bigr)[\star_2 H_{Z,2}]$
for $\star_1,\star_2\in\{\ast,!\}$,
we obtain
\[
 \Upsilon\bigl(
 \Xi_{f}^{(a)}(\nbigm)[\star H]
 \bigr)
 =\Upsilon\bigl(\Xi_f^{(a)}(\nbigm)\bigr)[\star H].
\]
Hence,
$\Upsilon(\nbigm[\star H])$
is equal to the cohomology of the following complex:
\[
 \Upsilon\bigl(\psi^{(1)}_f(\nbigm)\bigr)[\star H]
 \lrarr
 \Upsilon\bigl(
 \Xi^{(0)}_f(\nbigm)\bigr)[\star H]
 \oplus
 \Upsilon\bigl(
 \phi^{(0)}_f(\nbigm)\bigr)[\star H]
 \lrarr
 \Upsilon\bigl(
 \psi^{(0)}_f(\nbigm)\bigr)[\star H].
\]
It is equal to $\Upsilon(\nbigm)[\star H]$.
\hfill\qed

\begin{cor}
\mbox{{}}
 \begin{itemize}
 \item  For any $\nbigm\in\nbigc(X;H)$,
there uniquely exists a Malgrange extension
	$\Upsilon(\nbigm)$ of $\nbigm$.
	This induces a functor from
	$\nbigc(X;H)$ to
	$\gbigc(X;H)$.
	Let $\gbigc_{\Malg}(X;H)$
	denote the essential image.
  \item
       For any $\nbigm_0\in\nbigc(X)$,
we have
       $\Upsilon(\nbigm_0(\ast H))=\Upsilon(\nbigm_0)(\ast H)$.
In other words,
       $\gbigc_{\Malg}(X;H)$
       is the essential image of
       the natural functor
       $\gbigc_{\Malg}(X)\lrarr\gbigc(X;H)$.
 \end{itemize}
\end{cor}
\pf
Let $\nbigm_0\in\nbigc(X)$.
We have the Malgrange extension
$\Upsilon(\nbigm_0)$,
and
$\Upsilon(\nbigm_0[\ast H])=
\Upsilon(\nbigm_0)[\ast H]$ holds.
Because
$\nbigm_0[\ast H]_{|\gbigx\setminus\gbigx^0}
=\nbigm_0(\ast H)_{|\gbigx\setminus\gbigx^0}$,
$\Upsilon(\nbigm_0)(\ast H)$
is the Malgrange extension of
$\nbigm_0(\ast H)$.
Then, the claims of the corollary immediately follows.
\hfill\qed

\begin{cor}
\label{cor;21.6.22.17}
Let $F:X\lrarr Y$ be a projective morphism of complex manifolds.
Let $H_Y$ be a hypersurface of $Y$,
and we set $H_X:=F^{-1}(H_Y)$.
\begin{itemize}
 \item We obtain the induced functors
       $F^{j}_{\dagger}:\gbigc_{\Malg}(X;H_X)\lrarr \gbigc_{\Malg}(Y;H_Y)$.
 \item If $F$ induces an isomorphism $X\setminus H_X\simeq Y\setminus H_Y$,
       $F^0_{\dagger}:\gbigc_{\Malg}(X;H_X)\lrarr
       \gbigc_{\Malg}(Y;H_Y)$ is an equivalence.
       We have $F^0_{\dagger}(\gbigm)=F_{\ast}(\gbigm)$
       for $\gbigm\in\gbigc_{\Malg}(X;H_X)$.
       The quasi inverse is given by
       $\gbign\longmapsto F^{\ast}(\gbign)$
       for $\gbign\in\gbigc_{\Malg}(Y;H_Y)$.
 \item Suppose that $F_{|X\setminus H_X}$ is a closed embedding
       of $X\setminus H_X$ into $Y\setminus H_Y$.
       We set
       $\gbigc_{F(X),\Malg}(Y;H_Y)=
       \gbigc_{\Malg}(Y;H_Y)\cap\gbigc_{F(X)(Y;H_Y)}$.
       Then, 
       $F_{\dagger}^0:\gbigc_{\Malg}(X;H_X)\lrarr
       \gbigc_{F(X),\Malg}(Y;H_Y)$
       is an equivalence.
\end{itemize}
\end{cor}
\pf
The first claim is clear.
Because the functors
$\Upsilon:\nbigc(X;H_X)\lrarr\gbigc_{\Malg}(X;H_X)$
and
$\Upsilon:\nbigc(Y;H_Y)\lrarr\gbigc_{\Malg}(Y;H_Y)$,
we obtain the second claim from
Proposition \ref{prop;21.6.22.16}
and Proposition \ref{prop;21.3.13.10}.
Similarly, we obtain the third claim from
Proposition \ref{prop;21.4.9.21}
and Proposition \ref{prop;21.3.13.10}.
\hfill\qed

\begin{prop}
\label{prop;21.4.14.12} 
Let $g$ be any meromorphic function on $(X,H)$.
Let $\nbigm\in\nbigc(X;H)$.
\begin{itemize}
 \item
      $\Upsilon\bigl(\Pi^{a,b}_{g,\star}(\nbigm)\bigr)
      =\Pi^{a,b}_{g,\star}(\Upsilon(\nbigm))$
      $(\star=\ast,!)$.
 \item
 We have
$\Upsilon(\Xi^{(a)}_g(\nbigm))
=\Xi^{(a)}_g(\Upsilon(\nbigm))$,
$\Upsilon(\psi^{(a)}_g(\nbigm))
=\psi^{(a)}_g(\Upsilon(\nbigm))$
and 
$\Upsilon(\phi^{(a)}_g(\nbigm))
=\phi^{(a)}_g(\Upsilon(\nbigm))$.
\end{itemize}
In other words,
for any $\gbigm\in\gbigc_{\Malg}(X;H)$,
the induced objects 
$\Pi^{a,b}_{g,\star}(\gbigm)$,
 $\Xi^{(a)}_g(\gbigm)$,
 $\psi^{(a)}_g(\gbigm)$ and
 $\phi^{(a)}_g(\gbigm)$
of $\gbigc(X;H)$ 
are contained in $\gbigc_{\Malg}(X;H)$. 
\end{prop}
\pf
We use an induction of the dimension of the support of $\nbigm$.
It is enough to check the claim locally around any point of $X$.
Therefore, we may assume the existence of
$f$, $Z$, $\varphi$ and $\nbigv$
as in the proof of Proposition \ref{prop;21.2.18.1}.
Then, $\Pi^{a,b}_{g,\star}(\nbigm)$ is naturally isomorphic to
the cohomology of the following complex:
\[
 \Pi^{a,b}_{g,\star}(\psi^{(1)}_f(\nbigm))
 \lrarr
 \Pi^{a,b}_{g,\star}(\Xi^{(0)}_f(\nbigm))
 \oplus
 \Pi^{a,b}_{g,\star}(\phi^{(0)}_f(\nbigm))
 \lrarr
 \Pi^{a,b}_{g,\star}(\psi^{(0)}_f(\nbigm)).
\]
We obtain the first claim by using an argument
similar to the proof of Proposition \ref{prop;21.2.18.1}.
The second claim follows from the first claim.
\hfill\qed

\begin{prop}
\label{prop;21.4.14.20}
 For $\nbigm\in\nbigc(X;H)$,
we have
$\DD_{X(\ast H)}(\Upsilon(\nbigm))
=\Upsilon(\DD_{X(\ast H)}(\nbigm))$.
In other words,
for any $\gbigm\in\gbigc_{\Malg}(X;H)$,
$\DD_{X(\ast H)}(\gbigm)$
is an object of $\gbigc_{\Malg}(X;H)$. 
\end{prop}
\pf
We may assume that there exists
$f$, $Z$, $\nbigv$ and $\varphi$ as above.
We obtain the description of
$\DD\Upsilon(\nbigm)$
as the cohomology of the following complex:
\[
 \DD\bigl(
 \Upsilon(\psi^{(0)}_f(\nbigm) )
 \bigr)
\lrarr
 \DD\bigl(
 \Upsilon(\Xi^{(0)}_f(\nbigm) )
 \bigr)
\oplus
 \DD\bigl(
 \Upsilon(\phi^{(0)}_f(\nbigm) )
 \bigr)
 \lrarr
 \DD\bigl(
 \Upsilon(\psi^{(1)}_f(\nbigm) )
 \bigr).
\]
Because
$\DD\psi^{(i)}_f(\nbigm)
\simeq
\psi^{(1-i)}_f(\nbigm)$,
$\DD\phi^{(0)}_f(\nbigm)
\simeq
\phi^{(0)}_f(\nbigm)$,
and 
$\DD\Xi^{(0)}_f(\nbigm)
\simeq
\Xi^{(0)}_f(\nbigm)$,
it is enough to prove the claim for
$\DD\Xi^{(0)}_f(\nbigm)$.
We have
$\varphi_{\dagger}\DD\bigl(
\Upsilon\Xi_{f_Z}^{(0)}(\nbigm)
\bigr)
\simeq
 \DD\varphi_{\dagger}\bigl(
 \Upsilon\Xi_{f_Z}^{(0)}(\nbigm)
 \bigr)$.
Hence, it is enough prove the claim for
$\DD\bigl(
\Upsilon\Xi_{f_Z}^{(0)}(\nbigm)
\bigr)$
which follows from Lemma \ref{lem;21.3.23.20}.
\hfill\qed

\begin{prop}
 \label{prop;21.4.14.22}
Let $X_i$ $(i=1,2)$ be complex manifolds
with a hypersurface $H_i$.
For $\nbigm_i\in\nbigc(X_i;H_i)$ $(i=1,2)$,
we have
$\Upsilon(\nbigm_1\boxtimes\nbigm_2)
=\Upsilon(\nbigm_1)\boxtimes\Upsilon(\nbigm_2)$.
Namely, for $\gbigm_i\in\gbigc_{\Malg}(X_i;H_i)$ $(i=1,2)$,
$\gbigm_1\boxtimes\gbigm_2$
is an object of $\gbigc_{\Malg}(X_1\times X_2;H)$,
where $H=(H_1\times X_2)\cup (X_1\times H_2)$.
\end{prop}
\pf
We explain only an outline.
We may assume that
$f_i$, $Z_i$, $\nbigv_i$, $\varphi_i$ for $\nbigm_i$ $(i=1,2)$.
We set $f_{Z_i}:=\varphi_i^{\ast}(f_{i})$.
By the induction on the support of the dimension
of the support of $\nbigm_1$,
it is enough to prove the claim for
$\Xi^{(0)}_{f_{Z_i}}(\nbigv_i)\boxtimes \nbigm_2$.
By the induction on the support of the dimension
of the support of $\nbigm_2$,
it is enough to prove the claim for
$\Xi^{(0)}_{f_{Z_1}}(\nbigv_1)\boxtimes
 \Xi^{(0)}_{f_{Z_2}}(\nbigv_2)$.
Then, the claim follows from Lemma \ref{lem;21.3.23.11}.
\hfill\qed

\begin{prop}
\label{prop;21.6.29.12}
Let $F:X\lrarr Y$ be a morphism of complex manifolds.
Let $H_Y$ be a hypersurface of $Y$.
Let $\nbigm\in\nbigc(Y;H_Y)$.
Suppose that $F$ is strictly non-characteristic for $\nbigm$.
Then,
we have
$F^{\ast}(\Upsilon(\nbigm))=\Upsilon(F^{\ast}(\nbigm))$.
Namely,
for any $\gbigm\in\gbigc_{\Malg}(Y;H_Y)$
such that $F$ is non-characteristic for $\gbigm$,
$F^{\ast}(\gbigm)$ is an object of $\gbigc_{\Malg}(X;H_X)$
where $H_X=F^{-1}(H_Y)$. 
\end{prop}
\pf
It is enough to check the claim locally around any point of $X$.
If $F$ is a projection,
it is reduced to the functoriality
with respect to external products.
If $F$ is a closed embedding,
it is reduced to Proposition \ref{prop;21.2.18.1}.
(See the proof of Theorem \ref{thm;21.4.7.50}).
\hfill\qed

\begin{prop}
\label{prop;21.7.25.1}
In the situation of {\rm\S\ref{subsection;21.6.22.2}}
and {\rm\S\ref{subsection;21.6.22.20}},
for any $\nbigm\in\nbigc(Y;H_Y)$,
we have
\[
(\lefttop{T}f^{\star})^i(\Upsilon(\nbigm))
\simeq
\Upsilon\bigl(\lefttop{T}f^{\star}(\nbigm)\bigr)\quad
 (\star=!,\ast).
\]
Similarly, in the situation of {\rm\S\ref{subsection;21.6.22.11}}
and {\rm\S\ref{subsection;21.6.22.18}},
for $\gbigm_i\in\nbigc(X;H)$,
we have 
\[
  \Upsilon\bigl(\nbigh^k(\nbigm_1\otimes^{\star}\nbigm_2)\bigr)
\simeq
 \nbigh^k\bigl(
 \Upsilon(\nbigm_1)\otimes^{\star}
 \Upsilon(\nbigm_2)
 \bigr).
\]
\end{prop}
\pf
The first claim follows from
Proposition \ref{prop;21.2.18.1}
and Corollary \ref{cor;21.6.22.17}.
The second claim follows from
Proposition \ref{prop;21.2.18.1},
Proposition \ref{prop;21.4.14.22}
and Corollary \ref{cor;21.6.22.17}.
\hfill\qed

\subsubsection{Example}
\label{subsection;22.7.27.10}

Let $f$ be a meromorphic function on $(X,H)$.
Let $\gbigl(f)$ denote the $\gbigrtilde_{X(\ast H)}$-module
obtained as
$\gbigl(f)=\nbigo_{\gbigx}\bigl(\ast(\gbigx^{\infty}\cup\gbigh)\bigr)$
with the meromorphic integrable connection
$d+d(\lambda^{-1}f)$.
We have
$\gbigl(f)_{|\nbigx}\in\nbigc(X;H)$
as explained in {\rm\cite[\S3.2]{Mochizuki-GKZ}}.
It is easy to check that
$\gbigl(f)$ is strongly regular along $\gbigx^{\infty}$.
Hence, $\gbigl(f)\in\gbigc_{\Malg}(X;H)$.

Let $(M,F)$ be a good filtered $\nbigd_X$-module
underlying a mixed Hodge module.
Let $\Rtilde_F(M)$ denote the $\nbigrtilde_X$-module
obtained as the analytification of the Rees module of $(M,F)$.
(See \S\ref{subsection;22.7.28.1}.)
It naturally extends to
the $\gbigrtilde_X$-module $\gbigr_F(M)$
such that
$(\gbigrtilde_M)(\ast\gbigx^0)
=p_{1,X}^{\ast}(M)(\ast(\gbigx^0\cup\gbigx^{\infty}))$.
Then, 
we have $\Rtilde_F(M)\in\nbigc(M)\in\nbigc(X)$
as explained in \cite[\S13.5]{Mochizuki-MTM}.
Clearly $\gbigr_F(M)$ is strongly regular along $\gbigx^{\infty}$.
Hence, $\gbigr_F(M)\in \gbigc_{\Malg}(X)$.

Starting from these basic objects,
we can construct many examples
by using the functoriality.
For example, we have the following proposition.

\begin{prop}
$\gbigl(f)\otimes \gbigr_F(M)$
is an object of $\gbigc_{\Malg}(X;H)$.
\end{prop}
\pf
It follows from
Proposition \ref{prop;21.4.14.22}
and Proposition \ref{prop;21.6.29.12}.
\hfill\qed

\section{Rescalable objects and the associated irregular Hodge filtrations}

\subsection{Preliminary}

\subsubsection{Rescaling of $\nbigrtilde$-modules}

Let 
$\varphi_1:\cnum_{\lambda}\times\cnum_{\tau}^{\ast}\lrarr
\cnum_{\lambda}$
be defined by $\varphi_1(\lambda,\tau)=\lambda\tau^{-1}$.
For a complex manifold $X$,
we set $\lefttop{\tau}X:=\cnum_{\tau}\times X$
and $\lefttop{\tau}\nbigx:=\cnum_{\lambda}\times \lefttop{\tau}X$
by following \cite{Sabbah-irregular-Hodge}.
In general, for any complex analytic subset $Z\subset X$,
we set $\lefttop{\tau}Z=\cnum_{\tau}\times Z\subset \lefttop{\tau}X$.
We also set
$\lefttop{\tau}X_0:=\{0\}\times X\subset \lefttop{\tau}X$
and
$\lefttop{\tau}\nbigx_0:=\cnum_{\lambda}\times \lefttop{\tau}X_0$.
We obtain the induced morphism
$\varphi_1:
 \lefttop{\tau}\nbigx\setminus \lefttop{\tau}\nbigx_0
 \lrarr
 \nbigx$.

For an $\nbigo_{\nbigx}$-module $\nbigm$,
we obtain the
$\nbigo_{\lefttop{\tau}\nbigx\setminus\lefttop{\tau}\nbigx_0}$-module
$\lefttop{\tau}\nbigm=\varphi_1^{\ast}(\nbigm)$.
If $\nbigm$ is equipped with a meromorphic flat connection $\nabla$,
then $\lefttop{\tau}\nbigm$ is naturally
equipped with the induced meromorphic flat connection
$\varphi_1^{\ast}(\nabla)$.
Hence, if $\nbigm$ is an $\nbigrtilde_X$-module
then $\lefttop{\tau}\nbigm$ is naturally
an $\nbigrtilde_{\lefttop{\tau}X\setminus\lefttop{\tau}X_0}$-module.
(See \cite{Sabbah-irregular-Hodge}.)

\begin{lem}[\cite{Sabbah-irregular-Hodge}]
\label{lem;21.3.13.30}  
Let $\nbigm\in\nbigc(X)$.
\begin{itemize}
 \item For any projective morphism $F:X\lrarr Y$,
       we have
       $F^i_{\dagger}(\lefttop{\tau}\nbigm)
       =\lefttop{\tau}F^i_{\dagger}(\nbigm)$.
 \item For any hypersurface $H$ of $X$,
       the $\nbigrtilde_{\lefttop{\tau}X\setminus\lefttop{\tau}X_0}$-module
       $\lefttop{\tau}\nbigm$ is localizable along $H$,
       and we have
       $\lefttop{\tau}(\nbigm[\star H])
       =(\lefttop{\tau}\nbigm)[\star (\lefttop{\tau}H)]$.
 \item Let $f$ be any holomorphic function on $X$,
       and let $f_0$ denote the induced holomorphic function on $\lefttop{\tau}X$.
       Then, for any $a\leq b$,
       we have
       $\lefttop{\tau}\bigl(
       \Pi^{a,b}_{f,\star}(\nbigm)
       \bigr)
       =\Pi^{a,b}_{f_0,\star}(\lefttop{\tau}\nbigm)$
       $(\star=!,\ast)$,
       and
      $\lefttop{\tau}\bigl(
       \Pi^{a,b}_{f,\ast !}(\nbigm)
       \bigr)
       =\Pi^{a,b}_{f_0,\ast !}(\lefttop{\tau}\nbigm)$.
       Hence, we have
       $\lefttop{\tau}\bigl(
       \Xi^{(a)}_f(\nbigm)\bigr)
       =\Xi^{(a)}_{f_0}(\lefttop{\tau}\nbigm)$,
       $\lefttop{\tau}\bigl(
       \psi^{(a)}_f(\nbigm)\bigr)
       =\psi^{(a)}_{f_0}(\lefttop{\tau}\nbigm)$,
       and
       $\lefttop{\tau}\bigl(
       \phi^{(a)}_f(\nbigm)\bigr)
       =\phi^{(a)}_{f_0}(\lefttop{\tau}\nbigm)$.
\end{itemize}
\end{lem}
\pf
The first claim is clear by the constructions of
the direct image and the rescaling.
For the second claim,
we may assume that there exists a holomorphic function $f$
on $X$ such that $H=f^{-1}(0)$.
By using the first claim with the graphs of $f$ and $f_0$,
it is enough to consider the case where
$f$ is a coordinate function.
Let $\lefttop{f}V\nbigr_X\subset\nbigr_X$
denote the sheaf of subalgebras
generated by $\lambda \Theta_{X}(\log f)$ over $\nbigo_{\nbigx}$
and
$\lefttop{f}V\nbigrtilde_X=
\lefttop{f}V\nbigr_X\langle\lambda^2\del_{\lambda}\rangle$.
We use the notation
$\lefttop{f_0}V\nbigr_{\lefttop{\tau}X}$
and $\lefttop{f_0}V\nbigrtilde_{\lefttop{\tau}X}$ in similar meanings.
There exists a $V$-filtration $V_{\bullet}(\nbigm)$ of $\nbigm$
along $f$ by $\lefttop{f}V\nbigrtilde_X$-submodules
which are $\lefttop{f}V\nbigr_X$-coherent.
The induced $\nbigo_{\lefttop{\tau}\nbigx}$-submodules
$\lefttop{\tau}\bigl(
\lefttop{f}V_{a}(\nbigm)\bigr)$
are naturally $\lefttop{f_0}V\nbigrtilde_{\lefttop{\tau}X}$-submodules
which are coherent over $\lefttop{f}V\nbigr_{\lefttop{\tau}X}$.
It is easy to see that
the filtration
$\lefttop{\tau}\bigl(\lefttop{f}V_{\bullet}(\nbigm)\bigr)$
is a $V$-filtration of $\lefttop{\tau}\nbigm$ along $f_0$.
We have
$\lefttop{\tau}(\nbigm[!f])(\ast f)
=\lefttop{\tau}(\nbigm)(\ast f)$,
and the natural morphism
$\can:
\Gr^V_{-1}(\lefttop{\tau}\nbigm[!f])
\lrarr
\Gr^V_{0}(\lefttop{\tau}\nbigm[!f])$
is an isomorphism.
Hence, we have
$\lefttop{\tau}(\nbigm[!f])
=(\lefttop{\tau}\nbigm)[!f]$.
Similarly, we have
$\lefttop{\tau}(\nbigm[\ast f])
=(\lefttop{\tau}\nbigm)[\ast f]$.

By the construction,
we have
$\lefttop{\tau}\bigl(\Pi^{a,b}_{f}(\nbigm)\bigr)
=\Pi^{a,b}_{f_0}\bigl(\lefttop{\tau}\nbigm\bigr)$.
By the second claim,
we obtain 
$\lefttop{\tau}\bigl(
\Pi^{a,b}_{f,\star}(\nbigm)
\bigr)
=\Pi^{a,b}_{f_0,\star}(\lefttop{\tau}\nbigm)$
$(\star=!,\ast)$.
Because the rescaling is exact,
we obtain
$\lefttop{\tau}\bigl(
\Pi^{a,b}_{f,\ast !}(\nbigm)
\bigr)
=\Pi^{a,b}_{f_0,\ast !}(\lefttop{\tau}\nbigm)$.
It particularly implies 
$\lefttop{\tau}\bigl(
\Xi^{(a)}_f(\nbigm)\bigr)
=\Xi^{(a)}_{f_0}(\lefttop{\tau}\nbigm)$
and
$\lefttop{\tau}\bigl(
\psi^{(a)}_f(\nbigm)\bigr)
=\psi^{(a)}_{f_0}(\lefttop{\tau}\nbigm)$.
By the exactness of rescaling again,
we obtain
$\lefttop{\tau}\bigl(
\phi^{(a)}_f(\nbigm)\bigr)
=\phi^{(a)}_{f_0}(\lefttop{\tau}\nbigm)$.
\hfill\qed

\subsubsection{Rescaling of $\gbigrtilde$-modules}
\label{subsection;21.6.29.20}

Let $\varphi_0:\widetilde{\proj^1_{\lambda}\times\cnum_{\tau}}
\lrarr\proj^1_{\lambda}\times\cnum_{\tau}$
denote the blow up at $(\lambda,\tau)=(0,0)$.
The map $\varphi_1$ extends to a morphism
$\widetilde{\proj^1_{\lambda}\times\cnum_{\tau}}
\lrarr \proj^1_{\lambda}$,
which is also denoted by $\varphi_1$.

We set
$\lefttop{\tau}\gbigx:=\proj^1_{\lambda}\times\lefttop{\tau}X
=\proj^1_{\lambda}\times\cnum_{\tau}\times X$
and
$\widetilde{\lefttop{\tau}\gbigx}:=
\widetilde{(\proj^1_{\lambda}\times\cnum_{\tau})}\times X$.
We obtain the induced morphism
$\varphi_1:\widetilde{\lefttop{\tau}\gbigx}\lrarr\gbigx$.
We also have the natural morphism
$\varphi_0:\widetilde{\lefttop{\tau}\gbigx}
\lrarr \lefttop{\tau}\gbigx$.
We set 
$\lefttop{\tau}\gbigx_0:=\proj^1_{\lambda}\times \lefttop{\tau}X_0$,
$\gbigy:=\lefttop{\tau}\gbigx_0\cup\lefttop{\tau}\gbigx^{\infty}$
and 
$\gbigytilde=\varphi_0^{-1}(\gbigy)$.
The following lemma is obvious.
\begin{lem}
If $\gbigm$ is a good $\nbigo_{\gbigx}(\ast\gbigx^{\infty})$-module
in the sense of {\rm\cite[Definition 4.22]{kashiwara_text}},
then $\varphi_1^{\ast}(\gbigm)(\ast\gbigytilde)$
is a good $\nbigo_{\widetilde{\lefttop{\tau}\gbigx}}(\ast\gbigytilde)$-module.
If $\gbign$ is a good $\nbigo_{\lefttop{\tau}\gbigx}(\ast\gbigy)$-module,
then $\varphi_0^{\ast}(\gbign)$
is a good $\nbigo_{\widetilde{\lefttop{\tau}\gbigx}}(\ast\gbigytilde)$-module.
\hfill\qed
\end{lem}

\begin{lem}
\label{lem;21.3.29.2}
For any good $\nbigo_{\widetilde{\lefttop{\tau}\gbigx}}(\ast\gbigytilde)$-module
$\gbigm_1$,
we have $R^i\varphi_{0\ast}(\gbigm_1)=0$ $(i>0)$,
and $\varphi_{0\ast}(\gbigm_1)$
is a good $\nbigo_{\lefttop{\tau}\gbigx}(\ast\gbigy)$-module.
Moreover, the natural morphism
 $\varphi_0^{\ast}\bigl(
 \varphi_{0\ast}(\gbigm_1)
 \bigr)
 \lrarr\gbigm_1$ 
is an isomorphism.
\end{lem}
\pf
For any coherent $\nbigo_{\widetilde{\lefttop{\tau}\gbigx}}$-module $G$,
we have
$R\varphi_{0\ast}\bigl(
G(\ast\gbigytilde)
\bigr)
=R\varphi_{0\ast}(G)\bigl(\ast\gbigy\bigr)$.
Hence, it is easy to see that
$R^i\varphi_{0\ast}\bigl(
G(\ast\gbigytilde)
\bigr)=0$ $(i>0)$
and 
$\varphi_{0\ast}\bigl(
G(\ast\gbigytilde)
\bigr)$
is a good $\nbigo_{\lefttop{\tau}\gbigx}(\ast\gbigy)$-module.
The restriction of
\begin{equation}
\label{eq;21.3.29.1}
 \varphi_{0}^{\ast}\varphi_{0\ast}(G(\ast\gbigytilde))
\lrarr
 G(\ast\gbigytilde)
\end{equation}
to $\widetilde{\lefttop{\tau}\gbigx}\setminus\gbigytilde$
is an isomorphism,
and both
$\nbigo_{\widetilde{\lefttop{\tau}\gbigx}}(\ast\gbigytilde)$-modules
$\varphi_{0}^{\ast}\varphi_{0\ast}(G(\gbigytilde))$
and
$G(\ast\gbigytilde)$
are good.
Hence, (\ref{eq;21.3.29.1}) is an isomorphism.

Because $\gbigm_1$ is good,
there exists a directed family of coherent
$\nbigo_{\widetilde{\lefttop{\tau}\gbigx}}$-submodules
$G_i$ $(i\in\Lambda)$ of $\gbigm_1$
such that $\sum G_i=\gbigm_1$.
Because
$\varphi_{0\ast}(\gbigm_1)
=\sum_{i\in\Lambda} \varphi_{0\ast}(G_i(\ast\gbigytilde))$,
$\varphi_{0\ast}(\gbigm_1)$ is good.
Moreover, the natural morphism
$\varphi_0^{\ast}\varphi_{0\ast}(\gbigm_1)\lrarr\gbigm_1$
is an isomorphism
as in the case of $G(\ast\gbigytilde)$.

For each $i\in\Lambda$,
there exists a free $\nbigo_{\lefttop{\tau}\gbigx}(\ast\gbigy)$-module
$F'_i$ with an epimorphism
$F'_i\lrarr \varphi_{0\ast}(G_i)(\ast\gbigy)$.
We obtain an epimorphism
$\bigoplus_{i\in\Lambda}
\varphi_{0}^{\ast}F'_i\lrarr \gbigm_1$,
and the kernel is also good (see \cite[Proposition 4.23]{kashiwara_text}).
Hence, we can construct a complex $\nbigj^{\bullet}$
of free $\nbigo_{\widetilde{\lefttop{\tau}\gbigx}}(\ast\gbigytilde)$-modules
$\nbigj^{\bullet}$
with a quasi-isomorphism $\nbigj^{\bullet}\lrarr\gbigm_1$
such that $\nbigj^k=0$ $(k>0)$.
Because
$R^i\varphi_{0\ast}\nbigj^k=0$ $(i>0)$,
we obtain 
$R^i\varphi_{0\ast}(\gbigm_1)=0$ $(i>0)$.
\hfill\qed

\begin{cor}
The functors $\varphi_{0\ast}$ and $\varphi_0^{\ast}$
induce equivalences of the categories of
good $\nbigo_{\widetilde{\lefttop{\tau}\gbigx}}(\ast\gbigytilde)$-modules
and
good $\nbigo_{\lefttop{\tau}\gbigx}(\ast\gbigy)$-modules,
which are mutually quasi-inverse.
\hfill\qed
\end{cor}

\vspace{.1in}

For an $\nbigo_{\gbigx}(\ast\gbigx^{\infty})$-module $\gbigm$,
we obtain an
$\nbigo_{\lefttop{\tau}\gbigx}\Bigl(
\ast\gbigy
 \Bigr)$-module
\[
\lefttop{\tau}\gbigm:=
 \varphi_{0\ast}\Bigl(
 \varphi_1^{\ast}(\gbigm)\bigl(\ast\gbigytilde\bigr)
 \Bigr).
\]
If $\gbigm$ is equipped with a meromorphic flat connection
$\nabla$,
then 
$\lefttop{\tau}\gbigm$ is naturally equipped with
the induced meromorphic flat connection
$\lefttop{\tau}\nabla=\varphi_{0\ast}(\varphi_1^{\ast}(\nabla))$.
By this procedure,
if $\gbigm$ is an $\gbigrtilde_X$-module
then 
$\lefttop{\tau}\gbigm$ is naturally an 
$\gbigrtilde_{\lefttop{\tau}X}(\ast\tau)$-module.

\begin{lem}
\label{lem;22.7.29.1}
Let $\gbigm\in\gbigc(X)$.
\begin{itemize}
\item $\lefttop{\tau}\gbigm$ is coherent
      over $\gbigrtilde_{\lefttop{\tau}X}(\ast\tau)$.
      If $\gbigm\in\gbigc_{\Malg}(X)$,
      then $\lefttop{\tau}\gbigm$
      is coherent over $\pi^{\ast}\gbigr_X(\ast\tau)$,
      where $\pi:\lefttop{\tau}X\to X$
      denotes the projection.
 \item For any $\alpha\in\cnum$,
       the multiplication by $\lambda-\alpha\tau$
       on $\lefttop{\tau}\gbigm$
       is a monomorphism,
       and hence we have
$L\iota_{\lambda=\alpha\tau}^{\ast}(\lefttop{\tau}\gbigm)
       =\iota_{\lambda=\alpha\tau}^{\ast}(\lefttop{\tau}\gbigm)$,
       where
       $\iota_{\lambda=\alpha\tau}:\nbigx\lrarr \lefttop{\tau}\nbigx$
       denotes the morphism defined by
       $(\lambda,x)\longmapsto(\alpha\lambda,\lambda,x)$. 
 \item For any $\alpha\in\cnum$,
       we set $\gbigm^{\alpha}:=\iota_{\alpha}^{\ast}\gbigm$,
       where $\iota_{\alpha}:\{\alpha\}\times X\to\nbigx$
       denotes the inclusion.
       Then, there exists a natural isomorphism
      $\iota_{\lambda=\alpha\tau}^{\ast}(\lefttop{\tau}\gbigm)
      \simeq p_X^{\ast}(\gbigm^{\alpha})(\ast\lambda)$,
       where $p_X:\nbigx\lrarr X$ denotes the projection.
 \end{itemize}
\end{lem}
\pf
Because $\varphi_1^{\ast}(\gbigm)(\ast\gbigytilde)$
is coherent over
$\varphi_1^{\ast}(\gbigrtilde_X)(\ast\gbigytilde)$,
it is coherent over
$\varphi_0^{\ast}(\gbigrtilde_{\lefttop{\tau}X}(\ast\tau))$,
and hence 
$\lefttop{\tau}\gbigm$ is coherent
over $\gbigrtilde_{\lefttop{\tau}X}(\ast\tau)$.
If $\gbigm\in\gbigc_{\Malg}(X)$,
because $\gbigm$ is coherent over $\gbigr_X$,
$\varphi_1^{\ast}(\gbigm)(\ast\gbigytilde)$
is coherent over
$\varphi_1^{\ast}(\gbigr_X)(\ast\gbigytilde)
=\varphi_0^{\ast}(\pi^{\ast}(\gbigr_X)(\ast\tau))$.
Hence, we obtain the first claim.

Because $\gbigm$ is strict,
the multiplication by $\lambda-\alpha$ on $\gbigm$
is a monomorphism.
Hence,
the multiplication by
$\lambda-\alpha\tau=-\tau(\alpha-\lambda/\tau)$ on
$\varphi_1^{\ast}\gbigm(\ast\gbigytilde)$ is
a monomorphism.
Because
$\lefttop{\tau}\gbigm
=\varphi_{0\ast}\bigl(
\varphi_1^{\ast}\gbigm(\ast\gbigytilde)
\bigr)$,
we obtain the second claim.
There exists the natural isomorphism
as in Lemma \ref{lem;21.3.29.2}:
\begin{equation}
\label{eq;21.3.15.20}
\varphi_{0}^{\ast}
\bigl(
\lefttop{\tau}\gbigm
\bigr)
\lrarr
 \varphi_1^{\ast}\gbigm
 \bigl(\ast\gbigytilde \bigr).
\end{equation}
Let $\iota'_{\lambda=\alpha\tau}:\nbigx\lrarr
\widetilde{\lefttop{\tau}\gbigx}$
denote the morphism induced by
the embedding of the proper transform of
$\{(\lambda,\tau)\,|\,\lambda=\alpha\tau\}
\subset
\proj^1_{\lambda}\times\cnum_{\tau}$.
We have
\[
 \iota_{\lambda=\alpha\tau}^{\ast}
 \bigl(\lefttop{\tau}\gbigm\bigr)
=
 (\iota'_{\lambda=\alpha\tau})^{\ast}
 \varphi_0^{\ast}
 \bigl(\lefttop{\tau}\gbigm\bigr)
 =(\iota'_{\lambda=\alpha\tau})^{\ast}
 \Bigl(
  \varphi_1^{\ast}\bigl(
  \gbigm
  \bigr)(\ast\gbigytilde)
  \Bigr)
=p_X^{\ast}(\gbigm^{\alpha})(\ast\lambda)
\]
Thus, we obtain the third claim.
\hfill\qed

\begin{rem}
\label{rem;22.7.29.2}
In Lemma {\rm\ref{lem;22.7.29.1}},
$\gbigm^1$ and $\gbigm^0$
are often denoted by $\Xi_{\DR}(\gbigm)$,
and $\Xi_{\Dol}(\gbigm)$,
respectively. 
\hfill\qed
\end{rem}

The following proposition is clear
by the constructions of the direct image and the rescaling.

\begin{prop}
\label{prop;21.4.20.1}
Let $\gbigm\in\gbigc(X)$.
For any projective morphism $F:X\lrarr Y$,
we have
$F^i_{\dagger}
 \bigl(
 \lefttop{\tau}\gbigm
 \bigr)
=
 \lefttop{\tau}(F_{\dagger}^i\gbigm)$.
\hfill\qed
\end{prop}

\begin{rem}
Let $\nu:\cnum_{\lambda}\times\cnum_{\tau}\lrarr \cnum_{\lambda}\times\cnum_{\tau}$
be defined by $\nu(\lambda,\tau)=(\lambda\tau,\tau)$.
The induced morphism $\lefttop{\tau}\nbigx\lrarr \lefttop{\tau}\nbigx$
is also denoted by $\nu$.
Although there exists a natural morphism
$\nbigo_{\lefttop{\tau}\nbigx}(\ast\lefttop{\tau}\nbigx_0)
 \lrarr
 \nu_{\ast}\nbigo_{\lefttop{\tau}\nbigx}(\ast\lefttop{\tau}\nbigx_0)$,
it is not an isomorphism.
For example,
the holomorphic function $\exp(\lambda/\tau)$ is a section of 
$\nu_{\ast}\nbigo_{\lefttop{\tau}\nbigx}(\ast\lefttop{\tau}\nbigx_0)$
on $\lefttop{\tau}\nbigx$,
but not a section of
$\nbigo_{\lefttop{\tau}\nbigx}(\ast\lefttop{\tau}\nbigx_0)$.
In particular,
any coherent $\nbigo_{\lefttop{\tau}\nbigx}(\ast\lefttop{\tau}\nbigx_0)$-module
is not naturally
a $\nu_{\ast}\nbigo_{\lefttop{\tau}\nbigx}(\ast\lefttop{\tau}\nbigx_0)$-module.

Let $\pi:\lefttop{\tau}\nbigx\lrarr \nbigx$ denote the projection.
For $\nbigm\in\nbigc(X)$,
we obtain the $\nbigrtilde_{\lefttop{\tau}X}(\ast\tau)$-module
$\nu_{\ast}(\pi^{\ast}\nbigm)$.
It is much larger than
 $\lefttop{\tau}\Upsilon(\nbigm)_{|\lefttop{\tau}\nbigx}$,
 i.e.,
 \[
 \nu_{\ast}(\pi^{\ast}\nbigm)=
 \nu_{\ast}\bigl(
 \nbigo_{\lefttop{\tau}\nbigx}(\ast\lefttop{\tau}\nbigx_0)
 \bigr)
 \otimes_{\nbigo_{\lefttop{\tau}\nbigx}(\ast\lefttop{\tau}\nbigx_0)}
 \bigl(
 \lefttop{\tau}\Upsilon(\nbigm)_{|\lefttop{\tau}\nbigx}
 \bigr).
\]
Hence, in the complex analytic setting,
$\nu_{\ast}(\pi^{\ast}\nbigm)$ 
does not seem an appropriate object
in the study of irregular Hodge filtrations.
For example,
it is not coherent over $\nbigrtilde_{\lefttop{\tau}X}(\ast\tau)$,
and it cannot be strictly specializable along $\tau$. 

\hfill\qed
\end{rem}

\subsubsection{$\cnum^{\ast}$-actions}

We consider the $\cnum^{\ast}$-action
on $\proj^1_{\lambda}\times\cnum_{\tau}$
defined by
$a(\lambda,\tau)=(a\lambda,a\tau)$.
It induces a $\cnum^{\ast}$-action
on $\widetilde{\proj^1_{\lambda}\times\cnum_{\tau}}$.
We consider the trivial $\cnum^{\ast}$-action on
$\cnum_{\lambda}$.
We also consider the trivial $\cnum^{\ast}$-action on $X$.
They induce $\cnum^{\ast}$-actions on
$\lefttop{\tau}\gbigx$,
$\widetilde{\lefttop{\tau}\gbigx}$
and $\gbigx$.
The morphisms $\varphi_0$ and $\varphi_1$
are $\cnum^{\ast}$-equivariant.

\vspace{.1in}
Let $p_2:\cnum^{\ast}\times \lefttop{\tau}\gbigx\lrarr \lefttop{\tau}\gbigx$
denote the projection.
Let $\sigma:\cnum^{\ast}\times\lefttop{\tau}\gbigx\lrarr\lefttop{\tau}\gbigx$
denote the morphism induced by the $\cnum^{\ast}$-action.
Let $m:\cnum^{\ast}\times\cnum^{\ast}\lrarr\cnum^{\ast}$
be the multiplication.
Let $p_{2,3}:\cnum^{\ast}\times\cnum^{\ast}\times
\lefttop{\tau}\gbigx\lrarr\cnum^{\ast}\times\lefttop{\tau}\gbigx$
be defined by
$p_{2,3}(a_1,a_2,x)=(a_2,x)$.

An $\nbigo_{\lefttop{\tau}\gbigx} $-module
 $\nbige$ is called $\cnum^{\ast}$-equivariant
if it is equipped with an isomorphism
$\rho:
 p_2^{\ast}(\nbige)\simeq \sigma^{\ast}(\nbige)$
satisfying the cocycle condition
$(m\times\id_{\lefttop{\tau}\gbigx})^{\ast}\rho
 =(\id_{\cnum^{\ast}}\times\sigma)^{\ast}(\rho)
 \circ
 p_{23}^{\ast}(\rho)$.
(See \cite[\S9.10 and \S11.5]{Hotta-Takeuchi-Tanisaki}.)

 \begin{lem}
For any $\nbigo_{\gbigx}(\ast\gbigx^{\infty})$-module $\gbigm$,
 $\lefttop{\tau}\gbigm$
 is naturally $\cnum^{\ast}$-equivariant.
If $\gbigm$ is an $\gbigrtilde_X$-module
with a meromorphic flat connection $\nabla$,
then we obtain
$p_2^{\ast}(\lefttop{\tau}\nabla)
  =\sigma^{\ast}(\lefttop{\tau}\nabla)$
  under the isomorphism
  $p_2^{\ast}\bigl(\lefttop{\tau}\gbigm\bigr)
  \simeq
  \sigma^{\ast}\bigl(\lefttop{\tau}\gbigm\bigr)$.
 \end{lem}
\pf
Because the morphisms $\varphi_i$ $(i=0,1)$ are $\cnum^{\ast}$-equivariant,
$\lefttop{\tau}\gbigm=
\varphi_{0\ast}\varphi_1^{\ast}(\gbigm)$ is $\cnum^{\ast}$-equivariant.
We define the morphisms
$p_2:\cnum^{\ast}\times\widetilde{\lefttop{\tau}\gbigx}
\lrarr \widetilde{\lefttop{\tau}\gbigx}$
and
$\sigma:\cnum^{\ast}\times\widetilde{\lefttop{\tau}\gbigx}
\lrarr \widetilde{\lefttop{\tau}\gbigx}$
as in the case of $\lefttop{\tau}\gbigx$.
Then, we obtain
$\sigma^{\ast}(\varphi_1^{\ast}\nabla)
=p^{\ast}(\varphi_1^{\ast}\nabla)$,
from which we obtain 
  $p_2^{\ast}\bigl(\lefttop{\tau}\nabla\bigr)
  \simeq
  \sigma^{\ast}\bigl(\lefttop{\tau}\nabla\bigr)$.
\hfill\qed

\begin{lem} 
 \label{lem;21.3.15.20}
Let $H$ be a hypersurface of $X$.
Let $\gbigm\in\gbigc(X)(\ast H)$.
If the $\nbigrtilde_{\lefttop{\tau}X(\ast\lefttop{\tau}H)}(\ast\tau)$-module
 $\lefttop{\tau}\gbigm_{|\lefttop{\tau}\nbigx}$
 is strictly specializable along $\tau$,
then we have
 $p_2^{\ast}\bigl(
 V_{\bullet}(\lefttop{\tau}\gbigm_{|\lefttop{\tau}\nbigx})
 \bigr)
 =\sigma^{\ast}\bigl(
  V_{\bullet}(\lefttop{\tau}\gbigm_{|\lefttop{\tau}\nbigx})
 \bigr)$
under 
 $p_2^{\ast}\bigl(\lefttop{\tau}\gbigm\bigr)
  _{|\cnum^{\ast}\times\lefttop{\tau}\nbigx}
  \simeq
 \sigma^{\ast}\bigl(\lefttop{\tau}\gbigm\bigr)
 _{|\cnum^{\ast}\times\lefttop{\tau}\nbigx}$.
\end{lem}
\pf
Because both 
$p_2^{\ast}\bigl(
V_{\bullet}(\lefttop{\tau}\gbigm)
\bigr)$
and
$\sigma^{\ast}\bigl(
V_{\bullet}(\lefttop{\tau}\gbigm)
\bigr)$
are the $V$-filtrations of
$p_2^{\ast}\bigl(\lefttop{\tau}\gbigm\bigr)
_{|\cnum^{\ast}\times\lefttop{\tau}\nbigx}
\simeq
\sigma^{\ast}\bigl(\lefttop{\tau}\gbigm\bigr)
_{|\cnum^{\ast}\times\lefttop{\tau}\nbigx}$,
it follows from the uniqueness of the $V$-filtration.
\hfill\qed

\subsection{Rescalable objects}

We shall prove the following theorem
in \S\ref{subsection;21.3.13.100}
after the preliminaries in
\S\ref{subsection;21.3.13.3}--\ref{subsection;21.3.13.101}.

\begin{thm}
\label{thm;21.3.12.12}
Let $\nbigm\in\nbigc(X)$.
If there exists $\nbigm_0\in\nbigc(\lefttop{\tau}X)$
such that
 $\nbigm_{0|\lefttop{\tau}X\setminus\lefttop{\tau}X_0}
 =\lefttop{\tau}\nbigm$,
then we obtain
 $\lefttop{\tau}\Upsilon(\nbigm)=
 \Upsilon(\nbigm_0)(\ast\tau)$.
Moreover, the following holds.
\begin{itemize}
 \item For any hypersurface $H^{(1)}$ of $X$,
       we have $\lefttop{\tau}\Upsilon(\nbigm[\star H^{(1)}])
       =\Upsilon(\nbigm_0[\star \lefttop{\tau}H^{(1)}])(\ast\tau)$
       $(\star=!,\ast)$.
 \item Let $g$ be any holomorphic function on $X$,
       and let $g_0$ be the induced holomorphic function on $\lefttop{\tau}X$.
       Then, for any $a\leq b$,
       we have
       $\lefttop{\tau}\Upsilon\bigl(\Pi^{a,b}_{g,\star}(\nbigm)\bigr)
       =\Upsilon(\Pi^{a,b}_{g_0,\star}\nbigm_0)(\ast\tau)$
       $(\star=!,\ast)$
       and
       $\lefttop{\tau}\Upsilon\bigl(\Pi^{a,b}_{g,\ast !}(\nbigm)\bigr)
       =\Upsilon(\Pi^{a,b}_{g_0,\ast !}\nbigm_0)(\ast\tau)$.
       Hence, we have
\[
       \lefttop{\tau}\Upsilon(\Xi^{(a)}_{g}\nbigm)
       =\Upsilon(\Xi^{(a)}_{g_0}\nbigm_0)(\ast\tau),
       \quad
       \lefttop{\tau}\Upsilon(\psi^{(a)}_{g}\nbigm)
       =\Upsilon(\psi^{(a)}_{g_0}\nbigm_0)(\ast\tau),
       \quad
       \lefttop{\tau}\Upsilon(\phi^{(a)}_{g}\nbigm)
       =\Upsilon(\phi^{(a)}_{g_0}\nbigm_0)(\ast\tau).
\]
\end{itemize}
\end{thm}

Before proving Theorem \ref{thm;21.3.12.12},
we explain some consequences of the theorem.
Let $H$ be a hypersurface of $X$.
\begin{cor}
Let $\nbigm\in\nbigc(X;H)$.
If there exists
$\nbigm_0\in\nbigc(\lefttop{\tau}X;\lefttop{\tau}H)$
such that
$\nbigm_{0|\lefttop{\tau}X\setminus\lefttop{\tau}X_0}
=\lefttop{\tau}\nbigm$,
we obtain
$\lefttop{\tau}\Upsilon(\nbigm)=
\Upsilon(\nbigm_0)(\ast\tau)$.
\end{cor}
\pf
There exists $\nbigm_1\in\nbigc(X)$
such that $\nbigm_1(\ast H)=\nbigm$
and that $\nbigm_1[\ast H]=\nbigm_1$.
We have
$(\lefttop{\tau}\nbigm_1)[\ast \lefttop{\tau}H]
=\lefttop{\tau}\nbigm_1$
and that
$(\lefttop{\tau}\nbigm_1)(\ast \lefttop{\tau}H)
=\lefttop{\tau}\nbigm=
\nbigm_{0|\lefttop{\tau}X\setminus\lefttop{\tau}X_0}$.
We set
$\nbigm_2:=\nbigm_0[\ast\lefttop{\tau}H]$.
Then, we obtain
$\nbigm_{2|\lefttop{\tau}X\setminus \lefttop{\tau}X_0}
=\lefttop{\tau}\nbigm_1$.
By Theorem \ref{thm;21.3.12.12},
we obtain
$\Upsilon(\nbigm_2)(\ast\tau)
=\lefttop{\tau}\Upsilon(\nbigm_1)$,
which implies
$\lefttop{\tau}\Upsilon(\nbigm)=
\Upsilon(\nbigm_0)(\ast\tau)$.
\hfill\qed

\vspace{.1in}
We introduce the categories of rescalable objects.

\begin{df}  
Let $\nbigc_{\res}(X)\subset\nbigc(X)$
denote the full subcategory of
$\nbigm\in\nbigc(X)$ for which 
there exists $\nbigm_0\in\nbigc(\lefttop{\tau}X)$
such that
$\nbigm_{0|\lefttop{\tau}X\setminus\lefttop{\tau}X_0}
 =\lefttop{\tau}\nbigm$.
Let $\nbigc_{\res}(X;H)$
denote the essential image of
$\nbigc_{\res}(X)\lrarr\nbigc(X;H)$. 
Objects in $\nbigc_{\res}(X;H)$ are called
rescalable objects in $\nbigc(X;H)$.
\hfill\qed
\end{df}

\begin{df}
Let $\gbigc_{\res}(X)\subset\gbigc_{\Malg}(X)$
denote the full subcategory of
objects $\gbigm\in\gbigc_{\Malg}(X)$ such that
$\lefttop{\tau}\gbigm\in \gbigc(\lefttop{\tau}X;\lefttop{\tau}X_0)$.
Let $\gbigc_{\res}(X;H)$
denote the essential image of
$\gbigc_{\res}(X)\lrarr\gbigc(X;H)$.
Objects in  $\gbigc_{\res}(X;H)$
are called rescalable objects in $\gbigc(X;H)$.
\hfill\qed
\end{df}

The following lemma is clear by the theorem.
\begin{lem}
The restriction $\gbigm\longmapsto\gbigm_{|\nbigx}$
induces an equivalence
 $\gbigc_{\res}(X;H)
 \simeq
 \nbigc_{\res}(X;H)$.
 A quasi-inverse is given by
$\nbigm\longmapsto \Upsilon(\nbigm)$. 
\hfill\qed
\end{lem}

The categories
$\nbigc(\lefttop{\tau}X;\lefttop{\tau}X_0\cup \lefttop{\tau}H)$
and
$\gbigc(\lefttop{\tau}X;\lefttop{\tau}X_0\cup \lefttop{\tau}H)$
are also denoted by $\nbigc(\lefttop{\tau}X;\lefttop{\tau}H)(\ast\tau)$
and $\gbigc(\lefttop{\tau}X;\lefttop{\tau}H)(\ast\tau)$,
respectively.
We obtain the following corollary
by Theorem \ref{thm;21.3.12.12}.
\begin{cor}
For $\nbigm\in\nbigc(X;H)$,
the following three conditions are equivalent.
\begin{itemize}
 \item $\nbigm\in\nbigc_{\res}(X;H)$.
 \item $\lefttop{\tau}\Upsilon(\nbigm)_{|\lefttop{\tau}\nbigx}
       \in\nbigc(\lefttop{\tau}X;\lefttop{\tau}H)(\ast\tau)$.
 \item $\lefttop{\tau}\Upsilon(\nbigm)
       \in\gbigc(\lefttop{\tau}X;\lefttop{\tau}H)(\ast\tau)$.
\end{itemize}
For $\gbigm\in\gbigc_{\Malg}(X;H)$,
the following three conditions are equivalent.
\begin{itemize}
 \item $\gbigm\in\gbigc_{\res}(X;H)$.
 \item $\lefttop{\tau}\gbigm_{|\lefttop{\tau}\nbigx}
       \in\nbigc(\lefttop{\tau}X;\lefttop{\tau}H)(\ast\tau)$.
 \item $\lefttop{\tau}\gbigm
       \in\gbigc_{\Malg}(\lefttop{\tau}X;\lefttop{\tau}H)(\ast\tau)$.
       \hfill\qed
\end{itemize}
\end{cor}

\begin{example}
\label{example;22.7.27.11}
The object
$\gbigr_F(M)\otimes\gbigl(f)\in\gbigc_{\Malg}(X;H)$
in {\rm\S\ref{subsection;22.7.27.10}}
is contained in $\gbigc_{\res}(X;H)$.
Indeed,
we have
$\lefttop{\tau}\bigl(
 \gbigr_F(M)\otimes\gbigl(f)
 \bigr)
 =\bigl(
 \pi^{\ast}(\gbigr_F(M))\otimes\gbigl(\tau f)
 \bigr)(\ast\tau)
\in\gbigc_{\Malg}(\lefttop{\tau}X;\lefttop{\tau}H)$,
where $\pi:\lefttop{\tau}X\to X$ denotes the projection.
\hfill\qed
\end{example}

We note that the following lemma is implied in 
Theorem \ref{thm;21.3.12.12}.
\begin{lem}
For $\gbigm\in\gbigc_{\res}(X;H)$,
$\lefttop{\tau}\gbigm$
is the Malgrange extension of
$\lefttop{\tau}\gbigm_{|\lefttop{\tau}\nbigx}$.
\hfill\qed
\end{lem}

\subsubsection{Smooth case}
\label{subsection;21.3.13.3}

Let us study the first claim of
Theorem \ref{thm;21.3.12.12}
in the case
where $\nbigm$ is a smooth $\nbigrtilde_{X}$-module
with $H=\emptyset$.
We obtain the holomorphic vector bundle $E=\nbigm/\lambda\nbigm$ on $X$
with the Higgs field $\theta$ induced by
the action of $\nbigr_X$.
Let $\rho$ denote the endomorphism of
the Higgs bundle $(E,\theta)$
induced by the action of $\lambda^2\del_{\lambda}$.
We also obtain the holomorphic vector bundle
$\lefttop{\tau}E=\lefttop{\tau}\nbigm\big/\lambda\lefttop{\tau}\nbigm$
on $\lefttop{\tau}X\setminus\lefttop{\tau}X_0$
with the Higgs field $\lefttop{\tau}\theta$
induced by $\nbigr_{\lefttop{\tau}X\setminus\lefttop{\tau}X_0}$-action.
Let $\pi:\lefttop{\tau}X\setminus\lefttop{\tau}X_0\lrarr X$
denote the projection.
There exists a natural isomorphism
$\lefttop{\tau}E\simeq \pi^{\ast}(E)$
under which we have
$\lefttop{\tau}\theta=\tau\cdot\pi^{\ast}\theta
-\pi^{\ast}(\rho)d\tau$.

There exists a weight filtration $\nbigw$ of $\nbigm_0$
such that each $\Gr^{\nbigw}_j(\nbigm_0)$
underlies an integrable
polarizable pure twistor $\nbigd$-module of weight $j$.
It induces a filtration $\nbigw$ on
$\nbigm_{0|\lefttop{\tau}X\setminus\lefttop{\tau}X_0}
=\lefttop{\tau}\nbigm$
by smooth $\nbigrtilde_{\lefttop{\tau}X\setminus \lefttop{\tau}X_0}$-modules,
i.e., locally free $\nbigo_{\nbigx}$-modules.
As the restriction,
we obtain a filtration $W$ of
the Higgs bundle $(\lefttop{\tau}E,\lefttop{\tau}\theta)$
by Higgs subbundles.
Each Higgs bundle
$\Gr^{\nbigw}_j(\lefttop{\tau}E,\lefttop{\tau}\theta)$
has a pluri-harmonic metric $h_j$
which induces
the $\nbigr_{\lefttop{\tau}X\setminus\lefttop{\tau}X_0}$-module
$\Gr^{\nbigw}_j(\nbigm_{0|\lefttop{\tau}X\setminus\lefttop{\tau}X_0})$.

\begin{lem}
\label{lem;21.3.13.1}
 For any $\lambda\neq 0$,
the isomorphism
$\nbigm_{0|\{\lambda\}\times(\lefttop{\tau}X\setminus\lefttop{\tau}X_0)}
 \simeq
 \lefttop{\tau}\nbigm_{|\{\lambda\}\times(\lefttop{\tau}X\setminus\lefttop{\tau}X_0)}$
extends to
\[
 \Upsilon(\nbigm_0)(\ast\tau)_{|\{\lambda\}\times\lefttop{\tau}X}
\simeq
 \lefttop{\tau}\bigl(
  \Upsilon(\nbigm)
  \bigr)
  _{|\{\lambda\}\times\lefttop{\tau}X}. 
\]
\end{lem}
\pf
Because
$\lefttop{\tau}\theta=
\tau\cdot\pi^{\ast}\theta
-\pi^{\ast}(\rho)d\tau$,
the harmonic bundles
$(\Gr^{\nbigw}_j(\lefttop{\tau}E,\lefttop{\tau}\theta),h_j)$
are tame along $\tau=0$.
It implies that
$\Gr^{\nbigw}_j(\nbigm_0)(\ast\tau)_{|\{\lambda\}\times \lefttop{\tau}X}$
are regular singular meromorphic flat bundles
on $(\lefttop{\tau}X,\lefttop{\tau}X_0)$,
and hence
$\nbigm_0(\ast\tau)_{|\{\lambda\}\times\lefttop{\tau}X}$
is regular singular.
Because
$\lefttop{\tau}\bigl(
\Upsilon(\nbigm)\bigr)_{|\{\lambda\}\times\lefttop{\tau}X}$
is also regular singular,
we obtain the claim of Lemma \ref{lem;21.3.13.1}.
\hfill\qed

\begin{lem}
We obtain
$\Upsilon(\nbigm_0)(\ast\tau)_{|\lefttop{\tau}\nbigx}=
\lefttop{\tau}\bigl(
\Upsilon(\nbigm)
\bigr)_{|\lefttop{\tau}\nbigx}$.
\end{lem}
\pf
Both
$\Upsilon(\nbigm_0)(\ast\tau)_{|\lefttop{\tau}\nbigx}$
and
$\lefttop{\tau}\bigl(
\Upsilon(\nbigm)
\bigr)_{|\lefttop{\tau}\nbigx}$
are locally free
$\nbigo_{\lefttop{\tau}\nbigx}(\ast\lefttop{\tau}\nbigx_0)$-modules.
Let $\vecv_1$ and $\vecv_2$ be local frames
of 
$\Upsilon(\nbigm_0)(\ast\tau)_{|\lefttop{\tau}\nbigx}$
and
$\lefttop{\tau}\bigl(
\Upsilon(\nbigm)
\bigr)_{|\lefttop{\tau}\nbigx}$
on an open subset $\nbigu$ of $\lefttop{\tau}\nbigx$.
There exists a matrix $\nbiga$ determined by
$\vecv_{1|\nbigu\setminus\lefttop{\tau}\nbigx_0}
=\vecv_{2|\nbigu\setminus\lefttop{\tau}\nbigx_0}$.
Each entry $\nbiga_{i,j}$ of $\nbiga$ is a holomorphic function on
$\nbigu\setminus\lefttop{\tau}\nbigx_0$.
For any $\lambda\neq 0$,
the restriction of $\nbiga_{i,j}$ to
$(\{\lambda\}\times\lefttop{\tau}X)\cap\nbigu$ are meromorphic
along $\tau=0$.
Hence, we obtain that $\nbiga_{i,j}$ are meromorphic along $\tau=0$.
It implies the claim of Lemma 
\hfill\qed

\vspace{.1in}

Because the restrictions of both
$\Upsilon(\nbigm_0)$
and
$\lefttop{\tau}\Upsilon(\nbigm)$
to $\proj^1\times(\lefttop{\tau}X\setminus\lefttop{\tau}X_0)$
are regular singular,
we obtain the isomorphism
$\Upsilon(\nbigm_0)(\ast\tau)
\simeq
\lefttop{\tau}\bigl(
\Upsilon(\nbigm)
\bigr)$
on 
$\lefttop{\tau}\gbigx\setminus
(\{\infty\}\times \lefttop{\tau}X_0)$.
It extends to an isomorphism
on $\lefttop{\tau}\gbigx$
by the Hartogs theorem.
Thus, we obtain the first claim of Theorem \ref{thm;21.3.12.12}
in the case where $\nbigm$ is a smooth
$\nbigrtilde_{X}$-module.

\subsubsection{Localizations}
\label{subsection;21.3.13.101}

Let $H$ be a hypersurface of $X$.
Suppose that there exist
$\nbign\in\nbigc(X)$ and
$\nbign_0\in\nbigc(\lefttop{\tau}X)$
such that
$\lefttop{\tau}\bigl(
 \Upsilon(\nbign)(\ast H)
 \bigr)
 =\Upsilon(\nbign_0)\bigl(\ast (\lefttop{\tau}X_0\cup\lefttop{\tau}H)\bigr)$.
\begin{lem}
\label{lem;21.3.13.50}
 We have
 $\lefttop{\tau}\Upsilon\bigl(
 \nbign[\star H]
 \bigr)
 =\Upsilon(\nbign_0[\star \lefttop{\tau}H])(\ast\tau)$
 for $\star=!,\ast$
 as $\nbigrtilde_{\lefttop{\tau}X}(\ast\tau)$-modules.
\end{lem}
\pf
Let us explain the argument in the case $\star=!$.
The other case can be argued similarly.
We may assume
$\nbign=\nbign[!H]$
and
$\nbign_0=\nbign_0[!\lefttop{\tau}H]$.

It is enough to prove the claim locally around any point of $H$.
Hence, we may assume that there exists
a holomorphic function $g$ on $X$
such that $H=g^{-1}(0)$.
Let $g_0$ denote the induced holomorphic function
on $\lefttop{\tau}X$.
It is enough to consider the case where
$g$ is a coordinate function,
i.e.,
$X$ is an open subset of $X'\times\cnum_t$,
and $g=t$.
The function $g_0$ is also denoted by $t$.
Let $V\gbigrtilde_{X}\subset
\gbigrtilde_{X}$ denote the sheaf of subalgebras
generated by
$\lambda \Theta_{X}(\log t)$
and $\lambda^2\del_{\lambda}$.
We use the notation $V\gbigrtilde_{\lefttop{\tau}X}$
in a similar meaning.

The $\gbigrtilde_{\lefttop{\tau}X}$-module
$\Upsilon(\nbign_0)$
has the $V$-filtration
$V_{\bullet}(\Upsilon(\nbign_0))$ along $t$.
We have
\[
\Upsilon(\nbign_0)(\ast \tau)
=
\gbigrtilde_{\lefttop{\tau}X}(\ast\tau)
\otimes_{V\gbigrtilde_{\lefttop{\tau}X}(\ast\tau)}
V_{<0}\bigl(
\Upsilon(\nbign_0)
\bigr)(\ast \tau). 
\]
Note that
$V_{\bullet}(\Upsilon(\nbign_0))(\ast \tau)$
is the $V$-filtration
of the $\gbigrtilde_{\lefttop{\tau}X}(\ast\tau)$-module
$\Upsilon(\nbign_0)(\ast\tau)$ along $t$.

The $\gbigrtilde_{X}$-module
$\Upsilon(\nbign)$
has the $V$-filtration
$V_{\bullet}(\Upsilon(\nbign))$ along $t$.
We obtain
the $\nbigo_{\lefttop{\tau}\gbigx}(\ast\tau)$-submodules
$\lefttop{\tau}V_{a}(\Upsilon(\nbign))$
of $\lefttop{\tau}\Upsilon(\nbign)$.
We have
\[
\lefttop{\tau}\Upsilon(\nbign)
=
\gbigrtilde_{\lefttop{\tau}X}(\ast\tau)
\otimes_{V\gbigrtilde_{\lefttop{\tau}X}(\ast\tau)}
\lefttop{\tau}
\bigl(
V_{<0}\Upsilon(\nbign)
\bigr).
\]
Note that 
$\lefttop{\tau}V_{\bullet}(\Upsilon(\nbign))$
is the $V$-filtration of
the $\gbigrtilde_{\lefttop{\tau}X}(\ast\tau)$-module
$\lefttop{\tau}\Upsilon(\nbign)$.

Because
$\Bigl(
\lefttop{\tau}\Upsilon\bigl(
 \nbign
 \bigr)
 \Bigr)(\ast t)
=\lefttop{\tau}\Bigl(
 \Upsilon(\nbign)(\ast t)
 \Bigr)
=\Upsilon(\nbign_0)\bigl(\ast (\tau t)\bigr)
=\Bigl(
 \Upsilon(\nbign_0)(\ast\tau)
 \Bigr)(\ast t)$
as an $\gbigrtilde_{\lefttop{\tau}X}(\ast(\tau t))$-module,
we obtain
\[
V_{<0}\bigl(
\Upsilon(\nbign_0)(\ast \tau)
\bigr)
=
V_{<0}\bigl(
\bigl(\Upsilon(\nbign_0)(\ast \tau)\bigr)
(\ast t)
\bigr)
=
 V_{<0}\bigl(
\lefttop{\tau}\Upsilon(\nbign)(\ast t)
\bigr)
= V_{<0}\bigl(
\lefttop{\tau}\Upsilon(\nbign)
\bigr).
\]
Thus, we obtain Lemma \ref{lem;21.3.13.50}.
\hfill\qed

\subsubsection{Proof of Theorem \ref{thm;21.3.12.12}}
\label{subsection;21.3.13.100}

It is enough to check the claim locally around any point of $X$.
Therefore, we may assume the existence of
$f$, $Z$, $\varphi$ and $\nbigv$
as in the proof of Proposition \ref{prop;21.2.18.1}.
We obtain the description of
$\lefttop{\tau}\Upsilon(\nbigm)$
as the cohomology of the following complex of
$\gbigrtilde_{\lefttop{\tau}X}(\ast \tau)$-modules:
\begin{equation}
\label{eq;21.3.12.4}
\lefttop{\tau}\bigl(
 \Upsilon(\psi_f^{(1)}(\nbigm))
 \bigr)
 \lrarr
 \lefttop{\tau}\bigl(
 \Upsilon(\Xi_f^{(0)}(\nbigm))
 \bigr)
 \oplus
  \lefttop{\tau}\bigl(
  \Upsilon(\phi_f^{(0)}(\nbigm))
  \bigr)
  \lrarr
    \lefttop{\tau}\bigl(
    \Upsilon(\psi_f^{(0)}(\nbigm))
    \bigr).
\end{equation}
We set $f_0=f\circ\pi$.
We can describe $\Upsilon(\nbigm_0)(\ast\tau)$
as the cohomology of the following complex of
$\gbigrtilde_{\lefttop{\tau}X}(\ast \tau)$-modules:
\begin{equation}
\label{eq;21.3.13.2}
 \Upsilon(\psi_{f_0}^{(1)}(\nbigm_0))(\ast\tau)
 \lrarr
 \Upsilon(\Xi_{f_0}^{(0)}(\nbigm_0))(\ast\tau)
 \oplus
  \Upsilon(\phi_{f_0}^{(0)}(\nbigm_0))(\ast\tau)
  \lrarr
 \Upsilon(\psi_{f_0}^{(0)}(\nbigm_0))(\ast\tau).
\end{equation}

By Lemma \ref{lem;21.3.13.30},
there exist the following natural identifications
\[
\lefttop{\tau}\bigl(
\psi^{(a)}_f(\nbigm)\bigr)
=\psi^{(a)}_{f_0}(\lefttop{\tau}\nbigm)
=\psi^{(a)}_{f_0}(\nbigm_0)
_{|\lefttop{\tau}X\setminus\lefttop{\tau}X_0},
\]
\[
\lefttop{\tau}\bigl(
\phi^{(a)}_f(\nbigm)\bigr)
=\phi^{(a)}_{f_0}(\lefttop{\tau}\nbigm)
=\phi^{(a)}_{f_0}(\nbigm_0)
_{|\lefttop{\tau}X\setminus\lefttop{\tau}X_0},
\]
\[
\lefttop{\tau}\bigl(
 \Xi^{(a)}_f(\nbigm)\bigr)
=\Xi^{(a)}_{f_0}(\lefttop{\tau}\nbigm)
=\Xi^{(a)}_{f_0}(\nbigm_0)
_{|\lefttop{\tau}X\setminus\lefttop{\tau}X_0}.
\]
By the assumption of the induction,
we obtain
\[
 \lefttop{\tau}\Upsilon(\phi^{(a)}_f(\nbigm))
 =\Upsilon(\phi^{(a)}_{f_0}(\nbigm_0))(\ast\tau),\quad
 \lefttop{\tau}\Upsilon(\psi^{(a)}_f(\nbigm))
=\Upsilon(\psi^{(a)}_{f_0}(\nbigm_0))(\ast\tau).
\]
By the description of 
$\lefttop{\tau}\Upsilon(\nbigm)$
and $\Upsilon(\nbigm_0)(\ast\tau)$
as the cohomology of (\ref{eq;21.3.12.4})
and (\ref{eq;21.3.13.2}), respectively,
the proof of the first claim of Theorem \ref{thm;21.3.12.12}
is reduced to the following lemma.
\begin{lem}
\label{lem;21.3.12.11}
 $\lefttop{\tau}\Upsilon(\Xi^{(a)}_f(\nbigm))
 =\Upsilon(\Xi^{(a)}_{f_0}(\nbigm_0))(\ast\tau)$. 
\end{lem}
\pf
We set $W:=\varphi(Z)$.
According to \cite[Theorem 2.0.2]{Wlodarczyk},
there exists a projective morphism
$\rho:X'\lrarr X$
such that the following holds.
\begin{itemize}
 \item $\rho^{-1}(f^{-1}(0))$ is normal crossing in $X'$.
 \item The proper transform $Y$ of $W$ is smooth,
       and $Y$ intersects with $\rho^{-1}(f^{-1}(0))$
       in the normal crossing way.
\end{itemize}
We set $f':=f\circ\rho$. 
By the independence from the compactification in
\cite[Proposition 11.2.12]{Mochizuki-MTM},
there uniquely exists $\nbigm'\in\nbigc(X')$
such that
$\nbigm'=\nbigm'[\ast f']$
and that
$\rho_{\dagger}\nbigm'=\nbigm[\ast f]$.
The support of $\nbigm'$ is $Y$.
We set $H_Y:=Y\cap (f')^{-1}(0)$.
By Kashiwara equivalence 
for integrable mixed twistor $\nbigd$-modules
\cite[Proposition 7.2.8]{Mochizuki-MTM},
there uniquely exists
$\nbigm_{Y}\in\nbigc(Y)$
such that $\nbigm'$ is the direct image of $\nbigm_Y$.
Let $\rho_{Y}:Y\lrarr X$ be the induced morphism,
and we set $f_{Y}:=f\circ\rho_{Y}$.
We have
$\nbigm_{Y}[\ast f_{Y}]=\nbigm_{Y}$
and 
$\rho_{Y\dagger}(\nbigm_{Y})
=\nbigm[\ast f]$.

We can apply a similar procedure to $\nbigm_0$.
Namely, let $\rho_0:\lefttop{\tau}X'\lrarr \lefttop{\tau}X$
denote the induced morphism.
We set $f'_0:=f_0\circ\rho_0$.
There uniquely exists
$\nbigm'_0\in\nbigc(\lefttop{\tau}X')$
such that
$\nbigm'_0=\nbigm'_0[\ast f_0']$
and
$\rho_{0\dagger}(\nbigm'_0)=\nbigm_0[\ast f_0]$.
There uniquely exists
$\nbigm_{0,\lefttop{\tau}Y}\in\nbigc(\lefttop{\tau}Y)
\in\nbigc(\lefttop{\tau}Y)$
such that $\nbigm'_0$
is the direct image of
$\nbigm_{0,\lefttop{\tau}Y}$.
Let $\rho_{\lefttop{\tau}Y}:\lefttop{\tau}Y\lrarr \lefttop{\tau}X$
denote the induced morphism,
and we set
$f_{0,\lefttop{\tau}Y}:=f_0\circ \rho_{\lefttop{\tau}Y}$.
We have
$\nbigm_{0,\lefttop{\tau}Y}[\ast f_{0,\lefttop{\tau}Y}]
=\nbigm_{0,\lefttop{\tau}Y}$
and
$\rho_{\lefttop{\tau}Y\dagger}(\nbigm_{0,\lefttop{\tau}Y})
=\nbigm_0[\star f_0]$.

\begin{lem}
\label{lem;21.3.13.40}
 $\lefttop{\tau}\nbigm'
=\nbigm'_{0|\lefttop{\tau}X'\setminus\lefttop{\tau}X'_0}$.
\end{lem}
\pf
By Lemma \ref{lem;21.3.13.30},
we have
$\lefttop{\tau}(\nbigm[\ast f])
=(\lefttop{\tau}\nbigm)[\ast f_0]
=(\nbigm_0[\ast f_0])_{|\lefttop{\tau}X\setminus\lefttop{\tau}X_0}$.
Let $\rho_1:\lefttop{\tau}X'\setminus\lefttop{\tau}X'_0
\lrarr\lefttop{\tau}X\setminus\lefttop{\tau}X_0$
denote the restriction of $\rho_0$.
Because
$\rho_{1\dagger}(\nbigm'_{0|\lefttop{\tau}X'\setminus\lefttop{\tau}X'_0})
=
\lefttop{\tau}(\nbigm[\ast f])
=\rho_{1\dagger}(\lefttop{\tau}\nbigm')$,
we obtain
the claim of Lemma \ref{lem;21.3.13.40}.
\hfill\qed

\vspace{.1in}

By Lemma \ref{lem;21.3.13.40},
we obtain
$(\lefttop{\tau}\nbigm_{Y})_{|\lefttop{\tau}Y\setminus\lefttop{\tau}Y_0}
=(\nbigm_{0,\lefttop{\tau}Y})_{|\lefttop{\tau}Y\setminus\lefttop{\tau}Y_0}$.
Hence, we have
\begin{equation}
\label{eq;21.3.12.30}
 \lefttop{\tau}\Upsilon\bigl(\nbigm_{Y}(\ast f_{Y})\bigr)
  _{| \lefttop{\tau}\gbigy
\setminus
\bigl(
(\lefttop{\tau}\gbigy)^{\infty}
\cup
\lefttop{\tau}\gbigy_0
\bigr)
}
\simeq
\Upsilon(\nbigm_{0,\lefttop{\tau}Y})\bigl(\ast (\tau f_{0,\lefttop{\tau}Y})\bigr)
  _{| \lefttop{\tau}\gbigy
\setminus
\bigl(
(\lefttop{\tau}\gbigy)^{\infty}
\cup
\lefttop{\tau}\gbigy_0
\bigr)
}.
\end{equation}
By using the result
in the smooth case (see \S\ref{subsection;21.3.13.3}),
we obtain
\begin{equation}
 \label{eq;21.3.12.31}
 \lefttop{\tau}\Upsilon\bigl(\nbigm_{Y}(\ast f_{Y})\bigr)
  _{| \lefttop{\tau}\gbigy
\setminus
(\lefttop{\tau}\gbigh_{Y})^{\infty}
}
=
\Upsilon(\nbigm_{0,\lefttop{\tau}Y})\bigl(\ast (\tau f_{0,\lefttop{\tau}Y})\bigr)
  _{| \lefttop{\tau}\gbigy
\setminus
(\lefttop{\tau}\gbigh_{Y})^{\infty}
}.
\end{equation}
By the Hartogs theorem,
we obtain
\begin{equation}
 \label{eq;21.3.12.32}
 \lefttop{\tau}\Upsilon\bigl(\nbigm_{Y}(\ast f_{Y})\bigr)
=\Upsilon(\nbigm_{0,\lefttop{\tau}Y})\bigl(\ast (\tau f_{0,\lefttop{\tau}Y})\bigr).
\end{equation}
Therefore, we have
\[
\lefttop{\tau}\Bigl(
\Upsilon\bigl(
\Pi^{a,b}_{f_{Y}}(\nbigm_{Y})
\bigr)
\Bigr)
=
\Pi^{a,b}_{f_{0,\lefttop{\tau}Y}}
\Bigl(
\lefttop{\tau}\Upsilon(\nbigm_Y(\ast f_{Y}))
\Bigr)
=
\Pi^{a,b}_{f_{0,\lefttop{\tau}Y}}
\Bigl(
\Upsilon(\nbigm_{0,\lefttop{\tau}Y})\bigl(\ast(\tau f_{0,\lefttop{\tau}Y})\bigr)
\Bigr)
=
\Upsilon\Bigl(
\Pi^{a,b}_{f_{0,\lefttop{\tau}Y}}(\nbigm_{0,\lefttop{\tau}Y})
\Bigr)(\ast\tau).
\]
By Lemma \ref{lem;21.3.13.50},
we obtain
$\lefttop{\tau}
\Upsilon(\Pi^{a,b}_{f_Y,\star}\nbigm_Y)
=
\Upsilon\Bigl(
\Pi^{a,b}_{f_{0,\lefttop{\tau}Y},\star}(\nbigm_{0,\lefttop{\tau}Y})
\Bigr)(\ast\tau)$
for $\star=!,\ast$.
We obtain
\begin{equation}
\label{eq;21.3.12.100}
   \lefttop{\tau}\Bigl(
  \Xi_{f_{Y}}^{(a)}(\Upsilon(\nbigm_{Y}))
  \Bigr)
  \simeq
    \Bigl(
 \Upsilon\Xi^{(a)}_{f_{0,\lefttop{\tau}Y}}(\nbigm_{0,\lefttop{\tau}Y})
 \Bigr) (\ast\tau).
\end{equation}
As the direct image of
(\ref{eq;21.3.12.100}) via $\rho_{\lefttop{\tau}Y}$,
we obtain the desired isomorphism in Lemma \ref{lem;21.3.12.11}.
Thus, we obtain Lemma \ref{lem;21.3.12.11}
and the first claim of Theorem \ref{thm;21.3.12.12}.
\hfill\qed

\vspace{.1in}

As for the second claim of Theorem \ref{thm;21.3.12.12}
it is enough to check the claim locally around any point of $P\in X$.
Hence, we may assume that there exists a holomorphic function $g$
on $X$ such that $H=g^{-1}(0)$.
Then, the second claim of Theorem \ref{thm;21.3.12.12}
is reduced to Lemma \ref{lem;21.3.13.50}.
The third claim of Theorem \ref{thm;21.3.12.12}
is also reduced to Lemma \ref{lem;21.3.13.50}.
\hfill\qed

\subsection{Basic functorial properties}

\subsubsection{Direct image}

Let $F:X\lrarr Y$ be a projective morphism.
Let $H_Y$ be a hypersurface of $Y$.
We set $H_X:=F^{-1}(H_Y)$.

\begin{prop}
\label{prop;21.6.29.21}
The functors $F_{\dagger}^j:\gbigc(X;H_X)\lrarr \gbigc(Y;H_Y)$
induce 
$F_{\dagger}^j:\gbigc_{\res}(X;H_X)\lrarr\gbigc_{\res}(Y;H_Y)$. 
\begin{itemize}
 \item Suppose moreover that $F_{|X\setminus H_X}$ is a closed embedding of
       $X\setminus H_X$ into $Y\setminus H_Y$.
       We set
       $\gbigc_{F(X),\res}(Y)=
       \gbigc_{\res}(Y)\cap\gbigc_{F(X)}(Y)$.
       Then,
       $F_{\dagger}^0$ induces
       an equivalence
       $\gbigc_{\res}(X;H_X)\simeq
       \gbigc_{F(X),\res}(Y;H_Y)$.
\end{itemize}
\end{prop}
\pf
Let $\lefttop{\tau}F:\lefttop{\tau}X\lrarr\lefttop{\tau}Y$
denote the induced morphism.
Let $\gbigm\in\gbigc_{\res}(X;H_X)$.
Because 
$\lefttop{\tau}\gbigm
\in \gbigc(\lefttop{\tau}X;\lefttop{\tau}H)(\ast\tau)$,
we obtain
$(\lefttop{\tau}F)_{\dagger}^i\bigl(
\lefttop{\tau}\gbigm
\bigr)
\in \gbigc(\lefttop{\tau}Y;\lefttop{\tau}H_Y)(\ast\tau)$.
According to Proposition \ref{prop;21.4.20.1}
and Theorem \ref{thm;21.3.12.12}
$F_{\dagger}^i(\gbigm)$ is an object of
$\gbigc_{\res}(Y;H_Y)$.
Thus, we obtain the first claim.

Suppose moreover that
$F_{|X\setminus H_X}$ is a closed embedding of
$X\setminus H_X$ into $Y\setminus H_Y$.
Because
$F^0_{\dagger}:\gbigc(X;H_X)\simeq
\gbigc_{F(X)}(Y;H_Y)$,
it is enough to prove the essential surjectivity.
Let $\gbigm_Y\in\gbigc_{F(X),\res}(Y;H_Y)$.
There exists $\gbigm_X\in\gbigc(X;H_X)$
such that $F^0_{\dagger}(\gbigm_X)=\gbigm_Y$.
Because
$(\lefttop{\tau}F)^0_{\dagger}
\bigl(
\lefttop{\tau}\gbigm_X
\bigr)
=\lefttop{\tau}\gbigm_Y$,
we obtain
$\lefttop{\tau}\gbigm_X
\in\gbigc(\lefttop{\tau}X;\lefttop{\tau}H_X)(\ast\tau)$,
and hence
$\gbigm_X\in\gbigc(X;H_X)$.
\hfill\qed

\subsubsection{Localizability and Beilinson functors}

Let $X$ be a complex manifold with a hypersurface $H$.
\begin{prop}
\label{prop;21.6.29.24}
 Let $\gbigm\in\gbigc_{\res}(X;H)$.
\begin{itemize}
 \item For any hypersurface $H^{(1)}$ of $X$,
       $\gbigm[\star H^{(1)}]$ $(\star=!,\ast)$ are objects
       of $\gbigc_{\res}(X;H)$.
 \item For any meromorphic function $g$ of $(X,H)$,
       $\Pi^{a,b}_{g\star}(\gbigm)$,
       $\Pi^{a,b}_{g,\ast!}(\gbigm)$,
       $\Xi^{(a)}_g(\gbigm)$,
       $\psi^{(a)}_g(\gbigm)$
       and $\phi^{(0)}_g(\gbigm)$
       are objects  of $\gbigc_{\res}(X;H)$.
\end{itemize}
\end{prop}
\pf
There exists $\gbigm_1\in\gbigc_{\res}(X)$
such that $\gbigm_1(\ast H)=\gbigm$.
By Theorem \ref{thm;21.3.12.12},
$\gbigm_1[\star H^{(1)}]$
is an object of $\gbigc_{\res}(X)$.
Because
$\gbigm[\star H^{(1)}]
=(\gbigm_1[\star H^{(1)}])(\ast H)$,
we obtain
$\gbigm[\star H^{(1)}]\in\gbigc_{\res}(X;H)$.

Let us study the second claim.
Let $g_0$ be the meromorphic function on
$(\lefttop{\tau}X,\lefttop{\tau}H)$
induced by $g$.
There exist natural isomorphisms
$\lefttop{\tau}
\bigl(
\II_{\gbigx,g}^{a,b}\otimes \gbigm
\bigr)
\simeq
\II_{\lefttop{\tau}\gbigx,g_0}^{a,b}
\otimes
\lefttop{\tau}\gbigm$
of $\gbigrtilde
 _{\lefttop{\tau}\gbigx(\ast \lefttop{\tau}\gbigh)}
 \bigl(\ast(\tau g_0)\bigr)$-modules.
We set $H^{(2)}=|(g)_0|$.
We have $\lefttop{\tau}H^{(2)}=|(g_0)_0|$.
Note that
$\Pi^{a,b}_{g,\star}(\gbigm)
=\bigl(
\II_{\gbigx,g}^{a,b}\otimes \gbigm
\bigr)[\star H^{(2)}]$
and 
$\Pi^{a,b}_{g_0,\star}\bigl(\lefttop{\tau}\gbigm\bigr)
=\bigl(
\II_{\lefttop{\tau}\gbigx,g_0}^{a,b}\otimes \lefttop{\tau}\gbigm
\bigr)[\star \lefttop{\tau}H^{(2)}]$.
Then, the claim follows from Lemma \ref{lem;21.3.13.50}.
\hfill\qed

\begin{cor}
In the situation of {\rm\S\ref{subsection;21.6.22.2}}
and {\rm\S\ref{subsection;21.6.22.20}},
for any $\gbigm\in\gbigc_{\res}(Y;H_Y)$,
$(\lefttop{T}f^{\star})^i(\gbigm)$ $(\star=!,\ast)$
are objects of
$\gbigc_{\res}(X;H_X)$.
\hfill\qed
\end{cor}

\subsection{Regularity of rescalable objects along $\tau$}

Let $\pi:\lefttop{\tau}X\lrarr X$ denote the projection.
We shall prove the following theorem in \S\ref{subsection;21.4.16.13}.

\begin{thm}
\label{thm;21.3.15.2}
 For any $\gbigm\in\gbigc_{\res}(X)$,
$\lefttop{\tau}\gbigm$
is regular along $\tau$,
i.e.,
 each $V_{a}\bigl(
 \lefttop{\tau}\gbigm
 \bigr)$
is coherent over
$\pi^{\ast}\gbigr_X$.
\end{thm}

\begin{cor}
Let $H$ be a hypersurface of $X$.
For any $\gbigm\in\gbigc_{\res}(X;H)$,
$\lefttop{\tau}\gbigm$
is regular along $\tau$,
i.e.,
each $V_{a}\bigl(
\lefttop{\tau}\gbigm
\bigr)$
is coherent over
$\pi^{\ast}\gbigr_{X(\ast H)}$.
\hfill\qed
\end{cor}

As a complement to Theorem \ref{thm;21.3.15.2},
we shall prove the following proposition
in \S\ref{subsection;21.4.16.12}.
\begin{prop}
\label{prop;21.4.16.11}  
Let $\gbigm\in\gbigc_{\res}(X)$.
We have
$\bigcap_{a\in\real}V_a\bigl(\lefttop{\tau}\gbigm\bigr)=0$. 
Each $V_a\bigl(\lefttop{\tau}\gbigm\bigr)$
is flat over $\nbigo_{\proj^1_{\lambda}\times\cnum_{\tau}}$. 
\end{prop}

\subsubsection{Preliminary}

We set $\mu=\lambda^{-1}$.
There exists a $V$-filtration
$\lefttop{\infty}V\gbigm_{|\gbigx^{\circ}}$ along $\mu$.
Note that each $\lefttop{\infty}V_a\gbigm_{|\gbigx^{\circ}}$
is coherent over $\nbigd_{\gbigx^{\circ}/\mu}$.
We have the morphism
$\varphi_1^{\circ}:
\varphi_0^{-1}\bigl(
\lefttop{\tau}\gbigx^{\circ}
\bigr)
\lrarr
\gbigx^{\circ}$ induced by $\varphi_1$,
and the isomorphism
$\varphi_0^{\circ}:
\varphi_0^{-1}\bigl(
\lefttop{\tau}\gbigx^{\circ}
\bigr)
\simeq
\lefttop{\tau}\gbigx^{\circ}$.
We obtain
\[
V_a\Bigl(
 \lefttop{\tau}\gbigm_{|\lefttop{\tau}\gbigx^{\circ}}
\Bigr):=
\varphi_{0\ast}^{\circ}
 (\varphi_1^{\circ})^{\ast}
 \bigl(
 \lefttop{\infty}V_a\gbigm
 \bigr)(\ast\lefttop{\tau}\gbigx^{\infty})
 \subset
 \lefttop{\tau}\gbigm_{|\lefttop{\tau}\gbigx^{\circ}}.
\]
We can easily check that
$V_a\Bigl(
 \lefttop{\tau}\gbigm_{|\lefttop{\tau}\gbigx^{\circ}}
 \Bigr)$
are coherent over
$\pi^{\ast}(\gbigr_X)_{|\lefttop{\tau}\gbigx^{\circ}}$,
and that
$V_{\bullet}\Bigl(
\lefttop{\tau}\gbigm_{|\lefttop{\tau}\gbigx^{\circ}}
\Bigr)$
is a $V$-filtration of
$\lefttop{\tau}\gbigm_{|\lefttop{\tau}\gbigx^{\circ}}$
along $\tau=0$.
Hence, it is enough to check
the claim of Theorem \ref{thm;21.3.15.2}
for the restriction to $\lefttop{\tau}\nbigx$.

\subsubsection{Good-KMS case}

Let $X$ be a neighbourhood of $(0,\ldots,0)$ in $\cnum^n$.
Let $H=\bigcup_{i=1}^{\ell}\{z_i=0\}$.
Let $\nbigv$ be a good-KMS smooth $\nbigrtilde_{X(\ast H)}$-module
such that $\nbigv[\ast H]\in\nbigc_{\res}(X)$.
We set $\gbigv:=\Upsilon(\nbigv)$.
Let $I\sqcup J=\{1,\ldots,\ell\}$ be a decomposition.

\begin{prop}
\label{prop;21.3.15.10}
$\lefttop{\tau}\bigl(
\gbigv[\ast I!J]
\bigr)$ 
is regular along $\tau$.
\end{prop}
\pf
For any $e\in\seisuu_{>0}$,
let $\varphi_e:X^{(e)}\lrarr X$
be the morphism defined by
\[
 \varphi_e(\zeta_1,\ldots\zeta_n)
=(\zeta_1^{e},\ldots,\zeta_{\ell}^e,\zeta_{\ell+1},\ldots,\zeta_n).
\]
For an appropriate $e$,
$\nbigv^{(e)}:=\varphi_e^{\ast}(\nbigv)$ is
unramifiedly good-KMS.
We set $\gbigv^{(e)}:=\Upsilon(\nbigv^{(e)})$.
There exist natural monomorphisms
$\gbigv[\ast I!J]\lrarr
\varphi_{e\dagger}(\gbigv^{(e)}[\ast I!J])$,
which induce
\[
 \lefttop{\tau}\bigl(
 \gbigv[\ast I!J]\bigr)
 \lrarr
 \varphi_{e\dagger}
 \lefttop{\tau}\bigl(
 \gbigv^{(e)}[\ast I!J]
 \bigr).
\]
Therefore, it is enough to study the case
where $\nbigv$ is unramifiedly good.

Let $\nbigi$ be a good set of irregular values,
and suppose that $\nbigv$ is good-KMS over $\nbigi$.
We obtain $\nbigs(\nbigi)\subset\seisuu_{\leq 0}^{\ell}$
as in \S\ref{subsection;21.3.15.1}.
By setting $\Xtilde=\lefttop{\tau}X$,
we apply the construction in \S\ref{subsection;21.3.15.1}.
We obtain
$\Xtilde_{\nbigs(\nbigi)}$,
$\Htilde_{\nbigs(\nbigi)}$
and
$F_{\nbigs}$.
Recall that $F_{\nbigs(\nbigi)}^{\ast}(\vecnbigitilde)$ is
a good system of irregular values on
$(\Xtilde_{\nbigs(\nbigi)},\Htilde_{\nbigs(\nbigi)})$.

We set
$\gbigvtilde:=\lefttop{\tau}\gbigv$
on $(\gbigxtilde,\gbigxtilde_0\cup\gbightilde)$.
We obtain
$F^{\ast}(\gbigvtilde)$
on $(\gbigxtilde_{\nbigs(\nbigi)},\gbightilde_{\nbigs(\nbigi)})$.

\begin{lem}
$\nbigvtilde':=
 F^{\ast}_{\nbigs(\nbigi)}(\gbigvtilde)_{|\cnum_{\lambda}\times U_Q}$
is unramifiedly good-KMS over $F^{\ast}(\vecnbigitilde)$.
\end{lem}
\pf
There exists $\nbigm\in\nbigc(\lefttop{\tau}X)$
such that
$\nbigm(\ast \lefttop{\tau}H)(\ast\tau)
=\gbigvtilde_{|\lefttop{\tau}\nbigx}$.
There exists a filtration $\nbigw$ of $\nbigm$
such that $(\nbigm,\nbigw)$ underlies
an integrable mixed twistor $\nbigd$-module on
$\lefttop{\tau}X$.
It induces a filtration $\nbigw$
of $\nbigrtilde_{\lefttop{\tau}X(\ast\lefttop{\tau}H)}(\ast\tau)$-module
of $\gbigvtilde_{|\lefttop{\tau}\nbigx}$.
According to \cite[Theorem 11.1.2, Theorem 19.1.3]{Mochizuki-wild},
$F_{\nbigs(\nbigi)}^{\ast}
\bigl(
\Gr^{\nbigw}_{j}(\gbigvtilde_{|\lefttop{\tau}\nbigx})
\bigr)$
are unramifiedly good-KMS over
$F_{\nbigs(\nbigi)}^{\ast}(\vecnbigitilde)$.
If $C$ is a smooth curve in $\Xtilde_{\nbigs(\nbigi)}$
which intersects with the smooth part of
$\Htilde_{\nbigs(\nbigi)}$,
the restriction of
$F_{\nbigs(\nbigi)}^{\ast}(\gbigvtilde_{|\lefttop{\tau}\nbigx})$
to $C$ has an unramifiedly good KMS-structure
with which $\nbigw$ is compatible.
By \cite[Proposition 5.2.9]{Mochizuki-MTM},
we obtain that
$F_{\nbigs(\nbigi)}^{\ast}\bigl(
\gbigvtilde_{|\lefttop{\tau}\nbigx}
\bigr)$
is unramifiedly good-KMS over
$F_{\nbigs(\nbigi)}^{\ast}(\vecnbigitilde)$,
with which $\nbigw$ is compatible.
\hfill\qed

\vspace{.1in}
The rest of the argument is similar to
the proof of Proposition \ref{prop;21.1.23.11}.
We shall explain only an outline.

Let $\Htilde_{\nbigs(\nbigi)}=\bigcup_{j\in\Lambda(\nbigi)}\Htilde_{\nbigs(\nbigi),j}$
denote the irreducible decomposition.
For each $i=1,\ldots,\ell$,
we obtain
$\Lambda(\nbigi,i)\subset
\Lambda(\nbigi)$
determined by
$F_{\nbigs(\nbigi)}^{-1}(\Htilde_{i})
=\bigcup_{j\in\Lambda(\nbigi,i)}\Htilde_{\nbigs(\nbigi),j}$.
There exists $[\tau]\in\Lambda(\nbigi)$
such that
$\Htilde_{\nbigs(\nbigi),[\tau]}$
is the proper transform of $\Xtilde_0$.
Let $\veca(I,J)\in\real^{\Lambda(\nbigi)}$
be determined by
$\veca(I,J)_j=1$ $(j\in \{[\tau]\}\cup \bigcup_{i\in I}\Lambda(\nbigi,i))$
and
$\veca(I,J)_j=1-\epsilon$ $(j\in \bigcup_{i\in J}\Lambda(\nbigi,i))$
for a sufficiently small $\epsilon>0$.
Note that
$\nbigq_{\veca(I,J)}(\nbigvtilde')$
is independent of any sufficiently small $\epsilon>0$.

Let $V_{\Htilde_{\nbigs(\nbigi)}}\nbigrtilde_{\Xtilde_{\nbigs(\nbigi)}}
\subset
\nbigrtilde_{\Xtilde_{\nbigs(\nbigi)}}$
be the sheaf of subalgebras
generated by
$\lambda \Theta_{\Xtilde_{\nbigs(\nbigi)}}(\log \Htilde_{\nbigs(\nbigi)})$
and $\lambda^2\del_{\lambda}$.
We obtain the following
$V_{\Htilde_{\nbigs(\nbigi)}}\nbigrtilde_{\Xtilde_{\nbigs(\nbigi)}}$-submodule
\[
 \nbign(I,J):=
 V_{\Htilde_{\nbigs(\nbigi)}}\nbigrtilde_{\Xtilde_{\nbigs(\nbigi)}}\cdot
 \nbigq_{\veca(I,J)}(\nbigvtilde')
\subset\nbigvtilde'. 
\]
Let $V_{\Htilde\cup \Xtilde_0}\nbigrtilde_{\Xtilde}
\subset
\nbigrtilde_{\Xtilde}$
denote the sheaf of subalgebras generated by
$\lambda\Theta_{\Xtilde}(\log(\Htilde\cup\Xtilde_0))$
and $\lambda^2\del_{\lambda}$.
We obtain the $V_{\Htilde\cup\Xtilde_0}\nbigrtilde_{\Xtilde}$-submodule
\[
 F_{\nbigs(\nbigi)\ast}(\nbign(I,J))
 \subset
 \gbigvtilde_{|\nbigxtilde}.
\]
We obtain the following isomorphism
by using the argument in the proof of Lemma \ref{lem;21.3.9.10}:
\begin{equation}
 \nbigvtilde[\ast I!J]
 \simeq
 \nbigr_{\Xtilde}(\ast \tau)
 \otimes_{V_{\Htilde\cup\Xtilde_0}\nbigr_{\Xtilde}}
 F_{\nbigs(\nbigi)\ast}(\nbign(I,J)).
\end{equation}
Let $V\nbigrtilde_{\Xtilde}\subset \nbigrtilde_{\Xtilde}$
denote the sheaf of subalgebras generated by
$\lambda\Theta_{\Xtilde}(\log \tau)$
and $\lambda^2\del_{\lambda}$.
Let $\nbigl(I,J)\subset \nbigvtilde(\ast I!J)$
denote the $V\nbigrtilde_{\Xtilde}$-submodule
generated by the image of the natural morphism
$F_{\nbigs(\nbigi)\ast}\nbign(I,J)
\lrarr \nbigv(\ast I!J)$.
By using the argument in the proof of
Lemma \ref{lem;21.3.8.12}
and Lemma \ref{lem;21.3.8.11},
we obtain that $\nbigl(I,J)$
is coherent over $\pi^{\ast}\nbigr_X$.
Then, we obtain that
$V_a(\nbigv(\ast I!J))$
are coherent over $\pi^{\ast}\nbigr_{X}$
by using the argument in the proof of
Lemma \ref{lem;21.3.9.30}.
Thus, we obtain Proposition \ref{prop;21.3.15.10}.
\hfill\qed

\subsubsection{Proof of Theorem \ref{thm;21.3.15.2}}
\label{subsection;21.4.16.13}

We use the induction on the dimension of the support.
It is enough to prove the claim locally around any point of $X$.
Therefore, we may assume the existence of
$f$, $Z$, $\varphi$ and $\nbigv$
for $\gbigm_{|\nbigx}$
as in the proof of Proposition \ref{prop;21.2.18.1}.
Let $f_0$ be the holomorphic function on $\lefttop{\tau}X$
induced by $f$.
Then, $V_a(\lefttop{\tau}\gbigm)$ is described as the following complex:
\[
 V_a\bigl(
  \psi^{(1)}_{f_0}(\lefttop{\tau}\gbigm)
  \bigr)
  \lrarr
  V_a\bigl(\Xi^{(0)}_{f_0}(\lefttop{\tau}\gbigm)\bigr)
  \oplus
  V_a\bigl(\phi^{(0)}_{f_0}(\lefttop{\tau}\gbigm)\bigr)
  \lrarr
  V_a\bigl(
  \psi^{(0)}_{f_0}(\lefttop{\tau}\gbigm)
  \bigr).
\]
We set $f_Z:=\varphi^{\ast}(f)$.
Let $f_{Z,0}$ be the holomorphic function
on $\lefttop{\tau}Z$ induced by $f_Z$.
We set $\gbigv:=\Upsilon(\nbigv)$.
We have
$\varphi_{\dagger}\bigl(
 \Xi^{(0)}_{f_{Z,0}}(\lefttop{\tau}\gbigv)
 \bigr)
 \simeq
 \Xi^{(0)}_{f_0}(\lefttop{\tau}\gbigm)$.
Therefore, it is enough to prove the claim
for 
$\Xi^{(0)}_{f_{Z,0}}(\gbigv)$,
which follows from Proposition \ref{prop;21.3.15.10}.
\hfill\qed

\subsubsection{Proof of Proposition \ref{prop;21.4.16.11}}
\label{subsection;21.4.16.12}

We easily obtain
$\bigcap_{a\in\real}
V_a\bigl(\lefttop{\tau}\gbigm\bigr)=0$
on $\lefttop{\tau}\gbigx\setminus (\lefttop{\tau}\gbigx)^0$
from Proposition \ref{prop;21.4.16.2}.
We obtain 
$\bigcap_{a\in\real}
V_a\bigl(\lefttop{\tau}\gbigm\bigr)=0$
on $\lefttop{\tau}\gbigx$
from the strictness.

Let us study the flatness at
$(\lambda_0,0)\in\cnum_{\lambda}\times\cnum_{\tau}$.
Let $A$ denote the stalk of
$\nbigo_{\proj^1_{\lambda}\times\cnum_{\tau}}$
at $(\lambda_0,0)$.
Let $M$ denote the quotient of $A$
by the maximal ideal of $A$.
\begin{lem}
\label{lem;21.4.16.21}
 $\Tor^A_i\bigl(M,V_{a}\bigl(\lefttop{\tau}\gbigm\bigr)_{(\lambda_0,0)}
 \bigr)=0$
for $i\neq 0$.
\end{lem}
\pf
There exists the standard exact sequence
$0\lrarr A\stackrel{\kappa_1}{\lrarr} A^2
\stackrel{\kappa_2}{\lrarr} A\lrarr M\lrarr 0$,
where $\kappa_1(b)=((\lambda-\lambda_0)b,\tau b)$
and $\kappa_2(a_1,a_2)=\tau a_1-(\lambda-\lambda_0)a_2$.
Let us look at the complex
\[
 V_{a}\bigl(\lefttop{\tau}\gbigm\bigr)_{(\lambda_0,0)}
 \stackrel{\kappa_3}{\lrarr}
  V_{a}\bigl(\lefttop{\tau}\gbigm\bigr)^2_{(\lambda_0,0)}
  \stackrel{\kappa_4}{\lrarr}
 V_{a}\bigl(\lefttop{\tau}\gbigm\bigr)_{(\lambda_0,0)},
\]
where $\kappa_3$ and $\kappa_4$
are induced by $\kappa_1$ and $\kappa_2$,
respectively,
and the last 
$V_{a}\bigl(\lefttop{\tau}\gbigm\bigr)_{(\lambda_0,0)}$
sits at the degree $0$.
It is naturally quasi-isomorphic to
\[
  V_{a}\bigl(\lefttop{\tau}\gbigm\bigr)_{(\lambda_0,0)}
  \big/
  V_{a-1}\bigl(\lefttop{\tau}\gbigm\bigr)_{(\lambda_0,0)}
  \stackrel{\kappa_5}{\lrarr}
  V_{a}\bigl(\lefttop{\tau}\gbigm\bigr)_{(\lambda_0,0)}
  \big/
  V_{a-1}\bigl(\lefttop{\tau}\gbigm\bigr)_{(\lambda_0,0)},
\]
where $\kappa_5$ is induced by the multiplication of
$\lambda-\lambda_0$,
and 
the second
$V_{a}\bigl(\lefttop{\tau}\gbigm\bigr)_{(\lambda_0,0)}
 \big/
V_{a-1}\bigl(\lefttop{\tau}\gbigm\bigr)_{(\lambda_0,0)}$
sits at degree $0$.
Because $\Gr^V_b\bigl(\lefttop{\tau}\gbigm\bigr)$
$(a-1<b\leq a)$ are strict,
we obtain that $\kappa_5$ is a monomorphism.
Hence, we obtain the claim of the lemma.
\hfill\qed

\begin{lem}
\label{lem;21.4.16.20}
 For any $f\in A$,
the multiplication by $f$ on
$V_{a}\bigl(\lefttop{\tau}\gbigm\bigr)_{(\lambda_0,0)}$
is a monomorphism.
\end{lem}
\pf
There exists an expansion
$f=\sum_{j=0}^{\infty} f_j(\lambda)\tau^j$.
It is enough to consider the case where
$f_0(\lambda)$ is not constantly $0$.
Let $s\in
V_{a}\bigl(\lefttop{\tau}\gbigm\bigr)_{(\lambda_0,0)}$
be non-zero.
There exists $b\in\real$ such that
$s\in
V_{b}\bigl(\lefttop{\tau}\gbigm\bigr)_{(\lambda_0,0)}$
and that
the induced section
$[s]$ of
$\Gr^V_{b}\bigl(\lefttop{\tau}\gbigm\bigr)_{(\lambda_0,0)}$
is not $0$.
Because
$\Gr^V_{b}\bigl(\lefttop{\tau}\gbigm\bigr)$
is strict,
we obtain
$f_0[s]\neq 0$.
Hence, we obtain $fs\neq 0$.
\hfill\qed

\vspace{.1in}
Let $I$ be any ideal of $A$.
There exists a finite tuple $a_1,\ldots,a_m\in I$
which generates $I$.
Because $A$ is a unique factorization domain,
there exists a greatest common divisor $c$
of $a_1,\ldots,a_m$.
We set $I_0:=c\cdot A\supset I$.
There exists the exact sequence
$0\lrarr I_0/I\lrarr A/I\lrarr A/I_0\lrarr 0$.
By Lemma \ref{lem;21.4.16.20},
we have
$\Tor_i^A\Bigl(
A/I_0,
V_a\bigl(\lefttop{\tau}\gbigm\bigr)
\Bigr)=0$ for $i\neq 0$.
There exists a filtration $\nbigf$ on $I_0/I$ of finite length
such that $\Gr^{\nbigf}_j(I_0/I)\simeq M$.
Hence, we obtain
$\Tor_i^A\Bigl(I_0/I,
V_a\bigl(\lefttop{\tau}\gbigm\bigr)
\Bigr)=0$ for $i\neq 0$
by Lemma \ref{lem;21.4.16.21}.
Then, we obtain
$\Tor_i^A\Bigl(A/I,V_a\bigl(\lefttop{\tau}\gbigm\bigr)
\Bigr)=0$ $(i\neq 0)$.
It implies the flatness of 
$V_a\bigl(\lefttop{\tau}\gbigm\bigr)$ over $A$.
(See \cite{Matsumura}.)

We can apply a similar argument in the case
$(\lambda_0,\tau_0)\in\cnum_{\lambda}\times\cnum_{\tau}$
with $\tau\neq 0$,
and $(\infty,\tau_0)$.
Thus, we obtain Proposition \ref{prop;21.4.16.11}.
\hfill\qed

\subsubsection{Appendix: Strong regularity}

For a positive integer $m$,
let $\varphi_m:\cnum_{\tau_m}\times X\lrarr\lefttop{\tau}X$
be defined by
$\varphi_m(\lambda,\tau_m,x)=(\lambda,\tau_m^m,x)$.
For any open subset $U$ of $\lefttop{\tau}X$,
we set $U^{(m)}:=\varphi_m^{-1}(U)$.
Let $f$ be a coordinate function on $U^{(m)}$
defining $U^{(m)}\cap\{\tau_m=0\}$.
Let $\nbigr_{U^{(m)}/f}\subset\nbigr_{U^{(m)}}$
be the sheaf of subalgebras generated by
$\lambda p^{\ast}\Ker(df)$
over $\nbigo_{\cnum_{\lambda}\times U^{(m)}}$.
By using the arguments in Proposition \ref{prop;21.1.23.11},
we can prove the following proposition,
which amplifies Theorem \ref{thm;21.3.15.2}.

\begin{prop}
For any $\gbigm\in\gbigc_{\res}(X)$,
$\lefttop{\tau}\gbigm_{|\lefttop{\tau}\nbigx}$
is strongly regular along $\tau$ in the following sense.
\begin{itemize}
  \item For any $m$, $U$ and $f$ as above,
	$V_a(\varphi_m^{\ast}(\lefttop{\tau}\gbigm))$
	$(a\in\real)$
	are coherent over $\nbigr_{U^{(m)}/f}$.
	\hfill\qed
\end{itemize} 
\end{prop}

\subsection{Irregular Hodge filtrations associated with rescalable objects}

\subsubsection{Along $\lambda=\alpha\tau$ $(\alpha\neq 0)$} 

Let $\alpha\in\cnum^{\ast}$.
Let $\iota_{\lambda=\alpha\tau}:
\nbigx\lrarr \lefttop{\tau}\nbigx$
denote the morphism defined by
$\iota_{\alpha\lambda=\tau}(\lambda,x)
=(\alpha\lambda,\lambda,x)
\in\cnum_{\lambda}\times\cnum_{\tau}\times X
=\lefttop{\tau}\nbigx$.
Let $\gbigm\in\gbigc_{\res}(X)$.
By Theorem \ref{thm;21.3.15.2},
we obtain the $V\nbigrtilde_{\lefttop{\tau}X}$-module
$V_a(\lefttop{\tau}\gbigm)$
which is $\cnum^{\ast}$-equivariant
and coherent over $\pi^{\ast}\nbigr_{X}$.
By Lemma \ref{lem;21.3.15.20},
we obtain the $\cnum^{\ast}$-homogeneous coherent $\nbigr_{X}$-module
$\iota_{\lambda=\alpha\tau}^{\ast}V_a\bigl(
  \lefttop{\tau}\gbigm
  \bigr)$.
By Proposition \ref{prop;21.4.16.11},
we obtain
 $L\iota_{\lambda=\alpha\tau}^{\ast}V_a(\lefttop{\tau}\gbigm)
 =\iota_{\lambda=\alpha\tau}^{\ast}V_a(\lefttop{\tau}\gbigm)$.

\begin{lem}[\mbox{\cite[Lemma 2.21]{Sabbah-irregular-Hodge}}]
The multiplication by $\lambda-\alpha\tau$
on $\lefttop{\tau}\gbigm/V_a\lefttop{\tau}\gbigm$
is a monomorphism.
Hence,
 the natural morphism
 $\iota_{\lambda=\alpha\tau}^{\ast}V_a\bigl(
  \lefttop{\tau}\gbigm
      \bigr)
      \lrarr
       \iota_{\lambda=\alpha\tau}^{\ast}\bigl(
  \lefttop{\tau}\gbigm
      \bigr)$
is a monomorphism.
In particular, $\iota_{\lambda=\alpha\tau}^{\ast}V_a\bigl(
\lefttop{\tau}\gbigm
\bigr)$ is strict.
\end{lem}
\pf
The multiplication by $\lambda-\alpha\tau$
on $\Gr^V_b\bigl(\lefttop{\tau}\gbigm\bigr)$
is equal to the multiplication of $\lambda$,
which is a monomorphism.
It implies the claim of the lemma.
\hfill\qed

\vspace{.1in}

Let $\gbigm^{\alpha}$ be as in Lemma \ref{lem;22.7.29.1}.
We may naturally regard
$\gbigm^{\alpha}$
as a holonomic $\nbigd_{X}$-module.

\begin{prop}
\label{prop;21.6.29.20}
For $a\in\real$,
there uniquely exists
a coherent filtration
$F^{\irr}_{a+\bullet}$
of $\gbigm^{\alpha}$
such that
\[
 \Rtilde_{F^{\irr}_{a+\bullet}}(\gbigm^{\alpha})
=\iota_{\lambda=\alpha\tau}^{\ast}
V_a\bigl(\lefttop{\tau}\Upsilon(\gbigm)\bigr).
\]
They induce 
an $\real$-indexed coherent filtration $F^{\irr}_{\bullet}$
of $\gbigm^{\alpha}$.
 (See {\rm\S\ref{subsection;22.7.28.1}}
 for $\Rtilde_{F^{\irr}_{a+\bullet}}$.
See {\rm\S\ref{subsection;21.4.14.50}}
for the notion of $\real$-indexed good filtration.)
 \end{prop}
\pf
As in Lemma \ref{lem;22.7.29.1},
$\iota_{\lambda=\alpha\tau}^{\ast}(\lefttop{\tau}\gbigm)$
is isomorphic to
$p_{X}^{\ast}(\gbigm^{\alpha})$ on $\nbigx$
in a way compatible with the $\cnum^{\ast}$-action.
It naturally extends to
an isomorphism of $\gbigr_X$-modules
which is $\cnum^{\ast}$-equivariant.
Then, we have only to apply
Proposition \ref{prop;21.4.17.1}
and Proposition \ref{prop;21.4.14.51} below
to the coherent $\nbigr_X$-submodule
$\iota_{\lambda=\alpha\tau}^{\ast}V_a\bigl(
\lefttop{\tau}\gbigm
\bigr)
\subset
 \iota_{\lambda=\alpha\tau}^{\ast}
 \bigl(\lefttop{\tau}\gbigm\bigr)$. 
\hfill\qed

\vspace{.1in}

Let $\iota_{(\lambda,\tau)=(0,0)}:X\simeq (0,0)\times X
\lrarr \cnum_{\lambda}\times\cnum_{\tau}\times X
=\lefttop{\tau}\nbigx$
denote the inclusion.
Note that
\[
 \iota_{(\lambda,\tau)=(0,0)}^{\ast}
 V_a\bigl(\lefttop{\tau}\gbigm\bigr)
 \simeq
 \bigoplus_{j}
 \Gr_j^{F^{\irr}_{a+\bullet}}
 \gbigm^{\alpha}.
\]
On
$\Gr_j^{F^{\irr}_{a+\bullet}}\gbigm^{\alpha}$,
the natural $\cnum^{\ast}$-action
is given as the multiplication by
$c^j$ for $c\in\cnum^{\ast}$.
We also have the $\cnum^{\ast}$-equivariant filtrations
$F^{\irr}_{a+\bullet}$
on
$\iota_{\lambda=\alpha\tau}^{\ast}
V_a\bigl(\lefttop{\tau}\gbigm\bigr)$
such that
\[
 \Gr^{F^{\irr}_{a+\bullet}}
\iota_{\lambda=\alpha\tau}^{\ast}
V_a\bigl(\lefttop{\tau}\gbigm\bigr)
\simeq
 p_X^{\ast}\Bigl(
 \iota_{(\lambda,\tau)=(0,0)}^{\ast}
 V_a\bigl(\lefttop{\tau}\gbigm\bigr)
 \Bigr).
\]

\subsubsection{Along $\lambda=0$}

Let $\iota_{\lambda=0}:
\nbigx\lrarr \lefttop{\tau}\nbigx$
denote the morphism
defined by
$\iota_{\lambda=0}(\lambda,x)
=(0,\lambda,x)\in\cnum_{\lambda}\times\cnum_{\tau}\times X
=\lefttop{\tau}\nbigx$,
i.e.,
it is induced by the natural identification
$\nbigx\simeq \lefttop{\tau}X$.
We can identify 
$\iota_{\lambda=0}^{\ast}
 \pi^{\ast}(\nbigr_X)$
with 
$\Sym^{\bullet}\bigl(
\lambda p_{X}^{\ast}\Theta_X
\bigr)$
in the $\cnum^{\ast}$-equivariant way.
Hence,
$\iota_{\lambda=0}^{\ast}
V_a\bigl(
\lefttop{\tau}\gbigm
\bigr)$
is naturally a $\cnum^{\ast}$-equivariant
coherent $\Sym^{\bullet}\bigl(
\lambda p_X^{\ast}\Theta_X\bigr)$-module
on $\nbigx$.

Let $\gbigm^0$ be as in Lemma \ref{lem;22.7.29.1}.
It is naturally a coherent
$\Sym^{\bullet}\Theta_X$-module.
We obtain the following proposition from
Proposition \ref{prop;21.4.17.1}
and Proposition \ref{prop;21.4.14.51} below.
\begin{prop}
For each $a\in\real$,
there exists a coherent filtration
$F^{\irr}_{a+\bullet}$
of $\gbigm^0$
such that
\[
 \Rtilde_{F^{\irr}_{a+\bullet}}(\gbigm^0)
 \simeq
 \iota_{\lambda=0}^{\ast}
 V_a\bigl(
 \lefttop{\tau}\gbigm
 \bigr).
\]
They induce a good $\real$-coherent filtration
on $\gbigm^0$.
(See {\rm\S\ref{subsection;22.7.28.1}}
 and {\rm\ref{subsection;21.4.14.50}}
for coherent filtrations and $\Rtilde_{F^{\irr}}$
in this context.)
\hfill\qed
\end{prop}

There exists a natural isomorphism
\[
 \iota_{(\lambda,\tau)=(0,0)}^{\ast}
 V_a\bigl(\lefttop{\tau}\gbigm\bigr)
 \simeq
 \bigoplus_{j}
 \Gr_j^{F^{\irr}_{a+\bullet}}
 \gbigm^{0}.
\]

We also have the $\cnum^{\ast}$-equivariant filtrations
$F^{\irr}_{a+\bullet}$
on
$\iota_{\lambda=0}^{\ast}
V_a\bigl(\lefttop{\tau}\gbigm\bigr)$
such that
\[
 \Gr^{F^{\irr}_{a+\bullet}}
\iota_{\lambda=0}^{\ast}
V_a\bigl(\lefttop{\tau}\gbigm\bigr)
\simeq
 p_X^{\ast}\bigl(
 \iota_{(\lambda,\tau)=(0,0)}^{\ast}
 V_a\bigl(\lefttop{\tau}\gbigm\bigr)
 \bigr).
\]

\subsubsection{Along $\tau=0$}

Let $\iota_{\tau=0}:\nbigx=\lefttop{\tau}\nbigx_0
\lrarr\lefttop{\tau}\nbigx$
denote the inclusion.
We obtain the following $\nbigr_X$-module:
\[
 \iota_{\tau=0}^{\ast}V_a\bigl(
 \lefttop{\tau}\gbigm
 \bigr)
 =\iota_{\tau=0}^{-1}\Bigl(
V_a\bigl(
 \lefttop{\tau}\gbigm
 \bigr)
 \big/
V_{a-1}\bigl(
 \lefttop{\tau}\gbigm
 \bigr)
 \Bigr).
\]
For $a-1<b\leq a$, let
$V_b\Bigl(
  \iota_{\tau=0}^{\ast}V_a\bigl(
 \lefttop{\tau}\gbigm
 \bigr)\Bigr)$
 denote the image of
$\iota_{\tau=0}^{\ast}V_b\bigl(
 \lefttop{\tau}\gbigm
 \bigr)
 \lrarr
\iota_{\tau=0}^{\ast}V_a\bigl(
 \lefttop{\tau}\gbigm
 \bigr)$.
Then, there exists the natural isomorphism
\[
 \Gr^V\iota_{\tau=0}^{\ast}V_a\bigl(
 \lefttop{\tau}\gbigm
 \bigr)
 =\bigoplus_{a-1<b\leq a}
 \iota_{\tau=0}^{-1}\Bigl(
 \Gr^V_b\bigl(\lefttop{\tau}\gbigm\bigr)
 \Bigr).
\]

Let $\iota_{\infty}:\gbigx^{\infty}=\{\infty\}\times X\lrarr \gbigx$
denote the inclusion.
Let $\lefttop{\infty}V_{\bullet}\gbigm$
denote the $V$-filtration of
the $\nbigd_{\gbigx}$-module $\gbigm(\ast\gbigx^0)$
along $\lambda^{-1}$.
We may naturally regard 
$\iota_{\infty}^{\ast}\bigl(
 \lefttop{\infty}V_a\gbigm
 \bigr)$
 as a $\nbigd_X$-module.
For any $\lambda_0\in\cnum$,
let $\iota_{(\lambda,\tau)=(\lambda_0,0)}:
(\lambda_0,0)\times X\lrarr \lefttop{\tau}\nbigx$
denote the inclusion.
If $\lambda_0\neq 0$,
we may naturally regard
$\iota_{(\lambda,\tau)=(\lambda_0,0)}^{\ast}
V_a\bigl(\lefttop{\tau}\gbigm\bigr)$
as a $\nbigd_X$-module,
and there exists the following natural isomorphism of
$\nbigd_X$-modules:
\[
 \iota_{\infty}^{\ast}\bigl(
 \lefttop{\infty}V_a\gbigm
 \bigr)
 \simeq
 \iota_{(\lambda,\tau)=(\lambda_0,0)}^{\ast}
V_a\bigl(\lefttop{\tau}\gbigm\bigr).
\]
We obtain the following natural isomorphism of
$p_X^{\ast}\nbigd_X$-modules
which is $\cnum^{\ast}$-equivariant:
\begin{equation}
p_X^{\ast}\Bigl(
\iota_{\infty}^{\ast}\bigl(
 \lefttop{\infty}V_a\gbigm
 \bigr)
 \Bigr)(\ast\lambda)
 \simeq
 p_X^{\ast}\Bigl(
 \iota_{(\lambda,\tau)=(\lambda_0,0)}^{\ast}
V_a\bigl(\lefttop{\tau}\gbigm\bigr)
 \Bigr)(\ast\lambda)
 \simeq
 \iota_{\tau=0}^{\ast}
  \Bigl(
  V_a\bigl(\lefttop{\tau}\gbigm\bigr)
\Bigr)
(\ast\lambda).
\end{equation}

We obtain the following proposition from
Proposition \ref{prop;21.4.17.1}
and Proposition \ref{prop;21.4.14.51} below.
\begin{prop}
For each $a\in\real$,
there exist coherent filtrations $F^{\irr}_{a+\bullet}$
on $\iota_{\infty}^{\ast}V_a\gbigm$
and
$\iota_{(\lambda,\tau)=(\lambda_0,0)}^{\ast}
V_a\bigl(\lefttop{\tau}\gbigm\bigr)$ $(\lambda_0\neq 0)$
such that
\[
 \Rtilde_{F_{a+\bullet}^{\irr}}
 \bigl(
 \iota_{\infty}^{\ast}(\lefttop{\infty}V_a\gbigm)
 \bigr)
 \simeq
 \Rtilde_{F_{a+\bullet}^{\irr}}\bigl(
 \iota_{(\lambda,\tau)=(\lambda_0,0)}^{\ast}
V_a\bigl(\lefttop{\tau}\gbigm\bigr)
 \bigr)
 \simeq
  \iota_{\tau=0}^{\ast}
  \Bigl(
 V_a\bigl(\lefttop{\tau}\gbigm\bigr)
 \Bigr).
\]
\hfill\qed
\end{prop}

There exists a natural isomorphism
\[
 \iota_{(\lambda,\tau)=(0,0)}^{\ast}
 V_a\bigl(\lefttop{\tau}\gbigm\bigr)
 \simeq
 \bigoplus_{j}
 \Gr_j^{F_{a+\bullet}}
 \iota_{\infty}^{\ast}
 (\lefttop{\infty}V_a\gbigm)
 \simeq
  \bigoplus_{j}
 \Gr_j^{F_{a+\bullet}}
\iota_{(\lambda,\tau)=(\lambda_0,0)}^{\ast}
V_a\bigl(\lefttop{\tau}\gbigm\bigr).
\]
We also have the $\cnum^{\ast}$-equivariant filtration
$F^{\irr}_{a+\bullet}$
on
$\iota_{\tau=0}^{\ast}
V_a\bigl(\lefttop{\tau}\gbigm\bigr)$
such that
\[
 \Gr^{F^{\irr}_{a+\bullet}}
\iota_{\tau=0}^{\ast}
V_a\bigl(\lefttop{\tau}\gbigm\bigr)
\simeq
 p_X^{\ast}\Bigl(
 \iota_{(\lambda,\tau)=(0,0)}^{\ast}
 V_a\bigl(\lefttop{\tau}\gbigm\bigr)
 \Bigr).
\]

\subsubsection{A coherent filtration of
$V_a\lefttop{\tau}\gbigm$}

\begin{prop}
\label{prop;21.6.22.30}
 There uniquely exists a $\cnum^{\ast}$-invariant coherent filtration
$F_{\bullet}$ of $V_a\bigl(\lefttop{\tau}\gbigm\bigr)$
satisfying the following conditions.
\begin{itemize}
 \item $\iota^{\ast}_{\lambda=\alpha\tau}F_{\bullet}$
       is equal to
       $F^{\irr}_{a+\bullet}$
       on $\iota^{\ast}_{\lambda=\alpha\tau}
       V_a\bigl(\lefttop{\tau}\gbigm\bigr)$.
       Similar claims hold for
       $\iota^{\ast}_{\lambda=0}F_{\bullet}$
       and
       $\iota^{\ast}_{\tau=0}F_{\bullet}$.
 \item $\Gr^F\bigl(
       V_a\bigl(\lefttop{\tau}\gbigm\bigr)
       \bigr)
       =q^{\ast}
       \Bigl(
       \iota_{(\lambda,\tau)=(0,0)}^{\ast}
       V_a\bigl(\lefttop{\tau}\gbigm\bigr)
       \Bigr)$,
       where $q:\lefttop{\tau}\nbigx\lrarr X$
       denotes the projection.
\end{itemize}
\end{prop}
\pf
Let
$\varphi_0:\widetilde{\lefttop{\tau}\nbigx}\lrarr\lefttop{\tau}\nbigx$
be induced by the blow up
$\widetilde{\cnum_{\lambda}\times\cnum_{\tau}}\lrarr
\cnum_{\lambda}\times\cnum_{\tau}$
at $(\lambda,\tau)=(0,0)$.
We set $E:=\varphi_0^{-1}\bigl((0,0)\times X\bigr)$.
We obtain
\[
\varphi_0^{\ast}
 V_a\bigl(\lefttop{\tau}\gbigm\bigr)
\subset
\varphi_0^{\ast}
V_a\bigl(\lefttop{\tau}\gbigm\bigr)(\ast E).
\]

We may naturally identify
$\widetilde{\cnum_{\lambda}\times\cnum_{\tau}}$
as the total space of the line bundle $\nbigo_{\proj^1}(-1)$
on $\proj^1=\varphi_0^{-1}(0,0)$.
It induces
$\kappa:\widetilde{\cnum_{\lambda}\times\cnum_{\tau}}\lrarr
\varphi_0^{-1}(0,0)$.
Let 
$\kappa_X:\widetilde{\lefttop{\tau}\nbigx}\lrarr E$
denote the induced map.
On each open subset $U\subsetneq \varphi_0^{-1}(0,0)$,
there exists
a trivialization
\begin{equation}
\label{eq;21.4.17.2}
 \kappa^{-1}(U)\simeq
 \cnum\times U, 
\end{equation}
 which is $\cnum^{\ast}$-equivariant.
We obtain
$\kappa_X^{-1}(U\times X)
\simeq
 \cnum\times U\times X$.
Then, there exits a coherent filtration $F$
on
$\varphi_0^{\ast}\bigl(
V\bigl(\lefttop{\tau}\gbigm\bigr)
\bigr)_{|\{1\}\times U\times X}$
such that
the analytification of the Rees module
is $\cnum^{\ast}$-equivariantly isomorphic to
the restriction of
$\varphi_0^{\ast}\bigl(
V\bigl(\lefttop{\tau}\gbigm\bigr)
\bigr)$.
We obtain the $\cnum^{\ast}$-equivariant filtration
$F\Bigl(
\varphi_0^{\ast}\bigl(
V\bigl(\lefttop{\tau}\gbigm\bigr)
\bigr)_{|\kappa_X^{-1}(U\times X)}\Bigr)$
which is independent of
the choice of
a trivialization (\ref{eq;21.4.17.2}).
Hence, we obtain
the $\cnum^{\ast}$-equivariant filtration $F$ on 
$\varphi_0^{\ast}\bigl(
V\bigl(\lefttop{\tau}\gbigm\bigr)
\bigr)$.

Let $\qtilde:\widetilde{\lefttop{\tau}\nbigx}\lrarr X$
denote the projection.
By the construction,
we have
\begin{equation}
\label{eq;21.4.17.3}
 \Gr^F_j
 \varphi_0^{\ast}\bigl(
V\bigl(\lefttop{\tau}\gbigm\bigr)
\bigr)
=\qtilde^{\ast}
\Bigl(
 \iota_{(\lambda,\tau)=(0,0)}^{\ast}
 V\bigl(\lefttop{\tau}\gbigm\bigr)
\Bigr)_j.
\end{equation}
Then, it is easy to see that
\[
 R^k\varphi_{0\ast}
 F_j\Bigl(
 \varphi_0^{\ast}\bigl(
V\bigl(\lefttop{\tau}\gbigm\bigr)
\bigr)
 \Bigr)
=0\quad(k\neq 0).
\]
Thus, we obtain the filtration
\[
 F_j
 V\bigl(\lefttop{\tau}\gbigm\bigr)
 :=
 \varphi_{0\ast}
 F_j\Bigl(
 \varphi_0^{\ast}\bigl(
V\bigl(\lefttop{\tau}\gbigm\bigr)
\bigr)
 \Bigr).
\]
By the construction, the second condition is satisfied.
By the constructions of the filtrations,
we obtain
\begin{equation}
\label{eq;21.4.17.4}
 \iota_{\lambda=\alpha\tau}^{\ast}F_j
 \subset
 F_{a+j}^{\irr}
 \iota_{\lambda=\alpha\tau}^{\ast}
 V_{a}\lefttop{\tau}\gbigm.
\end{equation}
By using (\ref{eq;21.4.17.3}),
we obtain that
the equality holds in (\ref{eq;21.4.17.4}).
\hfill\qed

\subsubsection{Characteristic varieties}

Let $\gbigm\in\gbigc_{\res}(X)$.

\begin{cor}
$\Ch(\gbigm^{\alpha})$
is equal to the support of
$\nbigo_{T^{\ast}X}$-module
associated to 
$\iota_{(\lambda,\tau)=(0,0)}^{\ast}
V_a\bigl(\lefttop{\tau}\gbigm\bigr)$.
In particular,
$\Ch(\gbigm^{\alpha_1})=\Ch(\gbigm^{\alpha_2})$
for any $\alpha_1,\alpha_2\in\cnum$. 
\hfill\qed
\end{cor}

\begin{cor}
$\Ch(\gbigm_{|\nbigx})$
is equal to 
$\cnum\times\Ch(\Xi_{\DR}(\gbigm))$
in $\cnum\times T^{\ast}X$.
\end{cor}
\pf
The restriction of $F$ to $\cnum_{\lambda}\times\{1\}\times X$
induces a coherent filtration $F$ of
$\gbigm_{|\nbigx}$
such that
$\Gr^F(\gbigm_{|\nbigx})$ is isomorphic to
$p_{X}^{\ast}\bigl(
\iota_{(\lambda,\tau)=(0,0)}^{\ast}
V_{a}\bigl(\lefttop{\tau}\Upsilon(\nbigm)\bigr)
\bigr)$.
Then, the claim of the corollary follows.
\hfill\qed

\vspace{.1in}

We set $\mu=\lambda^{-1}$.
Let $\gbigm\in\gbigc_{\res}(X)$.
Let $\gbigr_{0,X}\subset \gbigr_X$
be the sheaf of subalgebras determined by
$\gbigr_{0,X|\nbigx}=\nbigr_X$
and
$\gbigr_{0,X|\gbigx^{\circ}}
=\nbigd_{\gbigx^{\circ}/\mu}$.
Note that $\gbigr_{0,X}$ is equipped with a natural filtration
$F(\gbigr_{0,X})$ by the orders of differential operators
such that
$\Gr^F(\gbigr_{0,X})\simeq
\Sym^{\bullet}\bigl(p_{1,X}^{\ast}(\Theta_X)\otimes\nbigo_{\proj^1}(-1)
\bigr)$.
For a coherent $\gbigr_{0,X}$-module $\gbign$,
we can naturally define the notion of good filtration,
from which we obtain
the characteristic variety
$\Ch(\gbign)\subset
p_{1,X}^{\ast}(T^{\ast}X)
\otimes\nbigo_{\proj^1}(1)$.
Because it is homogeneous,
we obtain
$\Ch(\gbign)\subset
p_{1,X}^{\ast}(T^{\ast}X)$.

\begin{cor}
$\Ch\bigl(\lefttop{\infty}V_a(\gbigm)\bigr)$
 is equal to
 $\proj^1\times \Ch(\Xi_{\DR}(\gbigm))$
in $\proj^1\times T^{\ast}X$. 
\hfill\qed
\end{cor}

We also have the following corollaries.

\begin{cor}
\[
 \Ch\bigl(
 \iota_{(\lambda,\tau)=(\lambda_0,0)}^{\ast}
 V_a(\lefttop{\tau}\gbigm)
 \bigr)
=\bigcup_{a-1<b\leq a}
 \Ch\bigl(
 \Gr^V_{b}(\lefttop{\tau}\gbigm)^{\alpha}
 \bigr)
=\Ch(\gbigm^{\alpha}).
\]
\hfill\qed
\end{cor}

\begin{cor}
$\Ch(V_a(\lefttop{\tau}\gbigm))$
is equal to
$\cnum_{\lambda}\times\cnum_{\tau}\times\Ch(\Xi_{\DR}\gbigm)$.
\hfill\qed
\end{cor}

\subsubsection{Meromorphic case}

Let $H$ be a hypersurface of $X$.
Let $\gbigm\in\gbigc_{\res}(X;H)$.
There exists $\gbigm_1\in\gbigc_{\res}(X)$
such that $\gbigm=\gbigm_1(\ast H)$.
For any $\alpha\in\cnum$,
we obtain the induced filtration
$F^{\irr}_{\bullet}(\gbigm^{\alpha})
=F^{\irr}_{\bullet}(\gbigm_1^{\alpha})(\ast H)$.
Similarly, we obtain the induced filtrations
$F^{\irr}_{\bullet}
\bigl(\iota_{\infty}^{\ast}(\lefttop{\infty}V_a\gbigm)\bigr)
=F^{\irr}_{\bullet}
\bigl(\iota_{\infty}^{\ast}(\lefttop{\infty}V_a\gbigm_1)\bigr)(\ast H)$.
We also obtain
$F_{\bullet}(V_a(\lefttop{\tau}\gbigm))
=F_{\bullet}(V_a(\lefttop{\tau}\gbigm_1))(\ast\lefttop{\tau}H)$
as the localization of the filtration in
Proposition \ref{prop;21.6.22.30}
which has a similar property.
The filtrations are independent of the choice of $\gbigm_1$.

\subsubsection{Some strictness properties}

We recall some strictness properties of the irregular Hodge filtrations
in \cite{Sabbah-irregular-Hodge} for the convenience.
Let $g:\gbigm_1\lrarr\gbigm_2$ be a morphism in
$\gbigc_{\res}(X;H)$.
Let $g^{\lambda}:\gbigmlambda_1\lrarr\gbigmlambda_2$,
$g_a^{\infty}:\iota_{\infty}^{\ast}(\lefttop{\infty}V_a(\gbigm_1))
\lrarr\iota_{\infty}^{\ast}(\lefttop{\infty}V_a(\gbigm_2))$
and 
$\lefttop{\tau}g:V_a(\lefttop{\tau}\gbigm_1)
 \lrarr V_a(\lefttop{\tau}\gbigm_2)$
denote the induced morphisms.
 
\begin{prop}[\cite{Sabbah-irregular-Hodge}]
We have
\[
 g^{\lambda}\bigl(F^{\irr}_{\bullet}\gbigm_1^{\lambda}\bigr)
       \subset
 F^{\irr}_{\bullet}\gbigm^{\lambda}_2,
 \quad
 g_a^{\infty}\bigl(
 F^{\irr}_{\bullet}\iota_{\infty}^{\ast}
 (\lefttop{\infty}V_a
  \gbigm_1^{\lambda})\bigr)
       \subset
 F^{\irr}_{\bullet}
 \iota_{\infty}^{\ast}(\lefttop{\infty}V_a\gbigm^{\lambda}_2),
 \quad
 \lefttop{\tau}g\bigl(F_{\bullet}
       V_a\lefttop{\tau}\gbigm_1\bigr)
       \subset
       F_{\bullet}
       V_a\lefttop{\tau}\gbigm_2.
\]
 If $\Cok(g)\in\gbigc_{\res}(X;H)$, 
 then $g^{\lambda}$ and $\lefttop{\tau}g$ are strict
 with respect to
 the irregular Hodge filtrations $F^{\irr}_{\bullet}$,
i.e., 
\[
 g^{\lambda}\bigl(
 F^{\irr}_{b}(\gbigmlambda_1)
 \bigr)
 =F^{\irr}_{b}(\gbigmlambda_2)
 \cap
 g^{\lambda}(\gbigmlambda_1),
 \quad
 g^{\infty}\bigl(
 F^{\irr}_{b}\iota_{\infty}^{\ast}(\lefttop{\infty}V_a(\gbigm_1))
 \bigr)
 =F^{\irr}_{b}\iota_{\infty}^{\ast}(\lefttop{\infty}V_a(\gbigm_2))
 \cap
 g^{\infty}\iota_{\infty}^{\ast}V_a(\gbigm_1),
\]
\[
 \lefttop{\tau}g\bigl(
 F_{n}(V_a\lefttop{\tau}\gbigmlambda_1)
 \bigr)
=F_{n}(V_a\lefttop{\tau}\gbigmlambda_2)
 \cap
 \lefttop{\tau}g(V_a\lefttop{\tau}\gbigmlambda_1)
\]
for any $b\in\real$.
\end{prop} 
\pf
The morphism $g$ induces the following
$\cnum^{\ast}$-equivariant morphism $\lefttop{\tau}g$
of $\pi^{\ast}\gbigr_{X(\ast H)}$-modules:
\begin{equation}
 \label{eq;21.3.15.30}
  V_a\lefttop{\tau}\gbigm_1
\lrarr
V_a\lefttop{\tau}\gbigm_2.
\end{equation}
Then, we obtain the first claim.
If $\Cok(g)\in\gbigc_{\res}(X;H)$,
then the cokernel of (\ref{eq;21.3.15.30})
is equal to
$V_a\lefttop{\tau}\Cok(g)$.
Then, we obtain the second claim.
\hfill\qed

\begin{prop}[\cite{Sabbah-irregular-Hodge}]
Let $f:X\lrarr Y$ be a projective morphism.
For any $\gbigm\in\gbigc_{\res}(X)$,
we have
 $f^j_{\dagger}(
 \iota_{\lambda=\alpha\tau}^{\ast}V_a(\lefttop{\tau}\gbigm))
 =\iota_{\lambda=\alpha\tau}^{\ast}
 V_a(\lefttop{\tau}(f^j_{\dagger}(\gbigm)))$.
In particular, we obtain the $E_1$-degeneration of
the spectral sequence associated with
the direct image of the filtered $\nbigd_X$-module
$(\gbigm^{\alpha},F^{\irr}_{a+\bullet})$ via $f$.
\hfill\qed
\end{prop}

Let $\gbigm\in\nbigc_{\res}(X)$.
Let $W$ be an increasing filtration of $\gbigm$
by $\nbigrtilde_X$-submodules indexed by $\seisuu$
such that the following holds.
\begin{itemize}
 \item $(\gbigm,W)_{|\nbigx}$ underlies
       an integrable mixed twistor $\nbigd$-module on $X$.
 \item There exists $\gbigm_0\in\gbigc(\lefttop{\tau}X)$
       with a filtration $W$ by
       $\gbigrtilde_{\lefttop{\tau}X}$-submodules
       such that
       $(\gbigm_0,W)_{|\lefttop{\tau}\nbigx}$ underlies
       an integrable mixed twistor $\nbigd$-module
       on $\lefttop{\tau}X$
       and that
       $(\gbigm_0,W)_{|\lefttop{\tau}\gbigx\setminus\lefttop{\tau}\gbigx_0}
       =\lefttop{\tau}(\gbigm,W)_{|\lefttop{\tau}\gbigx\setminus
        \lefttop{\tau}\gbigx_0}$.
\end{itemize}

\begin{prop}[\cite{Sabbah-irregular-Hodge}]
The irregular Hodge filtration
$F^{\irr}_{\bullet}$
on $W_j(\gbigm^{\alpha})/W_{k}(\gbigm^{\alpha})$
is equal to the filtration
induced by the irregular Hodge filtration 
$F^{\irr}_{\bullet}$ of $\gbigm^{\alpha}$. 
\end{prop}
\pf
We have
\[
V_{a}
\lefttop{\tau}\Upsilon\bigl(
\Gr^{\nbigw}_j\gbigm/\Gr^{\nbigw}_k\gbigm\bigr)
=
V_{a}
\lefttop{\tau}\Upsilon\bigl(
\Gr^{\nbigw}_j\gbigm\bigr)
\Big/
V_{a}
\lefttop{\tau}\Upsilon\bigl(
\Gr^{\nbigw}_k\gbigm\bigr).
\]
We obtain
\[
\iota_{\lambda=\alpha\tau}^{\ast}
V_{a}
\lefttop{\tau}\Upsilon\bigl(
\Gr^{\nbigw}_j\gbigm/\Gr^{\nbigw}_k\gbigm\bigr)
=
\iota_{\lambda=\alpha\tau}^{\ast}V_{a}
\lefttop{\tau}\Upsilon\bigl(
\Gr^{\nbigw}_j\gbigm\bigr)
\Big/
\iota_{\lambda=\alpha\tau}^{\ast}V_{a}
\lefttop{\tau}\Upsilon\bigl(
\Gr^{\nbigw}_k\gbigm\bigr).
\]
Then, the claim of the proposition follows.
\hfill\qed

\subsubsection{Strict specializability along
non-characteristic smooth hypersurfaces}

Suppose that $X$ is an open subset of $X_0\times\cnum_t$.
Let $\gbigm\in\gbigc_{\res}(X)$.
Suppose that $X_0\times\{0\}$
is non-characteristic to $\gbigm$.
\begin{prop}
For $\alpha\neq 0$,
the filtered $\nbigd$-module
$(\gbigm^{\alpha},F^{\irr}_{a+\bullet})$
is strictly specializable along $t$.
\end{prop}
\pf
Note that
$V_a\bigl(\lefttop{\tau}\Upsilon(\gbigm)\bigr)$
is coherent over
$\pi^{\ast}(\gbigr_{X/t})$.
Then, the claim follows.
\hfill\qed

\subsubsection{Appendix:
Coherent filtrations on $\nbigd$-modules
and Higgs sheaves}
\label{subsection;22.7.28.1}

We set $\nbiga_0:=\Sym^{\bullet}\Theta_X$.
It is equipped with a filtration
$F_j\nbiga_{0}=\bigoplus_{i\leq j} \Sym^i\Theta_X$ $(j\in\seisuu_{\geq 0})$.
We set $\nbiga_1:=\nbigd_X$.
Let $F_j\nbiga_1\subset\nbiga_1$ $(j\in\seisuu_{\geq 0})$
denote the subsheaves of differential operators
whose orders are less than $j$.
The sheaves of algebras
$\nbiga_i$ $(i=0,1)$ with $F$ are filtered rings,
i.e.,
(i) $\nbiga_i=\bigcup_nF_n(\nbiga_i)$,
(ii) $1\in F_0(\nbiga_i)$,
(iii) $F_j(\nbiga_i)\cdot F_k(\nbiga_i)\subset F_{j+k}(\nbiga_i)$,
(iv) $F_j(\nbiga_i)=0$ $(j<0)$.
(See \cite[A.1]{kashiwara_text}.)
We have the Rees algebras
$R_F(\nbiga_i)=\sum F_j(\nbiga_i)\lambda^j
\subset \nbiga[\lambda]$.
On $\gbigx$,
we set
$\gbiga_0:=
\Sym^{\bullet}(\lambda p_{1,X}^{\ast}\Theta_X)(\ast\gbigx^{\infty})$
and $\gbiga_1:=\gbigr_X$,
which are sheaves of algebras.
There exists a natural morphism
$p_{1,X}^{-1}R_F(\nbiga_i)\to \gbiga_i$
of sheaves of algebras.

Let $\mu:\cnum^{\ast}\times \gbigx\to\gbigx$
be the $\cnum^{\ast}$-action defined by
$\mu(a,(\lambda,x))=(a\lambda,x)$.
Let $p:\cnum^{\ast}\times\gbigx\to\gbigx$
denote the projection.
An $\nbigo_{\gbigx}(\ast\gbigx^{\infty})$-module $\gbigf$
is called $\cnum^{\ast}$-equivariant
if it is equipped with an isomorphism
$\rho:\mu^{\ast}\gbigf\simeq p^{\ast}\gbigf$
satisfying the cocycle condition.
If $\gbigf$ is $\cnum^{\ast}$-equivariant,
a subsheaf $\gbigf'$ is called $\cnum^{\ast}$-invariant
if $\mu^{\ast}(\gbigf')=p^{\ast}(\gbigf')$
under the isomorphism $\mu$.
It implies that $\gbigf'$ is $\cnum^{\ast}$-equivariant.

Let $M$ be a coherent $\nbiga_i$-module.
We set
$\gbigm:=p_{1,X}^{\ast}(M)\bigl(
\ast(\gbigx^0\cup\gbigx^{\infty})
\bigr)$.
It is naturally an $\gbiga_i(\ast\gbigx^0)$-module,
and $\cnum^{\ast}$-equivariant
as an $\nbigo_{\gbigx}(\ast\gbigx^{\infty})$-module.
Let $F$ be a coherent filtration of $M$.
(See \cite[Definition A.19]{kashiwara_text}
for the notion of coherent filtration.)
We obtain the Rees module
$R_F(M):=\bigoplus_{j\in\seisuu} F_j(M)\lambda^j$,
which is coherent over $R_F(\nbiga_i)$.
We set
\[
\gbigr_F(M)=\gbiga_i
\otimes_{p_{1,X}^{-1}R_F(\nbiga_i)}
R_F(M).
\]
It is naturally a coherent $\gbiga_i$-module,
and $\cnum^{\ast}$-equivariant
as an $\nbigo_{\gbigx}(\ast\gbigx^{\infty})$-module.
There exists a natural monomorphism
$\gbigr_F(M)\to\gbigm$
which induces
$\gbigr_F(M)(\ast\gbigx^0)\simeq\gbigm$.

\begin{rem}
The restriction $\gbigr_F(M)_{|\nbigx}$
is denoted by $\nbigrtilde_F(M)$.
\hfill\qed 
\end{rem}

\begin{lem}
Let $P$ be any point of $X$.
Then, $s\in M_P$ is contained in $F_jM_P$
if and only if $\lambda^jp_{1,X}^{\ast}(s)$
is contained in $\gbigr_F(M)_{(0,P)}$.
\end{lem}
\pf
The stalks of $\gbigr_F(M)$ and $\gbigm$
at $(0,P)\in\gbigx$
are
$\nbigo_{\gbigx,(0,P)}\otimes_{\nbigo_{X,P}[\lambda]}
R_F(M)_P$
and
$\nbigo_{\gbigx,(0,P)}\otimes_{\nbigo_{X,P}[\lambda]}
M_P[\lambda,\lambda^{-1}]$,
respectively.
Because the analytification is faithfully flat,
$\lambda^ks$ is $0$ in
$M_P[\lambda,\lambda^{-1}]/R_F(M)_P$
if and only if
$\lambda^kp_{1,X}^{\ast}(s)$ is $0$
in $\gbigm_{(0,P)}\big/\gbigr_F(M)_{(0,P)}$.
Then, the claim is clear.
\hfill\qed

\begin{prop}
\label{prop;21.4.17.1}
The above procedure induces an equivalence
between coherent filtrations $F$ on $M$
and coherent $\gbiga_i$-submodule
$\gbigf\subset\gbigm$
such that
(i) $\gbigf(\ast\gbigx^0)=\gbigm$,
(ii) $\gbigf$ is $\cnum^{\ast}$-invariant.
\end{prop}
\pf
Let $F$ be a coherent filtration of $M$.
For any $P\in X$ and $s$ of $M_P$,
we have
\begin{equation}
\label{eq;21.4.16.30}
 \min\bigl\{j\in\seisuu\,
 \big|\,
 \lambda^jp_{1,X}^{-1}(s)\in \gbigr_{F}(M)_{(0,P)}
 \bigr\}
 =\min\big\{
 s\in F_jM_P
 \}.
\end{equation}

Let $F$ and $F'$ be two coherent filtrations of $M$.
Suppose that
$\gbigr_{F'}(M)\subset\gbigr_{F}(M)$.
By (\ref{eq;21.4.16.30}),
we obtain $F'_j\subset F_j$ for any $j\in\seisuu$.
In particular,
if $\gbigr_{F'}(M)=\gbigr_{F}(M)$,
we obtain $F=F'$.

Let $\gbigf\subset\gbigm$
be a coherent $\gbigs_X$-submodule
such that
(i) $\gbigf(\ast\gbigx^0)=\gbigm$,
(ii) $\gbigf$ is $\cnum^{\ast}$-invariant.
Let $P\in X$
and $s\in M_P$.
There exists $j\in\seisuu$
such that
$\lambda^jp_{1,X}^{\ast}(s)\in\gbigf_{(0,P)}$.
\begin{lem}
 There exists the minimum of the set
$\bigl\{
 j\in\seisuu\,\big|\,
 \lambda^jp_{1,X}^{\ast}(s)
 \in\gbigf_{(0,P)}
 \bigr\}$
which we denote by 
$j_{\gbigf}(s)$.
\end{lem}
\pf
Suppose that $\lambda^jp_{1,X}^{\ast}(s)\in\gbigf_{(0,P)}$
for any $j\in\seisuu$.
There exists a neighbourhood $U$ of $P$ in $X$
on which $s$ is defined.
We set $\gbigu=\proj^1\times U\subset\gbigx$.
We set
$\gbign:=\sum_j \gbiga_i(\lambda^{j}s)\subset\gbigf_{|\gbigu}$.
By the Noetherian property of $\gbiga_i$,
we obtain that
$\gbign$ is a coherent $\gbiga_{i|\gbigu}$-module.
However,
for $N:=\nbiga_{i|U}\cdot s\subset M_{|U}$,
we obtain $\gbign=p_{1,U}^{\ast}(N)(\ast(\gbigu^0\cup\gbigu^{\infty}))$,
which is not coherent over $\gbiga_{i|\gbigu}$.
Thus, we obtain a contradiction.
\hfill\qed

\vspace{.1in}
For $P\in X$ and $j\in\seisuu$,
let $F_j(M_P)\subset M_P$ denote
the $\nbigo_{X,P}$-submodule of $s\in M_P$
such that $j_{\gbigf}(s)\leq j$.
They induce subsheaves $F_j(M)\subset M$.
We have
$F_j(\nbiga_i)\cdot F_k(M)\subset
F_{j+k}(M)$.
We obtain
the $\gbiga_i$-submodule
$\gbigr_F(M)\subset\gbigm$
associated with $F$.
By the construction,
we obtain
$\gbigr_F(M)\subset\gbigf$.

Let $(0,P)\in \gbigx^0$.
Let $U_P$ be a relatively compact neighbourhood of $P$ in $X$.
Let $U_0$ be a relatively compact neighbourhood of $0$ in $\cnum$.
Let $u$ be a section of $\gbigf/\lambda\gbigf$ on
$U_0\times U_P$.
We may assume that
there exist $s_1,\ldots,s_m$ of $M$ on $U_P$
and meromorphic functions $f_1,\ldots,f_m$ on $(U_0,0)\times U_P$
such that $\utilde=\sum_{i=1}^m f_ps_p$
is a section of $\gbigf$ on $U_0\times U_P$
which induces $u$.
There exist $k\in\seisuu$ such that $s_p\in F_k(M)$.
There exist the expansions
$f_p=\sum_j f_{p,j}\lambda^j$.
Note that
$\sum_{j\geq k+1}f_{p,j}\lambda^{j}s_p$
are sections of
$\lambda\gbigf$.
By replacing $f_p$ with $f_p-\sum_{j\geq k+1}f_{p,j}\lambda^j$,
we may assume that
$f_p=\sum_{j\leq k}f_{p,j}\lambda^j$.
Then, $\utilde$ is a section of
$\gbigf$ on $\proj^1\times U_P$.
Because $\gbigf$ is $\cnum^{\ast}$-invariant,
$a^{\ast}\utilde
=\sum_p\sum_j a^j\lambda^jf_{p,j}s_p$
are sections of $\gbigf$
for any $a\in\cnum^{\ast}$.
Hence, we obtain that
$\lambda^j \sum_pf_{p,j}s_p$
are sections of $\gbigm$,
which implies that
$\sum_pf_{p,j}s_p\in F_j(M)$.
Hence, $\utilde$ is a section of
$\gbigr_F(M)$.
We obtain that
the composition of
$\gbigr_F(M)\lrarr \gbigf\lrarr\gbigf/\lambda\gbigf$
is an epimorphism.
Then, we obtain
$\gbigr_F(M)\simeq \gbigf$.
We also obtain that
the Rees module of $(M,F)$ is coherent over
$R_F(\nbiga_i)$
from the $\gbiga_i$-coherence of $\gbigm$.
\hfill\qed

\subsubsection{Appendix: coherent filtrations indexed by $\real$}
\label{subsection;21.4.14.50}

We continue to use the notation in \S\ref{subsection;22.7.28.1}.
Let $M$ be a coherent $\nbiga_i$-module.
An $\real$-indexed coherent filtration $F$ of $M$
is an increasing tuple $F_a(M)$ $(a\in\real)$ of
coherent $\nbigo_X$-submodules
such that
the following holds.
\begin{itemize}
 \item For any $a\in\real$ and a compact subset $K$ of $X$,
       there exists $\epsilon>0$
       such that $F_a(M)_{|K}=F_{a+\epsilon}(M)_{|K}$.       
 \item $F_j\nbiga_i\cdot F_a(M)\subset F_{a+j}(M)$
       for any $j\in\seisuu_{\geq 0}$
       and $a\in\real$.
 \item For each $b\in\real$,
       $F_{b+\bullet}(M)=\{
       F_{b+n}(M)\,|\,n\in\seisuu\}$
       is a coherent filtration of $M$.
\end{itemize}
We obtain
an increasing sequence of $\gbiga_i$-submodules
$U_b(\gbigm)=\gbigr_{F_{b+\bullet}}(M)$ $(b\in\real)$
of $\gbigm$.
They satisfy the following conditions.
\begin{itemize}
 \item $U_b(\gbigm)(\ast\gbigx^0)=\gbigm$.
 \item $U_b(\gbigm)$
       are $\gbiga_i$-coherent
       and $\cnum^{\ast}$-equivariant.
 \item For each $b\in\real$
       and a compact subset $K$ of $X$,
       there exists $\epsilon>0$ such that
       $U_b(\gbigm)_{|K}
       =U_{b+\epsilon}(\gbigm)_{|K}$.
 \item $\lambda^n U_b(\gbigm)=U_{b-n}(\gbigm)$ for any $n\in\seisuu$.
\end{itemize}
We obtain the following proposition
from Proposition \ref{prop;21.4.17.1}.
\begin{prop}
\label{prop;21.4.14.51}
The above construction induces an equivalence between
$\real$-indexed coherent filtrations $F$ on $M$
and
$\real$-indexed increasing sequences $U_{\bullet}(\gbigm)$
of $\cnum^{\ast}$-equivariant
$\gbigr_X$-coherent submodules of 
$\gbigm$ satisfying the above conditions.
\hfill\qed
\end{prop}

\subsection{Duality}

Let $\gbigm\in\gbigc_{\res}(X)$.
We shall prove the following theorem
in \S\ref{subsection;21.4.14.60}.
\begin{thm}
\label{thm;21.3.25.20}
There exists a natural isomorphism
\begin{equation}
\label{eq;21.3.25.3}
 \lambda\cdot \lefttop{\tau}\bigl(
 \DD_X\gbigm\bigr)
\simeq
 \DD_{\lefttop{\tau}X(\ast \tau)}\bigl(
 \lefttop{\tau}\gbigm
 \bigr).
\end{equation}
In particular, $\DD_X\gbigm$ is an object of $\gbigc_{\res}(X)$.
\end{thm}

Before proving the theorem,
we state a consequence for the associated irregular Hodge filtrations
of the underlying $\nbigd$-module $\Xi_{\DR}(\gbigm)=\gbigm^1$.
According to Theorem \ref{thm;21.3.25.20}
and Proposition \ref{prop;21.3.30.3},
we have
\[
 V_{-1-a}\Bigl(
 \lefttop{\tau}\DD_X\gbigm
 \Bigr)
 \simeq
\nbigh^0\Bigl(
 \nrhom_{\pi^{\ast}\gbigr_X}
 \Bigl(
 V_{<a}\bigl(
 \lefttop{\tau}\gbigm\bigr),\,
 \pi^{\ast}\gbigr_X
 \otimes
 (\lambda^{d_X}\cdot\omega_{\gbigx}^{-1})
 \Bigr)[d_X]
\Bigr),
\]
and 
\[
 \nbigh^j\Bigl(
 \nrhom_{\pi^{\ast}\gbigr_X}
 \Bigl(
 V_{<a}\bigl(
 \lefttop{\tau}\gbigm\bigr),\,
 \pi^{\ast}\gbigr_X
 \otimes
 (\lambda^{d_X}\cdot\omega_{\gbigx}^{-1})
 \Bigr)[d_X]
 \Bigr)
 =0\quad (j\neq 0).
\]

\begin{cor}
For $\gbigm\in\gbigc_{\res}(X)$,
we obtain
\[
\iota_{\lambda=\tau}^{\ast}
 V_a\Bigl(
 \lefttop{\tau}\DD_X\gbigm
\Bigr)
 =\DD_X\Bigl(
 \iota_{\lambda=\tau}^{\ast}
 V_{<-1-a}\bigl(
 \lefttop{\tau}\gbigm
 \bigr)
 \Bigr).
\]
As a result,
we obtain the following duality
for the filtered $\nbigd_X$-modules:
\[
 R_{F^{\irr}_{a+\bullet}}
 \DDD_X\Xi_{\DR}(\gbigm)
 \simeq
 \DD_X\Bigl(
 R_{F^{\irr}_{-1-a-\epsilon+\bullet}}
 \Xi_{\DR}(\gbigm)
 \Bigr).
\]
Here, $\epsilon$ denotes any small positive number.
\hfill\qed
\end{cor}

\subsubsection{Proof of Theorem \ref{thm;21.3.25.20}}
\label{subsection;21.4.14.60}

Recall
$\gbigy=\lefttop{\tau}\gbigx_0
\cup
 \lefttop{\tau}\gbigx^{\infty}$
and
$\gbigytilde=\varphi_0^{-1}(\gbigy)$.

\begin{lem}
\label{lem;21.6.18.1}
For any
$\gbigrtilde_X$-module $\gbigm$
which is a coherent over $\gbigr_X$,
there exists a natural quasi-isomorphism of
$\varphi_1^{\ast}(\gbigr_X)
\langle \lambda^2\del_{\lambda},
  \deldel_{\tau}\rangle$-complexes:
\[
 \varphi_1^{\ast}\Bigl(
 \nrhom_{\gbigr_X}\bigl(\gbigm,
  \gbigr_X\otimes\omega^{-1}_{\gbigx}\bigr)
 \Bigr)(\ast\gbigytilde)
 \simeq
 \nrhom_{\varphi_1^{\ast}(\gbigr_X)(\ast\gbigytilde)}
 \Bigl(
 \varphi_1^{\ast}(\gbigm)(\ast\gbigytilde),
 \varphi_1^{\ast}\bigl(
 \gbigr_X\otimes\omega^{-1}_{\gbigx}
 \bigr)(\ast\gbigytilde)
 \Bigr).
\] 
\end{lem}
\pf
We set
\[
 \nbiga_0=
(\gbigr_X\otimes_{\nbigo_{\proj^1}}\gbigr_X)
 \langle\lambda^2\del_{\lambda}\rangle,
 \quad
 \nbiga_1=
 \Bigl(
 \varphi_1^{\ast}(\gbigr_X)(\ast\gbigytilde)
 \otimes_{\nbigo_{\widetilde{\proj^1\times\cnum}}}
 \varphi_1^{\ast}(\gbigr_X)(\ast\gbigytilde)
 \Bigr)
 \langle
  \lambda^2\del_{\lambda},
  \deldel_{\tau}
 \rangle.
\]
Let $\nbigg_0^{\bullet}$ be
an $\nbiga_0$-injective resolution
of $\gbigr_X\otimes\omega_{\gbigx}^{-1}$.
Let $\nbigg_1^{\bullet}$ be an
$\nbiga_1$-injective resolution
 of $\varphi_1^{\ast}(\gbigr_X\otimes\omega_{\gbigx}^{-1})(\ast\gbigytilde)$.
We set
\[
 \varphi_1^{\ast}(\nbigg_0^{\bullet})(\ast\gbigytilde):=
 \nbigo_{\widetilde{\lefttop{\tau}\gbigx}}(\ast\gbigytilde)
 \otimes^{\ell}_{\varphi_1^{-1}(\nbigo_{\gbigx})}
 \varphi_1^{-1}(\nbigg_0^{\bullet})
 \otimes^{r}_{\varphi_1^{-1}(\nbigo_{\gbigx})}
  \nbigo_{\widetilde{\lefttop{\tau}\gbigx}}(\ast\gbigytilde). 
\]
Because
the naturally induced morphism
$\varphi_1^{\ast}(\gbigr_X\otimes\omega_{\gbigx}^{-1})(\ast\gbigytilde)
 \lrarr
 \varphi_1^{\ast}(\nbigg_0^{\bullet})(\ast\gbigytilde)$
is a quasi-isomorphism of
$\nbiga_1$-complexes,
there exists a quasi-isomorphism
$\varphi_1^{\ast}(\nbigg_0^{\bullet})(\ast\gbigytilde)
 \lrarr
\nbigg_1^{\bullet}$ of 
$\nbiga_1$-resolutions of
$\varphi_1^{\ast}(\gbigr_X\otimes\omega_{\gbigx}^{-1})(\ast\gbigytilde)$.

We have the following natural morphisms
of $\varphi_1^{\ast}(\gbigr_X)
 \langle\lambda^2\del_{\lambda},\deldel_{\tau}
 \rangle$-complexes:
\begin{multline}
\label{eq;21.3.25.1}
 \varphi_1^{\ast}\Bigl(
 \nhom_{\gbigr_X}\bigl(\gbigm,
 \nbigg^{\bullet}_0
 \bigr)
 \Bigr)(\ast\gbigytilde)
 \lrarr
 \nhom_{\varphi_1^{\ast}\gbigr_X}\Bigl(
 \varphi_1^{\ast}(\gbigm),\,
 \varphi_1^{\ast}(\nbigg^{\bullet}_0)
 \Bigr)(\ast\gbigytilde)
 \\
 \lrarr
 \nhom_{\varphi_1^{\ast}\gbigr_X(\ast\gbigytilde)}\Bigl(
 \varphi_1^{\ast}(\gbigm)(\ast\gbigytilde),\,
 \varphi_1^{\ast}(\nbigg^{\bullet}_0)(\ast\gbigytilde)
 \Bigr)
 \lrarr
  \nhom_{\varphi_1^{\ast}\gbigr_X(\ast\gbigytilde)}\Bigl(
 \varphi_1^{\ast}(\gbigm)(\ast\gbigytilde),\,
 \nbigg_1
 \Bigr).
\end{multline}
We can check that the composition of
the morphisms in (\ref{eq;21.3.25.1})
is a quasi-isomorphism
by using a resolution of $\gbigm$
by $\gbigr_X$-free modules of finite rank.
\hfill\qed

\vspace{.1in}
Note that
$\varphi_1^{\ast}(\gbigr_X)(\ast\gbigytilde)
=\varphi_0^{\ast}\bigl(\pi^{\ast}(\gbigr_X)(\ast\tau)\bigr)$
and
$\varphi_1^{\ast}(\gbigr_X)(\ast\gbigytilde)
\langle\lambda^2\del_{\lambda},\deldel_{\tau}\rangle
=\varphi_0^{\ast}(\gbigrtilde_{\lefttop{\tau}X}(\ast\tau))$.

\begin{lem}
\label{lem;21.6.18.2}
For any $\varphi_0^{\ast}(\gbigrtilde_{\lefttop{\tau}X}(\ast\tau))$-module
$\gbigm$
which is good coherent over 
$\varphi_0^{\ast}\bigl(\pi^{\ast}(\gbigr_X)(\ast\tau)\bigr)$,
there exists a natural isomorphism
 \[
 \varphi_{0\ast}\Bigl(
 \nrhom_{\varphi_0^{\ast}\pi^{\ast}\gbigr_X(\ast\tau)}
 \Bigl(
 \gbigm,
 \varphi_0^{\ast}\bigl(
 \pi^{\ast}(\gbigr_X\otimes\omega_{\gbigx}^{-1})(\ast\tau)
 \bigr)
 \Bigr)
 \Bigr)
 \simeq
 \nrhom_{\pi^{\ast}\gbigr_X(\ast\tau)}
 \bigl(
 \varphi_{0\ast}\gbigm,\,
 \pi^{\ast}(\gbigr_X\otimes\omega_{\gbigx}^{-1})(\ast\tau)
 \bigr).
\]
\end{lem}
\pf
We set
\[
 \nbiga_2:=
 \Bigl(
 \pi^{\ast}\gbigr_X(\ast\tau)
 \otimes_{\nbigo_{\proj^1\times\cnum_{\tau}}}
  \pi^{\ast}\gbigr_X(\ast\tau)
\Bigr)
  \langle
  \lambda^2\del_{\lambda},\deldel_{\tau}
  \rangle.
\]
Let $\nbigg_2^{\bullet}$ be
an $\nbiga_2$-injective resolution of
$\pi^{\ast}(\gbigr_X\otimes\omega_{\gbigx}^{-1})(\ast\tau)$.
Because
$R\varphi_{0\ast}(
\varphi_{0}^{\ast}
\pi^{\ast}\gbigr_X(\ast\tau)
)
=\pi^{\ast}\gbigr_X(\ast\tau)$,
there exists a morphism of complexes
$\varphi_{0\ast}(\nbigg_1^{\bullet})
\lrarr
 \nbigg_2^{\bullet}$.
We obtain
\begin{multline}
\label{eq;21.3.25.2}
 \varphi_{0\ast}\Bigl(
 \nhom_{\varphi_0^{\ast}\pi^{\ast}(\gbigr_X)(\ast\tau)}
 \bigl(
 \gbigm,
 \nbigg_1^{\bullet}
 \bigr)
 \Bigr)
 \lrarr 
 \nhom_{\pi^{\ast}(\gbigr_X)(\ast\tau)}
 \bigl(
 \varphi_{0\ast}(\gbigm),
 \varphi_{0\ast}(\nbigg_1^{\bullet})
 \bigr)
 \\
 \lrarr
  \nhom_{\pi^{\ast}(\gbigr_X)(\ast\tau)}
 \bigl(
 \varphi_{0\ast}(\gbigm),
 \nbigg_2^{\bullet}
 \bigr).
\end{multline}
We can check that the composition of (\ref{eq;21.3.25.2})
is a quasi-isomorphism
by using a free resolution of $\gbigm$
by free $\varphi_1^{\ast}\pi^{\ast}(\gbigr_X)(\ast\tau)$-modules of
finite rank.
\hfill\qed

\vspace{.1in}

By Lemma \ref{lem;21.6.18.1} and
Lemma \ref{lem;21.6.18.2},
for $\gbigm\in\gbigc(X)$,
there exists a natural isomorphism
\begin{equation}
\lefttop{\tau}(\DD_X\gbigm)
 \simeq
\nrhom_{\pi^{\ast}\gbigr_X(\ast\tau)}
\Bigl(
\lefttop{\tau}\gbigm,\,
\lambda^{d_X}
\pi^{\ast}(\gbigr_{X}\otimes\omega_{\gbigx}^{-1})(\ast\tau)
\Bigr)[d_X].
\end{equation}
By a similar argument similar to the proof of
Proposition \ref{prop;21.3.23.10},
there exists a natural isomorphism
\[
 \nrhom_{\pi^{\ast}\gbigr_X(\ast\tau)}
\Bigl(
\lefttop{\tau}\gbigm,\,
\pi^{\ast}(\gbigr_{X}\otimes\omega_{\gbigx}^{-1})(\ast\tau)
\Bigr)[d_X]
\simeq
 \nrhom_{\gbigr_{\lefttop{\tau}X}(\ast\tau)}
\Bigl(
\lefttop{\tau}\gbigm,\,
 \gbigr_{\lefttop{\tau}X}(\ast\tau)
 \otimes\omega_{\lefttop{\tau}\gbigx}^{-1}
\Bigr)[d_X+1].
\]
Thus, we obtain 
the desired isomorphism (\ref{eq;21.3.25.3}),
and Theorem \ref{thm;21.3.25.20} is proved.
\hfill\qed

\subsection{External tensor product}

Let $\gbigm_i\in\gbigc(X_i)$ $(i=1,2)$.
We obtain the $\gbigrtilde_{\lefttop{\tau}X_i}(\ast\tau)$-modules
$\lefttop{\tau}\gbigm_i$.
Let $p_i:\lefttop{\tau}(X_1\times X_2)\lrarr \lefttop{\tau}X_i$
denote the projections.
\begin{lem}
$p_i^{\ast}\lefttop{\tau}\gbigm_i$ $(i=1,2)$
are non-characteristic,
and they are flat over
$\nbigo_{\proj^1_{\lambda}\times\cnum_{\tau}\times X_j}$ $(j\neq i)$.
\end{lem}
\pf
The non-characteristic property is clear.
Let $(\lambda,\tau,P_1,P_2)\in
\proj^1_{\lambda}\times \cnum_{\tau}\times X_1\times X_2$.
Note that $\lefttop{\tau}\gbigm_1$
is flat over $\nbigo_{\proj^1_{\lambda}\times\cnum_{\tau}}$.
Let $A_i$ denote the stalk of
$\nbigo_{\proj^1_{\lambda}\times\cnum_{\tau}\times X_i}$
at $(\lambda,\tau,P_i)$.
Let $A_{0}$ denote the stalk of
$\nbigo_{\proj^1_{\lambda}\times\cnum_{\tau}}$
at $(\lambda,\tau)$.
Let $\Atilde$ denote the stalk of
$\nbigo_{\proj^1_{\lambda}\times\cnum_{\tau}\times X_1\times X_2}$
at $(\lambda,\tau,P_1,P_2)$.
The extension
$A_1\otimes_{A_0}A_2\lrarr\Atilde$
is faithfully flat.
Because $(\lefttop{\tau}\gbigm_1)_{(\lambda,\tau,P_1)}$
is flat over $A_0$,
$A_2\otimes_{A_0}(\lefttop{\tau}\gbigm_1)_{(\lambda,\tau,P_1)}$
is flat over $A_2$.
Then, we obtain that
$\Atilde\otimes\bigl(
A_2\otimes_{A_0}(\lefttop{\tau}\gbigm_1)_{(\lambda,\tau,P_1)}
\bigr)$
is flat over $A_2$.
\hfill\qed

\vspace{.1in}

By the construction,
there exist the following natural isomorphisms:
\[
 \lefttop{\tau}(\gbigm_1\boxtimes\gbigm_2)=
 p_1^{\ast}\lefttop{\tau}\gbigm_1
 \otimes_{\nbigo_{\proj^1_{\lambda}\times \cnum_{\tau}\times X_1\times X_2}}
 p_2^{\ast}\lefttop{\tau}\gbigm_2
\simeq
 p_1^{\ast}\lefttop{\tau}\gbigm_1
 \otimes^L_{\nbigo_{\proj^1_{\lambda}\times \cnum_{\tau}\times X_1\times X_2}}
 p_2^{\ast}\lefttop{\tau}\gbigm_2
\]
Note that
$p_i^{\ast}V_{c}(\lefttop{\tau}\nbigm_i)$
are flat over
$\nbigo_{\proj^1_{\lambda}\times\cnum_{\tau}\times X_j}$ $(j\neq i)$
as in the case of $p_i^{\ast}(\lefttop{\tau}\gbigm_i)$.
Hence, the natural morphisms
\[
 p_1^{\ast}\bigl(
 V_{b_1}(\lefttop{\tau}\gbigm_1)
 \bigr)
\otimes
 p_2^{\ast}\bigl(
  V_{b_2}(\lefttop{\tau}\gbigm_2)
  \bigr)
\lrarr
  p_1^{\ast}\bigl(
 V_{b_1}(\lefttop{\tau}\gbigm_1)
 \bigr)
\otimes
 p_2^{\ast}\bigl(
  \lefttop{\tau}\gbigm_2
  \bigr)
\lrarr
   p_1^{\ast}(
 \lefttop{\tau}\gbigm_1)
\otimes
 p_2^{\ast}(
  \lefttop{\tau}\gbigm_2)
\]
are monomorphisms.
We also note that
$p_1^{\ast}\bigl(
 V_{b_1}(\lefttop{\tau}\gbigm_1)
 \bigr)
\otimes
 p_2^{\ast}\bigl(
  V_{b_2}(\lefttop{\tau}\gbigm_2)
  \bigr)
=p_1^{\ast}\bigl(
 V_{b_1-1}(\lefttop{\tau}\gbigm_1)
 \bigr)
\otimes
 p_2^{\ast}\bigl(
  V_{b_2+1}(\lefttop{\tau}\gbigm_2)
  \bigr)$.

\begin{thm}
\label{thm;21.3.26.11}
 If $\gbigm_i\in\gbigc_{\res}(X_i)$,
then
$\gbigm_1\boxtimes\gbigm_2\in\gbigc_{\res}(X_1\times X_2)$.
For the $V$-filtrations along $\tau$,
we have
\begin{equation}
\label{eq;21.3.26.2}
 V_a\Bigl(
 \lefttop{\tau}\bigl(\gbigm_1\boxtimes\gbigm_2\bigr)
 \Bigr)
 =\sum_{b_1+b_2\leq a}
 p_1^{\ast}V_{b_1}\bigl(
 \lefttop{\tau}\gbigm_1
 \bigr)
 \otimes
 p_2^{\ast}V_{b_2}\bigl(
 \lefttop{\tau}\gbigm_2
 \bigr).
\end{equation}
\end{thm}
\pf
There exist
$\gbigm_{0,i}\in\gbigc(\lefttop{\tau}X_{i})$
such that
$\gbigm_{0,i}(\ast\tau)
=\lefttop{\tau}\gbigm_i$.
We obtain
$\gbigm_0:=\gbigm_{0,1}\boxtimes\gbigm_{0,2}
\in\gbigc(\lefttop{\tau_1}X_1\times\lefttop{\tau_2}X_2)$.
We have
$\gbigm_0(\ast(\tau_1\tau_2))
=\lefttop{\tau_1}\gbigm_1
\boxtimes
\lefttop{\tau_2}\gbigm_2$.
There exist $V$-filtration
$V_{\bullet}(\gbigm_0)$ along $\tau_1-\tau_2$
as an $\gbigr_{\lefttop{\tau_1}X_1\times\lefttop{\tau_2}X_2}$-module.
It induces a $V$-filtration
$V_{\bullet}\gbigm_0(\ast(\tau_1\tau_2))$
of
$\gbigm_0(\ast\tau_1\tau_2)$
as an
$\gbigr_{\lefttop{\tau_1}X_1\times\lefttop{\tau_2}X_2}
(\ast(\tau_1\tau_2))$-module.

\begin{lem}
\begin{equation}
\label{eq;21.3.26.1}
 V_a(\gbigm_0(\ast\tau_1\tau_2))
 =\left\{
\begin{array}{ll}
 (\tau_1-\tau_2)^{-a-1}\gbigm_0(\ast(\tau_1\tau_2))
  & (a\in\seisuu_{\leq -1})\\
 0 & (\mbox{\rm otherwise}).
\end{array}
 \right.
\end{equation}
\end{lem}
\pf
Let $V'_{a}(\gbigm_0(\ast\tau_1\tau_2))$
be the filtration determined by the right hand side of
(\ref{eq;21.3.26.1}).
Note that $\gbigm_0(\ast(\tau_1\tau_2))$
is
$\pi_1^{\ast}(\gbigr_{X_1\times X_2})
(\ast\tau_1\tau_2)$-coherent,
where
$\pi_1:\lefttop{\tau_1}X_1\times\lefttop{\tau_2}X_{2}
\lrarr X_1\times X_2$ denotes the projection.
Hence, we can easily check that
$V'_{a}(\gbigm_0(\ast\tau_1\tau_2))$
is a $V$-filtration.
\hfill\qed

\vspace{.1in}

Let $E_{\tau_1=\tau_2}$ be the hypersurface
of $\lefttop{\tau_1}X_1\times\lefttop{\tau_2}X_2$
determined by $\tau_1=\tau_2$.
Let $\iota_1:
\lefttop{\tau}(X_1\times X_2)\lrarr
\lefttop{\tau_1}X_1\times
\lefttop{\tau_2}X_2$
be the morphism induced by
the diagonal embedding
$\cnum_{\tau}\lrarr \cnum_{\tau_1}\times\cnum_{\tau_2}$.
It induces $\lefttop{\tau}(X_1\times X_2)\simeq E_{\tau_1=\tau_2}$.
We obtain
\begin{multline}
\lambda\Bigl(
 \gbigm_0[\ast E_{\tau_1=\tau_2}]\big/\gbigm_0
 \Bigr)(\ast(\tau_1\tau_2))
\simeq
 \iota_{1\dagger}\Bigl(
 \Gr^V_{-1}(\gbigm_0)
 \Bigr)(\ast\tau_1\tau_2)
=\iota_{1\dagger}\Bigl(
 \iota_1^{\ast}\bigl(
 \gbigm_0(\ast\tau_1\tau_2)
 \bigr)
 \Bigr)
 \simeq
 \iota_{1\dagger}\Bigl(
 \lefttop{\tau}(\gbigm_1\boxtimes\gbigm_2)
 \Bigr).
\end{multline}
Hence,
$\lefttop{\tau}(\gbigm_1\boxtimes\gbigm_2)
\in\gbigc(\lefttop{\tau}(X_1\times X_2))(\ast\tau)$.

\begin{lem}
For $c_i<b_i$ $(i=1,2)$,
we have
\begin{equation}
\label{eq;21.3.26.4}
 \Bigl(
p_1^{\ast}V_{b_1}(\lefttop{\tau}\gbigm_1)
\otimes
p_2^{\ast}V_{c_2}(\lefttop{\tau}\gbigm_2)
\Bigr)
\cap
\Bigl(
p_1^{\ast}V_{c_1}(\lefttop{\tau}\gbigm_1)
\otimes
p_2^{\ast}V_{b_2}(\lefttop{\tau}\gbigm_2)
\Bigr)
=
p_1^{\ast}V_{c_1}(\lefttop{\tau}\gbigm_1)
\otimes
p_2^{\ast}V_{c_2}(\lefttop{\tau}\gbigm_2).
\end{equation} 
\end{lem}
\pf
Clearly, the right hand side of
(\ref{eq;21.3.26.4}) is contained in the left hand side.
We have
\begin{multline}
 p_1^{\ast}\bigl(
 V_{b_1}(\lefttop{\tau}\gbigm_1)
 \bigr)
\otimes
 p_2^{\ast}\bigl(
  V_{b_2}(\lefttop{\tau}\gbigm_2)
  \bigr)
  \Big/
 p_1^{\ast}\bigl(
 V_{b_1-1}(\lefttop{\tau}\gbigm_1)
 \bigr)
\otimes
 p_2^{\ast}\bigl(
  V_{b_2}(\lefttop{\tau}\gbigm_2)
 \bigr)
 \simeq \\
  p_1^{\ast}\Bigl(
 V_{b_1}(\lefttop{\tau}\gbigm_1)
 \big/
 V_{b_1-1}(\lefttop{\tau}\gbigm_1)
 \Bigr)
\otimes
 p_2^{\ast}\bigl(
 V_{b_2}(\lefttop{\tau}\gbigm_2)
  \bigr)
 \simeq
   p_1^{\ast}\Bigl(
 V_{b_1}(\lefttop{\tau}\gbigm_1)
 \big/
 V_{b_1-1}(\lefttop{\tau}\gbigm_1)
 \Bigr)
\otimes
 p_2^{\ast}\Bigl(
 V_{b_2}(\lefttop{\tau}\gbigm_2)
 \big/
  V_{b_2-1}(\lefttop{\tau}\gbigm_2)
 \Bigr).
\end{multline}
Note that
$V_{b_i}(\lefttop{\tau}\gbigm_i)
 \big/
 V_{b_i-1}(\lefttop{\tau}\gbigm_i)$
is naturally an $\nbigr_{X_i}$-module.
Hence, we have
\begin{multline}
\label{eq;21.6.18.10}
 p_1^{\ast}\bigl(
 V_{b_1}(\lefttop{\tau}\gbigm_1)
 \bigr)
\otimes
 p_2^{\ast}\bigl(
  V_{b_2}(\lefttop{\tau}\gbigm_2)
  \bigr)
  \Big/
 p_1^{\ast}\bigl(
 V_{b_1-1}(\lefttop{\tau}\gbigm_1)
 \bigr)
\otimes
 p_2^{\ast}\bigl(
  V_{b_2}(\lefttop{\tau}\gbigm_2)
  \bigr)
 \simeq \\
 \Bigl(
 V_{b_1}(\lefttop{\tau}\gbigm_1)
 \big/
 V_{b_1-1}(\lefttop{\tau}\gbigm_1)
 \Bigr)
\boxtimes
 \Bigl(
 V_{b_2}(\lefttop{\tau}\gbigm_2)
 \big/
 V_{b_2-1}(\lefttop{\tau}\gbigm_2)
 \Bigr).
 \end{multline}
We also note that there exist canonical splittings
\begin{equation}
\label{eq;21.6.18.11}
\Bigl(
  V_{b_i}(\lefttop{\tau}\gbigm_i)
 \big/
 V_{b_i-1}(\lefttop{\tau}\gbigm_i)
\Bigr)_{|\cnum^{\ast}\times X_i}
 \simeq
 \bigoplus_{b_i-1<a\leq b_i}
 \Gr^V_{a}(\lefttop{\tau}\gbigm_i)_{|\cnum^{\ast}\times X_i}.
\end{equation}

Let $b_1-1\leq c_1<b_1$ and $b_2-1\leq c_2<b_2$.
By (\ref{eq;21.6.18.10}) and (\ref{eq;21.6.18.11}),
we have
\begin{multline}
\label{eq;21.6.18.13}
\Bigl(
p_1^{\ast}V_{b_1}(\lefttop{\tau}\gbigm_1)
\otimes
p_2^{\ast}V_{c_2}(\lefttop{\tau}\gbigm_2)
\Bigr)_{|\cnum^{\ast}\times\lefttop{\tau}(X_1\times X_2)}
\cap
\Bigl(
p_1^{\ast}V_{c_1}(\lefttop{\tau}\gbigm_1)
\otimes
p_2^{\ast}V_{b_2}(\lefttop{\tau}\gbigm_2)
\Bigr)_{|\cnum^{\ast}\times\lefttop{\tau}(X_1\times X_2)}
= \\
\Bigl(
p_1^{\ast}V_{c_1}(\lefttop{\tau}\gbigm_1)
\otimes
p_2^{\ast}V_{c_2}(\lefttop{\tau}\gbigm_2)
\Bigr)_{|\cnum^{\ast}\times\lefttop{\tau}(X_1\times X_2)}.
\end{multline}
We have the following isomorphism
which is similar to (\ref{eq;21.6.18.10}):
\begin{multline}
 p_1^{\ast}\bigl(
 V_{b_1}(\lefttop{\tau}\gbigm_1)
 \bigr)
\otimes
 p_2^{\ast}\bigl(
  V_{b_2}(\lefttop{\tau}\gbigm_2)
  \bigr)
  \Big/
 p_1^{\ast}\bigl(
 V_{c_1}(\lefttop{\tau}\gbigm_1)
 \bigr)
\otimes
 p_2^{\ast}\bigl(
  V_{b_2}(\lefttop{\tau}\gbigm_2)
  \bigr)
 \simeq \\
 \Bigl(
 V_{b_1}(\lefttop{\tau}\gbigm_1)
 \big/
 V_{c_1}(\lefttop{\tau}\gbigm_1)
 \Bigr)
\boxtimes
 \Bigl(
 V_{b_2}(\lefttop{\tau}\gbigm_2)
 \big/
 V_{b_2-1}(\lefttop{\tau}\gbigm_2)
 \Bigr).
 \end{multline}
We also have
\begin{multline}
 p_1^{\ast}\bigl(
 V_{c_1}(\lefttop{\tau}\gbigm_1)
 \bigr)
\otimes
 p_2^{\ast}\bigl(
  V_{b_2}(\lefttop{\tau}\gbigm_2)
  \bigr)
  \Big/
 p_1^{\ast}\bigl(
 V_{c_1}(\lefttop{\tau}\gbigm_1)
 \bigr)
\otimes
 p_2^{\ast}\bigl(
  V_{c_2}(\lefttop{\tau}\gbigm_2)
  \bigr)
 \simeq \\
 \Bigl(
 V_{c_1}(\lefttop{\tau}\gbigm_1)
 \big/
 V_{c_1-1}(\lefttop{\tau}\gbigm_1)
 \Bigr)
\boxtimes
 \Bigl(
 V_{b_2}(\lefttop{\tau}\gbigm_2)
 \big/
 V_{c_2}(\lefttop{\tau}\gbigm_2)
 \Bigr).
\end{multline}
Hence, we obtain that
\[
  p_1^{\ast}\bigl(
 V_{b_1}(\lefttop{\tau}\gbigm_1)
 \bigr)
\otimes
 p_2^{\ast}\bigl(
  V_{b_2}(\lefttop{\tau}\gbigm_2)
  \bigr)
  \Big/
 p_1^{\ast}\bigl(
 V_{c_1}(\lefttop{\tau}\gbigm_1)
 \bigr)
\otimes
 p_2^{\ast}\bigl(
  V_{c_2}(\lefttop{\tau}\gbigm_2)
  \bigr)
\]
is strict.
Hence, we obtain (\ref{eq;21.3.26.4})
from (\ref{eq;21.6.18.13})
in the case $b_i-1<c_i<b_i$.
By an easy inductive argument,
we obtain (\ref{eq;21.3.26.4})
for any $c_i<b_i$.
\hfill\qed

\begin{lem}
\label{lem;21.3.26.10}
 \begin{multline}
 \Bigl(
 p_1^{\ast}V_{b_1}(\lefttop{\tau}\gbigm_1)
 \otimes
 p_2^{\ast}V_{b_2}(\lefttop{\tau}\gbigm_2)
 \Bigr)
 \cap
 \left(
 \sum_{c_1+c_2<b_1+b_2}
 p_1^{\ast}V_{c_1}(\lefttop{\tau}\gbigm_1)
 \otimes
 p_2^{\ast}V_{c_2}(\lefttop{\tau}\gbigm_2)
 \right)
 =\\
 \sum_{\substack{c_1+c_2<b_1+b_2\\ c_i\leq b_i}}
 p_1^{\ast}V_{c_1}(\lefttop{\tau}\gbigm_1)
 \otimes
 p_2^{\ast}V_{c_2}(\lefttop{\tau}\gbigm_2).
 \end{multline}
\end{lem}
\pf
Let $(c^{(k)}_1,c^{(k)}_2)$ $(k=1,\ldots,\ell)$
satisfy $c^{(k)}_i+c^{(k)}_2<b_1+b_2$.
Let $m_k$ be sections of
$p_1^{\ast}V_{c^{(k)}_1}(\lefttop{\tau}\gbigm_1)
 \otimes
 p_2^{\ast}V_{c^{(k)}_2}(\lefttop{\tau}\gbigm_2)$.
Suppose that
$m=\sum m_k$ is a section of
$p_1^{\ast}V_{b_1}(\lefttop{\tau}\gbigm_1)
 \otimes
 p_2^{\ast}V_{b_2}(\lefttop{\tau}\gbigm_2)$.

We may assume that
$c^{(1)}_1>c^{(k)}_1$ $(k\neq 1)$.
Suppose that $c^{(k)}_1>b_1$.
Let $d:=\max\{b_1,c^{(k)}_1\,\,k\neq 1\}$.
Then,
$m_1$ is contained in
\[
\Bigl(
 p_1^{\ast}V_{c^{(1)}_1}(\lefttop{\tau}\gbigm_1)
 \otimes
 p_2^{\ast}V_{c^{(1)}_2}
 (\lefttop{\tau}\gbigm_2)
 \Bigr)
\cap
\Bigl(
 p_1^{\ast}V_{d}(\lefttop{\tau}\gbigm_1)
 \otimes
 p_2^{\ast}V_{N}
 (\lefttop{\tau}\gbigm_2)
 \Bigr)
\]
for a large $N$.
Hence, $m_1$ is a section of
$p_1^{\ast}V_{d}(\lefttop{\tau}\gbigm_1)
  \otimes
   p_2^{\ast}V_{c^{(1)}_2}
 (\lefttop{\tau}\gbigm_2)$.
By an easy inductive argument,
we can prove that
there exist
$(c^{\prime(k)}_1,c^{\prime(k)}_2)\in\real^2$
$(k=1,\ldots,m')$
satisfying $c^{(\prime(k))}_1\leq b_1$
and $c^{\prime(k)}_1+c^{\prime(k)}_2<b_1+b_2$,
and that there exist sections
$m_k'$ of
$p_1^{\ast}V_{c^{\prime(k)}_1}
\otimes
p_2^{\ast}V_{c^{\prime(k)}_2}$
such that 
$m=\sum m'_k$.
By applying a similar argument to the second component,
we can obtain the claim of Lemma \ref{lem;21.3.26.10}.
\hfill\qed

\vspace{.1in}
By a similar argument, we can prove the following lemma.
\begin{lem}
\label{lem;21.6.18.20}
 For any $-1<b_1\leq 0$, we have
\[
 \Bigl(
 p_1^{\ast}V_{b_1}(\lefttop{\tau}\gbigm_1)
 \otimes
 p_2^{\ast}V_{b_2}(\lefttop{\tau}\gbigm_2)
 \Bigr)
\cap
  \left(
 \sum_{\substack{c_1+c_2=b_1+b_2\\
  -1<c_1<b_1}}
 p_1^{\ast}V_{c_1}(\lefttop{\tau}\gbigm_1)
 \otimes
 p_2^{\ast}V_{c_2}(\lefttop{\tau}\gbigm_2)
 \right)
\subset
 \sum_{\substack{c_1<b_1\\ c_2\leq b_2}}
 p_1^{\ast}V_{c_1}(\lefttop{\tau}\gbigm_1)
 \otimes
 p_2^{\ast}V_{c_2} (\lefttop{\tau}\gbigm_2).
\]
\hfill\qed
\end{lem}

Let $V'_a\Bigl(\lefttop{\tau}\bigl(\gbigm_1\boxtimes\gbigm_2\bigr)\Bigr)$
be defined by the right hand side of
(\ref{eq;21.3.26.2}).
By Lemma \ref{lem;21.3.26.10}
and Lemma \ref{lem;21.6.18.20},
we obtain
\begin{equation}
 \Gr^{V'}_a\Bigl(\lefttop{\tau}\bigl(\nbigm_1\boxtimes\nbigm_2\bigr)\Bigr)
  \simeq
  \bigoplus_{
  \substack{b_1+b_2=a\\ -1<b_1\leq 0}}
  \Gr^V_{b_1}\Bigl(
  \lefttop{\tau}\gbigm_1
  \Bigr)
  \boxtimes
  \Gr^V_{b_2}\Bigl(
  \lefttop{\tau}\gbigm_2
  \Bigr).
\end{equation}
In particular, it is strict.
Then, we can easily check that
$V'$ is a $V$-filtration,
and we obtain Theorem \ref{thm;21.3.26.11}.
\hfill\qed

\begin{cor}
\[
\iota_{\lambda=\tau}^{\ast}
 V_a\Bigl(
 \lefttop{\tau}\bigl(
 \gbigm_1\boxtimes\gbigm_2\bigr)
 \Bigr)
 =\sum_{b_1+b_2\leq a}
 \iota_{\lambda=\tau}^{\ast}\Bigl(
 V_{b_1}\bigl(
 \lefttop{\tau}\gbigm_1
 \bigr)
\Bigr)
\boxtimes
\iota_{\lambda=\tau}^{\ast}\Bigl(
 V_{b_2}\bigl(
 \lefttop{\tau}\gbigm_2
 \bigr)
 \Bigr).
\]
\end{cor}
\pf
Because
$\iota_{\lambda=\tau}^{\ast}
 \Bigl(
 p_1^{\ast}V_{b_1}
 \otimes
 p_2^{\ast}V_{b_2}
 \Bigr)
=\iota_{\lambda=\tau}^{\ast}
 p_1^{\ast}V_{b_1}
 \boxtimes
 \iota_{\lambda=\tau}^{\ast}
 p_2^{\ast}V_{b_2}$,
the claim follows from Theorem \ref{thm;21.3.26.11}.
\hfill\qed

\begin{cor}
\label{cor;21.6.29.22}
For $\gbigm_i\in\gbigc(X_i)$,
we have
\[
 F^{\irr}_{a}\Bigl(
 \Xi_{\DR}(\gbigm_1)\boxtimes
 \Xi_{\DR}(\gbigm_2)
 \Bigr)
=\sum_{b_1+b_2\leq a}
 F^{\irr}_{b_1}\Xi_{\DR}(\gbigm_1)\boxtimes
 F^{\irr}_{b_2}\Xi_{\DR}(\gbigm_2).
\]
\hfill\qed
\end{cor}

\begin{cor}
In the situation of {\rm\S\ref{subsection;21.6.22.11}}
and {\rm\S\ref{subsection;21.6.22.18}},
for $\gbigm_i\in\gbigc_{\res}(X;H)$,
$\nbigh^k(\gbigm_1\otimes^{\star}\gbigm_2)$
are objects of $\gbigc_{\res}(X;H)$.
\hfill\qed
\end{cor}

\subsection{Non-characteristic inverse image}

Let $f:X\lrarr Y$ be a morphism of complex manifolds.
Let $H_Y$ be a hypersurface of $Y$.
We set $H_X:=f^{-1}(H_Y)$.
Let $\gbigm\in\gbigc_{\res}(Y;H_Y)$.

\begin{prop}
\label{prop;21.6.29.23}
Suppose that $f$ is non-characteristic for $\gbigm$.
Then,  
$f^{\ast}(\gbigm)$
is an object of $\gbigc_{\res}(X;H_X)$.
Moreover, we have
\[
 F^{\irr}_{\bullet}f^{\ast}(\gbigm^{\lambda})
 =f^{\ast}\bigl(
 F^{\irr}_{\bullet}(\gbigm^{\lambda})
 \bigr),
\quad
 F^{\irr}_{\bullet}f^{\ast}
 \bigl(\iota_{\infty}^{\ast}(\lefttop{\infty}V_a\gbigm)\bigr)
=f^{\ast}\bigl(
 F^{\irr}_{\bullet}
 \bigl(\iota_{\infty}^{\ast}(\lefttop{\infty}V_a\gbigm^{\lambda})\bigr)
 \bigr),
\]
\[
 F^{\irr}_{\bullet}f^{\ast}(V_a\lefttop{\tau}\gbigm)
=f^{\ast}\bigl(
 F^{\irr}_{\bullet}(V_a\lefttop{\tau}\gbigm)
 \bigr).
\]
\end{prop}
\pf
The induced morphism
$f_0:\lefttop{\tau}X
\lrarr \lefttop{\tau}Y$
is non-characteristic for
$\lefttop{\tau}\gbigm
\in
\gbigc\bigl(\lefttop{\tau}Y;
 \lefttop{\tau}Y_0\cup
 \lefttop{\tau}H_Y
 \bigr)$.
We have
$f_0^{\ast}
 \bigl(
 \lefttop{\tau}\gbigm
 \bigr)
=\lefttop{\tau}(f^{\ast}(\gbigm))$.
Hence, we obtain the first claim of the proposition.

Let us prove the second claim.
It is enough to prove the case where
$f$ is a closed embedding.
We have only to prove the claim locally around
any point of $f(X)$.
It is enough to consider the case
where $f(X)$ is a smooth hypersurface of $Y$
defined by a coordinate function $t$.

Note that $t$ is non-characteristic for
$\Gr^V_a\bigl(\lefttop{\tau}\gbigm\bigr)$.
Hence, the multiplication by $t$ on 
$\Gr^V_a\bigl(\lefttop{\tau}\gbigm\bigr)$
is a monomorphism.
It follows that
$f_0^{\ast}\bigl(
 V_b\bigl(\lefttop{\tau}\gbigm \bigr)
 \bigr)
 \lrarr
 f_0^{\ast}\bigl(
 \lefttop{\tau}\gbigm
 \bigr)$
$(b\in\real)$
are monomorphisms,
and the images induce a $V$-filtration of
$f_0^{\ast}\bigl(
 \lefttop{\tau}\gbigm
 \bigr)$.
Then, we obtain the claim for the irregular Hodge filtrations.
\hfill\qed

\section{Irregular Hodge filtrations in the regular case}

\subsection{Strict specializability of
integrable mixed twistor $\nbigd$-modules in two directions}

\subsubsection{Strict $(t,\tau)$-specializability
of the underlying $\nbigrtilde$-modules}
\label{subsection;22.7.9.3}

Let $X$ be any complex manifold.
Let $Y$ be an open subset in
$\cnum_t\times\cnum_{\tau}\times X$.
Let $\pi:Y\lrarr X$ denote the projection.
We set
$\lefttop{t}V\nbigr_Y:=
 \pi^{\ast}\nbigr_X\langle t\deldel_t,\deldel_{\tau}\rangle
 \subset\nbigr_Y$
and 
$\lefttop{t}V\nbigrtilde_Y:=
 \lefttop{t}V\nbigr_Y\langle\lambda^2\del_{\lambda}\rangle
 \subset\nbigrtilde_Y$.
 Similarly,
we set 
 $\lefttop{\tau}V\nbigr_Y:=
 \pi^{\ast}\nbigr_X\langle \deldel_t,\tau\deldel_{\tau}\rangle
 \subset\nbigr_Y$
and 
$\lefttop{\tau}V\nbigrtilde_Y:=
 \lefttop{\tau}\nbigrtilde_Y\langle\lambda^2\del_{\lambda}\rangle
 \subset\nbigrtilde_Y$.
We also set
$\lefttop{t,\tau}V\nbigr_Y:=
 \pi^{\ast}\nbigr_X\langle t\deldel_t,\tau\deldel_{\tau}\rangle
 \subset\nbigr_Y$
and
$\lefttop{t,\tau}V\nbigrtilde_Y:=
\lefttop{t,\tau}V\nbigr_Y\langle\lambda^2\del_{\lambda}\rangle
\subset\nbigrtilde_Y$.

For $\nbigm\in\nbigc(Y)$,
there exist $V$-filtrations
$\lefttop{t}V_{\bullet}(\nbigm)$
(resp. $\lefttop{\tau}V_{\bullet}(\nbigm)$)
along $t$ (resp. $\tau$).

\begin{df}
\label{df;21.3.18.10}
 A $V$-filtration of $\nbigm$ along $(t,\tau)$
is a tuple of
$\lefttop{t,\tau}V\nbigrtilde_{Y}$-submodules
$\bigl\{
 \lefttop{t,\tau}V_{a,b}(\nbigm)\,\big|\,
 a,b\in\real
 \bigr\}$ of $\nbigm$
satisfying the following conditions.
\begin{enumerate}
 \item $\lefttop{t,\tau}V_{a,b}(\nbigm)$
       are coherent over $\lefttop{t,\tau}V\nbigr_Y$.
 \item
      $\lefttop{t,\tau}V_{a,b}(\nbigm)
      \subset
      \lefttop{t,\tau}V_{a',b'}(\nbigm)$
      for $a\leq a'$ and $b\leq b'$.
 \item $\lefttop{t,\tau}V_{a,b}(\nbigm)(\ast (t\tau))
       =\nbigm(\ast(t\tau))$.
 \item $\bigcup_{a}\lefttop{t,\tau}V_{a,b}(\nbigm)
       =\lefttop{\tau}V_b(\nbigm)$
       and
       $\bigcup_{b}\lefttop{t,\tau}V_{a,b}(\nbigm)
       =\lefttop{t}V_a(\nbigm)$.
 \item For any $(a,b)\in\real^2$ and any compact subset
       $K\subset \cnum_{\lambda}\times Y$,
       there exists $\epsilon>0$ such that
       $V_{a,b}(\nbigm)_{|K}=V_{a+\epsilon,b+\epsilon}(\nbigm)_{|K}$.
 \item $t\cdot \lefttop{t,\tau}V_{a,b}\subset \lefttop{t,\tau}V_{a-1,b}$
       and
       $\deldel_t\cdot \lefttop{t,\tau}V_{a,b}(\nbigm)\subset
       \lefttop{t,\tau}V_{a+1,b}(\nbigm)$.
 \item $\tau\cdot \lefttop{t,\tau}V_{a,b}(\nbigm)
       \subset \lefttop{t,\tau}V_{a,b-1}(\nbigm)$
       and
       $\deldel_{\tau}\cdot \lefttop{t,\tau}V_{a,b}(\nbigm)
        \subset \lefttop{t,\tau}V_{a,b+1}(\nbigm)$.
 \item For any $a,b\in\real$, the natural morphism
       $\lefttop{t,\tau}V_{a,b}(\nbigm)\big/
       \lefttop{t,\tau}V_{<a,b}(\nbigm)
       \lrarr
       \lefttop{t}\Gr^V_a(\nbigm)$
       induces an isomorphism
       \begin{equation}
	\label{eq;21.3.16.50}
       \lefttop{t,\tau}V_{a,b}(\nbigm)\big/
       \lefttop{t,\tau}V_{<a,b}(\nbigm)
       \simeq
       \lefttop{\tau}V_b\bigl(
       \lefttop{t}\Gr^V_a(\nbigm)
       \bigr).
       \end{equation}
 \item  For any $a,b\in\real$,
       the natural morphism
       $\lefttop{t,\tau}V_{a,b}\big/
       \lefttop{t,\tau}V_{a,<b}
       \lrarr
       \lefttop{\tau}\Gr^V_b(\nbigm)$
       induces an isomorphism
       \begin{equation}
	\label{eq;21.6.19.1}
       \lefttop{t,\tau}V_{a,b}\big/
       \lefttop{t,\tau}V_{a,<b}
       \simeq
       \lefttop{t}V_a\bigl(
       \lefttop{\tau}\Gr^V_b(\nbigm)
       \bigr).
       \end{equation}
 \end{enumerate}
Here, we set
$\lefttop{t,\tau}V_{<a,b}:=\sum_{c<a}\lefttop{t,\tau}V_{c,b}$
and
$\lefttop{t,\tau}V_{a,<b}:=\sum_{c<b}\lefttop{t,\tau}V_{a,c}$,
and 
$\lefttop{\tau}V_{\bullet}\bigl(
       \lefttop{t}\Gr^V_a(\nbigm)
       \bigr)$
(resp. $\lefttop{t}V_{\bullet}\bigl(
       \lefttop{\tau}\Gr^V_b(\nbigm)
       \bigr)$)
denotes the $V$-filtration of
$\lefttop{t}\Gr^V_a(\nbigm)\in\nbigc(Y\cap\{t=0\})$
(resp. $\lefttop{\tau}\Gr^V_b(\nbigm)\in\nbigc(Y\cap\{\tau=0\})$)
along $\tau$ (resp. $t$).

 If there exists a $V$-filtration of $\nbigm$ along $(t,\tau)$,
$\nbigm$ is called strictly $(t,\tau)$-specializable.
 \hfill\qed
\end{df}

\begin{lem}
\label{lem;21.3.18.20}
Let $\lefttop{t,\tau}V_{\bullet,\bullet}=
\{\lefttop{t,\tau}V_{a,b}(\nbigm)\,|\,a,b\in\real\}$
be a tuple of $\lefttop{t,\tau}V\nbigrtilde_Y$-submodules of $\nbigm$
satisfying the conditions {\rm 1--7} in Definition {\rm\ref{df;21.3.18.10}}.
If it satisfies the condition {\rm 8} (resp. {\rm 9}),
it also satisfies the condition {\rm 9} (resp. {\rm 8}),
i.e., $\lefttop{t,\tau}V_{\bullet,\bullet}$
is a $V$-filtration of $\nbigm$ along $(t,\tau)$.
\end{lem}
\pf
Suppose that the condition 8 is satisfied.
By (\ref{eq;21.3.16.50}),
we obtain
$\lefttop{t,\tau}V_{a,b}(\nbigm)\cap
 \lefttop{t}V_{<a}(\nbigm)=
 \lefttop{t,\tau}V_{<a,b}(\nbigm)$.
 Hence, for any $c<a$,
 we obtain 
 $\lefttop{t,\tau}V_{a,b}(\nbigm)\cap
 \lefttop{t}V_{c}(\nbigm)=
 \lefttop{t,\tau}V_{c,b}(\nbigm)$.
Let us prove that
$\lefttop{t,\tau}V_{a,b}(\nbigm)\cap
 \lefttop{\tau}V_{<b}(\nbigm)=
 \lefttop{t,\tau}V_{a,<b}(\nbigm)$.
Clearly, we have
$\lefttop{t,\tau}V_{a,b}(\nbigm)\cap
\lefttop{\tau}V_{<b}(\nbigm)\supset
\lefttop{t,\tau}V_{a,<b}(\nbigm)$.
Let $s$ be a local section of
$\lefttop{t,\tau}V_{a,b}(\nbigm)\cap
\lefttop{\tau}V_{<b}(\nbigm)$.
There exists $a'\geq a$ and $c<b$ such that
$s\in \lefttop{t,\tau}V_{a',c}(\nbigm)$.
Because
$\lefttop{t,\tau}V_{a',c}(\nbigm)
\cap \lefttop{t}V_{a}(\nbigm)
=\lefttop{t,\tau}V_{a,c}(\nbigm)$,
we obtain that $s$ is a local section of
$\lefttop{t,\tau}V_{a,c}
\subset
\lefttop{t,\tau}V_{a,<b}$.

As a result,
for any $b\in\real$,
we obtain monomorphisms
\begin{equation}
\label{eq;21.3.18.12}
\lefttop{t,\tau}V_{a,b}\big/
\lefttop{t,\tau}V_{a,<b}
\lrarr
 \lefttop{\tau}\Gr^V_{b}(\nbigm).
\end{equation}
Let
$\lefttop{t}V'_a(\lefttop{\tau}\Gr^V_b(\nbigm))$
denote the image of (\ref{eq;21.3.18.12}).
By the construction,
we obtain
\begin{equation}
 \label{eq.21.3.18.13}
\lefttop{t}\Gr^{V'}_a(\lefttop{\tau}\Gr^V_b(\nbigm))
\simeq 
 \lefttop{t,\tau}V_{a,b}\Big/
 \sum_{\substack{c\leq a\\ d\leq b\\ (c,d)\neq (a,b)}}
 \lefttop{t,\tau}V_{c,d}
 \simeq
 \lefttop{\tau}\Gr^V_b\lefttop{t}\Gr^V_a(\nbigm).
\end{equation}
Here, we obtain the second isomorphism
from the condition 8.
It implies that
$\lefttop{t}\Gr^{V'}_a(\lefttop{\tau}\Gr^V_b(\nbigm))$
is strict for any $a$.
It is easy to see that
$\lefttop{t}V'_a(\lefttop{\tau}\Gr^V_b(\nbigm))$
are coherent over
$\lefttop{t}V\nbigr_{Y\cap\{\tau=0\}}$.
By the construction,
we have
$t\bigl(
\lefttop{t}V'_a(\lefttop{\tau}\Gr^V_b(\nbigm))
\bigr)
\subset
\lefttop{t}V'_{a-1}(\lefttop{\tau}\Gr^V_b(\nbigm))$
and 
$t\bigl(
\lefttop{t}V'_a(\lefttop{\tau}\Gr^V_b(\nbigm))
\bigr)
\subset
\lefttop{t}V'_{a-1}(\lefttop{\tau}\Gr^V_b(\nbigm))$.
Moreover,
because of (\ref{eq.21.3.18.13}),
$-\deldel_tt-a\lambda$ are locally nilpotent on
$\lefttop{t}\Gr^{V'}_a(\lefttop{\tau}\Gr^V_b(\nbigm))$
for any $a\in\real$.
By the standard argument,
we can prove that
the filtration
$\lefttop{t}V'_{\bullet}(\lefttop{\tau}\Gr^V_b(\nbigm))$
is the $V$-filtration of
$\lefttop{\tau}\Gr^V_b(\nbigm)$.
(See \cite{mochi2, Sabbah-pure-twistor}.)
Thus, we obtain (\ref{eq;21.6.19.1}).
\hfill\qed

\begin{lem}
\label{lem;22.7.9.1}
Let $\lefttop{t,\tau}V_{\bullet,\bullet}$
be a $V$-filtration of $\nbigm$ along $(t,\tau)$.
For $a<0$ and $b\in\real$, we have
 $t\cdot\bigl(
 \lefttop{t,\tau}V_{a,b}(\nbigm)\bigr)
=\lefttop{t,\tau}V_{a-1,b}(\nbigm)$.
For $a>{-1}$ and $b\in\real$,
the morphism
$\deldel_t:
 \lefttop{t,\tau}V_{a,b}/\lefttop{t,\tau}V_{<a,b}
 \lrarr
 \lefttop{t,\tau}V_{a+1,b}/\lefttop{t,\tau}V_{<a+1,b}$
is an isomorphism of
$\lefttop{\tau}V\nbigrtilde_{\{t=0\}\cap Y}$-modules.
Similar clams also hold
when we exchange the roles of $t$ and $\tau$.
\end{lem}
\pf
If $a<0$,
we have
$t\cdot \lefttop{t}V_{a}(\nbigm)=\lefttop{t}V_{a-1}(\nbigm)$.
For $v\in \lefttop{t,\tau}V_{a-1,b}(\nbigm)$,
there exists $c\in\real$ and $v'\in \lefttop{t,\tau}V_{a,c}(\nbigm)$
such that $tv'=v$.
We obtain the induced local section  $[v']$ of
$\lefttop{t,\tau}V_{a,c}\big/\lefttop{t,\tau}V_{a,<c}
\simeq
\lefttop{t}V_{a}\lefttop{\tau}\Gr^V_c(\nbigm)$.
If $c>b$, we obtain $[tv']=0$ in
$\lefttop{t,\tau}V_{a-1,c}\big/\lefttop{t,\tau}V_{a-1,<c}
\simeq
\lefttop{t}V_{a-1}\lefttop{\tau}\Gr^V_c(\nbigm)$.
Because the multiplication of $t$ induces
an isomorphism of $\nbigo_{\nbigy}$-modules
$\lefttop{t}V_a\lefttop{\tau}\Gr^V_c(\nbigm)\simeq
\lefttop{t}V_{a-1}\lefttop{\tau}\Gr^V_c(\nbigm)$ $(a<0)$,
we obtain that $[v']=0$ in 
$\lefttop{t,\tau}V_{a,c}\big/\lefttop{t,\tau}V_{a,<c}$.
There exists a finite sequence
$b<c_m<\cdots<c_1<c_0=c$ such that
$\lefttop{\tau}\Gr^V_{c'}=0$ for $b< c'\leq c$
unless $c'\in\{c_i\}$.
Hence, we obtain $v'\in \lefttop{t,\tau}V_{a,b}$,
that is the first claim.

If $a>0$,
$\deldel_t$ induces an isomorphism
of $\nbigrtilde_{Y\cap\{t=0\}}$-modules
$\lefttop{t}\Gr^V_{a}(\nbigm)
 \lrarr
 \lefttop{t}\Gr^V_{a+1}(\nbigm)$.
Hence, the second claim follows from the isomorphisms (\ref{eq;21.3.16.50}).
\hfill\qed

\begin{lem}
\label{lem;22.7.9.2}
If $\lefttop{t,\tau}V_{\bullet,\bullet}$
 is a $V$-filtration of $\nbigm$ along $(t,\tau)$,
 the following holds.
\begin{itemize}
 \item If $a\geq a'$, we obtain
       $\lefttop{t,\tau}V_{a,b}\cap
       \lefttop{t}V_{a'}=\lefttop{t,\tau}V_{a',b}$.
       Similarly, if $b\geq b'$, we obtain
      $\lefttop{t,\tau}V_{a,b}\cap
       \lefttop{\tau}V_{b'}=\lefttop{t,\tau}V_{a,b'}$.
 \item
 $\lefttop{t,\tau}V_{a,b}
 =\lefttop{t}V_a\cap\lefttop{\tau}V_b$.       
\end{itemize}
\end{lem}
\pf
We have already observed the first claim
in the proof of Lemma \ref{lem;21.3.18.20}.
As for the second claim,
we clearly have
$\lefttop{t,\tau}V_{a,b}\subset
\lefttop{t}V_a\cap\lefttop{\tau}V_b$.
Let $v$ be a local section of
$\lefttop{t}V_a\cap\lefttop{\tau}V_b$.
There exists $a'\geq a$ and $b'\geq b$
such that
$v\in \lefttop{t,\tau}V_{a',b'}$.
Because
$\lefttop{t,\tau}V_{a',b'}\cap \lefttop{t}V_a\cap \lefttop{\tau}V_b
=\lefttop{t,\tau}V_{a,b'}
\cap
 \lefttop{\tau}V_{b}
=\lefttop{t,\tau}V_{a,b}$,
we obtain $v\in \lefttop{t,\tau}V_{a,b}$.
Thus, the second claim is proved.
\hfill\qed

\begin{cor}
If $\nbigm\in\nbigc(Y)$ is strictly $(t,\tau)$-specializable,
a $V$-filtration of $\nbigm$ along $(t,\tau)$
is uniquely determined.
\hfill\qed
\end{cor}

\begin{prop}
\label{prop;21.3.17.10}
 Let $F:\nbigm_1\lrarr\nbigm_2$ be a morphism
 in $\nbigc(Y)$
 such that
 $\Ker(F)$,
 $\Image(F)$ and $\Cok(F)$
 are also objects of $\nbigc(Y)$.
Assume that $\nbigm_i$ are strictly $(t,\tau)$-specializable.
\begin{itemize}
 \item $F$ is strict with respect to
       $\lefttop{t,\tau}V_{\bullet,\bullet}$,
       i.e.,
       $F(V_{a,b}(\nbigm_1))=\Image(F)\cap V_{a,b}(\nbigm_2)$
       for any $a,b\in\real$.
 \item $\Ker(F)$, $\Image(F)$ and $\Cok(F)$
       are also strictly $(t,\tau)$-specializable.
       The $V$-filtrations
       of $\Ker(F)$, $\Image(F)$ and $\Cok(F)$ along $(t,\tau)$
       are equal to the filtrations
       naturally induced by $\lefttop{t,\tau}V(\nbigm_i)$.
\end{itemize}
\end{prop}
\pf
It is enough to prove the claim locally around
any point of $\{(t,\tau)=(0,0)\}\cap Y$.
Recall that
$F$ is strict with respect to the $V$-filtrations $\lefttop{t}V$
of $\nbigm_i$ along $t$,
and the $V$-filtrations of $\Ker(F)$,
$\Image(F)$ and $\Cok(F)$ are the filtrations
induced by $\lefttop{t}V(\nbigm_i)$.
Similar claims hold for the $V$-filtrations along $\tau$.

Clearly,
$F(\lefttop{t,\tau}V_{a,b}(\nbigm_1))
\subset
\lefttop{t,\tau}V_{a,b}(\nbigm_2)$ holds.
The induced morphism
\begin{equation}
\label{eq;21.3.17.1}
 \lefttop{t,\tau}V_{a,b}(\nbigm_1)
 \big/
 \lefttop{t,\tau}V_{a,<b}(\nbigm_1)
 \lrarr
 \lefttop{t,\tau}V_{a,b}(\nbigm_2)
 \big/
 \lefttop{t,\tau}V_{a,<b}(\nbigm_2)
\end{equation}
is equal to the morphism induced by
$\lefttop{\tau}\Gr^V_b(F):
 \lefttop{\tau}\Gr^V_b(\nbigm_1)
 \lrarr
 \lefttop{\tau}\Gr^V_b(\nbigm_2)$
on the $\lefttop{t}V_a$-parts.
Hence, the kernel, the image and the cokernel of
(\ref{eq;21.3.17.1})
are equal to
$\lefttop{t}V_a\Ker\bigl(
  \lefttop{\tau}\Gr^V_b(F)\bigr)$,
$\lefttop{t}V_a\Image\bigl(
   \lefttop{\tau}\Gr^V_b(F)\bigr)$,
and
$\lefttop{t}V_a\Cok\bigl(
\lefttop{\tau}\Gr^V_b(F)\bigr)$,
respectively.
Then, it is standard to obtain the following equality
in $\lefttop{t,\tau}V_{a,b}(\nbigm_2)$
for any $c<b$:
\begin{equation}
\label{eq;22.7.23.1}
F(\lefttop{t,\tau}V_{a,b}(\nbigm_1))
+\lefttop{t,\tau}V_{a,c}(\nbigm_2)
=
\Bigl(
F(\lefttop{t}V_{a}(\nbigm_1))
\cap
\lefttop{t,\tau}V_{a,b}(\nbigm_2)
\Bigr)
+\lefttop{t,\tau}V_{a,c}(\nbigm_2).
\end{equation}

Fix $b_0<0$.
For any $n\in\seisuu_{\geq 0}$,
we set
\[
 \nbign_n:=
 \bigl\{
 s\in \lefttop{t,\tau}V_{a,b_0}(\nbigm_2)\,\big|\,
 \tau^n s\in F(\lefttop{t}V_a\nbigm_1)
 \bigr\}.
\]
We obtain an increasing sequence of
coherent $\lefttop{t,\tau}\nbigrtilde_Y$-submodules
$\nbign_n$ of $\lefttop{t,\tau}V_{a,b_0}(\nbigm_2)$.
Because $\lefttop{t,\tau}\nbigrtilde_Y$ is Noetherian,
there exists $n(0)$ such that
$\nbign_n=\nbign_{n(0)}$ for any $n\geq n(0)$.
It implies the following for any $p\in\seisuu_{\geq 0}$:
\[
 F(\lefttop{t}V_{a}(\nbigm_1))
 \cap
 \lefttop{t,\tau}V_{a,b_0-n(0)-p}(\nbigm_2)
 =\tau^p\Bigl(
 F(\lefttop{t}V_{a}(\nbigm_1))
 \cap
 \lefttop{t,\tau}V_{a,b_0-n(0)}(\nbigm_2)
 \Bigr).
\]
There exists $d\in\real$ such that
\[
  F(\lefttop{t}V_{a}(\nbigm_1))
 \cap
 \lefttop{t,\tau}V_{a,b_0-n(0)}(\nbigm_2)
 \subset
 F(\lefttop{t,\tau}V_{a,d}(\nbigm_1)).
\]
For any given $b$, we choose $j$
such that $d-j<b$.
We obtain
\begin{multline}
  F(\lefttop{t}V_{a}(\nbigm_1))
 \cap
 \lefttop{t,\tau}V_{a,b_0-n(0)-j}(\nbigm_2)
 \subset
 \tau^j
 \Bigl(
 F(\lefttop{t}V_{a}(\nbigm_1))
 \cap
 \lefttop{t,\tau}V_{a,b_0-n(0)}(\nbigm_2)
 \Bigr)
 \\
 \subset
 \tau^j F(\lefttop{t,\tau}V_{a,d}(\nbigm_1))
\subset F(\lefttop{t,\tau}V_{a,b}(\nbigm_1)).
\end{multline}
By considering the intersection of 
$F(\lefttop{t}V_{a}(\nbigm_1))$
with the both hand sides of (\ref{eq;22.7.23.1})
for $b$ and $c=b_0-n(0)-j$,
we obtain
$F(\lefttop{t,\tau}V_{a,b}(\nbigm_1))
=
F(\lefttop{t}V_{a}(\nbigm_1))
\cap
\lefttop{t,\tau}V_{a,b}(\nbigm_2)$.
Thus, we obtain that
$F:\lefttop{t}V_a(\nbigm_1)\lrarr
\lefttop{t}V_a(\nbigm_2)$ is strict
with respect to the filtrations
$\lefttop{t,\tau}V_{a,\bullet}(\nbigm_i)$,
and hence
\[
 F\Bigl(
 \lefttop{t,\tau}V_{a,b}(\nbigm_1)
 \Bigr)
 =F\Bigl(
  \lefttop{t}V_{a}(\nbigm_1)
 \Bigr)
 \cap
 \lefttop{t,\tau}V_{a,b}(\nbigm_2)
=F(\nbigm_1)
\cap \lefttop{t}V_a(\nbigm_2)
\cap
\lefttop{t,\tau}V_{a,b}(\nbigm_2)
=F(\nbigm_1)
\cap
\lefttop{t,\tau}V_{a,b}(\nbigm_2).
\]
Thus, we obtain the first claim.

We define the $\lefttop{t,\tau}V\nbigrtilde_{Y}$-modules
$\lefttop{t,\tau}V_{a,b}(\Ker F)$,
$\lefttop{t,\tau}V_{a,b}(\Image F)$
and
$\lefttop{t,\tau}V_{a,b}(\Cok F)$
as the kernel, the image and the cokernel of
$\lefttop{t,\tau}V_{a,b}(\nbigm_1)
\lrarr
\lefttop{t,\tau}V_{a,b}(\nbigm_2)$.
Because the induced morphism
$\lefttop{t}V_{a}(\nbigm_1)
\lrarr
 \lefttop{t}V_{a}(\nbigm_2)$
is strict with respect to the filtrations
$\lefttop{t,\tau}V_{a,\bullet}$,
the kernel, the image and the cokernel of
\[
 \lefttop{t,\tau}V_{a,b}(\nbigm_1)\big/
 \lefttop{t,\tau}V_{a,<b}(\nbigm_1)
\lrarr
 \lefttop{t,\tau}V_{a,b}(\nbigm_2)\big/
 \lefttop{t,\tau}V_{a,<b}(\nbigm_2)
\]
are isomorphic to
\[
\lefttop{t,\tau}V_{a,b}(\Ker(F))\big/
 \lefttop{t,\tau}V_{a,<b}(\Ker(F)),
\quad
 \lefttop{t,\tau}V_{a,b}(\Image(F))\big/
 \lefttop{t,\tau}V_{a,<b}(\Image(F)),
\quad
\lefttop{t,\tau}V_{a,b}(\Cok(F))\big/
 \lefttop{t,\tau}V_{a,<b}(\Cok(F)), 
\]
respectively.
Hence, we obtain
$\lefttop{t,\tau}V_{a,b}(\nbign)\big/
 \lefttop{t,\tau}V_{a,<b}(\nbign)
 \simeq
 \lefttop{t}V_a\lefttop{\tau}\Gr^V_b(\nbign)$
for 
$\nbign=\Ker(F),\,\Image(F),\,\Cok(F)$.
Similarly, we obtain
$\lefttop{t,\tau}V_{a,b}(\nbign)\big/
 \lefttop{t,\tau}V_{<a,b}(\nbign)
 \simeq
 \lefttop{\tau}V_b\lefttop{t}\Gr^V_a(\nbign)$
for 
$\nbign=\Ker(F),\,\Image(F),\,\Cok(F)$.
We can easily check the other conditions for
$\lefttop{t,\tau}V(\nbign)$ $(\nbign=\Ker(F),\Image(F),\Cok(F))$
to be $V$-filtrations along $(t,\tau)$.
\hfill\qed

\vspace{.1in}
We give some complements.
The following lemma is easy to see.
\begin{lem}
Suppose that the support of $\nbigm\in\nbigc(Y)$
is contained in $\{t=0\}$.
Then, $\nbigm$ is strictly $(t,\tau)$-specializable.
\hfill\qed
\end{lem}

The following lemma is similar to Lemma \ref{lem;22.7.9.1}.
\begin{lem}
Let $\nbigm\in\nbigc(Y)$.
If $\nbigm=\nbigm[\ast t]$,
then $t\cdot\bigl(\lefttop{t,\tau}V_{0,b}\nbigm\bigr)
=\lefttop{t,\tau}V_{-1,b}(\nbigm)$ holds
for any $b\in\real$.
If $\nbigm=\nbigm[!t]$,
then
$\deldel_t:
 \lefttop{t,\tau}V_{-1,b}(\nbigm)/
  \lefttop{t,\tau}V_{<-1,b}(\nbigm)
 \lrarr
  \lefttop{t,\tau}V_{0,b}(\nbigm)/
  \lefttop{t,\tau}V_{0,b}(\nbigm)$
are isomorphisms for any $b\in\real$.
\hfill\qed
\end{lem}

\subsubsection{Direct image by a projective morphism}
\label{subsection;22.7.9.10}

Let $f_X:X_1\lrarr X_2$ be a projective morphism of complex manifolds.
Let $U$ be an open subset in $\cnum_t\times \cnum_{\tau}$.
We set $Y_i:=U\times X_i$,
and let $f:Y_1\lrarr Y_2$ denote the induced morphism.

\begin{prop}
\label{prop;21.3.17.2}
Suppose that $\nbigm\in\nbigc(Y_1)$ is strictly
$(t,\tau)$-specializable.
Then,
$f_{\dagger}^j(\nbigm)$ are also
strictly $(t,\tau)$-specializable.
The $V$-filtrations of $f_{\dagger}^j(\nbigm)$
along $(t,\tau)$ are equal to the filtrations
induced by $\lefttop{t,\tau}V(\nbigm)$,
i.e.,
\begin{equation}
\label{eq;22.7.9.11}
 \lefttop{t,\tau}V_{a,b}(f_{\dagger}^j\nbigm)
=f_{\dagger}^j\bigl(
\lefttop{t,\tau}V_{a,b}(\nbigm)
 \bigr)
 :=R^jf_{\ast}
 \Bigl(
 \pi^{\ast}\nbigr_{X_2\larr X_1}
 \otimes^L_{\pi^{\ast}\nbigr_{X_1}}
 \lefttop{t,\tau}V_{a,b}(\nbigm)
 \Bigr).
\end{equation}
\end{prop}
\pf
Recall that
$\lefttop{t}V_af_{\dagger}^j(\nbigm)
=
f_{\dagger}^j(\lefttop{t}V_a\nbigm)$
and
$\lefttop{\tau}V_bf_{\dagger}^j(\nbigm)
=
f_{\dagger}^j(\lefttop{\tau}V_b\nbigm)$.
(See \cite{Sabbah-pure-twistor}.
See also \cite{Mebkhout-Sabbah}.)
Because Proposition \ref{prop;21.3.17.2} can be proved
similarly,
we explain only an outline of the proof.
Note that we also have
$\lefttop{t}V_af_{\dagger}^j(\lefttop{\tau}\Gr^V_b(\nbigm))
=f_{\dagger}^j(\lefttop{t}V_a\lefttop{\tau}\Gr^V_b(\nbigm))$
and 
$\lefttop{\tau}V_bf_{\dagger}^j(\lefttop{t}\Gr^V_a(\nbigm))
=f_{\dagger}^j(\lefttop{\tau}V_b\lefttop{t}\Gr^V_a(\nbigm))$.
In particular,
$f_{\dagger}^j(\lefttop{t}V_a\lefttop{\tau}\Gr^V_b(\nbigm))$
and 
$f_{\dagger}^j(\lefttop{\tau}V_b\lefttop{t}\Gr^V_a(\nbigm))$
are strict.
By an easy induction,
we can check that
$f_{\dagger}^j\Bigl(
 \lefttop{t,\tau}V_{a,b}(\nbigm)\big/
 \lefttop{t,\tau}V_{a,c}(\nbigm)
 \Bigr)$
 and
$f_{\dagger}^j\Bigl(
 \lefttop{t,\tau}V_{a,b}(\nbigm)\big/
 \lefttop{t,\tau}V_{d,b}(\nbigm)
 \Bigr)$
are strict for any $a>d$  and $b>c$.

We set
$\lefttop{\tau}V_b
 \Bigl(
 f_{\dagger}^j(\lefttop{t}V_a \nbigm)
 \Bigr)$
as the image of
the naturally induced morphism
$f_{\dagger}^j\bigl(
 \lefttop{t,\tau}V_{a,b}(\nbigm)
 \bigr)
 \lrarr
 f_{\dagger}^j( \lefttop{t}V_a\nbigm)$.
It is standard to prove that
$f_{\dagger}^j\bigl(
 \lefttop{t,\tau}V_{a,b}(\nbigm)
 \bigr)$
are $\lefttop{t,\tau}V\nbigr_{Y_2}$-coherent
(see \cite[\S4.7]{kashiwara_text}),
and that
$\lefttop{\tau}V_b
 \Bigl(
 f_{\dagger}^j(\lefttop{t}V_a\nbigm)
 \Bigr)$ 
 are $\lefttop{t,\tau}V\nbigr_{Y_2}$-coherent
(see \cite[Appendix]{kashiwara_text}).

Recall $\nbigy_2^0=\{0\}\times Y_2\subset\nbigy_2$.
It is also standard to see that
there exists a splitting
\[
 f_{\dagger}^j\Bigl(
 \lefttop{t}V_{a}(\nbigm)\big/
 \lefttop{t,\tau}V_{a,b}(\nbigm)
 \Bigr)_{|\nbigy_2\setminus\nbigy_2^0}
 \simeq
 \bigoplus_{b<c}
 f_{\dagger}^j\bigl(
 \lefttop{t}V_a\lefttop{\tau}\Gr^V_c(\nbigm)
 \bigr)_{|\nbigy_2\setminus\nbigy_2^0},
\]
where the actions of $-\deldel_tt-\lambda c$
on
$f_{\dagger}^j\bigl(
 \lefttop{t}V_a\lefttop{\tau}\Gr^V_c(\nbigm)
 \bigr)_{|\nbigy_2\setminus\nbigy_2^0}$
are nilpotent.

We set
$n:=\dim X-\dim Y$.
There exists the following epimorphism for any $b$:
\[
   f_{\dagger}^n\bigl(
 \lefttop{t}V_{a}(\nbigm)
 \bigr)
 \lrarr
 \Cok\Bigl(
 f_{\dagger}^n\bigl(
 \lefttop{t,\tau}V_{a,b}(\nbigm)
 \bigr)
\lrarr
    f_{\dagger}^n\bigl(
 \lefttop{t}V_{a}(\nbigm)
 \bigr)
\Bigr) \\
\simeq
f_{\dagger}^n\bigl(
 \lefttop{t}V_{a}(\nbigm)/
 \lefttop{t,\tau}V_{a,b}(\nbigm)
\bigr).
\]
Let $N$ denote the image of
the following morphism
of $\lefttop{t,\tau}V\nbigrtilde_{Y_2}$-modules:
\begin{equation}
\label{eq;21.1.30.2}
 f_{\dagger}^{n-1}\Bigl(
 \lefttop{t}V_{a}(\nbigm)
 \big/
 \lefttop{t,\tau}V_{a,b}(\nbigm)
 \Bigr)
 \lrarr
 f_{\dagger}^{n}\Bigl(
 \lefttop{t,\tau}V_{a,b}(\nbigm)
 \Bigr).
\end{equation}

\begin{lem}
\label{lem;21.3.16.20}
 $N=0$.
\end{lem}
\pf
For any $c\leq b$,
let $\lefttop{\tau}V_{c}f_{\dagger}^n(\lefttop{t,\tau}V_{a,b}(\nbigm))$
denote the image of
$f_{\dagger}^n(\lefttop{t,\tau}V_{a,c}(\nbigm))
\lrarr
f_{\dagger}^n(\lefttop{t,\tau}V_{a,b}(\nbigm))$.
We obtain the induced morphism
\begin{equation}
\label{eq;21.3.17.3}
 f_{\dagger}^{n-1}\Bigl(
 \lefttop{t}V_{a}(\nbigm)
 \big/
 \lefttop{t,\tau}V_{a,b}(\nbigm)
 \Bigr)
 \lrarr
 f_{\dagger}^{n}\Bigl(
 \lefttop{t,\tau}V_{a,b}(\nbigm)
 \Bigr)\Big/
\lefttop{\tau}V_{c}
 f_{\dagger}^{n}\Bigl(
 \lefttop{t,\tau}V_{a,b}(\nbigm)
 \Bigr)
 \simeq
 f_{\dagger}^n\Bigl(
 \lefttop{t,\tau}V_{a,b}(\nbigm)\big/
 \lefttop{t,\tau}V_{a,c}(\nbigm)
 \Bigr).
\end{equation}
As the restriction of (\ref{eq;21.3.17.3})
to $\nbigy_2\setminus\nbigy_2^0$,
we obtain
\begin{equation}
\label{eq;21.3.17.4}
 \bigoplus_{b'>b}
  f_{\dagger}^{n-1}\bigl(
   \lefttop{t}V_a\lefttop{\tau}\Gr^V_{b'}(\nbigm)
   \bigr)_{|\nbigy_2\setminus\nbigy_2^0}
   \lrarr
    \bigoplus_{c<b'\leq b}
  f_{\dagger}^n\bigl(
   \lefttop{t}V_a\lefttop{\tau}\Gr^V_{b'}(\nbigm)
   \bigr)_{|\nbigy_2\setminus\nbigy_2^0},
\end{equation}
and the actions of $-\deldel_tt-\lambda b'$
on 
$f_{\dagger}^j\bigl(
\lefttop{t}V_a\lefttop{\tau}\Gr^V_{b'}(\nbigm)
\bigr)_{|\nbigy_2\setminus\nbigy_2^0}$
are nilpotent.
Hence, we obtain that (\ref{eq;21.3.17.4}) is $0$.
Because
$f_{\dagger}^n\Bigl(
 \lefttop{t,\tau}V_{a,b}(\nbigm)\big/
 \lefttop{t,\tau}V_{a,c}(\nbigm)
 \Bigr)$
 is strict,
we obtain that (\ref{eq;21.3.17.3}) is $0$.
Therefore, 
$N$ is contained in
$\lefttop{\tau}V_{c}
 f_{\dagger}^{n}\Bigl(
 \lefttop{t,\tau}V_{a,b}(\nbigm)
 \Bigr)$ for any $c\in\real$.

Suppose that $N\neq 0$,
and we shall deduce a contradiction.
Note that $N$ is $\lefttop{t,\tau}V\nbigrtilde_{Y_2}$-coherent,
and hence it is a good $\nbigo_{\nbigy_2}$-module.
Let $N_0\subset N$ denote the kernel of $\tau$.
Because the support of $N$ is contained in
$\{\tau=0\}$,
we obtain $N_0\neq 0$.
Take $0\neq v\in N_0$.
For any $j\in\seisuu_{>0}$,
$N_0$ is contained in
$\lefttop{\tau}V_{-j-10}
 f_{\dagger}^n\Bigl(
 \lefttop{t,\tau}V_{a,b}(\nbigm)
 \Bigr)$,
 and hence 
there exists
$v^{(j)}\in
\lefttop{\tau}V_{<0}f_{\dagger}^n\Bigl(
 \lefttop{t,\tau}V_{a,b}(\nbigm)\Bigr)$
 such that
$\tau^jv^{(j)}=v$.
We consider the $\lefttop{t,\tau}V\nbigr_{Y_2}$-submodule
$N_1\subset N$ generated by
$v^{(j)}$ $(j=1,2,\ldots)$.
Because $N$ is $\lefttop{t,\tau}V\nbigr_{Y_2}$-coherent,
$N_1$ is also $\lefttop{t,\tau}V\nbigr_{Y_2}$-coherent.
\begin{lem}
For any $s\in N_1$,
there exists $L(s)\in\seisuu_{>0}$
such that
$\tau^{L(s)}s=0$. 
\end{lem}
\pf
We have an expression as a finite sum
$s=\sum_{j\leq m} P_jv^{(j)}$,
where
$m$ is a positive number,
and
$P_j\in \lefttop{t,\tau}V\nbigr_{Y_2}$.
We obtain $\tau^{m+1}s=0$.
\hfill\qed

\vspace{.1in}

Because $N_1$ is
$\lefttop{t,\tau}V\nbigr_{Y_2}$-coherent,
there exist
$u_1,\ldots,u_m$
which generates $N_1$ over
$\lefttop{t,\tau}V\nbigr_{Y_2}$.
Choose $L>\max\{L(u_i)\}$.
Then, we obtain that
$\tau^{L}s=0$ for any $s\in N_1$.
However, if $j>L$,
then $\tau^Lv^{(j)}\neq 0$.
Thus, we obtain a contradiction,
i.e.,
Lemma \ref{lem;21.3.16.20} is proved.
\hfill\qed

\vspace{.1in}

By Lemma \ref{lem;21.3.16.20},
for any $b$,
$f_{\dagger}^n\bigl(
 \lefttop{t,\tau}V_{a,b}\nbigm
 \bigr)
 \lrarr
  f_{\dagger}^n\bigl(
 \lefttop{t}V_{a}\nbigm
 \bigr)$
 is a monomorphism,
 and 
\[
 f_{\dagger}^{n-1}\bigl(
 \lefttop{t,\tau}V_{a,b}\nbigm
 \bigr)
\lrarr
 f_{\dagger}^{n-1}\bigl(
 \lefttop{t}V_{a}\nbigm
 \bigr)
 \lrarr
 f_{\dagger}^{n-1}\Bigl(
 \lefttop{t}V_{a}\nbigm\big/
 \lefttop{t,\tau}V_{a,b}\nbigm
 \Bigr)\lrarr 0
\]
is exact.
By using a similar argument with a descending induction,
we can prove that
\[
0\lrarr f_{\dagger}^j\bigl(
 \lefttop{t,\tau}V_{a,b}\nbigm
 \bigr)
 \lrarr
  f_{\dagger}^j\bigl(
 \lefttop{t}V_{a}\nbigm
 \bigr)
 \lrarr
   f_{\dagger}^j\Bigl(
   \lefttop{t}V_{a}\nbigm\big/
    \lefttop{t,\tau}V_{a,b}\nbigm
 \Bigr)
 \lrarr 0
\]
are exact for any $j$ and $b$.
In particular,
we obtain that
$f_{\dagger}^j\bigl(
 \lefttop{t,\tau}V_{a,b}\nbigm
 \bigr)
 \lrarr
 f_{\dagger}^j(\nbigm)$
 are monomorphisms.
We also obtain
\[
 f_{\dagger}^j\bigl(
 \lefttop{t,\tau}V_{a,b}\nbigm
 \bigr)
 \big/
 f_{\dagger}^j\bigl(
 \lefttop{t,\tau}V_{<a,b}\nbigm
 \bigr)
 \simeq
 \lefttop{\tau}V_b
 f_{\dagger}^j\bigl(
 \lefttop{t}\Gr^V_a(\nbigm)
 \bigr),
\]
\[
  f_{\dagger}^j\bigl(
 \lefttop{t,\tau}V_{a,b}\nbigm
 \bigr)
 \big/
 f_{\dagger}^j\bigl(
 \lefttop{t,\tau}V_{a,<b}\nbigm
 \bigr)
 \simeq
 \lefttop{t}V_a
 f_{\dagger}^j\bigl(
 \lefttop{\tau}\Gr^V_b(\nbigm)
 \bigr).
\]
Then, it is easy to check that
$f^j_{\dagger}\bigl(
\lefttop{t,\tau}V_{\bullet,\bullet}(\nbigm)
\bigr)$
are $V$-filtrations of
$f^j_{\dagger}\nbigm$
along $(t,\tau)$,
and we obtain Proposition \ref{prop;21.3.17.2}.
\hfill\qed

\subsubsection{Variant}
We use the notation in \S\ref{subsection;22.7.9.3}.

\begin{df}
A $V$-filtration of $\nbigm\in\nbigc(Y)(\ast\tau)$ along $(t,\tau)$
is a tuple
$\bigl\{
\lefttop{t,\tau}V_{a,b}(\nbigm)\,\big|\,
a,b\in\real
\bigr\}$ of
$\lefttop{t,\tau}V\nbigrtilde_{Y}$-submodules
satisfying the conditions in 
Definition {\rm\ref{df;21.3.18.10}}.
We say that $\nbigm$ is strictly $(t,\tau)$-specializable
if such a $V$-filtration exists. 
\hfill\qed
\end{df}

\begin{rem}
As in Lemma {\rm\ref{lem;21.3.18.20}},
if the conditions {\rm 1--7} are satisfied,
it is enough to impose either one of $8$ or $9$.
\hfill\qed
\end{rem}

The following lemma is similar to
Lemma \ref{lem;22.7.9.1} and Lemma \ref{lem;22.7.9.2}.
\begin{lem}
Let $\lefttop{t,\tau}V_{\bullet,\bullet}$
be a $V$-filtration of $\nbigm\in\nbigc(Y)(\ast\tau)$
along $(t,\tau)$.
\begin{itemize}
\item For $a<0$ and $b\in\real$, we have
      $t\cdot\left(
      \lefttop{t,\tau}V_{a,b}(\nbigm)\right)
=\lefttop{t,\tau}V_{a-1,b}(\nbigm)$.
For $a>{-1}$ and $b\in\real$,
the induced morphism
$\deldel_t:
 \lefttop{t,\tau}V_{a,b}/\lefttop{t,\tau}V_{<a,b}
 \lrarr
 \lefttop{t,\tau}V_{a+1,b}/\lefttop{t,\tau}V_{<a+1,b}$
is an isomorphism of
 $V\nbigrtilde_{\{t=0\}\cap Y}(\ast\tau)$-modules.
 \item
      For any $a\in\real$ and $b\in\real$, we have
$\tau\cdot \lefttop{t,\tau}V_{a,b}
=\lefttop{t,\tau}V_{a,b-1}$.
 \item $\lefttop{t,\tau}V_{a,b}=\lefttop{t}V_a\cap\lefttop{\tau}V_{b}$ holds.
       In particular,
       $\lefttop{t,\tau}V_{\bullet,\bullet}$
       is uniquely determined for
       a strictly $(t,\tau)$ specializable $\nbigm$.
 \item If $\nbigm=\nbigm[\ast t]$,
       we have
      $t\cdot\left(
      \lefttop{t,\tau}V_{0,b}(\nbigm)\right)
       =\lefttop{t,\tau}V_{-1,b}(\nbigm)$
       for any $b\in\real$.
       If $\nbigm=\nbigm[!t]$,
the induced morphism
$\deldel_t:
 \lefttop{t,\tau}V_{-1,b}/\lefttop{t,\tau}V_{<-1,b}
 \lrarr
 \lefttop{t,\tau}V_{0,b}/\lefttop{t,\tau}V_{<0,b}$
is an isomorphism of
       $V\nbigrtilde_{\{t=0\}\cap Y}(\ast\tau)$-modules
       for any $b\in\real$.       
\hfill\qed
\end{itemize}
\end{lem}

The following proposition is similar to Proposition \ref{prop;21.3.17.10}.
\begin{prop}
Let $F:\nbigm_1\lrarr\nbigm_2$ be a morphism
in $\nbigc(Y)(\ast\tau)$
such that
$\Ker(F)$,
$\Image(F)$ and $\Cok(F)$
are also objects of $\nbigc(Y)(\ast\tau)$.
Assume that $\nbigm_i$ are strictly $(t,\tau)$-specializable.
Then, $F$ is strict with respect to
$\lefttop{t,\tau}V_{\bullet,\bullet}$.
As a result,
$\Ker(F)$, $\Image(F)$ and $\Cok(F)$
are also strictly $(t,\tau)$-specializable,
and the $V$-filtrations
of $\Ker(F)$, $\Image(F)$ and $\Cok(F)$ along $(t,\tau)$
are equal to the filtrations
naturally induced by $\lefttop{t,\tau}V(\nbigm_i)$.
\hfill\qed
\end{prop}

\begin{lem}
If $\nbigm\in\nbigc(Y)$ is strictly $(t,\tau)$-specializable,
$\nbigm(\ast\tau)\in\nbigc(Y)(\ast\tau)$ 
is also strictly $(t,\tau)$-specializable.
\end{lem}
\pf
For $b\in\real$ and $a\in\real$,
we choose $n\in\seisuu_{\geq 0}$ such that $b-n<0$,
and we set
\[
\lefttop{t,\tau}V_{a,b}(\nbigm(\ast\tau))
=\tau^{-n}\Bigl(
\lefttop{t,\tau}V_{a,b-n}(\nbigm)
\Bigr)
\subset \nbigm(\ast\tau).
\]
It is independent of the choice of $n$.
It is easy to check the conditions {\rm 1--7} and {\rm 9}
in Definition \ref{df;21.3.18.10}.
\hfill\qed

\vspace{.1in}
The following proposition is similar to
Proposition \ref{prop;21.3.17.2}.

\begin{prop}
Let $f:Y_1\to Y_2$ be a morphism as in 
{\rm\S\ref{subsection;22.7.9.10}}.
If $\nbigm\in \nbigc(Y_1)(\ast\tau)$
is strictly $(t,\tau)$-specializable,
$f_{\dagger}^j(\nbigm)$ are also
strictly $(t,\tau)$-specializable,
and the $V$-filtration along $(t,\tau)$ of 
 $f_{\dagger}^j(\nbigm)$
are given as in {\rm(\ref{eq;22.7.9.11})}.
\hfill\qed 
\end{prop}

\subsubsection{$(g,\tau)$-specializability}
\label{subsection;21.3.17.11}

Let $Z$ be an open subset in $\cnum_{\tau}\times X$.
Let $g$ be any holomorphic function on $Z$.

\begin{df}
$\nbigm\in\nbigc(Z)$ (resp. $\nbigc(Z)(\ast\tau)$)
is called strictly $(g,\tau)$-specializable
if $\iota_{g\dagger}(\nbigm)\in\nbigc(\cnum_t\times Z)$
(resp. $\nbigc(\cnum_t\times Z)(\ast\tau)$)
is strictly $(t,\tau)$-specializable.
Here, $\iota_g:Z\to\cnum_t\times Z$ denotes the graph embedding.
\hfill\qed
\end{df}

We obtain the following proposition from Proposition \ref{prop;21.3.17.10}.
\begin{prop}
\label{prop;21.3.19.2}
 Let $F:\nbigm_1\lrarr\nbigm_2$ be a morphism
 in $\nbigc(Z)$ (resp. $\nbigc(Z)(\ast\tau)$)
such that $\Ker(F)$, $\Image(F)$ and $\Cok(F)$
 are also contained in $\nbigc(Z)$
 (resp. $\nbigc(Z)(\ast\tau)$).
If $\nbigm_i$ are strictly $(g,\tau)$-specializable,
then $\Ker(F)$, $\Image(F)$ and $\Cok(F)$
are also $(g,\tau)$-specializable.
\hfill\qed
\end{prop}

Let $f_X:X_1\lrarr X_2$ be a projective morphism of complex manifolds.
Let $U$ be an open subset in $\cnum_{\tau}$.
We set $Z_i=U\times X_i$.
Let $f:Z_1\lrarr Z_2$ denote the induced morphism.
Let $g$ be a holomorphic function on $Z_2$.
We obtain the following proposition from Proposition \ref{prop;21.3.17.2}.
\begin{prop}
\label{prop;21.3.19.3}
If $\nbigm\in\nbigc(Z_1)$
(resp. $\nbigc(Z_1)(\ast\tau)$)
is $(f^{\ast}(g),\tau)$-specializable,
then $f_{\dagger}^j\nbigm\in\nbigc(Z_2)$
(resp. $\nbigc(Z_2)(\ast\tau)$)
are strictly $(g,\tau)$-specializable.
\hfill\qed
\end{prop}

\subsubsection{Strictly $(g,\tau)$-specializable
integrable mixed twistor $\nbigd$-modules}

Let $Z$ be as in \S\ref{subsection;21.3.17.11}.

\begin{df}
An integrable mixed twistor $\nbigd$-module $\nbigt$
is called $(g,\tau)$-specializable
if the underlying $\nbigrtilde_Z$-modules
are $(g,\tau)$-specializable.
\hfill\qed
\end{df}

\begin{cor}
The full subcategory of $(g,\tau)$-specializable integrable
mixed twistor $\nbigd_Z$-modules
is abelian.
\hfill\qed
\end{cor}

\begin{prop}
Let $f:Z_1\lrarr Z_2$ be as in {\rm\S\ref{subsection;21.3.17.11}}.
If $\nbigt\in\MTM^{\integral}(Z_1)$
is strictly $(g,\tau)$-specializable,
$f^j_{\dagger}(\nbigt)$
are also strictly $(g,\tau)$-specializable.
\hfill\qed
\end{prop}

\subsubsection{Appendix: Specializability of $\nbigd$-modules in two directions}

Let $X$ be a complex manifold.
Let $Y$ be an open subset in $\cnum_t\times\cnum_{\tau}\times X$.
Let $\pi:Y\lrarr X$ denote the projection.
We set
$\lefttop{t,\tau}V\nbigd_Y:=
\pi^{\ast}(\nbigd_X)\langle t\del_t,\tau\del_{\tau}\rangle$.
We choose a total order $\leq_{\cnum}$ on $\cnum$
such that
(i) the restriction of $\leq_{\cnum}$ to $\seisuu$
is equal to the standard order,
(ii) $a_1\leq_{\cnum}a_2$
implies that
$a_1+n\leq_{\cnum}a_2+n$
for any $n\in\seisuu$.

Let $M$ be a holonomic $\nbigd_{Y}$-module.
Let $\lefttop{t}V$ and $\lefttop{\tau}V$ denote
the $V$-filtrations along $t$ and $\tau$, respectively,
with respect to the total order $\leq_{\cnum}$.
\begin{df}
A $V$-filtration of $M$
along $(t,\tau)$ is a tuple of
 $\lefttop{t,\tau}V\nbigd_Y$-coherent submodules
$\bigl\{\lefttop{t,\tau}V_{a,b}(M)
 \,\big|\,a,b\in\real \bigr\}$
of $M$ such that the following holds.
\begin{itemize}
 \item $\lefttop{t,\tau}V_{a,b}(M)\subset
       \lefttop{t,\tau}V_{a',b'}(M)$
       if $a\leq_{\cnum} a'$ and $b\leq_{\cnum} b'$.
 \item For any $(a,b)\in\real^2$ and for any compact subset $K\subset Y$,
      there exists $\epsilon>0$ such that
       $V_{a,b}(M)_{|K}=V_{a+\epsilon,b+\epsilon}(M)_{|K}$.       
 \item $\lefttop{t,\tau}V_{a,b}(M)(\ast(t\tau))
       =M(\ast (t\tau))$.
 \item $\bigcup_{a}\lefttop{t,\tau}V_{a,b}(M)=\lefttop{\tau}V_b(M)$
       and
       $\bigcup_{b}\lefttop{t,\tau}V_{a,b}(M)=\lefttop{t}V_a(M)$.
 \item $t\cdot \lefttop{t,\tau}V_{a,b}\subset \lefttop{t,\tau}V_{a-1,b}$
       and
       $\del_t\cdot \lefttop{t,\tau}V_{a,b}\subset \lefttop{t,\tau}V_{a+1,b}$.
 \item $\tau\cdot \lefttop{t,\tau}V_{a,b}\subset \lefttop{t,\tau}V_{a,b-1}$
       and
       $\del_{\tau}\cdot \lefttop{t,\tau}V_{a,b}
       \subset \lefttop{t,\tau}V_{a,b+1}$.
 \item The natural morphism
       $\lefttop{t,\tau}V_{a,b}\big/
       \lefttop{t,\tau}V_{<a,b}
       \lrarr
       \lefttop{t}\Gr^V_a(\nbigm)$
       induces an isomorphism
       \[
       \lefttop{t,\tau}V_{a,b}\big/
       \lefttop{t,\tau}V_{<a,b}
       \simeq
       \lefttop{\tau}V_b\bigl(
       \lefttop{t}\Gr^V_a(\nbigm)
       \bigr).
       \]
       Similarly,
       the natural morphism
       $\lefttop{t,\tau}V_{a,b}\big/
       \lefttop{t,\tau}V_{a,<b}
       \lrarr
       \lefttop{\tau}\Gr^V_b(\nbigm)$
       induces an isomorphism
       \[
       \lefttop{t,\tau}V_{a,b}\big/
       \lefttop{t,\tau}V_{a,<b}
       \simeq
       \lefttop{t}V_a\bigl(
       \lefttop{\tau}\Gr^V_b(\nbigm)
       \bigr).
       \]
\end{itemize}
If $M$ has a $V$-filtration along $(t,\tau)$,
we say that $M$ is $(t,\tau)$-specializable.
\hfill\qed
\end{df}

The following lemma is obvious.
\begin{lem}
Suppose that $M$ is $(t,\tau)$-specializable.
Then, $-\del_tt-a$ on
$\lefttop{t,\tau}V_{a,b}/\lefttop{t,\tau}V_{<a,b}$
is locally nilpotent,
and $-\del_{\tau}\tau-b$ on
$\lefttop{t,\tau}V_{a,b}/\lefttop{t,\tau}V_{a,<b}$
is locally nilpotent.
\hfill\qed 
\end{lem}

\begin{lem}
Let $\lefttop{t,\tau}V_{\bullet,\bullet}$
be $V$-filtration of $M$ along $(t,\tau)$.
For $a<_{\cnum}0$, we have
 $t\cdot \lefttop{t,\tau}V_{a,b}(M)=V_{a-1,b}(M)$.
For $a>_{\cnum}-1$,
the induced morphism
$\deldel_t:
 \lefttop{t,\tau}V_{a,b}/\lefttop{t,\tau}V_{<a,b}
 \lrarr
 \lefttop{t,\tau}V_{a+1,b}/\lefttop{t,\tau}V_{<a+1,b}$
is an isomorphism of
 $V\nbigd_{\{t=0\}\cap Y}$-modules.
 Similar clams also hold for $\tau$.
 \hfill\qed
\end{lem}

\begin{lem}
 If $\lefttop{t,\tau}V_{\bullet,\bullet}$
 is a $V$-filtration of $M$ along $(t,\tau)$,
 the following holds.
\begin{itemize}
 \item If $a\geq a'$, we obtain
       $\lefttop{t,\tau}V_{a,b}\cap
       \lefttop{t}V_{a'}=\lefttop{t,\tau}V_{a',b}$.
       Similarly, if $b\geq b'$, we obtain
      $\lefttop{t,\tau}V_{a,b}\cap
       \lefttop{\tau}V_{b'}=\lefttop{t,\tau}V_{a,b'}$.
 \item
 $\lefttop{t,\tau}V_{a,b}
 =\lefttop{t}V_a\cap\lefttop{\tau}V_b$.       
\end{itemize}
 As a result,
 if $M$ is $(t,\tau)$-specializable,
 a $V$-filtration of $M$ along $(t,\tau)$
depends only on the total order $\leq_{\cnum}$.
\hfill\qed
\end{lem}

The following proposition is similar to 
Proposition \ref{prop;21.3.17.10}.
\begin{prop}
Suppose that $M_i$ $(i=1,2)$ are $(t,\tau)$-specializable.
Let $F:M_1\lrarr M_2$ be a morphism of holonomic $\nbigd_Y$-modules.
\begin{itemize}
 \item $F$ is strictly compatible with $\lefttop{t,\tau}V(M_i)$,
       i.e.,
       $F(\lefttop{t,\tau}V_{a,b}(M_1))
       =F(M_1)\cap\lefttop{t,\tau}V_{a,b}(M_2)$.
 \item $\Ker(F)$, $\Image(F)$ and $\Cok(F)$
       are also specializable with respect to $(t,\tau)$.
       The $V$-filtration with respect to $(t,\tau)$
       on $\Ker(F)$, $\Image(F)$ and $\Cok(F)$
       are equal to the filtrations naturally induced 
       by $\lefttop{t,\tau}V(M_i)$.
\end{itemize}
As a result,
the full subcategory of
$(t,\tau)$-specializable
holonomic $\nbigd_Y$-modules
is abelian. 
\hfill\qed
\end{prop}

The following proposition is similar to 
Proposition \ref{prop;21.3.17.2}.
\begin{prop}
 \label{prop;21.3.16.30}
Let $f_X:X_1\lrarr X_2$
be a proper morphism of complex manifolds.
Let $U$ be an open subset of $\cnum_t\times\cnum_{\tau}$.
We set $Y_i:=U\times X_i$.
Let $f:Y_1\lrarr Y_2$ be the induced morphism.
 If $M$ is a $(t,\tau)$-specializable
 $\nbigd_{Y_1}$-module,
 then
 $f_{\dagger}^j(M)$
 are $(t,\tau)$-specializable
 $\nbigd_{Y_1}$-modules.
 The $V$-filtration of $f_{\dagger}^j(M)$
 along $(t,\tau)$
 is equal to the filtration
 naturally induced by
 $\lefttop{t,\tau}V_{\bullet,\bullet}(M)$,
 i.e.,
\[
 \lefttop{t,\tau}V_{a,b}(f_{\dagger}^jM)
 =f_{\dagger}^j\bigl(
 \lefttop{t,\tau}V_{a,b}(M)
 \bigr)
 :=
 R^jf_{\ast}\Bigl(
 \pi^{\ast}\nbigd_{X_2\larr X_1}\otimes^L_{\pi^{\ast}\nbigd_{X_1}}
 \lefttop{t,\tau}V_{a,b}(M)
 \Bigr).
\]
 Here, $\pi$ denotes the projection $Y_1\lrarr X_1$.
 \hfill\qed
\end{prop}

Let $Z$ be an open subset in $\cnum_{\tau}\times X$.
Let $g$ be any holomorphic function on $Z$.

\begin{df}
A holonomic $\nbigd$-module $M$ on $Z$
is called $(g,\tau)$-specializable
if $\iota_{g\dagger}(M)$
is $(t,\tau)$-specializable.
 Here, $\iota_g:\cnum_{\tau}\times
 X\lrarr \cnum_t\times\cnum_{\tau}\times X$
denotes the graph embedding.
\hfill\qed
\end{df}

\begin{prop}
The category of $(g,\tau)$-specializable
holonomic $\nbigd_{Z}$-modules
is abelian.
\hfill\qed
\end{prop}

 \begin{prop}
  Let $f_X:X_1\lrarr X_2$ be a proper morphism.
  Let $U$ be an open subset in $\cnum_{\tau}$.
We set $Z_i:=U\times X_i$,
and let $f:Z_1\lrarr Z_2$ denote the induced morphism.
Let $g$ be a holomorphic function on $Z_2$.
 If $M$ is a $(f^{\ast}(g),\tau)$-specializable
 holonomic $\nbigd_{Z_1}$-module,
 then $f_{\dagger}^jM$ are
 $(g,\tau)$-specializable holonomic $\nbigd_{Z_2}$-modules.
\hfill\qed
 \end{prop}

\subsection{Strict $(g_0,\tau)$-specializability
of regular rescalable objects}

Let $X$ be any complex manifold.
Let $g$ be a holomorphic function on $X$.
Let $g_0$ denote the induced holomorphic function on
$\lefttop{\tau}X$.
We shall prove the following theorem
in \S\ref{subsection;22.7.4.1}--\S\ref{subsection;22.7.4.2}.

\begin{thm}
\label{thm;21.3.17.20}
Let $\nbigm\in\nbigc_{\res}(X)$.
Suppose that $\Xi_{\DR}(\nbigm)$ is regular.
Then,
$\lefttop{\tau}\Upsilon(\nbigm)_{|\lefttop{\tau}\nbigx}$
is strictly $(g_0,\tau)$-specializable.
 Moreover,
\[
\lefttop{t,\tau}V_{a,b}\bigl(
 \lefttop{\tau}\Upsilon(\iota_{g\dagger}\nbigm)_{|\cnum_{\lambda}\times
 \cnum_t\times\times\cnum_{\tau}\times X}\bigr)
\]
are $\pi^{\ast}\nbigr_X$-coherent,
where $\iota_{g}:X\lrarr \cnum_t\times X$ denotes
the graph embedding,
and $\pi:\cnum_t\times\cnum_{\tau}\times X\lrarr X$
denotes the projection.
\end{thm}

\subsubsection{Regular integrable mixed twistor $\nbigd$-modules}
\label{subsection;22.7.4.1}

Let $X$ be any complex manifold.
An object $\nbigm\in\nbigc(X)$ is called regular
if $\Xi_{\DR}(\nbigm)$ is a regular holonomic $\nbigd_X$-module.
Let $\nbigt=(\nbigm',\nbigm'',C)$ be an integrable
mixed twistor $\nbigd_X$-module
with the weight filtration $W$.
\begin{lem}
If $\nbigm''$ is regular,
then $\nbigm'$ is also regular.
\end{lem}
\pf
There exist filtrations $W$ on $\nbigm'$ and $\nbigm''$
such that
$\Gr^W_w(\nbigt)=(\Gr^W_{-w}\nbigm',\Gr^W_{w}\nbigm'',\Gr^W_wC)$
are polarizable pure twistor $\nbigd$-modules.
In particular,
$\Gr^W_{-w}(\nbigm')\simeq\Gr^W_{w}(\nbigm'')$.
Hence, the regularity of $\nbigm''$
implies the regularity of $\nbigm'$.
\hfill\qed

\vspace{.1in}

Suppose that $\nbigm''$ is regular.
Let $P$ be any point of $X$
with a relatively compact neighbourhood $X_P$.
Let $\Gr^W_w(\nbigt_{|X_P})=\bigoplus_Z\nbigt_{w,Z}$
be the decomposition by the strict support.
For each $(w,Z)$, there exist
a holomorphic function $f$
and
a projective morphism of complex manifolds
$\varphi:\Ztilde\lrarr X_P$,
such that the following holds.
\begin{itemize}
 \item $Z\setminus f^{-1}(0)$ is a complex submanifold of $X_P$.
 \item $\varphi(\Ztilde)=Z$.
 \item $\Htilde=\varphi^{\ast}(f)^{-1}(0)$
       is a normal crossing hypersurface of $\Ztilde$.
 \item $\varphi$ induces
       $\Ztilde\setminus \Htilde\simeq Z\setminus f^{-1}(0)$.
 \item There exists
       a good-KMS $\nbigr_{\Ztilde(\ast\Htilde)}$-triple
       $\nbigttilde_{w,Z}$
       with an isomorphism
       $\varphi_{\dagger}\nbigttilde_{w,Z}
       \simeq
       \nbigt_{w,Z}(\ast f)$.
\end{itemize}
\begin{lem}
\label{lem;21.6.19.10}
 $\nbigttilde_{w,Z}$ is regular-KMS in the sense of
{\rm\cite[\S5.5]{Mochizuki-MTM}}.
\end{lem}
\pf
Let $\nbigm_{w,Z}'$ and $\nbigm_{w,Z}''$
denote the underlying $\nbigr_{\Ztilde(\ast \Htilde)}$-module.
Because $\nbigm''$ is assumed to be regular,
$\Xi_{\DR}(\nbigm_{w,Z}')$
and
$\Xi_{\DR}(\nbigm_{w,Z}'')$
are also regular.
Then, the claim easily follows.
\hfill\qed

\vspace{.1in}
Note that
$\nbigttilde_{w,Z|\Ztilde\setminus\Htilde}$
has a polarization $C_{w,Z}$
induced by a polarization of $\Gr^W_w(\nbigt)$.
The polarized pure twistor $\nbigd$-module
$(\nbigttilde_{w,Z|\Ztilde\setminus\Htilde},C_{w,Z})$
is induced by a harmonic bundle
$(E_{w,Z},\theta_{w,Z},h_{w,Z})$.
By Lemma \ref{lem;21.6.19.10},
the harmonic bundles are tame along $\Htilde$.

\subsubsection{Regular admissible variation of integrable
mixed twistor structure}

Let $X$ be a neighbourhood of $(0,\ldots,0)$
in $\cnum_{\tau}\times\cnum^n$.
We set $H=\{\tau=0\}\cup\bigcup_{i=1}^{\ell}\{z_i=0\}$.
Let $\nbigv$ be an integrable regular-KMS
smooth $\nbigrtilde_{X(\ast H)}$-module
underlying an admissible variation of integrable mixed twistor structure
on $(X,H)$.
Let $\nbigp_{\ast}\nbigv$ be the associated holomorphic family of
regular filtered $\lambda$-flat bundles on
$\cnum_{\lambda}\times (X,H)$
indexed by $\real\times\real^{\ell}$.
We also denote $\tau=z_0$.
We set
$\ellsitabar=\{1,\ldots,\ell\}$
and
$\zeroellsitabar=\{0\}\cup\ellsitabar$.

For any subset $I\subset\zeroellsitabar$,
let $q_I:\real^{\zeroellsitabar}\lrarr\real^I$
denote the projection.
For $\vecc\in\real^I$, we set
\[
 \lefttop{I}\nbigp_{\vecc}\nbigv
 =\bigcup_{\substack{\veca\in\real^{\zeroellsitabar}\\
 q_I(\veca)=\vecc }}
 \nbigp_{\veca}\nbigv
 \subset\nbigv.
\]
Let $\nbigr_{X,I}\subset\nbigr_X$
denote the sheaf of subalgebras generated by
$\deldel_i$ $(i\in I)$ over $\nbigo_{\nbigx}$.
Let $\lefttop{I}V\nbigr_{X,I}\subset\nbigr_{X,I}$
denote the sheaf of subalgebras generated by
$z_i\deldel_{i}$ for $i\in I$.

For $I\sqcup J\subset\zeroellsitabar$,
let $\veca(!I\ast J)\in\real^{I\sqcup J}$
be determined by
$a_i(!I\ast J)=1-\epsilon$ $(i\in I)$
and
$a_j(!I\ast J)=1$ $(j\in J)$.
We set
\[
 \nbigv[!I\ast J]:=
 \nbigr_{X,I\sqcup J}\otimes_{\lefttop{I\sqcup J}V\nbigr_{X,I\sqcup J}}
 \lefttop{I\sqcup J}\nbigp_{\veca(!I,\ast J)}\nbigv.
\]
For $L\subset (\zeroellsitabar)\setminus(I\sqcup J)$,
and for $\vecc\in\real^L$,
we set
\[
 \lefttop{L}V_{\vecc}
 \nbigv[!I\ast J]:=
 \nbigr_{X,I\sqcup J}\otimes_{\lefttop{I\sqcup J}V\nbigr_{X,I\sqcup J}}
 \lefttop{I\sqcup J\sqcup L}\nbigp_{(\veca(!I,\ast J),\vecc)}\nbigv,
\]
where $(\veca(!I,\ast J),\vecc)\in\real^{I\sqcup J\sqcup L}$
denote the element induced by
$\veca(!I,\ast J)$ and $\vecc$.

Let $I\sqcup J\sqcup K=\{1,\ldots,\ell\}$
be a decomposition.
We set
$\nbigv(!I\ast J!K):=\nbigv\bigl(
 !(I\sqcup K)\ast J
 \bigr)$
and
$\nbigv(!I\ast J\ast K):=\nbigv\bigl(
!I\ast(J\sqcup K)
\bigr)$.
For $b\in\real$,
we set
\[
 \lefttop{\tau}V_{b}\nbigv[!I\ast J\star K]:=
 \nbigr_{X,\ellsitabar}
 \otimes_{\lefttop{\ellsitabar}V\nbigr_{X,\ellsitabar}}
 \lefttop{\zeroellsitabar}\nbigp_{(b+1,\veca(!I,\ast J,\star K))}\nbigv.
\]
According to \cite{Mochizuki-MTM},
the natural morphism
$\lefttop{\tau}V_{b}\nbigv[!I\ast J\star K]
\lrarr
 \nbigv[!I\ast J\star K]$
is a monomorphism,
and the tuple
$\lefttop{\tau}V_{\bullet}\nbigv[!I\ast J\star K]$
is the $V$-filtration of $\nbigv[!I\ast J\star K]$
along $\tau$.

We set $g=\vecz^{\vecp}=\prod_{i\in K}z_i^{p_i}$
for $\vecp\in\seisuu_{>0}^{K}$.
We obtain the $\nbigrtilde_{\cnum_t\times X}(\ast\tau)$-modules
\[
 \iota_{g\dagger}\bigl(\nbigv[!I\ast J\star K]\bigr)
 =\iota_{g\ast}(\nbigv[!I\ast J\star K])
 \otimes\cnum[\deldel_t]
\]
on $\cnum_{\lambda}\times\bigl(\cnum_t\times X\bigr)$.
Let us recall the construction of a $V$-filtration
$U_{\bullet}(\iota_{g\dagger}(\nbigv[!I\ast J\star K]))$
along $t$.
For $\vecc\in\real^K$,
we set
\[
 \lefttop{K}V_{\vecc}\nbigv(!I\ast J)
 =\nbigr_{X,I\sqcup J}\otimes_{\lefttop{I\sqcup J}V\nbigr_{X,I\sqcup J}}
 \lefttop{\ellsitabar}\nbigp_{\veca(!I,\ast J),\vecc+\vecdelta_K}\nbigv.
\]
Here, $\vecdelta_K=(1,\ldots,1)\in\real^K$.
For $a\leq 0$,
we set
\[
 U_a\bigl(
 \iota_{g\dagger}(\nbigv[!I\ast J \ast K])
 \bigr)
=\nbigr_{X,K}
\Bigl(
 \lefttop{K}V_{a\vecp}\nbigv[!I\ast J]\otimes 1
\Bigr).
\]
For $a>0$, we set
\[
  U_a\bigl(
 \iota_{g\dagger}(\nbigv[!I\ast J \ast K])
 \bigr)
 =\sum_{\substack{c\leq 0\\ j\in\seisuu_{\geq 0} \\
 c+j\leq a}} \deldel_t^j U_c\bigl(
 \iota_{g\dagger}(\nbigv[!I\ast J\ast K])
 \bigr)
 \subset \iota_{g\dagger}(\nbigv[!I\ast J\ast K]).
\]
For $a<0$,
we set
\[
 U_a\bigl(
 \iota_{g\dagger}(\nbigv[!I\ast J!K])
 \bigr)
=\nbigr_{X,K}
\Bigl(
 \lefttop{K}V_{a\vecp}\nbigv[!I\ast J]\otimes 1
\Bigr).
\]
For $a\geq 0$, we set
\[
  U_a\bigl(
 \iota_{g\dagger}(\nbigv[!I\ast J!K])
 \bigr)
 =\sum_{\substack{c< 0\\ j\in\seisuu \\
 c+j\leq a}} \deldel_t^j U_c\bigl(
 \iota_{g\dagger}(\nbigv[!I\ast J!K])
 \bigr)
 \subset  \iota_{g\dagger}(\nbigv[!I\ast J!K]).
\]

Let us define
$\lefttop{t,\tau}V_{a,b}\bigl(
\iota_{g\dagger}(\nbigv[!I\ast J\star K])
\bigr)$.
For $(a,b)\in\real_{\leq 0}\times\real$,
we set
\begin{equation}
\label{eq;21.3.18.1}
 \lefttop{t,\tau}V_{a,b}(\iota_{g\dagger}\nbigv[!I\ast J\ast K])
=\nbigr_{X,K}
 \Bigl(
 \lefttop{0\sqcup K}V_{(b,a\vecp)}\nbigv[!I\ast J]\otimes 1
 \Bigr)
 \subset
\iota_{g\dagger}\nbigv[!I\ast J\ast K].
\end{equation}
Here, $(b,a\vecp)\in\real^{\zeroellsitabar}$
denote the element induced by
$b\in\real$ and $a\vecp\in\real^{\ellsitabar}$.
For $(a,b)\in\real_{>0}\times\real$,
we set
\begin{equation}
\label{eq;21.3.18.2}
 \lefttop{t,\tau}V_{a,b}(\iota_{g\dagger}\nbigv[!I\ast J\ast K])
 =\sum_{\substack{c\leq 0\\ c+j\leq a}}
 \deldel_t^j\cdot
 \lefttop{t,\tau}V_{c,b}(\iota_{g\dagger}\nbigv[!I\ast J\ast K]). 
\end{equation}
Similarly,
for $(a,b)\in\real_{<0}\times\real$,
we set
\begin{equation}
\label{eq;21.3.18.3}
 \lefttop{t,\tau}V_{a,b}(\iota_{g\dagger}\nbigv[!I\ast J!K])
=\nbigr_{X,K}
 \Bigl(
 \lefttop{0\sqcup K}V_{(b,a\vecp)}\nbigv[!I\ast J]\otimes 1
 \Bigr)
 \subset
\iota_{g\dagger}\nbigv[!I\ast J!K]. 
\end{equation}
For $(a,b)\in\real_{\geq 0}\times\real$,
we set
\begin{equation}
\label{eq;21.3.18.4}
 \lefttop{t,\tau}V_{a,b}(\iota_{g\dagger}\nbigv[!I\ast J!K])
 =\sum_{\substack{c\leq 0\\ c+j\leq a}}
 \deldel_t^j\cdot
 \lefttop{t,\tau}V_{c,b}(\iota_{g\dagger}\nbigv[!I\ast J!K]).
\end{equation}

\begin{prop}
\label{prop;21.3.19.1}
 $\iota_{g\dagger}\nbigv[!I\ast J\star K]$
are strictly $(t,\tau)$-specializable.
 The $V$-filtrations along $(t,\tau)$ are given as in
 {\rm(\ref{eq;21.3.18.1})},
 {\rm(\ref{eq;21.3.18.2})},
 {\rm(\ref{eq;21.3.18.3})}
 and
 {\rm(\ref{eq;21.3.18.4})}.
In particular,
 $\lefttop{t,\tau}V_{a,b}\bigl(
 \iota_{g\dagger}\nbigv[!I\ast J\star K]
 \bigr)$
are $\pi^{\ast}\nbigr_X$-coherent. 
\end{prop}
\pf
By the construction,
we can check that
$\lefttop{t,\tau}V_{a,b}$
are coherent over
$\nbigr_{X,\ellsitabar}$.
By the construction,
it is easy to check that
\[
 \lefttop{t,\tau}V_{a,b}\big/\lefttop{t,\tau}V_{a,<b}
 \simeq
 U_{a}\bigl(
 \iota_{g\dagger}
 \lefttop{\tau}\Gr^{V}_b(\nbigv)[\ast I!J\star K]
 \bigr),
\]
which is the $V$-filtration of
$\lefttop{\tau}\Gr^V_b\bigl(
\iota_{g\dagger}\nbigv[\ast I!J\star K]
\bigr)
=
\iota_{g\dagger}
 \lefttop{\tau}\Gr^{V}_b(\nbigv)[\ast I!J\star K]$.
Then, by using Lemma \ref{lem;21.3.18.20},
we can check that
$\lefttop{t,\tau}V_{\bullet,\bullet}
\bigl(\iota_{g\dagger}(\nbigv[\ast I!J\star K])\bigr)$
is a $V$-filtration along $(t,\tau)$.
Thus, we obtain
Proposition \ref{prop;21.3.19.1}
\hfill\qed

\subsubsection{An easy case}
\label{subsection;21.3.19.5}

Let us prove the claim of Theorem \ref{thm;21.3.17.20}
in the case $\Supp(\nbigm)\subset g^{-1}(0)$.
We obtain $\iota_{g\dagger}(\nbigm)\in\nbigc(\cnum_t\times X)$.
Because the support of
$\iota_{g\dagger}(\nbigm)$ is contained in $\{0\}\times X$,
there exists $\nbigm_1\in\nbigc(X)$
such that
$\iota_{g\dagger}(\nbigm)
=\iota_{0\dagger}(\nbigm_1)$,
where $\iota_0:X\simeq \{0\}\times X\subset \cnum_t\times X$.
We obtain
\begin{equation}
\lefttop{\tau}\Upsilon(\iota_{g\dagger}\nbigm)
=\lefttop{\tau}\Upsilon(\iota_{0\dagger}\nbigm_1)
=\iota_{1\dagger}\bigl(
\lefttop{\tau}\Upsilon(\nbigm_1)
\bigr),
\end{equation}
where
$\iota_{1}:\lefttop{\tau}X
=\cnum_{\tau}\times X
\simeq \cnum_{\tau}\times\{0\}\times X
\subset
\cnum_{\tau}\times \cnum_t\times X$
denote the natural inclusion.

There exists a $V$-filtration
$\lefttop{\tau}V_{\bullet}$
of $\lefttop{\tau}\Upsilon(\nbigm_1)_{|\lefttop{\tau}\nbigx}$.
There exists a natural isomorphism
\[
 \iota_{1\dagger}\bigl(
\lefttop{\tau}\Upsilon(\nbigm_1)
\bigr)
\simeq
\bigoplus_{j=0}^{\infty}
\iota_{1\ast}
\bigl(
\lefttop{\tau}\Upsilon(\nbigm_1)
\bigr)\cdot \deldel_t^j.
\]
We define
$\lefttop{t,\tau}V_{\bullet,\bullet}
\Bigl(
\iota_{1\dagger}\bigl(
\lefttop{\tau}\Upsilon(\nbigm_1)
\bigr)\Bigr)$
as follows:
\[
 \lefttop{t,\tau}V_{a,b}
 \Bigl(
\iota_{1\dagger}\bigl(
\lefttop{\tau}\Upsilon(\nbigm_1)
\bigr)_{|\lefttop{\tau}\nbigx}
\Bigr)
=\bigoplus_{0\leq j\leq a}
 \iota_{1\ast} \lefttop{\tau}V_{b}\bigl(
    \lefttop{\tau}\Upsilon(\nbigm)_{|\lefttop{\tau}\nbigx}
   \bigr)
 \deldel_t^j.
\]
It is easy to see that
$\lefttop{t,\tau}V$ is a $V$-filtration of
$\iota_{1\dagger}\bigl(
\lefttop{\tau}\Upsilon(\nbigm_1)
\bigr)$
along $(t,\tau)$,
and each
$\lefttop{t,\tau}V$ is
$\pi^{\ast}\nbigr_X$-coherent.

\subsubsection{Proof of Theorem \ref{thm;21.3.17.20}}
\label{subsection;22.7.4.2}

We use the induction on the dimension of the support.
It is enough to prove the claim locally around any point of $X$.
Therefore, we may assume the existence of
$f$, $Z$, $\varphi$ and $\nbigv$
as in the proof of Proposition \ref{prop;21.2.18.1}.
By the consideration in \S\ref{subsection;21.3.19.5},
we may assume that $Z\not\subset g^{-1}(0)$.
We may also assume that
$g^{-1}(0)\subset f^{-1}(0)$.
Then,
$\lefttop{\tau}\Upsilon(\nbigm)_{|\lefttop{\tau}\nbigx}$
is described as
the cohomology of the following complex:
\[
 \lefttop{\tau}\Upsilon(\psi^{(1)}_f(\nbigm))_{|\lefttop{\tau}\nbigx}
 \lrarr
 \lefttop{\tau}\Upsilon(\Xi^{(0)}_f(\nbigm))_{|\lefttop{\tau}\nbigx}
 \oplus
 \lefttop{\tau}\Upsilon(\phi^{(0)}_f(\nbigm))_{|\lefttop{\tau}\nbigx}
 \lrarr
 \lefttop{\tau}\Upsilon(\psi^{(0)}_f(\nbigm))_{|\lefttop{\tau}\nbigx}.
\]
By the assumption of the induction,
$\lefttop{\tau}\Upsilon(\psi^{(i)}_f(\nbigm))
_{|\lefttop{\tau}\nbigx}$
and
$\lefttop{\tau}\Upsilon(\phi^{(i)}_f(\nbigm))
_{|\lefttop{\tau}\nbigx}$
are strictly $(g_0,\tau)$-specializable.

Let $f_{0,Z}$ and $g_{0,Z}$ denote the holomorphic functions
on $\lefttop{\tau}Z$ induced by $f_Z$ and $g_Z$,
respectively.
By Proposition \ref{prop;21.3.19.1},
$\Pi^{a,b}_{f_{0,Z},\star}
 \bigl(
 \lefttop{\tau}\Upsilon(\nbigv)
 \bigr)$
are strictly $(g_{0,Z},\tau)$-specializable
for any $a\leq b$.
By Proposition \ref{prop;21.3.19.2},
$\Xi^{(0)}_{f_{0,Z}}
 \bigl(
 \nbigv
 \bigr)$
 are strictly $(g_{0,Z},\tau)$-specializable.
By Proposition \ref{prop;21.3.19.3},
we obtain that
$\Xi^{(0)}_{f_0}(\nbigm)$
is strictly $(g_{0,Z},\tau)$-specializable.
By Proposition \ref{prop;21.3.19.2} again,
$\lefttop{\tau}\Upsilon(\nbigm)_{|\lefttop{\tau}\nbigx}$
is strictly $(g_0,\tau)$-specializable.
We can prove the coherence of $V$-filtrations
similarly,
and the proof of Theorem \ref{thm;21.3.17.20}
is completed.
\hfill\qed

\subsection{Strict specializability of irregular Hodge filtrations
in the regular case}

Let $g$ be a holomorphic function on $X$.
The induced function on $\cnum_{\tau}\times X$ is denoted by $g_0$.

\begin{prop}
\label{prop;21.6.20.1}
Let $\nbigm\in\nbigc_{\res}(X)$.
Suppose the following conditions.
\begin{itemize}
 \item $\lefttop{\tau}\Upsilon(\nbigm)_{|\lefttop{\tau}\nbigx}$
is strictly $(g_0,\tau)$-specializable.
 \item
      $\lefttop{t,\tau}V_{a,b}\bigl(
      \iota_{g_0\dagger}\lefttop{\tau}\Upsilon(\nbigm)
      _{|\lefttop{\tau}\nbigx}
      \bigr)$
      are coherent over
      $\pi_1^{\ast}\bigl(
      \lefttop{t}V\nbigr_{\cnum_t\times\cnum_X}
      \bigr)$,
      where $\pi_1:\cnum_{\tau}\times\cnum_t\times X
      \lrarr \cnum_t\times X$ denotes the projection.
\end{itemize}
Then,
the filtered $\nbigd_X$-module
$(\Xi_{\DR}(\nbigm),F^{\irr}_{b+\bullet})$
is strictly specializable along $g$,
 i.e.,
 $\Rtilde_{F^{\irr}_{b+\bullet}}\Xi_{\DR}(\nbigm)$
is strictly specializable along $g$.
 Moreover, we have
\begin{equation}
\label{eq;22.7.9.31}
 \Gr^V_c
 \Rtilde_{F^{\irr}_{b+\bullet}}(\iota_{g\dagger}\Xi_{\DR}(\nbigm))
 \simeq
 \Rtilde_{F^{\irr}_{b+\bullet}}(\Gr^V_c(\iota_{g\dagger}\Xi_{\DR}(\nbigm)).
\end{equation}
\end{prop}
\pf
We have
$\iota_{g_0\dagger}(\lefttop{\tau}\Upsilon(\nbigm))
\simeq
 \lefttop{\tau}\Upsilon(\iota_{g\dagger}\nbigm)$.
We obtain the filtration
$\lefttop{t,\tau}V_{\bullet,b}\bigl(
 \lefttop{\tau}\Upsilon\bigl(
 \iota_{g\dagger}\nbigm
 \bigr)
 \bigr)$
 of
 $\lefttop{\tau}V_b\bigl(
 \lefttop{\tau}\Upsilon\bigl(
 \iota_{g\dagger}\nbigm
 \bigr)
 \bigr)$.
We have
\begin{equation}
\label{eq;22.7.4.3}
 \lefttop{t,\tau}V_{a,b}
 \bigl(
 \lefttop{\tau}\Upsilon\bigl(
 \iota_{g\dagger}\nbigm\bigr)
 \bigr)
 \big/
  \lefttop{t,\tau}V_{<a,b}
 \bigl(
 \lefttop{\tau}\Upsilon\bigl(
 \iota_{g\dagger}\nbigm\bigr)
 \bigr)
 \simeq
 \lefttop{\tau}V_b\Bigl(
  \lefttop{\tau}\Upsilon\bigl(
 \Gr^V_a\iota_{g\dagger}\nbigm\bigr)
 \Bigr).
\end{equation}
Because the multiplication of $\lambda-\tau$
on (\ref{eq;22.7.4.3}) is a monomorphism for any $a\in\real$,
we obtain that
the multiplication of $\lambda-\tau$
on 
$\lefttop{\tau}V_b\bigl(
 \lefttop{\tau}\Upsilon\bigl(
 \iota_{g\dagger}\nbigm
 \bigr)
 \bigr)
 \Big/
 \lefttop{t,\tau}V_{a,b}\bigl(
 \lefttop{\tau}\Upsilon\bigl(
 \iota_{g\dagger}\nbigm
 \bigr)
 \bigr)$
is a monomorphism.
Thus, we obtain coherent
$\lefttop{t}V\nbigr_{\cnum_t\times X}$-submodules:
\[
 \iota_{\lambda=\tau}^{\ast}
  \lefttop{t,\tau}V_{a,b}\bigl(
 \lefttop{\tau}\Upsilon\bigl(
 \iota_{g\dagger}\nbigm
 \bigr)
 \bigr)
\subset
  \iota_{\lambda=\tau}^{\ast}
  \lefttop{\tau}V_b\bigl(
 \lefttop{\tau}\Upsilon\bigl(
 \iota_{g\dagger}\nbigm
 \bigr)
 \bigr)
\]
Because
\begin{equation}
\label{eq;22.7.9.30}
  \iota_{\lambda=\tau}^{\ast}
  \lefttop{t,\tau}V_{a,b}\bigl(
 \lefttop{\tau}\Upsilon\bigl(
 \iota_{g\dagger}\nbigm
 \bigr)
 \bigr)
 \Big/
   \iota_{\lambda=\tau}^{\ast}
  \lefttop{t,\tau}V_{<a,b}\bigl(
 \lefttop{\tau}\Upsilon\bigl(
 \iota_{g\dagger}\nbigm
 \bigr)
 \bigr)
 \simeq
 \iota_{\lambda=\tau}^{\ast}
 \lefttop{\tau}V_b
 \bigl(
 \lefttop{\tau}\Upsilon(\Gr^V_a\iota_{g\dagger}\nbigm)
 \bigr)
\end{equation}
is strict,
we can check that
the $\nbigr_{\cnum_t\times X}$-module 
$\iota_{\lambda=\tau}^{\ast}
  \lefttop{\tau}V_{b}\bigl(
 \lefttop{\tau}\Upsilon\bigl(
 \iota_{g\dagger}\nbigm
 \bigr)
 \bigr)
 =\Rtilde_{F^{\irr}_{b+\bullet}}
 (\iota_{g\dagger}\Xi_{\DR}\nbigm)$
is strictly specializable along $t$,
and that
$\iota_{\lambda=\tau}^{\ast}
  \lefttop{t,\tau}V_{\bullet,b}\bigl(
 \lefttop{\tau}\Upsilon\bigl(
 \iota_{g\dagger}\nbigm
 \bigr)
 \bigr)$
is a $V$-filtration.
Hence, we obtain that
the filtered $\nbigd_X$-module
$(\Xi_{\DR}(\nbigm),F^{\irr}_{b+\bullet})$
is strictly specializable along $g$. 
We obtain (\ref{eq;22.7.9.31})
from (\ref{eq;22.7.9.30}).
\hfill\qed

\vspace{.1in}

We note that the $V$-filtration of
$\Rtilde_{F^{\irr}_{b+\bullet}}\iota_{g_{\dagger}}\Xi_{\DR}(M)$
is also obtained
in terms of the $V$-filtration and
the irregular Hodge filtration $F^{\irr}_{b+\bullet}$
of $\Xi_{\DR}(M)$.
For each $a\in\real$,
let $F^{\irr}_{b+\bullet}
 \lefttop{t}V_{a}
 \Bigl(
 \iota_{g\dagger}\Xi_{\DR}(\nbigm)
 \Bigr)$
denote the filtration on 
$\lefttop{t}V_{a}
 \Bigl(
 \iota_{g\dagger}\Xi_{\DR}(\nbigm)
 \Bigr)$
 induced by $F^{\irr}_{b+\bullet}$,
i.e.,
\[
 F^{\irr}_{b+\bullet}
\lefttop{t}V_{a}
\Bigl(
\iota_{g\dagger}\Xi_{\DR}(\nbigm)
\Bigr)
=\lefttop{t}V_{a}
 \Bigl(
 \iota_{g\dagger}\Xi_{\DR}(\nbigm)
 \Bigr)
 \cap
  F^{\irr}_{b+\bullet}
\Bigl(
\iota_{g\dagger}\Xi_{\DR}(\nbigm)
\Bigr).
\]
We obtain the $V\nbigr_{\cnum_t\times X}$-module
$\Rtilde_{F^{\irr}_{b+\bullet}}
 \Bigl(
 \lefttop{t}V_{a}
 \Bigl(
 \iota_{g\dagger}\Xi_{\DR}(\nbigm)
 \Bigr)
 \Bigr)$
obtained as the analytification of the Rees module of
$F^{\irr}_{b+\bullet}
 \lefttop{t}V_{a}
 \Bigl(
 \iota_{g\dagger}\Xi_{\DR}(\nbigm)
 \Bigr)$.
\begin{prop}
$\Rtilde_{F^{\irr}_{b+\bullet}}
 \Bigl(
 \lefttop{t}V_{a}
 \Bigl(
 \iota_{g\dagger}\Xi_{\DR}(\nbigm)
 \Bigr)
 \Bigr)
 =
 \iota_{\lambda=\tau}^{\ast}
  \lefttop{t,\tau}V_{a,b}\bigl(
 \lefttop{\tau}\Upsilon\bigl(
 \iota_{g\dagger}\nbigm
 \bigr)
 \bigr)$.
\end{prop}
\pf
By the construction,
$\Rtilde_{F^{\irr}_{b+\bullet}}
 \Bigl(
 \iota_{g\dagger}\Xi_{\DR}(\nbigm)
 \Bigr)
\Big/
 \Rtilde_{F^{\irr}_{b+\bullet}}
 \Bigl(
 \lefttop{t}V_{a}
 \Bigl(
 \iota_{g\dagger}\Xi_{\DR}(\nbigm)
 \Bigr)
 \Bigr)$
is strict.
Then, the claim follows from 
Lemma \ref{lem;21.3.19.10} below,
\hfill\qed

\begin{lem}
\label{lem;21.3.19.10}
 Let $(M,F)$ be a filtered $\nbigd_Y$-module.
Let $\Rtilde_F(M)$ denote the associated
$\nbigr_Y$-module. 
Let $\nbigg_i\subset \Rtilde_F(M)$ $(i=1,2)$
be $\nbigo_{\nbigy}$-submodules
 such that
 (i) $\nbigg_1(\ast\lambda)=\nbigg_2(\ast\lambda)$,
 (ii) $\Rtilde_F(M)/\nbigg_i$ are strict.
Then, $\nbigg_1=\nbigg_2$. 
\end{lem}
\pf
By the conditions,
the induced morphisms $\nbigg_1\lrarr \Rtilde_F(M)/\nbigg_2$
and  $\nbigg_2\lrarr \Rtilde_F(M)/\nbigg_1$
are $0$.
Hence, we obtain that $\nbigg_1=\nbigg_2$.
\hfill\qed

\begin{cor}
\label{cor;22.7.9.23}
Let $\nbigm$ be as in Proposition {\rm\ref{prop;21.6.20.1}}.
We set
$M:=\Xi_{\DR}(\nbigm)$
and
$M_1:=\iota_{g\dagger}M$.
Let $V_{\bullet}(M_1)$ denote the $V$-filtration of $M_1$
along $t$.
\begin{itemize}
 \item 
 For $c<0$,
we have 
\begin{equation}
\label{eq;22.7.9.20}
      F_{b+j}^{\irr}(M_1)
      \cap V_{c}(M_1)
      =\bigl(
      F_{b+j}^{\irr}(M_1)(\ast t)
      \bigr)
      \cap
      V_{c}(M_1).
\end{equation}
For $c>0$, we have
\begin{equation}
\label{eq;22.7.9.21}
      F_{b+j}^{\irr}(M_1)
      \cap V_{c}(M_1)
      =\sum_{\substack{a\leq 0,
      n\in\seisuu_{\geq 0} \\
      a+n\leq c}}
	\del_t^n
	\bigl(
	F_{b+j-n}^{\irr}(M_1)
	\cap V_{a}(M_1)
	\bigr).
\end{equation}
\item If moreover $M_1=M_1(\ast t)$,
      {\rm(\ref{eq;22.7.9.20})} holds
      for $c=0$.
\item If moreover $M_1=M_1(!t)$,
      we have
\begin{equation}
\label{eq;22.7.9.22}
       F_{b+j}^{\irr}(M_1)
      \cap V_{0}(M_1)
      =\sum_{\substack{a<0,
      n\in\seisuu_{\geq 0} \\
      a+n\leq 0}}
	\del_t^n
	\bigl(
	F_{b+j-n}^{\irr}(M_1)
	\cap V_{a}(M_1)
	\bigr).
\end{equation}
\end{itemize}
\end{cor}
\pf
Because 
the induced morphisms
$\deldel_t:
 \Gr^V_{c}(\Rtilde_{F^{\irr}_{b+\bullet}}(M_1))
 \lrarr
 \Gr^V_{c+1}(\Rtilde_{F^{\irr}_{b+\bullet}}(M_1))$
 are isomorphisms for $c>-1$,
 we obtain (\ref{eq;22.7.9.21}).

We set
$\Ftilde^{\irr}_{b+j}(M_1(\ast t)):=
 F^{\irr}_{b+j}(M_1)(\ast t)$ for any $j\in\seisuu$.
We also set
$\Ftilde^{\irr}_{b+j}V_c(M_1(\ast t)):=
V_c(M_1(\ast t))
\cap
\Ftilde^{\irr}_{b+j}(M_1(\ast t))$.
Note that
$\Rtilde_{\Ftilde^{\irr}_{b+\bullet}}(M_1(\ast t))
=\Rtilde_{F^{\irr}_{b+\bullet}}(M_1)(\ast t)$
is a strictly specializable $\nbigrtilde_{\cnum_t\times X}(\ast t)$-module,
and its $V$-filtration is
$\Rtilde_{\Ftilde^{\irr}_{b+\bullet}}(V_cM_1(\ast t))$ $(c\in\real)$.
Because 
$V_c\Rtilde_{F^{\irr}_{b+\bullet}}(M_1)
=V_c\Rtilde_{F^{\irr}_{b+\bullet}}(M_1(\ast t))$ for $c<0$,
we obtain (\ref{eq;22.7.9.20}).

We set $\nbigm_1=\iota_{g\dagger}\nbigm$. If $M_1=M_1(\ast t)$,
we obtain $\nbigm_1[\ast t]=\nbigm_1$,
and
$\lefttop{\tau}\Upsilon(\nbigm_1)[\ast t]
=\lefttop{\tau}\Upsilon(\nbigm_1)$.
It implies that
$\lefttop{t,\tau}V_{0,b}(\lefttop{\tau}\Upsilon(\nbigm_1)(\ast t))
=\lefttop{t,\tau}V_{0,b}(\lefttop{\tau}\Upsilon(\nbigm_1))$.
Hence, we obtain (\ref{eq;22.7.9.20}) for $c=0$.
If $M_1=M_1(!t)$,
the induced morphism
$\deldel_t:
 \Gr^V_{-1}(\Rtilde_{F^{\irr}_{b+\bullet}}(M_1))
 \lrarr
 \Gr^V_{0}(\Rtilde_{F^{\irr}_{b+\bullet}}(M_1))$
 is an isomorphism,
 and hence we obtain (\ref{eq;22.7.9.22}).
\hfill\qed

\begin{thm}
Suppose that $\nbigm\in\nbigc_{\res}(X)$ is regular. 
We set $M:=\Xi_{\DR}(\nbigm)$.
Then, for any open subset $U\subset X$
with a holomorphic function $g\in\nbigo(U)$,
the filtered $\nbigd$-module
$(M,F_{b+\bullet}^{\irr})_{|U}$ is strictly specializable along $g$.
 The $V$-filtration
and the irregular Hodge filtration $F^{\irr}_{b+\bullet}$
of $\iota_{g\dagger}M$ satisfy
the compatibility conditions in Corollary {\rm\ref{cor;22.7.9.23}}.
\end{thm}
\pf
It follows from Theorem \ref{thm;21.3.17.20},
Proposition \ref{prop;21.6.20.1}.
and Corollary {\rm\ref{cor;22.7.9.23}}.
\hfill\qed

\subsubsection{Complement}

Let $\nbigm\in\nbigc_{\res}(X)$ be regular. 
We set $M:=\Xi_{\DR}(\nbigm)$.
Let $g$ be a holomorphic function on $X$.
Let us explain that
$F^{\irr}_{b+\bullet}(M)$
is determined by
$F^{\irr}_{b+\bullet}(M)(\ast g)$
and
$F^{\irr}_{b+\bullet}(\phi_g^{(0)}M)$.
We set $\nbigm_1=\iota_{g\dagger}\nbigm$
and $M_1=\iota_{g\dagger}M$ on $\cnum_t\times X$.
 \begin{itemize}
  \item The filtered $\nbigd_{\cnum_t\times X}(\ast t)$-module
       $(M_1(\ast t),F^{\irr}_{b+\bullet})$
       is obtained as the push-forward of
       the filtered $\nbigd_X(\ast g)$-module
       $(M(\ast g),F^{\irr}_{b+\bullet}(M)(\ast g))$.
  \item $(\Pi^{a,b}_{t}(M_1(\ast t)),F^{\irr}_{b+\bullet})$
	are determined by
	the filtered $\nbigd_{\cnum_t\times X}(\ast t)$-module
	$(M_1(\ast t),F^{\irr}_{b+\bullet})$
       and
	the Hodge filtration of
	$\nbigd_{\cnum_t}(\ast t)$-module
       $\II^{a,b}_{t}$.
  \item The filtered $\nbigd_{\cnum_t\times X}$-modules
	$(\Pi^{a_1,a_2}_{g\star}(M_1),F^{\irr}_{b+\bullet})$
	$(\star=!,\ast)$
	are determined
	by $(\Pi^{a,b}_{t}(M_1(\ast t)),F^{\irr}_{b+\bullet})$
       according to Corollary \ref{cor;22.7.9.23}.
       Then,
       $(\Pi^{a_1,a_2}_{t\ast!}(M_1),F^{\irr}_{b+\bullet})$
       are naturally determined.
       In particular,
       $F^{\irr}_{b+\bullet}$ of
       $\Xi^{(0)}_{t}(M_1)$
       and $\psi^{(i)}_t(M_1)$ $(i=0,1)$
       are determined.
\item 
Hence, $(M_1,F^{\irr}_{b+\bullet})$
is isomorphic to the cohomology of
the following complex which is strictly compatible with
$F^{\irr}_{b+\bullet}$:
\[
 (\psi^{(1)}_t(M_1),F^{\irr}_{b+\bullet})
 \lrarr
 (\Xi^{(0)}_t(M_1),F^{\irr}_{b+\bullet})
 \oplus
  (\phi^{(0)}_t(M_1),F^{\irr}_{b+\bullet})
 \lrarr
  (\psi^{(0)}_t(M_1),F^{\irr}_{b+\bullet}).
\]     
  \item
In other words,  $(M,F^{\irr}_{b+\bullet})$
is isomorphic to the cohomology of
the following complex which is strictly compatible with
$F^{\irr}_{b+\bullet}$:
\[
 (\psi^{(1)}_g(M),F^{\irr}_{b+\bullet})
 \lrarr
 (\Xi^{(0)}_g(M),F^{\irr}_{b+\bullet})
 \oplus
  (\phi^{(0)}_g(M),F^{\irr}_{b+\bullet})
 \lrarr
  (\psi^{(0)}_g(M),F^{\irr}_{b+\bullet}).
\]     
 \end{itemize}

\section{Rescalability at $\infty$}
\label{section;21.7.3.1}

\subsection{Preliminary}

\subsubsection{Rescaling of $\nbigrtilde_X$-modules}

Let $\varphi_{3}:
\cnum_{\lambda}\times\cnum_{\xi}
\lrarr
\cnum_{\lambda}$
be the map defined by
$\varphi_3(\lambda,\xi)=\lambda\xi$.
We set
$\lefttop{\xi}X:=\cnum_{\xi}\times X$
and $\lefttop{\xi}X_0:=\{0\}\times X\subset\lefttop{\xi}X$.
The induced morphism
$\lefttop{\xi}\nbigx\lrarr \nbigx$
is also denoted by $\varphi_3$.

For any $\nbigo_{\nbigx}$-module $\nbigm$,
we obtain
the $\nbigo_{\lefttop{\xi}\nbigx}(\ast\lefttop{\xi}\nbigx_0)$-module
$\lefttop{\xi}\nbigm:=
\varphi_3^{\ast}(\nbigm)(\ast\xi)$.
If $\nbigm$ is equipped with a meromorphic flat connection $\nabla$,
$\lefttop{\xi}\nbigm$ is equipped with 
the induced meromorphic flat connection $\varphi_3^{\ast}\nabla$.
In this way,
if $\nbigm$ is an $\nbigrtilde_X$-module,
$\lefttop{\xi}\nbigm$ is naturally
an $\nbigrtilde_{\lefttop{\xi}X}(\ast\xi)$-module.
The following lemma is similar to Lemma \ref{lem;21.3.13.30}.
\begin{lem}[\cite{Sabbah-irregular-Hodge}]
Let $\nbigm\in\nbigc(X)$.
\begin{itemize}
 \item For any projective morphism $F:X\lrarr Y$,
       we have
       $F_{\dagger}(\lefttop{\xi}\nbigm)\simeq
       \lefttop{\xi}F_{\dagger}(\nbigm)$.
 \item For any hypersurface $H$ of $X$,
       the $\nbigrtilde_{\lefttop{\xi}X}(\ast\xi)$-module
       $\lefttop{\xi}\nbigm$ is localizable along
       $\lefttop{\xi}H$,
       and we have the isomorphism of
       $\nbigrtilde_{\lefttop{\xi}X}(\ast\xi)$-modules
       $\lefttop{\xi}(\nbigm[\star H])\simeq
       (\lefttop{\xi}\nbigm)[\star \lefttop{\xi}H]$.
 \item Let $f$ be any holomorphic function on $X$.
       Let $f_0$ denote the induced holomorphic function
       on $\lefttop{\xi}X$.
       For any $a\leq b$, we have
       $\lefttop{\xi}(\Pi^{a,b}_{f\star}(\nbigm))
       \simeq
       \Pi^{a,b}_{f_0\star}(\lefttop{\xi}\nbigm)$,
       and
       $\lefttop{\xi}(\Pi^{a,b}_{f,\ast !}(\nbigm))
       \simeq
       \Pi^{a,b}_{f_0,\ast !}(\lefttop{\xi}\nbigm)$.
       Hence, we have
       $\xi\bigl(
        \Xi^{(a)}_f(\nbigm)
       \bigr)
       =\Xi^{(a)}_{f_0}(\lefttop{\xi}\nbigm)$,
       $\lefttop{\xi}\bigl(
        \psi^{(a)}_f(\nbigm)
       \bigr)
       =\psi^{(a)}_{f_0}(\lefttop{\xi}\nbigm)$,
       and
       $\phi^{(0)}_f(\nbigm)
       =\phi^{(0)}_{f_0}(\lefttop{\xi}\nbigm)$.
       \hfill\qed
\end{itemize}
\end{lem}

\subsubsection{Rescaling of $\gbigrtilde_X$-modules}

Let $\varphi_4:\widetilde{\proj_{\lambda}^1\times\cnum_{\xi}}\lrarr
\proj_{\lambda}^1\times\cnum_{\xi}$
denote the blow up at $(\infty,0)$.
The morphism $\widetilde{\proj_{\lambda}^1\times\cnum_{\xi}}\lrarr
\proj^1_{\lambda}$
induced by $\varphi_3$
is also denoted by $\varphi_3$.
The induced morphisms
$\widetilde{\proj^1_{\lambda}}\times\cnum_{\xi}\times X
\lrarr \lefttop{\xi}\gbigx$
and
$\widetilde{\proj^1_{\lambda}}\times\cnum_{\xi}\times X
\lrarr \gbigx$
are also denoted by $\varphi_4$
and $\varphi_3$, respectively.

For an $\nbigo_{\gbigx}(\ast\gbigx^{\infty})$-module $\gbigm$,
we obtain
the $\nbigo_{\lefttop{\xi}\gbigx}\bigl(
\ast(\lefttop{\xi}\gbigx^{\infty}\cup\lefttop{\xi}\gbigx_0)
\bigr)$-module
$\lefttop{\xi}\gbigm:=
\varphi_{4\ast}
\varphi_3^{\ast}(\gbigm)
\bigl(
\ast(\lefttop{\xi}\gbigx^{\infty}\cup\lefttop{\xi}\gbigx_0)
\bigr)$.
If $\gbigm$ is an $\gbigrtilde_X$-module,
then $\lefttop{\xi}\gbigm$
is naturally an $\gbigrtilde_{\lefttop{\xi}X}(\ast\xi)$-module.

\subsection{Rescalability at $\infty$}

We shall prove the following theorem in
\S\ref{subsection;21.7.2.11}--\S\ref{subsection;21.7.25.10}.
\begin{thm}
\label{thm;21.7.2.20}
Let $\nbigm\in\nbigc(X)$.  
Suppose that there exists
$\nbigm_0\in\nbigc(\lefttop{\xi}X)$ such that
$\nbigm_{0|\lefttop{\xi}X\setminus\lefttop{\xi}X_0}
 =(\lefttop{\xi}\nbigm)
  _{|\lefttop{\xi}X\setminus\lefttop{\xi}X_0}$.
Then, we have $\nbigm_0(\ast\xi)=\lefttop{\xi}\nbigm$.
In particular, we obtain
$\lefttop{\xi}\nbigm\in\nbigc(\lefttop{\xi}X;\lefttop{\xi}X_0)$. 
\end{thm}

Let $\nbigc_{\res,\infty}(X)\subset\nbigc(X)$
denote the full subcategory of
$\nbigm\in\nbigc(X)$ such that
$\lefttop{\xi}\nbigm\in\nbigc(\lefttop{\xi}X;\lefttop{\xi}X_0)$.
We shall prove the following proposition in \S\ref{subsection;21.7.25.11}.
\begin{prop}
\label{prop;21.7.2.30}
Let $\nbigm\in\nbigc_{\res,\infty}(X)$.
Then,
$\Upsilon\bigl(\lefttop{\xi}\nbigm\bigr)
=\lefttop{\xi}\Upsilon(\nbigm)$. 
\end{prop}

Let $\gbigc_{\res,\infty}(X)\subset\gbigc(X)$
denote the essential image of
$\Upsilon:\nbigc_{\res,\infty}(X)\lrarr\gbigc(X)$.
According to Proposition \ref{prop;21.7.2.30},
$\gbigm\in\gbigc_{\Malg}(X)$ is an object of
$\gbigc_{\res,\infty}(X)$
if and only if
$\lefttop{\xi}\gbigm\in\gbigc_{\Malg}(\lefttop{\xi}X;\lefttop{\xi}X_0)$.

Let $H$ be a hypersurface of $X$.
Let $\nbigc_{\res,\infty}(X;H)$ denote the essential image of
$\nbigc_{\res,\infty}(X)\lrarr\nbigc(X;H)$.
Similarly,
let $\gbigc_{\res,\infty}(X;H)$ denote the essential image of
$\gbigc_{\res,\infty}(X)\lrarr\gbigc(X;H)$.

\begin{example}
$\gbigr_F(M)\otimes\gbigl(f)\in\gbigc_{\Malg}(X;H)$
in {\rm\S\ref{subsection;22.7.27.10}}
is an object of $\gbigc_{\res,\infty}(X;H)$.
(See Example {\rm\ref{example;22.7.27.11}}.)
Hence,
$\bigl(
\gbigr_F(M)\otimes\gbigl(f)
\bigr)[\star H]
\in\gbigc_{\res,\infty}(X)$.
\hfill\qed
\end{example}

\subsubsection{The $0$-dimensional pure case}
\label{subsection;22.7.5.2}

Let $\pt$ denote the one point set.
Let $\nbigm\in\nbigc(\pt)$.
We set $\cnum_{\xi}^{\ast}=\cnum_{\xi}\setminus\{0\}$.
Suppose that
there exists 
an integrable pure twistor $\nbigd$-module
$\nbigt_1=(\nbigm_1,\nbigm_1,C_1)$ of weight $0$ with
the integrable polarization $\nbigs_1=(\id,\id)$
on $\cnum_{\xi}$ such that
$\nbigm_{1|\cnum_{\xi}^{\ast}}=
\lefttop{\xi}\nbigm_{|\cnum_{\xi}^{\ast}}$.
We obtain the meromorphic flat bundle
$\nbigv=\nbigm_1\bigl(\ast(\lambda\xi)\bigr)$
on $(\cnum_{\lambda}\times\cnum_{\xi},\{\lambda\xi=0\})$.

Let $(E,\delbar_E,\theta,h)$ be a harmonic bundle
corresponding to $(\nbigt_1,\nbigs_1)_{|Y^{\circ}}$.
Because the Higgs bundle $(E,\delbar_E,\theta)$
is isomorphic to
$\nbigm_{1|\cnum_{\lambda}\times\cnum_{\xi}^{\ast}}
\big/\lambda\nbigm_{1|\cnum_{\lambda}\times\cnum_{\xi}^{\ast}}
\simeq
\lefttop{\xi}\nbigm_{|\cnum_{\lambda}\times\cnum_{\xi}^{\ast}}\big/
\lambda
\lefttop{\xi}\nbigm_{|\cnum_{\lambda}\times\cnum_{\xi}^{\ast}}$,
it is naturally $\cnum^{\ast}$-homogeneous.
The following lemma is clear.
\begin{lem}
\label{lem;21.7.2.1}
There exist a finite set
$S\subset\cnum$
and a decomposition
$(E,\delbar_E)=\bigoplus_{\alpha\in S}(E_{\alpha},\delbar_{E_{\alpha}})$
such that the following holds.
\begin{itemize}
 \item $\theta(E_{\alpha})\subset E_{\alpha}\otimes\Omega^1_{Y^{\circ}}$,
       and
       $\theta_{|E_{\alpha}}-\alpha d(\xi^{-1})\id_{E_{\alpha}}$
       is nilpotent.
       \hfill\qed
\end{itemize} 
\end{lem}

\begin{lem}
\label{lem;22.7.5.1}
$\nbigv_{|\cnum_{\lambda}^{\ast}\times\cnum_{\xi}}$
is unramifiedly good along $\xi=0$.
The good set of irregular values is  
$\{\alpha \lambda^{-1}\xi^{-1}\,|\,\alpha\in S\}$.
\end{lem}
\pf
It follows from
\cite[Theorem 11.1.2, Theorem 12.6.1]{Mochizuki-wild}.
\hfill\qed

\vspace{.1in}

Let $\rho_0,\rho_1:
\cnum^{\ast}\times(\cnum_{\lambda}\times\cnum_{\xi})\to
\cnum_{\lambda}\times\cnum_{\xi}$
be given by
$\rho_0(a,\lambda,\xi)=(\lambda,\xi)$
and $\rho_1(a,\lambda,\xi)=(a\lambda,a^{-1}\xi)$.
Because
$\nbigm_{1|\cnum_{\lambda}\times\cnum_{\xi}^{\ast}}
=\lefttop{\xi}\nbigm_{|\cnum_{\lambda}\times\cnum_{\xi}^{\ast}}$,
there exists the following isomorphism
satisfying the cocycle condition:
\begin{equation}
\label{eq;22.7.4.20}
 (\rho^{\ast}_0\nbigv)
 _{|\cnum^{\ast}\times(\cnum_{\lambda}\times\cnum^{\ast}_{\xi})}
 \simeq
 (\rho^{\ast}_1\nbigv)
 _{|\cnum^{\ast}\times(\cnum_{\lambda}\times\cnum^{\ast}_{\xi})}.
\end{equation}
\begin{lem}
The isomorphism {\rm(\ref{eq;22.7.4.20})} 
uniquely extends to an isomorphism
$(\rho^{\ast}_0\nbigv)
\simeq
(\rho^{\ast}_1\nbigv)$
on $\cnum^{\ast}\times(\cnum_{\lambda}\times \cnum_{\xi})$
satisfying the cocycle condition,
i.e., $\nbigv$ is $\cnum^{\ast}$-homogeneous.
\end{lem}
\pf
The restriction of 
the isomorphism {\rm(\ref{eq;22.7.4.20})}
to $\{1\}\times(\cnum_{\lambda}^{\ast}\times\cnum_{\xi}^{\ast})$
extends to an isomorphism
on $\{1\}\times(\cnum_{\lambda}^{\ast}\times\cnum_{\xi})$.
Both
$\rho_i^{\ast}(\nbigv)
_{|\cnum^{\ast}\times\cnum_{\lambda}^{\ast}\times\cnum_{\xi}}$
are unramifiedly good along $\xi=0$,
and the good set of irregular values are
$\{\alpha\lambda^{-1}\xi^{-1}\,|\,\alpha\in S\}$.
Hence, by \cite[Theorem 4.4.1]{Mochizuki-wild},
the isomorphism {\rm(\ref{eq;22.7.4.20})}
extends to an isomorphism
$(\rho^{\ast}_0\nbigv)
_{|\cnum^{\ast}\times((\cnum_{\lambda}\times\cnum_{\xi})\setminus\{(0,0)\})}
\simeq
(\rho^{\ast}_1\nbigv)
_{|\cnum^{\ast}\times((\cnum_{\lambda}\times\cnum_{\xi})\setminus\{(0,0)\})}$.
By the theorem of Hartogs,
we obtain the desired isomorphism.
\hfill\qed

\vspace{.1in}

\begin{lem}
\label{lem;22.7.4.21}
$\nbigv$ is an unramifiedly good meromorphic flat bundle
on $(\nbigy,\{\lambda\xi=0\})$.
The set of irregular values is
$\{\alpha \lambda^{-1}\xi^{-1}\,|\,\alpha\in S\}$. 
\end{lem}
\pf
By taking the Malgrange extension along
$\{\infty\}\times \cnum_{\xi}$,
we obtain the meromorphic flat bundle
$\Upsilon(\nbigv)$
on $\proj^1_{\lambda}\times\cnum_{\xi}$,
which is also $\cnum^{\ast}$-homogeneous.
Because $\nbigv_{|\cnum_{\lambda}\times\cnum_{\xi}^{\ast}}
\simeq
\lefttop{\xi}\nbigm(\ast\lambda)
_{|\cnum_{\lambda}\times\cnum_{\xi}^{\ast}}$
is constructed by rescaling,
$\Upsilon(\nbigv)$
naturally extends to a $\cnum^{\ast}$-homogeneous
meromorphic flat bundle $\nbigvtilde$
on $\proj^1_{\lambda}\times\proj^1_{\xi}$.

Let $\nbigvtilde^{\DM}$ denote the Deligne-Malgrange lattice
of $\nbigvtilde$,
which is a $\cnum^{\ast}$-invariant lattice of $\nbigvtilde$.
We set
$\nbigl:=\nbigvtilde^{\DM}
\bigl(\ast\bigl(
 (\{\infty\}\times\proj^1_{\xi})
 \cup(\proj^1_{\lambda}\times\{\infty\})
 \bigr)
\bigr)$.
There exists the $\cnum^{\ast}$-equivariant surjection
$H^0(\proj^1_{\lambda}\times \proj^1_{\xi},\nbigl)
\lrarr \nbigl_{|(0,0)}$.
We choose a base
$v_1,\ldots,v_r$ of $\nbigl_{|(0,0)}$
such that $a^{\ast}(v_i)=a^{\rho_i}v_i$ for some $\rho_i\in\seisuu$.
There exist
$\vtilde_i\in H^0(\proj^1_{\lambda}\times \proj^1_{\xi},\nbigl)$
such that
$\vtilde_{i|(0,0)}=v_i$
and that
$a^{\ast}(\vtilde_i)=a^{\rho_i}\vtilde_i$.
Note that
$\vtilde_1,\ldots,\vtilde_r$
is a frame of the restriction of $\nbigl$
on a $\cnum^{\ast}$-invariant Zariski open neighbourhood $\nbigu$ of
$\{\lambda\xi=0\}$.

We set
$\vtilde^{(1)}_i=\lambda^{-\rho_i}\vtilde_i$.
Then,
$\vtilde^{(1)}_1,\ldots,\vtilde^{(1)}_r$
is a $\cnum^{\ast}$-invariant frame of $\nbigv$  on $\nbigu$.
Let $\nbiga_{\lambda}\,d\lambda+\nbiga_{\xi}d\xi$
denote the connection form of $\nabla$
with respect to
the frame $\vtilde^{(1)}_1,\ldots,\vtilde^{(1)}_r$.
We have the naturally defined morphism
$\nbigu\lrarr \cnum$ induced by
$(\lambda,\xi)\longmapsto \lambda\xi$.
Let $U\subset\cnum$ denote the image
which is a neighbourhood of $0$.
We obtain the isomorphism
$(\nbigu\setminus\{\lambda\xi=0\})/\cnum^{\ast}
\simeq
 U\setminus\{0\}$.
There exists an $M_r(\cnum)$-valued holomorphic function
$B$ on $U\setminus\{0\}$
such that
$\nbiga_{\lambda}\,d\lambda
+\nbiga_{\xi}d\xi=B(\lambda\xi)\,d(\lambda\xi)$.
We obtain that
$\nbiga_{\xi}(\lambda,\xi)=B(\lambda\xi)\lambda$
and
$\nbiga_{\lambda}(\lambda,\xi)=B(\lambda\xi)\xi$.
By Lemma \ref{lem;22.7.5.1},
we obtain that
$\nbigv_{|\cnum_{\lambda}\times\cnum_{\xi}^{\ast}}$
is unramifiedly good along $\lambda=0$,
and that the set of the irregular values
along $\lambda=0$
is $\{\alpha\lambda^{-1}\xi^{-1}\,|\,\alpha\in S\}$.
Thus, we obtain the claim of Lemma \ref{lem;22.7.4.21}
from Proposition \ref{prop;21.1.22.2}.
\hfill\qed

\begin{lem}
\label{lem;22.7.4.22}
$\nbigm_{1}(\ast\xi)
=\lefttop{\xi}\nbigm$.
\end{lem}
\pf
Because
$\lefttop{\xi}\nbigm_{|\cnum_{\lambda}\times \cnum_{\xi}^{\ast}}(\ast\lambda)
=\nbigv_{|\cnum_{\lambda}\times\cnum_{\xi}^{\ast}}$,
we obtain that
$\lefttop{\xi}\nbigm_{|\cnum_{\lambda}\times\cnum_{\xi}^{\ast}}(\ast\lambda)$
is unramifiedly good along
$\{0\}\times\cnum_{\xi}^{\ast}$,
and the set of irregular values is
$\{\alpha\lambda^{-1}\xi^{-1}\,|\,\alpha\in S\}$.
By the construction,
$\lefttop{\xi}
\nbigm_{|\cnum^{\ast}_{\lambda}\times\cnum_{\xi}}$
is unramifiedly good along $\xi=0$,
and the set of irregular values is
$\{\alpha\lambda^{-1}\xi^{-1}\,|\,\alpha\in S\}$.
By Proposition \ref{prop;21.1.22.2},
we obtain that
$\lefttop{\xi}\nbigm$ is
unramifiedly good
on $(\cnum_{\lambda}\times\cnum_{\xi},\{\lambda\xi=0\})$.
We obtain that
the isomorphism
$\nbigv_{|\cnum\times\cnum_{\xi}^{\ast}}
\simeq
\lefttop{\xi}\nbigm(\ast\lambda)_{|\cnum\times\cnum_{\xi}^{\ast}}$
extends to an isomorphism
$\nbigv
\simeq
\lefttop{\xi}\nbigm(\ast\lambda)$.
(See \cite[Proposition 2.10]{Mochizuki-subanalytic2}.)
Thus, we obtain Lemma \ref{lem;22.7.4.22}.
\hfill\qed

\subsubsection{The $0$-dimensional mixed case}

Let $\nbigm\in \nbigc(\pt)$.
Suppose that
there exists 
$\nbigm_0\in\nbigc(\cnum_{\xi})$
such that
$\nbigm_{0|\cnum_{\xi}^{\ast}}
\simeq\lefttop{\xi}\nbigm_{|\cnum_{\xi}^{\ast}}$.
Here, we recall that
$\nbign_{|Z}$
means
$\nbign_{|\cnum_{\lambda}\times Z}$
for an $\nbigr_{Y}$-module $\nbign$
and $Z\subset Y$.

\begin{lem}
\label{lem;21.7.2.21}
 We have
$\nbigm_0(\ast\xi)\simeq\lefttop{\xi}\nbigm$.
\end{lem}
\pf
There exists a filtration $W$ of $\nbigm_0$
by $\nbigrtilde_Y$-submodules
such that each $\Gr^W_w(\nbigm_0)$
underlies an integrable polarizable
pure twistor $\nbigd$-module of weight $w$.
We obtain the induced filtration
$W$ on $\nbigm_0(\ast\xi)$.

Because $\nbigm_{0|\cnum_{\lambda}\times\{1\}}
\simeq \nbigm$,
we obtain the induced filtration $W$ on $\nbigm$.
We have
$\lefttop{\xi}(W_w\nbigm)_{|\cnum_{\xi}^{\ast}}
\simeq
W_w\nbigm_{0|\cnum_{\xi}^{\ast}}$
and
$\lefttop{\xi}(\Gr^W_w\nbigm)_{|\cnum_{\xi}^{\ast}}
\simeq
\Gr^W_w\nbigm_{0|\cnum_{\xi}^{\ast}}$.
By Lemma \ref{lem;22.7.4.22},
we obtain
$\lefttop{\xi}\bigl(
\Gr^W_w\nbigm
\bigr)
\simeq
\Gr^W_w(\nbigm_0(\ast\xi))$.
Moreover,
each $\lefttop{\xi}\bigl(
\Gr^W_w\nbigm
\bigr)(\ast(\{0\}\times \cnum_{\xi}))$
is an unramifiedly good meromorphic flat bundle.
By \cite[Proposition 2.2.13]{Mochizuki-wild},
both
$\lefttop{\xi}\nbigm\bigl(\ast(\{0\}\times\cnum_{\xi})\bigr)$
and
$\nbigm_0(\ast\xi)\bigl(\ast(\{0\}\times\cnum_{\xi})\bigr)$
are unramifiedly good,
and the set of irregular values are the same,
and of the form
$\{\alpha \lambda^{-1}\xi^{-1}\,|\,\alpha\in S\}$
for a finite subset $S\subset\cnum$.
Then, the isomorphism
on $\cnum_{\lambda}\times\cnum_{\xi}^{\ast}$
extends to an isomorphism on $\cnum_{\lambda}\times\cnum_{\xi}$.
(See \cite[Proposition 2.10]{Mochizuki-subanalytic2}.)
\hfill\qed

\subsubsection{Smooth case}

Let $\nbigm\in\nbigc(X)$.
We assume that $\nbigm$ is a smooth $\nbigrtilde_X$-module.
Suppose that there exist
$\nbigm_0\in\nbigc(\lefttop{\xi}X)$
such that
$\nbigm_{0|\lefttop{\xi}X\setminus \lefttop{\xi}X_0}\simeq
\lefttop{\xi}\nbigm_{|\lefttop{\xi}X\setminus \lefttop{\xi}X_0}$.
\begin{lem}
\label{lem;21.7.2.23}
 We have
$\nbigm_0(\ast\xi)\simeq\lefttop{\xi}\nbigm$.
\end{lem}
\pf
For any point $P\in X$,
the isomorphism
$\nbigm_{0|\cnum_{\xi}^{\ast}\times P}\simeq
\lefttop{\xi}\nbigm_{|\cnum_{\xi}^{\ast}\times P}$
extends to an isomorphism
$\nbigm_{0|\cnum_{\xi}\times P}\simeq
\lefttop{\xi}\nbigm_{|\cnum_{\xi}\times P}$
by Lemma \ref{lem;21.7.2.21}.
Then, we obtain the claim of the lemma.
\hfill\qed

\subsubsection{The general case}
\label{subsection;21.7.25.10}
Let us prove the claim of Theorem \ref{thm;21.7.2.20}
in the general case.
We use an induction on the dimension of the support.
It is enough to prove the claim locally around any point of $X$.
We may assume the existence of 
$f$, $Z$, $\varphi$ and $\nbigv$
as in the proof of Proposition \ref{prop;21.2.18.1}.
We obtain the description of
$\lefttop{\xi}\nbigm$
as the cohomology of the following complex
of $\nbigrtilde_{\lefttop{\xi}X}(\ast\xi)$-modules:
\[
\lefttop{\xi}\bigl(
 \psi^{(1)}_f(\nbigm)
\bigr)
\lrarr
 \lefttop{\xi}\bigl(
 \phi^{(0)}_f(\nbigm)
 \bigr)
 \oplus
\lefttop{\xi}\bigl(
 \Xi^{(0)}_f(\nbigm)
 \bigr)
 \lrarr
 \lefttop{\xi}\bigl(
 \psi^{(0)}_f(\nbigm)
\bigr).
\]
We also have the description of
$\nbigm_0(\ast\xi)$
as the cohomology of the following complex of
$\nbigrtilde_{\lefttop{\xi}X}(\ast\xi)$-modules:
\[
 \psi^{(1)}_{f_0}(\nbigm_0(\ast\xi))
 \lrarr
 \Xi^{(0)}_{f_0}(\nbigm_0(\ast\xi))
 \oplus
 \phi^{(0)}_{f_0}(\nbigm_0(\ast\xi))
 \lrarr
 \psi^{(0)}_{f_0}(\nbigm_0(\ast\xi)).
\]
By the assumption of the induction,
we have
$\lefttop{\xi}\bigl(
 \psi^{(a)}_f(\nbigm)
 \bigr)
 = \psi^{(a)}_{f_0}(\nbigm_0(\ast\xi))$
and
$\lefttop{\xi}\bigl(
 \phi^{(0)}_f(\nbigm)
 \bigr)
 =\phi^{(0)}_{f_0}(\nbigm_0(\ast\xi))$.

Let $\varphi_0:\lefttop{\xi}Z\lrarr \lefttop{\xi}X$
denote the induced morphism.
We set $H=f^{-1}(0)$
and $H_Z=\varphi^{-1}(H)\subset Z$.
There exists
$\nbigv_0\in
 \nbigc\bigl(\lefttop{\xi}Z;
\lefttop{\xi}H_Z\cup\lefttop{\xi}Z_0\bigr)$
such that
$\varphi_{0\dagger}(\nbigv_0)
=\nbigm_{0}(\ast (\lefttop{\xi}H\cup \lefttop{\xi}X_0))$.
By using Lemma \ref{lem;21.7.2.23},
we obtain that
$\lefttop{\xi}\nbigv_{|\lefttop{\xi}Z\setminus\lefttop{\xi}H}
\simeq
 \nbigv_{0|\lefttop{\xi}Z\setminus\lefttop{\xi}H}$.
By the construction,
we have
$\lefttop{\xi}\nbigv_{|\lefttop{\xi}Z\setminus\lefttop{\xi}Z_0}
\simeq
 \nbigv_{0|\lefttop{\xi}Z\setminus\lefttop{\xi}Z_0}$.
By using the Hartogs theorem,
we obtain that
$\lefttop{\xi}\nbigv
\simeq
\nbigv_0$.
Let $f_1:=\varphi_0^{-1}(f_0)$.
We obtain
$\lefttop{\xi}\Xi^{(a)}_{f_1}(\nbigv)
\simeq 
\Xi^{(a)}_{f_1}(\lefttop{\xi}\nbigv)
\simeq
\Xi^{(a)}_{f_1}(\nbigv_0)$.
It follows that
$\lefttop{\xi}\bigl(
\Xi^{(0)}_f(\nbigm)
\bigr)
\simeq
\Xi^{(0)}_{f_0}(\nbigm_0(\ast\xi))$.
Hence, we obtain
$\lefttop{\xi}\nbigm\simeq
\nbigm_{0}(\ast\xi)$.
Thus, we obtain Theorem \ref{thm;21.7.2.20}.
\hfill\qed

\subsubsection{Proof of Proposition \ref{prop;21.7.2.30}}
\label{subsection;21.7.25.11}

We use an induction on the dimension of the support
as in \S\ref{subsection;21.7.25.10}.
It is enough to prove the claim locally around any point of $X$.
We may assume the existence of 
$f$, $Z$, $\varphi$ and $\nbigv$
as in the proof of Proposition \ref{prop;21.2.18.1}.
By using the standard argument,
it is enough to prove
$\Upsilon(\nbigv_0)\simeq
\lefttop{\xi}\Upsilon(\nbigv)$.
There is the isomorphism of the restrictions
$\Upsilon(\nbigv_0)_{|\lefttop{\xi}\gbigz\setminus
\lefttop{\xi}\gbigz_0^{\infty}}
\simeq
\lefttop{\xi}\Upsilon(\nbigv)_{|\lefttop{\xi}\gbigz\setminus
\lefttop{\xi}\gbigz_0^{\infty}}$.
By the Hartogs theorem, it extends to
the desired isomorphism.
Thus, we obtain Proposition \ref{prop;21.7.2.30}.
\hfill\qed

\subsubsection{Appendix: $S^1$-homogeneous harmonic bundles on $\cnum^{\ast}$}
\label{subsection;21.7.2.11}

We refine some statements in \S\ref{subsection;22.7.5.2}.
We set $Y=\proj^1_{\xi}$
and $Y^{\circ}=Y\setminus\{0,\infty\}$.
We consider the $S^1$-action on $Y$
given by $a\bullet\xi=a^{-1}\xi$.

Let $(E,\delbar_E,\theta,h)$ be
a harmonic bundle on $Y^{\circ}$
which is $S^1$-homogeneous in the sense that
(i) $(E,\delbar_E)$ is $S^{1}$-equivariant,
(ii) $a^{\ast}\theta=a\cdot\theta$ for any $a\in S^1$,
(iii) $h$ is $S^1$-invariant.
Lemma \ref{lem;21.7.2.1} holds by
the conditions (i) and (ii).
Let $\nbige$ be the associated
$\nbigr_{Y^{\circ}}$-module.
According to \cite{Mochizuki-wild},
it extends to the $\nbigr_{Y}(\ast\{0,\infty\})$-module
$\nbigq\nbige$,
and
there exists a polarizable pure twistor $\nbigd$-module
$\nbigt=(\nbigm,\nbigm,C)$ of weight $0$
with the polarization $(\id_{\nbigm},\id_{\nbigm})$
on $Y$
such that
$\nbigm(\ast\{0,\infty\})=\nbigq\nbige$.
Moreover,
if there exists a polarizable pure twistor $\nbigd$-module
$\nbigt_1=(\nbigm_1,\nbigm_1,C_1)$
of weight $0$ with the polarization $(\id_{\nbigm_1},\id_{\nbigm_1})$
on $Y\setminus\{\infty\}$
such that
$\nbigm_{1|Y\setminus\{0,\infty\}}=\nbige$,
then we obtain
$\nbigm_1(\ast\{0\})=\nbigq\nbige_{|Y\setminus\{\infty\}}$.
Note that
$\nbigq\nbige$ is unramifiedly good along $\xi=0$
as $\nbigr_{Y\setminus\{\infty\}}(\ast 0)$-module
whose set of the irregular values is
$\{\alpha \xi^{-1}\,|\,\alpha\in S\}$,
where $S$ is the set in Lemma \ref{lem;21.7.2.1}.

As explained in \cite[\S3.2.1]{Mochizuki-Toda-lattice-II},
$\nbige$ is naturally
an $\nbigrtilde_{Y^{\circ}}$-module,
and $\nbigq\nbige$ is naturally
an $\nbigrtilde_{Y}(\ast\{0,\infty\})$-module.
It extends to
$\gbigrtilde_Y(\ast\{0,\infty\})$-module
$\nbigq\nbigetilde$.
Moreover, it is naturally $\cnum^{\ast}$-homogeneous
with respect to the $\cnum^{\ast}$-action
on $\proj^1\times Y$ given by
$a(\lambda,\xi)=(a\lambda,a^{-1}\xi)$.

For $\lambda\in\proj^1$,
we set $\gbigy^{\lambda}=\{\lambda\}\times Y\subset
\proj^1\times Y=\gbigy$.
Note that we obtain a meromorphic flat bundle
$(\nbigv,\nabla)=(\nbigq\nbigetilde(\ast \gbigy^0),\nabla)$
on $\bigl(\gbigy,
 (\proj^1_{\lambda}\times\{0,\infty\})
 \cup
 \gbigy^0\cup\gbigy^{\infty}
 \bigr)$.
We set $Z=\{(\infty,0),(0,\infty)\}\subset\proj^1\times Y$.
The following proposition is similar to
Lemma \ref{lem;22.7.4.21}.

\begin{prop}
\label{prop;21.7.2.1}
 The meromorphic flat bundle
$(\nbigv,\nabla)
 _{|\gbigy\setminus Z}$
is unramifiedly good.
It is regular along
$\gbigy^{\infty}\cup(\{\infty\}\times\proj^1)$.
The set of the irregular values
along $\gbigy^0\cup(\{0\}\times\proj^1)$
is
$\{\alpha \lambda^{-1}\xi^{-1}\,|\,\alpha\in S\}$. 
\hfill\qed
\end{prop}

\begin{cor}
$\nbige_{|\cnum_{\lambda}\times \{1\}}(\ast\lambda)$
is unramifiedly good
whose set of irregular values is
$\{\alpha \lambda^{-1}\,|\,\alpha\in S\}$. 
Moreover, we have
$\nbigq\nbige
=\lefttop{\xi}\bigl(\nbige_{|\cnum_{\lambda}\times \{1\}}\bigr)$.
\hfill\qed
\end{cor}

We obtain the following corollary.
\begin{cor}
\label{cor;22.7.4.10}
$\lefttop{\xi}(\nbige_{|\cnum_{\lambda}\times\{1\}})$
is an object of
$\nbigc(Y;\{0,\infty\})$.
\hfill\qed
\end{cor}

Let $\nbigm\in\nbigc(\pt)$.
Suppose that
there exists 
an integrable pure twistor $\nbigd$-module
$\nbigt_1=(\nbigm_1,\nbigm_1,C_1)$ of weight $0$ with
the integrable polarization $\nbigs_1=(\id,\id)$
on $Y^{\circ}$
such that 
$\nbigm_{1}=
\lefttop{\xi}\nbigm_{|Y^{\circ}}$.
Let $(E,\delbar_E,\theta,h)$ be a harmonic bundle
corresponding to $(\nbigt_1,\nbigs_1)$.
It is easy to see that
$(E,\delbar_E,\theta)$ satisfies
the conditions (i) and (ii).
\begin{lem}
$h$ is $S^1$-invariant.
\end{lem}
\pf
We have the corresponding
polarized integrable variation of pure twistor structure
$(\nbigv^{\sankaku},S^{\sankaku})$
of weight $0$ on $Y^{\circ}$.
It is $S^1$-equivariant
with respect to the $S^1$-action
on $\proj^1_{\lambda}\times Y^{\circ}$
given by $a\bullet(\lambda,\xi)=(a\lambda,a^{-1}\xi)$.
For any $a\in S^1$,
we have the isomorphism
$a^{\ast}:
(\nbigv^{\sankaku},S^{\sankaku})_{|\proj^1\times\{\xi\}}
\simeq
(\nbigv^{\sankaku},S^{\sankaku})_{|\proj^1\times \{a\xi\}}$.
Note that $E_{|\xi}$ is identified with
$H^0(\proj^1,\nbigv_{|\proj^1\times\{\xi\}})$,
and that $h_{|\xi}$ is equal
to the Hermitian pairing induced by
$S^{\sankaku}_{|\proj^1\times\{\xi\}}$.
Hence, we obtain the $S^1$-invariance
of the harmonic metrics.
\hfill\qed

\subsection{Basic functorial properties}

The following proposition is similar to 
Proposition \ref{prop;21.6.29.21}.
\begin{prop}
Let $F:X\lrarr Y$ be a projective morphism.
Let $H_Y$ be a hypersurface of $Y$,
and we put $H_X:=F^{-1}(H_Y)$.
\begin{itemize}
 \item  The functors $F^j_{\dagger}:\gbigc(X;H_X)\lrarr\gbigc(Y;H_Y)$
	induce
	$F_{\dagger}^j:\gbigc_{\res,\infty}(X;H_X)
	\lrarr\gbigc_{\res,\infty}(Y;H_Y)$.  
 \item Suppose moreover that $F_{|X\setminus H_X}$ is a closed embedding of
       $X\setminus H_X$ into $Y\setminus H_Y$.
       We set
       $\gbigc_{F(X),\res,\infty}(Y)=
       \gbigc_{\res,\infty}(Y)\cap\gbigc_{F(X)}(Y)$.
       Then,
       $F_{\dagger}^0$ induces
       an equivalence
       $\gbigc_{\res,\infty}(X;H_X)\simeq
       \gbigc_{F(X),\res,\infty}(Y;H_Y)$.
       \hfill\qed
\end{itemize}
\end{prop}

The following proposition is similar to
Proposition \ref{prop;21.6.29.24}
and
Theorem \ref{thm;21.3.25.20}.
\begin{prop}
Let $X$ be a complex manifold with a hypersurface $H$.
 Let $\gbigm\in\gbigc_{\res,\infty}(X;H)$.
\begin{itemize}
 \item For any hypersurface $H^{(1)}$ of $X$,
       $\gbigm[\star H^{(1)}]$ $(\star=!,\ast)$ are objects
       of $\gbigc_{\res,\infty}(X;H)$.
 \item Let $g$ be any meromorphic function on $(X,H)$.
       Then, $\Pi^{a,b}_{g\star}(\gbigm)$,
       $\Pi^{a,b}_{g,\ast!}(\gbigm)$,
       $\Xi^{(a)}_g(\gbigm)$,
       $\psi^{(a)}_g(\gbigm)$
       and $\phi^{(0)}_g(\gbigm)$
       are objects  of $\gbigc_{\res,\infty}(X;H)$.
 \item
There exists a natural isomorphism
$\lambda\cdot \lefttop{\xi}\bigl(
 \DD_X\gbigm\bigr)
\simeq
 \DD_{\lefttop{\xi}X(\ast \xi)}\bigl(
 \lefttop{\xi}\gbigm
 \bigr)$.
In particular, $\DD_X\gbigm$ is an object of $\gbigc_{\res,\infty}(X)$.      
       \hfill\qed
\end{itemize}
\end{prop}

\begin{cor}
In the situation of {\rm\S\ref{subsection;21.6.22.2}}
and {\rm\S\ref{subsection;21.6.22.20}},
for any $\gbigm\in\gbigc_{\res,\infty}(Y;H_Y)$,
$(\lefttop{T}f^{\star})^i(\gbigm)$ $(\star=!,\ast)$
are objects of
$\gbigc_{\res,\infty}(X;H_X)$.
\hfill\qed
\end{cor}

The following proposition is similar to
Theorem \ref{thm;21.3.26.11}.
\begin{prop}
Let $X_i$ $(i=1,2)$ be complex manifolds.
If $\gbigm_i\in\gbigc_{\res,\infty}(X_i)$,
then
$\gbigm_1\boxtimes\gbigm_2\in\gbigc_{\res,\infty}(X_1\times X_2)$.
\hfill\qed
\end{prop}

\begin{prop}
Let $f:X\lrarr Y$ be a morphism of complex manifolds.
Let $H_Y$ be a hypersurface of $Y$,
and we set $H_X=f^{-1}(H_Y)$.
Let $\gbigm\in\gbigc_{\res,\infty}(Y;H_Y)$.
If $f$ is strictly non-characteristic for $\gbigm$,
$f^{\ast}(\gbigm)$ is an object of $\gbigc_{\res,\infty}(X;H_X)$.
\end{prop}
\pf
We set $\nbigm_1:=\gbigm_{|\nbigy\setminus \nbigh_Y}$.
Let $F(\nbigm_1)$ is a good filtration of $\nbigm_1$
as an $\nbigr_{Y\setminus H_Y}$-module.
It induces a good filtration
$F(\lefttop{\xi}\nbigm_1)$
of $\lefttop{\xi}\nbigm_1$
as an $\nbigr_{\lefttop{\xi}Y\setminus
(\lefttop{\xi}H\cup\lefttop{\xi}Y_0)}$-module.
Hence,
the induced morphism
$\lefttop{\xi}f:\lefttop{\xi}X\lrarr \lefttop{\xi}Y$
is strictly non-characteristic for
$\lefttop{\xi}\gbigm\in
\gbigc(\lefttop{\xi}Y;\lefttop{\xi}Y_0\cup \lefttop{\xi}H_Y)$.
Then, we easily obtain the claim of the proposition.
\hfill\qed

\section{Partial Fourier transforms and $\cnum^{\ast}$-homogeneous objects}

\subsection{Preliminary}
\label{subsection;22.7.27.12}

Let $Z$ be a quasi-projective complex manifold.
There exists an algebraic Zariski open embedding $\iota_Z:Z\to \Zbar$
to a projective manifold $\Zbar$
such that $H_{Z,\infty}:=\Zbar\setminus Z$
is a hypersurface.
We set
\[
\gbigctilde(Z\times X):=
\gbigc_{\Malg}(\Zbar\times X;H_{Z,\infty}\times X).
\]
Let $\iota_Z':Z\to\Zbar'$
be an algebraic Zariski open embedding 
to a projective manifold $\Zbar'$
such that $H'_{Z,\infty}:=\Zbar'\setminus Z$
is a hypersurface.
There exists
an algebraic Zariski open embedding 
$\iota''_Z:Z\to \Zbar''$
to a projective manifold $\Zbar''$
with morphisms
$\varphi_1:\Zbar''\to \Zbar$
and
$\varphi_2:\Zbar''\to \Zbar'$
such that
(i) $H''_{Z,\infty}:=\Zbar''\setminus Z$ is a hypersurface,
(ii) $\varphi_i$ induces the identity on $Z$.
We obtain the following equivalences:
\[
\begin{CD}
 \gbigc_{\Malg}(\Zbar\times X;H_{Z,\infty}\times X)
 @<{\varphi_{1\dagger}}<<
 \gbigc_{\Malg}(\Zbar''\times X;H''_{Z}\times X)
 @>{\varphi_{2\dagger}}>>
 \gbigc_{\Malg}(\Zbar'\times X;H'_{Z}\times X).
\end{CD}
\]
We have the inverse $\varphi_i^{\ast}$ of $\varphi_{i\dagger}$.
Hence, $\gbigctilde(Z\times X)$ is independent
of a projective completion.

Let $f:Z_1\to Z_2$ be an algebraic morphism of
complex quasi-projective manifolds.
It extends to a morphism of complex projective manifolds
$\overline{f}:\overline{Z}_1\to \overline{Z}_2$
such that
(i) $Z_i$ are Zariski open subsets of $\overline{Z}_i$,
(ii) $H_{Z_i}:=\overline{Z}_i\setminus Z_i$ are hypersurfaces,
(iii) $\overline{f}_{|Z_1}$ induces $f$.
For any
$\gbigm\in \gbigctilde(Z_1\times X)$,
we set
\[
 \lefttop{T}f_{\star}^j(\gbigm):=
 f^j_{\dagger}
 \bigl(
 \gbigm[\star (H_{Z_1}\times X)]
 \bigr)(\ast (H_{Z_2}\times X))
\in \gbigctilde(Z_2\times X).
\]

\subsection{Partial Fourier transforms}
\label{subsection;21.7.25.2}

Let $\cnum^m_{z_1,\ldots,z_m}$ denote
the space $\cnum^m$ with the standard coordinate system $(z_1,\ldots,z_m)$.
Let $X$ be a complex manifold.
Let $p_{\tau,X}:\cnum^2_{t,\tau}\times X\lrarr
\cnum_{\tau}\times X$
and $p_{\tau,X}:\cnum^2_{t,\tau}\times X\lrarr
\cnum_{\tau}\times X$ denote the projections.
We have objects
$\gbigl(t\tau)$ and $\gbigl(-t\tau)$
of $\gbigctilde(\cnum^2_{t,\tau})$
(see \S\ref{subsection;22.7.27.10}).
The pull back of $\gbigl(\pm t\tau)$ by
$\cnum^2_{t,\tau}\times
X\lrarr\cnum^2_{t,\tau}$
are also denoted by $\gbigl(\pm t\tau)$.

Let $\gbigm\in\gbigctilde(\cnum\times X)$.
We obtain
$p_{t,X}^{\ast}(\gbigm)\otimes\gbigl(t\tau)
\in\gbigctilde(\cnum^2_{t,\tau}\times X)$.
Let $H_{\infty,t,X}:=\{\infty\}\times\cnum_{\tau}\times X$.
We obtain the following morphisms
in $\gbigctilde(\proj^1_t\times\cnum_{\tau}\times X)$:
\begin{equation}
\label{eq;21.7.4.1}
\bigl(
 p_{t,X}^{\ast}(\gbigm)
 \otimes\gbigl(\pm t\tau)
 \bigr)
[!H_{\infty,t,X}]
\lrarr
\bigl(
 p_{t,X}^{\ast}(\gbigm)
 \otimes\gbigl(\pm t\tau)
 \bigr)
[\ast H_{\infty,t,X}].
\end{equation}

\begin{lem}
\label{lem;21.7.5.1}
The morphisms {\rm(\ref{eq;21.7.4.1})}
are isomorphisms.
Moreover, the following natural morphism is an isomorphism:
\begin{equation}
\label{eq;21.7.4.3}
 \bigl(
 p_{t,X}^{\ast}(\gbigm)
 \otimes\gbigl(\pm t\tau)
 \bigr)
[\star H_{\infty,t,X}]
\lrarr
  p_{t,X}^{\ast}(\gbigm)
 \otimes\gbigl(\pm t\tau).
\end{equation}
\end{lem}
\pf
We explain the proof in the case of
$\gbigl(t\tau)$.
There exist $V$-filtrations
$V_a\Bigl(
\bigl(
 p_{t,X}^{\ast}(\gbigm)
 \otimes\gbigl(t\tau)
 \bigr)
[\star H_{\infty,t,X}]
 \Bigr)$
along $t^{-1}$. 
For $a<0$, we have
\[
 V_a\Bigl(
\bigl(
 p_{t,X}^{\ast}(\gbigm)
 \otimes\gbigl(t\tau)
 \bigr)
[!H_{\infty,t,X}]
 \Bigr)
=  V_a\Bigl(
\bigl(
 p_{t,X}^{\ast}(\gbigm)
 \otimes\gbigl(t\tau)
 \bigr)
[\ast H_{\infty,t,X}]
 \Bigr)
\subset
 p_{t,X}^{\ast}(\gbigm)
 \otimes\gbigl(t\tau).
\]
Note that the following holds for any $a\in\real$:
\[
\deldel_{\tau}\cdot
 V_a\Bigl(
\bigl(
 p_{t,X}^{\ast}(\gbigm)
 \otimes\gbigl(t\tau)
 \bigr)
[\star H_{\infty,t,X}]
 \Bigr) 
 \subset
  V_a\Bigl(
\bigl(
 p_{t,X}^{\ast}(\gbigm)
 \otimes\gbigl(t\tau)
 \bigr)
[\star H_{\infty,t,X}]
 \Bigr).
\]

Let $U_X$ be an open subset of $X$.
Let $U_{t,\infty}$ be an open neighbourhood of $\infty$
in $\proj^1_{\infty}$.
Let $s$ be a section of $\gbigm$ on $U_{t,\infty}\times U_X$.
There exists $N\in\seisuu_{>0}$ such that
$t^{-N}p_{t,X}^{\ast}(s)$
is a section of
$V_{<0}\Bigl(
\bigl(
 p_{t,X}^{\ast}(\gbigm)
 \otimes\gbigl(t\tau)
 \bigr)
 [\star H_{\infty,t,X}]
\Bigr)$.
Hence,  for any $m\in\seisuu_{\geq 0}$,
$t^mp_{t,X}^{\ast}(s)
=\deldel_{\tau}^{N+m}(t^{-N}p_{t,X}^{\ast}(s))$
is a section of
$V_{<0}\Bigl(
\bigl(
 p_{t,X}^{\ast}(\gbigm)
 \otimes\gbigl(t\tau)
 \bigr)
 [\star H_{\infty,t,X}]
\Bigr)$.
We obtain
$p_{t,X}^{\ast}(\gbigm)
\otimes\gbigl(t\tau)
=V_{<0}\Bigl(
\bigl(
p_{t,X}^{\ast}(\gbigm)
\otimes\gbigl(t\tau)
\bigr)
[\star H_{\infty,t,X}]
\Bigr)$.
Then, the claim of Lemma  \ref{lem;21.7.5.1} follows.
\hfill\qed

\begin{lem}
We have
$\lefttop{T}(p_{\tau,X})_{\star}^i
 \bigl(
 p_{t,X}^{\ast}(\gbigm)
 \otimes\gbigl(\pm t\tau)
 \bigr)
 =0$
for $i\neq 0$.
Moreover,
the natural morphism
\[
 \lefttop{T}(p_{\tau,X})_{!}^i
 \bigl(
 p_{t,X}^{\ast}(\gbigm)
 \otimes\gbigl(\pm t\tau)
 \bigr)
 \lrarr
  \lefttop{T}(p_{\tau,X})_{\ast}^i
 \bigl(
 p_{t,X}^{\ast}(\gbigm)
 \otimes\gbigl(\pm t\tau)
\bigr)
\]
is an isomorphism.
\end{lem}
\pf
By (\ref{eq;21.7.4.3}),
we obtain
the vanishing in the case $i\geq 1$.
By the duality, we obtain the vanishing in the case $i\leq -1$.
Because (\ref{eq;21.7.4.1}) is an isomorphism,
we obtain the second claim.
\hfill\qed

\vspace{.1in}
For any $\gbigm\in\gbigctilde(\cnum_t\times X)$,
we define the objects $\FT_+(\gbigm)$ and $\FT_-(\gbigm)$
of $\gbigctilde(\cnum_{\tau}\times X)$
as follows:
\[
 \FT_{\pm}\bigl(\gbigm\bigr)
 =
 (p_{\tau,X})_{\dagger}^0\bigl(
  p_{t,X}^{\ast}(\gbigm)
 \otimes\gbigl(\pm t\tau)
 \bigr)
 \simeq
 \lefttop{T}
 (p_{\tau,X})_{\star}^0\bigl(
  p_{t,X}^{\ast}(\gbigm)
 \otimes\gbigl(\pm t\tau)
 \bigr)\quad(\star=!,\ast)
\]

\subsubsection{Partially algebraic description and the inversion}

Let $\pi_{\lambda,t}:\proj^1_{\lambda}\times\cnum_{t}\times X\lrarr X$
denote the projection.
We set
\[
\gbiga_{X,t}:=(\pi_{\lambda,t})_{\ast}\Bigl(
\gbigrtilde_{\proj^1_t\times X}
\bigl(\ast(\{\infty\}\times X)\bigr)
\Bigr)
\subset
\nbigd_X[\lambda,t]
\langle
\del_{\lambda},\del_{t}
\rangle.
\]
It is the sheaf of subalgebras
of $\nbigd_X[\lambda,t]
\langle
\del_{\lambda},\del_{t}
\rangle$
generated by
$\lambda\Theta_X$,
$\lambda^2\del_{\lambda}$
and $\deldel_t=\lambda\del_t$
over $\nbigo_X[\lambda,t]$.

Let $\gbign$ be an $\gbiga_{X,t}$-module.
We define the $\gbiga_{X,\tau}$-modules
$\gbigf_{+}(\gbign)$ and $\gbigf_-(\gbign)$ as follows.
\begin{itemize}
 \item We set
       $\gbigf_{\pm}(\gbign):=\lambda^{-1}\gbign$
       as $\nbigo_X[\lambda]$-modules.
 \item  For sections $s\in\gbigf_{\pm}(\gbign)$,
	we define       
       $\deldel_{\tau}\bullet s=\pm ts$,
	$\tau\bullet s=\mp \deldel_{t}s$
	and
       $\lambda^2\del_{\lambda}\bullet s:=
       \lambda^2\del_{\lambda}(s)
       \pm\lambda\del_t(ts)$.
\end{itemize}

For any $\gbigm\in\gbigctilde(\cnum_t\times X)$,
we obtain
the $\gbiga_{X,t}$-module
$(\pi_{\lambda,t})_{\ast}(\gbigm)$.
We can reconstruct $\gbigm$ from
$(\pi_{\lambda,t})_{\ast}(\gbigm)$
as the analytification.
By Lemma \ref{lem;21.7.5.1} and the standard computation,
we obtain the following lemma.

\begin{lem}
For $\gbigm\in\gbigctilde(\cnum_t\times X)$,
we have
$(\pi_{\lambda,\tau})_{\ast}\FT_{\pm}(\gbigm)
 \simeq
 \gbigf_{\pm}\bigl(
 (\pi_{\lambda,t})_{\ast}(\gbigm)
  \bigr)$.
\hfill\qed
\end{lem}

\begin{prop}
For $\gbigm\in\gbigctilde(\cnum\times X)$,
there exists the natural isomorphism
$\FT_{\mp}\circ\FT_{\pm}(\gbigm)
\simeq\lambda^{-1}\gbigm$.
\end{prop}
\pf
We have the natural isomorphism
$\gbigf_{\mp}\bigl(
\gbigf_{\pm}\bigl((\pi_{\lambda,t})_{\ast}(\gbigm)\bigr)
\bigr)
\simeq
\lambda^{-2}(\pi_{\lambda,t})_{\ast}\gbigm$
as $\nbigo_X[\lambda]$-modules.
The naturally induced actions of $\deldel_t$ and $t$ are the same
on the both sides.
We note that
$\lambda^2\del_{\lambda}+\lambda\del_tt
 -\lambda t\del_t
 =\lambda^2\del_{\lambda}+\lambda$.
Hence,
$\lambda^{-2}\gbigm\simeq \lambda^{-1}\gbigm$
given by $\lambda^{-2}(s)\longmapsto \lambda^{-1}(s)$
induces an isomorphism of
$\gbiga_{X,t}$-modules
$\gbigf_{\mp}(\gbigf_{\pm}(\gbigm))\simeq \lambda^{-1}\gbigm$.
\hfill\qed

\subsubsection{Equivalences of some subcategories}

Let $\star\in\{!,\ast\}$.
Let $\gbigctilde(\cnum\times X,[\star 0])
\subset\gbigctilde(\cnum\times X)$
denote the full subcategories of
$\gbigm\in\gbigctilde(\cnum\times X)$
such that $\gbigm[\star(\{0\}\times X)]=\gbigm$.
Let $\gbigctilde(\cnum\times X)_{\star}
\subset\gbigctilde(\cnum\times X)$
denote the full subcategories of
$\gbigm\in\gbigctilde(\cnum\times X)$
such that
$\lefttop{T}(\pi_X)^j_{\star}\bigl(
\gbigm
\bigr)=0$ $(j\in\seisuu)$,
where $\pi_X:\proj^1_t\times X\lrarr X$ denote the projection.
The following lemma is clear.
\begin{lem}
Let $(\star_1,\star_2)=(!,\ast),(\ast,!)$.
We have
$\DD\gbigctilde(\cnum\times X,[\star_1 0])
 \subset \gbigctilde(\cnum\times X,[\star_2 0])$
and  
 $\DD\gbigctilde(\cnum\times X)_{\star_1}
 \subset \gbigctilde(\cnum\times X)_{\star_2}$.
\hfill\qed
\end{lem}

\begin{prop}
\label{prop;21.7.11.10}
$\FT_{\pm}$ induce equivalences
$\FT_{\pm}:
 \gbigctilde(\cnum_t\times X,[\star 0])
 \simeq
 \gbigctilde(\cnum_{\tau}\times X)_{\star}$.
\end{prop}
\pf
We explain the proof in the case $\star=\ast$ and $\FT_+$.
The other case can be argued similarly.
Let $\gbigm\in\gbigctilde(\cnum_t\times X)$.
Let $p_X:\cnum^2_{t,\tau}\times X\lrarr X$
denote the projection.
We obtain
\begin{equation}
\label{eq;21.7.12.1}
\lefttop{T}(\pi_{X})_{\ast}^j\Bigl(
 \FT_{+}(\gbigm)
\Bigr)
=\lefttop{T}(p_{X})_{\ast}^j\Bigl(
 p_{t,X}^{\ast}(\gbigm)
 \otimes\gbigl(t\tau)
 \Bigr).
\end{equation}

Let $\Delta:\cnum^2_{t,\tau}\times X
\lrarr \cnum^3_{t_1,t_2,\tau}\times X$
be the morphism defined by
$\Delta(t,\tau,x)=(t,t,\tau,x)$.
Let
$\Htilde_{\Delta}\subset
\cnum^3_{t_1,t_2,\tau}\times X$
denote the image of $\Delta$.
Let
$p_i:\cnum^3_{t_1,t_2,\tau}\times X
\lrarr\cnum^2_{t,\tau}\times X$
denote the morphism
defined by
$p_i(t_1,t_2,\tau,x)=(t_i,\tau,x)$.
We consider the following complex
in
$\gbigctilde(\cnum^3_{t_1,t_2,\tau}\times X)$:
\begin{equation}
\label{eq;21.7.5.2}
\Bigl(
 (p_{t,X}\circ p_1)^{\ast}(\gbigm)
 \otimes
 p_2^{\ast}\gbigl(t\tau) \Bigr)
 \lrarr
 \Bigl(
 (p_{t,X}\circ p_1)^{\ast}(\gbigm)
 \otimes
 p_2^{\ast}\gbigl(t\tau)
  \Bigr)[\ast\Htilde_{\Delta}],
\end{equation}
where the second term sits in the degree $0$.
By Proposition \ref{prop;21.6.22.12},
the complex is quasi-isomorphic to
\[
 \Delta_{\dagger}
\Bigl(
  \lambda^{-1} p_{t,X}^{\ast}(\gbigm)
  \otimes\gbigl(t\tau)
\Bigr).
\]
Let
$p_{1,2}:\cnum^3_{t_1,t_2,\tau}\times X
\lrarr\cnum^2_{t_1,t_2}\times X$
denote the projection.
Let $\iota_{0,2}:
\cnum_t\times X\lrarr
\cnum^2_{t_1,t_2}\times X$
denote the inclusion defined by
$\iota_{0,2}(t,x)=(t,0,x)$.
By Lemma \ref{lem;22.7.5.10} below,
we have
\[
\lefttop{T}
(p_{1,2})^j_{\ast}\Bigl(
(p_{t,X}\circ p_1)^{\ast}(\gbigm)
 \otimes
 p_{2}^{\ast}\gbigl(t\tau)
 \Bigr)
 =\left\{
   \begin{array}{ll}
 \iota_{0,2\dagger} \lambda^{-1} p_{t,X}^{\ast}(\gbigm)
  & (j=0)\\
 0 & (j\neq 0).
   \end{array}
 \right.
\]
Hence, we obtain the following complex
in $\gbigctilde(\cnum^2_{t_1,t_2}\times X)$
by applying $\lefttop{T}(p_{1,2})_{\ast}$
to the complex (\ref{eq;21.7.5.2}):
\begin{equation}
\label{eq;21.7.5.3}
 \iota_{0,2\dagger} (\lambda^{-1} \gbigm)
\lrarr
\iota_{0,2\dagger} (\lambda^{-1} \gbigm)
 [\ast H_{\Delta}],
\end{equation}
where
$H_{\Delta}$ denotes the image of
$\cnum_t\times X\lrarr\cnum^2_{t_1,t_2}\times X$
defined by $\Delta(t,x)=(t,t,x)$.
By \cite[Proposition 11.2.7]{Mochizuki-MTM},
the complex (\ref{eq;21.7.5.3}) is identified with the following complex:
\begin{equation}
 \iota_{0,2\dagger}(\lambda^{-1}\gbigm)
  \lrarr
 \iota_{0,2\dagger}\bigl(\lambda^{-1}\gbigm[\ast (\{0\}\times X)]\bigr).
\end{equation}
Note that the support of the kernel and the cokernel
are contained in
$\{(0,0)\}\times X$.
Therefore,
(\ref{eq;21.7.12.1})
is naturally identified with
the $j$-th cohomology of the following complex
in $\gbigctilde(\cnum_t\times X)$:
\[
 \lambda^{-1}\gbigm
 \lrarr\lambda^{-1}\gbigm[\ast (\{0\}\times X)].
\]
Hence,
$\gbigm\in \gbigctilde(\cnum\times X,[\ast 0])$
if and only if
$\FT_+(\gbigm)\in\gbigctilde(\cnum\times X)_{\ast}$.
\hfill\qed

\begin{lem}
\label{lem;22.7.5.10}
For $\gbign\in \gbigc_{\Malg}(X)$,
 we have
$\lefttop{T}(p_{t,X})_{\ast}^i
 \Bigl(
 \gbigl(t\tau)\boxtimes\gbign
 \Bigr)=0$ $(i\neq 0)$,
and 
$\lefttop{T}(p_{t,X})^0_{\ast}
 \Bigl(
 \gbigl(t\tau)\boxtimes\gbign
 \Bigr)
 \simeq
 \iota_{0\dagger}(\lambda^{-1}\gbign)$,
where $\iota_0:X\to\cnum_t\times X$
is given by $\iota_0(x)=(0,x)$.
\end{lem}
\pf
Though this is essentially contained in
\cite[Proposition 3.22]{Mochizuki-GKZ},
we give an indication of the proof.
(See also Lemma \ref{lem;22.7.26.11} below.)
It is easy to see that
$\lefttop{T}(p_{t,X})^i_{\ast}
\bigl(
\gbigl(t\tau)\boxtimes\gbign
\bigr)=0$ $(i\neq 0)$
and the support of 
$\lefttop{T}(p_{t,X})^0_{\ast}
 \bigl(
 \gbigl(t\tau)\boxtimes\gbign
 \bigr)$
 is contained in $\{0\}\times X$.
By the partial algebraic description
we may identify
$\lefttop{T}(p_{t,X})^0_{\ast}
 \bigl(
 \gbigl(t\tau)\boxtimes\gbign
 \bigr)$
as the cokernel of 
$\del_t+\lambda^{-1}\tau:
\gbign[t,\tau]\to\gbign[t,\tau]\lambda^{-1}$,
and the action of
$\deldel_{\tau}=\lambda\del_{\tau}$
and $\lambda^{2}\del_{\lambda}$
are given as follows:
\[
 \deldel_{\tau}\bullet(t^jm\lambda^{-1})
 =t^{j+1}m\lambda^{-1},
 \quad
 \lambda^2\del_{\lambda}\bullet
 \bigl(t^jm\lambda^{-1}\bigr)
 =\lambda(t^j\lambda\del_{\lambda}(m)\lambda^{-1}
 +\lambda jt^jm\lambda^{-1}).
\]
Hence, we can identify 
$\lefttop{T}(p_{t,X})^0_{\ast}
 \bigl(
 \gbigl(t\tau)\boxtimes\gbign
 \bigr)$
with $\gbign[\deldel_{\tau}]$
by 
$t^jm\lambda^{-1}\longmapsto \deldel_{\tau}^jm$.
It is identified with
$\iota_{0\dagger}(\lambda^{-1}\gbign)
\simeq
\iota_{0\ast}(\lambda^{-1}\gbign (d\tau/\lambda)^{-1})
\otimes\cnum[\deldel_{\tau}]$.
\hfill\qed

\subsubsection{Basic functoriality}

\begin{prop}
\label{prop;21.7.14.10}
Let $F:X\lrarr Y$ be a projective morphism of complex manifolds.
For any $\gbigm\in\gbigctilde(\cnum\times X)$,
there exist natural isomorphisms
\begin{equation}
\label{eq;21.7.14.1}
\FT_{\pm}(\id_{\cnum}\times F)^j_{\dagger}(\gbigm)\simeq
(\id_{\cnum}\times F)^j_{\dagger}\FT_{\pm}(\gbigm).
\end{equation}
\end{prop}
\pf
It is easy to see
$(\id_{\cnum^2_{t,\tau}}\times F)_{\dagger}^j
(p_{t,X}^{\ast}\gbigm)
 \simeq
 p_{t,Y}^{\ast}\bigl(
 (\id_{\cnum_t}\times F)_{\dagger}^j(\gbigm)
 \bigr)$.
Then, we obtain
\[
 (\id_{\cnum^2_{t,\tau}}\times F)_{\dagger}^j
  (p_{t,X}^{\ast}\gbigm\otimes\gbigl(\pm t\tau))
 \simeq
 p_{t,Y}^{\ast}\bigl(
 (\id_{\cnum_t}\times F)_{\dagger}^j(\gbigm)
 \bigr)\otimes\gbigl(\pm t\tau).
\] 
(See \cite[Lemma 2.1.16]{Mochizuki-MTM}, for example.)
Note that
$p_{\tau,Y}\circ (\id_{\cnum^2_{t,\tau}}\times F)
=(\id_{\cnum_{\tau}}\times F)\circ p_{\tau, X}$.
We also have
$(p_{\tau,Y})_{\dagger}^k\bigl(
(\id_{\cnum^2_{t,\tau}}\times F)_{\dagger}^j
(p_{t,X}^{\ast}\gbigm\otimes\gbigl(\pm t\tau))
\bigr)=0$ $(k\neq 0)$
and
$(\id_{\cnum_{\tau}}\times F)^j_{\dagger}
\bigl(
 (p_{\tau,X})_{\dagger}^k
 (p_{t,X}^{\ast}\gbigm\otimes\gbigl(\pm t\tau))
\bigr)=0$ $(k\neq 0)$.
Hence, we obtain (\ref{eq;21.7.14.1}).
\hfill\qed

\begin{cor}
\label{cor;21.7.14.11}
A projective morphism of complex manifolds
$F:X\lrarr Y$ induces
$(\id\times F)_{\dagger}^j:\gbigctilde(\cnum\times X,[\star 0])\lrarr
\gbigctilde(\cnum\times Y,[\star 0])$
and 
$(\id\times F)_{\dagger}^j:\gbigctilde(\cnum\times X)_{\star}\lrarr
\gbigctilde(\cnum\times Y)_{\star}$
for $\star=!,\ast$.
\end{cor}
\pf
Because
$(\id\times F)^j_{\dagger}(\gbigm[\star (0\times X)])
=(\id\times F)^j_{\dagger}(\gbigm)[\star (0\times Y)]$,
we obtain
$(\id\times F)^j_{\dagger}:
\gbigctilde(\cnum\times X,[\star 0])
\lrarr
\gbigctilde(\cnum\times Y,[\star 0])$.
By using (\ref{eq;21.7.14.1}),
we obtain
$(\id\times F)^j_{\dagger}:
\gbigctilde(\cnum\times X)_{\star}
\lrarr
\gbigctilde(\cnum\times Y)_{\star}$.
\hfill\qed

\begin{prop}
\label{prop;21.7.13.20}
Let $\gbigm\in\gbigctilde(\cnum\times X)$.
\begin{itemize}
 \item For any hypersurface $H$ of $X$,
       there exist natural isomorphisms
\[
       \FT_{\pm}(\gbigm[\star (\cnum\times H)])
       \simeq
       \FT_{\pm}(\gbigm)[\star (\cnum\times H)].
\]

 \item Let $f$ be any holomorphic function on $X$.
       Let $\ftilde$ denote the induced holomorphic function
       on $\proj^1\times X$.
       There exist natural isomorphisms
       $\FT_{\pm}\Pi^{a,b}_{\ftilde,\star}\bigl(
       \gbigm
       \bigr)
       \simeq
       \Pi^{a,b}_{\ftilde,\star}\bigl(
       \FT_{\pm}\gbigm
       \bigr)$ $(\star=!,\ast)$
       and 
       $\FT_{\pm}\Pi^{a,b}_{\ftilde,\ast !}\bigl(
       \gbigm
       \bigr)
       \simeq
       \Pi^{a,b}_{\ftilde,\ast !}\bigl(
       \FT_{\pm}\gbigm
       \bigr)$.
       In particular,
       there exist natural isomorphisms
       $\FT_{\pm}\Xi^{(a)}_{\ftilde}(\gbigm)
       \simeq
       \Xi^{(a)}_{\ftilde}(\FT_{\pm}(\gbigm))$
       and 
       $\FT_{\pm}\psi^{(a)}_{\ftilde}(\gbigm)
       \simeq
       \psi^{(a)}_{\ftilde}(\FT_{\pm}(\gbigm))$.
       We also obtain
       $\FT_{\pm}\phi^{(0)}_{\ftilde}(\gbigm)
       \simeq
       \phi^{(0)}_{\ftilde}\FT_{\pm}(\gbigm)$.
\end{itemize}
\end{prop}
\pf
We obtain the first claim from the following isomorphisms:
\begin{multline}
\lefttop{T}(p_{\tau,X})^0_{\star}
 \Bigl(\bigl(
 p_{t,X}^{\ast}(\gbigm[\star(\cnum_t^1\times H)])
 \otimes\gbigl(\pm t\tau)
 \bigr)
 \Bigr)
 \simeq
\lefttop{T}(p_{\tau,X})_{\star}^0
 \Bigl(\bigl(
 p_{t,X}^{\ast}(\gbigm)
 \otimes\gbigl(\pm t\tau)
 \bigr)
 \bigl[
 \star(\cnum^2_{t,\tau}\times H)
 \bigr]
 \Bigr) \\
 \simeq
 \lefttop{T}(p_{\tau,X})_{\star}^0
 \bigl(
 p_{t,X}^{\ast}(\gbigm)
 \otimes\gbigl(\pm t\tau)
 \bigr)
 [\star(\cnum_{\tau}\times H)].
\end{multline}
We obtain 
$\FT_{\pm}\Pi^{a,b}_{\ftilde,\star}(\gbigm)
\simeq
\Pi^{a,b}_{\ftilde,\star}\FT_{\pm}(\gbigm)$
similarly.
For any $P\in X$,
there exist a sufficiently large number $N(P)$
and a neighbourhood $X_P$  of $P$ in $X$,
such that 
$\Pi^{a,b}_{\ftilde,\ast !}(\gbigm)$ on $\proj^1_t\times X_P$
is canonically isomorphic to
the cokernel of
$\Pi^{b,N(P)}_{\ftilde,!}(\gbigm)
\lrarr
\Pi^{a,N(P)}_{\ftilde,\ast}(\gbigm)$.
It is easy to see that
$\FT_{\pm}\Pi^{a,b}_{\ftilde,\ast !}(\gbigm)$ on $\proj^1_{\tau}\times X_P$
are isomorphic to
the cokernel of
$\FT_{\pm}\Pi^{b,N(P)}_{\ftilde,!}(\gbigm)
\lrarr
\FT_{\pm}\Pi^{a,N(P)}_{\ftilde,\ast}(\gbigm)$,
which is naturally identified with the cokernel of
the morphism
$\Pi^{b,N(P)}_{\ftilde,!}(\FT_{\pm}\gbigm)
\lrarr
\Pi^{a,N(P)}_{\ftilde,\ast}(\FT_{\pm}\gbigm)$,
Hence, we obtain
$\FT_{\pm}\Pi^{a,b}_{\ftilde,\ast !}(\gbigm)
\simeq
\Pi^{a,b}_{\ftilde,\ast !}(\FT_{\pm}\gbigm)$.
In particular, we obtain the isomorphisms
for $\psi_{\ftilde}^{(a)}$ and $\Xi^{(a)}_{\ftilde}$.
Because $\phi^{(0)}_{\ftilde}(\gbigm)$
is the cohomology of the complex
\[
\gbigm[!\ftilde]
\lrarr
\Xi^{(0)}_{\ftilde}(\gbigm)
\oplus
\gbigm
\lrarr
\gbigm[\ast\ftilde],
\]
we obtain $\FT_{\pm}(\gbigm)$
as the cohomology of the complex
\[
\FT_{\pm}(\gbigm[!\ftilde])
\lrarr
\FT_{\pm}(\Xi^{(0)}_{\ftilde}(\gbigm))
\oplus
\FT_{\pm}(\gbigm)
\lrarr
\FT_{\pm}(\gbigm[\ast\ftilde]).
\]
It is naturally identified with
\[
 (\FT_{\pm}\gbigm)[!\ftilde]
\lrarr
\Xi^{(0)}_{\ftilde}(\FT_{\pm}\gbigm)
\oplus
\FT_{\pm}(\gbigm)
\lrarr
(\FT_{\pm}\gbigm)[\ast\ftilde].
\]
Hence, we obtain the isomorphism for $\phi_{\ftilde}^{(0)}$.
\hfill\qed

\begin{cor}
In the situation of {\rm\S\ref{subsection;21.6.22.2}}
and {\rm\S\ref{subsection;21.6.22.20}},
for any $\gbigm\in\gbigctilde(\cnum\times Y)$,
there exist natural isomorphisms 
$\FT_{\pm}(\lefttop{T}(\id\times f)^{\star})^i(\gbigm)
 \simeq
 (\lefttop{T}(\id\times f)^{\star})^i\FT_{\pm}(\gbigm)$.
\hfill\qed
\end{cor}

We have the following standard compatibility
of the duality and the partial Fourier transform.
\begin{prop}
For $\gbigm\in\gbigctilde(\cnum\times X)$,
there exist natural isomorphisms
$\lambda \FT_{\pm}(\DD\gbigm)
\simeq
 \DD(\FT_{\mp}(\gbigm))$.
\end{prop}
\pf
Because
$\bigl(
 \lambda p_{t,X}^{\ast}(\DD \gbigm)
 \otimes\gbigl(\mp t\tau)
 \bigr)
 [!H_{\infty,t,X}]
\simeq
\DD\Bigl(
 \bigl(
 p_{t,X}^{\ast}(\gbigm)
 \otimes
 \gbigl(\pm t\tau)
 \bigr)[\ast H_{\infty,t,X}]
 \Bigr)$
in $\gbigctilde(\proj^1_t\times\cnum_{\tau}\times X)$
as in \cite[Proposition 13.3.5, Proposition 13.3.6,
Proposition 13.3.9]{Mochizuki-MTM},
we obtain the claim of the proposition.
\hfill\qed

\subsubsection{Convolution and tensor product}
\label{subsection;21.7.14.30}

Let $\Delta:\cnum_{t}\lrarr \cnum^2_{t_1,t_2}$
denote the diagonal embedding.
Let $\mu:\cnum^2_{t_1,t_2}\to \cnum$
be the morphism defined by
$\mu(t_1,t_2)=t_1+t_2$.
Let $X_i$ be complex manifolds.
We set $X=X_1\times X_2$.
For $\gbigm_i\in \gbigctilde(\cnum_{t_i}\times X_i)$ $(i=1,2)$,
we obtain
$\gbigm_1\boxtimes\gbigm_2
\in\gbigctilde(\cnum^2_{t_1,t_2}\times X)$.
We define
\[
\nbigh^j\Bigl(
 \gbigm_1\langle\boxtimes,\otimes\rangle^{\star}
 \gbigm_2
 \Bigr)
 =
 \bigl(
 \lefttop{T}(\Delta\times\id_{X})^{\star}
 \bigr)^j
 (\gbigm_1\boxtimes\gbigm_2)
\in\gbigctilde(\cnum\times X),
\]
\[
 \nbigh^j\Bigl(
 \gbigm_1\langle\boxtimes,\ast\rangle_{\star}
 \gbigm_2
 \Bigr)
 =
 \lefttop{T}(\mu_{\dagger})_{\star}^j (\gbigm_1\boxtimes\gbigm_2)
\in\gbigctilde(\cnum\times X).
\]
It is standard that
the partial Fourier transforms
exchange the convolution product and the tensor product
as in the following proposition.
\begin{prop}
\label{prop;21.7.13.10}
There exist the following natural isomorphisms.
\begin{equation}
\label{eq;21.7.13.1}
 \nbigh^j\Bigl(
 \FT_{\pm}(\gbigm_1)\langle\boxtimes,\otimes\rangle^!
 \FT_{\pm}(\gbigm_2)
 \Bigr)
 \simeq
 \lambda^{-1}\FT_{\pm}
\Bigl(
 \nbigh^{j-1}\bigl(
 \gbigm_1\langle\boxtimes,\ast\rangle_{\ast}
 \gbigm_2
 \bigr)
 \Bigr),
\end{equation}
\begin{equation}
\label{eq;21.7.13.2}
 \nbigh^j\Bigl(
 \FT_{\pm}(\gbigm_1)\langle\boxtimes,\otimes\rangle^{\ast}
 \FT_{\pm}(\gbigm_2)
 \Bigr)
 \simeq
 \FT_{\pm}
 \Bigl(
 \nbigh^{j+1} 
 \bigl(
 \gbigm_1\langle\boxtimes,\ast\rangle_{!}
 \gbigm_2
 \bigr)
 \Bigr).
\end{equation}
\end{prop}
\pf
We study (\ref{eq;21.7.13.1})
for $\FT_+$.
We indicate only an outline.
There exists the following fiber square:
\begin{equation}
  \begin{CD}
  \cnum^3_{t_1,t_2,\tau}\times X
  \times X
   @>{\Delta_1}>>
  \cnum^4_{t_1,t_2,\tau_1,\tau_2}\times X
  \\
  @V{p_1}VV @V{p_2}VV\\
\cnum_{\tau}\times X
  @>{\Delta_2}>>
\cnum^2_{\tau_1,\tau_2}\times X.
  \end{CD}
\end{equation}
Here, $\Delta_i$
are induced by the diagonal embedding
$\cnum_{\tau}\to\cnum^2_{\tau_1,\tau_2}$,
and $p_i$ denote the projections.
Let $p_{t_i,X_i}:\cnum^4_{t_1,t_2,\tau_1,\tau_2}
\lrarr \cnum_{t_i}\times X_i$
denote the projections.
We set
\[ 
\gbigmtilde:=
 p_{t_1,X_1}^{\ast}(\gbigm_1)
 \otimes
 p_{t_2,X_2}^{\ast}(\gbigm_2)
 \otimes
 \gbigl\bigl(t_1\tau_1+t_2\tau_2\bigr)
\in
 \gbigctilde
 \bigl(
 \cnum^4_{t_1,t_2,\tau_1,\tau_2}\times X
  \bigr).
\]
Let $H_{\tau_1=\tau_2}^{(1)}$
denote the hypersurface of
$\cnum^4_{t_1,t_2,\tau_1,\tau_2}\times X$
defined by $\tau_1=\tau_2$.
We obtain the following complex
in
$\gbigctilde(\cnum^4_{t_1,t_2,\tau_1,\tau_2}\times X)$:
\begin{equation}
\label{eq;21.7.13.3}
 \gbigmtilde\lrarr\gbigmtilde[\ast H^{(1)}_{\tau_1=\tau_2}].
\end{equation}
Here, the first term sits in the degree $0$.
Let $H^{(2)}_{\tau_1=\tau_2}$
denote the hypersurface of
$\cnum^2_{\tau_1,\tau_2}\times X$
defined by $\tau_1=\tau_2$.
We obtain the following complex
in $\gbigctilde(\cnum^2_{\tau_1,\tau_2}\times X)$
by applying $\lefttop{T}(p_2)_{\ast}$
to (\ref{eq;21.7.13.3}):
\begin{equation}
\label{eq;21.7.13.4}
 \FT_{+}(\gbigm_1)\boxtimes
 \FT_+(\gbigm_2)
 \lrarr
 \Bigl(
 \FT_{+}(\gbigm_1)\boxtimes
 \FT_+(\gbigm_2)
 \Bigr)[\ast H^{(2)}_{\tau_1=\tau_2}].
\end{equation}
The $j$-th cohomology of
(\ref{eq;21.7.13.4})
is quasi-isomorphic to
$\Delta_{2\dagger}\nbigh^j\bigl(
\FT_{+}(\gbigm_1)\langle\boxtimes,\otimes\rangle^!
\FT_+(\gbigm_2)
\bigr)$.

The projections
$\cnum^3_{t_1,t_2,\tau}\times X
\to\cnum_{t_i}\times X_i$
are also denoted by $p_{t_i,X_i}$.
We obtain
\[
 \gbigmtilde^{(1)}:=
 p_{t_1,X_1}^{\ast}(\gbigm_1)
 \otimes
 p_{t_2,X_2}^{\ast}(\gbigm_2)
 \otimes
 \gbigl((t_1+t_2)\tau)
 \in
 \gbigctilde(\cnum^3_{t_1,t_2,\tau}\times X).
\]
Similarly to Proposition \ref{prop;21.6.22.12},
the complex (\ref{eq;21.7.13.3}) is naturally quasi-isomorphic to
$\Delta_{1\dagger}\bigl(
\lambda^{-1}\gbigmtilde^{(1)}
\bigr)[-1]$.
Hence,
we obtain
\[
 \Delta_{2\dagger}\bigl(
 p_{1\dagger}^{j-1}
 \lambda^{-1}\gbigmtilde^{(1)}
 \bigr)
 \simeq
 \Delta_{2\dagger}
\nbigh^j\bigl(
\FT_{+}(\gbigm_1)\langle\boxtimes,\otimes\rangle^!
\FT_+(\gbigm_2)
\bigr).
\]
It is easy to see that
there exists a natural isomorphism:
\[
 p_{1\dagger}^{j}
 \gbigmtilde^{(1)}
 \simeq
 \FT_{+}\Bigl(
 \nbigh^j\bigl(
 \gbigm_1
 \langle\boxtimes,\ast\rangle_{\ast}
 \gbigm_2
 \bigr)
 \Bigr).
\]
Thus, we obtain
(\ref{eq;21.7.13.1})
for $\FT_+$.
The other cases can be argued similarly.
\hfill\qed

\vspace{.1in}

Let $\nbigo\in\gbigctilde(\cnum)$
denote the object naturally induced by
$\nbigo_{\cnum_{\lambda}\times\cnum}$
with the exterior derivative.
We have
$\nbigo[\star 0]\in \gbigctilde(\cnum)$.
Let $\gbigm\in\gbigctilde(\cnum\times X)$.
Similarly to Corollary \ref{cor;21.6.22.13},
it is easy to obtain
\[
\nbigh^j
\bigl(
\gbigm\langle\boxtimes,\otimes\rangle^!
 \nbigo[\ast 0]
 \bigr)
 \simeq
 \left\{
\begin{array}{ll}
 \lambda^{-1}\gbigm[\ast (0\times X)]& (j=1)\\
 0 & (j\neq 1),
\end{array}
 \right.
\quad\quad
 \nbigh^{j}\bigl(
 \gbigm\langle\boxtimes,\otimes\rangle^{\ast}
 \nbigo[!0]
 \bigr)
\simeq
 \left\{
\begin{array}{ll}
 \gbigm[!(0\times X)]& (j=-1)\\
 0 & (j\neq -1).
\end{array}
 \right.
\]
\begin{lem}
We have
$\FT_{\pm}(\nbigo[\ast 0])
 \simeq
 \lambda^{-1}\nbigo[!0]$
and
$\FT_{\pm}(\nbigo[!0])
 \simeq
 \nbigo[\ast 0]$.
\end{lem}
\pf
It follows from 
\cite[Proposition 2.23]{Mochizuki-GKZ}.
See also the proof of Lemma \ref{lem;22.7.26.11}
below.
\hfill\qed

\vspace{.1in}

By Proposition \ref{prop;21.7.13.10},
we have
$\nbigh^j\bigl(
 \gbigm\langle\boxtimes,\ast\rangle_{\ast}
 \nbigo[!0]
 \bigr)=0$
and 
$\nbigh^j\bigl(
 \gbigm\langle\boxtimes,\ast\rangle_{!}
 \nbigo[\ast 0]
 \bigr)=0$
unless $j=0$.
Moreover,
if $\gbigm\in\gbigctilde(\cnum\times X)_{\ast}$,
we obtain
$\nbigh^0\bigl(
 \gbigm\langle\boxtimes,\ast\rangle_{\ast}
 \nbigo[!0]
 \bigr)=\gbigm$,
and if 
$\gbigm\in\gbigctilde(\cnum\times X)_{!}$,
we obtain
$\nbigh^0\bigl(
 \gbigm\langle\boxtimes,\ast\rangle_{!}
 \nbigo[\ast 0]
 \bigr)=\gbigm$.

For $\gbigm\in\gbigctilde(\cnum\times X)$,
we define
\[
 \ttP_{\ast}(\gbigm):=
 \nbigh^0\bigl(
 \gbigm\langle\boxtimes,\ast\rangle_{\ast}\nbigo[!0]
 \bigr)
\in\gbigctilde(\cnum\times X)_{\ast},
 \quad\quad
 \ttP_{!}(\gbigm):=
 \nbigh^0\bigl(
 \gbigm\langle\boxtimes,\ast\rangle_{!}\nbigo[\ast 0]
 \bigr)
 \in\gbigctilde(\cnum\times X)_!.
\]
For $\star=!,\ast$,
there exist natural isomorphisms
$\ttP_{\star}\circ\ttP_{\star}(\gbigm)\simeq \ttP_{\star}(\gbigm)$.
If $\gbigm\in\gbigctilde(\cnum\times X)_{\star}$,
there exists a natural isomorphism
$\gbigm\simeq\ttP_{\star}(\gbigm)$.

\subsection{$\cnum^{\ast}$-homogeneity}

\subsubsection{Preliminary}

Let $Y$ be a complex manifold with a hypersurface $H$.
Let $\mu_1:\cnum^{\ast}\times Y\lrarr Y$ be a $\cnum^{\ast}$-action on $Y$,
which preserves $H$.
Let $\mu_2:\cnum^{\ast}\times\proj^1_{\lambda}\lrarr\proj^1_{\lambda}$
be a $\cnum^{\ast}$-action such that
$\mu_2(a,\lambda)=a^{m}\lambda$ for a non-zero integer $m$.
Let $\mu:\cnum^{\ast}\times \proj^1_{\lambda}\times Y\lrarr \proj^1_{\lambda}\times Y$
be the $\cnum^{\ast}$-action induced by $\mu_1$ and $\mu_2$.

Let $p:\cnum^{\ast}\times\proj^1_{\lambda}\times Y
 \lrarr
 \proj^1_{\lambda}\times Y$
denote the projection. 
We say that $\gbigm\in\gbigc(Y;H)$ is
$\cnum^{\ast}$-homogeneous with respect to $\mu$
if there exists an isomorphism
$\mu^{\ast}(\gbigm)
\simeq
p^{\ast}(\gbigm)$
of $\gbigrtilde_{\cnum^{\ast}\times Y}(\ast(\cnum^{\ast}\times H))$
satisfying the naturally defined cocycle condition.
(See \cite[\S3.6]{Mochizuki-GKZ}.)

\subsubsection{$\cnum^{\ast}$-homogeneity and Fourier transforms}

For $\vecn=(n_1,n_2)\in\seisuu^2_{\geq 0}$
such that $\gcd(n_1,n_2)=1$ and $n_1\geq n_2$,
we define the $\cnum^{\ast}$-action $\rho_{\vecn}$ on
$\proj^1_{\lambda}\times\proj^1_{t}\times X$
by
$\rho_{\vecn}(a)(\lambda,t,x)=(a^{n_1}\lambda,a^{n_2}t,x)$.
\begin{rem}
Though we may also consider the action $\rho_{\vecn}$
for any $\vecn\in\seisuu^2$,
we restrict ourselves to the case $n_1\geq n_2\geq 0$,
for simplicity.
We are mainly interested in the cases
$\vecn=(1,1),(1,0)$.
\hfill\qed
\end{rem}

Let $\gbigctilde_{\vecn}(\cnum_t\times X)
\subset\gbigctilde(\cnum_t\times X)$
denote the full subcategory of
$\gbigm\in\gbigctilde(\cnum_t\times X)$
which are $\cnum^{\ast}$-homogeneous with respect to $\rho_{\vecn}$.
Let $\gbigctilde_{\vecn}(\cnum_t\times X,[\star 0])
=\gbigctilde(\cnum_t\times X,[\star 0])
\cap
\gbigctilde_{\vecn}(\cnum_t\times X)$
and 
$\gbigctilde_{\vecn}(\cnum_t\times X)_{\star}
=\gbigctilde(\cnum_t\times X)_{\star}
\cap
\gbigctilde_{\vecn}(\cnum_t\times X)$.

\begin{lem}
\label{lem;21.7.11.11}
For $\vecn=(n_1,n_2)$, we define $\FT(\vecn):=(n_1,n_1-n_2)$.
Then, $\FT_{\pm}$ induce equivalences
\[
\gbigctilde_{\vecn}(\cnum_t\times X)
 \simeq
 \gbigctilde_{\FT(\vecn)}(\cnum_{\tau}\times X),
 \quad
\gbigctilde_{\vecn}(\cnum_t\times X,[\star 0])
 \simeq
 \gbigctilde_{\FT(\vecn)}(\cnum_{\tau}\times X)_{\star}.
\]
\end{lem}
\pf
Let us consider the action $\rhotilde_{\vecn}$
on $\proj^1_{\lambda}\times\proj^1_{t}\times\proj^1_{\tau}\times X$
defined by
$\rhotilde_{\vecn}(a)(\lambda,t,\tau,x)=
(a^{n_1}\lambda,a^{n_2}t,a^{n_1-n_2}\tau,x)$.
Then,
$\id_{\proj^1_{\lambda}}\times p_{t,X}:
\proj^1_{\lambda}\times
\cnum^2_{t,\tau}
\times X
\lrarr
\proj^1_{\lambda}\times
\cnum_{t}\times X$
is $\cnum^{\ast}$-equivariant with respect to
$\rhotilde_{\vecn}$ and $\rho_{\vecn}$,
and
$\id_{\proj^1_{\lambda}}\times p_{\tau,X}:
\proj^1_{\lambda}\times
\cnum^2_{t,\tau}\times X
\lrarr
\proj^1_{\lambda}\times
\cnum_{\tau}\times X$
is $\cnum^{\ast}$-equivariant with respect to
$\rhotilde_{\vecn}$ and $\rho_{\FT(\vecn)}$.
Note that
$\gbigl(\pm t\tau)$ are $\cnum^{\ast}$-homogeneous
with respect to $\rhotilde_{\vecn}$.
Hence, 
$p_{t,X}^{\ast}(\gbigm)\otimes\gbigl(t\tau)$
is naturally $\cnum^{\ast}$-homogeneous
with respect to $\rhotilde_{\vecn}$.
Then, the claim follows.
\hfill\qed

\begin{rem}
Because
$\FT_+$ induces
an equivalence
 $\gbigctilde_{\vecn}(\cnum_t\times X,[\star 0])\simeq
 \gbigctilde_{\vecn}(\cnum_{\tau}\times X)_{\star}$,
$\FT_-$ induces
an equivalence
$\gbigctilde_{\vecn}(\cnum_{\tau}\times X)_{\star}\simeq
\gbigctilde_{\vecn}(\cnum_{t}\times X,[\star 0])$.
\hfill\qed
\end{rem}

\begin{cor}
$\FT_{\pm}$ induce equivalences
$\gbigctilde_{(1,1)}(\cnum_t\times X)
 \simeq
 \gbigctilde_{(1,0)}(\cnum_t\times X)$.
Moreover, we obtain
$\gbigctilde_{(1,1)}(\cnum_{\tau}\times X,[\star 0])
\simeq
\gbigctilde_{(1,0)}(\cnum_t\times X)_{\star}$.
\hfill\qed
\end{cor}

\subsubsection{Functoriality}
\label{subsection;21.7.15.1}

We mention some basic functorial properties.

\begin{lem}
A projective morphism of complex manifolds
$F:X\lrarr Y$ induces
$F_{\dagger}^j:\gbigctilde_{\vecn}(\cnum\times X)\lrarr
 \gbigctilde_{\vecn}(\cnum\times Y)$,
$F_{\dagger}^j:\gbigctilde_{\vecn}(\cnum\times X)_{\star}
 \lrarr
 \gbigctilde_{\vecn}(\cnum\times Y)_{\star}$
and
$F_{\dagger}^j:\gbigctilde_{\vecn}(\cnum\times X,[\star 0])
 \lrarr
 \gbigctilde_{\vecn}(\cnum\times Y,[\star 0])$.
\hfill\qed
\end{lem}

\begin{lem}
Let $\gbigm\in\gbigctilde_{\vecn}(\cnum\times X)$.
For any hypersurface $H$ of $X$,
$\gbigm[\star (\cnum\times H)]$ $(\star=!,\ast)$
are also objects of
$\gbigctilde_{\vecn}(\cnum\times X)$.
Let $\star=!,\ast$.
For $\gbigm\in\gbigctilde_{\vecn}(\cnum\times X,[\star 0])$,
we obtain
$\gbigm[\star(\cnum\times H)]
\in \gbigctilde_{\vecn}(\cnum\times X,[\star 0])$.
Similarly,
for $\gbigm\in\gbigctilde_{\vecn}(\cnum\times X)_{\star}$,
we obtain
$\gbigm[\star (\cnum\times H)]
 \in\gbigctilde_{\vecn}(\cnum\times X)_{\star}$.
\hfill\qed
\end{lem}

\begin{cor}
In the situation of {\rm\S\ref{subsection;21.6.22.2}}
and {\rm\S\ref{subsection;21.6.22.20}},
for any $\gbigm\in\gbigctilde_{\vecn}(\cnum\times Y)$,
$(\lefttop{T}(\id\times f)^{\star})^i(\gbigm)$
are objects of
$\gbigctilde_{\vecn}(\cnum\times X)$.
\hfill\qed  
\end{cor}

\begin{lem}
\label{lem;21.7.15.2}
For any holomorphic function $f$ on $X$,
let $\ftilde$ denote the induced holomorphic function
on $\proj^1\times X$.
Then, $\Pi^{a,b}_{\ftilde,\star}(\gbigm)$
and
$\Pi^{a,b}_{\ftilde,\ast!}(\gbigm)$
are objects of
$\gbigctilde_{\vecn}(\cnum\times X)$.
As a result,
$\Xi^{(a)}_{\ftilde}(\gbigm)$,
$\psi^{(a)}_{\ftilde}(\gbigm)$
and $\phi^{(0)}_{\ftilde}(\gbigm)$
are objects of
$\gbigctilde_{\vecn}(\cnum\times X)$. 
\hfill\qed
\end{lem}

\begin{lem}
For $\gbigm\in\gbigctilde_{\vecn}(\cnum\times X)$,
$\DD(\gbigm)$ is an object of
$\gbigctilde_{\vecn}(\cnum\times X)$.
Let $(\star_1,\star_2)=(!,\ast),(\ast,!)$.
If $\gbigm\in\gbigctilde_{\vecn}(\cnum\times X,[\star_1 0])$,
 $\DD\gbigm$ is an object of
 $\gbigctilde_{\vecn}(\cnum\times,[\star_20])$.
If $\gbigm\in\gbigctilde_{\vecn}(\cnum\times X)_{\star_1}$,
$\DD\gbigm$ is an object of
$\gbigctilde_{\vecn}(\cnum\times X)_{\star_2}$.
\hfill\qed
\end{lem}

We set $X=X_1\times X_2$.
\begin{lem}
For $\gbigm_i\in\gbigctilde_{\vecn}(\cnum\times X_i)$,
$\nbigh^j\bigl(
 \gbigm_1\langle\boxtimes,\otimes\rangle^{\star}
 \gbigm_2
\bigr)$
and 
$\nbigh^j\bigl(
 \gbigm_1\langle\boxtimes,\ast\rangle_{\star}
 \gbigm_2
\bigr)$
 are objects of
$\gbigctilde_{\vecn}(\cnum\times X)$.
\hfill\qed
\end{lem}

\begin{lem}
For any $\gbigm_i\in \gbigctilde_{(1,1)}(\cnum\times X_i,[\ast 0])$,
we have the vanishing
$\nbigh^j(\gbigm_1\langle\boxtimes,\otimes\rangle^!\gbigm_2)=0$
unless $j=1$,
and 
$\nbigh^1(\gbigm_1\langle\boxtimes,\otimes\rangle^!\gbigm_2)$
is an object of
$\gbigctilde_{(1,1)}(\cnum\times X,[\ast 0])$.
For any $\gbigm_i\in \gbigctilde_{(1,1)}(\cnum\times X_i,[!0])$,
we have the vanishing
$\nbigh^j(\gbigm_1\langle\boxtimes,\otimes\rangle^{\ast}\gbigm_2)=0$
unless $j=-1$,
and
$\nbigh^{-1}(\gbigm_1\langle\boxtimes,\otimes\rangle^{\ast}\gbigm_2)$
is an object of
$\gbigctilde_{(1,1)}(\cnum\times X,[! 0])$.
\end{lem}
\pf
Let us explain the first claim.
Note that
$\nbigh^j(\gbigm_1\langle\boxtimes,\otimes\rangle^!\gbigm_2)$
are objects in
$\gbigctilde_{(1,1)}(\cnum\times X,[\ast 0])$.
We obtain
$\nbigh^j(\gbigm_1\langle\boxtimes,\otimes\rangle^!\gbigm_2)(\ast
(0\times X))=0$ unless $j=1$
by using Corollary \ref{cor;21.4.21.1}.
Then, the first claim of the lemma follows.
The second claim can be obtained similarly.
\hfill\qed

\subsection{Rescalable objects along $\{0,\infty\}$}

\subsubsection{Rescalable objects and $\cnum^{\ast}$-homogeneous objects}

We set $\gbigc_{\res,0,\infty}(X)\subset\gbigc(X)$
denote the full subcategory of
$\gbigm\in\gbigc(X)$ such that
$\gbigm\in\gbigc_{\res}(X)$ and $\gbigm\in\gbigc_{\res,\infty}(X)$.

\begin{example}
Let $\gbigr_F(M)\otimes\gbigl(f)\in\gbigc_{\Malg}(X;H)$
be as in {\rm\S\ref{subsection;22.7.27.10}}.
Then, 
$\bigl(
\gbigr_F(M)\otimes\gbigl(f)
\bigr)[\star H]$
are objects of $\gbigc_{\res,0,\infty}(X)$.
\hfill\qed
\end{example}

For any $\gbigm\in\gbigc_{\res,0,\infty}(X)$,
we can naturally regard $\lefttop{\tau}\gbigm$
as an object of
$\gbigctilde(\cnum^{\ast}_{\tau}\times X)$.
We obtain
\[
\ttQ_{\star}(\gbigm):=
\lefttop{\tau}\gbigm[\star (\{0\}\times X)]\in
 \gbigctilde_{(1,1)}(\cnum_{\tau}\times X,[\star 0]). 
\]

\begin{prop}
\label{prop;21.7.11.31}
The procedure induces equivalences
$\ttQ_{\star}:\gbigc_{\res,0,\infty}(X)\simeq
\gbigctilde_{(1,1)}(\cnum_{\tau}\times X,[\star 0])$.
\end{prop}
\pf
It is easy to see that the functor is fully faithful.
Let $\gbign\in\gbigctilde_{(1,1)}(\cnum_{\tau}\times X,[\star 0])$.
By the $\cnum^{\ast}$-homogeneity with respect to $\rho_{1,1}$,
we obtain that
$\{1\}\times X\subset \cnum_{t}\times X$ is strictly 
non-characteristic for $\gbign$.
Let $\iota_1:\{1\}\times X\lrarr\cnum_t\times X$
denote the inclusion.
We set $\gbigm=\iota_1^{\ast}\gbign\in\gbigc_{\Malg}(X)$.
By the $\cnum^{\ast}$-homogeneity,
there exists an isomorphism
$\lefttop{\tau}(\gbigm_{|\nbigx})
\simeq
 \gbign_{|\cnum_{\lambda}\times\cnum_{\tau}^{\ast}\times X}$.
By Theorem \ref{thm;21.3.12.12} and Theorem \ref{thm;21.7.2.20},
we obtain $\gbign(\ast\{0,\infty\})=\lefttop{\tau}\gbigm$
in $\gbigctilde(\cnum_{\tau}^{\ast}\times X)$,
which implies that
$\gbign\simeq \ttQ_{\star}(\gbigm)$.
\hfill\qed

\vspace{.1in}
We state some basic functorial property,
which are easy to prove.

\begin{prop}
Let $\gbigm\in\gbigc_{\res,0,\infty}(X)$.
\begin{itemize}
 \item 
Let $F:X\lrarr Y$ be a projective morphism of complex manifolds.
Then,
$F^j_{\dagger}(\gbigm)$ are objects of
$\gbigc_{\res,0,\infty}(Y)$,
and there exist natural isomorphisms
$\ttQ_{\star}\bigl(
F^j_{\dagger}(\gbigm)
\bigr)
\simeq
(\id\times F)^j_{\dagger}\bigl(
\ttQ_{\star}(\gbigm)
\bigr)$.
 \item Let $H$ be a hypersurface of $X$.
       Then, for any $\star=!,\ast$,
       $\gbigm[\star H]$
       are objects of $\gbigc_{\res,0,\infty}(X)$,
       and there exist isomorphisms
       $\ttQ_{\star}(\gbigm[\star H])
       \simeq
       \bigl(
       \ttQ_{\star}(\gbigm)[\star (\proj^1\times H)]
       \bigr)$.
 \item Let $f$ and $\ftilde$ be as in Proposition {\rm\ref{prop;21.7.13.20}}.
       Then, $\Pi^{a,b}_{f,\star}(\gbigm)$ $(\star=!,\ast)$
       and $\Pi^{a,b}_{f,\ast!}(\gbigm)$
       are objects of
       $\gbigc_{\res,0,\infty}(X)$,
       and there exist the following natural isomorphisms
       for any $\star'=!,\ast$:
\[
       \ttQ_{\star'}\Pi^{a,b}_{f,\star}(\gbigm)
       \simeq
       \bigl(
       \Pi^{a,b}_{\ftilde,\star}
       \ttQ_{\star'}(\gbigm)
       \bigr)[\star'(0\times X)]\,\,(\star=!,\ast),
       \quad
       \ttQ_{\star'}\Pi^{a,b}_{f,\ast!}(\gbigm)
       \simeq
       \bigl(
       \Pi^{a,b}_{\ftilde,\ast !}
       \ttQ_{\star'}(\gbigm)
       \bigr)[\star'(0\times X)].
\]       
       As a result,
       $\Xi^{(a)}_{f}(\gbigm)$,
       $\psi^{(a)}_f(\gbigm)$
       and $\phi^{(0)}_f(\gbigm)$
       are objects of
       $\gbigc_{\res,0,\infty}(X)$,
       and there exist natural isomorphisms
       for $\star'=!,\ast$:
\[
       \ttQ_{\star'}\Xi^{(a)}_f(\gbigm)
       \simeq
       \bigl(
       \Xi^{(a)}_{\ftilde}\ttQ_{\star'}(\gbigm)
       \bigr)[\star' (0\times X)],
       \quad
       \ttQ_{\star'}\psi^{(a)}_f(\gbigm)
       \simeq
       \bigl(
       \psi^{(a)}_{\ftilde}\ttQ_{\star'}(\gbigm)
       \bigr)[\star'(0\times X)],
\]
\[
       \ttQ_{\star'}\phi^{(0)}_f(\gbigm)
       \simeq
       \bigl(
       \phi^{(0)}_{\ftilde}\ttQ_{\star'}(\gbigm)
       \bigr)[\star'(0\times X)].
\]
\hfill\qed
\end{itemize}
\end{prop}

\begin{cor}
In the situation of {\rm\S\ref{subsection;21.6.22.2}}
and {\rm\S\ref{subsection;21.6.22.20}},
 for any $\gbigm\in\gbigc_{\res,0,\infty}(Y)$
 and for any $\star=!,\ast$,
$(\lefttop{T}f^{\star})^j(\gbigm)$
are objects of
$\gbigc_{\res,0,\infty}(X)$,
and there exist natural isomorphisms
$\ttQ_{\star'}(\lefttop{T}f^{\star})^j(\gbigm)
 \simeq
 \Bigl(
 \bigl(
 \lefttop{T}(\id\times f)^{\star}
 \bigr)^j
 \ttQ_{\star'}(\gbigm)
 \Bigr)[\star'(0\times X)]$ $(\star'=!,\ast)$.
\hfill\qed
\end{cor}

\begin{prop}
Let $(\star_1,\star_2)=(!,\ast),(\ast,!)$.
For $\gbigm\in\gbigc_{\res,0,\infty}(X)$,
$\DD\gbigm$ is an object of
$\gbigc_{\res,0,\infty}(X)$,
and there exist natural isomorphisms
$\lambda \ttQ_{\star_1}(\DD\gbigm)\simeq
 \DD\ttQ_{\star_2}(\gbigm)$.
\hfill\qed
\end{prop}

Let $X=X_1\times X_2$.
\begin{prop}
For $\gbigm_i\in\gbigc_{\res,0,\infty}(X_i)$,
$\gbigm_1\boxtimes\gbigm_2$
is an object of $\gbigc_{\res,0,\infty}(X)$,
and there exist natural isomorphisms
\[
 \lambda^{-1}\ttQ_{\ast}(\gbigm_1\boxtimes\gbigm_2)
 \simeq
 \nbigh^1\bigl(
 \ttQ_{\ast}(\gbigm_1)
 \langle
 \boxtimes,\otimes
 \rangle^!
 \ttQ_{\ast}(\gbigm_2)
 \bigr),
 \quad
 \ttQ_{!}(\gbigm_1\boxtimes\gbigm_2)
 \simeq
 \nbigh^{-1}\bigl(
 \ttQ_{!}(\gbigm_1)
 \langle
 \boxtimes,\otimes
 \rangle^{\ast}
 \ttQ_{!}(\gbigm_2)
 \bigr).
\]
\hfill\qed
\end{prop}

\subsubsection{Filtered $\nbigd$-modules and $\cnum^{\ast}$-homogeneous objects}
\label{subsection;21.7.25.3}

Let $\nbigctilde_{\Hod}(\cnum\times X)$ denote the category of
filtered $\nbigd_{\proj^1\times X}(\ast(\{\infty\}\times X))$-modules
$(M,F)$ 
such that
$\gbigr_F(M)\in\gbigc_{\Malg}(\proj^1\times X;\{\infty\}\times X)$.
(See \S\ref{subsection;22.7.27.10}
for $\gbigr_F(M)$.)
According to \cite[Proposition 3.29]{Mochizuki-GKZ},
it induces an equivalence
$\nbigctilde_{\Hod}(\cnum\times X)
\to
\gbigctilde_{(1,0)}(\cnum\times X)$.
For $\gbigm\in\gbigctilde_{(1,0)}(\cnum\times X)$,
we obtain the $\nbigd_X$-module $\Xi_{\DR}(\gbigm)$,
on which there exists a unique filtration $F$
such that
$\gbigr_F(\Xi_{\DR}(\gbigm))=\gbigm$.
Thus, we obtain
$\ttR:\gbigctilde_{(1,0)}(\cnum\times X)
\lrarr \nbigctilde_{\Hod}(\cnum\times X)$
by
$\ttR(\gbigm)=(\Xi_{\DR}(\gbigm),F)$,
which is a quasi-inverse of the above equivalence.
Let $\nbigctilde_{\Hod}(\cnum\times X)_{\star}
\subset\nbigctilde_{\Hod}(\cnum\times X)$
denote the full subcategory of
$(M,F)\in\nbigctilde_{\Hod}(\cnum\times X)$
such that
$\pi_{X\dagger}^j\bigl(
 M(\star(\infty\times X))
\bigr)=0$
$(j\in\seisuu)$.
We have the equivalences
$\ttR:
\gbigctilde_{(1,0)}(\cnum\times X)_{\star}\simeq
\nbigctilde_{\Hod}(\cnum\times X)_{\star}$.
As a consequence of the previous results,
we obtain the following theorem.

\begin{thm}
\label{thm;21.7.14.2}
There exist equivalences
$\ttA_{\pm,\star}:=\ttR\circ\FT_{\pm}\circ\ttQ_{\star}:
 \gbigc_{\res,0,\infty}(X)
 \simeq
 \nbigctilde_{\Hod}(\cnum\times X)_{\star}$
 $(\star=!,\ast)$.
\hfill\qed
\end{thm}

The equivalences in Theorem \ref{thm;21.7.14.2}
are compatible with basic functors.
Let
$\ttP_{\star}:\nbigctilde_{\Hod}(\cnum\times X)
\lrarr
\nbigctilde_{\Hod}(\cnum\times X)_{\star}$
denote the equivalence 
induced by
$\ttP_{\star}:\gbigctilde_{(1,0)}(\cnum\times X)
\lrarr
\gbigctilde_{(1,0)}(\cnum\times X)_{\star}$
and the equivalence $\ttR$.
The following proposition is easy to see.

\begin{prop}
\label{prop;21.7.26.10}
Let $\gbigm\in\gbigc_{\res,0,\infty}(X)$.
\begin{itemize}
 \item
      Let $F:X\lrarr Y$ be a projective morphism of complex manifolds.
      There exist natural isomorphisms
      $\ttA_{\pm,\star}F_{\dagger}^j(\gbigm)
      \simeq
      (\id\times F)^j_{\dagger}(\ttA_{\pm,\star}(\gbigm))$.
 \item Let $H$ be a hypersurface of $X$.
       For any $\star,\star'\in\{!,\ast\}$,
       there exist natural isomorphisms
       $\ttA_{\pm,\star'}(\gbigm[\star H])
       \simeq
        \ttP_{\star'}\bigl(
       \ttA_{\pm,\star'}(\gbigm)[\star (\cnum\times H)]
       \bigr)$.
 \item Let $f$ and $\ftilde$ be as in Proposition {\rm\ref{prop;21.7.13.20}}.
       For $\star,\star'\in\{!,\ast\}$,
       there exist the following natural isomorphisms:
 \[
       \ttA_{\pm,\star}\Pi^{a,b}_{f,\star'}(\gbigm)
       \simeq
       \ttP_{\star}
       \bigl(
       \Pi^{a,b}_{\ftilde,\star'}
       \ttA_{\pm,\star}(\gbigm)
       \bigr),
       \quad
       \ttA_{\pm,\star}\Pi^{a,b}_{f,\ast!}(\gbigm)
       \simeq
       \ttP_{\star}
       \bigl(
       \Pi^{a,b}_{\ftilde,\ast !}
       \ttA_{\pm,\star}(\gbigm)
       \bigr).
 \]
       As a result,
       there exist the following natural isomorphisms
       for $\star=!,\ast$:
\[
       \ttA_{\pm,\star}\Xi^{(a)}_f(\gbigm)
       \simeq
       \ttP_{\star}
       \bigl(
       \Xi^{(a)}_{\ftilde}\ttA_{\pm,\star}(\gbigm)
       \bigr),
       \quad
       \ttA_{\pm,\star}\psi^{(a)}_f(\gbigm)
       \simeq
       \ttP_{\star}
       \bigl(
       \psi^{(a)}_{\ftilde}\ttA_{\pm,\star}(\gbigm)
       \bigr),
       \quad
       \ttA_{\pm,\star}\phi^{(0)}_f(\gbigm)
       \simeq
       \ttP_{\star}
       \bigl(
       \phi^{(0)}_{\ftilde}\ttA_{\pm,\star}(\gbigm)
       \bigr).
\]
\hfill\qed
\end{itemize} 
\end{prop}

\begin{cor}
\label{cor;21.7.26.11}
In the situation of {\rm\S\ref{subsection;21.6.22.2}}
and {\rm\S\ref{subsection;21.6.22.20}},
for any $\gbigm\in\gbigc_{\res,0,\infty}(Y)$,
there exist natural isomorphisms
$\ttA_{\pm,\star'}(\lefttop{T}f^{\star})^j(\gbigm)
 \simeq
 \ttP_{\star'}
 \Bigl(
 \bigl(
 \lefttop{T}(\id\times f)^{\star}
 \bigr)^j
 \ttA_{\pm,\star'}(\gbigm)
 \Bigr)$.
\hfill\qed
\end{cor}

For a filtered $\nbigd$-module
$(M,F_{\bullet})$,
let $\ttS_i(M,F)$ denote the filtered $\nbigd$-module $(M,F')$
determined by $F'_j(M)=F_{j-i}(M)$.
Note that $\gbigr_{F'}(M)=\lambda^i\gbigr_F(M)$.
We obtain the following propositions.

\begin{prop}
Let $(\star_1,\star_2)=(!,\ast),(\ast,!)$.
For any $\gbigm\in\gbigc_{\res,0,\infty}(X)$,
there exist natural isomorphisms
$\ttA_{\pm,\star_1}(\DD\gbigm)\simeq
\ttS_{-2}\DDD\ttA_{\pm,\star_2}(\gbigm)$.
\hfill\qed
\end{prop}

\begin{prop}
For $\gbigm_i\in\gbigc_{\res,0,\infty}(X_i)$,
there exist natural isomorphisms
\[
 \ttS_{-1}\ttA_{\pm,\star}(\gbigm_1\boxtimes\gbigm_2)
 \simeq
 \nbigh^0\bigl(
 \ttA_{\pm,\star}(\gbigm_1)
 \langle
 \boxtimes,\ast
 \rangle_{\star}
 \ttA_{\pm,\star}(\gbigm_2)
 \bigr)
 \quad(\star=!,\ast).
\]
\hfill\qed
\end{prop}

\subsubsection{Quasi-inverse functors}
\label{subsection;21.7.26.12}

Let $\iota_1:X\lrarr \proj^1\times X$
be the morphism defined by
$\iota_1(x)=(1,x)$.
We define the functors
$\ttB_{\pm}:\nbigctilde_{\Hod}(\cnum\times X)
\lrarr\gbigc_{\res,0,\infty}(X)$
by
\[
 \ttB_{\pm}(M,F)
 =
 \iota_1^{\ast}\FT_{\pm}(\gbigr_F(M)).
\]
We obtain
$\ttB_{\pm,\star}:\nbigctilde_{\Hod}(\cnum\times X)_{\star}
\lrarr\gbigc_{\res,0,\infty}(X)$
as the restrictions of $\ttB_{\pm}$.
By a formal computation,
we have
\[
 \ttB_{\mp,\star}\circ\ttA_{\pm,\star}(\gbigm)
 \simeq
 \lambda^{-1}\gbigm.
\]
Therefore, $\ttB_{\mp,\star}$ are quasi-inverse
functors of $\ttA_{\pm,\star}$.
By the construction, we have
\[
 \ttB_{\pm,\star}(M,F)
=
 \pi_{X\dagger}^0
 \bigl(
 \gbigr_F(M)\otimes\gbigl(t)
 \bigr).
\]
We have the compatibility of $\ttB_{\pm,\star}$
with the other functors as in
Proposition \ref{prop;21.7.26.10}
and Corollary \ref{cor;21.7.26.11}.
We have the following compatibility of
$\ttB_{\pm,\star}$ with the duality functor,
for $(\star_1,\star_2)=(!,\ast),(\ast,!)$:
\[
 \ttB_{\pm,\star_1}(\DDD(M,F))
 \simeq
 \DD\ttB_{\mp,\star_2}(M,F).
\]
We have the following compatibility of
$\ttB_{\pm,\star}$
with the tensor and convolution products:
\[
 \ttB_{\pm,\star}
 \bigl(
 (M_1,F)\langle \boxtimes,\star\rangle_{\star}
 (M_2,F)
 \bigr)
 \simeq
 \ttB_{\pm,\star}(M_1,F)
 \boxtimes
 \ttB_{\pm,\star}(M_2,F).
\]

\section{Rescalable integrable mixed twistor $\nbigd$-modules}

\subsection{Preliminary}

\subsubsection{Some categories}

Let $X$ be a complex manifold with a hypersurface $H$.
We set $d_X:=\dim X$.
An $\nbigrtilde_{X(\ast H)}$-triple
is a pair of $\nbigrtilde_{X(\ast H)}$-modules
$\nbigm'$ and $\nbigm''$
with an integrable sesqui-linear pairing $C$
of $\nbigm'$ and $\nbigm''$.
(See \cite[\S2.1]{Mochizuki-MTM}.)
Such a triple is denoted by $\nbigt=(\nbigm',\nbigm'',C)$.

An integrable mixed twistor $\nbigd$-module on $X$
consists of
an $\nbigrtilde_X$-triple $\nbigt$
with a weight filtration $W$
satisfying some conditions.
(See \cite[\S7.2]{Mochizuki-MTM}.)
We shall often omit to denote the weight filtration $W$.
Let $\MTM^{\integral}(X)$
denote the category of integrable mixed twistor $\nbigd$-modules
on $X$.

A sesqui-linear pairing $C$ of
$\nbigrtilde_X$-modules $\nbigm'$ an $\nbigm''$
uniquely induces
a sesqui-linear pairing $C(\ast H)$
of the $\nbigrtilde_{X(\ast H)}$-modules
$\nbigm'(\ast H)$ and $\nbigm''(\ast H)$.
For any $\nbigrtilde_X$-triple $\nbigt=(\nbigm',\nbigm'',C)$,
the $\nbigrtilde_{X(\ast H)}$-triple
$(\nbigm'(\ast H),\nbigm''(\ast H),C(\ast H))$
is denoted by $\nbigt(\ast H)$.
This induces a functor from
$\MTM^{\integral}(X)$
to the category of filtered $\nbigrtilde_{X(\ast H)}$-modules.
The essential image is denoted by $\MTM^{\integral}(X;H)$.

Let $\nbigt=(\nbigm',\nbigm'',C)\in\MTM^{\integral}(X;H)$.
There exist $\nbigrtilde_X$-modules
$\nbigm'[\star H]$ and $\nbigm''[\star H]$
$(\star=!,\ast)$,
and $C$ uniquely induces the sesqui-linear pairings
$C_{\ast}$ of $\nbigm'[!H]$ and $\nbigm''[\ast H]$,
and
$C_!$ of $\nbigm'[\ast H]$ and $\nbigm''[!H]$.
We obtain the $\nbigrtilde_{X}$-triples
$\nbigt[\ast H]=(\nbigm'[!H],\nbigm''[\ast H],C_{\ast})$
and
$\nbigt[!H]=(\nbigm'[\ast H],\nbigm''[!H],C_{!})$.
There exist naturally induced weight filtrations $W$
on $\nbigt[\star H]$ $(\star=!,\ast)$
with which
$\nbigt[\star H]\in\MTM^{\integral}(X)$.
(See \cite[\S11.2]{Mochizuki-MTM}.)

Let $\MTM^{\integral}(X,[\star H])\subset\MTM^{\integral}(X)$
$(\star=!,\ast)$
denote the full subcategories of
objects $(\nbigt,W)\in\MTM^{\integral}(X)$
such that
$(\nbigt,W)[\star H]=(\nbigt,W)$.
Then, the induced functor
$\MTM^{\integral}(X,[\star H])\lrarr \MTM^{\integral}(X;H)$
is an equivalence.
A quasi-inverse 
$\MTM^{\integral}(X;H)\lrarr \MTM^{\integral}(X,[\star H])$
is given by
$(\nbigt,W)\longmapsto
(\nbigt,W)[\star H]$.

If $H=H_1\cup H_2$ with $\dim (H_1\cap H_2)<\dim H$,
we have the functors
$\MTM^{\integral}(X;H)
\lrarr \MTM^{\integral}(X;H_1)$
induced by 
$(\nbigt,W)\longmapsto
\bigl(
(\nbigt,W)[\star H]
\bigr)(\ast H_1)
=:(\nbigt,W)[\star H_2]$.

For any object $\nbigt=(\nbigm',\nbigm'',C)\in\MTM^{\integral}(X;H)$,
we naturally regard
$\Xi_{\DR}(\nbigt)=\nbigm''/(\lambda-1)\nbigm''$
as a $\nbigd_{X(\ast H)}$-module,
called the underlying $\nbigd_{X(\ast H)}$-module of $\nbigt$.
Note that the functor $\Xi_{\DR}$
from $\MTM^{\integral}(X)$
to the category of $\nbigd_{X(\ast H)}$-module
is faithful
as remarked in \cite[Remark 7.2.9]{Mochizuki-MTM}.

Let $f:X\to Y$ be a projective morphism.
For $\nbigt=(\nbigm',\nbigm'',C)\in\MTM^{\integral}(X)$,
the $\nbigrtilde_Y$-triple
$f_{\dagger}^j(\nbigt)
=(f_{\dagger}^{-j}(\nbigm'),f_{\dagger}^j(\nbigm''),
f_{\dagger}C)$
is obtained as in \cite{Sabbah-pure-twistor}.
There exists an induced weight filtration $W$
on $f_{\dagger}^j(\nbigt)$
with which
$f_{\dagger}^j(\nbigt)$
are objects in $\MTMtilde^{\integral}(Y)$.
It induces a cohomological functor
from $\MTM^{\integral}(X)$ to $\MTM^{\integral}(Y)$.
(See \cite[\S7.2]{Mochizuki-MTM}.)

\subsubsection{Real structure}

Let $j:\cnum_{\lambda}\lrarr\cnum_{\lambda}$
be the morphism defined by $j(\lambda)=-\lambda$.
The induced morphism
$\nbigx\lrarr\nbigx$ is also denoted by $j$.
Let $\nbigt=(\nbigm',\nbigm'',C)\in\MTM^{\integral}(X;H)$.
Recall that we have
the induced objects 
$j^{\ast}(\nbigt)=(j^{\ast}\nbigm',j^{\ast}\nbigm'',j^{\ast}C)
\in\MTM^{\integral}(X;H)$
and 
$\DD\nbigt=(\DD\nbigm',\DD\nbigm'',\DD C)
\in \MTM^{\integral}(X;H)$.
(See \cite[\S13]{Mochizuki-MTM} for $\DD\nbigt$.)
Recall that we have the Hermitian adjoint
$\nbigt^{\ast}=(\nbigm'',\nbigm',C^{\ast})$
of $\nbigt$
(see \cite{Sabbah-pure-twistor} or \cite[\S2.1]{Mochizuki-MTM}).
We set
$\gammatilde^{\ast}(\nbigt)
=\DD(j^{\ast}\nbigt^{\ast})$.
Recall that
a real structure of $\nbigt\in\MTM^{\integral}(X;H)$
is an isomorphism
$\kappa:\gammatilde^{\ast}\nbigt\simeq\nbigt$
satisfying
$\gammatilde^{\ast}(\kappa)\circ\kappa=\id$
(see \cite[\S13.4]{Mochizuki-MTM}).
Let $\MTM^{\integral}((X;H),\real)$
denote the category of
objects $(\nbigt,W)\in\MTM^{\integral}(X;H)$
equipped with a real structure $\kappa$.

\subsubsection{Integrable mixed twistor $\nbigd$-modules
associated with functions}
\label{subsection;22.7.26.5}

Let $f$ be a meromorphic function on $(X,H)$.
Let $\nbigl(f)$ be the $\nbigrtilde_{X(\ast H)}$-module
obtained as
$\nbigl(f)=\nbigo_{\nbigx}(\ast H)$
with the meromorphic integrable connection $d+d(\lambda^{-1}f)$.
There exists a naturally defined sesqui-linear pairing
$C(f)$ of the pair
of the $\nbigrtilde_{X(\ast H)}$-modules
$(\nbigl(f),\nbigl(f))$.
The tuple
\[
 \nbigt(f)=(\lambda^{d_X}\nbigl(f),\nbigl(f),C(f))
\]
is an object of
$\MTM^{\integral}(X;H)$,
which is pure of weight $d_X$.
(See \cite[\S3.2]{Mochizuki-GKZ}.)
We also have the following smooth $\nbigr_{X(\ast H)}$-triple
\[
 \nbigt_{\sm}(f)=(\nbigl(f),\nbigl(f),C(f)).
\]

We denote $\nbigt(0)$ by $\nbigu_X(d_X,0)$.
(See \cite[\S2.1.8.1]{Mochizuki-MTM}
for more general $\nbigu_X(a,b)$.)
We have
$\nbigt(f)=\nbigu_X(d_X,0)\otimes\nbigt_{\sm}(f)$.
There exists the real structure
$\kappa:\gammatilde^{\ast}\nbigu_X(d_X,0)
\simeq \nbigu_X(d_X,0)$
as in
\cite[\S13.4.2.1]{Mochizuki-MTM} and \cite[\S3.3.1]{Mochizuki-GKZ}.
We have the natural isomorphism
$j^{\ast}\nbigu_X(d_X,0)\simeq\nbigu_X(d_X,0)$.
We also have the sesqui-linear duality
$\nbigu_X(d_X,0)^{\ast}\simeq
\nbigu_X(d_X,0)\otimes\newTate(d_X)$
as in \cite[\S2.1.8.1]{Mochizuki-GKZ}.
We obtain the isomorphism
$\DD\nbigu_X(d_X,0)\simeq\nbigu_X(d_X,0)\otimes\newTate(d_X)$.
Because the smooth $\nbigr_{X(\ast H)}$-triple
$\nbigt_{\sm}(f)$ has the real structure
as in \cite[\S2.1.7.2]{Mochizuki-MTM},
$\nbigt(f)$ has the real structure
by \cite[Proposition 13.4.6]{Mochizuki-MTM}.
(See also \cite[\S3.3]{Mochizuki-GKZ}.)
Similarly, we have the natural isomorphisms
$j^{\ast}\nbigt(f)\simeq\nbigt(f)$,
$\nbigt(f)^{\ast}\simeq\nbigt(f)\otimes\newTate(d_X)$
and
$\DD\nbigt(f)\simeq\nbigt(f)\otimes\newTate(d_X)$.

Let $Y$ be any smooth hypersurface of $X$
such that $\dim (H\cap Y)<\dim Y$.
We set $H_Y:=H\cap Y$ and $f_Y:=f_{|Y}$.
Let $\iota_Y:Y\to X$ denote the inclusion.
As explained in \cite[Proposition 3.21]{Mochizuki-GKZ},
there exist the following natural exact sequences
in $\MTM^{\integral}((X;H),\real)$:
\begin{equation}
\label{eq;22.7.26.1}
 0\lrarr \nbigt(f)
 \lrarr \nbigt(f)[\ast Y]
 \lrarr \iota_{Y\dagger}(\nbigt(f_Y))\otimes\newTate(-1)
 \lrarr 0,
\end{equation}
\begin{equation}
\label{eq;22.7.26.2}
 0\lrarr
 \iota_{Y\dagger}\nbigt(f_Y)
 \lrarr
 \nbigt(f)[!Y]
 \lrarr\nbigt(f)\lrarr 0.
\end{equation}

\subsubsection{Notation}

Let $Z$ be a quasi-projective complex manifold.
There exists an algebraic Zariski open embedding $\iota_Z:Z\to \Zbar$
to a projective manifold $\Zbar$
such that $H_{Z,\infty}:=\Zbar\setminus Z$
is a hypersurface.
We set
\[
 \MTMtilde^{\integral}(Z\times X):=
 \MTM^{\integral}(\Zbar\times X;H_{Z,\infty}\times X).
\]
If we consider the real structure,
the category is denoted by 
$\MTMtilde^{\integral}(Z\times X,\real)$.
We set
$\nbigu_{Z}(d_Z,0):=
\nbigu_{\Zbar}(d_Z,0)(\ast H_{Z,\infty})
\in\MTM^{\integral}\bigl(Z,\real\bigr)$.

Let $\iota_Z':Z\to\Zbar'$
be an algebraic Zariski open embedding 
to a projective manifold $\Zbar'$
such that $H'_{Z,\infty}:=\Zbar'\setminus Z$
is a hypersurface.
Then, there exists
an algebraic Zariski open embedding 
$\iota''_Z:Z\to \Zbar''$
to a projective manifold $\Zbar''$
with morphisms
$\varphi_1:\Zbar''\to \Zbar$
and
$\varphi_2:\Zbar''\to \Zbar'$
as in \S\ref{subsection;22.7.27.12}.
We obtain the following equivalences:
\[
\begin{CD}
 \MTM^{\integral}(\Zbar\times X;H_{Z,\infty}\times X)
 @<{\varphi_{1\dagger}}<<
 \MTM^{\integral}(\Zbar''\times X;H''_{Z}\times X)
 @>{\varphi_{2\dagger}}>>
 \MTM^{\integral}(\Zbar'\times X;H'_{Z}\times X).
\end{CD}
\]
We have the inverse $\varphi_i^{\ast}$ of $\varphi_{i\dagger}$.
(See \cite[\S2.1]{Mochizuki-MTM}
for $\varphi_i^{\ast}$ in this situation.)
Hence, $\MTMtilde^{\integral}(Z\times X)$ is independent
of a projective completion.

Let $f:Z_1\to Z_2$ be an algebraic morphism of
complex quasi-projective manifolds.
It extends to a morphism of complex projective manifolds
$\overline{f}:\overline{Z}_1\to \overline{Z}_2$
as in \S\ref{subsection;22.7.27.12}.
For any
$\nbigt\in \MTMtilde^{\integral}(Z_1\times X)$,
we set
\[
 \lefttop{T}f_{\star}^j(\nbigt):=
 f^j_{\dagger}
 \bigl(
 \nbigt[\star (H_{Z_1}\times X)]
 \bigr)(\ast (H_{Z_2}\times X))
\in \MTMtilde(Z_2\times X).
\]

\subsubsection{Some lemmas}

We have $\nbigt(t\tau)\in\MTMtilde^{\integral}(\cnum^2_{t,\tau},\real)$
associated with the function $t\tau$
as in \S\ref{subsection;22.7.26.5}.
\begin{lem}
\label{lem;22.7.26.11}
Let
$p:\cnum^2_{t,\tau}\to \cnum_t$ 
denote the projection. 
Let $\iota_1:\{0\}\lrarr\cnum_t$ denote inclusion.
We have 
$\lefttop{T}p^j_{\ast}(\nbigt(t\tau))=0$ $(j\neq 0)$
and there exists an isomorphism
$\lefttop{T}p^0_{\ast}(\nbigt(t\tau))
\simeq
\iota_{1\dagger}\newTate(-1)$
in $\MTMtilde^{\integral}(\cnum_t,\real)$.
\end{lem}
\pf
Though this is a special case of \cite[Proposition 3.22]{Mochizuki-GKZ},
we explain an outline of the proof in this easy case
for the convenience of the readers.
We set $H_{\tau=0}=\cnum_t\times\{0\}\subset\cnum^2_{t,\tau}$.
By calculating the underling $\nbigd$-modules,
it is easy to check that
$\lefttop{T}p^j_{\ast}(\nbigt)=0$ $(j\neq 0)$
for
$\nbigt=\nbigt(t\tau),\,\nbigt(t\tau)[\ast H_{\tau=0}],\,
\nbigt(t\tau)[\ast H_{\tau=0}]\big/\nbigt(t\tau)$.
By (\ref{eq;22.7.26.1}),
we obtain
\[
 \lefttop{T}p^0_{\ast}\bigl(
 \nbigt(t\tau)[\ast H_{\tau=0}]\big/\nbigt(t\tau)
 \bigr)
\simeq \nbigu_{\cnum_t}(1,0)\otimes\newTate(-1).
\]
There exists the natural morphism
$\rho:\lefttop{T}p^0_{\ast}\bigl(
 \nbigt(t\tau)[\ast H_{\tau=0}]
 \bigr)[!0]
\lrarr
 \lefttop{T}p^0_{\ast}\bigl(
 \nbigt(t\tau)[\ast H_{\tau=0}]
 \bigr)$
in $\MTMtilde^{\integral}(\cnum_t)$.
By the standard computation of the Fourier transform,
we can check that 
$\Xi_{\DR}(\rho)$ is an isomorphism.
Hence, $\rho$ is an isomorphism.
It is also easy to check that the support of
$\lefttop{T}p^0_{\ast}\bigl(
 \nbigt(t\tau)
 \bigr)$
is contained in $\{0\}$.
Hence, we obtain
\[
  \lefttop{T}p^0_{\ast}\bigl(
 \nbigt(t\tau)[\ast H_{\tau=0}]
 \bigr)
 \simeq
 \lefttop{T}p^0_{\ast}\bigl(
 \nbigt(t\tau)[\ast H_{\tau=0}]
 \bigr)[!0]
 \simeq
 \lefttop{T}p^0_{\ast}\bigl(
 \nbigt(t\tau)[\ast H_{\tau=0}]\big/\nbigt(t\tau)
 \bigr)[!0]
 \simeq
 \nbigu_{\cnum_t}(1,0)[!0]\otimes\newTate(-1).
\]
By (\ref{eq;22.7.26.2}),
we obtain 
$\lefttop{T}p^0_{\ast}(\nbigt(t\tau))
\simeq
\iota_{1\dagger}\newTate(-1)$.
\hfill\qed

\vspace{.1in}

Let $\nbigt\in\MTMtilde^{\integral}(\cnum_t\times\cnum^m\times X)$.
We obtain
$\nbigt\boxtimes\nbigu_{\cnum_{\tau}}(1,0)
\in\MTMtilde^{\integral}(\cnum^2_{t,\tau}\times\cnum^m\times X)$.
Let $\iota:\cnum_t\times\cnum^m\times X\lrarr
\cnum^2_{t,\tau}\times\cnum^m\times X$
denote the inclusion defined by
$\iota(t,z_1,\ldots,z_m,x)=\iota(t,t,z_1,\ldots,z_m,x)$.
Let $H_{t=\tau}$ denote the image.

\begin{lem}
\label{lem;22.7.26.10}
 There exists the following exact sequence
in $\MTMtilde^{\integral}(\cnum^2_{t,\tau}\times\cnum^m\times X)$:
\begin{equation}
\label{eq;22.7.26.3}
 0\lrarr
 \nbigt\boxtimes\nbigu_{\cnum_{\tau}}(1,0)
 \lrarr
 \bigl(
 \nbigt\boxtimes\nbigu_{\cnum_{\tau}}(1,0)
 \bigr)[\ast H_{t=\tau}]
 \lrarr
 \iota_{\dagger}
 \bigl(
 \nbigt
 \bigr)\otimes\newTate(-1)
 \lrarr 0.
\end{equation}
If $\nbigt$ is equipped with a real structure,
{\rm (\ref{eq;22.7.26.3})} is compatible with
 the induced real structures.
\end{lem}
\pf
Let $\iota_0:\cnum_s\times\cnum^m\times X\lrarr
\cnum^2_{s,\sigma}\times\cnum^m\times X$
denote the embedding defined by
$\iota_0(s,z_1,\ldots,z_m,x)=(s,0,z_1,\ldots,z_m,x)$.
Let $H_{\sigma=0}$ denote the image.
We have the following exact sequence
in $\MTMtilde^{\integral}(\cnum^2_{s,\sigma}\times\cnum^m\times X)$:
\begin{equation}
\label{eq;22.7.26.4}
 0\lrarr
 \nbigt\boxtimes\nbigu_{\cnum_{\sigma}}(1,0)
 \lrarr
 \nbigt\boxtimes\bigl(
 \nbigu_{\cnum_{\sigma}}(1,0)[\ast H_{\sigma=0}]
 \bigr)
 \lrarr
 \iota_{0\dagger}
 \bigl(
 \nbigt
 \bigr)\otimes\newTate(-1)
 \lrarr 0.
\end{equation}
We obtain (\ref{eq;22.7.26.3})
from (\ref{eq;22.7.26.4})
and the isomorphism
$\cnum^2_{t,\tau}\simeq
\cnum^2_{s,\sigma}$
defined by
$\rho(t,\tau)
=\rho(t,\tau-t)$.
\hfill\qed

\subsection{Partial Fourier transforms}
\label{subsection;21.7.16.10}
Let $p_{\tau}:\cnum^2_{t,\tau}\times X\to\cnum_{\tau}\times X$
denote the projection.
For $\nbigt\in\MTMtilde^{\integral}(\cnum_t\times X)$,
we define
\[
\FT_{\pm}(\nbigt):=
  \lefttop{T}(p_{\tau})^0_{\ast}
 \Bigl(
 \bigl(
 \nbigu_{\cnum_{\tau}}(1,0)\boxtimes \nbigt
 \bigr)
 \otimes
 \nbigt_{\sm}(\pm t\tau)
 \Bigr)
\in \MTMtilde(\cnum_{\tau}\times X).
\]
By Lemma \ref{lem;21.7.5.1},
we have
$\FT_{\pm}(\nbigt)\simeq
 \lefttop{T}(p_{\tau})^0_{!}
 \Bigl(
 \bigl(
 \nbigu_{\cnum_{\tau}}(1,0)\boxtimes \nbigt
 \bigr)
 \otimes
 \nbigt_{\sm}(\pm t\tau)
 \Bigr)$.

\subsubsection{Inversion formula}

\begin{prop}
\label{prop;21.7.12.12}
There exist natural isomorphisms
$\FT_{\mp}\circ\FT_{\pm}(\nbigt)
 \simeq\nbigt\otimes\newTate(-1)$.
\end{prop}
\pf
Let us study
$\FT_-\circ\FT_+(\nbigt)$.
The other case can be argued similarly.
Let $p_{\sigma}:\cnum^3_{t,\tau,\sigma}\times X\to\cnum_{\sigma}\times X$
denote the projection.
We have
\[
 \FT_{-}\circ\FT_+(\nbigt)
 =
 \lefttop{T}(p_{\sigma})_{\ast}^0\Bigl(
  \bigl(\nbigt\boxtimes
   \nbigu_{\cnum_{\tau}\times\cnum_{\sigma}}(2,0)\bigr)
  \otimes\nbigt_{\sm}((t-\sigma)\tau)
 \Bigr).
\]
Let $H^{(1)}_{t_1=t_2}$ denote the hypersurface of
$\cnum^4_{t_1,t_2,\tau,\sigma}\times X$
defined by $t_1=t_2$.
We consider the following complex
in $\MTMtilde^{\integral}(\cnum^4_{t_1,t_2,\tau,\sigma}\times X)$:
\begin{equation}
\label{eq;21.7.11.2}
 \bigl(
 \nbigt\boxtimes\nbigt((t_2-\sigma)\tau)
 \bigr)
 \otimes\newTate(1)
\lrarr 
\bigl(
 \nbigt
 \boxtimes\nbigt((t_2-\sigma)\tau)
 \bigr)[\ast H^{(1)}_{t_1=t_2}]
 \otimes\newTate(1).
\end{equation}
Here, the second term sits in the degree $0$.
Let $\Delta^{(1)}_{t_1,t_2}:
\cnum^3_{t,\tau,\sigma}\times X
\lrarr
\cnum^4_{t_1,t_2,\tau,\sigma}\times X$
denote the morphism
defined by $\Delta^{(1)}_{t_1,t_2}(t,\tau,\sigma,x)
=(t,t,\tau,\sigma,x)$.
By Lemma \ref{lem;22.7.26.10},
the complex (\ref{eq;21.7.11.2})
is quasi-isomorphic to
\begin{equation}
\label{eq;22.7.26.13}
 (\Delta^{(1)}_{t_1,t_2})_{\dagger}\Bigl(
  \bigl(\nbigt\boxtimes
  \nbigu_{\cnum_{\tau}\times\cnum_{\sigma}}(2,0)
  \bigr)
  \otimes\nbigt_{\sm}((t-\sigma)\tau)
 \Bigr).
\end{equation}

Let $p_{t_1,t_2,\sigma}:\cnum^4_{t_1,t_2,\tau,\sigma}\times X
\lrarr
\cnum^3_{t_1,t_2,\sigma}\times X$
denote the projection.
Let $\Delta^{(2)}_{t_2,\sigma}:\cnum^2_{t_1,\sigma}\times X
 \lrarr
 \cnum^3_{t_1,t_2,\sigma}\times X$
 be the embedding defined by
 $\Delta^{(2)}_{t_2,\sigma}(t_1,\sigma,x)=(t_1,\sigma,\sigma,x)$.
By Lemma \ref{lem;22.7.26.11},
we obtain
\[
 \lefttop{T}(p_{t_1,t_2,\sigma})_{\ast}^0
 \Bigl(
 \bigl(
  \nbigt\boxtimes
  \nbigt((t_2-\sigma)\tau)
 \bigr)
 \otimes\newTate(1)
 \Bigr)
\simeq
 (\Delta^{(2)}_{t_2,\sigma})_{\dagger}\bigl(
\nbigt\boxtimes
 \nbigu_{\cnum_{\sigma}}(1,0)
 \bigr).
\]
Let $H^{(2)}_{t_1=\sigma}$ be the hypersurface of
$\cnum^2_{t_1,\sigma}\times X$
defined by $t_1=\sigma$.
Applying $\lefttop{T}(p_{t_1,t_2,\sigma})^0_{\ast}$
to the complex (\ref{eq;21.7.11.2}),
we obtain
\begin{equation}
\label{eq;21.7.11.4}
 (\Delta^{(2)}_{t_2,\sigma})_{\dagger}
 \bigl(
 \nbigt\boxtimes
 \nbigu_{\cnum_{\sigma}}(1,0)
\bigr)
 \lrarr
 (\Delta^{(2)}_{t_2,\sigma})_{\dagger}
  \Bigl(
\bigl(
 \nbigt\boxtimes
 \nbigu_{\cnum_{\sigma}}(1,0)
 \bigr)
 [\ast  H^{(2)}_{t_1=\sigma}]
 \Bigr).
\end{equation}
The push-forward of (\ref{eq;21.7.11.4})
by the morphism forgetting $t_2$
is naturally identified with
the following complex 
in $\MTMtilde^{\integral}(\cnum^2_{t_1,\sigma}\times X)$:
\begin{equation}
\label{eq;21.7.11.5}
 \nbigt\boxtimes
 \nbigu_{\cnum_{\sigma}}(1,0)
 \lrarr
\bigl(
 \nbigt\boxtimes
 \nbigu_{\cnum_{\sigma}}(1,0)
 \bigr)
 [\ast H^{(2)}_{t_1=\sigma}].
\end{equation}
Let $\Delta^{(3)}_{t_1,\sigma}:\cnum_{\sigma}\times X\lrarr
\cnum^2_{t_1,\sigma}\times X$
denote the morphism induced by
the diagonal embedding
$\cnum_{\sigma}\to\cnum^2_{t_1,\sigma}$.
By Lemma \ref{lem;22.7.26.11},
the complex (\ref{eq;21.7.11.5})
is naturally quasi-isomorphic to
$(\Delta^{(3)}_{t_1,\sigma})_{\dagger}\bigl(
 \nbigt\otimes\newTate(-1)
 \bigr)$.
Because the push-forward of (\ref{eq;22.7.26.13})
by the projection forgetting $(t_2,\tau)$
is quasi-isomorphic to
$(\Delta^{(3)}_{t_1,\sigma})_{\dagger}\bigl(
\FT_-\circ\FT_+(\nbigt)
\bigr)$,
we obtain the claim of Proposition
\ref{prop;21.7.12.12}.
\hfill\qed

\subsubsection{Some subcategories}
\label{subsection;21.7.16.11}

Let $\star\in\{!,\ast\}$.
Let $\MTMtilde^{\integral}(\cnum\times X,[\star 0])
\subset\MTMtilde^{\integral}(\cnum\times X)$
denote the full subcategory of
the objects $\nbigt\in\MTMtilde^{\integral}(\cnum\times X)$
such that $\nbigt[\star (0\times X)]\simeq\nbigt$.
Let 
$\MTMtilde^{\integral}(\cnum\times X)_{\star}
\subset\MTMtilde^{\integral}(\cnum\times X)$
denote the full subcategory of
the objects $\nbigt\in\MTMtilde^{\integral}(\cnum\times X)$
such that
$\lefttop{T}(\pi_{X})_{\star}^j(\nbigt)=0$ $(j\in\seisuu)$,
where $\pi_X:\cnum\times X\to X$ denotes the projection.

\begin{prop}
\label{prop;21.7.12.20}
$\FT_{\pm}$ induce equivalences
$\MTMtilde^{\integral}(\cnum_t\times X,[\star 0])
 \simeq
 \MTMtilde^{\integral}(\cnum_t\times X)_{\star}$.
\end{prop}
\pf
It follows from Proposition \ref{prop;21.7.11.10}.
\hfill\qed

\subsubsection{Functoriality}

Let $F:X\lrarr Y$ be any projective morphism of complex manifolds.
We obtain the following proposition as in the case of 
Proposition \ref{prop;21.7.14.10}
and Corollary \ref{cor;21.7.14.11}.
\begin{prop}
For $\nbigt\in\MTMtilde^{\integral}(\cnum\times X)$,
we have the following natural isomorphisms:
\[
 \FT_{\pm}\Bigl(
 F_{\dagger}^j(\nbigt)
 \Bigr)
 \simeq
 F_{\dagger}^j\FT_{\pm}(\nbigt).
\] 
We obtain 
$F_{\dagger}^j:
\MTMtilde^{\integral}(\cnum\times X,[\star 0])
\lrarr
\MTMtilde^{\integral}(\cnum\times Y,[\star 0])$
and
$F_{\dagger}^j:
\MTMtilde^{\integral}(\cnum\times X)_{\star}
\lrarr
\MTMtilde^{\integral}(\cnum\times Y)_{\star}$
for $\star=!,\ast$.
\hfill\qed
\end{prop}

We obtain the following proposition
as in the case of Proposition \ref{prop;21.7.13.20}.
\begin{prop}
\label{prop;21.7.16.12}
Let $\nbigt\in\MTMtilde^{\integral}(\cnum\times X)$.
\begin{itemize}
 \item For any hypersurface $H$ of $X$,
       there exist natural isomorphisms
\[
       \FT_{\pm}(\nbigt[\star (\proj^1\times H)])
       \simeq
       \FT_{\pm}(\nbigt)[\star (\proj^1\times H)].
\]

 \item 
       Let $f$ and $\ftilde$ be as in Proposition
       {\rm\ref{prop;21.7.13.20}}.
       There exist natural isomorphisms
       $\FT_{\pm}\Pi^{a,b}_{\ftilde,\star}\bigl(
       \nbigt
       \bigr)
       \simeq
       \Pi^{a,b}_{\ftilde,\star}\bigl(
       \FT_{\pm}\nbigt
       \bigr)$ $(\star=!,\ast)$
       and 
       $\FT_{\pm}\Pi^{a,b}_{\ftilde,\ast !}\bigl(
       \nbigt
       \bigr)
       \simeq
       \Pi^{a,b}_{\ftilde,\ast !}\bigl(
       \FT_{\pm}\nbigt
       \bigr)$.
       In particular,
       we have natural isomorphisms
       $\FT_{\pm}\Xi^{(a)}_{\ftilde}(\nbigt)
       \simeq
       \Xi^{(a)}_{\ftilde}(\FT_{\pm}(\nbigt))$
       and 
       $\FT_{\pm}\psi^{(a)}_{\ftilde}(\nbigt)
       \simeq
       \psi^{(a)}_{\ftilde}(\FT_{\pm}(\nbigt))$.
       We also obtain
       $\FT_{\pm}\phi^{(0)}_{\ftilde}(\nbigt)
       \simeq
       \phi^{(0)}_{\ftilde}\FT_{\pm}(\nbigt)$.
       \hfill\qed
\end{itemize}
\end{prop}

\begin{cor}
 In the situation of {\rm\S\ref{subsection;21.6.22.2}}
and {\rm\S\ref{subsection;21.6.22.20}},
for any $\nbigt\in\MTMtilde^{\integral}(\cnum\times Y)$,
there exist natural isomorphisms 
$\FT_{\pm}(\lefttop{T}(\id\times f)^{\star})^i(\nbigt)
 \simeq
 (\lefttop{T}(\id\times f)^{\star})^i\FT_{\pm}(\nbigt)$.
\hfill\qed
\end{cor}

\begin{prop}
\label{prop;21.7.15.3}
 For $\nbigt\in\MTMtilde^{\integral}(\cnum\times X)$,
there exist natural isomorphisms
$\FT_{\pm}(\DD\nbigt)\otimes\newTate(1)
\simeq
 \DD(\FT_{\mp}(\nbigt))$.
There also exist natural isomorphisms:
\[
 \FT_{\pm}(\nbigt^{\ast})
 \otimes\newTate(1)
 \simeq
 \FT_{\pm}(\nbigt)^{\ast},
 \quad
 \FT_{\pm}(j^{\ast}\nbigt)
 \simeq
 j^{\ast}\FT_{\mp}(\nbigt),
 \quad
 \FT_{\pm}(\gammatilde^{\ast}\nbigt)
 \simeq
 \gammatilde^{\ast}\FT_{\pm}(\nbigt).
\]
\end{prop}
\pf
We use the notation in \S\ref{subsection;21.7.25.2}.
As in \cite[Proposition 13.3.5,
Proposition 13.3.6, Proposition 13.3.9]{Mochizuki-MTM},
we have
\begin{multline}
 \DD\Bigl(
 \bigl(
 \nbigt\boxtimes\nbigu_{\cnum}(1,0)
 \otimes\nbigt_{\sm}(t\tau)\bigr)[\ast H_{\infty,t,X}]
 \Bigr)
 \simeq
 \DD\Bigl(
 \nbigt\boxtimes\nbigu_{\cnum}(1,0)
 \otimes\nbigt_{\sm}(t\tau)
 \Bigr)[! H_{\infty,t,X}] 
 \simeq \\
 \Bigl(
 \DD\bigl(\nbigt\boxtimes\nbigu_{\cnum}(1,0)
 \bigr)
 \otimes\nbigt_{\sm}(-t\tau)
 \Bigr)[! H_{\infty,t,X}] 
 \simeq
 \Bigl(
 \DD(\nbigt)\boxtimes
 \DD(\nbigu_{\cnum}(1,0))
 \otimes\nbigt_{\sm}(-t\tau)
\Bigr)[! H_{\infty,t,X}] 
 \\
 \simeq
 \Bigl(
 \DD(\nbigt)\boxtimes
  \nbigu_{\cnum}(1,0)
  \otimes\nbigt_{\sm}(-t\tau)
 \Bigr)[! H_{\infty,t,X}]
  \otimes\newTate(1).
\end{multline}
Then, we obtain the first claim.
We obtain the other isomorphisms similarly.
\hfill\qed

\subsubsection{Convolution and tensor product}

We use the notation in \S\ref{subsection;21.7.14.30}.
For $\nbigt_i\in\MTMtilde^{\integral}(\cnum_{t_i}\times X_i)$
$(i=1,2)$,
we obtain
$\nbigt_1\boxtimes\nbigt_2
\in \MTMtilde^{\integral}(\cnum^2_{t_1,t_2}\times X)$.
We define
\begin{equation}
\nbigh^j\Bigl(
 \nbigt_1\langle\boxtimes,\otimes\rangle^{\star}
 \nbigt_2
 \Bigr)
:=
 \lefttop{T}\bigl(
 (\Delta\times\id_{X})^{\star}
 \bigr)^j
 (\nbigt_1\boxtimes\nbigt_2)
\in\MTMtilde^{\integral}(\cnum\times X),
\end{equation}
\begin{equation}
 \nbigh^j\Bigl(
 \nbigt_1\langle\boxtimes,\ast\rangle_{\star}
 \nbigt_2
 \Bigr)
:=
\lefttop{T}\mu_{\star}^j
 (\nbigt_1\boxtimes\nbigt_2)
\in\MTMtilde^{\integral}(\cnum\times X).
\end{equation}

\begin{prop}
\label{prop;21.7.14.40}
 There exist the following natural isomorphisms.
\begin{equation}
\label{eq;21.7.14.30}
 \nbigh^j\Bigl(
 \FT_{\pm}(\nbigt_1)\langle\boxtimes,\otimes\rangle^!
 \FT_{\pm}(\nbigt_2)
 \Bigr)
 \simeq
 \FT_{\pm}
\Bigl(
 \nbigh^{j-1}\bigl(
 \nbigt_1\langle\boxtimes,\ast\rangle_{\ast}
 \nbigt_2
 \bigr)
 \Bigr)
 \otimes\newTate(-1),
\end{equation}
\begin{equation}
 \nbigh^j\Bigl(
 \FT_{\pm}(\nbigt_1)\langle\boxtimes,\otimes\rangle^{\ast}
 \FT_{\pm}(\nbigt_2)
 \Bigr)
 \simeq
 \FT_{\pm}
 \Bigl(
 \nbigh^{j+1} 
 \bigl(
 \nbigt_1\langle\boxtimes,\ast\rangle_{!}
 \nbigt_2
 \bigr)
 \Bigr).
\end{equation}
\end{prop}
\pf
We explain an outline of the proof for
(\ref{eq;21.7.14.30}) and $\FT_+$:
\[
 \nbigttilde:=
\bigl(
\nbigt_1\boxtimes
\nbigt_2\boxtimes\nbigu_{\cnum^2_{\tau_1,\tau_2}}(2,0)
\bigr)
\otimes\nbigt_{\sm}(t_1\tau_1+t_2\tau_2)
\in
\MTMtilde^{\integral}(
\cnum^4_{t_1,t_2,\tau_1,\tau_2}\times X).
\]
Let $H^{(1)}_{\tau_1=\tau_2}$ be the hypersurface of
$\cnum^4_{t_1,t_2,\tau_1,\tau_2}\times X$
determined by $\tau_1=\tau_2$.
We obtain the following complex
in
$\MTMtilde^{\integral}(\cnum^4_{t_1,t_2,\tau_1,\tau_2}\times X)$:
\begin{equation}
 \label{eq;21.7.14.31}
 \nbigttilde\lrarr\nbigttilde[\ast H^{(1)}_{\tau_1=\tau_2}].
\end{equation}
Here, the first term sits in the degree $0$.
Let $H^{(2)}_{\tau_1=\tau_2}$ denote the hypersurface of
$\cnum^2_{\tau_1,\tau_2}\times X$
determined by $\tau_1=\tau_2$.
We obtain the following complex as the push-forward of
(\ref{eq;21.7.14.31})
by the projection forgetting $(t_1,t_2)$:
\begin{equation}
\label{eq;21.7.14.32}
 \Bigl(
 \FT_{+}(\nbigt_1)\boxtimes
 \FT_+(\nbigt_2)
 \Bigr)
 \lrarr
 \Bigl(
 \FT_{+}(\nbigt_1)\boxtimes
 \FT_+(\nbigt_2)
 \Bigr)[\ast H^{(2)}_{\tau_1=\tau_2}].
\end{equation}
Let $\Delta_2:\cnum_{\tau}\times X\to\cnum^2_{\tau_1,\tau_2}\times X$
denote the morphism induced by the diagonal embedding
$\cnum_{\tau}\to\cnum^2_{\tau_1,\tau_2}$.
The $j$-th cohomology of
(\ref{eq;21.7.14.32})
is quasi-isomorphic to
$\Delta_{2\dagger}\nbigh^j\bigl(
\FT_{+}(\nbigt_1)\langle\boxtimes,\otimes\rangle^!
\FT_+(\nbigt_2)
\bigr)$.

Let $\Delta_1:
\cnum^3_{t_1,t_2,\tau}\times X\to
\cnum^4_{t_1,t_2,\tau_1,\tau_2}\times X$
be the morphism induced by the diagonal embedding
$\cnum_{\tau}\to\cnum^2_{\tau_1,\tau_2}$.
We set
\[
 \nbigttilde^{(1)}:=
 \bigl(\nbigt_1\boxtimes\nbigt_2
 \boxtimes\nbigu_{\cnum_{\tau}}(1,0)\bigr)
 \otimes
 \nbigt_{\sm}((t_1+t_2)\tau)
 \in
 \MTMtilde^{\integral}\bigl(
 \cnum^3_{t_1,t_2,\tau}
 \times X
 \bigr).
\]
By (\ref{eq;22.7.26.1}),
the complex (\ref{eq;21.7.14.31}) is naturally quasi-isomorphic to
$\Delta_{1\dagger}\bigl(
\nbigttilde^{(1)}
\otimes\newTate(-1)
\bigr)[-1]$.
Let $p_1:\cnum^3_{t_1,t_2,\tau}\times X\to\cnum_{\tau}\times X$
denote the projection.
We obtain
\[
 \Delta_{2\dagger}\bigl(
 \lefttop{T}(p_{1})_{\ast}^{j-1}
 \nbigttilde^{(1)}
 \otimes\newTate(-1)
 \bigr)
 \simeq
 \Delta_{2\dagger}
\nbigh^j\bigl(
\FT_{+}(\nbigt_1)\langle\boxtimes,\otimes\rangle^!
\FT_+(\nbigt_2)
\bigr).
\]
Because there exists a natural isomorphism
\[
 \lefttop{T}(p_{1})_{\ast}^{j}
 \nbigttilde^{(1)}
 \simeq
 \FT_{+}\Bigl(
 \nbigh^j\bigl(
 \nbigt_1
 \langle\boxtimes,\ast\rangle_{\ast}
 \nbigt_2
 \bigr)
 \Bigr),
\]
we obtain
(\ref{eq;21.7.14.30})
for $\FT_+$.
The other cases can be argued similarly.
\hfill\qed

\vspace{.1in}
For $\nbigt\in\MTMtilde^{\integral}(\cnum\times X)$,
it is easy to see
\[
 \nbigh^j\Bigl(
 \nbigt\langle\boxtimes,\otimes\rangle^!
 \bigl(\nbigu_{\cnum}(1,0)[\ast 0]\bigr)
 \Bigr)
 \simeq
 \left\{
\begin{array}{ll}
 \nbigt[\ast (0\times X)]\otimes\newTate(-1)& (j=1)\\
 0 & (j\neq 1),
\end{array}
 \right.
\]
\[
 \nbigh^j\Bigl(
 \nbigt\langle\boxtimes,\otimes\rangle^{\ast}
 \bigl(\nbigu_{\cnum}(1,0)[!0]\bigr)
 \Bigr)
 \simeq
 \left\{
\begin{array}{ll}
 \nbigt[!(0\times X)] & (j=-1)\\
 0 & (j\neq -1).
\end{array}
 \right.
\]
By \cite[Proposition 3.22, Corollary 3.23]{Mochizuki-GKZ},
(see the proof of Lemma \ref{lem;22.7.26.11}),
there exist the following isomorphisms
in $\MTMtilde^{\integral}(\cnum,\real)$:
\[
 \FT_{\pm}(\nbigu_{\cnum}(1,0)[\ast 0])
 \simeq
 \nbigu_{\cnum}(1,0)[!0]\otimes\newTate(-1),
 \quad\quad
 \FT_{\pm}(\nbigu_{\cnum}(1,0)[!0])
 \simeq
 \nbigu_{\cnum}(1,0)[\ast 0].
\]
By Proposition \ref{prop;21.7.14.40},
we have
$\nbigh^j\bigl(
\nbigt\langle\boxtimes,\ast\rangle_{\ast}
(\nbigu_{\cnum}(1,0)[!0])
\bigr)=0$
and 
$\nbigh^j\bigl(
\nbigt\langle\boxtimes,\ast\rangle_{!}
(\nbigu_{\cnum}(1,0)[\ast 0])
\bigr)=0$
unless $j=0$.
Moreover,
if $\nbigt\in\MTMtilde^{\integral}(\cnum\times X)_{\ast}$,
we obtain
$\nbigh^0\bigl(
 \nbigt\langle\boxtimes,\ast\rangle_{\ast}
 (\nbigu_{\cnum}(1,0)[!0])
 \bigr)=\nbigt$,
and if 
$\nbigt\in\MTMtilde^{\integral}(\cnum\times X)_{!}$,
we obtain
$\nbigh^0\bigl(
 \nbigt\langle\boxtimes,\ast\rangle_{!}
 (\nbigu_{\cnum}(1,0)[\ast 0])
 \bigr)=\nbigt$.

For $\nbigt\in\MTMtilde^{\integral}(\cnum\times X)$,
we define
\begin{equation}
\label{eq;21.7.16.21}
\left\{
\begin{array}{l}
 \ttP_{\ast}(\nbigt):=
 \nbigh^0\bigl(
 \nbigt\langle\boxtimes,\ast\rangle_{\ast}
 (\nbigu_{\cnum}(1,0)[!0])
 \bigr)
\in\MTMtilde^{\integral}(\cnum\times X)_{\ast},\\
 \ttP_{!}(\nbigt):=
 \nbigh^0\bigl(
 \nbigt\langle\boxtimes,\ast\rangle_{!}
 (\nbigu_{\cnum}(1,0)[\ast 0])
 \bigr)
 \in\MTMtilde^{\integral}(\cnum\times X)_!.
\end{array}
\right.
\end{equation}
For $\star=!,\ast$,
there exist natural isomorphisms
$\ttP_{\star}\circ\ttP_{\star}(\nbigt)\simeq \ttP_{\star}(\nbigt)$.
If $\nbigt\in\MTMtilde^{\integral}(\cnum\times X)_{\star}$,
there exists a natural isomorphism
$\nbigt\simeq\ttP_{\star}(\nbigt)$.

\subsubsection{Real structure}

By the compatibility of $\gammatilde^{\ast}$ and $\FT_{\pm}$
in Proposition \ref{prop;21.7.15.3},
we obtain the functors
\begin{equation}
\label{eq;21.7.12.13}
\FT_{\pm}:
\MTMtilde^{\integral}(\cnum_t\times X,\real)
\lrarr
\MTMtilde^{\integral}(\cnum_{\tau}\times X,\real).
\end{equation}
The isomorphisms in Proposition \ref{prop;21.7.12.12}
are compatible with $\real$-structures.
The functors (\ref{eq;21.7.12.13})
are equivalent.

Let $\star=!,\ast$.
Let 
$\MTMtilde^{\integral}(\cnum\times X,[\star 0],\real)
\subset
\MTMtilde^{\integral}(\cnum\times X,\real)$
denote the full subcategory of objects $\nbigt$
such that $\nbigt[\star (0\times X)]=\nbigt$.
Let 
$\MTMtilde^{\integral}(\cnum\times X,\real)_{\star}
\subset
\MTMtilde^{\integral}(\cnum\times X,\real)$
denote the full subcategory of objects
$\nbigt$ such that
$\lefttop{T}(\pi_{X})^j_{\star}\bigl(
\nbigt\bigr)=0$ $(j\in\seisuu)$,
where $\pi_X:\cnum\times X\to X$ denotes the projection.
Proposition \ref{prop;21.7.12.20}
is enhanced with real structure as follows.
\begin{prop}
\label{prop;21.7.12.21}
$\FT_{\pm}$ induce equivalences
$\MTMtilde^{\integral}(\cnum_t\times X,[\star 0],\real)
 \simeq
 \MTMtilde^{\integral}(\cnum_t\times X,\real)_{\star}$.
\hfill\qed
\end{prop}

\subsubsection{Filtrations}
\label{subsection;21.7.16.100}

Though any object $\nbigt$ of
$\MTMtilde^{\integral}(\cnum\times X,[\star 0])$
is equipped with the weight filtration $W$,
it is also natural to consider the filtration
$\Wtilde^{(1)}(\nbigt)$ determined by
$\Wtilde^{(1)}_j(\nbigt)=
\bigl(W_j(\nbigt)\bigr)[\star 0]$.
The cokernel of
$\Wtilde^{(1)}_j(\nbigt)\lrarr\Wtilde^{(1)}_{k}(\nbigt)$
is naturally an object of
$\MTMtilde^{\integral}(\cnum\times X,[\star 0])$.
Similarly,
for any $\nbigt\in\MTMtilde^{\integral}(\cnum\times X)_{\star}$,
we obtain the filtration
$\Wtilde^{(2)}_j(\nbigt):=
\ttP_{\star}(W_j\nbigt)$.
The following lemma is easy to check.
\begin{prop}
For $\nbigt\in\MTMtilde^{\integral}(\cnum\times X)_{\star}$
such that $\Gr^{\Wtilde^{(2)}}_j(\nbigt)=0$ unless $j=m$,
$\Gr^{\Wtilde^{(1)}}_{j}(\FT_{\pm}(\nbigt))=0$ holds unless $j=m+1$.
For $\nbigt\in\MTMtilde^{\integral}(\cnum\times X,[\star 0])$
such that $\Gr^{\Wtilde^{(1)}}_j(\nbigt)=0$ unless $j=m$,
we obtain $\Gr^{\Wtilde^{(2)}}_{j}(\FT_{\pm}(\nbigt))=0$
unless $j=m+1$.
\hfill\qed
\end{prop}

\subsection{$\cnum^{\ast}$-homogeneity}

Let $\vecn\in\seisuu^2_{\geq 0}$
such that $\gcd(n_1,n_2)=1$ and $n_1\geq n_2$.
Let $\MTMtilde^{\integral}_{\vecn}(\cnum\times X)$
denote the full subcategory of
$(\nbigt,W)=((\nbigm',\nbigm'',C),W)\in\MTMtilde^{\integral}(\cnum\times X)$
such that
$\Upsilon(\nbigm'),\Upsilon(\nbigm'')\in\gbigctilde_{\vecn}(\cnum\times X)$.
We obtain the full subcategories
$\MTMtilde^{\integral}_{\vecn}(\cnum\times X)_{\star}$
and
$\MTMtilde^{\integral}_{\vecn}(\cnum\times X,[\star 0])$
in natural ways.

\begin{prop}
\label{prop;21.7.12.11}
$\FT_{\pm}$ induce equivalences
\[
 \MTMtilde^{\integral}_{\vecn}(\cnum_t\times X)
 \simeq
 \MTMtilde^{\integral}_{\FT(\vecn)}(\cnum_{\tau}\times X),
 \quad
\MTMtilde^{\integral}_{\vecn}(\cnum_t\times X)_{\star}
 \simeq
 \MTMtilde^{\integral}_{\FT(\vecn)}(\cnum_{\tau}\times X,[\star 0]).
\]
\end{prop}
\pf
It follows from Lemma \ref{lem;21.7.11.11}.
\hfill\qed

\begin{cor}
\label{cor;21.7.12.10}
$\FT_{\pm}$ induce equivalences
$\MTMtilde^{\integral}_{(1,0)}(\cnum_t\times X)
 \simeq
 \MTMtilde^{\integral}_{(1,1)}(\cnum_{\tau}\times X)$.
Moreover, we obtain
$\MTMtilde^{\integral}_{(1,0)}(\cnum_t\times X)_{\star}
 \simeq
 \MTMtilde^{\integral}_{(1,1)}(\cnum_{\tau}\times X,[\star 0])$.
\hfill\qed
\end{cor}

\subsubsection{Functoriality}

As in \S\ref{subsection;21.7.15.1},
we mention the functorial properties
which are easy to prove.

\begin{lem}
A projective morphism of complex manifolds
$F:X\lrarr Y$ induces
 $(\id\times F)_{\dagger}^j:
 \MTMtilde^{\integral}_{\vecn}(\cnum\times X)\lrarr
 \MTMtilde^{\integral}_{\vecn}(\cnum\times Y)$,
 $(\id\times F)_{\dagger}^j:
 \MTMtilde^{\integral}_{\vecn}(\cnum\times X)_{\star}
 \lrarr
 \MTMtilde^{\integral}_{\vecn}(\cnum\times Y)_{\star}$
and
 $(\id\times F)_{\dagger}^j:
 \MTMtilde^{\integral}_{\vecn}(\cnum\times X,[\star 0])
 \lrarr
 \MTMtilde^{\integral}_{\vecn}(\cnum\times Y,[\star 0])$.
\hfill\qed
\end{lem}

\begin{lem}
Let $\nbigt\in\MTMtilde^{\integral}_{\vecn}(\cnum\times X)$.
\begin{itemize}
 \item For any hypersurface $H$ of $X$,
$\nbigt[\star (\proj^1\times H)]$ $(\star=!,\ast)$
are also objects of
$\MTMtilde^{\integral}_{\vecn}(\cnum\times X)$.
 \item Let $f$ and $\ftilde$ be as in Lemma {\rm\ref{lem;21.7.15.2}}.
       Then, $\Pi^{a,b}_{\ftilde,\star}(\nbigt)$
       and
       $\Pi^{a,b}_{\ftilde,\ast!}(\nbigt)$
       are objects of
       $\MTMtilde^{\integral}_{\vecn}(\cnum\times X)$.
       As a result,
       $\Xi^{(a)}_{\ftilde}(\nbigt)$,
       $\psi^{(a)}_{\ftilde}(\nbigt)$
       and $\phi^{(0)}_{\ftilde}(\nbigt)$
       are objects of
       $\MTMtilde^{\integral}_{\vecn}(\cnum\times X)$.
       \hfill\qed
\end{itemize}
\end{lem}

\begin{cor}
In the situation of {\rm\S\ref{subsection;21.6.22.2}}
and {\rm\S\ref{subsection;21.6.22.20}},
for any $\nbigt\in\MTMtilde^{\integral}_{\vecn}(\cnum\times Y)$,
$(\lefttop{T}(\id\times f)^{\star})^i(\nbigt)$
are objects of
$\MTMtilde^{\integral}_{\vecn}(\cnum\times X)$.
\hfill\qed  
\end{cor}

\begin{lem}
For $\nbigt\in\MTMtilde^{\integral}_{\vecn}(\cnum\times X)$,
$j^{\ast}(\nbigt)$
and $\gammatilde^{\ast}(\nbigt)$
are objects of
$\MTMtilde^{\integral}_{\vecn}(\cnum\times X)$.
Let $\star=!,\ast$.
If $\nbigt\in\MTMtilde_{\vecn}(\cnum\times X,[\star 0])$,
$j^{\ast}(\nbigt)$
and $\gammatilde^{\ast}(\nbigt)$
are objects of
$\MTMtilde^{\integral}_{\vecn}(\cnum\times,[\star 0])$.
If $\nbigt\in\MTMtilde^{\integral}_{\vecn}(\cnum\times X)_{\star}$,
$j^{\ast}(\nbigt)$ and $\gammatilde^{\ast}(\nbigt)$ are
objects of
$\MTMtilde^{\integral}_{\vecn}(\cnum\times X)_{\star}$.
\hfill\qed
\end{lem}

\begin{lem}
For $\nbigt\in\MTMtilde^{\integral}_{\vecn}(\cnum\times X)$,
$\nbigt^{\ast}$,
and $\DD(\nbigt)$
are objects of
$\MTMtilde^{\integral}_{\vecn}(\cnum\times X)$.
Let $(\star_1,\star_2)=(!,\ast),(\ast,!)$.
If $\nbigt\in\MTMtilde^{\integral}_{\vecn}(\cnum\times X,[\star_1 0])$,
$\nbigt^{\ast}$ and
$\DD\nbigt$ are objects of
$\MTMtilde^{\integral}_{\vecn}(\cnum\times,[\star_20])$.
If $\nbigt\in\MTMtilde^{\integral}_{\vecn}(\cnum\times X)_{\star_1}$,
$\nbigt^{\ast}$ and $\DD\nbigt$ are objects of
$\MTMtilde^{\integral}_{\vecn}(\cnum\times X)_{\star_2}$.
\hfill\qed
\end{lem}

We set $X=X_1\times X_2$.
\begin{lem}
For $\nbigt_i\in\MTMtilde^{\integral}_{\vecn}(\cnum\times X_i)$,
$\nbigh^j\bigl(
 \nbigt_1\langle\boxtimes,\otimes\rangle^{\star}
 \nbigt_2
\bigr)$
and 
$\nbigh^j\bigl(
 \nbigt_1\langle\boxtimes,\ast\rangle_{\star}
 \nbigt_2
\bigr)$
 are objects of
$\MTMtilde^{\integral}_{\vecn}(\cnum\times X)$.
\hfill\qed
\end{lem}

\begin{lem}
For any $\nbigt_i\in \MTMtilde^{\integral}_{(1,1)}(\cnum\times X_i,[\ast 0])$,
we have the vanishing
$\nbigh^j(\nbigt_1\langle\boxtimes,\otimes\rangle^!\nbigt_2)=0$
unless $j=1$,
and 
$\nbigh^1(\nbigt_1\langle\boxtimes,\otimes\rangle^!\nbigt_2)$
is an object of
$\MTMtilde^{\integral}_{(1,1)}(\cnum\times X,[\ast 0])$.
For any $\nbigt_i\in \MTMtilde^{\integral}_{(1,1)}(\cnum\times X_i,[!0])$,
we have the vanishing
$\nbigh^j(\nbigt_1\langle\boxtimes,\otimes\rangle^{\ast}\nbigt_2)=0$
unless $j=-1$,
$\nbigh^{-1}(\nbigt_1\langle\boxtimes,\otimes\rangle^{\ast}\nbigt_2)$
is an object of
$\MTMtilde^{\integral}_{(1,1)}(\cnum\times X,[! 0])$.
\hfill\qed
\end{lem}

\subsubsection{Real structure}
\label{subsection;21.7.16.30}

Let
$\MTMtilde^{\integral}_{\vecn}(\cnum\times X,\real)
\subset
\MTMtilde^{\integral}(\cnum\times X,\real)$
denote the full subcategory of
the objects
$(\nbigt,W,\kappa)$
such that $(\nbigt,W)\in\MTMtilde^{\integral}_{\vecn}(\cnum\times X)$.
For $(\nbigt,W,\kappa)=((\nbigm',\nbigm'',C),W,\kappa)$,
$\nbigm'$ is isomorphic to
$j^{\ast}\DD(\nbigm'')$.
Hence, the condition is equivalent to
that $\Upsilon(\nbigm')\in\gbigctilde_{\vecn}(\cnum\times X)$.

We obtain
$\MTMtilde^{\integral}_{\vecn}(\cnum\times X,[\star 0],\real)
=\MTMtilde^{\integral}_{\vecn}(\cnum\times X,\real)
\cap
\MTMtilde^{\integral}(\cnum\times X,[\star 0],\real)$
and
$\MTMtilde^{\integral}_{\vecn}(\cnum\times X,\real)_{\star}
=\MTMtilde^{\integral}_{\vecn}(\cnum\times X,\real)
\cap
\MTMtilde^{\integral}(\cnum\times X,\real)_{\star}$.
Proposition \ref{prop;21.7.12.11}
and Corollary \ref{cor;21.7.12.10}
are enhanced as follows.
\begin{prop}
 \label{prop;21.7.12.22}
$\FT_{\pm}$ induce equivalences
\[
 \MTMtilde^{\integral}_{\vecn}(\cnum_t\times X,\real)
 \simeq
 \MTMtilde^{\integral}_{\FT(\vecn)}(\cnum_{\tau}\times X,\real),
 \quad
\MTMtilde^{\integral}_{\vecn}(\cnum_t\times X,\real)_{\star}
 \simeq
 \MTMtilde^{\integral}_{\FT(\vecn)}(\cnum_{\tau}\times X,[\star 0],\real).
\]
In particular,
we obtain
\[
 \MTMtilde^{\integral}_{(1,0)}(\cnum_t\times X,\real)
 \simeq
 \MTMtilde^{\integral}_{(1,1)}(\cnum_{\tau}\times X,\real),
\quad
 \MTMtilde^{\integral}_{(1,0)}(\cnum_t\times X,\real)_{\star}
 \simeq
 \MTMtilde^{\integral}_{(1,1)}(\cnum_{\tau}\times X,[\star 0],\real).
\]
\hfill\qed
\end{prop}

\subsection{Rescalability and exponential $\real$-Hodge modules}

\subsubsection{Rescalable mixed twistor $\nbigd$-modules}
\label{subsection;21.7.16.51}

Let $\nbigt\in\MTMtilde^{\integral}_{(1,1)}(\cnum\times X,\real)$.
Let $\iota_1:\{1\}\times X\lrarr \proj^1\times X$
denote the inclusion.
There exists 
$\ttU(\nbigt)\in\MTM^{\integral}(X,\real)$
such that
$\iota_{1\dagger}\ttU(\nbigt)
=\psi^{(1)}_{t-1}(\nbigt)$.
Because $\iota_1$ is strictly non-characteristic for $\nbigt$,
we obtain the following exact sequence:
\begin{equation}
\label{eq;21.7.15.10}
 0\lrarr
 \iota_{1\dagger}\ttU(\nbigt)
 \lrarr
 \nbigt[!(1\times X)]
 \lrarr
 \nbigt\lrarr 0.
\end{equation}
We also have
\begin{equation}
\label{eq;21.7.15.11}
 0\lrarr
\nbigt
 \lrarr
 \nbigt[\ast(1\times X)]
 \lrarr
 \iota_{1\dagger}\ttU(\nbigt)
 \otimes\newTate(-1)
 \lrarr 0.
\end{equation}

We obtain the functor
$\ttU:\MTMtilde^{\integral}_{(1,1)}(\cnum\times X,\real)
\lrarr\MTM^{\integral}(X,\real)$.
As the restriction to full subcategories,
we obtain 
\begin{equation}
\label{eq;21.7.11.30}
 \ttU:
\MTMtilde^{\integral}_{(1,1)}(\cnum\times X,[\star 0],\real)
\lrarr
\MTM^{\integral}(X,\real).
\end{equation}

\begin{prop}
\label{prop;21.7.12.30}
The functor {\rm(\ref{eq;21.7.11.30})} is
fully faithful.
\end{prop}
\pf
By Proposition \ref{prop;21.7.11.31},
we obtain that the functor (\ref{eq;21.7.11.30})
is faithful.
Let us prove that (\ref{eq;21.7.11.30}) is full.
Let $\nbigt_i\in\MTMtilde^{\integral}_{(1,1)}
(\cnum\times X,[\star 0],\real)$ $(i=1,2)$.
There exist $Z_i\subset X$
such that $\proj^1\times Z_i$ are the supports of $\nbigt_i$.
We set $d:=\max\{\dim Z_1,\dim Z_2\}$.
We use an induction on $d$.
If $d=0$, the claim is obvious.
Suppose that we have already proved the claim in the case $d-1$.
It is enough to prove the claim locally around any point
$P$ of $X$.
Let $g:\ttU(\nbigt_1)
\lrarr\ttU(\nbigt_2)$
be a morphism in
$\MTM^{\integral}(X,\real)$.
Let $\nbigt_i=(\nbigm'_i,\nbigm''_i,C_i)$.
Note that
\[
\ttU(\nbigt_i)
=\Bigl(
 \lambda^{-1}\iota_1^{\ast}(\nbigm'_i),\,
 \iota_1^{\ast}(\nbigm''_i),\Ctilde_i
 \Bigr),
\]
where $\Ctilde_i$ denote the induced sesqui-linear pairing.
We have the morphisms
$g':
\lambda^{-1}\iota_1^{\ast}(\nbigm'_2)\lrarr
\lambda^{-1}\iota_1^{\ast}(\nbigm'_1)$
and 
$g'':
\iota_1^{\ast}(\nbigm''_1)\lrarr
\iota_1^{\ast}(\nbigm''_2)$
which determine $g$.
By Proposition \ref{prop;21.7.11.31},
we have the unique morphisms
$\gtilde':\nbigm'_2\lrarr\nbigm'_1$
and
$\gtilde'':\nbigm''_1\lrarr\nbigm''_2$
which induce $g'$ and $g''$, respectively.
We have only to check that
they define a morphism
$\nbigt_1\lrarr\nbigt_2$
compatible with the $\real$-structures.
By the standard argument,
we can reduce the issue
to the study the case where
$\nbigm'_i$ and $\nbigm''_i$ are
locally free $\nbigo_{\nbigx}$-modules,
which we can easily check.
\hfill\qed

\begin{df}
\label{df;21.7.16.50}
Let
$\MTM^{\integral}_{\res}(X,\real)$
denote the essential image of
the functor {\rm(\ref{eq;21.7.11.30})},
which is independent of $\star=!,\ast$.
Objects of $\MTM^{\integral}_{\res}(X,\real)$
are called rescalable mixed twistor $\nbigd$-modules
on $X$.
\hfill\qed
\end{df}

We state some functorial properties
which are easy to prove.

\begin{prop}
Let $F:X\lrarr Y$ be a projective morphism of complex manifolds.
For 
$\nbigt\in\MTMtilde^{\integral}_{(1,1)}(\cnum\times X,[\star 0],\real)$,
there exist natural isomorphisms
$F_{\dagger}^j\bigl(
 \ttU(\nbigt)
 \bigr)
 \simeq
 \ttU((\id\times F)_{\dagger}^j\nbigt)$
in $\MTMtilde^{\integral}(Y,\real)$.
In particular, we obtain
$F_{\dagger}^j:
\MTM^{\integral}_{\res}(X,\real)
\lrarr
\MTM^{\integral}_{\res}(Y,\real)$.
\hfill\qed
\end{prop}

\begin{prop}
Let $H$ be any hypersurface of $X$.
For any
$\nbigt\in\MTMtilde^{\integral}_{(1,1)}(\cnum\times X,[\star 0],\real)$,
we have
\[
 \ttU
 \bigl(
  (\nbigt[\star' (\proj^1\times H)])
 \bigr)
=\ttU
 \bigl(
  (\nbigt[\star' (\proj^1\times H)])[\star(0\times X)]
 \bigr)
 \simeq
  \ttU(\nbigt)[\star'H]
 \quad (\star'=!,\ast).
\] 
In particular,
for any $\nbigt_1\in\MTM^{\integral}_{\res}(X,\real)$,
$\nbigt_1[\star'H]$ $(\star'=!,\ast)$
are objects of 
$\MTM^{\integral}_{\res}(X,\real)$.
\hfill\qed
\end{prop}

\begin{prop}
Let $f$ and $\ftilde$ be as in Proposition {\rm\ref{prop;21.7.13.20}}.
Let $\nbigt\in\MTMtilde^{\integral}_{(1,1)}(\cnum\times X,[\star 0],\real)$.
Then, there exist natural isomorphisms
\[
 \ttU
 \bigl(
 \Pi^{a,b}_{\ftilde,\star'}(\nbigt)[\star(0\times X)]
 \bigr)
 \simeq
 \ttU
 \bigl(
 \Pi^{a,b}_{\ftilde,\star'}(\nbigt)
 \bigr)
 \simeq
 \Pi^{a,b}_{f,\star'}
 \ttU(\nbigt),
\]
\[
 \ttU
 \bigl(
 \Pi^{a,b}_{\ftilde,\ast !}(\nbigt)[\star(0\times X)]
 \bigr)
 \simeq
 \ttU
 \bigl(
 \Pi^{a,b}_{\ftilde,\ast !}(\nbigt)
 \bigr)
 \simeq
 \Pi^{a,b}_{f,\ast !}
 \ttU(\nbigt).
\] 
As a result, 
there exist natural isomorphisms
\[
 \ttU\bigl(
 \Xi^{(a)}_{\ftilde}(\nbigt)[\star(0\times X)]
 \bigr)
 \simeq
 \ttU\bigl(
 \Xi^{(a)}_{\ftilde}(\nbigt)
 \bigr)
  \simeq
 \ttU\bigl(
 \Xi^{(a)}_{\ftilde}(\nbigt)
 \bigr)
 \simeq
 \Xi^{(a)}_{f}
 \ttU(\nbigt),
\]
\[
 \ttU\bigl(
 \psi^{(a)}_{\ftilde}(\nbigt)[\star(0\times X)]
 \bigr)
 \simeq
 \ttU\bigl(
 \psi^{(a)}_{\ftilde}(\nbigt)
 \bigr)
 \simeq
 \ttU\bigl(
 \psi^{(a)}_{\ftilde}(\nbigt)
 \bigr)
 \simeq
 \psi^{(a)}_{f}
 \ttU(\nbigt),
\]
\[
 \ttU\bigl(
 \phi^{(0)}_{\ftilde}(\nbigt)[\star(0\times X)]
 \bigr)
 \simeq
 \ttU\bigl(
 \phi^{(0)}_{\ftilde}(\nbigt)
 \bigr) 
 \simeq
 \phi^{(0)}_{f}
 \ttU(\nbigt).
\]
Therefore,
for any $\nbigt_1\in\MTM^{\integral}_{\res}(X,\real)$,
$\Pi^{a,b}_{f,\star}(\nbigt_1)$,
$\Pi^{a,b}_{f,\ast !}(\nbigt_1)$,
$\Xi^{(a)}_f(\nbigt_1)$,
$\psi^{(a)}_f(\nbigt_1)$
and $\phi^{(0)}_f(\nbigt)$
are objects of $\MTM^{\integral}_{\res}(X,\real)$.
\hfill\qed
\end{prop}

\begin{cor}
In the situation of {\rm\S\ref{subsection;21.6.22.2}}
and {\rm\S\ref{subsection;21.6.22.20}},
for any $\nbigt_1\in\MTM^{\integral}_{\res}(Y,\real)$
and for any $\star=!,\ast$,
$(\lefttop{T}f^{\star})^j(\nbigt_1)$
are objects of
$\MTM^{\integral}_{\res}(X,\real)$.
For any $\nbigt\in
\MTMtilde^{\integral}_{(1,1)}(\cnum\times Y,[\star' 0],\real)$,
there exist natural isomorphisms
\[
 \ttU\bigl(
 (\lefttop{T}(\id\times f)^{\star})^j(\nbigt)[\star' (0\times Y)]
 \bigr) 
 \simeq
 \ttU\bigl(
 (\lefttop{T}(\id\times f)^{\star})^j(\nbigt)
 \bigr) 
 \simeq
 (\lefttop{T}f^{\star})^j
 \ttU(\nbigt).
\]
\hfill\qed
\end{cor}

\begin{prop}
Let $\star=!,\ast$.
For $\nbigt\in
\MTMtilde^{\integral}_{(1,1)}(\cnum\times X,[\star 0],\real)$,
there exist isomorphisms
\[
 \ttU
 j^{\ast}(\nbigt)
 \simeq
 j^{\ast}
 \ttU(\nbigt),
 \quad
 \ttU
 \gammatilde^{\ast}(\nbigt)
 \simeq
 \gammatilde^{\ast}
 \ttU(\nbigt).
\]
Therefore,
for any $\nbigt_1\in\MTM^{\integral}_{\res}(X,\real)$,
$j^{\ast}\nbigt_1$
and $\gammatilde^{\ast}\nbigt_1$
are objects of $\MTM^{\integral}_{\res}(X,\real)$.
\hfill\qed
\end{prop}

By using the exact sequences (\ref{eq;21.7.15.10})
and (\ref{eq;21.7.15.11}),
we obtain the following.
\begin{prop}
Let $\star=!,\ast$.
For $\nbigt\in
\MTMtilde^{\integral}_{(1,1)}(\cnum\times X,[\star 0],\real)$,
there exist isomorphisms
\[
 \bigl(
 \ttU(\nbigt)
 \bigr)^{\ast}
 \simeq
 \ttU
 (\nbigt^{\ast})\otimes\newTate(-1),
 \quad
 \DD\bigl(
 \ttU(\nbigt)
 \bigr)
 \simeq
 \ttU\bigl(
 \DD(\nbigt)
\bigr)
\otimes\newTate(-1).
\]
Therefore,
for any $\nbigt_1\in\MTM^{\integral}_{\res}(X,\real)$,
$\nbigt_1^{\ast}$
and $\DD\nbigt_1$
are objects of $\MTM^{\integral}_{\res}(X,\real)$.
\hfill\qed
\end{prop}

Let $X=X_1\times X_2$.
\begin{prop}
\label{prop;21.7.16.60}
For $\nbigt_i\in\MTMtilde^{\integral}_{(1,1)}(\cnum\times X_i,\real)$,
there exist the following natural isomorphisms:
\[
 \ttU\Bigl(
 \nbigh^1\bigl(
 \nbigt_1\langle\boxtimes,\otimes\rangle^!\nbigt_2
 \bigr)
 \Bigr)
 \simeq
 \bigl(
 \ttU(\nbigt_1)\boxtimes\ttU(\nbigt_2)
 \bigr)
 \otimes\newTate(-1),
\quad\quad
 \ttU\Bigl(
 \nbigh^{-1}\bigl(
 \nbigt_1\langle\boxtimes,\otimes\rangle^{\ast}\nbigt_2
 \bigr)
  \Bigr)
 \simeq
 \ttU(\nbigt_1)\boxtimes\ttU(\nbigt_2).
\]
As a result,
for any $\nbigt_i\in\MTM^{\integral}_{\res}(X_i,\real)$,
$\nbigt_1\boxtimes\nbigt_2$
is an object of 
$\MTM^{\integral}_{\res}(X,\real)$. 
\end{prop}
\pf
Let us study the first isomorphism.
Let $H_{t_1=t_2}$ be the hypersurface of
$\cnum^2_{t_1,t_2}\times X$
defined by $t_1=t_2$.
Let $H_{t_i=1}$ be the hypersurfaces of
$\cnum^2_{t_1,t_2}\times X$
defined by $t_i=1$.
We also use the notation in the proof of Lemma \ref{lem;21.4.8.2}.
We set $\Gamma_0:=\{0,1,2\}$,
and let $\gbigk_0:\Gamma_0\lrarr
\nbigs(\cnum^2_{t_1,t_2}\times X)$
be defined by
$\gbigk_0(0)=H_{t_1=t_2}$,
$\gbigk_0(1)=H_{t_1=1}$,
and $\gbigk_0(2)=H_{t_2=1}$.
We put $\Gamma_1=\{0,2\}$
and $\Gamma_2=\{1,2\}$.
Let $\gbigk_i$ $(i=1,2)$
denote the restriction of $\gbigk_0$
to $\Gamma_i$.
We obtain the complexes
$C^{\bullet}(\nbigt_1\boxtimes\nbigt_2,\Gamma_i,\gbigk_i)$
$(i=0,1,2)$
in $\MTMtilde^{\integral}(\cnum^2_{t_1,t_2}\times X,\real)$.
There exist the natural morphisms of complexes
\begin{equation}
\label{eq;21.7.15.20}
\begin{CD}
 C^{\bullet}(\nbigt_1\boxtimes\nbigt_2,\Gamma_1,\gbigk_1)
 @<<<
 C^{\bullet}(\nbigt_1\boxtimes\nbigt_2,\Gamma_0,\gbigk_0)
 @>>>
 C^{\bullet}(\nbigt_1\boxtimes\nbigt_2,\Gamma_2,\gbigk_2).
\end{CD}
\end{equation}
Because
$B(\Gamma_i,\gbigk_i)=\{(1,1)\}\times X$,
the morphisms in (\ref{eq;21.7.15.20})
are quasi-isomorphisms.
It is easy to see that
the $j$-th cohomology of the complexes is $0$
unless $j=2$.
The second cohomology of
$C^{\bullet}(\nbigt_1\boxtimes\nbigt_2,\Gamma_1,\gbigk_1)$
is naturally isomorphic to
$\ttU\bigl(
\nbigh^1(\nbigt_1\langle\boxtimes,\otimes\rangle^!\nbigt_2)
\otimes\newTate(-1)
\bigr)$.
The second cohomology of
$C^{\bullet}(\nbigt_1\boxtimes\nbigt_2,\Gamma_2,\gbigk_2)$
is naturally isomorphic to
$\bigl(
\ttU(\nbigt_1)\otimes\newTate(-1)
\bigr)
\boxtimes
\bigl(
\ttU(\nbigt_2)\otimes\newTate(-1)
\bigr)$.
Thus, we obtain the first isomorphism.
We can obtain the second isomorphism similarly.
\hfill\qed

\subsubsection{Comparison with exponential $\real$-Hodge modules}
\label{subsection;21.7.16.71}

Let $\MHMtilde(\cnum\times X,\real)$
denote the category of mixed $\real$-Hodge modules
on $\cnum\times X$
which are extendable to mixed $\real$-Hodge modules
on $\proj^1\times X$.
A mixed $\real$-Hodge module on $\cnum\times X$
is denoted as a tuple of
a regular holonomic $\nbigd_{\cnum\times X}$-module $M$,
a Hodge filtration $F$ on $M$,
a $\real$-perverse sheaf $P_{\real}$
with $P_{\real}\simeq \DR(M)$,
and a weight filtration $W$ on $P_{\real}$
satisfying the conditions for mixed Hodge modules.
By the Rees construction
(see \cite[\S13.5]{Mochizuki-MTM}),
we obtain the functor
\[
 \ttV:\MHMtilde(\cnum\times X,\real)
 \lrarr
 \MTMtilde^{\integral}_{(1,0)}(\cnum\times X,\real).
\]
According to \cite[Theorem 3.35]{Mochizuki-GKZ},
$\ttV$ is equivalent.
For $\star=!,\ast$,
let
$\MHMtilde(\cnum\times X,\real)_{\star}
\subset
\MHMtilde(\cnum\times X,\real)$
denote the full subcategory
corresponding to
$\MTMtilde^{\integral}_{(1,0)}(\cnum\times X,\real)_{\star}$,
i.e.,
mixed $\real$-Hodge modules $(M,F,P_{\real})$ on $\cnum\times X$
such that
(i) they are extendable across $\{\infty\}\times X$,
(ii) $R\pi_{X\star}(P_{\real})=0$.
We have the equivalence
$\ttV:\MHMtilde(\cnum\times X,\real)_{\star}
\simeq
\MTMtilde^{\integral}_{(1,0)}(\cnum\times X,\real)_{\star}$.

\begin{rem}
\label{rem;21.7.16.70}
As an analogue of exponential mixed Hodge structure in
{\rm\cite{Kontsevich-Soibelman-exponential-Hodge}},
we call
$\MHMtilde(\cnum\times X,\real)_{\ast}$
the category of 
exponential $\real$-Hodge modules on $X$.
\hfill\qed
\end{rem}

We obtain the following functor:
\[
 \ttB_{\pm}:=\ttU\circ\FT_{\pm}\circ\ttV:
 \MHMtilde(\cnum\times X,\real)
 \lrarr
 \MTM^{\integral}_{\res}(X,\real).
\]
The restriction of $\ttB_{\pm}$
to $\MHMtilde(\cnum\times X,\real)_{\star}$
are denoted by $\ttB_{\pm,\star}$.
We obtain the following theorem
as a consequence of Proposition \ref{prop;21.7.12.30}.

\begin{thm}
\label{thm;21.7.16.3}
The functors $\ttB_{\pm,\star}$ are equivalent.
\hfill\qed
\end{thm}

We have the compatibility of
the equivalence in Theorem \ref{thm;21.7.16.3}
with some basic functors.

\begin{prop}
Let $(M,F,P_{\real},W)\in\MHMtilde(\cnum\times X,\real)$.
\begin{itemize}
 \item
      Let $f:X\lrarr Y$ be a projective morphism of complex manifolds.
      There exist natural isomorphism
      $\ttB_{\pm}(\id\times f)_{\dagger}^j(M,F,P_{\real},W)
      \simeq
      f^j_{\dagger}(\ttB_{\pm}(M,F,P_{\real},W))$.
 \item Let $H$ be a hypersurface of $X$.
       For any $\star=!,\ast$,
       there exist natural isomorphisms
\[
       \ttB_{\pm}\bigl(
       (M,F,P_{\real},W)[\star (\proj^1\times H)]
       \bigr)
       \simeq
       \ttB_{\pm}(M,F,P_{\real},W)[\star H]. 
\]
 \item Let $f$ and $\ftilde$ be as in Proposition {\rm\ref{prop;21.7.13.20}}.
       There exist the following natural isomorphisms:
 \[
       \ttB_{\pm}
       \Bigl(
       \Pi^{a,b}_{\ftilde,\star}(M,F,P_{\real},W)
       \Bigr)
       \simeq
       \Pi^{a,b}_{f,\star}
       \ttB_{\pm}(M,F,P_{\real},W) \quad(\star=!,\ast),
\]
\[
       \ttB_{\pm}
       \Bigl(
       \Pi^{a,b}_{\ftilde,\ast!}(M,F,P_{\real},W)
       \Bigr)
       \simeq
       \Pi^{a,b}_{f,\ast !}
       \ttB_{\pm}(M,F,P_{\real},W).
 \]
       As a result,
       there exist the following natural isomorphisms:
\[
       \ttB_{\pm}
        \Xi^{(a)}_{\ftilde}(M,F,P_{\real},W)
       \simeq
       \Xi^{(a)}_{f}\ttB_{\pm}(M,F,P_{\real},W),
\quad
       \ttB_{\pm}
       \psi^{(a)}_{\ftilde}(M,F,P_{\real},W)
       \simeq
       \psi^{(a)}_{f}\ttB_{\pm}(M,F,P_{\real},W),
\]
\[
       \ttB_{\pm}\phi^{(0)}_f(M,F,P_{\real},W)
       \simeq
       \phi^{(0)}_{\ftilde}\ttB_{\pm}(M,F,P_{\real},W).
\]
\hfill\qed
\end{itemize} 
\end{prop}

\begin{cor}
In the situation of {\rm\S\ref{subsection;21.6.22.2}}
and {\rm\S\ref{subsection;21.6.22.20}},
for any $(M,F,P_{\real},W)\in\MHMtilde(\cnum\times X,\real)$,
there exist natural isomorphisms
 $\ttB_{\pm}
 (\lefttop{T}(\id\times f)^{\star})^j(M,F,P_{\real},W)
 \simeq
 \bigl(
 \lefttop{T}f^{\star}
 \bigr)^j
 \ttB_{\pm}(M,F,P_{\real},W)$
 $(\star=!,\ast)$.
\hfill\qed
\end{cor}

\begin{prop}
For any $(M,F,P_{\real},W)\in\MHMtilde(\cnum\times X,\real)$,
there exist natural isomorphisms
\[
\ttB_{\pm}\bigl(
\DDD(M,F,P_{\real},W)
\bigr)
\simeq
\DD(\ttB_{\mp}(M,F,P_{\real},W)).
\]
For $(M_i,F,P_{i\real},W)\in \MHMtilde(\cnum\times X_i)$,
there exist the following natural isomorphisms:
\[
 \ttB_{\pm}\bigl(
 (M_1,F,P_{1\real},W)
 \langle\boxtimes,\ast\rangle_{\star}
 (M_2,F,P_{2\real},W)
 \bigr)
 \simeq
 \ttB_{\pm}
 (M_1,F,P_{1\real},W)
 \boxtimes
 \ttB_{\pm}
 (M_2,F,P_{2\real},W).
\]
\hfill\qed 
\end{prop}

Let $\nbigt=\ttV(M,F,P_{\real},W)
\in\MTMtilde^{\integral}_{(1,0)}(\cnum_t\times X,\real)$.
Let us describe
$\ttU\circ\FT_+(\nbigt)$ more directly.
Let $H_{\tau=1}$ denote the hypersurface of
$\cnum^2_{t,\tau}\times X$
determined by $\tau=1$.
We have the following morphism on
in $\MTMtilde^{\integral}(\cnum^2_{t,\tau}\times X,\real)$:
\begin{equation}
\label{eq;21.7.16.1}
 \Bigl(
 \bigl(
 \nbigt\boxtimes\nbigu_{\cnum_{\tau}}(1,0)
 \bigr)\otimes\nbigt_{\sm}(t\tau)
 \Bigr)[!H_{\tau=1}]
\lrarr
 \bigl(
 \nbigt\boxtimes\nbigu_{\cnum_{\tau}}(1,0)
 \bigr)\otimes\nbigt_{\sm}(t\tau).
\end{equation}
Let $p_{\tau,X}:\cnum^2_{t,\tau}\times X\to\cnum_{\tau}\times X$
denote the projection.
Applying
$\lefttop{T}(p_{\tau,X})_!$
to (\ref{eq;21.7.16.1}),
we obtain the following morphism
in $\MTMtilde^{\integral}(\cnum_{\tau}\times X,\real)$:
\begin{equation}
\label{eq;21.7.16.2}
\FT_+(\nbigt)[!(1\times X)]
\lrarr
\FT_+(\nbigt).
\end{equation}
Let $\iota_{1}:X\lrarr \cnum_{\tau}\times X$
denote the inclusion
defined by $\iota_1(x)=(1,x)$.
Then, $\iota_{1\dagger}\ttU\circ\FT_+(\nbigt)$
is naturally isomorphic to the kernel of
(\ref{eq;21.7.16.2}).
Let $\iota_1^{(1)}:\cnum_t\times X\lrarr
\cnum^2_{t,\tau}\times X$
denote the inclusion defined by
$\iota_1^{(1)}(t,x)=(t,1,x)$.
Then, the kernel of
(\ref{eq;21.7.16.1})
is naturally isomorphic to
$\iota^{(1)}_{1\dagger}\bigl(
 \nbigt\otimes\nbigt_{\sm}(t)
 \bigr)$.
Hence, we obtain
\begin{equation}
\label{eq;21.7.16.4}
 \ttU\circ\FT_+(\nbigt)
 \simeq
\lefttop{T}
 (\pi_{X})^0_!\bigl(
 \nbigt\otimes\nbigt_{\sm}(t)
 \bigr)
 \simeq
 \lefttop{T}
 (\pi_{X})^0_{\ast}\bigl(
 \nbigt\otimes\nbigt_{\sm}(t)
 \bigr).
\end{equation}

\begin{rem}
Any object of $\MTM^{\integral}_{\res}(X,\real)$
is described as {\rm(\ref{eq;21.7.16.4})}
for some $\nbigt=\ttV(M,F,P_{\real},W)$
by Theorem {\rm\ref{thm;21.7.16.3}}.
\hfill\qed
\end{rem}

Let
$\ttP_{\star}:\MHMtilde(\cnum\times X,\real)
\lrarr
\MHMtilde(\cnum\times X,\real)_{\star}$
denote the functor
induced by
$\ttP_{\star}:\MTMtilde^{\integral}_{(1,0)}(\cnum\times X,\real)
\lrarr
\MTMtilde^{\integral}_{(1,0)}(\cnum\times X,\real)_{\star}$
and the equivalence $\ttV$.
We have
$\ttB_{\pm}=\ttB_{\pm,\star}\circ\ttP_{\star}$.
For any $(M,F,P_{\real},W)\in \MHMtilde(\cnum\times X,\real)_{\star}$,
we have the filtration $\Wtilde^{(2)}$
as in \S\ref{subsection;21.7.16.100},
i.e.,
$\Wtilde^{(2)}_j(M,F,P_{\real})
=\ttP_{\star}(W_j(M,F,P_{\real}))$.
The following holds.
\begin{prop}
For $(M,F,P_{\real},W)\in\MHMtilde(\cnum\times X,\real)_{\star}$
such that
$\Gr^{\Wtilde^{(2)}}_j(M,F,P_{\real})=0$ unless $j=m$,
we have
$\Gr^{W}_j\ttB_{\pm,\star}(M,F,P_{\real},W)=0$ unless $j=m$.
\hfill\qed
\end{prop}

\subsubsection{Examples}

Let $\nbigt\in\MTM^{\integral}(X,\real)$
be the integrable mixed twistor $\nbigd$-module
with a real structure,
obtained as the analytification of
the Rees construction of a mixed $\real$-Hodge module.
Let $H$ be a hypersurface of $X$.
Let $f$ be a meromorphic function on $(X,H)$.
We obtain
$\nbigt\otimes\nbigt_{\sm}(f)\in
\MTM^{\integral}((X;H),\real)$.
We obtain
$(\nbigt\otimes\nbigt_{\sm}(f))[\ast H]
\in\MTM^{\integral}(X,\real)$.

\begin{prop}
\label{prop;22.7.6.1}
$(\nbigt\otimes\nbigt_{\sm}(f))[\ast H]$
is a rescalable mixed twistor $\nbigd$-module on $X$.
\end{prop}
\pf
We obtain
\[
\nbigttilde:=
\Bigl(
\bigl(
\nbigt\boxtimes\nbigu_{\cnum}(1,0)
\bigr)\otimes\nbigt_{\sm}(tf)
\Bigr)[\ast((\proj^1\times H)\cup(\{0\}\times X) )]
\in\MTM^{\integral}(\cnum\times X,[\ast 0],\real).
\]
It is easy to see that
$\nbigttilde\in\MTM_{(1,1)}
\bigl(\cnum\times X,[\ast 0],\real\bigr)$
and that
$\ttU(\nbigttilde)=(\nbigt\otimes\nbigt_{\sm}(f))[\ast H]$.
Hence,
$(\nbigt\otimes\nbigt_{\sm}(f))[\ast H]$
is a rescalable mixed twistor $\nbigd$-module.
\hfill\qed

\begin{rem}
The irregular Hodge filtration of objects
in Proposition {\rm\ref{prop;22.7.6.1}}
was first studied in {\rm\cite{Sabbah-Yu}}.
\hfill\qed
\end{rem}

\section{Algebraic case}
\label{section;21.6.28.1}

\subsection{Preliminary}

\subsubsection{Analytification}

Let $Y$ be a complex quasi-projective manifold with the complex topology.
Let $Y^{\alg}$ denote the scheme with the Zariski topology whose analytification is $Y$.
Let $\nbigo^{\alg}_Y$ denote the structure sheaf of $Y^{\alg}$.
There exists the naturally defined morphism of
the ringed spaces
$\rho:(Y,\nbigo_Y)\lrarr(Y^{\alg},\nbigo_{Y}^{\alg})$.
For any $\nbigo^{\alg}_Y$-module $\nbign$ on $Y^{\alg}$,
let $\nbign^{\an}$ denote its analytification on $Y$,
i.e., $\nbign^{\an}=\rho^{\ast}(\nbign)=
\nbigo_{Y}\otimes_{\rho^{-1}\nbigo^{\alg}_{Y}}\rho^{-1}(\nbign)$.
According to \cite[Proposition 10]{Serre-GAGA} (see also \cite{Danilov}),
the analytification functor is exact and faithful.
If $\nbign$ is a coherent $\nbigo^{\alg}_Y$-module,
then $\nbign^{\an}$ is a coherent $\nbigo_Y$-module.
If moreover $Y$ is projective,
according to \cite[Theorem 3]{Serre-GAGA},
the analytification functor induces
an equivalence between the categories of
coherent $\nbigo_Y^{\alg}$-modules
and coherent $\nbigo_Y$-modules.
Let us add minor complements for the convenience of our arguments.

\begin{lem}
If $\nbign$ is a quasi-coherent $\nbigo^{\alg}_Y$-module,
$\nbign^{\an}$ is a good $\nbigo_Y$-module
in the sense of
{\rm\cite[Definition 4.22]{kashiwara_text}}.
\end{lem}
\pf
There exists a directed family $\{\nbigq_i\}$
of coherent $\nbigo^{\alg}_Y$-submodules of $\nbign$
such that $\sum\nbigq_i=\nbign$.
(For example, see \cite[II. Exercise 5.15 (e)]{Hartshorne}.)
The natural morphisms
$\nbigq_i^{\an}\lrarr\nbign^{\an}$
are monomorphisms,
and we have $\nbign^{\an}=\sum\nbigq_i^{\an}$.
Hence, $\nbign^{\an}$ is good.
\hfill\qed

\begin{lem}
\label{lem;21.6.23.10}
If $Y$ is projective,
then the analytification induces
an equivalence between the category of
quasi-coherent $\nbigo_Y^{\alg}$-modules
and the category of good $\nbigo_Y$-modules.
\end{lem}
\pf
Let $\nbigm$ be a good $\nbigo_Y$-module.
There exists a directed family $\{\nbigg_i\}$ of
coherent $\nbigo_X$-submodules of $\nbigm$
such that $\sum \nbigg_i=\nbigm$.
Let $\nbigm'$ denote the kernel of the natural morphism
$\bigoplus\nbigg_i\lrarr\nbigm$,
which is also good (see \cite[Proposition 4.23]{kashiwara_text}).
There exists a directed family $\{\nbigg_j'\}$
of coherent $\nbigo_X$-submodules of $\nbigm'$
such that $\sum\nbigg_j'=\nbigm'$.
Let $\varphi:\bigoplus\nbigg_j'\lrarr\bigoplus\nbigg_i$
denote the naturally induced morphism.
We have $\Cok(\varphi)\simeq\nbigm$.
There exist coherent $\nbigo_{Y}^{\alg}$-modules
$\nbigg_i^{\alg}$ and
$\nbigg_j^{\prime\alg}$
whose analytifications are
$\nbigg_i$ and $\nbigg_j'$, respectively.
By the construction of the analytification,
$\bigl(
\bigoplus\nbigg^{\alg}_i\bigr)^{\an}$
and
$\bigl(
\bigoplus\nbigg^{\prime\alg}_j\bigr)^{\an}$
are naturally isomorphic to
$\bigoplus\nbigg_i$
and $\bigoplus\nbigg_j$, respectively.
Moreover, there exists a unique morphism
$\varphi^{\alg}:
\bigoplus \nbigg_j^{\prime\alg}
\lrarr
\bigoplus\nbigg_i^{\alg}$
which induces $\varphi$.
We obtain
$\bigl(\Cok(\varphi^{\alg})\bigr)^{\an}\simeq\nbigm$.

Let $\nbign_a$ $(a=1,2)$ be quasi-coherent $\nbigo^{\alg}_Y$-modules.
Let $f:\nbign_1^{\an}\lrarr\nbign_2^{\an}$
be an $\nbigo_Y$-morphism.
Let $\{\nbigg_{a,i}\}$ be directed families of
coherent sheaves of $\nbign_a$ such that
$\sum \nbigg_{a,i}=\nbign_a$.
For each $i$, there exists $k(i)$
such that
$f(\nbigg_{1,i}^{\an})\subset \nbigg_{2,k(i)}^{\an}$.
Hence, there exists $f^{\alg}_i:\nbigg_{1,i}\lrarr\nbign_2$
which induces the restriction
$f_{|\nbigg_{1,i}^{\an}}$.
Let 
$\psi:\bigoplus \nbigg_{1,i}\lrarr \nbign_{2}$
denote the morphism induced by $f_i^{\alg}$.
Let $\nbign_1'$ denote the kernel of
the natural morphism
$\bigoplus \nbigg_{1,i}\lrarr\nbign_1$.
Because the induced morphism
$\psi^{\an}:
\bigl(\bigoplus\nbigg_{1,i}\bigr)^{\an}\lrarr\nbign_2^{\an}$
factors through $\nbign_1^{\an}$,
the restriction of $\psi^{\an}$ to
$(\nbign_1')^{\an}$ is $0$.
It implies that
the restriction of $\psi$ to
$\nbign_1'$ is $0$,
i.e.,
$\psi$ induces a morphism
$f^{\alg}:\nbign_1\lrarr\nbign_2$.
By the construction,
the analytification of $f^{\alg}$ is equal to $f$.
\hfill\qed

\vspace{.1in}

Suppose that $Y$ is a Zariski open subset of
a smooth complex projective variety $\Ybar$
such that $H=\Ybar\setminus Y$ is a hypersurface.
Let $\iota_Y:Y^{\alg}\lrarr \Ybar^{\alg}$
denote the inclusion.
We may naturally regard
$\iota_{Y\ast}\nbigo^{\alg}_Y$
as the sheaf of algebraic meromorphic functions
on $(\Ybar^{\alg},H^{\alg})$.
For any $\nbigo^{\alg}_Y$-module $\nbign$,
we obtain $\iota_{Y\ast}\nbigo^{\alg}_{Y}$-module $\iota_{Y\ast}(\nbign)$,
which induces an $\nbigo_{\Ybar}(\ast H)$-module
$\iota_{Y\ast}(\nbign)^{\an}$.

\begin{lem}
The above procedure induces an equivalence
between the category of
coherent $\nbigo^{\alg}_Y$-modules
and the category of coherent $\nbigo_{\Ybar}(\ast H)$-modules
which is good as an $\nbigo_{\Ybar}$-module.
It also induces an equivalence between the category of
quasi-coherent $\nbigo^{\alg}_Y$-modules
and the category of $\nbigo_{\Ybar}(\ast H)$-modules
which is good as an $\nbigo_{\Ybar}$-module.
\end{lem}
\pf
Let $\nbigm$ be an $\nbigo_{\Ybar}(\ast H)$-module which
is a good $\nbigo_{\Ybar}$-module.
There exists a quasi-coherent $\nbigo^{\alg}_{\Ybar}$-module
$\nbigm^{\alg}$ such that $(\nbigm^{\alg})^{\an}\simeq\nbigm$.
It is easy to see that
the natural morphism
$\nbigm^{\alg}\lrarr
 \iota_{Y\ast}\iota_Y^{\ast}\nbigm^{\alg}$
is an isomorphism.
Thus, we obtain the essential surjectivity in the second claim.
 
Suppose moreover that $\nbigm$ is coherent over $\nbigo_{\Ybar}(\ast H)$-module.
Let us observe that there exists a coherent $\nbigo_{\Ybar}$-submodule
$\nbigm_0\subset\nbigm$ such that
$\nbigm_0(\ast H)=\nbigm$.
Let $\{\nbigg_i\}$ be a directed family of
coherent $\nbigo_{\Ybar}$-submodules of $\nbigm$ such that
$\sum\nbigg_i=\nbigm$.
We have
$\sum\nbigg_i(\ast H)=\nbigm$.
Because $\nbigo_{\Ybar}(\ast H)$ is Noetherian,
there exists $i_0$ such that 
$\nbigg_{i_0}(\ast H)=\nbigm$.
There exists a coherent $\nbigo_{\Ybar}^{\alg}$-module
$\nbigg_{i_0}^{\alg}$ such that
$(\nbigg_{i_0}^{\alg})^{\an}\simeq\nbigg_{i_0}$.
Then, we obtain
$(\iota_{Y\ast}\iota_{Y}^{\ast}(\nbigg_{i_0}^{\alg}))^{\an}
=\nbigm$.
Thus, we obtain the essential surjectivity in the first claim.

Let $\nbign_i$ $(i=1,2)$ be quasi-coherent $\nbigo^{\alg}_Y$-modules.
Let $f:(\iota_{Y\ast}\nbign_1)^{\an}\lrarr(\iota_{Y\ast}\nbign_2)^{\an}$
be an $\nbigo_{\Ybar}(\ast H)$-homomorphism.
There exists a unique
$\iota_{Y\ast}\nbigo^{\alg}_Y$-homomorphism
$f^{\alg}:\iota_{Y\ast}\nbign_1\lrarr\iota_{Y\ast}\nbign_2$
such that $(f^{\alg})^{\an}=f$,
which is induced by
$\iota_Y^{\ast}(f^{\alg}):\nbign_1\lrarr\nbign_2$.
\hfill\qed

\subsubsection{Algebraic $\nbigrtilde$-modules}

Let $X$ be a smooth quasi-projective variety.
Let $X^{\alg}$ denote the associated scheme with the Zariski topology.
Let $p_X:\cnum_{\lambda}\times X\lrarr X$ denote the projection.
Let $\Theta_X^{\alg}$ denote the algebraic tangent sheaf.
Let $\nbigd^{\alg}_{\cnum_{\lambda}\times X}$
denote the sheaf of algebraic differential operators on
$\cnum_{\lambda}\times X$.
Let $\nbigrtilde^{\alg}_X\subset\nbigd^{\alg}_{\cnum_{\lambda}\times X}$
denote the sheaf of subalgebras generated by
$\lambda p_X^{\ast}\Theta_X^{\alg}$
and $\lambda^2\del_{\lambda}$
over $\nbigo^{\alg}_{\cnum_{\lambda}\times X}$.

Let $\Xbar$ be a smooth projective manifold
with an open embedding $X\lrarr \Xbar$
such  that $D_{\infty}=\Xbar\setminus X$ is a hypersurface.
Any $\nbigrtilde^{\alg}_X$-module $\gbigm$
naturally induces
an $\gbigrtilde_{\Xbar}(\ast D_{\infty})$-module $\gbigm^{\an}$.

\begin{lem}
\label{lem;21.6.23.1}
 The above procedure induces an equivalence
between the category of coherent
$\nbigrtilde^{\alg}_{X}$-modules
and the category of good coherent
$\gbigrtilde_{\Xbar}(\ast D_{\infty})$-modules. 
\end{lem}
\pf
Let $\gbigm$ be a good coherent
$\gbigrtilde_{\Xbar}(\ast D_{\infty})$-module.
There exists a directed family
$\{\gbigg_i\}$ of coherent $\nbigo_{\gbigxbar}$-submodules
of $\gbigm$
such that $\sum\gbigg_i=\gbigm$.
There exists $i_0$ such that
the induced $\gbigrtilde_{\Xbar}(\ast D_{\infty})$-morphism
$\gbigrtilde_{\Xbar}(\ast\Dbar_{\infty})
\otimes\gbigg_{i_0}\lrarr\gbigm$
is an epimorphism.
We set $\gbigb_1:=\gbigg_{i_0}$.
Let $\gbigk$ denote the kernel of
$\gbigrtilde_{\Xbar}(\ast\Dbar_{\infty})
\otimes\gbigb_1\lrarr\gbigm$.
Similarly, there exists a coherent
$\nbigo_{\gbigxbar}$-module
$\gbigb_2$
with a monomorphism of $\nbigo_{\gbigxbar}$-modules
$\varphi:\gbigb_2\lrarr\gbigk$
such that
the induced $\gbigrtilde_{\Xbar}(\ast\Dbar_{\infty})$-morphism
$\gbigrtilde_{\Xbar}(\ast D_{\infty})
\otimes\gbigb_2
\lrarr \gbigk$
is an epimorphism.
There exist $\nbigo_{\gbigxbar}$-modules
$\gbigb_i^{\alg}$
whose analytifications are isomorphic to
$\gbigb_i$.
There exists a morphism
of $\nbigo_{\nbigx^{\alg}}$-modules
$\varphi^{\alg}:
(\gbigb_2^{\alg})_{|\nbigx^{\alg}}
\lrarr
 \nbigrtilde_X^{\alg}
 \otimes
 (\gbigb_1^{\alg})_{|\nbigx^{\alg}}$
which induces $\varphi$.
The analytification of $\Cok(\varphi^{\alg})$
is naturally isomorphic to $\Cok(\varphi)$.

Let $\gbigm_i$ be $\nbigrtilde^{\alg}_X$-modules.
Let $f:\gbigm_1^{\an}\lrarr\gbigm_2^{\an}$
be a morphism of
$\gbigrtilde_{\gbigxbar}(\ast D_{\infty})$-modules.
There exists an $\nbigo_{\nbigx^{\alg}}$-morphism
$f^{\alg}:\gbigm_1\lrarr\gbigm_2$
which induces $f$.
Let $U$ be any affine open subset of $X$.
Let $v$ be a section of $\Theta^{\alg}_U$.
We obtain an $\nbigo^{\alg}_{U}$-homomorphism
$v\circ f^{\alg}-f^{\alg}\circ v:
\gbigm_{1|U^{\alg}}
\lrarr
\gbigm_{2|U^{\alg}}$.
Because its analytification is $0$,
we obtain $\lambda v\circ f^{\alg}-f^{\alg}\circ (\lambda v)=0$.
It implies that $f^{\alg}$ is a morphism of
$\nbigr^{\alg}_X$-modules.
Similarly,
we can prove that
$\lambda^2\del_{\lambda}\circ f^{\alg}
-f^{\alg}\circ (\lambda^2\del_{\lambda})=0$.
Thus, we obtain that $f^{\alg}$
is a morphism of $\nbigrtilde^{\alg}_X$-modules.
\hfill\qed

\subsection{Algebraic $\nbigrtilde_X$-modules underlying
integrable mixed twistor $\nbigd$-modules}

Let $X$ be a quasi-projective manifold
with a smooth projective compactification $\Xbar$.
Set $D_{\infty}=\Xbar\setminus X$.
Let $\nbigc^{\alg}(X)$ denote the category of
$\nbigrtilde^{\alg}_{X}$-modules $\gbigm$ such that
$\gbigm^{\an}\in \gbigc_{\Malg}(\Xbar;D_{\infty})$.
The condition is independent of the choice of $\Xbar$.
By definition, we have the following.
\begin{lem}
Let $U$ be a Zariski open subset of $X$.
For any $\gbigm_U\in \nbigc^{\alg}_{\res}(U)$
there exists $\gbigm_X\in\nbigc^{\alg}_{\res}(X)$
such that
$\gbigm_{X|U}=\gbigm_U$. 
\hfill\qed
\end{lem}

\subsubsection{Direct image by projective morphisms}

Let $\rho:X\lrarr Y$ be a projective morphism.
Let $\omega_X^{\alg}$ denote the canonical sheaf of $X^{\alg}$.
We set $\omega_{\nbigx}^{\alg}=\lambda^{-\dim X}p_{X}^{\ast}\omega^{\alg}_X$.
We define
$\nbigr^{\alg}_{Y\larr X}:=
\rho^{\ast}(\nbigr^{\alg}_Y)\otimes (\omega^{\alg}_{\nbigx})^{-1}$.
For a coherent $\nbigrtilde^{\alg}_X$-module $\nbigm$,
we define
$\rho_{\dagger}^i(\nbigm)$
as the $i$-th cohomology sheaf of
$R\rho_{\ast}\Bigl(
\nbigr^{\alg}_{Y\larr X}\otimes^L_{\nbigr_X}\nbigm
\Bigr)$.
\begin{lem}
There exist natural isomorphisms
$\rho_{\dagger}^i(\nbigm)^{\an}
\simeq
\rho_{\dagger}^i(\nbigm^{\an})$.
\end{lem}
\pf
We can obtain a natural transform
$\rho_{\dagger}^i(\nbigm)^{\an}
\lrarr
\rho_{\dagger}^i(\nbigm^{\an})$
in a natural way.
We can check that it is an isomorphism
by using a resolution of
$\nbigg^{\bullet}\lrarr \nbigm$
such that each $\nbigg^i$ is of the form
$\nbigrtilde^{\alg}_X\otimes \nbigj^i$
for a coherent $\nbigo^{\alg}_{\cnum\times X}$-module $\nbigj^i$.
\hfill\qed

\vspace{.1in}
We obtain the following.

\begin{prop}
Let $\rho:X_1\lrarr X_2$ be a projective morphism.
For any $\gbigm\in\nbigc^{\alg}(X_1)$,
$\rho^i_{\dagger}(\gbigm)$
are objects of $\nbigc^{\alg}(X_2)$.
\hfill\qed
\end{prop}

\subsubsection{Strict specializability}

Suppose that $X=\cnum\times X_0$ for a complex algebraic manifold $X_0$.
Let $t$ denote the standard coordinate of $\cnum$.
Let $\pi$ denote the projection $X\lrarr X_0$.
Let $V\nbigr^{\alg}_X\subset\nbigr^{\alg}_X$ denote
the sheaf of subalgebras generated by
$\pi^{\ast}\nbigr^{\alg}_{X_0}$ and $\lambda t\del_t$.
We set $V\nbigrtilde^{\alg}_X=
V\nbigr_X^{\alg}\langle\lambda^2\del_{\lambda}\rangle$.
Any $V\nbigrtilde_X^{\alg}$-module $\gbigm$
induces
$V\gbigrtilde_{\Xbar(\ast D_{\infty})}$-module
$\gbigm^{\an}$.
The following lemma is similar to
Lemma \ref{lem;21.6.23.1}.

\begin{lem}
The above procedure induces an equivalence
between the category of
$V\nbigrtilde^{\alg}_X$-modules which are coherent over $V\nbigr^{\alg}_X$,
and the category of
$V\gbigrtilde_{\Xbar(\ast D_{\infty})}$-modules
which are coherent over $V\gbigr_{\Xbar(\ast D_{\infty})}$.
\hfill\qed
\end{lem}

Let $\gbigm$ be a coherent $\nbigrtilde^{\alg}_X$-module
or a coherent $\nbigrtilde^{\alg}_X(\ast t)$-module.
We define the notion of $V$-filtrations for $\gbigm$
as a filtration $V_{\bullet}(\gbigm)$ indexed by $\real$
such that
(i) each $V_a(\gbigm)$ is
a $V\nbigrtilde^{\alg}_X$-module
which is coherent over $V\nbigr^{\alg}_X$,
(ii) they satisfy the conditions in \S\ref{subsection;21.4.12.40}.
If $\gbigm$ has a $V$-filtration,
it is called strictly specializable along $t$.
It is called regular along $t$
if moreover $V_a(\gbigm)$ is coherent over
$\pi^{\ast}\nbigr^{\alg}_{X_0}$.
If $\gbigm$ has a $V$-filtration,
we define
$\gbigm[!t]=
\nbigrtilde^{\alg}_X\otimes_{V\nbigrtilde^{\alg}_X}V_{<0}\gbigm$
and
$\gbigm[\ast t]=
\nbigrtilde^{\alg}_X\otimes_{V\nbigrtilde^{\alg}_X}V_{0}\gbigm$.
We have the following lemma.
\begin{lem}
Let $\nbigm$ be a coherent
$\nbigrtilde^{\alg}_X$-module.
\begin{itemize}
 \item $\nbigm$ is strictly specializable along $t$
if and only if
$\nbigm^{\an}$ has a $V$-filtration in the sense of
Lemma {\rm\ref{lem;21.3.13.21}}
       such that each $V_a(\nbigm^{\an})$ is coherent over $V\gbigr_{X(\ast H)}$.
       In the case,
       we have
       $\nbigm[\star t]^{\an}\simeq\nbigm^{\an}[\star t]$.      
 \item $\nbigm$ is regular along $t$ if and only if
       $\nbigm^{\an}$ is regular along $t$.
\end{itemize} 
Similar claims hold for $\nbigm(\ast t)$.
\hfill\qed
\end{lem}

\begin{cor}
Let $\gbigm\in\nbigc^{\alg}(X)$. 
Then, $\gbigm$ is strictly specializable along $t$,
and $\gbigm[\star t]$ are objects in $\nbigc^{\alg}(X)$. 
(Note Proposition {\rm\ref{prop;21.6.26.2}}.)
\hfill\qed
\end{cor}

\subsubsection{Strict specializability
and Beilinson functors along an algebraic function}

Let $f$ be an algebraic function on a smooth projective variety $X$.
An $\nbigrtilde^{\alg}_{X}$-module
(resp. $\nbigrtilde^{\alg}_X(\ast f)$-module)
$\gbigm$
is called strictly specializable along $f$
if the $\gbigrtilde_{\Xbar(\ast D_{\infty})}$-module
(resp. $\gbigrtilde_{\Xbar(\ast D_{\infty})}(\ast|(f)_0|)$-module)
$\gbigm^{\an}$ is strictly specializable along $f$
modulo $D_{\infty}^{\an}$.
We say that $\gbigm$ is localizable along $f$
if moreover there exist
$\nbigrtilde^{\alg}_X$-modules $\gbigm[\star f]$
such that
$\gbigm[\star f]^{\an}
\simeq
\gbigm^{\an}[\star f]$.
Note that the conditions are stated
in terms of the push-forward by
the graph embedding of $f$ in an algebraic way.

\begin{prop}
Let $\gbigm\in\nbigc^{\alg}(X)$.
Then, $\gbigm$ is strictly specializable and localizable along $f$,
 and $\gbigm[\star f]$ are objects in $\nbigc^{\alg}(X)$.
\hfill\qed
\end{prop}

We define the $\nbigrtilde^{\alg}_X$-module
$(\II_f^{a,b})^{\alg}$ on $\nbigx^{\alg}$
as in the case of $\II_f^{a,b}$.
For $\nbigrtilde^{\alg}_X(\ast f)$-module $\gbigm$,
we define
$\Pi^{a,b}_f(\gbigm):=(\II_f^{a,b})^{\alg}\otimes\gbigm$.
Suppose that $\Pi^{a,b}_f(\gbigm)$ are strictly specializable
and localizable along $f$ for any $a,b$.
Then, we put
$\Pi^{a,b}_{f\star}(\gbigm)
=\Pi^{a,b}_f(\gbigm)[\star f]$,
and we define
\[
\Pi^{a,b}_{f\ast!}(\gbigm):=
\varprojlim_{N\to\infty}
\Cok\Bigl(
\Pi^{b,N}_{f!}(\gbigm)\lrarr
\Pi^{a,N}_{f\ast}(\gbigm)
\Bigr).
\]
In particular, we define
$\psi_f^{(a)}(\gbigm)=\Pi^{a,a}_{f\ast!}(\gbigm)$,
and 
$\Xi^{(a)}_f(\gbigm)
=\Pi^{a,a+1}_{f,\ast!}(\gbigm)$.
We define $\phi^{(0)}_f(\gbigm)$
as the cohomology of the induced complex
\[
 \gbigm[!f]
 \lrarr
 \Xi^{(0)}_f(\gbigm)
 \oplus\gbigm
 \lrarr
 \gbigm[\ast f].
\]
We can recover $\gbigm$
as the cohomology of the following complex
as in the analytic case:
\[
 \psi^{(1)}_f(\gbigm)
 \lrarr
 \phi^{(0)}_f(\gbigm)
 \oplus
 \Xi^{(0)}_f(\gbigm)
 \lrarr
 \psi^{(0)}_f(\gbigm).
\]

\begin{prop}
Let $\gbigm\in\nbigc^{\alg}(X)$.
The $\nbigrtilde^{\alg}_X(\ast f)$-modules
$\Pi^{a,b}_f(\gbigm)$
are strictly specializable and localizable along $f$,
and 
$\Pi^{a,b}_{f\star}(\gbigm)
=\Pi^{a,b}_f(\gbigm)[\star f]$
are objects of $\nbigc^{\alg}(X)$. 
In particular, we obtain
the objects
$\psi^{(a)}_f(\gbigm)$,
$\Xi^{(a)}_f(\gbigm)$
and $\phi^{(0)}_f(\gbigm)$
of $\nbigc^{\alg}(X)$.
\hfill\qed
\end{prop}

\subsubsection{Localizability along hypersurfaces}

Let $H$ be an algebraic hypersurface of $X$.
Let $\Hbar$ denote the closure of $H$ in $\Xbar$.
Let $\gbigm\in\nbigrtilde^{\alg}_{X(\ast H)}$-module.
It is called localizable along $H$
if $\gbigm^{\an}$ is localizable along $\Hbar$
modulo $D_{\infty}$.
If $\gbigm$ is localizable along $H$,
there exists $\gbigm[\star H]$
such that
$\gbigm[\star H]^{\an}\simeq
\gbigm^{\an}[\star \Hbar]$.
For any Zariski open subset $U\subset X$
with an algebraic function $f\in\nbigo^{\alg}_X(U)$
such that $f^{-1}(0)=H\cap U$,
we have
$\nbigm_{|U}[\star f]=(\nbigm[\star H])_{|U}$.

\begin{prop}
Any $\nbigm\in\nbigc^{\alg}$ is localizable along $H$.
\hfill\qed
\end{prop}

\subsubsection{Some other functoriality}

For any $\nbigrtilde^{\alg}_X$-module $\gbigm$,
we define
$\DD_X(\gbigm)=
\lambda^{d_X}
\nrhom_{\nbigr^{\alg}_X}(\gbigm,\nbigr^{\alg}_{X}\otimes(\omega^{\alg}_X)^{-1})[d_X]$
which is naturally an object
in the derived category of $\nbigrtilde^{\alg}_X$-modules.
There exists a natural isomorphism
$\DD_X(\gbigm)^{\an}
\simeq
 \DD_{\Xbar(\ast D_{\infty})}
  (\gbigm^{\an})$.
\begin{prop}
If $\gbigm\in\nbigc^{\alg}(X)$,
then $\DD_X(\gbigm)\in\nbigc^{\alg}(X)$.
\hfill\qed
\end{prop}

For any 
$\nbigrtilde^{\alg}_{X_i}$-modules $\gbigm_i$ $(i=1,2)$,
we naturally obtain
the $\nbigrtilde^{\alg}_{X_1\times X_2}$-module
$\gbigm_1\boxtimes\gbigm_2$.
There exists a natural isomorphism
$(\gbigm_1\boxtimes\gbigm_2)^{\an}
\simeq
 \gbigm_1^{\an}\boxtimes\gbigm_2^{\an}$.

\begin{prop}
If $\gbigm_i\in\nbigc^{\alg}(X_i)$,
 $\gbigm_1\boxtimes\gbigm_2$ is an object in
 $\nbigc^{\alg}(X_1\times X_2)$.
\hfill\qed 
\end{prop}

Let $f:X\lrarr Y$ be any morphism of complex algebraic manifolds.
Let $\gbigm\in\nbigc^{\alg}(Y)$.
If $f$ is non-characteristic for $\gbigm$,
we define
\[
(\lefttop{T}f^!)^i(\gbigm)=
\left\{
\begin{array}{ll}
 0 & (i\neq \dim Y-\dim X)\\
 \lambda^{-\dim Y+\dim X}f^{\ast}(\gbigm)
  & (i=\dim Y-\dim X)
\end{array}
\right.
\]
\[
(\lefttop{T}f^{\ast})^i(\gbigm)=
\left\{
\begin{array}{ll}
 0 & (i\neq \dim X-\dim Y)\\
  f^{\ast}(\gbigm)
  & (i=\dim X-\dim Y)
\end{array}
\right.
\]
If $f$ is a closed immersion,
we define
$(\lefttop{T}f^{\star})^i(\gbigm)$
in the way explained in \S\ref{subsection;21.6.22.2}.
Note that the additional assumption in \S\ref{subsection;21.6.22.2}
is always satisfied in the algebraic case.
The definitions are compatible
if $f$ is a closed immersion and non-characteristic for $\gbigm$.
For a general morphism $f$,
we define $(\lefttop{T}f^{\star})^i(\gbigm)$
by using the decomposition of $f$
into the closed embedding and the projection.
We have
$(\lefttop{T}f^{\star})^i(\gbigm)^{\an}
\simeq
(\lefttop{T}f^{\star})^i(\gbigm^{\an})$.

\subsection{Rescalable objects}

Let $X$ be a quasi-projective manifold.
For any $\nbigrtilde^{\alg}_X$-module $\gbigm$,
we obtain
$\nbigrtilde^{\alg}_{\lefttop{\tau}X\setminus\lefttop{\tau}X_0}$-module
$\lefttop{\tau}\gbigm$
in the natural way.
Let $\nbigc_{\res}^{\alg}(X)\subset\nbigc^{\alg}(X)$
denote the subcategory of $\gbigm\in\nbigc^{\alg}(X)$
such that
$\lefttop{\tau}\gbigm\in\nbigc^{\alg}(\lefttop{\tau}X\setminus\lefttop{\tau}X_0)$.
The condition is independent of the choice of $\Xbar$.

\begin{rem}
In the above definition,
$(\lefttop{\tau}\gbigm)^{\an}$
is an object of 
$\nbigc\bigl(\proj^1\times \Xbar,
 (\proj^1\times D_{\infty})
 \cup (\{0,\infty\}\times\Xbar)
 \bigr)$,
which is a stronger condition
than
$\gbigm^{\an}\in\nbigc_{\res}(\Xbar;D_{\infty})$. 
\hfill\qed
\end{rem}

\begin{prop}
Let $\gbigm\in\nbigc^{\alg}_{\res}(X)$.
\begin{itemize}
 \item
      For any hypersurface $H$ of $X$,
      $\gbigm[\star H]$ are objects of
      $\nbigc^{\alg}_{\res}(X)$ $(\star=!,\ast)$.
 \item  For any algebraic function $g$ on $X$,
$\Xi^{(a)}_g(\gbigm)$,
$\psi_g^{(a)}(\gbigm)$
and $\phi_g^{(0)}(\gbigm)$
are objects in $\nbigc^{\alg}_{\res}(X)$.
 \item For any morphism $f:Y\lrarr X$,
       $(\lefttop{T}f^{\star})^i(\gbigm)$
       are objects of
       $\nbigc^{\alg}_{\res}(Y)$.       
\end{itemize}
\end{prop}
\pf
The first claim follows from
the natural isomorphisms
$(\lefttop{\tau}\gbigm)[\star \lefttop{\tau}H]
\simeq
\lefttop{\tau}\bigl(
 \gbigm[\star H]
 \bigr)$.
We obtain the other claims similarly.
\hfill\qed

\subsubsection{Irregular Hodge filtrations}

Let $\gbigm\in\nbigc^{\alg}_{\res}(X)$.
Because $(\lefttop{\tau}\gbigm)^{\an}$ is regular along $\tau$,
each $V_a(\lefttop{\tau}\gbigm)$
is coherent over $\pi^{\ast}\nbigr^{\alg}_X$.

\begin{prop}
There exists a filtration
$F^{\irr}_{\bullet}(V_a(\lefttop{\tau}\gbigm))$
by coherent $\nbigo^{\alg}_{\lefttop{\tau}\nbigx}$-submodules
such that
\[
F^{\irr}_{\bullet}(V_a(\lefttop{\tau}\gbigm))^{\an}_{|\lefttop{\tau}\nbigx}
=F^{\irr}_{\bullet}(V_a(\lefttop{\tau}\gbigm^{\an})_{|\lefttop{\tau}\nbigx}).
\]
\end{prop}
\pf
Let $X_1:=\proj^1_{\tau}\times \Xbar$.
We set
$H_1:=(\proj^1_{\tau}\times H)$
and
$D_{1,\infty}=(\proj^1_{\tau}\times D_{\infty})\cup
 (\{\infty\}\times X)$.
By the construction,
the $V\gbigrtilde_{\lefttop{\tau}\Xbar}
\bigl(\ast (\lefttop{\tau}H\cup
 \lefttop{\tau}D_{\infty})\bigr)$-module
$V_a(\lefttop{\tau}\gbigm^{\an})$
naturally extends to
a good coherent
$V\gbigrtilde_{X_1}(\ast(H_1\cup D_{1,\infty}))$-module
$\bigl(
V_a(\lefttop{\tau}\gbigm)
\bigr)^{\an}$.
The $\cnum^{\ast}$-equivariant filtration $F^{\irr}$ of
$V_a(\lefttop{\tau}\gbigm^{\an})_{|\lefttop{\tau}\nbigx}$
by coherent $\nbigo_{\lefttop{\tau}\nbigx}$-modules
naturally extends to
a $\cnum^{\ast}$-equivariant filtration of 
$\bigl(
V_a(\lefttop{\tau}\gbigm)
\bigr)^{\an}$
by good coherent
$\nbigo_{\gbigx_1}
\bigl(\ast(\gbigh_1\cup\gbigd_{1,\infty}\cup \gbigx_1^{\infty})
\bigr)$-modules,
denoted by $F^{\irr}$.
By Lemma \ref{lem;21.6.23.10},
there exists a filtration $F^{\irr}$
of $V_a(\lefttop{\tau}\gbigm)$
by $\nbigo^{\alg}_{\lefttop{\tau}\nbigx}(\ast\lefttop{\tau}\nbigh)$-coherent
submodules,
which induces
the filtration $F^{\irr}$ of $V_a(\lefttop{\tau}\gbigm)^{\an}$.
\hfill\qed

\vspace{.1in}
For any $\lambda\in\cnum$,
let $\iota_{\lambda}:X\simeq \{\lambda\}\times X\lrarr \nbigx$
denote the inclusion.
We set $\gbigm^{\lambda}:=\iota_{\lambda}^{\ast}(\gbigm)$.

\begin{cor}
There exists a filtration $F^{\irr}$ of
$\gbigm^{\lambda}$ by
coherent $\nbigo^{\alg}_X$-submodules
such that
$\bigl(
F^{\irr}_{\bullet}(\gbigm^{\lambda})\bigr)^{\an}
=F^{\irr}_{\bullet}\bigl(
(\gbigm^{\lambda})^{\an}
\bigr)$.
\hfill\qed
\end{cor}

Let $\iota_{(\lambda,\tau)=(1,0)}$
denote the inclusion
$X\simeq (1,0)\times X\lrarr \lefttop{\tau}\nbigx$.
Note that
$\iota_{(\lambda,\tau)=(1,0)}^{\ast}V_a(\lefttop{\tau}\gbigm)^{\an}$
is naturally isomorphic to
$\iota_{\infty}^{\ast}V_a(\gbigm^{\an})$.

\begin{cor}
There exists a filtration $F$
of $\iota_{(1,0)}^{\ast}V_a(\lefttop{\tau}\gbigm)$
 by coherent $\nbigo^{\alg}_X$-submodules
such that 
\[
 F^{\irr}_{\bullet}
 \iota_{(\lambda,\tau)=(1,0)}^{\ast}V_a(\lefttop{\tau}\gbigm)^{\an}
=F^{\irr}_{\bullet}
 \iota_{\infty}^{\ast}V_a(\gbigm^{\an}).
\]
 \hfill\qed
\end{cor}

\begin{prop}
Let $f:\gbigm_1\lrarr \gbigm_2$ be a morphism
in $\nbigc^{\alg}_{\res}(X)$
such that $\Cok(f)\in\nbigc_{\res}^{\alg}(X)$.
Then, the induced morphisms
$V_a(\lefttop{\tau}\gbigm_1)\lrarr V_a(\lefttop{\tau}\gbigm_2)$,
$\gbigm_1^{\lambda}\lrarr\gbigm_2^{\lambda}$ 
and
 $\iota_{(\lambda,\tau)=(1,0)}^{\ast}
 V_a(\lefttop{\tau}\gbigm_1)\lrarr
 \iota_{(\lambda,\tau)=(1,0)}^{\ast}
 V_a(\lefttop{\tau}\gbigm_2)$ 
are strictly compatible with $F^{\irr}$.
\hfill\qed
\end{prop}

\begin{prop}\mbox{{}}
\begin{itemize}
 \item 
Let $\rho:X_1\lrarr X_2$ be a projective morphism.
For any $\gbigm\in\nbigc^{\alg}_{\res}(X_1)$,
$\rho^i_{\dagger}(\gbigm)$
are objects of $\nbigc^{\alg}_{\res}(X_2)$.
Moreover, we have
$\rho^i_{\dagger}(\Rtilde_{F^{\irr}_{a+\bullet}}(\Xi_{\DR}(\gbigm)))
=\Rtilde_{F^{\irr}_{a+\bullet}}(\rho^i_{\dagger}\Xi_{\DR}(\gbigm))$,
where $\Rtilde_{F^{\irr}_{a+\bullet}}$
denote the Rees construction associated with
the filtration $F^{\irr}_{a+\bullet}$.
\item
For any $\gbigm\in\nbigc^{\alg}_{\res}(X)$,
$\DD_X(\gbigm)$ is an object of $\nbigc^{\alg}_{\res}(X)$.
Moreover, for any $\lambda\neq 0$,
we have
\[
 \DD_X(\Rtilde_{F^{\irr}_{a+\bullet}}\gbigm^{\lambda})
 =\Rtilde_{F^{\irr}_{<-1-a+\bullet}}\DDD_X(\gbigm^{\lambda}).
\]
 \item
For $\gbigm_i\in\nbigc^{\alg}_{\res}(X_i)$,
$\gbigm_1\boxtimes\gbigm_2\in\nbigc^{\alg}_{\res}(X_1\times X_2)$.
 We have
 \[
F^{\irr}_a\Xi_{\DR}(\gbigm_1\boxtimes\gbigm_2)
=\sum_{b_1+b_2\leq a}
 F^{\irr}_{b_1}\Xi_{\DR}(\gbigm_1)
 \boxtimes
F^{\irr}_{b_2}\Xi_{\DR}(\gbigm_2).
 \]   
 \item
Let $\rho:X\lrarr Y$ be an algebraic morphism of
quasi-projective manifolds.
Let $\gbigm\in\nbigc^{\alg}_{\res}(Y)$.
If $\rho$ is non-characteristic for $\gbigm$,
$\rho^{\ast}(\gbigm)\in\nbigc^{\alg}_{\res}(X)$.
Moreover,
$F^{\irr}_{\bullet}\rho^{\ast}\Xi_{\DR}(\gbigm)
=\rho^{\ast}F^{\irr}_{\bullet}\Xi_{\DR}(\gbigm)$. 
\hfill\qed
\end{itemize}  
\end{prop}

\subsubsection{Examples}

Let us mention two examples.
Starting from these examples,
we can construct many examples of objects
in $\nbigc^{\alg}_{\res}(X)$
by using the functoriality.

\begin{prop}
Let $f$ be an algebraic function on $X$.
Let $\nbigl^{\alg}(f)$ be the $\nbigr^{\alg}_X$-module
obtained as $\nbigo_{\nbigx}^{\alg}$ with $d+d(\lambda^{-1}f)$.
Then, $\nbigl^{\alg}(f)\in \nbigc^{\alg}_{\res}(X)$.
\hfill\qed
\end{prop}

\begin{prop}
Let $(M,F)$ be a filtered $\nbigd$-module
underlying a mixed Hodge module
on $X$ which is extendable
to a mixed Hodge module on a projective completion of $X$. 
Then,
$\Rtilde_F(M)\in\nbigc^{\alg}_{\res}(X)$. 
\hfill\qed
\end{prop}

\subsection{Rescalable mixed twistor $\nbigd$-modules}

\subsubsection{Integrable mixed twistor $\nbigd$-modules}

Let $\MTM^{\integral}(X)^{\alg}:=\MTM(\Xbar;D_{\infty})$.
For any $(\nbigt,W)=((\nbigm',\nbigm'',C),W)\in\MTM^{\integral}(X)^{\alg}$,
the underlying $\nbigrtilde_{\Xbar(\ast D_{\infty})}$-modules
$\nbigm'$ and $\nbigm''$ extends to
$\gbigrtilde_{\Xbar(\ast D_{\infty})}$-modules
$\Upsilon(\nbigm')$ and $\Upsilon(\nbigm'')$.
The extensions are compatible with
the functors for integrable mixed twistor $\nbigd$-modules.
Hence, we may regard $\nbigm'$ and $\nbigm''$
as $\nbigrtilde^{\alg}_X$-modules from the beginning.

As explained in \cite[\S14]{Mochizuki-MTM},
there are $6$-operations for integrable mixed twistor $\nbigd$-modules
in the algebraic setting.

\subsubsection{Partial Fourier transforms}

As in \S\ref{subsection;21.7.16.10},
for any $\nbigt\in\MTM^{\integral}(\cnum_t\times X)^{\alg}$,
we define
\[
 \FT_{\pm}(\nbigt):=
 \lefttop{T}(p_{\tau,X})^0_{!}\bigl(
 (\lefttop{T}p_{t,X}^{\ast})^1(\nbigt)
 \otimes\nbigt^{\alg}_{\sm}(\pm t\tau)
 \bigr)
 \simeq
 \lefttop{T}(p_{\tau,X})^0_{\ast}\bigl(
 (\lefttop{T}p_{t,X}^{\ast})^1(\nbigt)
 \otimes\nbigt^{\alg}_{\sm}(\pm t\tau)
 \bigr)
\in \MTM^{\integral}(\cnum\times X)^{\alg}.
\]

The following proposition is similar to
Proposition \ref{prop;21.7.12.12}.
\begin{prop}
There exist natural isomorphisms
$\FT_{\mp}\circ\FT_{\pm}(\nbigt)\simeq\nbigt\otimes\newTate(-1)$.
\hfill\qed
\end{prop}

We define the subcategories
$\MTM^{\integral}(\cnum\times X,[\star 0])^{\alg}$
and
$\MTM^{\integral}(\cnum\times X)^{\alg}_{\star}$
of $\MTM^{\integral}(\cnum\times X)^{\alg}$
as in \S\ref{subsection;21.7.16.11}.
The following proposition is similar to
Proposition \ref{prop;21.7.12.20}.
\begin{prop}
 $\FT_{\pm}$ induce equivalences
$\MTM^{\integral}(\cnum_t\times X,[\star 0])^{\alg}
 \simeq
 \MTM^{\integral}(\cnum_t\times X)^{\alg}_{\star}$.
\hfill\qed
\end{prop}

Let $F:X\lrarr Y$ be an algebraic morphism of
quasi-projective manifolds.
For $\star=!,\ast$,
we have
$(\id\times F)_{\star}:
D^b(\MTM^{\integral}(\cnum\times X)^{\alg})
\lrarr
D^b(\MTM^{\integral}(\cnum\times Y)^{\alg})$.
The following proposition is easy to see.
\begin{prop}
For any
$\nbigt^{\bullet}\in
D^b(\MTM^{\integral}(\cnum\times X)^{\alg})$,
there exist natural isomorphisms
\[
 \FT_{\pm}\circ
 \lefttop{T}(\id\times F)_{\star}(\nbigt^{\bullet})
 \simeq
 \lefttop{T}(\id\times F)_{\star}\circ\FT_{\pm}(\nbigt^{\bullet}).
\]

\hfill\qed
\end{prop}

The following proposition is similar to
Proposition \ref{prop;21.7.16.12}.
\begin{prop}
\label{prop;21.7.16.22}
 Let $\nbigt\in\MTM^{\integral}(\cnum\times X)^{\alg}$.
\begin{itemize}
 \item For any hypersurface $H$ of $X$,
       there exist natural isomorphisms
\[
       \FT_{\pm}(\nbigt[\star (\proj^1\times H)])
       \simeq
       \FT_{\pm}(\nbigt)[\star (\proj^1\times H)].
\]
 \item Let $f$ be an algebraic function on $X$.
       Let $\ftilde$ be the induced function on $\cnum\times X$.
       There exist natural isomorphisms
       $\FT_{\pm}\Pi^{a,b}_{\ftilde,\star}\bigl(
       \nbigt
       \bigr)
       \simeq
       \Pi^{a,b}_{\ftilde,\star}\bigl(
       \FT_{\pm}\nbigt
       \bigr)$ $(\star=!,\ast)$
       and 
       $\FT_{\pm}\Pi^{a,b}_{\ftilde,\ast !}\bigl(
       \nbigt
       \bigr)
       \simeq
       \Pi^{a,b}_{\ftilde,\ast !}\bigl(
       \FT_{\pm}\nbigt
       \bigr)$.
       In particular,
       we have natural isomorphisms
       $\FT_{\pm}\Xi^{(a)}_{\ftilde}(\nbigt)
       \simeq
       \Xi^{(a)}_{\ftilde}(\FT_{\pm}(\nbigt))$
       and 
       $\FT_{\pm}\psi^{(a)}_{\ftilde}(\nbigt)
       \simeq
       \psi^{(a)}_{\ftilde}(\FT_{\pm}(\nbigt))$.
       We also obtain
       $\FT_{\pm}\phi^{(0)}_{\ftilde}(\nbigt)
       \simeq
       \phi^{(0)}_{\ftilde}\FT_{\pm}(\nbigt)$.
       \hfill\qed
\end{itemize}
\end{prop}

Let $F:X\lrarr Y$ be an algebraic morphism of
quasi-projective manifolds.
For $\star=!,\ast$,
we have
$\lefttop{T}(\id\times F)^{\star}:
D^b(\MTM^{\integral}(\cnum\times Y)^{\alg})
\lrarr
D^b(\MTM^{\integral}(\cnum\times X)^{\alg})$.

\begin{cor}
For any $\nbigt^{\bullet}\in
D^b(\MTM^{\integral}(\cnum\times Y)^{\alg})$,
there exist natural isomorphisms 
\[
 \FT_{\pm}\circ
 \lefttop{T}(\id\times F)^{\star}(\nbigt^{\bullet})
 \simeq
 \lefttop{T}(\id\times F)^{\star}\circ\FT_{\pm}(\nbigt^{\bullet}).
\]
\hfill\qed
\end{cor}

The following proposition is similar to
Proposition \ref{prop;21.7.15.3}.
\begin{prop}
For $\nbigt\in\MTM^{\integral}(\cnum\times X)^{\alg}$,
there exist natural isomorphisms
\[
 \FT_{\pm}(\DD\nbigt)\otimes\newTate(1)
\simeq
 \DD(\FT_{\mp}(\nbigt)),
 \quad
 \FT_{\pm}(\nbigt^{\ast})
 \otimes\newTate(1)
 \simeq
 \FT_{\pm}(\nbigt)^{\ast}. 
\]
There also exist natural isomorphisms:
\[
 \FT_{\pm}(j^{\ast}\nbigt)
 \simeq
 j^{\ast}\FT_{\mp}(\nbigt),
 \quad
 \FT_{\pm}(\gammatilde^{\ast}\nbigt)
 \simeq
 \gammatilde^{\ast}\FT_{\pm}(\nbigt).
\]
\hfill\qed
\end{prop}

We use the notation in \S\ref{subsection;21.7.14.30}.
For $\nbigt_i^{\bullet}\in
D^b(\MTM^{\integral}(\cnum_{t_i}\times X_i)^{\alg})$,
we obtain
$\nbigt_1^{\bullet}\boxtimes\nbigt_2^{\bullet}
\in D^b\bigl(
\MTM^{\integral}(\cnum_{t_1}\times\cnum_{t_2}\times X)^{\alg}
\bigr)$.
We define
\[
 \nbigt_1^{\bullet}\langle
 \boxtimes,\otimes\rangle
 \rangle^{\star}\nbigt_2^{\bullet}
 :=
 \lefttop{T}(\Delta\times\id_X)^{\star}
 (\nbigt_1^{\bullet}\boxtimes\nbigt_2^{\bullet})
 \in
  D^b\bigl(
  \MTM^{\integral}(\cnum\times X)^{\alg}
  \bigr).
\]
We also define
\[
 \nbigt_1^{\bullet}
 \langle\boxtimes,\ast\rangle_{\star}
  \nbigt_2^{\bullet}:=
 \mu_{\star}
 \bigl(
 \nbigt_1^{\bullet}\boxtimes\nbigt_2^{\bullet}
 \bigr)
 \in
 D^b\bigl(
 \MTM^{\integral}(\cnum\times X)^{\alg}
 \bigr).
\]
The following proposition is similar to
Proposition \ref{prop;21.7.14.40}.
\begin{prop}
There exist the following natural isomorphisms:
\[
 \FT_{\pm}(\nbigt^{\bullet}_1)
 \langle\boxtimes,\otimes
 \rangle^!
 \FT_{\pm}(\nbigt^{\bullet}_2)
 \simeq
 \Bigl(
 \FT_{\pm}\bigl(
 \nbigt_1^{\bullet}\langle
 \boxtimes,\ast\rangle_{\ast}
 \nbigt_2^{\bullet}
 \bigr)
 \otimes\newTate(-1)
 \Bigr)
 [-1].
\]
\[
 \FT_{\pm}(\nbigt^{\bullet}_1)
 \langle\boxtimes,\otimes
 \rangle^{\ast}
 \FT_{\pm}(\nbigt^{\bullet}_2)
 \simeq
 \FT_{\pm}\bigl(
 \nbigt_1^{\bullet}\langle
 \boxtimes,\ast\rangle_{!}
 \nbigt_2^{\bullet}
 \bigr)[1].
\]
\hfill\qed
\end{prop}

As in (\ref{eq;21.7.16.21}),
for $\nbigt\in \MTM^{\integral}(\cnum\times X)^{\alg}$,
we define
\[
\left\{
\begin{array}{l}
 \ttP_{\ast}(\nbigt):=
  \nbigt\langle\boxtimes,\ast\rangle_{\ast}
  \bigl(
  \nbigu_{\cnum}(1,0)[!0]
  \bigr)
\in \MTM^{\integral}(\cnum\times X)^{\alg}_{\ast}
 \\
 \ttP_{!}(\nbigt):=
  \nbigt\langle\boxtimes,\ast\rangle_{!}
\bigl(
\nbigu_{\cnum}(1,0)[\ast 0]
\bigr)
 \in\MTM^{\integral}(\cnum\times X)^{\alg}_!.
\end{array}
\right.
\]
For $\star=!,\ast$,
there exist natural isomorphisms
$\ttP_{\star}\circ\ttP_{\star}(\nbigt)\simeq \ttP_{\star}(\nbigt)$.
If $\nbigt\in\MTM^{\integral}(\cnum\times X)^{\alg}_{\star}$,
there exists a natural isomorphism
$\nbigt\simeq\ttP_{\star}(\nbigt)$.
We obtain the induced functors
\[
\ttP_{\star}:
D^b(\MTM^{\integral}(\cnum\times X)^{\alg})
\lrarr
D^b(\MTM^{\integral}(\cnum\times X)^{\alg}_{\star}).
\]

As in (\ref{eq;21.7.12.13}),
we also obtain
\[
\FT_{\pm}:
\MTM^{\integral}(\cnum_t\times X,\real)^{\alg}
\lrarr
\MTM^{\integral}(\cnum_{\tau}\times X,\real)^{\alg}.
\]
We obtain equivalences
\[
 \FT_{\pm}:
\MTM^{\integral}(\cnum_t\times X,[\star 0],\real)^{\alg}
\simeq
\MTM^{\integral}(\cnum_{\tau}\times X,\real)^{\alg}_{\star}.
\]

\subsubsection{$\cnum^{\ast}$-homogeneity}

For $\nbigrtilde^{\alg}_X$-modules,
the $\cnum^{\ast}$-homogeneity is naturally defined
as in the case of $\nbigrtilde_X$-modules.
For $\vecn\in\seisuu^2_{\geq 0}$
such that $\gcd(n_1,n_2)=1$ and $n_1\geq n_2$,
let $\rho_{\vecn}$ denote 
the $\cnum^{\ast}$-action on
$\cnum_{\lambda}\times\cnum_t\times X$
by $a(\lambda,t,x)=(a^{n_1}\lambda,a^{n_2}t,x)$.
Let $\MTM^{\integral}_{\vecn}(\cnum\times X)^{\alg}$
denote the full subcategory of 
$(\nbigt,W)=((\nbigm',\nbigm'',C),W)\in
\MTM^{\integral}(\cnum\times X)^{\alg}$
such that
$\nbigm'$ and $\nbigm''$
are homogeneous with respect to $\rho_{\vecn}$.
Clearly,
$\FT_{\pm}$ induce the following equivalences:
\[
 \FT_{\pm}:
  \MTM^{\integral}_{\vecn}(\cnum_t\times X)^{\alg}
 \simeq
 \MTM^{\integral}_{\FT(\vecn)}(\cnum_{\tau}\times X)^{\alg},
 \quad
 \MTM^{\integral}_{\vecn}(\cnum_t\times X)_{\star}^{\alg}
 \simeq
 \MTM^{\integral}_{\FT(\vecn)}(\cnum_{\tau}\times X,[\star 0])^{\alg}.
\]

Let $F:X\lrarr Y$ be an algebraic morphism of
quasi-projective manifolds.
For $\star=!,\ast$, we obtain
\[
 \lefttop{T}(\id\times F)_{\star}:
 D^b\bigl(
 \MTM_{\vecn}(\cnum\times X)^{\alg}
 \bigr)
 \lrarr
 D^b\bigl(
 \MTM_{\vecn}(\cnum\times Y)^{\alg}
 \bigr),
\]
\[
 \lefttop{T}(\id\times F)_{\star}:
 D^b\bigl(
 \MTM_{\vecn}(\cnum\times X,[\star 0])^{\alg}
 \bigr)
 \lrarr
 D^b\bigl(
 \MTM_{\vecn}(\cnum\times Y,[\star 0])^{\alg}
 \bigr),
\]
\[
 \lefttop{T}(\id\times F)_{\star}:
 D^b\bigl(
 \MTM_{\vecn}(\cnum\times X)^{\alg}_{\star}
 \bigr)
 \lrarr
 D^b\bigl(
 \MTM_{\vecn}(\cnum\times Y)^{\alg}_{\star}
 \bigr).
\]
We also obtain
\[
 \lefttop{T}(\id\times F)^{\star}:
 D^b\bigl(
 \MTM_{\vecn}(\cnum\times Y)^{\alg}
 \bigr)
 \lrarr
 D^b\bigl(
 \MTM_{\vecn}(\cnum\times X)^{\alg}
 \bigr).
\]
For $(\star_1,\star_2)=(!,\ast),(\ast,!)$,
we have
\[
 \lefttop{T}(\id\times F)^{\star_1}:
 D^b\bigl(
 \MTM_{\vecn}(\cnum\times Y,[\star_2 0])^{\alg}
 \bigr)
 \lrarr
 D^b\bigl(
 \MTM_{\vecn}(\cnum\times X,[\star_2 0])^{\alg}
 \bigr),
\]
\[
 \lefttop{T}(\id\times F)^{\star_1}:
 D^b\bigl(
 \MTM_{\vecn}(\cnum\times Y)^{\alg}_{\star_2}
 \bigr)
 \lrarr
 D^b\bigl(
 \MTM_{\vecn}(\cnum\times X)^{\alg}_{\star_2}
 \bigr).
\]

\begin{lem}
Let $\nbigt\in\MTM^{\integral}_{\vecn}(\cnum\times X)^{\alg}$.
\begin{itemize}
 \item For any hypersurface $H$ of $X$,
$\nbigt[\star (\proj^1\times H)]$ $(\star=!,\ast)$
are also objects of
$\MTM^{\integral}_{\vecn}(\cnum\times X)^{\alg}$.
 \item Let $f$ and $\ftilde$ be
       as in Proposition {\rm\ref{prop;21.7.16.22}}.
       Then, $\Pi^{a,b}_{\ftilde,\star}(\nbigt)$
       and
       $\Pi^{a,b}_{\ftilde,\ast!}(\nbigt)$
       are objects of
       $\MTM^{\integral}_{\vecn}(\cnum\times X)^{\alg}$.
       As a result,
       $\Xi^{(a)}_{\ftilde}(\nbigt)$,
       $\psi^{(a)}_{\ftilde}(\nbigt)$
       and $\phi^{(0)}_{\ftilde}(\nbigt)$
       are objects of
       $\MTM^{\integral}_{\vecn}(\cnum\times X)^{\alg}$.
       \hfill\qed
\end{itemize}
\end{lem}

\begin{lem}
For $\nbigt\in\MTM^{\integral}_{\vecn}(\cnum\times X)^{\alg}$,
$j^{\ast}(\nbigt)$
and $\gammatilde^{\ast}(\nbigt)$
are objects of
$\MTM^{\integral}_{\vecn}(\cnum\times X)^{\alg}$.
Let $\star=!,\ast$.
If $\nbigt\in\MTM_{\vecn}(\cnum\times X,[\star 0])^{\alg}$.
$j^{\ast}(\nbigt)$
and $\gammatilde^{\ast}(\nbigt)$
are objects of
$\MTM^{\integral}_{\vecn}(\cnum\times,[\star 0])^{\alg}$.
If $\nbigt\in\MTM^{\integral}_{\vecn}(\cnum\times X)^{\alg}_{\star}$,
$j^{\ast}(\nbigt)$ and $\gammatilde^{\ast}(\nbigt)$ are
objects of
$\MTM^{\integral}_{\vecn}(\cnum\times X)_{\star}^{\alg}$.
\hfill\qed
\end{lem}

\begin{lem}
For $\nbigt\in\MTM^{\integral}_{\vecn}(\cnum\times X)^{\alg}$,
$\nbigt^{\ast}$,
and $\DD(\nbigt)$
are objects of
$\MTM^{\integral}_{\vecn}(\cnum\times X)^{\alg}$.
Let $(\star_1,\star_2)=(!,\ast),(\ast,!)$.
If $\nbigt\in\MTM^{\integral}_{\vecn}(\cnum\times X,[\star_1 0])^{\alg}$,
$\nbigt^{\ast}$ and
$\DD\nbigt$ are objects of
$\MTM^{\integral}_{\vecn}(\cnum\times,[\star_20])^{\alg}$.
If $\nbigt\in\MTM^{\integral}_{\vecn}(\cnum\times X)_{\star_1}^{\alg}$,
$\nbigt^{\ast}$ and $\DD\nbigt$ are objects of
$\MTM^{\integral}_{\vecn}(\cnum\times X)_{\star_2}^{\alg}$.
\hfill\qed
\end{lem}

We set $X=X_1\times X_2$.
\begin{lem}
 For $\nbigt_i\in
  D^b(\MTM^{\integral}_{\vecn}(\cnum\times X_i)^{\alg})$,
$\nbigt_1\langle\boxtimes,\otimes\rangle^{\star}
 \nbigt_2$
and 
$\nbigt_1\langle\boxtimes,\ast\rangle_{\star}
 \nbigt_2$
 are objects of
 $D^b\bigl(
 \MTM^{\integral}_{\vecn}(\cnum\times X)^{\alg}
 \bigr)$.
\hfill\qed
\end{lem}

\begin{lem}
For any
$\nbigt_i\in \MTM^{\integral}_{(1,1)}(\cnum\times X_i,[\ast 0])^{\alg}$,
$\nbigt_1\langle\boxtimes,\otimes\rangle^!\nbigt_2[1]$
is an object of
$\MTM^{\integral}_{(1,1)}(\cnum\times X,[\ast 0])^{\alg}$.
For any
$\nbigt_i\in \MTM^{\integral}_{(1,1)}(\cnum\times X_i,[!0])^{\alg}$,
$\nbigt_1\langle\boxtimes,\otimes\rangle^{\ast}\nbigt_2[-1]$
is an object of
$\MTM^{\integral}_{(1,1)}(\cnum\times X,[! 0])^{\alg}$.
\hfill\qed
\end{lem}

Let
$\MTM^{\integral}_{\vecn}(\cnum\times X,\real)^{\alg}
\subset
\MTM^{\integral}(\cnum\times X,\real)^{\alg}$
denote the full subcategory of
the objects
$(\nbigt,W,\kappa)$
such that $(\nbigt,W)\in\MTM^{\integral}_{\vecn}(\cnum\times X)^{\alg}$.
We obtain the full subcategories
$\MTM^{\integral}_{\vecn}(\cnum\times X,[\star 0],\real)$
and 
$\MTM^{\integral}_{\vecn}(\cnum\times X,\real)_{\star}$
as in \S\ref{subsection;21.7.16.30}.
We obtain the following equivalences:
\[
 \MTM^{\integral}_{\vecn}(\cnum_t\times X,\real)^{\alg}
 \simeq
 \MTM^{\integral}_{\FT(\vecn)}(\cnum_{\tau}\times X,\real)^{\alg},
\]
\[
\MTM^{\integral}_{\vecn}(\cnum_t\times X,\real)_{\star}^{\alg}
 \simeq
 \MTM^{\integral}_{\FT(\vecn)}(\cnum_{\tau}\times X,[\star 0],\real)^{\alg}.
\]

\subsubsection{Rescalable mixed twistor $\nbigd$-modules}

Let $\nbigt\in\MTM^{\integral}_{(1,1)}(\cnum\times X,\real)^{\alg}$.
Let $\iota_1:\{1\}\times X\lrarr \proj^1\times X$
denote the inclusion.
There exists 
$\ttU(\nbigt)\in\MTM^{\integral}(X,\real)^{\alg}$
such that
$\iota_{1\dagger}\ttU(\nbigt)
=\psi^{(1)}_{t-1}(\nbigt)$.
We obtain the functor
$\ttU:\MTM^{\integral}_{(1,1)}(\cnum\times X,\real)^{\alg}
\lrarr\MTM^{\integral}(X,\real)^{\alg}$.
As the restriction to full subcategories,
we obtain 
\begin{equation}
\label{eq;21.7.16.40}
 \ttU:
\MTM^{\integral}_{(1,1)}(\cnum\times X,[\star 0],\real)^{\alg}
\lrarr
\MTM^{\integral}(X,\real)^{\alg}.
\end{equation}
We obtain the following proposition from
Proposition \ref{prop;21.7.12.30}.
\begin{prop}
\label{prop;21.7.16.41}
The functor {\rm(\ref{eq;21.7.16.40})} is
fully faithful.
\hfill\qed
\end{prop}

The following definition is similar to
Definition \ref{df;21.7.16.50}.
\begin{df}
Let
$\MTM^{\integral}_{\res}(X,\real)^{\alg}$
denote the essential image of
the functor {\rm(\ref{eq;21.7.16.40})},
which is independent of $\star=!,\ast$.
Objects of $\MTM^{\integral}_{\res}(X,\real)^{\alg}$
are called rescalable mixed twistor $\nbigd$-modules
on $X^{\alg}$.
\hfill\qed
\end{df}

We state some functorial properties
as in \S\ref{subsection;21.7.16.51}.
\begin{prop}
Let $F:X\lrarr Y$ be an algebraic morphism of
quasi-projective manifolds.
For 
$\nbigt^{\bullet}\in
D^b(\MTM^{\integral}_{(1,1)}(\cnum\times X,[\star 0],\real)^{\alg})$,
there exist natural isomorphisms
\[
 F_{\star}\bigl(
 \ttU(\nbigt^{\bullet})
 \bigr)
 \simeq
 \ttU\bigl(
 \bigl(
 (\id\times F)_{\star}\nbigt^{\bullet}
\bigr)
 \bigr)
\]
in $D^b(\MTM^{\integral}(Y,\real)^{\alg})$.
In particular, we obtain
$F_{\star}:
D^b(\MTM^{\integral}_{\res}(X,\real)^{\alg})
\lrarr
D^b(\MTM^{\integral}_{\res}(Y,\real)^{\alg})$.
\hfill\qed
\end{prop}

\begin{prop}
Let $H$ be any hypersurface of $X$.
For any
$\nbigt\in\MTM^{\integral}_{(1,1)}(\cnum\times X,[\star 0],\real)^{\integral}$,
we have
\[
 \ttU
 \bigl(
  (\nbigt[\star' (\proj^1\times H)])
 \bigr)
=\ttU
 \bigl(
  (\nbigt[\star' (\proj^1\times H)])[\star(0\times X)]
 \bigr)
 \simeq
  \ttU(\nbigt)[\star'H]
 \quad (\star'=!,\ast).
\] 
In particular,
for any $\nbigt_1\in\MTM^{\integral}_{\res}(X,\real)^{\alg}$,
$\nbigt_1[\star'H]$ $(\star'=!,\ast)$
are objects of 
$\MTM^{\integral}_{\res}(X,\real)^{\alg}$.
\hfill\qed
\end{prop}

\begin{prop}
Let $f$ and $\ftilde$ be as in
as in Proposition {\rm\ref{prop;21.7.16.22}}.
Let $\nbigt\in
\MTM^{\integral}_{(1,1)}(\cnum\times X,[\star 0],\real)^{\alg}$.
Then, there exist natural isomorphisms
\[
 \ttU
 \bigl(
 \Pi^{a,b}_{\ftilde,\star'}(\nbigt)[\star(0\times X)]
 \bigr)
 \simeq
 \ttU
 \bigl(
 \Pi^{a,b}_{\ftilde,\star'}(\nbigt)
 \bigr)
 \simeq
 \Pi^{a,b}_{f,\star'}
 \ttU(\nbigt),
\]
\[
 \ttU
 \bigl(
 \Pi^{a,b}_{\ftilde,\ast !}(\nbigt)[\star(0\times X)]
 \bigr)
 \simeq
 \ttU
 \bigl(
 \Pi^{a,b}_{\ftilde,\ast !}(\nbigt)
 \bigr)
 \simeq
 \Pi^{a,b}_{f,\ast !}
 \ttU(\nbigt).
\] 
As a result, 
there exist natural isomorphisms
\[
 \ttU\bigl(
 \Xi^{(a)}_{\ftilde}(\nbigt)[\star(0\times X)]
 \bigr)
 \simeq
 \ttU\bigl(
 \Xi^{(a)}_{\ftilde}(\nbigt)
 \bigr)
  \simeq
 \ttU\bigl(
 \Xi^{(a)}_{\ftilde}(\nbigt)
 \bigr)
 \simeq
 \Xi^{(a)}_{f}
 \ttU(\nbigt),
\]
\[
 \ttU\bigl(
 \psi^{(a)}_{\ftilde}(\nbigt)[\star(0\times X)]
 \bigr)
 \simeq
 \ttU\bigl(
 \psi^{(a)}_{\ftilde}(\nbigt)
 \bigr)
 \simeq
 \ttU\bigl(
 \psi^{(a)}_{\ftilde}(\nbigt)
 \bigr)
 \simeq
 \psi^{(a)}_{f}
 \ttU(\nbigt),
\]
\[
 \ttU\bigl(
 \phi^{(0)}_{\ftilde}(\nbigt)[\star(0\times X)]
 \bigr)
 \simeq
 \ttU\bigl(
 \phi^{(0)}_{\ftilde}(\nbigt)
 \bigr) 
 \simeq
 \phi^{(0)}_{f}
 \ttU(\nbigt).
\]
Therefore,
for any $\nbigt_1\in\MTM^{\integral}_{\res}(X,\real)^{\alg}$,
$\Pi^{a,b}_{f,\star}(\nbigt_1)$,
$\Pi^{a,b}_{f,\ast !}(\nbigt_1)$,
$\Xi^{(a)}_f(\nbigt_1)$,
$\psi^{(a)}_f(\nbigt_1)$
and $\phi^{(0)}_f(\nbigt)$
are objects of $\MTM^{\integral}_{\res}(X,\real)^{\alg}$.
\hfill\qed
\end{prop}

\begin{prop}
Let $\star=!,\ast$.
For $\nbigt\in
\MTM^{\integral}_{(1,1)}(\cnum\times X,[\star 0],\real)^{\alg}$,
there exist isomorphisms
\[
 \ttU
 j^{\ast}(\nbigt)
 \simeq
 j^{\ast}
 \ttU(\nbigt),
 \quad
 \ttU
 \gammatilde^{\ast}(\nbigt)
 \simeq
 \gammatilde^{\ast}
 \ttU(\nbigt).
\]
Therefore,
for any $\nbigt_1\in\MTM^{\integral}_{\res}(X,\real)^{\alg}$,
$j^{\ast}\nbigt_1$
and $\gammatilde^{\ast}\nbigt_1$
are objects of $\MTM^{\integral}_{\res}(X,\real)^{\alg}$.
\hfill\qed
\end{prop}

\begin{prop}
Let $\star=!,\ast$.
For $\nbigt\in
\MTM^{\integral}_{(1,1)}(\cnum\times X,[\star 0],\real)^{\alg}$,
there exist isomorphisms
\[
 \bigl(
 \ttU(\nbigt)
 \bigr)^{\ast}
 \simeq
 \ttU
 (\nbigt^{\ast})\otimes\newTate(-1),
 \quad
 \DD\bigl(
 \ttU(\nbigt)
 \bigr)
 \simeq
 \ttU\bigl(
 \DD(\nbigt)
\bigr)
\otimes\newTate(-1).
\]
Therefore,
for any $\nbigt_1\in\MTM^{\integral}_{\res}(X,\real)^{\alg}$,
$\nbigt_1^{\ast}$
and $\DD\nbigt_1$
are objects of $\MTM^{\integral}_{\res}(X,\real)^{\alg}$.
\hfill\qed
\end{prop}

Let $X=X_1\times X_2$.
The following proposition is similar to
Proposition \ref{prop;21.7.16.60}.
\begin{prop}
For $\nbigt_i\in\MTM^{\integral}_{(1,1)}(\cnum\times X_i)^{\alg}$,
there exist the following natural isomorphisms:
\[
 \ttU\Bigl(
 \bigl(
 \nbigt_1\langle\boxtimes,\otimes\rangle^!\nbigt_2
 \bigr)[1]
 \Bigr)
 \simeq
 \bigl(
 \ttU(\nbigt_1)\boxtimes\ttU(\nbigt_2)
 \bigr)
 \otimes\newTate(-1),
\quad\quad
 \ttU\Bigl(
 \bigl(
 \nbigt_1\langle\boxtimes,\otimes\rangle^{\ast}\nbigt_2
 \bigr)[-1]
  \Bigr)
 \simeq
 \ttU(\nbigt_1)\boxtimes\ttU(\nbigt_2).
\]
As a result,
for any $\nbigt_i\in\MTM^{\integral}_{\res}(X_i,\real)^{\alg}$,
$\nbigt_1\boxtimes\nbigt_2$
is an object of 
$\MTM^{\integral}_{\res}(X,\real)^{\alg}$.
\hfill\qed
\end{prop}

\subsubsection{Comparison with exponential mixed $\real$-Hodge modules}

Let $\MHM(\cnum\times X,\real)^{\alg}$
denote the category of algebraic mixed $\real$-Hodge modules
on $(\cnum\times X)^{\alg}$.
As in \S\ref{subsection;21.7.16.71},
by the Rees construction,
we obtain an equivalence
\[
 \ttV:\MHM(\cnum\times X,\real)^{\alg}
 \lrarr
 \MTM^{\integral}_{(1,0)}(\cnum\times X,\real)^{\alg}.
\]
Let
$\MHM(\cnum\times X,\real)_{\star}^{\alg}
\subset
\MHM(\cnum\times X,\real)^{\alg}$
denote the full subcategory
corresponding to
$\MTM^{\integral}_{(1,0)}(\cnum\times X,\real)_{\star}^{\alg}$,
i.e.,
algebraic mixed $\real$-Hodge modules
$(M,F,P_{\real})$ on $\cnum\times X$
such that
$R\pi_{X\star}(P_{\real})=0$.
We have the equivalence
$\ttV:\MHM(\cnum\times X,\real)_{\star}^{\alg}
\simeq
\MTM^{\integral}_{(1,0)}(\cnum\times X,\real)_{\star}^{\alg}$.
As in Remark \ref{rem;21.7.16.70},
we call $\MHM(\cnum\times X,\real)_{\ast}^{\alg}$
the category of 
algebraic exponential mixed $\real$-Hodge modules on $X$.

We obtain the following functors:
\[
 \ttB_{\pm}:=\ttU\circ\FT_{\pm}\circ\ttV:
 \MHM(\cnum\times X,\real)^{\alg}
 \lrarr
 \MTM^{\integral}_{\res}(X,\real)^{\alg}.
\]
As in (\ref{eq;21.7.16.4}),
we have
\[
 \ttB_{\pm}(M,F,P_{\real},W)
=\pi_{X!}\bigl(
 \ttV(M,F,P_{\real},W)\otimes\nbigt^{\alg}_{\sm}(\pm t)
 \bigr)
 \simeq
 \pi_{X\ast}\bigl(
 \ttV(M,F,P_{\real},W)\otimes\nbigt^{\alg}_{\sm}(\pm t)
 \bigr).
\]
Let $\ttB_{\pm,\star}$
denote the restriction of $\ttB_{\pm}$
to 
$\MHM(\cnum\times X,\real)_{\star}^{\alg}$.
The following theorem is similar to Theorem \ref{thm;21.7.16.3}.
\begin{thm}
The functors $\ttB_{\pm,\star}$ are equivalent.
\hfill\qed
\end{thm}

We have the compatibility with $\ttB_{\pm}$
and the other basic functors,
which are easy to prove.
We just state such compatibility in the following.

\begin{prop}
Let $f:X\lrarr Y$ be an algebraic morphism of quasi-projective manifolds.
For any object
$(M,F,P_{\real},W)\in D^b(\MHM(\cnum\times X,\real)^{\alg})$,
there exist the following isomorphisms:
\[
 \ttB_{\pm}(\id\times f)_{\star}(M,F,P_{\real},W)
 \simeq
 f_{\star}\ttB_{\pm}(M,F,P_{\real},W).
\]
For any $(M,F,P_{\real},W)\in D^b(\MHM(\cnum\times Y,\real)^{\alg})$,
there exist the following isomorphisms:
\[
 \ttB_{\pm}(\id\times f)^{\star}(M,F,P_{\real},W)
 \simeq
 \lefttop{T}f^{\star}\ttB_{\pm}(M,F,P_{\real},W).
\] 
\hfill\qed
\end{prop}

\begin{prop}
Let $(M,F,P_{\real},W)\in\MHM(\cnum\times X,\real)^{\alg}$.
\begin{itemize}
 \item Let $H$ be a hypersurface of $X$.
       There exist natural isomorphisms
\[
       \ttB_{\pm}\Bigl(
       (M,F,P_{\real},W)[\star (\proj^1\times H)]
       \Bigr)
       \simeq
       \ttB_{\pm}(M,F,P_{\real},W)[\star H]. 
\]
 \item Let $f$ and $\ftilde$ be as in
       Proposition {\rm\ref{prop;21.7.16.22}}.
       There exist the following natural isomorphisms:
 \[
       \ttB_{\pm}
       \Pi^{a,b}_{\ftilde,\star}(M,F,P_{\real},W)
       \simeq
       \Pi^{a,b}_{f,\star}
       \ttB_{\pm}(M,F,P_{\real},W),
\]
\[
       \ttB_{\pm}
       \Pi^{a,b}_{\ftilde,\ast!}(M,F,P_{\real},W)
       \simeq
       \Pi^{a,b}_{f,\ast !}
       \ttB_{\pm}(M,F,P_{\real},W).
 \]
       As a result,
       there exist the following natural isomorphisms:
\[
       \ttB_{\pm}
        \Xi^{(a)}_{\ftilde}(M,F,P_{\real},W)
       \simeq
       \Xi^{(a)}_{f}\ttB_{\pm}(M,F,P_{\real},W),
\quad
       \ttB_{\pm}
       \psi^{(a)}_{\ftilde}(M,F,P_{\real},W)
       \simeq
       \psi^{(a)}_{f}\ttB_{\pm}(M,F,P_{\real},W),
\]
\[
       \ttB_{\pm}
       \phi^{(0)}_f(M,F,P_{\real},W)
       \simeq
       \phi^{(0)}_{\ftilde}\ttB_{\pm}(M,F,P_{\real},W).
\]
\hfill\qed
\end{itemize} 
\end{prop}

\begin{prop}
For $(M,F,P_{\real},W)\in\MHM(\cnum\times X,\real)^{\alg}$,
 there exist natural isomorphisms
\[
\ttB_{\pm}\bigl(
\DD(M,F,P_{\real},W)
\bigr)
\simeq
\DD(\ttB_{\pm}(M,F,P_{\real},W)).
\] 
For $(M_i,F,P_{i\real},W)\in \MHM(\cnum\times X_i)^{\alg}$,
there exist the following natural isomorphisms:
\[
 \ttB_{\pm}\bigl(
 (M_1,F,P_{1\real},W)
 \langle\boxtimes,\ast\rangle_{\star}
 (M_2,F,P_{2\real},W)
 \bigr)
 \simeq
 \ttB_{\pm}
 (M_1,F,P_{1\real},W)
 \boxtimes
 \ttB_{\pm}
 (M_2,F,P_{2\real},W).
\]
\hfill\qed 
\end{prop}

Let
$\ttP_{\star}:\MHM(\cnum\times X,\real)^{\alg}
\lrarr
\MHM(\cnum\times X,\real)_{\star}^{\alg}$
denote the functor
induced by
$\ttP_{\star}:\MTM^{\integral}_{(1,0)}(\cnum\times X,\real)
\lrarr
\MTM^{\integral}_{(1,0)}(\cnum\times X,\real)_{\star}$
and the equivalence $\ttV$.
For any $(M,F,P_{\real},W)\in \MHM(\cnum\times X,\real)_{\star}^{\alg}$,
we have the filtration $\Wtilde^{(2)}$
as in \S\ref{subsection;21.7.16.100},
i.e.,
$\Wtilde^{(2)}_j(M,F,P_{\real})
=\ttP_{\star}(W_j(M,F,P_{\real}))$.
The following holds.
\begin{prop}
For $(M,F,P_{\real},W)\in\MHM(\cnum\times X,\real)_{\star}^{\alg}$
such that
$\Gr^{\Wtilde^{(2)}}_j(M,F,P_{\real})=0$ unless $j=m$,
we have
$\Gr^{W}_j\ttB_{\pm,\star}(M,F,P_{\real},W)=0$ unless $j=m$.
\hfill\qed
\end{prop}

\noindent
{\em Address\\
Research Institute for Mathematical Sciences,
Kyoto University,
Kyoto 606-8502, Japan\\
takuro@kurims.kyoto-u.ac.jp
}
	
\end{document}

%% file: notation.tex
\newcommand{\nbiga}{\mathcal{A}}
\newcommand{\nbigb}{\mathcal{B}}
\newcommand{\nbigc}{\mathcal{C}}
\newcommand{\nbigd}{\mathcal{D}}
\newcommand{\nbige}{\mathcal{E}}
\newcommand{\nbigf}{\mathcal{F}}
\newcommand{\nbigg}{\mathcal{G}}
\newcommand{\nbigh}{\mathcal{H}}
\newcommand{\nbigi}{\mathcal{I}}
\newcommand{\nbigj}{\mathcal{J}}

\newcommand{\nbigl}{\mathcal{L}}
\newcommand{\nbigm}{\mathcal{M}}
\newcommand{\nbign}{\mathcal{N}}
\newcommand{\nbigo}{\mathcal{O}}
\newcommand{\nbigp}{\mathcal{P}}
\newcommand{\nbigq}{\mathcal{Q}}
\newcommand{\nbigr}{\mathcal{R}}
\newcommand{\nbigs}{\mathcal{S}}
\newcommand{\nbigt}{\mathcal{T}}
\newcommand{\nbigu}{\mathcal{U}}
\newcommand{\nbigv}{\mathcal{V}}
\newcommand{\nbigw}{\mathcal{W}}
\newcommand{\nbigx}{\mathcal{X}}
\newcommand{\nbigy}{\mathcal{Y}}
\newcommand{\nbigz}{\mathcal{Z}}

\newcommand{\proj}{\mathbb{P}}
\newcommand{\seisuu}{{\mathbb Z}}

\newcommand{\cnum}{{\mathbb C}}
\newcommand{\real}{{\mathbb R}}

\newcommand{\DD}{\mathbb{D}}

\newcommand{\II}{\mathbb{I}}

\newcommand{\gbiga}{\mathfrak A}
\newcommand{\gbigb}{\mathfrak B}
\newcommand{\gbigc}{\mathfrak C}
\newcommand{\gbigd}{\mathfrak D}

\newcommand{\gbigf}{\mathfrak F}
\newcommand{\gbigg}{\mathfrak G}
\newcommand{\gbigh}{\mathfrak H}

\newcommand{\gbigk}{\mathfrak K}
\newcommand{\gbigl}{\mathfrak L}
\newcommand{\gbigm}{\mathfrak M}
\newcommand{\gbign}{\mathfrak N}

\newcommand{\gbigr}{\mathfrak R}
\newcommand{\gbigs}{\mathfrak S}

\newcommand{\gbigu}{\mathfrak U}
\newcommand{\gbigv}{\mathfrak V}

\newcommand{\gbigx}{\mathfrak X}
\newcommand{\gbigy}{\mathfrak Y}
\newcommand{\gbigz}{\mathfrak Z}

\newcommand{\gminia}{\mathfrak a}
\newcommand{\gminib}{\mathfrak b}

\newcommand{\gminik}{\mathfrak k}

\newcommand{\vecv}{{\pmb v}}

\newcommand{\veca}{{\pmb a}}
\newcommand{\vecb}{{\pmb b}}

\newcommand{\vecdelta}{{\pmb \delta}}

\newcommand{\vecc}{{\pmb c}}

\newcommand{\vecm}{{\pmb m}}

\newcommand{\vecn}{{\pmb n}}
\newcommand{\vecp}{{\pmb p}}

\newcommand{\vecz}{{\pmb z}}

\newcommand{\larr}{\leftarrow}
\newcommand{\llarr}{\longleftarrow}

\newcommand{\rarr}{\rightarrow}
\newcommand{\lrarr}{\longrightarrow}

\newcommand{\pf}{{\bf Proof}\hspace{.1in}}
\newcommand{\qed}{\mbox{\rule{1.2mm}{3mm}}}

\def\Cok{\mathop{\rm Cok}\nolimits}

\def\Image{\mathop{\rm Im}\nolimits}

\def\Gr{\mathop{\rm Gr}\nolimits}

\def\Ker{\mathop{\rm Ker}\nolimits}

\def\Gr{\mathop{\rm Gr}\nolimits}
\def\Sym{\mathop{\rm Sym}\nolimits}

\def\Res{\mathop{\rm Res}\nolimits}

\def\ord{\mathop{\rm ord}\nolimits}

\def\can{\mathop{\rm can}\nolimits}
\def\var{\mathop{\rm var}\nolimits}
\def\id{\mathop{\rm id}\nolimits}

\def\gcd{\mathop{\rm g.c.d.}\nolimits}

\def\Supp{\mathop{\rm Supp}\nolimits}

\newcommand{\del}{\partial}
\newcommand{\delbar}{\overline{\del}}

\newcommand{\nhom}{{\mathcal Hom}}

\newcommand{\sankaku}{\triangle}

\newcommand{\Hbar}{\overline{H}}

\newcommand{\lefttop}[1]{{}^{#1}\!}

\def\irr{\mathop{\rm irr}\nolimits}
\def\reg{\mathop{\rm reg}\nolimits}

\newcommand{\deldel}{\eth}

\newcommand{\rhotilde}{\widetilde{\rho}}

\newcommand{\nablatilde}{\widetilde{\nabla}}

\newcommand{\nbigetilde}{\widetilde{\nbige}}

\newcommand{\vtilde}{\widetilde{v}}

\newcommand{\ftilde}{\widetilde{f}}
\newcommand{\utilde}{\widetilde{u}}
\newcommand{\Ftilde}{\widetilde{F}}

\newcommand{\stilde}{\widetilde{s}}

\newcommand{\Rtilde}{\widetilde{R}}

\newcommand{\Itilde}{\widetilde{I}}

\newcommand{\nbigvtilde}{\widetilde{\nbigv}}

\newcommand{\ellsitabar}{\underline{\ell}}
\newcommand{\Utilde}{\widetilde{U}}

\newcommand{\Xtilde}{\widetilde{X}}
\newcommand{\vecy}{\pmb y}

\newcommand{\nbigrtilde}{\widetilde{\nbigr}}

\newcommand{\nbigmtilde}{\widetilde{\nbigm}}

\def\ord{\mathop{\rm ord}\nolimits}

\def\MTM{\mathop{\rm MTM}\nolimits}
\def\MHM{\mathop{\rm MHM}\nolimits}

\def\DR{\mathop{\rm DR}\nolimits}

\def\Ch{\mathop{\rm Ch}\nolimits}

\def\res{\mathop{\rm res}\nolimits}

\def\integral{\mathop{\rm int}\nolimits}
\def\pt{\mathop{\rm pt}\nolimits}

\def\sm{\mathop{\rm sm}\nolimits}

\def\DM{\mathop{\rm DM}\nolimits}
\def\alg{\mathop{\rm alg}\nolimits}
\def\an{\mathop{\rm an}\nolimits}
\def\Tor{\mathop{\rm Tor}\nolimits}
\def\Malg{\mathop{\rm Malg}\nolimits}
\def\FT{\mathop{\rm FT}\nolimits}
\def\Hod{\mathop{\rm Hod}\nolimits}
\def\Dol{\mathop{\rm Dol}\nolimits}

\newcommand{\Dbar}{\overline{D}}

\newcommand{\Wtilde}{\widetilde{W}}
\newcommand{\Ptilde}{\widetilde{P}}

\newcommand{\gtilde}{\widetilde{g}}
\newcommand{\Ztilde}{\widetilde{Z}}
\newcommand{\nbigttilde}{\widetilde{\nbigt}}

\newcommand{\Ctilde}{\widetilde{C}}

\newcommand{\gbigctilde}{\widetilde{\gbigc}}

\newcommand{\Atilde}{\widetilde{A}}

\newcommand{\gbigrtilde}{\widetilde{\gbigr}}
\newcommand{\gbigmtilde}{\widetilde{\gbigm}}

\newcommand{\Ytilde}{\widetilde{Y}}
\newcommand{\Ybar}{\overline{Y}}

\newcommand{\ztilde}{\widetilde{z}}

\newcommand{\nbigxtilde}{\widetilde{\nbigx}}

\newcommand{\gbigvtilde}{\widetilde{\gbigv}}

\newcommand{\nbigitilde}{\widetilde{\nbigi}}

\newcommand{\vecnbigi}{{\pmb \nbigi}}

\newcommand{\Lambdatilde}{\widetilde{\Lambda}}

\newcommand{\nrhom}{R{\mathcal Hom}}
\newcommand{\DDD}{\pmb D}

\newcommand{\Hhat}{\widehat{H}}
\newcommand{\Htilde}{\widetilde{H}}

\newcommand{\Xbar}{\overline{X}}

\newcommand{\gammatilde}{\widetilde{\gamma}}
\newcommand{\Gammatilde}{\widetilde{\Gamma}}

\newcommand{\qtilde}{\widetilde{q}}

\newcommand{\nbigctilde}{\widetilde{\nbigc}}

\newcommand{\newTate}{\pmb{T}}

\newcommand{\vecnbigitilde}{\boldsymbol\nbigitilde}
\newcommand{\Qhat}{\widehat{Q}}

\newcommand{\Mhat}{\widehat{M}}
\newcommand{\Uhat}{\widehat{U}}
\newcommand{\gbigxtilde}{\widetilde{\gbigx}}
\newcommand{\gbightilde}{\widetilde{\gbigh}}
\newcommand{\zeroellsitabar}{0\sqcup\underline{\ell}}
\newcommand{\Jtilde}{\widetilde{J}}
\newcommand{\gbigytilde}{\widetilde{\gbigy}}

\newcommand{\gbigktilde}{\widetilde{\gbigk}}
\newcommand{\gbigmlambda}{\gbigm^{\lambda}}
\newcommand{\gbigxbar}{\overline{\gbigx}}

\newcommand{\MTMtilde}{\widetilde{\MTM}}
\newcommand{\MHMtilde}{\widetilde{\MHM}}
\newcommand{\ttP}{\texttt{P}}
\newcommand{\ttQ}{\texttt{Q}}
\newcommand{\ttR}{\texttt{R}}
\newcommand{\ttA}{\texttt{A}}
\newcommand{\ttU}{\texttt{U}}
\newcommand{\ttV}{\texttt{V}}
\newcommand{\ttB}{\texttt{B}}
\newcommand{\ttS}{\texttt{S}}
\newcommand{\Zbar}{\overline{Z}}

%% file: new_theorem.tex
\newtheorem{thm}{Theorem}[section]
\newtheorem{cor}[thm]{Corollary}

\newtheorem{rem}[thm]{Remark}
\newtheorem{lem}[thm]{Lemma}
\newtheorem{prop}[thm]{Proposition}
\newtheorem{df}[thm]{Definition}
\newtheorem{example}[thm]{Example}
\newtheorem{condition}[thm]{Condition}

\newtheorem{notation}[thm]{Notation}

%% file: 1-19.bbl
\begin{thebibliography}{99}

\bibitem{beilinson2}
A. Beilinson,
{\em How to glue perverse sheaves},
in; {\em $K$-theory, arithmetic and geometry
 (Moscow, 1984--1986)}, 
Lecture Notes in Math., {\bf 1289},
Springer, Berlin, (1987),
42--51.

\bibitem{bbd}
A. Beilinson, J. Bernstein, P. Deligne,
{\em Faisceaux pervers},
Analysis and topology on singular spaces, I (Luminy, 1981), 
Ast\'{e}risque, {\bf 100},
(1982), 5--171.
	
\bibitem{CDS}	
	A. Casta\~{n}o Dom\'{\i}nguez, 
	C. Sevenheck,
	{\em Irregular Hodge filtration of
	some confluent hypergeometric systems},
	J. Inst. Math. Jussieu {\bf 20} (2021),
	627--668.
\bibitem{CDSR}
	A. Casta\~{n}o Dom\'{\i}nguez,
	T. Reichelt,
	C. Sevenheck,
	{\em Examples of hypergeometric twistor D-modules,}
	Algebra Number Theory {\bf 13} (2019), 1415--1442. 
	
\bibitem{Chen-Yu}
	K.-C. Chen,
	J.-D. Yu,
	{\em The K\"{u}nneth formula for
	the twisted de Rham and Higgs cohomologies,}
	SIGMA Symmetry Integrability Geom. Methods Appl. {\bf 14} (2018),
	Paper No. 055
	
\bibitem{Danilov}
	V. I. Danilov,
	{\em Cohomology of algebraic varieties},
	Algebraic geometry, II, 1--125, 255--262,
	Encyclopaedia Math. Sci., {\bf 35}, Springer, Berlin, 1996
\bibitem{Deligne-Malgrange-Ramis}
P. Deligne,
B. Malgrange,
J-P. Ramis,
{\em Singularit\'es Irr\'{e}guli\`{e}res},
Documents Math\'ematiques {\bf 5},
Soci\'et\'e Math\'ematique de France,
2007.
	
\bibitem{Esnault-Sabbah-Yu}
H. Esnault, 
C. Sabbah, 
J.-D. Yu,
(with an appendix by M. Saito),
{\em $E_1$-degeneration of the irregular Hodge filtration},
J. reine angew. Math. (2015), doi:10.1515/crelle-2014-0118.

 \bibitem{Fedorov}
	 R. Fedorov,
	 {\em Variations of Hodge structures for
	 hypergeometric differential operators
	 and parabolic Higgs bundles,}
	 International Mathematics Research Notices,
	 Volume 2018, Issue 18 (2018), 5583--5608. 
	
\bibitem{Fresan-Sabbah-Yu}
J. Fres\'{a}n, C. Sabbah, J.-D. Yu,
{\em Hodge theory of Kloosterman connections},
Duke Math. J. {\bf 171} (2022), 1649--1747.
	
\bibitem{Hartshorne}	
R. Hartshorne,
{\em Algebraic geometry},
	Graduate Texts in Mathematics, {\bf 52}.
	Springer-Verlag, New York-Heidelberg, 1977.

 \bibitem{Hertling}
	 C, Hertling,
	 {\em $tt^{\ast}$ geometry, Frobenius manifolds,
	 their connections, and the construction for singularities},
	 J. Reine Angew. Math. {\bf 555} (2003), 77--161. 
	
\bibitem{Hertling-Sevenheck}
	C. Hertling,
	C. Sevenheck,
	{\em Nilpotent orbits of a generalization of Hodge structures.}
	J. Reine Angew. Math. {\bf 609} (2007), 23--80.

\bibitem{Hertling-Sevenheck2}	
	C. Hertling,
	C. Sevenheck,
	{\em Limits of families of Brieskorn lattices and
	compactified classifying spaces},
	Adv. Math. {\bf 223} (2010), 1155--1224.
	
\bibitem{Hotta-Takeuchi-Tanisaki}
	R. Hotta, K. Takeuchi, T. Tanisaki,
	{\em D-modules, perverse sheaves, and representation theory},
	Progress in Mathematics, {\bf 236}.
	Birkh\"{a}user Boston, Inc., Boston, MA, 2008. 
	
 \bibitem{kashiwara_text}
M. Kashiwara,
{\em $\nbigd$-modules and microlocal calculus},
Translations of Mathematical Monographs, 217. 
Iwanami Series in Modern Mathematics,
American Mathematical Society, 
2003

\bibitem{Kashiwara-Schapira}
M. Kashiwara and P. Schapira,
{\em Sheaves on manifolds},
Springer-Verlag, Berlin, 1990

	 
\bibitem{Kontsevich-Soibelman-exponential-Hodge}
	M. Kontsevich,
	Y. Soibelman,
	{\em Cohomological Hall algebra,
	exponential Hodge structures and
	motivic Donaldson-Thomas invariants},
	Commun. Number Theory Phys. {\bf 5} (2011), 231--352. 
\bibitem{malgrange}
B. Malgrange,
{\em Connexions m\'eromorphies $2$, Le r\'eseau canonique},
 Invent. Math. {\bf 124}, (1996) 367--387.


 \bibitem{malgrange-holonomic-D-modules}
B. Malgrange,	 
	 {\em On irregular holonomic D-modules},
	 in {\'{E}l\'{e}ments de la th\'{e}orie des
	 syst\`{e}mes diff\'{e}rentiels g\'{e}om\'{e}triques},
	 391--410, S\'{e}min. Congr., {\bf 8},
	 Soc. Math. France, Paris, 2004.

\bibitem{Martin}
	N. Martin,
	{\em Middle multiplicative convolution and hypergeometric equations.}
	J. Singul. {\bf 23} (2021), 194--204.
\bibitem{Matsumura}
	H. Matsumura,
	{\em Commutative ring theory.}
	Translated from the Japanese by M. Reid.
	Cambridge Studies in Advanced Mathematics, {\bf 8}.
	Cambridge University Press, Cambridge, 1986.
\bibitem{Mebkhout}
Z. Mebkhout,
	{\em Le formalisme des six op\'{e}rations de Grothendieck
	pour les $\nbigd_X$-modules coh\'{e}rents},
	With supplementary material by the author and L. Narv\'{a}ez Macarro.
	Hermann, Paris, 1989.
\bibitem{Mebkhout-Sabbah}
	Z. Mebkhout, C. Sabbah,
	{\S III.4 $\nbigd$-modules et cycles \'{e}vanescents}
	in \cite{Mebkhout}, 201--239.

\bibitem{mochi2}
T. Mochizuki,
{\em Asymptotic behaviour of tame harmonic bundles
and an application to pure twistor $\nbigd$-modules I, II},
Mem. AMS. {\bf 185}, (2007).
	
\bibitem{Mochizuki-DM-lattice}
	T. Mochizuki,
	{\em On Deligne-Malgrange lattices, resolution of turning points
	and harmonic bundles},
	Ann. Inst. Fourier (Grenoble) {\bf 59} (2009),
	2819--2837.

\bibitem{Mochizuki-Stokes}
	T. Mochizuki,
	{\em The Stokes structure of a good meromorphic flat bundle},
	J. Inst. Math. Jussieu {\bf 10} (2011), 675--712.

\bibitem{mochi8}
T. Mochizuki,
{\em Asymptotic behaviour of variation of
pure polarized TERP structures},
Publ. Res. Inst. Math. Sci. {\bf 47} (2011), 419--534.
	
\bibitem{Mochizuki-wild}
T. Mochizuki,
{\em Wild harmonic bundles and 
 wild pure twistor $\nbigd$-modules},
Ast\'{e}risque {\bf 340}, (2011)
	
\bibitem{Mochizuki-Toda-lattice-II}
T. Mochizuki,
{\em Harmonic bundles and Toda lattices with opposite sign II},
	Comm. Math. Phys. {\bf 328} (2014), 1159--1198. 
	
\bibitem{Mochizuki-MTM}
T. Mochizuki,
{\em Mixed twistor $\nbigd$-modules},
Springer, 2015.

\bibitem{Mochizuki-Kontsevich-complexes}
	T. Mochizuki,
	{\em A twistor approach to the Kontsevich complexes},
	Manuscripta Math. {\bf 157} (2018), 193--231.

\bibitem{Mochizuki-subanalytic2}
	T. Mochizuki,
	{\em Curve test for enhanced ind-sheaves and holonomic D-modules, II},
	to appear in Ann. l'ENS.
\bibitem{Mochizuki-GKZ}
	T. Mochizuki,
	{\em Twistor property of GKZ-hypergeometric systems}
	 arXiv:1501.04146, version 2.

\bibitem{Patrikis-Taylor}
	S. Patrikis,
	R. Taylor,
	{\em Automorphy and irreducibility of
	some $\ell$-adic representations.}
	Compos. Math. {\bf 151} (2015), 207--229.
	
\bibitem{Sabbah-pure-twistor}
	 C. Sabbah, {\em Polarizable twistor $\nbigd$-modules},
	 Ast\'{e}risque {\bf 300} (2005).

\bibitem{Sabbah-wild-twistor}
	C. Sabbah,
	{\em Wild twistor $\nbigd$-modules},
	in {\em Algebraic analysis and around},
	293--353, Adv. Stud. Pure Math., {\bf 54},
	Math. Soc. Japan, Tokyo, 2009. 

 \bibitem{Sabbah-Fourier-LaplaceII}
	 C. Sabbah,
	 {\em Fourier-Laplace transform of a variation of
	 polarized complex Hodge structure, II},
	 in
	 {\em New developments in algebraic geometry,
	 integrable systems and mirror symmetry (RIMS, Kyoto, 2008)},
	 289--347,
	 Adv. Stud. Pure Math., {\bf 59},
	 Math. Soc. Japan, Tokyo, 2010. 
 \bibitem{Sabbah-irregular-Hodge}
	C. Sabbah,
	{\em Irregular Hodge theory},
	With the collaboration of Jeng-Daw Yu.
	M\'{e}m. Soc. Math. Fr. (N.S.) {\bf 156} (2018).

\bibitem{Sabbah-Yu}
C. Sabbah,
J.-D. Yu.
{\em On the irregular Hodge filtration
of exponentially twisted mixed Hodge modules},
Forum Math. Sigma {\bf 3} (2015), 71 pp. 

\bibitem{Sabbah-Yu2}
	C. Sabbah,
	J.-D. Yu,
	{\em Irregular Hodge numbers of
	confluent hypergeometric differential equations},
	\'{E}pijournal G\'{e}om. Alg\'{e}brique {\bf 3}
	(2019).
	
\bibitem{saito1}
M. Saito,
{\em Modules de Hodge polarisables},
Publ. RIMS., {\bf 24},
(1988), 849--995.

\bibitem{saito2}
M. Saito,
{\em Mixed Hodge modules},
Publ. RIMS., {\bf 26}, (1990),
221--333.

\bibitem{Takahiro-Saito}
	T. Saito,
	{\em A description of monodromic mixed Hodge modules},
	J. Reine Angew. Math. {\bf 786} (2022), 107--153.

\bibitem{Takahiro-Saito2}
	T. Saito,
	The Hodge filtrations of monodromic mixed Hodge modules
	and the irregular Hodge filtrations,
	arXiv:2204.13381
	
\bibitem{Serre-GAGA}
J. P. Serre,
{\em G\'{e}om\'{e}trie alg\'{e}brique et g\'{e}om\'{e}trie analytique},
Ann. Inst. Fourier (Grenoble) {\bf 6} (1956), 1--42.

\bibitem{s3}
C. Simpson,
{\it Mixed twistor structures},
math.AG/9705006.
	
\bibitem{Siu-extension}
Y. T. Siu,
{\em Techniques of extension of analytic objects},
 Vol. 8. Marcel Dekker, Inc., New York, 1974.
	
\bibitem{Wlodarczyk}
J. W{\l}odarczyk, 
{\em Resolution of singularities of analytic spaces},
Proceedings of G\"{o}kova Geometry-Topology Conference 2008, 
G\"{o}kova Geometry/Topology Conference (GGT), G\"{o}kova, (2009),
31--63.

\bibitem{Yu}
J.-D. Yu,
{\em Irregular Hodge filtration on twisted de Rham cohomology},
Manuscripta Math. {\bf 144}(1–2) (2014), 99--133.

	
\end{thebibliography}
